\documentclass[draft,tbtags]{crmm-l}
\makeindex

\newdimen\secnumwidth
\makeatletter
\renewcommand{\tocsection}[3]{%
\indentlabel{\@ifnotempty{#2}{\ignorespaces#1\space \hbox to\secnumwidth{#2.\hfil}\quad}}#3}
\makeatother

\usepackage{xspace}
\usepackage[utf8x]{inputenc}
\usepackage[T1]{fontenc}
\usepackage{lmodern}
\usepackage[american,frenchb]{babel}
\let\textenglish\emph
\newcommand{\textfrench}[1]{\foreignlanguage{frenchb}{#1}}
\makeatletter
\@namedef{subjclassname@2010}{Classification mathématique \textup{2010}}

\def\@part[#1]#2{%
  \ifnum \c@secnumdepth >-2\relax \refstepcounter{part}%
    \addcontentsline{toc}{part}{\partname.
        \protect\enspace\protect\noindent#1}%
  \else
    \addcontentsline{toc}{part}{#1}\fi
  \begingroup\centering
  \ifnum \c@secnumdepth >-2\relax
       {\fontsize{\@xviipt}{22}\bfseries
         \partname} \vskip 20\p@ \fi
  \fontsize{\@xxpt}{25}\bfseries
      #1\vfil\vfil\endgroup \newpage\thispagestyle{empty}}
\addto\captionsfrench{%
\renewcommand{\partname}{{\ifcase\value{part}\or Première\or
 Deuxième\or Troisième\or Quatrième\or Cinquième\or
 Sixième\or Septième\or Huitième\or Neuvième\or
 Dixième\or Onzième\or Douzième\or Treizième\or
 Quatorzième\or Quinzième\or Seizième\or Dix-septième\or
 Dix-huitième\or Dix-neuvième\or Vingtième\fi}\space partie}

}
\makeatother

\usepackage{upref}
\usepackage{mathtools}
\mathtoolsset{mathic}

\usepackage{amssymb}

\DeclareFontFamily{U}{mathx}{\hyphenchar\font45}
\DeclareFontShape{U}{mathx}{m}{n}{
      <5> <6> <7> <8> <9> <10>
      <10.95> <12> <14.4> <17.28> <20.74> <24.88>
      mathx10
      }{}
\DeclareSymbolFont{mathx}{U}{mathx}{m}{n}
\DeclareFontSubstitution{U}{mathx}{m}{n}
\DeclareMathAccent{\widehat}{0}{mathx}{"70}
\DeclareMathAccent{\widecheck}{0}{mathx}{"71}
\DeclareMathAccent{\widetilde}{0}{mathx}{"72}
\DeclareMathAccent{\widearrow}{0}{mathx}{"74}
\DeclareMathAccent{\wideparen}{0}{mathx}{"75}

\usepackage[abbrev,initials]{amsrefs}
\makeatletter
\renewcommand{\BibLabel}{\hfil\hyper@anchorstart{cite.\CurrenptBib}\relax\thebib\hyper@anchorend.}%
\renewcommand{\PrintNames@a}[4]{%
    \PrintSeries{\name}
        {#1}
        {}{ et \set@othername}
        {,}{ \set@othername}
        {}{ et \set@othername}
        {#2}{#4}{#3}%
}

\renewcommand{\eprintpages}[1]{%
    #1\IfEmptyBibField{eprint}{}{\IfEmptyBibField{journal}{ p.}{}}%
}

\BibSpec{collection.article}{%
    +{}  {\PrintAuthors}                {author}
    +{,} { \textit}                     {title}
    +{.} { }                            {part}
    +{:} { \textit}                     {subtitle}
    +{,} { \PrintContributions}         {contribution}
    +{,} { \PrintConference}            {conference}
    +{}  {\PrintBook}                   {book}
    +{,} { }                            {booktitle}
    +{,} { \PrintDateB}                 {date}
    +{,} { p.~}                         {pages}
    +{,} { }                            {status}
    +{,} { \PrintDOI}                   {doi}
    +{,} { disponible à \eprint}        {eprint}
    +{}  { \parenthesize}               {language}
    +{}  { \PrintTranslation}           {translation}
    +{;} { \PrintReprint}               {reprint}
    +{.} { }                            {note}
    +{.} {}                             {transition}
    +{}  {\SentenceSpace \PrintReviews} {review}
}
\BibSpec{nameLE}{
    +{}{}{initials}
    +{}{\IfEmptyBibField{initials}{}{ }\textsc}{surname}
    +{}{ }{jr}
}
\BibSpec{nameBE}{
    +{}{\textsc}{surname}
    +{}{ }{initials}
}
\BibSpec{nameinverted}{
    +{} {\textsc}  {surname}
    +{,}{ } {initials}
    +{,}{ } {jr}
}
\makeatother

\usepackage[newitem,newenum,neverdecrease]{paralist}
\makeatletter
\if@plflushright

\else

\fi
\renewcommand{\@asparaenum@}{%
\expandafter\list\csname label\@enumctr\endcsname{%
\usecounter{\@enumctr}%
\labelwidth\z@
\labelsep.5em
\leftmargin\z@
\parsep\parskip
\itemsep\z@
\topsep\listisep
\partopsep\parskip
\itemindent\parindent
\advance\itemindent\labelsep
\def\makelabel##1{\upshape ##1}}}
\def\@asparaitem@{%
\expandafter\list\csname\@itemitem\endcsname{%
\labelwidth\z@
\labelsep.5em
\leftmargin\z@
\parsep\parskip
\itemsep\z@
\topsep\listisep
\partopsep\parskip
\itemindent\parindent
\advance\itemindent\labelsep
\def\makelabel##1{##1}}}
\makeatother

\renewcommand{\thesection}{\thechapter.\arabic{section}}

\newtheorem{proposition}{Proposition}[section] 
\newtheorem{lemme}[proposition]{Lemme}
\newtheorem{theoreme}[proposition]{Théorème}
\newtheorem{corollaire}[proposition]{Corollaire}

\theoremstyle{remark}
\newtheorem{remarque}[proposition]{Remarque}

\numberwithin{equation}{proposition}
\renewcommand{\theequation}{\arabic{equation}}

\makeatletter
\renewenvironment{subequations}{%
  \refstepcounter{equation}%
  \protected@edef\theparentequation{\theequation}%
  \setcounter{parentequation}{\value{equation}}%
  \setcounter{equation}{0}%
  \def\theequation{\theparentequation)(\roman{equation}}%
  \ignorespaces
}{%
  \setcounter{equation}{\value{parentequation}}%
  \ignorespacesafterend
}
\makeatother

\DeclareMathOperator{\Sup}{Sup} 
\DeclareMathOperator{\Ad}{Ad}
\DeclareMathOperator{\Aut}{Aut}
\DeclareMathOperator{\End}{End}
\DeclareMathOperator{\GL}{GL}
\DeclareMathOperator{\Hom}{Hom}
\DeclareMathOperator{\Int}{Int}
\DeclareMathOperator{\Isom}{Isom}
\DeclareMathOperator{\Lie}{Lie}
\DeclareMathOperator{\Norm}{Norm}
\DeclareMathOperator{\Out}{Out}
\DeclareMathOperator{\Pol}{\mathbf{Pol}}
\DeclareMathOperator{\SL}{SL}
\DeclareMathOperator{\Support}{Support}
\DeclareMathOperator{\Vect}{Vect}
\DeclareMathOperator{\ad}{ad}
\DeclareMathOperator{\tr}{trace}
\DeclareMathOperator{\vol}{vol}

\newcommand{\AAA}{\mathbf{A}}
\newcommand{\AM}{\Math{A}}
\newcommand{\AOO}{A_{00}}
\newcommand{\AO}{A_0}
\newcommand{\CM}{\Math{C}}
\newcommand{\FF}{\mathcali{F}}
\newcommand{\FPQ}{F_{\PO}^{\Q}}
\newcommand{\GG}{\mathcali{G}}
\newcommand{\GM}{$(\G,\M)$}
\newcommand{\Gadef}{G(\adef)}
\newcommand{\G}{G}
\newcommand{\HG}{\HH_{\G}}
\newcommand{\HH}{\mathbf{H}}
\newcommand{\HO}{\HH_0}
\newcommand{\HPO}{\HH_{\PO}}
\newcommand{\HP}{\HH_P}
\newcommand{\HQ}{\HH_{\Q}}
\newcommand{\JJJ}{\mathbf{J}}
\newcommand{\JJ}{\td{J}}
\newcommand{\KG}{\K_{\G}}
\newcommand{\K}{\mathbf{K}}
\newcommand{\Levi}{\mathcali{L}}
\newcommand{\MM}{\mathcali{M}}
\newcommand{\MOO}{\M_{00}}
\newcommand{\MO}{\M_0}
\newcommand{\Mint}{\mathbf{M}}
\newcommand{\Ms}{\M_s}
\newcommand{\M}{M}
\newcommand{\NM}{\Math{N}}
\newcommand{\NO}{\N_0}
\newcommand{\NQp}{\N_{\Qp}}
\newcommand{\Nadef}{N(\adef)}
\newcommand{\N}{N}
\newcommand{\POO}{P_{00}}
\newcommand{\PO}{P_0}
\newcommand{\Parab}{\mathcali{P}}
\newcommand{\QM}{\Math{Q}}
\newcommand{\Qdo}{\Q_0}
\newcommand{\Qm}{\Q^{-}}
\newcommand{\QpQ}{\mathcali{Q}}
\newcommand{\Qp}{\Q^{+}}
\newcommand{\Qquad}[1]{\qquad\hbox{#1}\qquad}
\newcommand{\Quad}[1]{\quad\hbox{#1}\quad}
\newcommand{\Q}{Q}
\newcommand{\RM}{\Math{R}}
\newcommand{\Rac}{\mathcali{R}}
\newcommand{\Rm}{\R^{-}}
\newcommand{\R}{R}
\newcommand{\Sieg}{\boldsymbol{\goth{S}}}
\newcommand{\Spcusp}{S'_{\mathrm{cusp}}} 
\newcommand{\Sp}{S'}  
\newcommand{\TK}{T_0}
\newcommand{\TO}{S_0}
\newcommand{\Tsieg}{\T_{\G}}
\newcommand{\Tsigma}{\widetilde{\sigma}}
\newcommand{\T}{T}
\newcommand{\UU}{\mathcali{U}}
\newcommand{\U}{U}
\newcommand{\XGa}{\X_{\G}^a}
\newcommand{\XG}{\X_{\G}}
\newcommand{\XMQ}{\X_{\M_{\Q}}}
\newcommand{\XM}{\X_{\M}}
\newcommand{\XPG}{\X_{P,\G}}
\newcommand{\XP}{\X_P}
\newcommand{\XQG}{\X_{\Q,\G}}
\newcommand{\XQ}{\X_{\Q}}
\newcommand{\XSG}{\X_{S,\G}}
\newcommand{\XS}{\X_{S}}
\newcommand{\XX}{\mathcali{X}}
\newcommand{\XpQG}{\X_{\pQ,\G}}
\newcommand{\X}{\mathbf{X}}
\newcommand{\YG}{\mathbf{Y}_{\G}}
\newcommand{\YP}{\mathbf{Y}_P}
\newcommand{\YQdo}{\mathbf{Y}_{\Qdo}}
\newcommand{\YQ}{\mathbf{Y}_{\Q}}
\newcommand{\YY}{\mathcali{Y}}
\newcommand{\ZM}{\Math{Z}}
\newcommand{\aT}{\HO(a)-T}
\newcommand{\autom}{\mathcali{A}}
\newcommand{\base}{\mathcal{B}}
\newcommand{\cad}{c'est-à-dire\xspace}
\newcommand{\cbf}{\mathbf{c}}
\newcommand{\centdi}{\cent_{\delti}}
\newcommand{\centd}{\cent_{\delta}}
\newcommand{\cent}{I}
\newcommand{\cgh}{\goth{h}}
\newcommand{\ctc}{\mathcal{C}_{c}}
\newcommand{\cty}{\mathcal{C}^{\infty}}
\newcommand{\ctyc}{\mathcal{C}^{\infty}_{c}}
\newcommand{\dPO}{\dist_{\PO}}
\newcommand{\delti}{\delta_s}
\newcommand{\delto}{\delta_0}
\newcommand{\demisomO}{\demisom_0}
\newcommand{\demisomP}{\demisom_P}
\newcommand{\demisomQ}{\demisom_{\Q}}
\newcommand{\dist}{\mathbf{d}}
\newcommand{\ebf}{\mathbf{e}}
\newcommand{\extb}{g}
\newcommand{\extc}{h}
\newcommand{\ext}{f}
\newcommand{\fff}{f}
\newcommand{\frakg}{\goth{g}} 
\newcommand{\gAG}{\gA_{\G}}
\newcommand{\gAO}{\goth{A}_0}
\newcommand{\gAtG}{\gA_{\tG}}
\newcommand{\gA}{\goth{A}}
\newcommand{\gHM}{\gH_{\M}}
\newcommand{\gH}{\goth{H}}
\newcommand{\gO}{\goth{O}}
\newcommand{\gSP}{\gS_{\PO}}
\newcommand{\gStOm}{\gS_{t,\Omega}}
\newcommand{\gS}{\boldsymbol{\goth{S}}}
\newcommand{\gY}{\goth{Y}}
\newcommand{\gZ}{\goth{Z}}
\newcommand{\gaPQ}{\ga_P^{\Q}}
\newcommand{\gaP}{\ga_P}
\newcommand{\gao}{\goth{a}_0}
\newcommand{\gasns}{(\gasn)^{*}}
\newcommand{\gasn}{\ga_{\M}^\G}
\newcommand{\gatPO}{\ga_{\tPO}}
\newcommand{\gatP}{\ga_{\tP}}
\newcommand{\gats}{\ga_{\tnL}^{\tG}}
\newcommand{\ga}{\goth{a}}
\newcommand{\gb}{\goth{b}}
\newcommand{\gc}{\goth{c}}
\newcommand{\gk}{\goth{k}}
\newcommand{\gmomega}{\boldsymbol{\omega}}
\newcommand{\gnQ}{\goth{n}}
\newcommand{\gn}{\goth{n}}
\newcommand{\go}{\goth{o}}
\newcommand{\gu}{\goth{u}}
\newcommand{\hDelta}{\widehat\Delta}
\newcommand{\hff}{g}
\newcommand{\htautPQ}{\htau_{\tP}^{\tQ}}
\newcommand{\htautP}{\htau_{\tP}}
\newcommand{\htau}{\widehat\tau}
\newcommand{\hvf}{\widehat{\varphi}}
\newcommand{\ie}{c.-à-d.\xspace}
\newcommand{\ighMst}{\ima\gH_{\M}^{*}}
\newcommand{\ima}{\mathrm{i}}
\newcommand{\jtG}{j(\tG)}
\newcommand{\jtL}{j(\tL)}
\newcommand{\jtQ}{j(\tQ)}
\newcommand{\mG}{m}
\newcommand{\mm}{m}
\newcommand{\moins}{^{-1}}
\newcommand{\mtG}{m}
\newcommand{\nL}{L}
\newcommand{\oS}{S_0}
\newcommand{\orb}{\mathcali{O}}
\newcommand{\pQ}{\Q'}
\newcommand{\pni}{\par\noindent}
\newcommand{\ptf}{\;.}
\newcommand{\reg}{\boldsymbol{\rho}}
\newcommand{\sQR}{\sigma_{\Q}^{\R}}
\newcommand{\so}{s_0}
\newcommand{\tGM}{$(\tG,\tM)$}
\newcommand{\tGadef}{\tG(\adef)}
\newcommand{\tGamma}{\widetilde{\Gamma}}
\newcommand{\tG}{\td{\G}}
\newcommand{\tHO}{{\widetilde\HH_0}}
\newcommand{\tL}{\widetilde{L}}
\newcommand{\tMO}{\tM_0}
\newcommand{\tMQp}{\tM_{\Qp}}
\newcommand{\tM}{\widetilde{\M}}
\newcommand{\tPO}{\tP_0}
\newcommand{\tP}{\widetilde{P}}
\newcommand{\tQm}{{\tQ}^{-}}
\newcommand{\tQp}{{\tQ}^{+}}
\newcommand{\tQ}{\widetilde{\Q}}
\newcommand{\tRm}{{\tR}^{-}}
\newcommand{\tR}{\widetilde{\R}}
\newcommand{\talpha}{\widetilde{\alpha}}
\newcommand{\tbeta}{\widetilde{\beta}}
\newcommand{\tdM}{{}_tM}
\newcommand{\tdS}{{}_tS}
\newcommand{\tff}{\widetilde{\fff}}
\newcommand{\tgaQR}{\tga_{\Q}^{\R}}
\newcommand{\tga}{\widetilde{\ga}}
\newcommand{\thetomoins}{\theta_0\moins}
\newcommand{\theto}{\theta_0}
\newcommand{\tnL}{\widetilde{L}}
\newcommand{\tpi}{\td{\pi}}
\newcommand{\treg}{\widetilde{\reg}}
\newcommand{\tronc}{\boldsymbol{\Lambda}}
\newcommand{\tsQR}{\Tsigma_{\Q}^{\R}}
\newcommand{\tsuQuR}{\Tsigma_{\uQ}^{\uR}}
\newcommand{\ttmu}{\theto(\mu)}
\newcommand{\tus}{\widetilde{\us}}
\newcommand{\tvedQRt}{\tved(Q,R;t)}
\newcommand{\tvedQR}{\tved(Q,R)}
\newcommand{\tved}{\widetilde{\eta}}
\newcommand{\tve}{\widetilde{\ve}}
\newcommand{\tvpi}{\widetilde{\vpi}}
\newcommand{\tweyl}{\widetilde{\weyl}}
\newcommand{\uQ}{\Q_1}
\newcommand{\uR}{\R_1}
\newcommand{\uS}{S_1}
\newcommand{\us}{u}
\newcommand{\weyl}{\mathbf{W}}
\newcommand{\xT}{\HO(x)-T}
\newcommand{\xixT}{\HO(\xix)-T}
\newcommand{\xix}{\xi\,x}
\newcommand{\yy}{\mathbf{y}}
\newcommand{\zz}{\mathbf{z}}

\let\AdG=\Ad
\let\AdtG=\Ad
\let\Gpi=\pi
\let\LL=\Lambda
\let\Math=\mathbb
\let\Qbeta=\lambda
\let\WPhi=\Phi
\let\WPsi=\Psi

\let\adef=\ade
\let\bs=\backslash
\let\chm=\kappa
\let\demisom=\rho
\let\exta=\ext
\let\goth=\mathfrak 
\let\mathcali=\mathcal
\let\parag=\S
\let\tAd=\Ad
\let\td=\widetilde 
\let\tmu=\lambda
\let\ve=\varepsilon
\let\vf=\varphi
\let\vpi=\varpi

\newcommand{\RE}{\operatorname{Re}}

\newcommand{\chapintro}[2]{\subsubsection*{Chapitre \textup{#1}. #2}}
\newcommand{\partintro}[2]{\setcounter{part}{#1}\subsection*{\partname. #2}}
\newcommand{\newindex}[2]{\index{#1}}
\let\newnot=\newindex

\newcommand{\dd}{\mathop{}\!\mathrm{d}}
\newcommand{\ee}{\mathrm{e}}
\newcommand{\cf}{\emph{cf}.~}
\newcommand{\Cf}{\emph{Cf}.~}
\newcommand{\hyph}{\nobreakdash-\hspace{0pt}\relax}
\newcommand{\ldot}{\mathbin{.}}
\let\rest=\rvert
\newcommand{\nsubset}{\not\subset}
\newcommand{\nsupset}{\not\supset}
\newcommand{\nequiv}{\not\equiv}

\let\widebar=\overline

\addto\extrasfrench{\hyphenation{con-ver-gence con-nus cham-bre con-ver-gent}}

\addto\extrasamerican{%
\renewcommand{\chapintro}[2]{\subsubsection*{Chapter \textup{#1}. #2}}
\renewcommand{\partintro}[2]{\subsection*{Part #1. #2}}
}

\DeclareMathOperator{\meas}{meas}

\begin{document}
\frontmatter
\title{La formule des traces tordue\\
d'après le\\
\textenglish{Friday Morning Seminar}}

\author{Jean-Pierre Labesse}
\address{Institut de Mathématiques de Luminy, Campus de Luminy, Case 907, 13288 Marseille Cedex 9, France}
\email{labesse@iml.univ-mrs.fr}

\author{Jean-Loup Waldspurger}
\address{Institut de mathématiques de Jussieu, 2, place Jussieu,\newline 75005 Paris, France}
\email{waldspur@math.jussieu.fr}

\subjclass[2010]{primaire 11F72; secondaire 11R56, 20G35}
\maketitle

{\def\mathversion#1{}
\tableofcontents}

\chapter*{La formule des traces tordue: Avant-propos}

\begin{otherlanguage}{american}
The trace formula as an algebraic idea is as old as representation theory itself and can be
regarded as a form of the Frobenius reciprocity theorem. Suppose that $G$ is a finite group
and $\Gamma$ a subgroup. Consider the representation $\rho$ of $G$ on the functions on $\Gamma\backslash G$,
$g\colon \phi\mapsto \rho(g)\phi$, where $\rho(g)\phi(h)=\phi(hg)$. It is the representation
induced from the trivial representation of $\Gamma$ and the action of
a function on $G$ on the space of $\rho$
is
\[
f\colon \phi\mapsto\rho(f)\phi=\phi',\quad \phi'(h)=\frac{1}{\lvert G\rvert}\sum_G\phi(hg)f(g).
\]
The trace of $\rho(f)$ is readily calculated as
\[
\frac{1}{\lvert G\rvert}\sum_{\{(g,g_1)\,\vert\,g_1g=\gamma g_1,\gamma\in\Gamma\}}f(g)=
\frac{1}{\lvert\Gamma\rvert}\int_{\Gamma\backslash G}\sum_\Gamma f(g_1^{-1}\gamma g_1)\dd g_1,
\]
where the measure on $\Gamma\backslash G$ is implicitly normalized to $1$ by a factor $\lvert\Gamma\backslash G\rvert$
and $g_1$ runs over $\Gamma\backslash G$. The integral on the right can be expressed in various ways,
as a sum over conjugacy classes in $\Gamma$,
\[
\sum_{\{\gamma\}}\int_{\Gamma_\gamma\backslash G}f(g_1^{-1}\gamma g_1)\dd g_1,
\]
where $\Gamma_\gamma$ is the centralizer of $\gamma$ in $\Gamma$, or as
\begin{equation}
\sum_{\{\gamma\}}\meas(\Gamma_\gamma\backslash G_\gamma)\int_{G_\gamma\backslash G}f(g_1^{-1}\gamma g_1)\dd g_1.
\label{lang_eq1}
\end{equation}

We apply these simple formulas to the function $f(g)=\chi_\pi(g)$ given 
by the character $\chi_\pi$ of an irreducible representation $\pi$ of $G$. By the orthogonality
relations for matrix coefficients $\operatorname{tr}\bigl(\rho(f)\bigr)$
is the number of times $\pi$ is contained in $\rho$. The present function $f$ is a class
function, so that the formula \eqref{lang_eq1} for this trace reduces~to
\[
\sum_{\{\gamma\}}\meas(\Gamma_\gamma\backslash G)\chi_\pi(\gamma)=\frac{1}{\lvert\Gamma\rvert}\sum_\gamma\chi_\pi(\gamma),
\]
which is the number of times the trivial representation is contained in the restriction of $\pi$ to $\Gamma$.
In other words, we arrive at Frobenius reciprocity for the trivial representation of $\Gamma$.

One of the achievements of Selberg, in some regards perhaps the major achievement of
his career, was to recognize that there
was not only a formula similar to \eqref{lang_eq1} for discrete groups $\Gamma$ with compact quotient but also
one, the trace formula, for discrete subgroups of rank one. He, himself, was 
primarily concerned with subgroups of $\SL(2,\mathbb{R})$,
but the principles are similar for all groups of rank one for which the usual reduction theory is valid.
The spectral theory for the invariant differential operators is similar to the classical
spectral theory for a second-order differential equation on
a half-line: a one-dimensional continuous spectrum, empty if $\Gamma\backslash G$
is compact, together with a discrete spectrum. This is a classical theory
with which Selberg was more than familiar. If the quotient is not compact, the eigenfunctions
for the continuous spectrum
are constructed by the analytic continuation of Eisenstein series, a topic initiated by
Maa{\ss} and Roelcke, the central difficulty being resolved by Selberg. 

The rank-one theory together with the reduction theory for general arithmetic groups suggested
a more general theory, but there were difficulties, often misunderstood, even underestimated,
by commentators with limited familiarity with the methods used for their solution. Even the
general reduction theory developed in the nineteenth century by Eisenstein, H.~J.~S.~Smith,
Minkowski, Hermite and others and rescued, I am tempted to suggest, from oblivion by
C.~L.~Siegel in the twentieth, with the last proofs being provided by Borel and Harish-Chandra,
is sorely in need of a competent historical description. It is not surprising that,
in the fifties, when Siegel was still 
alive, still at the Institute for Advanced Study, Selberg, also at the Institute
and like Siegel an analytic
number theorist, was influenced not only by Maa{\ss} and Roelcke but also by Siegel.
Selberg wrote little and, in my experience, talked little about the sources of his
knowledge, so that it is difficult to understand why he made so little progress
with the theory in higher rank. The few notes he left suggest, although a final
judgement will have to await the closer examination of them by D.~Hejhal and others,
that in essence he made none and that he failed to understand subsequent developments.
Since he was certainly a strong mathematician, this is puzzling.

After years of unsystematic reflection I have concluded that the failure may have lain
in his lack of a clear understanding of the algebraic theory of semisimple groups
and, as a consequence, of the reduction theory of arithmetic groups. Specifically,
he failed to understand the notion of a cusp form, as it appeared in papers
by Godement and Harish-Chandra and in Gelfand's Stockholm lecture, and of the related 
decomposition of the spectrum according
to classes of parabolic subgroups. This is the clue to the general theory of
Eisenstein series. This lack may have been a consequence of an independent style
and a refusal, at least in his later years, of systematic study. I am not certain.

The general proofs of the analytic continuation of the
Eisenstein series and the description of the spectrum associated
to those of each type certainly incorporate basic ideas from those for rank-one groups
but they demand in addition not only a mastery of the theory of harmonic analysis
on reductive groups as created by Harish-Chandra but also the introduction of an
appropriate inductive structure and the solution of a number of specific
problems. 

Since the number of mathematicians with the necessary analytic experience
and the necessary understanding of the theory of semisimple groups is
limited, the general theory of Eisenstein series 
is not familiar to a large group
of mathematicians. Although this theory is necessary for the trace
formula for a general reductive group, even the first step toward the formula, the initial truncation of
the kernel of 
\begin{equation}
\rho(f)=\int f(g)\rho(g)\colon \phi\mapsto\phi',\quad \phi'(h)=\int\phi(hg)f(g)\dd g,\label{lang_eq2}
\end{equation}
for a smooth function with compact support, which is what permits the development of
the trace formula, was by no means obvious to me. I tried, unsuccessfully to
take it. It was taken, after much reflection\,---\,about two years I would suggest\,---\,by Arthur, 
who followed it by years of
effort and many papers, most of which neither I nor many other specialists have
yet digested. I attempted a summary some time ago in an article
\emph{The trace formula and its applications\textup: An introduction
to the work of James Arthur} that appeared in the Canadian Bulletin of Mathematics. 

In my view, although there is still much to do in the way of
integrating it with classical analytic number theory and classical
algebraic number theory, which will not be easy, the trace formula itself is the key to the construction
of a theory of automorphic forms in which functoriality and reciprocity appear
in their full generality. Reciprocity will demand, of course, more.

Fortunately the twisted trace formula, within or without endoscopy as the case requires, offers
an alternative, the proof
of whose consequences is less demanding from both an analytic and a number theoretic standpoint 
than those of the trace formula itself, but that offers none the less
substantial rewards. One of the earliest such applications was to cyclic base change for $\GL(2)$
and some special
cases of the Artin conjecture, which were later put to spectacular use by A. Wiles in the proof of
Fermat's theorem. More recently, Arthur has, in the book \emph{The endoscopic classification of representations\textup: orthogonal
and symplectic groups} systematically developed the 
twisted trace formula and twisted
endoscopy for $\GL(n)$ and applied it to establish important cases of functoriality. Its appearance will be,
in my view, 
an event of major importance in the theory of automorphic forms, and the book itself an opportunity 
for specialists, many of them mired in specific and limited techniques, to grasp the analytic
possibilities of the modern theory of automorphic forms.
\end{otherlanguage}

\bigskip
Selon son titre, c'est cette formule tordue qui est le sujet du présent livre de Jean-Pierre Labesse et
Jean-Loup Waldspurger. Pour la formule tordue, l'objet principal n'est pas l'opérateur $\rho(f)$
de la formule \eqref{lang_eq2}; il est plutôt, s'il est permis de m'exprimer d'une façon un peu plus simple
que celle du livre, $\rho(f)\circ\theta$, où $\theta$ est un automorphisme extérieur d'ordre fini du groupe
$G(\mathbb{A})$. Autant que je sache, c'est H.~Jacquet qui au début des années~80
avait proposé qu'une formule des traces pour de tels opérateurs pourrait résoudre quelques questions
qui l'intéressaient alors. J'avais essayé de trouver une telle formule, mais j'ai échoué. 

Toutefois, un peu plus tard, Y.~Flicker soumit aux \textenglish{Annals of Mathematics} un article qui contenait
une formule tordue convaincante.
Comme je voulais moi aussi le faire, il avait, mais avec apparemment 
plus de succès, trouvé une modification 
du noyau tronqué de Arthur. Malheureusement en examinant la démonstration qu'il avait proposée, je n'ai trouvé
qu'un argument sans fondement. Il n'était que le
simulacre d'une démonstration. Autant que je sache,
l'auteur n'a jamais réussi à proposer une autre démonstration
plus solide. Par la suite, même assez tôt, lors d'une année thématique
à l'\textenglish{Institute for Advanced Study} sur les formes automorphes, j'ai proposé à Jean-Pierre Labesse
que nous offrions avec quelques collègues un séminaire dans lequel nous essaierions de trouver une
démonstration de la formule donnée par Flicker. C'est ce séminaire qui est devenu le \textenglish{Morning
Seminar}, pour lequel les notes des divers conférenciers en autant qu'elles étaient disponibles
ont été distribuées aux auditeurs. Elles sont toujours disponibles en ligne. 

Ce fut une année très exigeante pour tous les participants et quelques conférenciers
n'ont pas réussi à rédiger des notes; d'autres les ont rédigées hâtivement et, pour la
plupart, nous ne sommes pas retournés réfléchir à ce que nous avions écrit lors du séminaire.
Là le sujet s'expose à des dangers car 
Arthur, en rédigeant son livre, avait suffisamment de pain sur la planche 
qu'il devait forcément tenir pour acquises les conclusions du \textenglish{Morning Seminar}. Même s'il n'y avait pas
d'autres raisons pour leur savoir gré, à cause de celle-ci seule nous devons beaucoup à
Jean-Pierre Labesse et à Jean-Loup Waldspurger d'avoir repris la formule tordue et de l'avoir assise 
sur des fondations sûres.

Ils observent d'ailleurs,
qu'ils se sont délestés de quelques arguments que j'avais esquissés dans mes notes
et qu'ils les ont remplacés par d'autres. Quoique j'aie beaucoup de confiance en le
mathématicien que j'étais, bien plus que dans celui que je suis, il est vrai que
vers la fin du séminaire j'étais éreinté et manquais de temps. De toutes façons
je n'ai pas eu depuis envie de reprendre ces arguments. Il est certainement
possible que j'avais trop simplifié les choses.

Les auteurs expliquent eux-mêmes dans la préface la genèse de leur entreprise et il n'y a pas
besoin de répéter ce qu'ils écrivent. Je voudrais toutefois ajouter un sentiment de reconnaissance
personelle, d'abord à Jean-Pierre Labesse et Jean-Loup Waldspurger et ensuite à 
l'\textenglish{American Mathematical Society} et surtout à son président Eric Friedlander d'avoir
accepté de publier ce texte en français, ce qui à ma surprise n'était pas si
évident. 

Il me semble que les mathématiciens se doivent de bien comprendre les conséquences
de leur utilisation toujours croissante et maintenant presque universelle de l'anglais,
ou plutôt de l'américain car nous sommes tous, comme mathématiciens, à la remorque
des américains. Cette pratique est accompagnée d'une ignorance d'autres langues de sorte
que les mathématiciens sont coupés à maints égards non seulement des mathématiques des
dix-septième et dix-huitième siècles mais aussi, et à un degré bien plus grave,
des mathématiques allemandes, françaises, russes ou italiennes du dix-neuvième siècle
et de la première moitié du vingtième siècle. On ne trouve pas facilement dans 
des textes 
contemporains tout
ce que nous avons hérité
des grands mathématiciens de ces périodes.
Ce déshéritement est peut-être moins grave pour les domaines qui n'ont apparu que récemment
qu'il l'est pour la théorie des formes automorphes qui est une fusion et une continuation
de sujets comme la théorie algébrique
des nombres, la géométrie algébrique, la théorie des groupes et leurs
représentations, ou la théorie
analytique des nombres, donc des sujets bien enracinés
dans l'histoire des mathématiques et où les façons de travailler ont beaucoup 
changé depuis le temps de, disons, Weierstrass et Jacobi ou Dedekind
et Frobenius et le présent. 

L'utilisation d'une langue unique et les moyens de communications et de voyager contemporains permettent
la formation de petites cellules de mathématiciens qui ne communiquent guère entre elles. C'est
une pratique que ceux et celles qui veulent contribuer d'une façon sérieuse à la théorie
des formes automorphes ne peuvent pas se permettre. Même ce livre ou le livre de Arthur ne sont,
à mon avis, au moins
en partie, qu'une préparation pour des travaux nettement plus difficiles. Il s'agit toutefois avec les deux
d'un apprentissage que l'on ne peut pas contourner.
D'ailleurs notre devoir est de résoudre les problèmes actuels, ou au moins d'essayer de résoudre ceux-ci,
et non pas ceux de l'avenir et, en plus, de ne pas s'adonner trop aux rêves. 
Il faut donc reconnaître lesquels sont abordables à présent et essayer de surmonter les difficultés
concrètes qu'ils posent. Le mérite de Jacquet, Labesse-\hspace{0pt}Waldspurger et Arthur,
au moins en ce qui concerne la formule tordue, est en partie
de l'avoir fait.

La formule des traces elle-même -- donc la formule introduite par Selberg pour les groupes
de rang $1$ mais dont
la forme générale a été créée par Arthur -- sera, à mon avis, dans sa forme
stable -- donc accompagnée par une théorie d'endoscopie, pour laquelle le lemme fondamental, démontré
dans une suite de travaux de divers mathématiciens dont le plus connu est celui de Ngô Bau Châu,
est une composante critique -- l'outil essentiel pour établir les théorèmes de base
de la théorie des formes ou représentations automorphes, mais seulement lorsque on aura réussi
à l'accompagner avec des formes nouvelles de la théorie analytique des nombres et
de la théorie des corps de classes, les deux fusionées comme elles étaient jadis --
par exemple, dans le rapport \emph{Klassenkörperbericht}
de Hasse -- mais portant aussi sur des extensions non abéliennes,
de sorte que la fusion n'est pas évidente. 
Il s'agit certainement de projets pour l'avenir. Pour l'instant,
pour convaincre des jeunes mathématiciens de la valeur des conjectures sur la fonctorialité
et la réciprocité dans le cadre des formes automorphes et de celle de la formule des traces
elle-même on a besoin -- sinon pour la formule stable au moins pour la théorie plus
générale de la formule tordue, stable ou non --
de moyens à portée de main et de théorèmes dont la
valeur est indiscutable. Ces moyens et ces théorèmes se trouveront dans deux livres.

Pour des théorèmes importants sur la fonctorialité qui
découlent de la formule tordue, on peut consulter le livre de Arthur; pour les fondations
d'un traitement
rigoureux et complet de la formule tordue, c'est ce livre de Labesse-\hspace{0pt}Waldspurger qu'il
faut consulter. Les deux maîtrisés, le jeune mathématicien peut, s'il le veut, entamer lui-même
les problèmes rattachés à la formule simple et décrits, en partie, dans un essai
\emph{A prologue to \textup{``}Functoriality and Reciprocity\textup{''}, Part~\textup{I}},  que j'ai
rédigé récemment et qui paraîtra bientôt.

Ceci dit, il me semble que toute conséquence laissée de côté le livre 
de Labesse-\hspace*{0pt}Waldspurger serait en soi
un plaisir de lire et peut servir à introduire un novice à la théorie analytique
des formes automorphes et à la formule des traces générale.

\aufm{Robert Langlands}
\chapter*{Préface}

{\renewcommand{\thesection}{\arabic{section}}
\section{La genèse du texte}

La formule des traces pour un groupe réductif connexe sur un corps de caractéristique zéro
est due à James Arthur. On renvoie à \cite{Aintro} pour une introduction et une bibliographie
complète. 

Le cas tordu a fait l'objet du \textenglish{Fri\-day Morn\-ing Sem\-i\-nar} à l'\textenglish{In\-sti\-tute for Ad\-vanced Stud\-y} (IAS) de
Princeton en 1983-1984, souvent cité 
dans la littérature sous le nom de \textenglish{Morn\-ing Sem\-i\-nar on the Trace For\-mu\-la}.
Lors de ce séminaire, les exposés ont été présentés par Laurent Clozel,
Jean-Pierre Labesse et Robert Langlands.
Les exposés~1, 2, 6, 7, 8 et~15 de Langlands ainsi que les exposés~3, 4, 5, 9, 12 et~13 de Labesse
ont donné lieu à des notes, rédigées et distribuées au fur et à mesure. 
Les exposés~10, 11 et~14 de Clozel n'ont pas été rédigés.
Ces notes, citées \cite{MS} dans la suite,
sont accessibles sur la page web de Langlands à l'IAS\@.
Toutefois, ayant été rédigées dans l'urgence, elles laissent
à désirer sur de nombreux points.

Notre ambition est de donner, en nous basant pour l'essentiel sur les notes de \cite{MS},
une version complète de la preuve de la formule des traces 
dans le cas tordu dans sa version primitive, \cad non invariante.
Ce travail s'inscrit dans le projet de l'équipe parisienne animée par L.~Clozel 
et J.-L.~Waldspurger pour rédiger la variante tordue de la formule
des traces et de sa stabilisation, outil indispensable sur lequel se fondent les travaux récents
d'Arthur sur les groupes classiques. En effet, ceux-ci reposent sur
la stabilisation de la formule des traces pour $\GL(n)$ tordu
par l'automorphisme $x\mapsto {}^tx\moins$.

Cette rédaction a dans un premier temps été menée en collaboration entre Laurent Clozel et Jean-Pierre Labesse.
On doit savoir gré à Clozel d'avoir accepté de tenter cette aventure où Labesse craignait de s'engager seul,
même si, en définitive, cette collaboration s'est interrompue
et si c'est Jean-Loup Waldspurger qui a collaboré pour la fin de ce travail.
Il convient de dire que Clozel a écrit un premier jet pour certaines sections,
relu diverses versions préliminaires et participé à de
nombreuses discussions qui ont permis de progresser dans la compréhension
de points obscurs. Qu'il en soit ici remercié.

Nous devons bien entendu remercier tout particulièrement R.~P.~Langlands de nous avoir permis d'utiliser
les notes du séminaire de Princeton \cite{MS} et singulièrement le texte de son dernier exposé
qui contient une esquisse des parties les plus originales et les plus difficiles de la 
preuve dans le cas tordu. Ce texte a été notre guide, même si nous avons dû nous
en écarter en certains points.

\section{Contenu des divers chapitres}

Nous allons maintenant décrire brièvement le contenu des parties et chapitres.
Les deux premières parties sont le plus souvent une simple réexposition du contenu de
{\cite{ATFI}, \cite{ATFII}} et, partiellement, \cite{Arinvo} avec quelques compléments pour les adapter au cas tordu.
Comme dans \cite{MS}, mais de manière plus systématique,
nous avons préféré réexposer ces articles plutôt que de renvoyer à la littérature car, avec le temps, la structure
des preuves est apparue plus clairement et il est désormais possible de les présenter dans un ordre
plus naturel et plus facile à suivre pour le lecteur; au surplus cela rend l'extension au cas tordu transparente.

Dans la troisième partie, la torsion joue un rôle plus important, en compliquant quelque peu les preuves de convergence,
mais là encore, comme dans \cite{MS}, nous suivons de près {\cite{ATFI} et \cite{ATFII}}. 
Les trois premières parties couvrent les exposés~1 à~14 de \cite{MS}.

La quatrième partie, qui donne l'extension au cas tordu
de \cite{AeisI} et \cite{AeisII}, reprend pour l'essentiel le contenu de \cite{MS}*{Lecture~15}
à ceci près que nous avons dû nous en écarter quelque peu pour le calcul de certains termes.
 Dans cette partie, la torsion joue un rôle essentiel en introduisant des termes
qui étaient absents ou négligeables dans le cas classique (\ie non tordu) et dont l'étude est très délicate.

\partintro{1}{Géométrie et combinatoire}

Cette partie contient trois chapitres sur la géométrie 
des groupes et espaces tordus ainsi que sur la combinatoire des cônes et convexes
associés aux systèmes de racines. Sauf naturellement dans le chapitre~\ref{esptor}, 
qui introduit les espaces tordus, la torsion ne joue guère de rôle.
 Mais, faute de référence commode et comportant des preuves complètes ainsi que 
pour convaincre le lecteur que l'extension au cas tordu était facile, il nous a souvent 
paru nécessaire d'exposer en détail le cas classique.

\chapintro{1}{Racines et convexes}

Nous rappelons tout d'abord la construction des espaces vectoriels $\gaPQ$
associés aux paires de sous-groupes paraboliques $P\subset\Q$ d'un groupe réductif $\G$
défini sur un corps de nombres $F$,
ainsi que la propriété fondamentale pour la combinatoire des cônes associés aux 
racines: à savoir le fait que les bases $\Delta_P^\Q$ sont obtuses. 
Puis nous rappelons quelques propriétés, élémentaires et classiques, 
des éléments et des sous-ensembles des groupes de Weyl, qui interviennent fréquemment
en particulier via la décomposition de Bruhat.
Nous en fournissons des preuves lorsque nous ne connaissons pas de références commodes.
Ensuite nous donnons des énoncés 
concernant les familles de cônes et de convexes attachées aux systèmes de racines
et leur relation avec les \GM-familles. Nous reprenons pour l'essentiel les preuves
données dans \cite{MS}*{Lecture~13} où on voit que beaucoup d'énoncés combinatoires sont des conséquences de
la simple identité matricielle $\tau\htau=\htau\tau=1$ (\cf proposition~\ref{totoun}).

Nous n'avons pas toujours repris les preuves classiques. De plus, certains énoncés semblent
nouveaux, quoique implicites chez Arthur ou Langlands;
c'est par exemple le cas des lemmes~\ref{soma} et \ref{decomp}.
La preuve des propriétés~\ref{lisse} et~\ref{decompGM}
des \GM-familles à partir de la combinatoire des cônes,
via la transformée de Fourier est inspirée par le traitement de la combinatoire dans \cite{MS}*{Lecture~15}. 
La clef en est l'énoncé de globalisation~\ref{exist} qui lui aussi semble nouveau. 

\chapintro{2}{Espaces tordus}

Pour l'étude de la formule des traces tordue,
il est commode d'utiliser le langage des espaces tordus introduit dans \cite{Ltw}
(certains préfèrent parler de groupes tordus). Nous en rappelons la définition.
Notre cadre, celui des espaces tordus, est une variante légèrement plus générale
du cadre utilisé dans le \textenglish{Morn\-ing Sem\-i\-nar} et repris
par Arthur dans divers articles ultérieurs. C'est, aux notations près, le cadre de \cite{KS}.

L'extension au cas tordu de la combinatoire des cônes associés aux poids et racines
est immédiate en observant 
que la seule propriété des systèmes de racines utilisée par cette combinatoire 
dans le cas usuel (non tordu), est que les racines simples forment une base obtuse; or 
la généralisation au cas tordu de cette propriété est elle aussi immédiate.

 Dans la section \ref{volconvpol} on montre que le volume de certains convexes peut
 se calculer au moyen de polynômes du premier degré en chaque variable
 pour un choix astucieux des variables paramétrant le convexe.
 Ceci montre que dualement des termes définis au moyen de
 certaines \GM-familles peuvent s'exprimer
 au moyen de produits de dérivées du premier ordre
 en chacune de ces variables.
 Ce résultat, de nature combinatoire, dû à Arthur dans le cas non tordu et que nous généralisons,
 peut être vu comme un cas très simple de résultats plus généraux de Finis et
 Lapid \cite{FL}.
 Ceci permet d'étendre au cas tordu les techniques d'Arthur n\'ecessaires
 pour la preuve du théorème~\ref{mainthc}.
}

On introduit ensuite la fonction caractéristique de cône
$\tsQR$ qui joue un rôle essentiel dans \og l'identité fondamentale\fg (qui fait l'objet du chapitre~\ref{ch8}). 
Sa définition est légèrement plus subtile que pour
son analogue non tordu $\sQR$.

Le chapitre se conclut par diverses inégalités liées à la géométrie de cônes 
qui elles sont spécifiques au cas tordu
(en particulier le lemme~\ref{bigron} qui provient de \cite{MS}*{Lecture~15}).

\chapintro{3}{Théorie de la réduction}

Ce chapitre contient essentiellement
la définition et les propriétés de la fonction $\HO$
sur les groupes adéliques ainsi que des rappels sur la théorie de la réduction.
Il s'agit, là encore, de propriétés très classiques ne faisant pas intervenir la torsion;
de fait, la torsion n'intervient que très peu dans tout ce chapitre. 

La fonction $\HO$, qui se définit via la décomposition d'Iwasawa,
fait le lien entre la géométrie du groupe et celle des espaces vectoriels
associés aux racines. Les lemmes du paragraphe~\ref{hown},
qui permettent le contrôle de $\HO(wn)$ lorsque $n$ est dans l'unipotent et $w$ dans le groupe de Weyl,
sont pour l'essentiel empruntés à \cite{MS}*{Lecture~6}). Ces lemmes jouent un rôle
important dans de nombreuses estimations, en particulier dans le chapitre suivant.
La partition de la section~\ref{unepartXG} et les estimées de la section~\ref{sec3.7} sont empruntées
à \cite{ATFI} (voir aussi \cite{MS}*{Lectures~3 et~4}).

\partintro{2}{Théorie spectrale, troncatures et noyaux}

Cette partie est pour l'essentiel un exposé de résultats classiques 
sur l'opérateur de troncature et la décomposition spectrale de l'espace des formes automorphes,
qu'il était nécessaire de rappeler au moins pour introduire les notations.
La torsion ne joue encore ici qu'un rôle accessoire.
Toutefois quelques nouveautés apparaissent ici où là. 

\chapintro{4}{L'opérateur de troncature}

Ce chapitre rappelle des faits bien connus, dus à Arthur, sur l'opérateur de troncature.
La torsion n'intervient pas du tout ici. On suit pour l'essentiel l'exposé~6 de Langlands
\cite{MS}*{Lecture~6} qui soi-même s'inspire du contenu du premier paragraphe
de l'article d'Arthur \cite{ATFII}. 

Le résultat technique le plus important de ce chapitre est le lemme~\ref{troncnul} 
qui reprend \cite{ATFII}*{Lemma~1.1}.
Les arguments de la preuve de ce lemme, essentiel pour la suite, 
semblent légèrement incomplets dans \cite{ATFII}. 
En effet, Arthur y utilise l'analogue de notre lemme~\ref{wn} mais sous une forme forte: \cad avec $c=0$. 
Cette forme forte est prouvée dans les notes de Langlands pour les groupes de Chevalley 
avec un choix optimal du sous-groupe compact maximal \cite{MS}*{Lemma~6.3}.
Mais il ne semble pas possible d'établir cette forme forte en toute généralité.
Fort heureusement, la preuve donnée par Langlands dans \cite{MS},
et que nous reprenons, montre que la forme forte du lemme~\ref{wn} n'est pas indispensable
pour prouver le lemme~\ref{troncnul}.
Pour le reste les arguments sont dus à Arthur.

Le second résultat technique important est la proposition~\ref{rapdec}
qui reproduit le lemme~6.6 de \cite{MS}*{Lecture~6} lui-même emprunté à \cite{ATFII}*{Lemma~1.4}.
Les arguments sont rappelés pour la commodité du lecteur.

\chapintro{5}{Formes automorphes et produits scalaires}

Après un bref rappel des résultats dus à Langlands sur le prolongement méromorphe 
des opérateurs d'entrelacement et des séries d'Eisenstein, on donne une preuve
simple de la formule, également due à Langlands, 
pour le produit scalaire de deux séries d'Eisenstein tronquées,
provenant de fonctions cuspidales,
au moyen de la \GM-famille spectrale. La preuve, donnée ici, est celle qui est esquissée
dans \cite{MS}*{Lecture~12}; elle est beaucoup plus directe et élémentaire que celle rédigée
par Arthur dans \cite{ATFII}. Dans le cas où les fonctions ne sont plus cuspidales, on ne
dispose alors que d'une formule asymptotique.
Le passage du cas cuspidal au cas non cuspidal est, lui, dû à Arthur. Nous nous contentons
de citer le résultat et nous renvoyons à la littérature pour sa preuve.

\chapintro{6}{Le noyau intégral}  

On introduit dans le cas tordu le noyau de la formule des traces
et on en donne des estimées. On rappelle la factorisation de Dixmier-Malliavin que nous
substituons dans diverses preuves
à l'argument de paramétrix utilisé par Arthur, qui lui est emprunté à Duflo-Labesse.

\chapintro{7}{Décomposition spectrale} 

La décomposition spectrale pour le noyau joue bien évidemment un rôle essentiel dans le développement
spectral de la formule des traces. La décomposition spectrale, due à Langlands, est brièvement rappelée. 
Puis on donne des estimées pour la décomposition spectrale du noyau.

\partintro{3}{La formule des traces grossière}

L'adjectif grossier se veut la traduction de \textenglish{coarse} utilisé dans \cite{MS}.
Dans cette partie on introduit tout d'abord l'identité fondamentale (c'est la \textenglish{ba\-sic i\-den\-ti\-ty} de \cite{MS}) qui donne naissance aux
développements géométrique et spectral de la formule des traces.
Puis on étudie le développement géométrique sous sa forme grossière mais aussi fine
(quoique très rapidement).
Ensuite on donne le développement spectral sous sa forme grossière
(\textenglish{coarse spec\-tral ex\-pan\-sion}). Ceci permet de prouver une première forme
de la formule des traces ainsi que les propriétés formelles
des termes des développements grossiers de deux membres de cette identité. 

\chapintro{8}{Formule des traces\textup: état zéro}

Ce chapitre contient la preuve de l'identité fondamentale
qui est le point de départ de la formule des traces.
On établit l'égalité de deux variantes tronquées pour la restriction à la diagonale
du noyau. L'une se prête bien au développement géométrique, \cad suivant les classes de conjugaison,
l'autre au développement spectral. Dans le séminaire de Princeton
une première forme de l'identité fondamentale est établie dans 
\cite{MS}*{Lecture~2} puis, une variante est donnée beaucoup plus tard dans \cite{MS}*{Lecture~9};
c'est cette variante qui s'avère être la bonne et qui est donnée ici à la proposition~\ref{fonfon}.
Il s'agit d'une simple identité combinatoire.

\chapintro{9}{Développement géométrique}

Ce chapitre est consacré à la décomposition 
suivant les classes de conjugaison de l'intégrale sur la diagonale du noyau après troncature \og géométrique\fg.
Le théorème~\ref{geoconv} établit la convergence du développement géométrique grossier
(\cad du développement
suivant les classes de conjugaison des parties quasi semi-simples).
 C'est une adaptation facile des argument de \cite{ATFI}.
On suit pour cela \cite{MS}*{Lectures~3 et~4}. On continue ce chapitre en donnant 
l'expression des termes associés aux classes de conjugaison semi-simples
au moyen d'intégrales orbitales pondérées suivant \cite{MS}*{Lecture~5}
repris et développé dans \cite{MS}*{Lecture~9}. 

Un dernier et bref paragraphe est consacré au développement géométrique fin
(\textenglish{fine $\go$-ex\-pan\-sion}). Il nous a paru suffisant de
renvoyer à la littérature pour le traitement des termes non semi-simples.
D'ailleurs, il n'y a rien concernant ces termes dans \cite{MS}.
En effet le traitement de ces termes n'a été fait, par Arthur, qu'après le \textenglish{Morn\-ing Sem\-i\-nar}.
Comme ceci a été rédigé par Arthur en y incluant le cas tordu 
(quoique dans un cadre légèrement plus restrictif
que le cas général traité par ailleurs dans notre texte),
il ne nous a pas paru nécessaire d'en reprendre la rédaction.

\chapintro{10}{Développement spectral grossier}

La décomposition spectrale suivant les \og données cuspidales\fg induit
le développement grossier (appelé \textenglish{coarse spec\-tral ex\-pan\-sion} dans \cite{MS}). La preuve
de sa convergence suit celle donnée par Langlands dans \cite{MS}
qui avait fait l'objet des exposés~7 et~8, preuve
qui est elle même inspirée de \cite{ATFII}, quoique la torsion induise quelques complications
techniques. La principale différence entre le cas classique et le cas tordu est
que (avec les notations du théorème~\ref{chiconv}) dans le cas tordu, le développement spectral fait intervenir 
une combinaison linéaire de termes indexés par des paires de sous-groupes
paraboliques standard $\Q\subset\R$:
\[
\int_{\YQdo} \tsQR (\xT) \tronc_1^{T,\Q } K_{\Q,\delto,\chi}(x,x)\dd x
\]
pouvant donner des contributions non triviales alors que dans le cas classique seul le terme
\[ 
\int_{\XG} \tronc_1^{T} K_{\chi}(x,x)\dd x
\]
correspondant au cas $\Qdo=\Q=\R=\G$, est non nul (pour $\T$ assez régulier).

\chapintro{11}{Formule des traces\textup: propriétés formelles}

Les termes des développements géométriques et spectraux \og grossiers\fg (appelés
\textenglish{coarse ex\-pan\-sions} dans \cite{MS}) ont des propriétés formelles
remarquables. La propriété essentielle est que l'on obtient, de façon asymptotique,
des polynômes en la variable de troncature $\T$. Les preuves dans le cas tordu
sont une adaptation immédiate des preuves données par Arthur dans \cite{Arinvo} pour le cas 
classique. Nous suivons ici \cite{MS}*{Lecture 13}.

Au total, les trois premières parties fournissent une preuve complète de
la variante tordue de l'ensemble des résultats d'Arthur contenus dans
{\cite{ATFI}  et \cite{ATFII}} ainsi qu'une partie des résultats de \cite{Arinvo} 
(essentiellement ceux concernant les \GM-familles et les propriétés formelles des termes de la formule des traces).

\partintro{4}{Forme explicite des termes spectraux}

Cette partie, la plus difficile et la plus originale de tout l'ensemble,
est consacrée à l'extension au cas tordu des résultats des articles
\cite{AeisI} et \cite{AeisII} d'Arthur. La difficulté nouvelle, par rapport au cas traité par Arthur,
 provient de la nécessité de prendre en compte des termes
attachés à des couples $\Q\subset\R$ avec $\Q\neq\G$ évoqués ci-dessus. 
L'analyse de leur comportement est beaucoup plus délicate. 
 
L'étude du développement spectral
de ces termes utilise le calcul du produit scalaire
de séries d'Eisenstein tronquées qui peut être fait explicitement,
au moins dans le cas où on part de séries d'Eisenstein construites à partir de
fonctions cuspidales, en se ramenant, moyennant une inversion d'intégrale,
 au calcul classique et rappelé ci-dessus (\cf chapitre~\ref{ch5}).
 On obtient alors une expression au moyen de \GM-familles spectrales
 généralisant le cas classique.
Toutefois, pour les termes attachés à des couples $\Q\subset\R$ avec $\theto(\Q)\ne\Q$,
le calcul auquel on est naturellement amené suppose, pour être convergent, d'avoir auparavant déplacé
le contour d'intégration en dehors du domaine naturel des variables spectrales (\cad qu'elles ne sont plus imaginaires pures),
du moins pour une partie d'entre elles. Cela se fait sans grosses difficultés.
Mais, pour achever la combinatoire il convient, calcul fait,
 de revenir ensuite au domaine naturel pour les variables spectrales. 
Il faut donc déplacer des contours d'intégration dans des intégrales
faisant intervenir des \GM-familles. C'est la démarche proposée par Langlands dans \cite{MS}*{Lecture~15}. 
Cela suppose des estimées sur les opérateurs
d'entrelacements et leur dérivées que nous n'avons pas su obtenir.

Une méthode ne supposant pas de déplacement de contour, mais 
très délicate du point de vue combinatoire et analytique, découverte par Waldspurger,
a permis de résoudre la question. On se ramène en définitive à l'expression donnée par Langlands.
 
\chapintro{12}{Introduction d'une fonction $B$}

Il s'agit d'adapter au cas tordu une technique due à Arthur et développée dans \cite{AeisI}. 
L'introduction d'une fonction $B$ à support compact dans l'expression spectrale
pour les termes évoqués ci-dessus va permettre de pallier l'absence d'estimées uniformes
de certains développements spectraux. Comme dit plus haut le traitement
des termes attachés aux couples $\Q\subset\R$ avec $\Q\ne\G$ 
est en général beaucoup plus difficile que le cas $\Q=\G$ traité par Arthur. 
La fonction $B$ apparaît le plus souvent dans les calculs via sa transformée de Fourier. 
Celle-ci n'est pas à support compact, mais seulement à décroissance rapide, 
ce qui pose de délicats problèmes de convergence. 
Pour les traiter, on a besoin de majorations plus fines que dans les paragraphes précédents.

\chapintro{13}{Calcul de $A^T(B)$}

Ce chapitre peut être vu comme l'analogue tordu de la seconde partie de \cite{AeisI}.
Il s'agit, entre autre, de tenir compte du caractère asymptotique des expressions
en termes de \GM-familles obtenues par le calcul de produit scalaire
dans le cas où les séries d'Eisenstein ne sont pas construites à partir de fonctions cuspidales.
Ici encore une difficulté nouvelle provient des termes avec $\theto(\Q)\ne\Q$.
La démarche empruntée ici fournit au total une expression plus simple que celle
obtenue par Arthur dans le cas non tordu. Elles diffèrent par des termes asymptotiquement
petits; il en résulte qu'une des étapes combinatoires de \cite{AeisII}*{section~3}
se trouve ainsi déjà prise en compte.

\chapintro{14}{Formules explicites}

Ce chapitre exploite l'analyse faite dans les deux chapitres précédents
pour obtenir dans le cas tordu, l'analogue des formules obtenues
par Arthur dans \cite{AeisII} donnant l'expression explicite des termes spectraux
de la formule des traces. 

La section~\ref{sec14.1} s'inspire du traitement proposé
par Langlands dans \cite{MS}*{Lecture~15} mais en utilisant de façon systématique
la globalisation des \GM-familles 
ce qui rend plus transparent l'argumentaire combinatoire et simplifie considérablement tant 
cette combinatoire que les notations.

L'objet de la section~\ref{sec14.2} est de débarrasser les divers termes de la fonction
auxiliaire $B$ en la faisant tendre vers~$1$. Pour cela il convient d'établir la convergence absolue
de ces termes. Arthur utilise deux arguments:
\begin{enumerate}
\item il\hspace*{-.2pt} suppose\hspace*{-.2pt} l'existence\hspace*{-.2pt} d'une\hspace*{-.2pt} normalisation\hspace*{-.2pt} des\hspace*{-.2pt} opérateurs d'entrelacements;
\item il montre que les termes à contrôler, qui font intervenir
 des \GM-familles définis au moyen des facteurs de normalisation, peuvent s'exprimer
 comme une combinaison linéaire de produits de dérivées du premier ordre
 en certaines variables. Ceci permet la réduction à un problème en rang un.
\end{enumerate}

L'existence d'une normalisation a été
établie pour la première fois par Langlands dans \cite{MS}*{Lecture~15}. Cette
normalisation a depuis été reprise par Arthur. N'ayant rien à ajouter,
 nous nous contentons de citer Arthur \cite{Ainter} pour cette normalisation ainsi que \cite{AeisII}
 pour la fin de la preuve, à un détail près qui fait l'objet de la section \ref{volconvpol}.

 Une dernière section reformule le développement spectral en exploitant 
 la convergence absolue due à Finis, Lapid et Müller.

\begin{otherlanguage}{american}
\chapter*{Foreword}
\setcounter{section}{0}

{\renewcommand{\thesection}{\arabic{section}}
\section{The genesis of the paper}

The trace formula for an arbitrary connected reductive group over a number field is due
 to James Arthur. We refer the reader to \cite{Aintro} for an introduction and a complete
 bibliography.
 
The twisted case was the subject of the \emph{Friday Morning Seminar} at the Institute for Advanced Study
(Princeton) during the academic year 1983--1984, often quoted
in the literature as \emph{Morning Seminar on the Trace Formula}.
During this seminar lectures were given by Laurent Clozel,
Jean-Pierre Labesse and Robert Langlands.
Notes for Lectures~1, 2, 6, 7, 8, and~15 by Langlands and Lectures~3, 4, 5, 9, 12, and~13 by Labesse
were written up and made available to the audience a few days after each lecture.
Lectures~10, 11, and~14 by Clozel were never written up.
The lecture notes, quoted \cite{MS} in the sequel,
are available on Langlands webpage.
But, having been written quite hastily they contain quite a few errors,
and in addition some proofs are not complete.

Our ambition is to give, following \cite{MS},
a complete proof of the twisted trace formula
in its primitive version, i.e., its noninvariant form.
This is a part of the project of the Parisian team led by L.~Clozel 
and J.-L.~Waldspurger whose aim is to give a complete proof of the
stablilization of the twisted trace formula
which is the basic tool for Arthur's book on classical groups. 
In fact it relies on the stabilization of the trace formula for $\GL(n)$ twisted
by the automorphism $x\mapsto {}^tx\moins$.

At the beginning, this book was a collaboration between
 Laurent Clozel and Jean-Pierre Labesse.
We are grateful to Laurent Clozel to have agreed to try
 this adventure where Labesse was afraid to embark alone,
even if, at some point, this collaboration was stopped and 
Jean-Loup Waldspurger helped to finish the job.
It is fair to say that Clozel has written the first draft for some sections,
read many preliminary versions and helped clear up many obscure points.
We thank him very much for this.

We are glad to thank R.~P.~Langlands who allowed us to use the IAS lecture notes
 \cite{MS} and in particular the notes of his Lecture~15 which contains
a sketch of the most original and most difficult part of the proof in the twisted case.
 These notes were our guide even if we had to follow a slightly different path at some point.

\section{Contents of the chapters}

We shall now describe briefly the contents of the various parts and chapters.
The first two parts are often simply a rewriting of the contents of 
\cite{ATFI},  \cite{ATFII} and part of \cite{Arinvo}
 with a few additions to fit with the twisted case.
 As in \cite{MS}, but in a more systematic way,
we preferred to repeat the arguments rather than to refer to the literature since 
by now the structure of the proofs is much better understood and it is possible
to present them in a more natural order, one that is easier to follow for the reader;
moreover this makes the extension to the twisted case more or less obvious.

In the third part, the twisting plays a more important role, as it makes the proofs for the convergences slightly more complicated, but again, as in \cite{MS}, we follow closely {\cite{ATFI} and \cite{ATFII}}. 
The first three parts cover Lectures~1 to~14 in \cite{MS}.

The fourth part, which extends \cite{AeisI} and \cite{AeisII} to the twisted case, 
is mainly based on \cite{MS}*{Lecture~15}
 except that we have used a slightly different approach for the computation of some terms.
In this part, the twisting plays an essential role by introducing terms
that are absent or negligible in the classical (i.e., nontwisted) case
and whose study is quite subtle. 

\partintro{1}{Geometry and combinatorics}

This part contains three chapters on the geometry of groups and twisted spaces,
and also on the combinatorics of cones and convex sets attached to root systems.

Except, of course, in Chapter~\ref{esptor} the twisting plays barely no role.
But by lack of a convenient reference with complete proofs and 
also to convince the reader that the extension to the twisted case was easy
we thought it better to give a detailed account for the classical case.

\chapintro{1}{Roots and convex sets}

We recall first the construction of vector spaces $\gaPQ$
attached to pairs of parabolic subgroups $P\subset\Q$ in a reductive group $\G$
over a number field $F$,
and the fundamental property for the combinatorics of cones attached to 
root systems: namely that basis $\Delta_P^\Q$ are obtuse. 
Then we recall some classic and elementary
 properties of elements and subsets of Weyl groups that arise quite often in particular
when using Bruhat decomposition.
We give proofs when we don't know of any easily accessible reference.
Then we give statements for family of cones and convex sets 
attached to root systems and their relation to \GM-families. Most of the time we
follow the proofs given in \cite{MS}*{Lecture~13} where it is shown that
many of these statements are consequences of the simple matrix identity
 $\tau\htau=\htau\tau=1$ (cf.~Proposition~\ref{totoun}).

Some of our proofs do not follow the classical patterns.
 Moreover some of the statements seem to be new, albeit implicit
 in Langlands or Arthur.
This is, for example, the case of Lemmas~\ref{soma} and~\ref{decomp}.
The proof of Lemmas~\ref{lisse} and~\ref{decompGM}
for \GM-families using,
via a Fourier transform, the combinatorics of cones
is inspired by the use of this combinatorics in \cite{MS}*{Lecture~15}. 
The key is the globalization statement~\ref{exist} that also seems new. 

\chapintro{2}{Twisted spaces}

To study the twisted trace formula 
it is convenient to use the language of twisted spaces introduced in \cite{Ltw}
(some prefer to speak of twisted groups). We recall the definition.
Our setting, twisted spaces, is a slightly more general variant of the setting used in 
the Morning Seminar and again by Arthur in the papers written after it.
Up to notation this is the setting used in \cite{KS}.

The extension to the twisted case of the combinatorics of cones associated to roots and weights
is immediate when observing that it only relies on the following fact: simple roots
define an obtuse basis and that this is also the case for the corresponding basis 
that arise in the twisted case.

In section \ref{volconvpol} we show that the volume of certain convex sets can
be computed in terms of polynomials of degree one in each variable
for a good choice of the parameters defining the convex set.
This shows that dually some terms defined by certain \GM-families
can be computed in terms of products of derivatives of first order
in each of these variables. This combinatorial result, due to Arthur in the
nontwisted case, can be seen as a very simple instance of quite general results
of Finis and Lapid \cite{FL}.
This allows to extend to the twisted case Arthur's techniques in order to prove
Theorem~\ref{mainthc}.}

We then introduce $\tsQR$ the characteristic function of another cone
 which is the key object in the ``basic identity'' (treated in Chapter~\ref{ch8}).
 Its definition is slightly more subtle 
 than that of its classical analogue.
 
 The chapter ends with various inequalities stemming from the geometry of cones
 specific of the twisted case (in particular Lemma~\ref{bigron} 
 borrowed from \cite{MS}*{Lecture~15}).

\chapintro{3}{Reduction theory}

This chapter is mainly concerned with the
definition and the properties of function $\HO$
on adelic groups and other classical statements of reduction theory.
Here again, twisting plays little if any role in the whole chapter.

The function $\HO$, defined via Iwasawa decomposition,
allows to relate the geometry on the group and that of vector spaces
attached to roots.

The lemmas in Section~\ref{hown},
that allow to control $\HO(wn)$ when $n$ belong to the unipotent and $w$ is in the Weyl group,
are borrowed from \cite{MS}*{Lecture~6}. These lemmas are quite important in order
to obtain estimates, in particular in the next chapter.
The partition in Section~\ref{unepartXG} and the estimates in Section~\ref{sec3.7} are borrowed from \cite{ATFI} (see also \cite{MS}*{Lectures~3 and~4}).

\partintro{2}{Spectral theory, truncation and kernels}

This part is mainly a report on classical results about the truncation operator
and the spectral decomposition of the space of automorphic forms, which
we had to recall at least to fix notation. Here again twisting plays little role.
Nevertheless some new features show up.

\chapintro{4}{The truncation operator}

This chapter recalls well known facts, due to Arthur, on truncation operators.
The twisting is completely absent here. We follow Lecture~6 by Langlands \cite{MS} which is itself inspired by Arthur's paper \cite{ATFII}. 

The main technical result of the chapter is Lemma~\ref{troncnul} 
which is nothing but \cite{ATFII}*{Lemma~1.1}.
The proof of this key lemma seems slightly incomplete in \cite{ATFII}. 
In fact Arthur makes use of an analogue of our \ref{wn} but in a stronger form,
i.e., with $c=0$. This strong form is established by Langlands for Chevalley groups
under an optimal choice for the maximal compact subgroup \cite{MS}*{Lemma~6.3},
but such a strong form is not likely to hold in general.
Fortunately the proof given by Langlands in \cite{MS}, which we reproduce here, shows that
such a strong form of Lemma~\ref{wn} is not necessary to establish Lemma~\ref{troncnul}.
The remaining arguments are due to Arthur.

The second important technical result is Proposition~\ref{rapdec}
which reproduces Lemma~6.6 in \cite{MS}*{Lecture~6} which is itself borrowed from
\cite{ATFII}*{Lemma 1.4}. The proof is recalled for the convenience of the reader.

\chapintro{5}{Automorphic forms and scalar products}

We first recall briefly the results due to Langlands on analytic continuation of intertwining operators
and Eisenstein series, then we give a simple proof of the formula, also due to Langlands, for the
scalar product of truncated Eisenstein series built from cuspidal functions
in terms of the spectral \GM-family. The proof given here was outlined in \cite{MS}*{Lecture~12}; 
it is much more direct and elementary than the proof given by Arthur in \cite{ATFII}.
When dealing with the noncuspidal case one gets only an asymptotic formula.
The extension to the noncuspidal case is due to Arthur. We only quote the result and we
refer the reader to the literature for a proof.

\chapintro{6}{The integral kernel} 

We introduce the kernel that is used in the trace formula in the twisted case
and we give estimates. We recall Dixmier--Malliavin's factorization
which we substitute for the parametrix argument used by Arthur, which itself
was borrowed from Duflo--Labesse.

\chapintro{7}{Spectral deecomposition} 

The spectral decomposition for the kernel plays naturally a key role in the spectral
expansion of the trace formula. The spectral decomposition, due to Langlands, is briefly recalled.
Then we give estimates for the spectral decomposition of the kernel.

\partintro{3}{The coarse trace formula}

The word ``coarse'' used in \cite{MS} is to be translated by ``\textfrench{gros\-sier}'' in French.
In this part we first give the basic identity (``\textfrench{iden\-ti\-té fon\-da\-men\-tale}'' in French)
which gives rise to the geometric and spectral expansions of the trace formula.
Then we study the coarse geometric expansion and the fine geometric expansion as well
 (although very briefly). Next we give the coarse spectral expansion.
 This yields a first form of the trace formula, and we establish the formal properties
 of the various terms in the coarse expansions of both sides of this identity.

\chapintro{8}{Zero state of the trace formula}

This chapter contains the proof for the basic identity which is the starting point
for the trace formula. We establish the equality between two truncated variants
of the restriction to the diagonal of the kernel. One side yields the geometric expansion
by expanding it according to conjugacy classes, the other one yields the spectral expansion.
During the Princeton seminar a first form of the basic identity was established in \cite{MS}*{Lecture~2}
then a variant of it was given quite later in \cite{MS}*{Lecture~9};
this variant turned out to be the right one and is given here in Proposition~\ref{fonfon}.
This is merely a combinatorial identity.

\chapintro{9}{Geometric expansion}

This chapter deals with the decomposition along conjugacy classes of the
integral over the diagonal of the kernel after geometric truncation.
Theorem~\ref{geoconv} establishes the convergence of the coarse geometric expansion
(i.e., the expansion along the conjugacy classes of the quasisemisimple parts).
This is an easy adaptation of arguments in \cite{ATFI}.
We follow \cite{MS}*{Lectures~3 and~4}. The chapter goes on with the expression of terms
defined by semisimple conjugacy classes by means of weighted orbital integrals, following
 \cite{MS}*{Lecture 5} and again but with more details in \cite{MS}*{Lecture 9}. 

A last and short paragraph deals with the fine geometric expansions (``fine $\go$-expansion'').
We thought it enough to refer the reader to the literature for the study of nonsemisimple
elements. We observe that these terms are not studied in \cite{MS}. 
In fact they were treated by Arthur only after the Morning Seminar.
Since this was written up by Arthur including the twisted case
(although in a setting slightly more restrictive than here)
we thought it worthless to rewrite it here.

\chapintro{10}{The coarse spectral expansion}

The spectral expansion according to cuspidal data induces the coarse spectral expansion.
The proof of its convergence follows Lectures~7 and~8 by Langlands in \cite{MS}
which in turn are inspired by \cite{ATFII}, although a few technical difficulties arise from
the twisting. The main difference between the classical and the twisted case
is that (with notation from Theorem~\ref{chiconv}) the twisted spectral expansion
is a linear combination of terms indexed by pairs of standard parabolic subgroups
 $\Q\subset\R$:
\[
\int_{\YQdo} \tsQR (\xT) \tronc_1^{T,\Q } K_{\Q,\delto,\chi}(x,x)\dd x
\]
that can be nontrivial, while in the classical case there is only one term
\[ 
\int_{\XG} \tronc_1^{T} K_{\chi}(x,x)\dd x
\]
corresponding to the case $\Qdo=\Q=\R=\G$,
that is nonvanishing (for $T$ regular enough).

\chapintro{11}{Trace formula\textup: formal properties}

The terms in the geometric and spectral coarse expansions have remarkable
formal properties. Their main property is that one gets, asymptotically
polynomials in $T$, the truncation variable. The proofs in the twisted case are
an easy adaptation of the proofs given by Arthur in \cite{Arinvo} in the classical case. We follow here 
\cite{MS}*{Lecture~13}.

Altogether the first three parts give a complete account for the twisted variant
of the results contained in {\cite{ATFI} and  \cite{ATFII}} 
and a part of the results in \cite{Arinvo} as well (mainly those concerning
the \GM-families and the formal properties of terms in the trace formula).

\partintro{4}{Explicit form for spectral terms}

This part, which is the most original and difficult part of this work, deals with
the extension to the twisted case of Arthur's results in {\cite{AeisI}  and \cite{AeisII}}. 
 The new feature is that one has to take care of terms
indexed by pairs $\Q\subset\R$ with $\Q\ne\G$ alluded to above.
Their study is much more delicate.

The study of the spectral expansion of such terms uses the scalar product
of truncated Eisenstein series that can be computed exactly, at least
when dealing with Eisenstein series constructed from cuspidal functions,
 by means of the classical calculation recalled above
(cf. Chapter~\ref{ch5}) after an inversion of integrals. One thus gets an expression in terms of spectral \GM-families
generalizing the classical case. But for terms attached to pairs
$\Q\subset\R$ with $\theto(\Q)\ne\Q$ this natural formal computation
makes sense only after shifting the contour integral for some of the spectral variables
away from their natural domain (i.e., these are no longer purely imaginary).
This can be done without much trouble.
But now, to go on with the combinatorics, one needs to come back to the natural domain
for these spectral variables. One is thus led to shift back the contour integral
for terms that now involve spectral \GM-families. This is what is suggested by Langlands
in \cite{MS}*{Lecture~15}. To do this one needs estimates on intertwining operators
and their derivatives which we do not know how to get.

 A method that does not involve any shifting of contour, but that is quite delicate
 from the combinatorial and analytic point of view, has been devised by Waldspurger
 and allows to solve this problem. One is eventually led to the expression given by Langlands.
 
\chapintro{12}{Introduction of a function $B$}

One has to adapt a technique due to Arthur and developed in \cite{AeisI}. 
The introduction of a function $B$ with compact support in the spectral expansion
for terms described above is a remedy to the lack of uniform estimates
of certain spectral expansions. As said above the terms attached to pairs
$\Q\subset\R$ with $\Q\ne\G$ is in general much more difficult than the case
$Q=G$ treated by Arthur. The functions $B$ appears often via its Fourier transform
which is not compactly supported but only rapidly decreasing, and this
generates quite delicate convergence problems. To handle them one needs
to refine some of the estimates already used in previous sections.

\chapintro{13}{Computation of $A^T(B)$}

This chapter can be seen as the twisted analogue of the second part of \cite{AeisI}.
One aim is to take care of the asymptotic nature of the formulas involving
the \GM-families coming from the computation of the scalar product when
Eisenstein series are not constructed from cuspidal functions.
Here again a new difficulty comes from terms where $\theto(\Q)\ne\Q$.
The way this difficulty is solved here yields eventually an expression
which is simpler than the one obtained by Arthur in the nontwisted case.
They differ by asymptotically small terms; as a consequence one of
the combinatorial steps in \cite{AeisII}*{Section~3} is already taken care of.

\chapintro{14}{Explicit formulas}

Based on the analysis made in the two preceding chapters, we
establish the twisted analogue of the formulas obtained by Arthur in
\cite{AeisII} giving an explicit expression for spectral terms in the
trace formula.

Section~\ref{sec14.1} is inspired by what Langlands does in
\cite{MS}*{Lecture~15} but the systematic use of the globalization
of \GM-families makes the combinatorial argument much more transparent
since it simplifies the combinatorics and the notation as well.

In Section~\ref{sec14.2} we get rid of the auxiliary function $B$ by letting it
tend to 1. To do this one needs to prove the absolute convergence
of these terms. Arthur uses two kind of arguments:
\begin{enumerate}
\item he assumes known a normalization of the intertwining operators;
\item he shows then that the terms to control, that contain \GM-families,
can be expressed as a linear combination of products of first order derivatives
in certain variables. This allows to reduce to a problem in rank one.
\end{enumerate}

The existence of a normalization was first established by Langlands
in \cite{MS}*{Lecture~15}. This normalization later appeared in Arthur's work.
Since we have nothing more to say we simply refer to \cite{Ainter} for
the normalization and to \cite{AeisII} for the end of the proof up to a detail
which occupies section \ref{volconvpol}.

In a last section we reformulate the result by
making use of the absolute convergence of the spectral expansion
due to Finis, Lapid and Müller.

\end{otherlanguage}

\setcounter{part}{0}
\mainmatter
\part{Géométrie et combinatoire}

\chapter{Racines et convexes}\label{ch1}

\section[Les espaces $\ga_P$]{\mathversion{bold}Les espaces $\ga_P$}\label{vectrac}

\newindex{aP@$\gaP$}{espavect}%
Soit $F$ un corps de nombres. On note $\adef$ l'anneau des adèles de $F$. Soit $\G$
\newindex{G@$\G$}{Groupe}%
un groupe linéaire algébrique connexe défini sur $F$. On note $X_F(\G)$
\newindex{XF@$X_F$}{caracteeres}%
le groupe des caractères rationnels de $\G$ et on pose
\[
\ga_\G=\Hom(X_F(\G),\RM)\ptf
\]
C'est un espace vectoriel sur $\RM$. On note $a_\G$ sa dimension. On dispose alors d'un homomorphisme
\[
\HG\colon\Gadef \to\ga_\G
\]
défini par
\[
\HG(x)= \{\chi\mapsto\log\lvert\chi(x)\rvert\}\ptf
\]
\newindex{HG@$\HG$}{harish-map}%
On notera
\[
\Gadef^1
\]
\newindex{G(A)1@$\Gadef^1$}{harish-map-ker}%
le noyau de cette application.
Considérons une décomposition de Levi: $\G=L\N$ où $\N$ est le radical unipotent de $\G$. L'homomorphisme $\HG$ est trivial sur $\N(\adef)$ ainsi que sur $L_{\mathrm{der}}(\adef)$ où $L_{\mathrm{der}}$ est le groupe dérivé du sous-groupe de Levi $L$.

Supposons maintenant que $\G$ est réductif. On note $Z_\G$ son centre.
Soit $\G_{\QM}$ la restriction des scalaires de $F$ à $\QM$ de $\G$.
On note $\gA_\G$ la composante neutre du groupe des points réels
du $\QM$\hyph tore déployé maximal du centre
de $\G_\QM$. Donc $\gA_\G$
\newindex{AG@$\gA_\G$}{AG}%
est un sous-groupe de Lie
connexe de $Z_\infty$:
\[
\gA_\G\subset Z_\infty\coloneqq Z_\G(F\otimes\RM)\subset Z_\G(\adef)\subset\Gadef\ptf
\]
Par restriction de $\HG$ à $\gA_\G$ on obtient un homomorphisme:
\[
\gA_\G\to\ga_\G
\]
qui est un isomorphisme.
On peut alors interpréter $\ga_\G$ comme l'algèbre de Lie
de $\gA_\G$ et voir $\ga_\G$ comme une sous-algèbre de
l'algèbre de Lie $\frakg$ de $\G_\QM(\RM)$.
On notera
\[
a=\ee^H
\]
le point de $\gA_\G$ d'image $H$ dans $\ga_\G$. On observera que l'application naturelle
\[
\Gadef^1\to\gA_\G\bs\Gadef
\]
est un isomorphisme.

Soit $P=\M\N$ un sous-groupe parabolique de $\G$ où $\M$ est un sous-groupe de Levi de $P$,
définis sur $F$; on observe que
$X_F(P)=X_F(\M)$ et donc
\[
\ga_P=\ga_\M\ptf
\]
On choisit un sous-groupe parabolique minimal $\PO$ et un sous-groupe de Levi $\MO$
de $\PO$. Le groupe $\MO$ est fixé une fois pour toutes dans la suite de ce texte.
On pose
\[
\gao\coloneqq \ga_{\PO}=\ga_{\MO}\ptf
\]
Dans toute la suite, nous ne considèrerons que des sous-groupes paraboliques
$P$ semi-standard \cad contenant $\MO$.
Le sous-groupe de Levi $\M$ est déterminé par $P$ et la condition
$\MO\subset\M$ et
on écrira parfois $\gA_P$ pour $\gA_{\M}$. On prendra garde que
\[
\gA_P\coloneqq \gA_\M
\]
n'est pas central dans $P(\adef)$. Toutefois $\gA_P\N(\adef)$ est un
sous-groupe distingué dans $P(\adef)$.
L'inclusion $P\subset\G$
induit une inclusion
\[
X_F(\G)\subset X_F(P)
\]
et donc une surjection
\[
\ga_P\to\ga_\G
\] dont le noyau sera noté $\ga_P^\G$.
Compte tenu des isomorphismes
\[
\gA_\G\simeq\ga_\G \Qquad{et} \gA_{\M}\simeq\ga_{\M}=\ga_P
\]
et de l'inclusion
\[
\gA_{\G}\subset\gA_{\M}
\] 
on obtient une section $\ga_\G\to\ga_P$ de la surjection $\ga_P\to\ga_\G$ et donc une décomposition
\[
\ga_P=\ga_\G\oplus\ga_P^\G\ptf
\]
Plus généralement, soient $P\subset\Q$ deux sous-groupes paraboliques de $\G$.
En utilisant que $\HQ$ est trivial sur le radical unipotent de $\Q$
ce qui précède fournit une décomposition
\[
\ga_P=\ga_\Q\oplus\gaPQ\ptf
\]
On pose
\[
a_P^\Q=\dim\gaPQ\ptf
\]
\newindex{aPQ@$a_P^Q$}{dimensions}%

On considère $\gao$ comme l'algèbre de Lie de $\gA_{\MO}$ et donc
comme une sous-algèbre de $\frakg$.
La forme de Killing induit un produit scalaire sur $\gao^\G$.
Ceci définit une structure euclidienne sur son dual et
on notera $\langle\alpha,\beta\rangle $ le produit scalaire de deux éléments du dual.
Cette structure euclidienne définit également une mesure sur $\gao^\G$.
Plus généralement on dispose ainsi de mesures canoniques sur les espaces $\ga_P^\G$
vus comme quotients $\gao^\G/\gao^P$.

\section{Sous-groupes paraboliques et bases de racines}

On suppose désormais $\G$ réductif.
On dispose des racines attachées au couple $(\MO,\G)$.
Ce sont des formes linéaires sur $\gao$ nulles sur $\ga_\G$.
L'ensemble de leurs restrictions à $\gao^\G$ est un système de racines,
non réduit en général.
On a choisi un sous-groupe parabolique minimal $\PO$;
on dispose donc de la notion de racines positives
et d'une base de racines simples notée $\Delta_{\PO}^\G$.
On notera $\Rac^\G$
\newindex{R@$\Rac$}{systracines}%
le système de racines réduit formé des racines
$\beta$ telles que $\beta/2$ ne soit pas une racine.
On écrira souvent $\Rac$ pour $\Rac^\G$.
C'est le système de racines réduit admettant $\Delta_{\PO}^\G$ comme base.

Si $P$ est standard de sous-groupe de Levi $\M$,
on notera $\Delta_{\PO}^P$ la base des racines simples
pour le couple $(\MO,\M)$.
C'est une base du dual de $\gao^P=\gao^\M$. On considèrera les
éléments de $\Delta_{\PO}^P$ comme des formes linéaires
sur $\gao$ ou $\gao^P$ suivant les besoins. En particulier
on peut voir $\Delta_{\PO}^P$ comme
un sous-ensemble de $\Delta_{\PO}^\G$.
La combinatoire utilisera de façon systématique le fait bien connu suivant:

\begin{lemme}\label{bijpar}
L'application
\[
P\mapsto\Delta_{\PO}^P
\]
est une bijection entre l'ensemble des sous-groupes paraboliques standard de $\G$
et l'ensemble des parties de $\Delta_{\PO}^\G$.
\end{lemme}

On dispose également de la base des coracines ${\check\Delta}_{\PO}^P$
dans $\gao^P$; on notera $\hDelta_{\PO}^P$ la base duale
de la base des coracines.
Lorsque le groupe est déployé $\hDelta_{\PO}^\G$
est l'ensemble des poids dominants fondamentaux du groupe dérivé.
De plus, $\hDelta_{\PO}^P$ est l'ensemble des restrictions non nulles
des $\vpi\in\hDelta_{\PO}^\G$ au sous-espace $\gao^P$.
Plus généralement, soient $P\subset\Q$ deux sous-groupes paraboliques standard.
On note $\Delta_P^\Q$
\newindex{DeltaPQ@$\Delta_P^Q$}{deltapq}
l'ensemble des restrictions non nulles des
éléments de $\Delta_{\PO}^\Q$ au sous-espace $\ga_P$.
Cet ensemble de formes linéaires est une base du dual $(\gaPQ)^*$ de $\gaPQ$
et $\ga_\Q$ s'identifie avec
le sous-espace de $\ga_P$ intersection
des noyaux des $\alpha\in\Delta_P^\Q$.
On prendra garde toutefois qu'en général $\Delta_P^\Q$ n'est pas la base d'un système
de racines. On notera $\hDelta_P^\Q$ le sous ensemble des $\vpi\in\hDelta_{\PO}^\Q$
\newindex{DeltaPQ@$\protect\hDelta_P^Q$}{hdeltapq}%
nuls sur $\gao^P$. On prolonge les éléments de $\Delta_P^\Q$ et $\hDelta_P^\Q$ en des formes linéaires
sur $\gao$ en les composant avec la projection
\[
\gao\to\gaPQ\ptf
\]
On écrira parfois $\Delta_P$ pour $\Delta_P^\G$ ainsi que
$\hDelta_P$ pour $\hDelta_P^\G$.
On observera que les bases $\Delta_P^\Q$ et $\hDelta_P^\Q$
sont indépendantes du choix du sous-groupe parabolique minimal $\PO\subset P$.
On peut donc définir de telles bases pour toute paire
de sous-groupes paraboliques $P\subset\Q$ sans les supposer standard.

\begin{lemme}\label{inclus}
Si $P\subset\Q\subset\R$ sont trois sous-groupes paraboliques\textup, alors on a les inclusions
\[
\Delta_P^\Q\subset\Delta_P^\R\Qquad{et}\hDelta_\Q^\R\subset\hDelta_P^\R
\]
et il existe un sous-groupe parabolique $S$ tel que $P\subset S\subset\R$ et
\[
\Delta_P^S=\Delta_P^\R-\Delta_P^\Q\Qquad{et}\hDelta_S^\R=\hDelta_P^\R-\hDelta_\Q^\R\ptf
\]
\end{lemme}

\begin{proof}
La première assertion est claire; la seconde résulte du lemme~\ref{bijpar}.
\end{proof}

\begin{lemme}\label{binome}
Soient $P\subset\R$ deux sous-groupes paraboliques.
\[
\sum_{P\subset\Q\subset \R}(-1)^{a_P-a_\Q}=
\begin{cases}
1&\text{si $P=\R$,}\\
0& \text{sinon.}
\end{cases}
\]
\end{lemme}

\begin{proof}
Il suffit d'observer que d'après le lemme~\ref{bijpar}
la famille des sous-groupes paraboliques $\Q$ entre $P$ et $\R$ est en bijection avec la
famille des
sous-ensembles de $\Delta_P^\R$ puis d'invoquer la formule du binôme.
\end{proof}

On observera que si on identifie $\gaPQ$ et son dual au moyen de la structure euclidienne canonique,
les éléments de $\hDelta_P^\Q$
sont colinéaires aux éléments de la base duale de la base $\Delta_P^\Q$ et plus précisément ne
diffèrent que par des scalaires rationnels strictement positifs.
Dans la combinatoire des cônes,
les longueurs des vecteurs des bases ne jouent aucun rôle; seuls les angles importent.
On pourrait donc remplacer partout les $\hDelta_P^\Q$ par la base duale de $\Delta_P^\Q$,
mais il reste commode de penser aux éléments de $\hDelta_P^\Q$ comme des
restrictions de poids.

Les angles seront contrôlés via les deux lemmes bien connus ci-dessous.
Ils sont au cœur de la combinatoire qui commande toute la suite.

\begin{lemme}\label{obtua}
Considérons un espace vectoriel euclidien $V$ de dimension finie\textup,
muni d'une base obtuse $\Delta$ \cad que pour $\alpha\ne\beta$ dans $\Delta$
\[
\langle\alpha,\beta\rangle \le 0\ptf
\]
Soit $\Delta^1$ une partie de $\Delta$.
On note $\Delta_1$ la projection de $\Delta-\Delta^1$ sur l'orthogonal $V_1$
de $\Delta^1$. Alors $\Delta_1$ est une base obtuse de $V_1$.
\end{lemme}

\begin{proof}
Considérons trois
vecteurs distincts $\alpha$, $\beta$ et $\gamma$ appartenant à $\Delta$.
La projection $\bar\alpha$ de $\alpha $ sur l'orthogonal
de $\gamma$ s'écrit:
\[
\bar\alpha=\alpha-c_\alpha\,\gamma\qquad\text{avec 
$c_\alpha =\frac{\langle\alpha,\gamma\rangle}{\langle\gamma,\gamma\rangle}$}\ptf
\]
Mais, si $\bar\beta$ est la projection de $\beta$ sur l'orthogonal de
$\gamma$ on a
\[
\langle\beta,\gamma\rangle=\langle\bar\beta,\gamma\rangle
\]
et donc
\[
\langle\bar\alpha,\bar\beta\rangle =\langle\alpha,\beta\rangle
-\frac{\langle\alpha,\gamma\rangle \langle\beta,\gamma\rangle }{\langle\gamma,\gamma\rangle }
\le\langle\alpha,\beta\rangle\le0\ptf
\]
Le lemme résulte de cette remarque par récurrence sur le cardinal de $\Delta^1$.
\end{proof}

\begin{lemme}\label{obtub}
Considérons un espace vectoriel euclidien $V$ de dimension finie\textup,
muni d'une base obtuse $\Delta$.
Alors la base duale $\hDelta$ est une base aigüe\footnote{Le rédacteur principal a choisi d'écrire aigüe plutôt que aiguë (nonobstant la préférence du second rédacteur pour cette graphie traditionnelle),
suivant en cela les récentes recommandations du Conseil supérieur de la langue française.}
de $V$\textup: le produit scalaire $\langle\vpi,\vpi'\rangle $ est positif ou nul pour tout $\vpi$ et $\vpi'$ dans
$\hDelta$.
\end{lemme}

\begin{proof}
Soient $\vpi$ et $\vpi'$ deux vecteurs distincts dans la base duale $\hDelta$.
On désigne par $\alpha$ et $\alpha'$ les éléments de $\Delta$ correspondant à $\vpi$ et $\vpi'$.
Notons $\Delta^1$ le complémentaire de
$\{\alpha,\alpha'\}$ dans $\Delta$
et $V_1$ l'orthogonal de $\Delta^1$.
On observe que $\vpi$ et $\vpi'$ forment une base de $V_1$.
Notons enfin $\bar\alpha$ et $\bar\alpha'$
les projections de $\alpha$ et $\alpha'$ sur $V_1$.
D'après le lemme~\ref{obtua} l'ensemble $\{\bar\alpha,\bar\alpha'\}$ est une base
obtuse de $V_1$. Comme c'est la base duale de la base $\{\vpi,\vpi'\}$,
on est ramené à prouver le lemme en dimension~2, ce qui est élémentaire.
\end{proof}

\begin{lemme}\label{obtu}\footnote{Langlands donne un énoncé légèrement plus fort dans \cite{Eisen}*{Lemme~2.9, p.~20}.}
L'ensemble $\Delta_P^\Q$ est une base obtuse du dual de $\gaPQ$ et la base duale $\hDelta_P^\Q$ est aigüe.
\end{lemme}

\begin{proof} Supposons que $P\subset\Q$ sont deux sous-groupes paraboliques standard.
On sait que $\Delta_{\PO}^\Q$ est une base obtuse du dual de $\gao^\Q$
pour la structure euclidienne induite par la forme de Killing.
Il résulte alors du lemme~\ref{obtua} que $\Delta_P^\Q$,
qui est la projection de $\Delta_{\PO}^\Q-\Delta_{\PO}^P$ sur l'orthogonal
$\ga_P^*$ de $\ga_{\PO}^P$ est aussi obtuse.
De même la base des coracines est obtuse.
Maintenant le dual d'une base obtuse est une base aigüe
d'après le lemme~\ref{obtub}.
\end{proof}

\begin{lemme}\label{recone}
On suppose que $P\subset\Q$ sont deux sous-groupes paraboliques standard et on considère $H\in\gao$ tel que
\[
\alpha(H)>0\quad\forall\alpha\in\Delta_P^\Q \Qquad{et}\vpi(H)\le0\quad\forall\vpi\in\hDelta_{\PO}^P\ptf
\]
Alors
\[
\gamma(H)>0\qquad\forall\gamma\in\Delta_{\PO}^\Q-\Delta_{\PO}^P\ptf
\]
\end{lemme}

\begin{proof}
Par hypothèse, si on note $H_\Q$ la projection de $H$ sur $\ga_\Q$, on a
\[
H=\sum_{\vpi\in\hDelta_P^\Q} a_\vpi\vpi^\vee +\sum_{\beta\in\hDelta_{\PO}^P} b_\beta\beta^\vee+H_\Q
\]
avec $a_\vpi>0$ et $b_\beta\le0$.
Mais, pour $\gamma\in\Delta_{\PO}^\Q-\Delta_{\PO}^P$ on a $\gamma(\beta^\vee)\le0$
d'après le lemme~\ref{obtu}. Il reste à observer que $\gamma(H_\Q)=0$ et que
puisque $\gamma\notin\Delta_{\PO}^P$ alors $\gamma(\vpi^\vee)=1$
pour l'un des
\[
\vpi\in\hDelta_P^\Q\subset\hDelta_{\PO}^\Q
\]
alors que $\gamma(\vpi^\vee)=0$ pour tous les autres $\vpi^\vee$.
\end{proof}

\begin{lemme}\label{precone}
Soit $P$ et $\Q$ deux sous-groupes paraboliques. Si $\alpha(X)>0$ pour tout $\alpha\in\Delta_P^\Q$ on a
\[
\vpi(X)>0\qquad\text{pour tout $\vpi\in\hDelta_P^\Q$} \ptf
\]
\end{lemme}

\begin{proof}
Il suffit de montrer que tout $\vpi\in\hDelta_P^\Q$ peut s'écrire
\[
\vpi=\sum_{\alpha\in\Delta_P^\Q} c_{\alpha}\alpha
\qquad\text{avec $c_\alpha\ge0$ pour tout $\alpha\in\Delta_P^\Q$}\ptf
\]
Ceci résulte de ce que $\Delta_P^\Q$ et $\hDelta_P^\Q$
sont deux bases de $\ga_P^\Q$ et de ce que, pour tout $\alpha\in\Delta_P^\Q$,
si $\vpi_\alpha^\vee$ est l'élément de la base duale
correspondant à $\alpha$, alors
\[
c_\alpha=\vpi(\vpi_\alpha^\vee)\ge0
\]
car, d'après le lemme~\ref{obtu},
$\vpi_\alpha^\vee$ appartient à une base aigüe.
\end{proof}

\begin{lemme}\label{trecone}
Soient $P\subset\Q\subset\R$ trois sous-groupes paraboliques.
Supposons $\alpha(X)>0$ pour tout $\alpha\in\Delta_P^\R$.
Considérons $\bar\alpha\in\Delta_\Q^\R$ projection de
$\alpha\in\Delta_P^\R$ sur $\ga_\Q^\R$. Alors
\[
\bar\alpha(X)\ge\alpha(X) >0 \ptf
\]
\end{lemme}

\begin{proof}
On peut écrire $\bar\alpha$ sous la forme
\[
\bar\alpha=\alpha+\sum_{\beta\in\Delta_P^\Q}c_\beta\vpi_\beta
\qquad\text{avec $\vpi_\beta\in\hDelta_P^\Q$}
\]
et on doit avoir $\langle\beta,\bar\alpha\rangle =0$ pour $\beta\in\Delta_P^\Q$. Mais
\[
\langle\beta,\bar\alpha\rangle =\langle\alpha,\beta\rangle +c_\beta=0
\]
implique $c_\beta\ge0$ puisque $\Delta_P^\R$ est obtuse. On en déduit, compte tenu du lemme~\ref{precone}, que
\[
\bar\alpha(X)=\alpha(X)+\sum_{\beta\in\Delta_P^\Q}c_\beta\vpi_\beta(X)\ge
\alpha(X)\ptf\qedhere
\]
\end{proof}

Un élément $X\in\gao$ sera dit \og positif régulier\fg ou simplement \og régulier\fg
si
\[
\alpha(X)>0\qquad\forall\alpha\in\Delta_{\PO}^\G\ptf
\]
Nous utiliserons aussi la variante suivante:
on introduit le nombre
\[
\dPO(X)=\inf_{\alpha\in\Delta_{\PO}}\alpha(X)\ptf
\]
\newindex{dpo@$\dPO$}{distance T}%
Alors, $X$ est régulier si $\dPO(X)>0$.

\begin{lemme}\label{dist}
Soit $P$ un sous-groupe parabolique standard et considérons
$\bar\alpha\in\Delta_P$ qui est la projection de $\alpha\in\Delta_{\PO}$.
Soit $X\in\gao$ régulier.
On a
\[
\bar\alpha(X)\ge\alpha(X)\ge\dPO(X)\ptf
\]
\end{lemme}

\begin{proof}
C'est une conséquence immédiate du lemme~\ref{trecone}.
\end{proof}

\section{Géométrie et groupe de Weyl}

Le quotient du normalisateur de $\MO$ dans $\G$ par $\MO$
est le groupe de Weyl de $\G$ et sera noté $\weyl^\G$
ou simplement $\weyl$.
\newindex{W@$\weyl$}{groupes de Weyl}%
Si $P$ est un sous-groupe parabolique (semi-standard) de sous groupe de Levi $\M$
on notera souvent $\weyl^P$ au lieu de $\weyl^\M$ le groupe de Weyl de $\M$.
On notera $\ell(s)$ la longueur de $s\in\weyl$.

Soit $s\in\weyl$; on définit un sous-ensemble de $\Rac$ par
\[
\Rac(s)=\{\beta\in\Rac\mid \beta>0\text{ et }s(\beta)<0\}\ptf
\]
\newindex{R(s)@$\Rac(s)$}{racines s}%
Plus généralement, pour $s$ et $t$ dans $\weyl$ on pose
\[
\Rac(s,t)=\{\beta\in\Rac\mid t(\beta)>0\text{ et }s(\beta)<0\}\ptf
\]
\newindex{R(s,t)@$\Rac(s,t)$}{racines s t}%
On remarquera que
$\beta\mapsto t\beta$ induit une bijection $\Rac(s,t)\to\Rac(st\moins)$.

\begin{lemme}\label{bki} Considérons
$s=s_\alpha u$ avec $\ell(s)=\ell(u)+1$
où $s_\alpha$ est la symétrie définie par rapport à la racine simple $\alpha$\textup; alors
\[
\Rac(s)=\Rac(u)\cup\{\gamma\}
\]
avec $\gamma=u\moins(\alpha)$\textup; en particulier le cardinal de $\Rac(s)$ est la longueur de $s$.
Plus généralement\textup, posons $v=s t\moins $ et
supposons que $v=s_\alpha w$ avec $\ell(v)=\ell(w)+1$
où $s_\alpha$ est la symétrie définie par rapport à une racine simple $\alpha$.
Posons\textup, comme ci-dessus, $u=s_\alpha\moins s$ et $\gamma=u\moins\alpha$ alors
\[
\Rac(s,t)=\Rac(u,t)\cup\{\gamma\}\ptf
\]
\end{lemme}

\begin{proof}
La première assertion est un résultat classique que l'on trouve par exemple dans
\cite{Bki}*{Chapitre~VI, \parag 1, \no 6, Corollaire~2, p.~158}.
Pour le cas général on invoque les bijections $\Rac(u,t)\to\Rac(w)$
et $\Rac(s,t)\to\Rac(v)$ induites par $\beta\mapsto t\beta$
puis on remarque que $t(u\moins\alpha)=w\moins\alpha$.
\end{proof}

\begin{lemme}\label{swp}
Soit $P=\M\N$ un sous-groupe parabolique standard
et soit $s\in\weyl$.
Supposons que les racines $\beta\in\Rac(s)$
sont combinaison de racines simples $\alpha\in\Delta_{\PO}^P$ alors
$s$ appartient à $\weyl^P$.
\end{lemme}

\begin{proof}
Cette assertion s'obtient par récurrence sur la longueur de $s$. C'est clair pour $\ell(s)=0$.
Maintenant supposons que $s=s_\alpha\,t$ avec $\ell(s)=\ell(t)+1$
où $s_\alpha$ est la symétrie définie par rapport à la racine simple $\alpha$.
On a vu au lemme~\ref{bki} que
\[
\Rac(s)=\Rac(t)\cup\{\gamma\}
\]
avec $\gamma=t\moins(\alpha)$.
Par hypothèse de récurrence $t$ appartient au groupe de Weyl de $M$, le sous-groupe de Levi de $P$.
Comme $\gamma$ ne fait intervenir que des racines simples dans $\Delta_{\PO}^P$
et comme
\[
\alpha=t(\gamma)
\]
on a $\alpha\in\Delta_{\PO}^P$.
Donc $s_\alpha$ appartient aussi au groupe de Weyl de $M$ ainsi que $s=s_\alpha t$.
\end{proof}

Les lemmes suivants sont également classiques mais leur démonstration est souvent laissée en exercice\footnote{\Cf par exemple \cite{Bki}*{Chapitre~IV, \parag1, Exercice~3, p.~37}}.
Faute de référence commode, nous en donnons des preuves pour le confort du lecteur.

\begin{lemme}\label{echange}
Soit $P$ un sous-groupe parabolique standard.
Toute classe dans $\weyl/\weyl^P$ possède un unique représentant $s$ de longueur minimale
et\textup, pour tout $t\in\weyl^P$\textup, on a
\[
\ell(st)=\ell(s)+\ell(t)\ptf
\]
\end{lemme}

\begin{proof}
Soit $s$ un élément de longueur minimale dans sa classe et soit
$t\in\weyl^P$. Considérons des décompositions réduites de $s$ et $t$:
\[
s=s_1\dotsm s_p \Qquad{et} t=t_1\dotsm t_q\ptf
\]
Alors ou bien $\ell(st)=p+q$
ou bien il existe un plus petit indice $0\le r<q$ tel que
\[
\ell(s\,t_1\dotsm t_r)=p+r
\] et
\[
\ell(st_1\dotsm t_{r+1})=p+r-1\ptf
\]
Il résulte alors de la \og condition d'échange\fg
(\cf \cite{Bki}*{Chapitre~IV, \parag~1, Proposition~4, p.~15}) que l'on a soit\footnote{Dans ce qui suit la notation $\hat t_i$ signifie que $t_i$ est omis.}
\[
s\,t_1\cdots t_{r+1}=s\,t_1\dotsm \hat t_i\dotsm t_r
\]
ce qui est impossible puisque $t_1\cdots t_{r+1}$ est une décomposition réduite,
soit
\[
s\,t_1\dotsm t_{r+1}=s_1\dotsm \hat s_i \dotsm s_p\, t_1\dotsm t_r
\]
ce qui contredit la minimalité de la longueur de $s$ dans sa classe.
\end{proof}

Le lemme~\ref{echange} admet la généralisation suivante:

\begin{lemme}\label{echangebis}
Soit $P$ et $\Q$ deux sous-groupes paraboliques standard.
Toute classe dans $\weyl^{P}\bs\weyl/\weyl^{\Q}$ possède un unique représentant de longueur minimale.
\end{lemme}

\begin{proof}
Soient $s$ et $\sigma $ deux éléments de longueur minimale dans la même classe. On a donc
$\sigma =ust$ avec $u\in\weyl^{P}$ et
$t\in\weyl^{\Q}$ et supposons de plus que $t$ et $u$ sont choisis de sorte que $q=\ell(t)$ soit minimal.
Considérons des décompositions réduites de $s$, $t$ et $u$:
\[
s=s_1\dotsm s_p,\qquad t=t_1\dotsm t_q\Qquad{et}u=u_r\dotsm u_1\ptf
\]
Comme $s$ est minimal dans sa double classe il résulte du lemme~\ref{echange} que
\[
\ell(st)=\ell(s)+\ell(t)=p+q\ptf
\]
Si nous supposons $r\ge1$, il existe un indice $k<r$ tel que
\[
\ell(u_{k}\dotsm u_1\,st)=p+q+k
\]
et
\[
\ell(u_{k+1}\dotsm u_1\,st)=\ell(u_k\dotsm u_1\,st)-1\ptf
\]
Il résulte alors de la \og condition d'échange\fg que, puisque $u_{k+1}\dotsm u_1$ est une
décomposition réduite, alors on a
\[
u_{k+1}\dotsm u_1\,st=u_{k}\dotsm u_1s't'
\]
avec soit
\[
s'\,t'=s'\,t=s_1\dotsm \hat s_i \dotsm s_p \,t
\]
ce qui contredit la minimalité de la longueur de $s$ dans sa classe,
soit, $q\ge1$ et
\[
s'\,t'=s\,t'=s\,t_1\dotsm \hat t_i\dotsm t_q
\]
et donc si on pose
\[
u'=u_r\dotsm\hat u_{k+1}\dotsm u_1
\]
on aura
\[
\sigma =ust=u'st'
\]
avec $\ell(t')=\ell(t)-1$
ce qui contredit la minimalité de $q$. On a donc $r=0$ et
comme $\ell(\sigma)=\ell(s)$ on aura aussi $q=0$
et $\sigma=s$.
\end{proof}

\begin{lemme}\label{weylg}
Soit $P$ un sous-groupe parabolique standard.
Tout $s\in\weyl/\weyl^{P}$ admet un unique représentant\textup, encore noté
$s$\textup, dans $\weyl$ satisfaisant l'une des conditions
équivalentes suivantes\textup:
\begin{enumerate}[(i)]
\item $s$ est de longueur minimale
dans sa classe à gauche modulo $\weyl^{P}$;
\item $s\alpha>0$ pour toute $\alpha\in\Delta_{\PO}^{P}$.
\end{enumerate}
\end{lemme}

\begin{proof}\hspace*{-\labelsep}\footnote{Ce lemme est la proposition~3.9 de \cite{BT2}.}\hspace*{\labelsep}%
D'après le lemme~\ref{echange}, dans toute classe à gauche modulo $\weyl^P$
il existe un unique élément $s\in\weyl$ de longueur minimum
et la longueur de $s$ est le nombre de racines $\beta$ positives,
appartenant au système de racines réduit $\Rac$,
telles que $s(\beta)$ soit négatif (\cf lemme~\ref{bki}).
Considérons cet élément $s$ et supposons qu'il existe une racine $\alpha\in\Delta_{\PO}^{P}$
avec $s\alpha<0$; comme la symétrie $s_\alpha$, relative à cette racine simple,
ne change pas le signe des racines positives autres que
$\alpha$ on en déduit que
\[
\ell(ss_\alpha)=\ell(s)-1
\]
ce qui contredit la minimalité de $\ell(s)$.
La condition~(i) implique donc~(ii). Maintenant on observe que, puisque $\weyl^P$ agit trivialement sur $\ga_P$,
le signe de $t(\beta)$ pour $\beta\in\Rac$
est indépendant de $t\in\weyl^P$ si
la projection de $\beta$ sur $\ga_P$ est non nulle.
La condition~(ii), lorsqu'elle est réalisée, permet donc de minimiser la longueur de~$s$.
\end{proof}

Soient $P$ et $\Q$ deux sous-groupes paraboliques (semi-standard). On note
\[
\weyl(\ga_P,\ga_\Q)
\]
\newindex{W(ap,aQ)@$\weyl(\ga_P,\ga_\Q)$}{groupes de Weyl a p q}%
l'ensemble des restrictions à $\ga_P$ des
$s\in\weyl$ tels que
\[
s(\ga_P)=\ga_\Q\ptf
\]
C'est un sous-ensemble de $\weyl/\weyl^{P}\ptf$
On dit que deux sous-groupes paraboliques standard
$P$ et $\Q$ sont associés si $\weyl(\ga_P,\ga_\Q)$ est non vide.

\begin{lemme}\label{weylga}
Supposons $P$ et $\Q$ standard.
Tout $s\in\weyl(\ga_P,\ga_\Q)$ admet un unique représentant\textup, encore noté
$s$\textup, dans $\weyl$ satisfaisant l'une des conditions
équivalentes suivantes\textup:
\begin{enumerate}[(i)]
\item $s$ est de longueur minimale
dans sa classe à gauche modulo $\weyl^P$;
\item$s$ est de longueur minimale
dans sa classe à droite modulo $\weyl^\Q$;
\item $s\alpha>0$ pour toute $\alpha\in\Delta_{\PO}^{P}$;
\item $s\moins\alpha>0$ pour toute $\alpha\in\Delta_{\PO}^{\Q}$;
\item $s(\Delta_{\PO}^{P})=\Delta_{\PO}^{\Q}$.
\end{enumerate}
\end{lemme}

\begin{proof} L'équivalence de~(i) et~(ii) est claire.
L'existence de $s$ et l'équivalence de~(i) et~(iii) est l'objet du lemme~\ref{weylg}.
L'équivalence de~(ii) et~(iv) se démontre de même en changeant $s$ en $s\moins$.
La condition~(iii) nous dit que $s(\Delta_{\PO}^{P})$, qui est une base pour
le système de racines du sous-groupe de Levi $\M_\Q$ de $\Q$,
est formé de racines positives; c'est donc $\Delta_{\PO}^{\Q}$.
Donc~(iii) implique~(v). Maintenant~(v) implique évidemment~(iii) et~(iv).
\end{proof}

Soient maintenant $P$ et $\R$ deux sous-groupes paraboliques (semi-standard). On note
\[
\weyl(\ga_P,\R)
\]
\newindex{W(ap,R)@$\weyl(\ga_P,\R)$}{WeylPR}%
l'ensemble des doubles classes dans
\[
\weyl^{\R}\bs\weyl/\weyl^P
\]
formées d'éléments $s\in\weyl$ tels que $s(\ga_P)\supset\ga_\R$.
Le lemme~\ref{weylga} admet la généralisation suivante
(\cf \cite{MW}*{assertion~(3), p.~93}):\footnote{Ce lemme fait partie des exercices laissés à la lectrice dans \cite{MW}.}

\begin{lemme}\label{weylgab}
Si $P$ et $\R$ sont standard\textup,
l'ensemble $\weyl(\ga_P,\R)$ est
en bijection avec l'ensemble des
$s\in\weyl$ tels que
\begin{enumerate}[(i)]
\item $s(\ga_P)\supset\ga_\R$;
\item $s\moins\alpha>0$ pour toute $\alpha\in\Delta_{\PO}^\R$.
\end{enumerate}
Pour un tel $s$ on a
\[
s(\Delta_{\PO}^P)=\Delta_{\PO}^\Q\subset\Delta_{\PO}^\R
\]
où $\Q$ est un sous-groupe parabolique standard dans $\R$. L'ensemble
$\weyl(\ga_P,\R)$ est en bijection avec l'union disjointe des quotients
\[
\weyl^\R(\ga_\Q,\ga_\Q)\bs\weyl(\ga_P,\ga_\Q)
\]
où $\Q$ parcourt les sous-groupes paraboliques standard dans $\R$\textup, modulo $\M_\R$\hyph association.
\end{lemme}

\begin{proof}
La condition $s(\ga_P)\supset\ga_\R$ équivaut à dire que
le sous-groupe de Levi $\M_\R$ de $\R$ contient $s(\M)$.
L'ensemble $\weyl(\ga_P,\R)$ est donc formé de
doubles classes d'éléments $s$ tels que
\[
s(\weyl^P)\subset\weyl^\R\ptf
\]
C'est donc aussi le sous-ensemble des classes
dans $\weyl^{\R}\bs\weyl^\G$ d'éléments vérifiant~(i).
Maintenant, d'après le lemme~\ref{weylg},
tout élément de $\weyl(\ga_P,\R)$ admet un unique représentant
vérifiant~(ii), à savoir l'élément de longueur minimale dans
sa classe. Avec ce choix de $s$ l'ensemble de racines
$s\moins(\Delta_{\PO}^\R)$ est une base du système de racines de $s\moins(\M_\R)$
formé de racines positives pour l'ordre induit par l'ordre sur les racines de $\G$.
L'ensemble $\Delta_{\PO}^P$ est inclus dans l'ensemble
des racines de $s\moins(\M_\R)$.
Comme les racines dans $\Delta_{\PO}^P$ sont des racines simples (pour $\G$)
elles sont \emph{a fortiori} simples dans le système de racines de $s\moins(\M_\R)$ avec l'ordre induit
et donc
\[
\Delta_{\PO}^P\subset s\moins(\Delta_{\PO}^\R)
\]
ce qui équivaut à
\[
s(\Delta_{\PO}^P)\subset\Delta_{\PO}^\R\ptf
\]
Donc $s(\Delta_{\PO}^P)$ est une base pour les racines d'un sous-groupe de Levi
standard $\Q$. La dernière assertion en résulte.
\end{proof}

\section{Chambres et facettes}\label{chambresfaces}

Soit $\M$ un sous-groupe de Levi et $\Q$ un sous-groupe parabolique
contenant $\M$. L'ensemble des sous-groupes paraboliques $P\subset\Q$ et admettant $\M$
comme sous-groupe de Levi sera noté
\[
\Parab^\Q(\M)\ptf
\]
\newindex{PQ(M)@$\Parab^Q(M)$}{ensemble de prarab}%
L'ensemble des sous-groupes paraboliques $P$ avec
$\M\subset P\subset\Q$
sera noté
\[
\FF^\Q(\M)\ptf
\]
\newindex{FQ(M)@$\FF^Q(M)$}{sur-ensemble de prarab}%
Enfin, on désigne par
\[
\Levi^\Q(\M)
\]
\newindex{LQ(M)@$\Levi^Q(M)$}{levig}%
l'ensemble des sous-groupes de Levi $L$
avec $\M\subset L\subset\Q$. On omettra souvent l'exposant $\Q$ lorsque $\Q=\G$.

Soit $P$ un sous-groupe parabolique standard et soit $\M$ son sous-groupe de Levi.
On note
\[
\weyl^\G(\ga_{\M})\Qquad{ou simplement} \weyl(\ga_{\M})
\]
\newindex{W(aM)@$\weyl(\ga_\M)$}{WeylaM}%
l'union (disjointe) de tous les $\weyl(\ga_P,\ga_\R)$
où $\R$ est un sous-groupe parabolique standard de $\G$.
Cet ensemble est en bijection avec un sous-ensemble du quotient
$\weyl^\G/\weyl^{\M}$ que l'on identifie à un sous-ensemble de $\weyl^\G$
en choisissant le représentant de longueur minimale.
D'après le lemme~\ref{weylga},
$\weyl(\ga_{\M})$ est l'ensemble des $s\in\weyl^\G$
tels que $s(\Delta_{\PO}^{P})\subset\Delta_{\PO}$.
On notera
\[
\weyl^\G(\M)\Qquad{ou simplement}\weyl(\M)
\]
le groupe quotient du groupe $N^\weyl(\M)$ des $s\in\weyl$
tels que $s(\M)=\M$ par le sous-groupe $\weyl^\M$. On peut identifier
$\weyl(\M)$ avec un sous groupe de $\weyl$: à tout $s\in N^\weyl(\M)$
on associe $\bar s$ le représentant de longueur minimale
dans la classe
\[
s\weyl^\M=\weyl^\M s=\weyl^\M s\weyl^\M\ptf
\]
Une variante de la preuve du lemme~\ref{echangebis} montre que $s\mapsto\bar s$
est un homomorphisme $N^\weyl(\M)\to\weyl$ dont le noyau est $\weyl^\M$.
Son image est donc isomorphe à $\weyl(\M)$.
On notera $n(\M)$
le cardinal de $\weyl(\ga_\M)$ et $w(\M)$ le cardinal de $\weyl(\M)$.
Alors $n(\M)/w(\M)$, le cardinal du quotient $\weyl(\ga_\M)/\weyl(\M)$,
est le nombre de sous-groupes paraboliques standard associés à $P$.

On dispose dans $\ga_\M$ des chambres de Weyl complémentaires
des hyperplans définis par les racines $\beta$ dans $\Rac$
qui ne sont pas identiquement nulles sur à $\ga_\M$.
Soit $s\in\weyl(\ga_{\M})$, on note $\Delta(\M,s)$
\newindex{Delta(M,s)@$\Delta(\M,s)$}{deltams}%
l'ensemble
des projections sur $\ga_\M^*$ des $\alpha\in
s\moins(\Delta_{\PO})$ qui ne sont pas identiquement nulles sur $\ga_\M$.
On définit une chambre $C_\M(s)$ dans $\ga_\M$ par les inégalités
\[
\alpha(H)>0\qquad\text{pour $\alpha\in\Delta(\M,s)$}\ptf
\]
Plus généralement on dispose des facettes lorsqu'on remplace les inégalités par
des égalités pour les $\alpha$ appartenant au sous-ensemble de racines
associé à un sous-groupe de Levi contenant $\M$ (\cf \cite{Bki}*{Chapitre~V, \parag1}).

\begin{lemme}\label{weylbija}
Il y a une bijection naturelle entre les trois ensembles suivants
\begin{enumerate}[(i)]
\item L'ensemble $C_\M$ des chambres de Weyl dans $\ga_\M$
\item L'ensemble $\Parab(\M)$ des sous-groupes paraboliques $P$ admettant $\M$
comme sous-groupe de Levi
\item L'ensemble $\weyl(\ga_\M)$.
\end{enumerate}
\end{lemme}

\begin{proof}
Considérons $s\in\weyl(\ga_\M)$. On lui associe la chambre $C_\M(s)$.
On remarque que $C_\M(s)$ est une facette d'une chambre de Weyl dans $\gao$
et que les chambres qui admettent $C_\M(s)$ comme facette forment une orbite sous $\weyl^\M$.
L'application
\[
\weyl(\ga_\M)\to C_\M
\]
est donc injective. Maintenant, toute chambre dans $\ga_\M$ est une
facette d'une chambre dans $\gao$.
Comme le groupe de Weyl est simplement transitif sur l'ensemble
des chambres dans $\gao$ l'application ci-dessus est surjective.
Enfin à $\Delta(\M,s)$ on associe le
sous-groupe parabolique $\Q_s\in\Parab(\M)$
tel que
\[
\Delta_{\Q_s}=\Delta(\M,s)\ptf
\]
On en déduit la bijection entre~(ii) et~(iii).
\end{proof}

\begin{lemme}\label{weylbijb}
Il y a une bijection naturelle
entre les deux ensembles suivants
\begin{enumerate}[(i)]
\item Les facettes dans $\ga_\M$
\item L'ensemble $\FF(\M)$ des sous-groupes paraboliques $P$ contenant $\M$.
\end{enumerate}
\end{lemme}

\begin{proof}
On observe que
\[
\FF(\M)=\bigcup_{L\in\Levi(\M)}\Parab(L)
\]
et l'assertion résulte alors du lemme~\ref{weylbija}.
\end{proof}

Soit $s\in\weyl(\ga_{\M})$ alors $s(\M)$
est le sous-groupe de Levi d'un sous-groupe parabolique standard
que nous noterons $\R_s$.
Soit $\R^s$ le sous-groupe parabolique standard
dont le sous-groupe de Levi admet comme racines simples les $\alpha\in\Delta_{\PO}$
avec $s\moins\alpha>0$.
On pose
\[
\Q_s=s\moins(\R_s)\Qquad{et}\Q^s=s\moins(\R^s)\ptf
\]
\newindex{Qs@$\Q_s$}{qsqs}%
\newindex{Qs@$\Q^s$}{qsqs}%
En particulier $\Q_s$ est le sous-groupe parabolique dans $\Parab(\M)$
associé à $s$ par le lemme~\ref{weylbija}.
Il résulte du lemme~\ref{weylga} que $\R_s\subset\R^s$ et donc $\Q_s\subset\Q^s$.
On pose
\[
\FF_s(\M)=\{\Q\mid \Q_s\subset\Q\subset\Q^s\}\ptf
\]

\begin{lemme}\label{soma}
L'ensemble $\FF(\M)$ des sous-groupes paraboliques $\Q$ contenant $\M$\textup,
est l'union disjointe\textup, indexée par $s\in\weyl(\ga_{\M})$\textup, des
$\FF_s(\M)$.
En d'autres termes\textup, si $f$ est une fonction sur $\FF(\M)$ on a\textup:
\begin{equation}
\sum_{\Q\in\FF(\M)}f(\Q)=
\sum_{s\in\weyl(\ga_{\M})}\,\sum_{\Q\in\FF_s(\M)}f(\Q )\ptf\label{eq1.1}
\end{equation}
En particulier\textup, si $f(\Q)$ ne dépend que de $s$ lorsque $\Q\in\FF_s(\M)$\textup, alors
\begin{equation}
\sum_{\Q\in\FF(\M)}(-1)^{a_\Q-a_P}f(\Q)=f(\widebar{P})\label{eq1.2}
\end{equation}
où $\widebar{P}$ est le sous-groupe parabolique
qui correspond à l'opposé de la chambre de Weyl positive dans $\ga_\M$.
\end{lemme}

\begin{proof}
Étant donné $\Q$ on lui associe l'élément $s\in\weyl$ de longueur minimale
tel que $s(\Q)$ soit standard. C'est un élément de $\weyl(\ga_\M)$.
On a alors
\[
Q_s\subset\Q\subset\Q^s\ptf
\]
L'assertion~\eqref{eq1.1} est ainsi établie.
Maintenant, on observe que $\Q_s=\Q^s$ si et seulement si $s$ est
l'élément qui est associé, par le lemme~\ref{weylbija}, à l'opposé de la chambre positive dans $\ga_\M$.
L'assertion~\eqref{eq1.2} résulte alors de~\eqref{eq1.1} et du lemme~\ref{binome}.
\end{proof}

Nous utiliserons aussi la variante suivante:

\begin{lemme}\label{somb}
Soit $P$ un sous-groupe parabolique standard
de sous-groupe de Levi $\M$.
Soit $g(s,\R)$ une fonction dépendant de $s\in\weyl$
et d'un sous-groupe parabolique standard $\R$ et telle que
\[
g(ts,\R)=g(s,\R)\qquad\text{pour $t\in\weyl^\R$}\ptf
\]
On a
\[
\sum_{\PO\subset\R}\,\sum_{s\in\weyl(\ga_P,\R)}g(s,\R)=
\sum_S\sum_{s\in\weyl(\ga_P,\ga_S)}\,\sum_{s\moins(\R)\in\FF_s(\M)} g(s,\R)
\]
la somme en $S$ au second membre portant sur les paraboliques standard associés à~$P$.
\end{lemme}

\begin{proof}
On observe que tout $\Q\in\FF(\M)$ s'écrit d'une façon et d'une seule
sous la forme $\Q= s\moins(\R)$ avec $s\in\weyl(\ga_Q,P)$
et on rappelle que $\weyl(\ga_\M)$ est l'union disjointe des $\weyl(\ga_P,\ga_S)$
où $S$ décrit l'ensemble des sous-groupes paraboliques standard.
Le lemme résulte alors de l'assertion~\eqref{eq1.1} du lemme~\ref{soma}.
\end{proof}

\section{Familles orthogonales}\label{famorth}

On appelle famille orthogonale la donnée d'une fonction sur le groupe de Weyl
\[
\XX\colon s\mapsto X_s
\]
à valeurs dans $\gao$, telle que si $s=s_\alpha\,t$
où $s_\alpha$ est la symétrie définie par rapport à la racine simple $\alpha$ alors
\[
X_t- X_s= b_\gamma(s,t)\,\gamma^\vee
\]
avec $b_\gamma(s,t)\in\CM$ et $\gamma=t\moins(\alpha)$. On observera que
\[
\Rac(s,t)=\{\gamma\}
\Quad{et}\Rac(t,s)=\{-\gamma\}
\Qquad{et donc} b_\gamma(s,t)=b_\gamma(t,s)
\]
ne dépend que de la paire $\{s,t\}$.
En d'autres termes
si $s$ et $t$ définissent des chambres
adjacentes, alors $X_t-X_s$ est orthogonal au mur
séparant les deux chambres.
On dit qu'une famille orthogonale est régulière si
$b_\gamma(s,t)>0$ pour toute paire $\{s,t\}$ définissant des chambres adjacentes.

\begin{lemme}\label{ws}
Soit $\XX$ une famille orthogonale.
Soient $s$ et $t$ deux éléments du groupe de Weyl\textup,
il existe des scalaires $b_\beta(s,t)$\textup, dépendant du choix d'une décomposition réduite de
$v=st\moins$\textup, tels que
\[
X_t-X_s=\sum_{\beta\in\Rac(s,t)}b_\beta(s,t)\,\beta^\vee\ptf
\]
Si $\XX$ est régulière les coefficients $b_\beta(s,t)$
sont strictement positifs.
\end{lemme}

\begin{proof}
La preuve se fait par récurrence sur la longueur de $v=s t\moins $.
Si $s=t$ l'assertion est triviale. Maintenant supposons $v=s_\alpha w$ avec $\ell(v)=\ell(w)+1$
et $s_\alpha$ la symétrie définie par rapport à une racine simple $\alpha$.
Posons $u=s_\alpha s$ et donc $ut\moins=w$.
Par hypothèse de récurrence, et compte tenu du cas
particulier de la longueur 1 où l'assertion n'est autre que la propriété de définition
des familles orthogonales, on a
\[
X_t-X_s=(X_t-X_u)+(X_u-X_s)=
\sum_{\beta\in\Rac(u,t)}b_\beta(u,t)\,\beta^\vee
+b_\gamma(s,u)\,\gamma^\vee
\]
avec $\gamma=u\moins(\alpha)$. On posera
\[
b_\beta(s,t)=b_\beta(u,t)\quad\text{pour $\beta\in \Rac(u,t)$}
\Qquad{et} b_\gamma(s,t)=b_\gamma(s,u)\ptf
\]
Pour conclure on observe que, d'après le lemme~\ref{bki}, on a
\[
\Rac(s,t)=\Rac(u,t)\cup\{\gamma\}\ptf
\]
On aurait aussi pu procéder ainsi: si $v=s_{n}\dotsm s_{1}$ est une décomposition réduite de $v=st\moins$
et si on pose
\[
v_i=s_i\dotsm s_{1},\quad t_i=v_it
\Qquad{et} \Rac(t_{i},t_{i-1})=\{\beta_i\}
\]
on a $t_0=t$ et $t_n=s$ et donc
\[
X_t-X_s=\sum_{i=1}^{i=n} X_{t_{i-1}}-X_{t_{i}}=\sum_{i=1}^{i=n} b_{\beta_i}(t_{i},t_{i-1})\beta_i^\vee\ptf
\]
Enfin il convient d'observer que $\Rac(s,t)$ est l'ensemble de ces $\beta_i$.
\end{proof}

Un premier exemple de famille orthogonale est fourni par le lemme suivant:

\begin{lemme}\label{wT}
Soit $s$ un élément du groupe de Weyl et soit
$T\in\gao$. Alors pour chaque $\beta\in \Rac(s)$ il existe une
forme linéaire
\[
T\mapsto c_\beta(s,T)
\]
strictement positive sur la chambre de Weyl positive\textup,
de sorte que
\[
(1-s\moins)T=
T-s\moins(T)=\sum_{\beta\in \Rac(s)}c_\beta(s,T)\,\beta^\vee\ptf
\]
De plus on a
\[
c_\beta(s,T)\ge \dPO(T)\ptf
\]
En particulier $s\mapsto s\moins\T$ est une famille orthogonale.
Elle est régulière si $\T$ est régulier.
\end{lemme}

\begin{proof}
Supposons $s=s_\alpha t$ avec $l(s)=l(t)+1$
et $s_\alpha$ la symétrie définie par rapport à la racine simple $\alpha$.
On observe que
\[
t\moins T-s\moins T=t\moins(T-s_\alpha\moins T)=\alpha(T)\,t\moins(\alpha^\vee)
\]
Le lemme résulte alors du lemme~\ref{ws} pour
\[
s\mapsto s\moins T\ptf\qedhere
\]
\end{proof}

Au lieu des éléments du groupe de Weyl
on pourra utiliser les sous-groupes paraboliques minimaux pour indexer les
éléments d'une famille orthogonale: à tout $s\in\weyl$ on associe
le sous-groupe parabolique minimal $P=s\moins\PO$ et on écrira $X_P$ pour $X_s$.
Ceci a l'avantage de fournir une indexation indépendante du choix de $\PO$.
La condition d'orthogonalité pour la famille $\XX$ est équivalente
à demander que si $P$ et $\Q$ sont deux sous-groupes paraboliques
minimaux adjacents alors $X_P$ et $X_\Q$ ont la même projection $X_\R$ sur le mur
séparant les chambres associées à $P$ et $\Q$ \cad sur $\ga_\R$ où
$\R$ est le sous-groupe parabolique engendré par $P$ et $\Q$.
Plus généralement si $\R$ est un sous-groupe parabolique et si
$P\subset\R$
on note $X_\R$ la projection de $X_P$ sur $\ga_\R$. Cette projection est indépendante
du choix de $P$.

Nous aurons besoin de la généralisation suivante: si $\M$ est un sous-groupe de Levi,
on appellera famille $\M$\hyph orthogonale la donnée d'une famille d'éléments $X_\Q\in\ga_\Q$
pour chaque $\Q\in\FF(\M)$ telle que si $P\in\FF^\Q(\M)$
alors la projection de $X_P$ sur $\ga_\Q$ soit $X_\Q$. Il suffit bien entendu de se donner
les $X_P$ pour $P\in\Parab(\M)$.
\'Etant donnés une famille orthogonale $\XX$
et $\M$ un sous-groupe de Levi standard, on définit une
famille $\M$\hyph orthogonale
en considérant les projections sur $\ga_\M$
des $X_s$ pour $s\in\weyl(\ga_\M)$ (\cf lemme~\ref{weylbija}).

On notera $\gHM$ l'espace vectoriel des familles $\M$\hyph orthogonales.
C'est la limite projective des $\ga_P$
sur l'ensemble des $P\in\FF(\M)$ muni de l'ordre inverse de celui défini par l'inclusion:
\[
\gHM=\varprojlim_{P\in\FF(\M)}\ga_P
\]
avec pour flèches les projections
\[
\pi_{P,\Q}\colon \ga_P\to\ga_\Q
\]
lorsque $P\subset\Q$.
On dispose donc de projections
\[
\pi_P\colon \gHM\to\ga_P
\]
indexées par les sous-groupes paraboliques $P\in\FF(\M)$
et il est facile de voir qu'elles sont surjectives (par exemple en utilisant le lemme~\ref{wT}).

\section{Enveloppes convexes de familles orthogonales}\label{envconvfamorth}

Nous rappelons maintenant des résultats
établis par Arthur dans \cite{ArInvent} et
qui généralisent des lemmes combinatoires
empruntés aux travaux de Langlands sur les séries d'Eisenstein.

Soit $\XX=\{X_s\}$ une famille orthogonale dans $\gao$ et soit $\M$ un sous-groupe de Levi
standard.
Soient $s$ et $t$ dans $\weyl$ tels que les chambres $C(s)$ et $C(t)$
dans $\gao$ associées à $s$ et $t$ admettent
des facettes définissant la même chambre dans $\ga_\M$ et donc
diffèrent par un élément de $\weyl^\M$. Maintenant le lemme~\ref{ws} montre que $X_t-X_s$ est
orthogonal à $\ga_\M$ et donc que $X_s$ et $X_t$
ont la même projection sur $\ga_\M$.

\begin{lemme}\label{adja}
Soient deux chambres $C(s)$ et $C(t)$
dans $\gao$ associées à $s$ et $t$ définissant des chambres adjacentes $C_\M(s)$ et $C_\M(t)$ dans
$\ga_\M$. Alors la projection de $X_t-X_s$ sur $\ga_\M$ est
orthogonale au mur séparant les chambres.
\end{lemme}

\begin{proof}
Les deux chambres $C_\M(s)$ et $C_\M(t)$ dans $\ga_\M$ étant adjacentes il existe une forme linéaire sur cet espace,
unique à un scalaire près, qui est positive sur l'une et négative sur l'autre.
Soit $\lambda$ une telle forme linéaire séparant $C_\M(s)$ et $C_\M(t)$.
On sait que (lemme~\ref{ws})
\[
X_t-X_s=\sum_{\beta\in\Rac(s,t)}b_\beta(s,t)\,\beta^\vee
\]
où $\Rac(s,t)$ est l'ensemble des racines
telles que $t(\beta)>0$ et $s(\beta)<0$.
On a donc $\bar\beta=c_\beta\lambda$ si $\bar\beta$ est la projection de $\beta$ sur $\ga_\M$
pour tout $\beta\in\Rac(s,t)$.
\end{proof}

Soit $\chm$ dans une chambre de $\ga_\M^\G$.
On définit $\phi_{\M,s}^\chm$
\newindex{phiMskappa@$\phi_{M,s}^\chm$}{phiMsKappa}
comme la fonction caractéristique des $H\in\gao$ tels que
\[
\vpi_\alpha(H)\le0\quad\text{si $\alpha(\chm)>0$}
\]
et
\[
\vpi_\alpha(H)>0\quad\text{si $\alpha(\chm)<0$}
\]
où les $\alpha$ parcourent $\Delta(\M,s)$ (cet ensemble est introduit juste avant le lemme~\ref{weylbija}).
On note $a(s,\chm)$ le nombre de $\alpha\in\Delta(\M,s)$ \newindex{ask@$a(s,\chm)$}{asssk}%
avec $\alpha(\chm)<0$\footnote{Nous avons suivi Langlands \cite{MS}*{Lecture~15} et Arthur en définissant les fonctions $\phi_{\M,s}^\chm$
et les nombres $a(s,\chm)$ au moyen d'un paramètre $\chm$ dans une chambre. Mais bien entendu,
seul le choix de la chambre importe pour ces définitions.}.
On introduit alors%
\newindex{GammaMQ(H,X,kappa)@$\Gamma_\M^\Q(H,\XX,\chm)$}{gammaMQkappa}%
\[
\Gamma_\M(H,\XX,\chm)=\sum_{s\in\weyl(\ga_{\M})}(-1)^{a(s,\chm)}\phi_{\M,s}^\chm(H-X_s)\ptf
\]

\begin{lemme}\label{supportconva}
Supposons que, pour tout $s$ on ait
\[
\vpi_\alpha(H-X_s)\le0
\]
pour tout $\alpha\in\Delta(\M,s)$. Alors $\Gamma_\M(H,\XX,\chm)=1$.
\end{lemme}

\begin{proof}
Le point $\chm$ étant fixé dans une chambre, soit $s\in\weyl(\ga_{\M})$ l'élément du groupe de Weyl associé. On a
\[
(-1)^{a(s,\chm)}\phi_{\M,s}^\chm(H-X_s)=1
\]
alors que $\phi_{\M,t}^\chm(H-X_t)=0$ si $s\ne t$.
\end{proof}

\begin{lemme}\label{supportconvb}
Soit $s\in\weyl(\ga_{\M})$.
Supposons $\XX$ régulière et soit $\chm\in C_\M(s)$.
Alors
\begin{enumerate}[(i)]
\item Le point $X_s$ appartient au support de
$\Gamma_\M(\bullet,\XX,\chm)$
\item Si $H$ appartient au support de
$\Gamma_\M(\bullet,\XX,\chm)$ on a
\[
\langle\chm,H-X_s\rangle \le0\ptf
\]
\end{enumerate}
\end{lemme}

\begin{proof} On observe que
\[
\phi_{\M,s}^\chm(0)=1\ptf
\]
Par ailleurs on sait d'après le lemme~\ref{ws} que puisque $\YY$ est régulière $X_s-X_t$ est une combinaison à coefficients
strictement positifs des coracines $\beta^\vee$ telles que $s(\beta)>0$ et $t(\beta)<0$.
L'ensemble $\Rac(t,s)$ de ces racines est non vide
si $t \ne s$; elles sont positives ou nulles sur la chambre $C_\M(s)$. Elles ne sont pas
toutes nulles sur cette chambre si $s$ et $t$ définissent des chambres distinctes
dans $\ga_\M$.
On aura donc $\vpi_\alpha(\beta^\vee)>0$ pour au moins un $\alpha\in\Delta(\M,s)$
et un $\beta\in\Rac(t,s)$, ce qui implique
\[
\vpi_\alpha(X_s-X_t)>0
\]
alors que $\alpha(\chm)>0$ et donc
\[
\phi_{\M,t}^\chm(X_s-X_t)=0\ptf
\]
On a ainsi établi~(i).
Pour prouver~(ii) on remarque que
si $\Gamma_\M(H,\XX,\chm)$ est non nul, alors il y a au moins un $t$ tel que
\[
\phi_{\M,t}^\chm(H-X_t)=1
\]
et donc
\begin{equation}
\langle\chm,H-X_t\rangle =\sum \alpha(\chm)\vpi_\alpha(H-X_t)\le0\ptf\label{eq1.1a}
\end{equation}
Pour conclure, il reste à observer comme ci-dessus que $X_t-X_s$ est une combinaison à coefficients
positifs des coracines $\beta^\vee$ telles que $t(\beta)>0$ et $s(\beta)<0$
et donc on aura $\vpi_\alpha(\beta^\vee)\le0$ pour toute $\alpha\in\Delta(\M,s)$ alors que $\alpha(\chm)>0$
ce qui implique
\begin{equation}
\langle\chm,X_t-X_s\rangle \le0\ptf\label{eq1.2a}
\end{equation}
et la conjonction de~\eqref{eq1.1a} et~\eqref{eq1.2a} implique~(ii).
\end{proof}

\begin{lemme}\label{indep}\footnote{Cet énoncé se trouve déjà pour l'essentiel dans \cite{Bould}*{\parag8, p.~246}.}
La fonction $\Gamma_\M(H,\XX,\chm)$ est indépendante de $\chm$.
\end{lemme}

\begin{proof}
Il suffit de montrer que si
$\sigma$ et $\tau$ définissent des chambres adjacentes et si $\chm_\sigma$ et $\chm_\tau$
appartiennent aux chambres $C_\M(\sigma)$ et $C_\M(\tau)$ respectivement alors
\[
\Gamma_\M(H,\XX,\chm_\sigma)=\Gamma_\M(H,\XX,\chm_\tau)\ptf
\]
Examinons les différents termes. Tout d'abord on observe que
\[
(-1)^{a(s,\chm_\sigma)}\phi_{\M,s}^{\chm_\sigma}= (-1)^{a(s,\chm_\tau)}\phi_{\M,s}^{\chm_\tau}
\]
si
\[ \alpha(\chm_\sigma)\alpha(\chm_\tau)>0\qquad\forall\alpha\in\Delta(\M,s)\ptf
\]
Reste à examiner les contributions des autres $s$.
Soit $\lambda$ une forme linéaire définissant le mur entre les deux chambres $C_\M(\sigma)$ et $C_\M(\tau)$
et soient $s$ et $t$ deux éléments de $\weyl$ associés à deux chambres adjacentes
$C_\M(s)$ et $C_\M(t)$ séparées par ce mur. On observe qu'il existe dans $\Delta(\M,s)$ et dans
$\Delta(\M,t)$ une unique racine proportionnelle à $\lambda$. On voit que
\[
\xi_s(H)=\phi_{\M,s}^{\chm_\sigma}(H-X_s)+\phi_{\M,s}^{\chm_\tau}(H-X_s)
\]
est la fonction caractéristique des $H$ tels que
\[
\vpi_\alpha(H-X_s)\le0\qquad\text{si $\alpha(\chm_\sigma)>0$}
\]
et
\[
\vpi_\alpha(H-X_s)>0\qquad\text{si $\alpha(\chm_\sigma)<0$}
\]
pour les $\alpha\in\Delta(\M,s)$ qui ne sont pas proportionnels à $\lambda$.
Ce sont ceux pour lesquels on a $\alpha(\chm_\sigma)\alpha(\chm_\tau)>0$.
Mais pour chaque tel $\alpha\in\Delta(\M,s)$ il existe un unique $\tbeta\in\Delta(t)$
ayant la même projection (non nulle) sur le mur. Pour un tel couple on a
$\vpi_\alpha=\vpi_\beta$ et ils sont orthogonaux à $\lambda$. On a alors
\[
\vpi_\alpha(H-X_s)=\vpi_\beta(H-X_t)
\]
En effet, la projection de $X_s-X_t$ est proportionnelle à $\lambda^\vee$ d'après le lemme~\ref{adja}.
On a donc
\[
\xi_s(H)=\xi_t(H)
\]
et l'assertion en résulte.
\end{proof}

\begin{proposition}\label{envconv}
Supposons que la famille orthogonale $\XX$ est régulière.
Alors\textup, la fonction
\[
H\mapsto \Gamma_\M(H,\XX,\chm)
\]
est la fonction caractéristique de l'ensemble des $H$ dont la projection sur
$\ga_\M^\G$ appartient à l'enveloppe convexe des projections des $X_s$
avec $s\in\weyl(\ga_{\M})$.
\end{proposition}

\begin{proof}
Il résulte des lemmes~\ref{supportconvb}(ii) et~\ref{indep} que le support de cette fonction
est contenu dans l'ensemble des $H$ tels que
pour tout $s\in\weyl(\ga_{\M})$ on ait
\[
\vpi_\alpha(H-X_s)\le0
\]
pour tout $\alpha\in\Delta(\M,s)$.
Il résulte alors des lemmes~\ref{supportconva} et~\ref{indep}
que c'est la fonction caractéristique de cet ensemble, qui est un convexe fermé.
Maintenant, si $H$ n'appartient
pas à l'enveloppe convexe des projections des $X_s$ il existe un cône ouvert non vide
dans $\ga_\M^\G$
tel que pour $\chm$ dans ce cône on ait
\[
\langle T,H-X_s\rangle >0
\]
pour tout $s\in\weyl(\ga_{\M})$ et donc, d'après
le lemme~\ref{supportconvb}(ii), $H$ n'appartient pas au support.
Mais par ailleurs le lemme~\ref{supportconvb}(i)
montre que les projections des $X_s$ avec $s\in\weyl(\ga_{\M})$
appartiennent à ce support.
\end{proof}

\section{Combinatoire des cônes}
\label{conepart}

Soient $P\subset\Q$ deux sous-groupes paraboliques (semi-standard).
On note $\tau_P^\Q$
\newindex{tauPQ@$\tau_P^Q$}{topq}%
\newindex{tauPQ@$\protect\htau_P^Q$}{topq}%
la fonction caractéristique du cône ouvert dans $\gao$
défini par $\Delta_P^\Q$:
\[
\tau_P^\Q(H)=1 \iff\alpha(H)>0\quad \forall \alpha\in\Delta_P^\Q
\]
et on note $\htau_P^\Q$ la fonction caractéristique du cône ouvert défini par $\hDelta_P^\Q$:
\[
\htau_P^\Q(H)=1 \iff\varpi(H)>0\quad \forall\varpi\in\hDelta_P^\Q\ptf
\]
Lorsque $\Q=\G$ nous écrirons souvent $\tau_P$ au lieu de $\tau_P^\G$.
On observera que la valeur de $\tau_P^\Q(H)$ et de $\htau_P^\Q(H)$ ne dépendent que
de la projection $H_P^\Q$ de $H$ sur le sous-espace vectoriel~$\gaPQ$

\begin{lemme}\label{cone}
Si $P\subset\Q$\textup, on a
\[
\tau_P^\Q\le \htau_P^\Q\ptf
\]
Plus généralement\textup, si $P\subset\Q\subset\R$\textup,
\[
\tau_P^\Q\htau_\Q^\R\le\htau_P^\Q\htau_\Q^\R\le\htau_P^\R\ptf
\]
\end{lemme}

\begin{proof}
La première assertion n'est autre que le lemme~\ref{precone}.
Montrons la seconde assertion.
D'après ce qui vient d'être montré on a
\[
\tau_P^\Q\htau_\Q^\R\le\htau_P^\Q\htau_\Q^\R\ptf
\]
Il suffit maintenant de montrer que
\[
\htau_P^\Q\htau_\Q^\R\le\htau_P^\R\ptf
\]
Pour cela on observe que si $H$ est tel que
\[
\htau_P^\Q(H)\htau_\Q^\R(H)=1
\]
on a
\[
\vpi(H)>0\qquad\forall\vpi\in\hDelta_\Q^\R
\]
et d'autre part
\[
\widebar\vpi(H)>0\qquad\forall\vpi\in\hDelta_P^\R-\hDelta_\Q^\R
\]
où $\widebar\vpi$ est la projection orthogonale de $\vpi$ sur le dual de $\gaPQ$
\cad que
\[
\vpi=\widebar\vpi+\sum_{\alpha\in\Delta_P^\R-\Delta_P^\Q}\lambda_\alpha\,\vpi_\alpha
\qquad\text{avec $\lambda_\alpha+\widebar\vpi(\alpha^\vee)=0$}\ptf
\]
Mais, on a vu dans la preuve de la première assertion du lemme que
\[
\widebar\vpi\in\hDelta_P^\Q
\]
est une combinaison
à coefficients positifs des racines dans $\Delta_P^\Q$
et comme $\Delta_P^\R$ est une base obtuse on a $\widebar\vpi(\alpha^\vee)\le0$ et donc
$\lambda_\alpha\ge0$. Il en résulte que, comme escompté,
$\vpi(H)>0$ pour $\vpi\in\hDelta_P^\R-\hDelta_\Q^\R$.
\end{proof}

On définit une matrice dont les coefficients sont indexés par des paires
de sous-groupes paraboliques standard
\[
\tau=(\tau_{P,\Q})
\]
avec
\[
\tau_{P,\Q}=\begin{cases}
(-1)^{a^{}_P}\tau_P^\Q&\text{si $P\subset\Q$}\\
0& \text{sinon}
\end{cases}
\]
et de même
\[
\htau=(\htau_{P,\Q})
\]
avec
\[
\htau_{P,\Q}=\begin{cases}
(-1)^{a^{}_P}\htau_P^\Q&\text{si $P\subset\Q$}\\
0&\text{sinon.}
\end{cases}
\]

Un des clefs essentielles pour la combinatoire des sous-groupes paraboliques,
est la proposition suivante, appelée parfois \og Lemme combinatoire de Langlands\fg
(\cf \cite{MS}*{Corollary~13.1.2}):

\begin{proposition}\label{totoun}
Les matrices $\tau$ et $\htau$ sont inverses l'une de l'autre.
\end{proposition}

\begin{proof}
Il suffit de montrer que
\[
\sum_{P\subset\Q\subset\R}(-1)^{a_P-a_\Q}\tau_P^\Q\,\htau_\Q^\R
=\begin{cases}
1& \text{si $P=\R$}\\
0& \text{sinon.}
\end{cases}
\]

Considérons $H\in\gao$. Il existe un unique sous-groupe parabolique $S$
avec $P\subset S\subset\R$ tel que $\alpha\in\Delta_P^S$ équivaut à $\alpha(H)>0$
pour $\alpha\in\Delta_P^\R$. De même il existe un unique sous-groupe parabolique $T$
avec $P\subset T\subset\R$ tel que $\vpi\in\hDelta_T^\R$ équivaut à $\vpi(H)>0$
pour $\vpi\in\hDelta_P^\R$. On a alors
\[
\sum_{P\subset\Q\subset\R}(-1)^{a_P-a_\Q}\tau_P^\Q(H)\,\htau_\Q^\R(H)=
\sum_{T\subset\Q\subset S}(-1)^{a_P-a_\Q}
\]
et d'après le lemme~\ref{binome}
cette dernière expression est nulle sauf si $\T=S$ auquel cas elle vaut $(-1)^{a_P-a_S}$.
Maintenant, si elle est non nulle cela implique que $\alpha(H)>0$ pour $\alpha\in\Delta_P^S$
et $\vpi(H)>0$ pour $\vpi\in\hDelta_S^\R$. Compte tenu du lemme~\ref{cone}
il en résulte que $\vpi(H)>0$ pour tout
$\vpi\in\hDelta_P^\R$. Or on doit avoir $\vpi(H)\le0$ pour $\vpi\in\hDelta_P^\R-\hDelta_S^\R$
donc $S=T=P$. Par définition de $S$ on a $\alpha(H)\le0$
pour tout $\alpha\in\Delta_P^\R$, puisque $S=P$.
Par ailleurs $\vpi(H)>0$ pour tout $\vpi\in\hDelta_P^\R$
(car $T=P$). Ces deux propriétés sont contradictoires sauf si $P=\R$.
\end{proof}

On considère trois sous-groupes paraboliques $P\subset\Q\subset\R$ et on pose
\[
\phi_P^{\Q,\R} =\sum_{P\subset S\subset\Q}(-1)^{a_S-a_\Q}\,\htau_S^\R\ptf
\]
\newindex{phiPQR@$\phi_P^{\Q,\R}$}{phipqr}

\begin{lemme}\label{prepar}
La fonction
$\phi_P^{\Q,\R}$ est la fonction caractéristique des $H\in\gao$ tels que $\vpi(H)\le0$
pour tous les $\vpi\in\hDelta_P^\R-\hDelta_\Q^\R$ et $\vpi(H)>0$
pour tous les $\vpi\in\hDelta_\Q^\R$.
\end{lemme}

\begin{proof}
Considérons $H\in\gao$ et soit $P'$ le sous-groupe parabolique
avec $P\subset P'\subset\Q$ tel que $\vpi\in\hDelta_{P'}^\R$ équivaut à $\vpi\in\hDelta_P^\R$
et $\vpi(H)>0$. Donc
\[
\phi_P^{\Q,\R}(H) =\sum_{P'\subset S\subset\Q}(-1)^{a_S-a_\Q}=\begin{cases}
1&\text{si $P'=\Q$}\\
0& \text{sinon}
\end{cases}
\]
d'après le lemme~\ref{binome}. En particulier
\[
\phi_P^{\Q,\R}(H)\ne0\iff P'=\Q\ptf\qedhere
\]
\end{proof}

Nous allons relier les fonctions $\phi_P^{\Q,\R}$ et les fonctions
$\phi_{\M,s}^\chm$ introduites dans la section~\ref{envconvfamorth},
au moyen des sous-groupes paraboliques $\Q_s$, $\Q^s$ et de l'ensemble $\FF_s(\M)$
introduits dans le lemme~\ref{soma}.

\begin{lemme}\label{preparb}
Soit $\M$ un sous-groupe de Levi standard.
Lorsque $\chm$ appartient à la chambre de Weyl positive on a
\[
\phi_{\M,s}^\chm (H)=\phi_{\Q_s}^{\Q^s,\G} (H)
\]
et
\[
(-1)^{a(s,\chm)}\phi_{\M,s}^\chm (H)=
\sum_{\Q\in\FF_s(\M)}(-1)^{a_\Q-a_\G}\,\htau_\Q(H)\ptf
\]
\end{lemme}

\begin{proof}
On observe tout d'abord que, pour $\chm$ régulier, par définition de $\phi_{\M,s}^\chm$ et
d'après le lemme~\ref{prepar},
\[
\phi_{\M,s}^\chm (H)=\phi_{\Q_s}^{\Q^s,\G} (H)=\sum_{\Q\in\FF_s(\M)}
(-1)^{a_\Q-a_{\Q^s}}\,\htau_\Q(H)\ptf
\]
De plus
\[
a(s,\chm)=a_{\Q^s}-a_\G\ptf\qedhere
\]
\end{proof}

On notera $\phi_P^\Q$ la fonction $\phi_P^{\Q,\Q}$.
\newindex{phiPQ@$\phi_P^Q$}{phipq}%
Les fonctions $\phi_P^\Q$ donnent naissance à des partitions.

\begin{lemme}\label{partition}
Deux sous-groupes paraboliques $P$ et $\R$ étant fixés on a
\[
\sum_{\{\Q\mid P\subset\Q\subset\R\}}\phi_P^\Q\,\tau_\Q^\R\equiv1\ptf
\]
\end{lemme}

\begin{proof}
On a
\[
\sum_{P\subset\Q\subset\R}\phi_P^\Q\,\tau_\Q^\R=
\sum_{P\subset S\subset\Q\subset\R}(-1)^{a_S-a_\Q}\htau_S^\Q
\,\tau_\Q^\R=\sum_{P\subset S\subset\Q\subset\R}\htau_{S,\Q}\,\tau_{\Q,\R}
\]
Mais, d'après la proposition~\ref{totoun} on a $\htau\tau=1$ et donc, en notant
$\delta_{S,\R}$ le symbole de Kronecker, on a
\[
\sum_{P\subset\Q\subset\R}\phi_P^\Q\,\tau_\Q^\R=
\sum_{P\subset S\subset\R}\delta_{S,\R}=1\ptf\qedhere
\]
\end{proof}

\section{Cônes et convexes}

Soient $H$ et $X$ deux éléments de $\gao$ et
$P\subset\Q\subset\R$ trois sous-groupes paraboliques.
On notera $C(P,\Q,\R,X)$ l'ensemble des $H$ qui vérifient les inégalités
suivantes
\[
\alpha(H)>0\quad\forall \alpha\in\Delta_P^\Q\Qquad{et}
\alpha(H)\le0\quad\forall \alpha\in\Delta_P^\R-\Delta_P^\Q
\]
ainsi que
\[
\vpi(H-X)>0\quad\forall \vpi\in\hDelta_\Q^\R\Qquad{et}
\vpi(H-X)\le0\quad\forall \vpi\in\hDelta_P^\R-\hDelta_\Q^\R
\]
On remarque que pour $P$, $\R$ et $X$ fixés, les $C(P,\Q,\R,X)$ sont disjoints.

\begin{lemme}\label{convcomp}
L'ensemble
$C(P,\Q,\R,X)$ est un convexe dont la projection dans $\ga_P^\R$ est
d'adhérence compacte. Plus précisément\textup, il existe $c>0$ tel que l'on ait
\[
\lVert H_P^\R\rVert \le c \lVert X_P^\R\lVert \qquad\text{pour $H\in C(P,\Q,\R,X)$}\ptf
\]
De plus, si $P$ est standard et si $X$ est dans l'adhérence de la chambre de Weyl positive (en particulier si $X$ est régulier), $C(P,\Q,\R,X)$ est vide sauf si $\Q=\R$. Enfin si $X=0$ alors $C(P,\Q,\R,X)$ est vide si $P\ne\R$.
\end{lemme}

\begin{proof}
D'après le lemme~\ref{inclus} il existe un sous-groupe parabolique $S$ tel que
\[
\Delta_P^S=\Delta_P^\R-\Delta_P^\Q
\Qquad{et}\hDelta_S^\R=\hDelta_P^\R-\hDelta_\Q^\R\ptf
\]
L'espace $\ga_S^\R$, dont le dual a pour base $\hDelta_S^\R\subset\hDelta_P^\R$,
est l'orthogonal des $\beta\in\Delta_P^S$.
L'espace $\ga_\Q^\R$ qui est l'orthogonal des $\alpha\in\Delta_P^\Q$
a un dual qui admet pour base $\hDelta_\Q^\R$.
Les sous-espaces $\ga_S^\R$ et $\ga_\Q^\R$ ne sont pas en général orthogonaux
mais on a cependant une décomposition en somme directe:
\[
\ga_P^\R=\ga_S^\R\oplus\ga_\Q^\R\ptf
\]
Il suffit de prouver que les projections orthogonales de $C(P,\Q,\R,X)$
sur $\ga_S^\R$ et $\ga_\Q^\R$ sont relativement compactes.
Considérons
\[ H\in C(P,\Q,\R,X)
\]
et notons $H_1$ et $X_1$ (resp. $H_2$ et $X_2$)
les projections orthogonales sur $\ga_S^\R$ (resp. $\ga_\Q^\R$)
de $H$ et $X$.
On a par hypothèse
\[
\vpi(H_1-X_1)=\vpi(H-X)\le0\quad\forall\vpi\in\hDelta_S^\R=\hDelta_P^\R-\hDelta_\Q^\R
\]
et donc
\[
\vpi(H_1)\le \vpi(X_1)\quad\forall\vpi\in\hDelta_S^\R\ptf
\]
Soit $\bar\alpha\in\Delta_S^\R$ la projection de $\alpha\in\Delta_P^\Q$ sur
$\ga_S^\R$.
On a donc $\beta^\vee(\bar\alpha)=0$ pour $\beta\in\Delta_P^S$.
On dispose d'une bijection $\beta\mapsto\vpi_\beta$ entre $\Delta_P^S$ et $\hDelta_\Q^\R$.
Les $\vpi_\beta\in\hDelta_\Q^\R$ forment une base de $\ga_\Q^\R$ qui est un supplémentaire de $\ga_S^\R$
dans $\ga_P^\R$.
On peut donc écrire $\bar\alpha$ sous la forme
\[
\bar\alpha=\alpha+\sum_{\beta\in\Delta_P^S}\mu_\beta
\vpi_\beta\qquad\text{avec $\mu_\beta=-\alpha(\beta^\vee)\ge0$}\ptf
\]
On rappelle que par hypothèse $\alpha(H)>0$ pour $\alpha\in\Delta_P^\Q$
et $\vpi_\beta(H-X)>0$ pour $\beta\in\Delta_P^S$ puisque
dans ce cas on a $\vpi_\beta\in\hDelta_\Q^\R$.
On en déduit que
\[
\bar\alpha(H_1)=\bar\alpha(H)=\bar\alpha(H-X)+\bar\alpha(X_1)
\ge\alpha(H-X)+\bar\alpha(X_1) > \alpha(X_1-X)
\]
pour tout $\alpha\in\Delta_P^\Q$.
Il résulte du lemme~\ref{cone} qu'il existe une constante $c(X)$ avec $\vpi(H_1)>c(X)$
pour tout $\vpi\in\hDelta_S^\R$, soit compte tenu de ce qui précède,
\[
c(X)\le \vpi(H_1)\le \vpi(X_1)\quad\forall\vpi\in\hDelta_S^\R\ptf
\]
Donc $H_1$ reste dans un compact de $\ga_S^\R$.
La discussion pour $H_2$ est analogue et on obtient
que l'on a simultanément
\[
\vpi(H_2)>\vpi(X_2)\quad \forall\vpi\in\hDelta_\Q^\R
\Qquad{et}
\bar\alpha(H_2)\le \alpha(X_2-X)\quad\forall
\bar\alpha\in\Delta_\Q^\R \ptf
\]
On conclut comme ci-desus que $H_2$ reste dans un compact.
De plus, si $X$ est dans l'adhérence de la chambre de Weyl positive, on a
\[
\bar\alpha(H_2)\le\bar\alpha(X_2)\qquad\forall\bar\alpha\in\Delta_\Q^\R
\]
ce qui, d'après le lemme~\ref{cone}, implique
\[
\vpi(H_2)\le\vpi(X_2)\qquad \forall\vpi\in\hDelta_\Q^\R\ptf
\]
Comme, par ailleurs, on a vu que
\[
\vpi(H_2)>\vpi(X_2)\qquad \forall\vpi\in\hDelta_\Q^\R
\]
ces d'inégalités sont incompatibles si $\Q\ne\R$ et
la seconde assertion du lemme en découle. Dans le cas $X=0$
on obtient de plus des inégalités
\[
\vpi(H_1)\le0\qquad \forall\vpi\in\hDelta_S^\R\Qquad{et}
\bar\alpha(H_1)>0\qquad\forall\bar\alpha\in\Delta_S^\R
\]
sur le sous-espace $\ga_S^\R$ incompatibles, d'après le lemme~\ref{cone},
s'il est non nul.
\end{proof}

On considère la matrice
$\Gamma(H,X)=\{\Gamma_{P,\Q}(H,X)\}$
définie par
\[
\Gamma(H,X)=\tau(H)\htau(H-X)
\]
et on pose
\[
\Gamma_P^\Q(H,X)=(-1)^{a_P-a_\Q}\Gamma_{P,\Q}(H,X)\ptf
\]
\newindex{GammaPQ(H,X)@$\Gamma_P^\Q(H,X)$}{gammaPR}%
On a donc
\[
\Gamma_P^\R(H,X)=\sum_{P\subset\Q\subset\R}(-1)^{a_\Q-a_\R}\,\tau_P^\Q(H)\htau_\Q^\R(H-X)\ptf
\]

\begin{lemme}\label{GammaHX}
On a
\begin{gather}
\tau_P^\R(H)=\sum_{P\subset\Q\subset\R}\Gamma_P^\Q(H,X)\tau_\Q^\R(H-X)\label{eq1.1b}\\
\htau_P^\R(H-X)=\sum_{P\subset\Q\subset\R}(-1)^{a_\Q-a_\R}\,
\htau_P^\Q(H)\Gamma_\Q^\R(H,X)\label{eq1.2b}
\end{gather}
et
\begin{equation}
\Gamma_P^\R(H,X+Y)=\sum_{P\subset\Q\subset\R}\Gamma_P^\Q(H,X)\Gamma_\Q^\R(H-X,Y)\ptf\label{eq1.3b}
\end{equation}
\end{lemme}

\begin{proof}
En effet, d'après la proposition~\ref{totoun}, on a
les égalités matricielles
\[
\tau(H)=\Gamma(H,X)\tau(H-X),\quad\htau(H-X)=\htau(H)\Gamma(H,X)
\]
et
\[
\Gamma(H,X+Y)=\Gamma(H,X)\Gamma(H-X,Y) \ptf\qedhere
\]
\end{proof}

\begin{lemme}\label{Gammacar}
La fonction
\[
H\mapsto\Gamma_P^\R(H,X)
\]
est combinaison à coefficients
\[
(-1)^{a_\Q-a_\R}
\]
des fonctions caractéristiques
des ensembles $C(P,\Q,\R,X)$. En particulier la projection
de son support dans $\ga_P^\R$ est compacte. Plus précisément\textup, il existe $c>0$ tel que l'on ait
\[
\lVert H_P^\R\rVert\le c\lVert X_P^\R\rVert
\]
lorsque $H$ appartient à ce support.
Lorsque $P$ est standard et
$X$ est régulier\textup, c'est la fonction caractéristique des $H$ tels que
$\alpha(H)>0$ et $\vpi_\alpha(H-X)\le0$ pour tous les $\alpha\in\Delta_P^\R$\textup:
\[
\Gamma_P^\R(H,X)=\tau_P^\R(H)\phi_P^\R(H-X)\ptf
\]
\end{lemme}

\begin{proof}
Fixons $H$ et $X$. Soit $S$ le plus grand sous-groupe parabolique
tel que $\tau_P^S(H)=1$ et $\T$ le plus petit sous-groupe parabolique
tel que $\htau_\T^\R(H-X)=1$; alors, comme dans la proposition~\ref{totoun}, on voit que
la somme sur $\Q$
\[
\Gamma_P^\R(H,X)=\sum_{\T\subset\Q\subset S}(-1)^{a_\Q-a_\R}
\]
est nulle sauf si $S=\T$ et, dans ce cas, un tel $H$ appartient à $C(P,\Q,\R,X)$.
La compacité résulte alors du lemme~\ref{convcomp}.
Lorsque $X$ est régulier il résulte également du lemme~\ref{convcomp} que
$\Gamma_P^\R(H,X)$ est la fonction caractéristique de $C(P,\R,\R,X)$
ce qui établit la dernière assertion.
\end{proof}

Toujours sous l'hypothèse $X$ régulier, l'égalité~\eqref{eq1.1b} du lemme~\ref{GammaHX}
peut donc s'interpréter
comme une partition du cône associé à $\tau_P^\R$ en produits de cônes
par des convexes relativement compacts. C'est une variante du lemme~\ref{partition}.

Soit maintenant $\M$ un sous-groupe de Levi et $\Q$ un sous-groupe parabolique
contenant $\M$. On rappelle que $\FF^\Q(\M)$ est l'ensemble des sous-groupes paraboliques de
$\Q$ qui contiennent $\M$.
On associe à une famille $\M$\hyph orthogonale $\XX=\{X_P\}$ la fonction suivante:
\[
\Gamma_\M^\Q(H,\XX)=\sum_{P\in\FF^\Q(\M)}(-1)^{a_P-a_\Q}\,
\htau_P^\Q(H-X_P)\ptf
\]
\newindex{GammaMQ(H,X)@$\Gamma_\M^\Q(H,\XX)$}{gammaMQ}%
On notera $\delta_\M^\Q$ la fonction caractéristique du sous-espace des $H\in\gao$
tels que
\[
H_\M^\Q=0
\]
où $H_\M^\Q$
est la projection de $H$ sur $\ga_\M^\Q$; autrement dit, $\delta_\M^\Q$ est
la fonction caractéristique
du sous-espace $\gao^\M\oplus\ga_\Q$.

\begin{lemme}\label{decomp}
On a les identités suivantes\textup:
\begin{gather}
\sum_{\Q\in\FF^\R(\M)}\delta_\M^\Q(H)\tau_\Q^\R(H)\equiv1\label{eq1.1c}\\
\Gamma_\M^\R(H,\XX)=\sum_{\Q\in\FF^\R(\M)}\delta_\M^\Q(H)\Gamma_\Q^\R(H,X_\Q)
\label{eq1.2c}\\
\sum_{\Q\in\FF^\R(\M)}\Gamma_\M^\Q(H,\XX)\,\tau_\Q^\R(H-X_\Q)\equiv 1\ptf
\label{eq1.3c}
\end{gather}
\end{lemme}

\begin{proof}
L'assertion~\eqref{eq1.1c} est l'écriture, au moyen de fonctions caractéristiques,
de la décomposition de l'espace vectoriel $\ga_\M^\R$
en chambres et facettes attachées aux divers sous-groupes paraboliques.
Maintenant, par définition de $\Gamma_\Q^\R$ on voit que
\[
\sum_{\Q\in\FF^\R(\M)} \delta_\M^\Q(H)\Gamma_\Q^\R(H,X_\Q)
= \mspace{-5mu}\sum_{\Q\in\FF^\R(\M)}\,\sum_{\Q\subset P\subset\R} (-1)^{a_P-a_\R}\delta_\M^\Q(H)
\tau_\Q^P(H)\htau_P^\R(H-X_P)
\]
qui, d'après~\eqref{eq1.1c} est encore égal à
\[
\sum_{P\in\FF^\R(\M)}(-1)^{a_P-a_\R}\htau_P^\R(H-X_P) =\Gamma_\M^\R(H,\XX)
\ptf
\]
Ceci établit l'assertion~\eqref{eq1.2c}. On en déduit que
\[
\sum_{\Q\in\FF^\R(\M)}\Gamma_\M^\Q(H,\XX)\,\tau_\Q^\R(H-X_\Q)=
\sum_{P\in\FF^\R(\M)}\delta_\M^P(H)\Gamma_P^\Q(H,X_S)\tau_\Q^\R(H-X_\Q)
\]
qui est encore égal à
\[
\sum_{P\in\FF^\R(\M)}\delta_\M^P(H)\sum_{P\subset \U\subset\Q}(-1)^{a_\U-a_\Q}
\tau_P^\U(H)\htau_\U^\Q(H-X_\U)\tau_\Q^\R(H-X_\Q)
\]
Mais, compte tenu de la proposition~\ref{totoun},
\[
\sum_\Q(-1)^{a_\U-a_\Q}\,
\htau_\U^\Q(H-X_\U)\tau_\Q^\R(H-X_\Q)
\]
est nul sauf si $U=\R$ auquel cas cette somme vaut identiquement~1.
On obtient donc
\[
\sum_{\Q\in\FF^\R(\M)}\Gamma_\M^\Q(H,\XX)\,\tau_\Q^\R(H-X_\Q)=
\sum_{P\in\FF^\R(\M)}\delta_\M^P(H)\tau_P^\R(H)
\]
et l'assertion~\eqref{eq1.3c} résulte alors de~\eqref{eq1.1c}.
\end{proof}

\begin{corollaire}\label{compact}
La fonction
\[
H\mapsto \Gamma_\M^\Q(H,\XX)
\]
a un support dont la projection sur $\ga_\M^\Q$ est
compacte.
Plus précisément\textup, il existe $c>0$ tel que l'on ait
\[
\lVert H_\M^\Q\rVert \le c\sup_{P\in\Parab^\Q(\M)}\lVert X_P^\Q\rVert
\]
lorsque $H$ appartient à ce support.
\end{corollaire}

\begin{proof}
C'est une conséquence immédiate de l'équation~\eqref{eq1.2c} du lemme~\ref{decomp} et du lemme~\ref{Gammacar}.
\end{proof}

On pourra observer que d'après le lemme~\ref{convcomp}
\[
\Gamma_\Q^\R(H,0)\equiv
\begin{cases}
0&\text{si $\Q\ne\R$}\\
1&\text{sinon}
\end{cases}
\]
et donc, compte tenu de l'équation~\eqref{eq1.2c} du lemme~\ref{decomp},
on a
\[
\delta_\M^\Q(H)=\Gamma_\M^\Q(H,0)\ptf
\]

\begin{lemme}\label{convol}
Soient $\XX$ et $\YY$ deux familles orthogonales.
On a
\[
\Gamma_\M^\R(H,\XX+\YY)=\sum_{\Q\in\FF^\R(\M)}
\Gamma_\M^\Q(H,\XX)\Gamma_\Q^\R(H-X_\Q, Y_\Q)
\]
\end{lemme}

\begin{proof}
Par définition de $\Gamma_\M^\R$ on a:
\[
\Gamma_\M^\R(H,\XX+\YY)=\sum_{P\in\FF^\R(\M)}(-1)^{a_P- a_\R}\,\htau_P^\R(H-X_P-Y_P)
\]
soit encore, au vu de~l'égalité~\eqref{eq1.2b} du lemme~\ref{GammaHX},
\[
\sum_{P\in\FF^\R(\M)}\sum_{P\subset\Q\subset\R}
(-1)^{a_P- a_\Q}\,\htau_P^\Q(H-X_P)\Gamma_\Q^\R(H-X_\Q,Y_\Q)
\ptf
\]
Ce qui peut se récrire
\[
\sum_{\Q\in\FF^\R(\M)}\sum_{P\in\FF^\Q(\M)}
(-1)^{a_P- a_\Q}\,\htau_P^\Q(H-X_P)\Gamma_\Q^\R(H-X_\Q,Y_\Q)
\ptf
\] On conclut en utilisant la définition de $\Gamma_\M^\Q$.
\end{proof}

Pour alléger les notations nous nous limiterons dans ce qui suit au cas $\Q=\G$.
Soit $\M$ un sous-groupe de Levi standard.
On a introduit et étudié plus haut
(\cf lemme~\ref{supportconva} à proposition~\ref{envconv}), pour $\chm$ en dehors des murs dans $\ga_\M$,
des nombres et des fonctions
\[
a(s,\chm),\qquad \phi_{\M,s}^\chm\Qquad{et} \Gamma_\M(H,\XX,\chm)\ptf
\]
Pour alléger les notations on écrira $a(s)$ et $\phi_{\M,s}$
pour $a(s,\chm)$ et $\phi_{\M,s}^\chm$
lorsque $\chm$ est dans la chambre positive.
On a vu au lemme~\ref{indep}
que la fonction $\Gamma_\M(H,\XX,\chm)$ était indépendante de $\chm$.
Nous allons de plus montrer qu'elle coïncide avec $\Gamma_\M(H,\XX)$.

\begin{proposition}\label{preconv}
Soit $\XX$ une famille orthogonale. Avec les notations du lemme~\ref{soma}\textup, la fonction $\Gamma_\M(H,\XX)$ vérifie l'identité\textup:
\begin{equation}
\Gamma_\M(H,\XX)=\sum_{s\in\weyl(\ga_{\M})}\,\sum_{\Q\in\FF_s(\M)}(-1)^{a_\Q-a_\G}\,\htau_\Q(H-X_s )\label{eq1.1d}
\end{equation}
ainsi que
\begin{equation}
\Gamma_\M(H,\XX)=\sum_{s\in\weyl(\ga_{\M})}(-1)^{a(s)}\phi_{\M,s}(H-X_s)\ptf\label{eq1.2d}
\end{equation}
Plus généralement\textup, on a
\begin{equation}
\Gamma_\M(H,\XX)=\Gamma_\M(H,\XX,\chm)\label{eq1.3d}
\end{equation}
pour tout $\chm$ en dehors de murs. Lorsque la famille orthogonale $\XX$ est régulière la fonction
\[
H\mapsto \Gamma_\M(H,\XX)
\]
est la fonction caractéristique de l'ensemble des $H$ dont la projection sur $\ga_\M^\G$ appartient à l'enveloppe convexe des $X_P$ pour $P\in\Parab(\M)$.
\end{proposition}

\begin{proof}
Par définition
\[
\Gamma_\M(H,\XX)=\sum_{\Q\in\FF(\M)}(-1)^{a_\Q-a_\G}\,
\htau_\Q(H-X_\Q)\ptf
\]
Considérons $\Q\in\FF_s(\M)$ pour $s\in\weyl(\ga_\M)$. On observe que $X_\Q$ est la projection de
\[
X_s=X_{\Q_s}
\]
sur $\ga_\Q$. L'équation~\eqref{eq1.1d} résulte alors de la première assertion du lemme~\ref{soma}.
Maintenant on rappelle que par définition
\[
\Gamma_\M(H,\XX,\chm)=\sum_{s\in\weyl(\ga_{\M})}(-1)^{a(s,\chm)}\phi_{\M,s}^\chm(H-X_s)\ptf
\]
Supposons que $\chm$ appartient à la chambre de Weyl positive. On observe alors que, d'après lemme~\ref{preparb}, on a:
\[
(-1)^{a(s,\chm)}\phi_{\M,s}^\chm (H)=(-1)^{a(s)}\phi_{\Q_s}^{\Q^s,\G} (H)= \sum_{\Q\in\FF_s(\M)}(-1)^{a_\Q-a_\G}\,\htau_\Q(H)\ptf
\]
En invoquant l'équation~\eqref{eq1.1d} on en déduit que l'égalité~\eqref{eq1.3d} est vraie pour $\chm$ dans la chambre positive, ce qui établit~\eqref{eq1.2d} et, compte tenu du lemme~\ref{indep}, l'égalité~\eqref{eq1.3d} est encore vraie pour tout $\chm$ en dehors de murs. La dernière assertion résulte alors de la proposition~\ref{envconv}.
\end{proof}

\begin{corollaire}\label{diff}
Soit $\XX$ une famille orthogonale telle que $sX_s$ soit régulier pour tout $s\in\weyl$. En particulier elle est régulière. Soient $L$ et $\M$ deux sous-groupes de Levi avec $\M\subset L$. La différence
\[
\Gamma_L(H,\XX)-\Gamma_\M(H,\XX)
\]
est soit nulle soit égale à~$1$. Il existe une constante $c>0$ telle que si la différence est non nulle on a
\[
\lVert H\rVert \ge c\inf_{P\in\Parab(\M)}\lVert X_P\rVert\ptf
\]
\end{corollaire}

\begin{proof}
D'après la proposition~\ref{preconv} le support de la fonction
\[
H\mapsto \Gamma_\M(H,\XX)
\]
est l'ensemble des $H$ dont la projection sur $\ga_\M^\G$ appartient à l'enveloppe convexe des $X_P$. Comme $\ga_L\subset\ga_\M$, ce support est inclus dans le support de
\[
H\mapsto \Gamma_L(H,\XX)\ptf
\]
Maintenant si la différence est non nulle $H$ est en dehors du support de $\Gamma_\M(H,\XX)$ et donc en dehors de l'enveloppe convexe des $X_P$. La seconde assertion en résulte.
\end{proof}

\section{Cônes et convexes: version duale}

Calculons d'abord les transformées de Laplace des fonctions caractéristiques de cônes. On choisit de manière cohérente des mesures de Haar sur les divers espaces vectoriels $\gaPQ$, par exemple celles définies au moyen des formes de Killing. Soit $\Lambda$ une forme linéaire sur $\gao$. On pose, pour $P$ standard et $\Lambda$ régulier\footnote{Arthur suivi en cela par Langlands \cite{MS}*{Lecture~15} introduisent des fonctions $\hat\theta_P^\Q$ qui sont les inverses de nos $\hat\epsilon_P^\Q$ etc. Dans cet article nous réservons la lettre $\theta$ pour désigner un automorphisme de $\G$.}
\[
\hat\epsilon_P^\Q(\Lambda)=\int_{\gaPQ}\tau_P^\Q(H)\, \ee^{-\Lambda(H)}\dd H\ptf
\]
\newindex{epsilonPQ@$\hat\epsilon_P^Q$}{hepsilonpq}%
On a
\[
\hat\epsilon_P^\Q(\Lambda)=\vol(\hDelta_P^\Q)\prod_{\vpi\in\hDelta_P^\Q} \Lambda (\vpi^\vee)\moins
\]
où $\vol(\hDelta_P^\Q)$ est le volume du parallélépipède engendré par $\hDelta_P^\Q$:

\[ \vol(\hDelta_P^\Q)=\vol\bigl(\ga_P^\Q/\ZM(\hDelta_P^\Q)\bigr)
\]
où $\ZM(\Delta_P^\Q)$ désigne le réseau engendré par $\hDelta_P^\Q$. L'intégrale, convergente si $\Lambda$ est régulier, admet donc un prolongement méromorphe. De même on pose
\[
\epsilon_P^\Q(\Lambda)=\int_{\gaPQ}\htau_P^\Q(H)\, \ee^{-\Lambda(H)}\dd H
\]
\newindex{epsilonPQ@$\epsilon_P^Q$}{epsilonpq}%
et on a
\[
\epsilon_P^\Q(\Lambda)=\vol(\Delta_P^\Q)
\prod_{\alpha\in\Delta_P^\Q} \Lambda (\alpha^\vee)\moins
\]
où
\[
\vol(\Delta_P^\Q)=\vol\bigl(\ga_P^\Q/\ZM(\Delta_P^\Q)\bigr)\ptf
\]
Dans ce qui suit, afin de minimiser le nombre de signes moins dans les exponentielles, on utilisera plutôt les transformées anti-Laplace obtenues en changeant $\Lambda$ en $-\Lambda$. On observe que
\[
\epsilon_P^\Q(-\Lambda)=(-1)^{a_P-a_\Q}\epsilon_P^\Q(\Lambda)\ptf
\]
On a une formule similaire pour $\hat\epsilon_P^\Q$.

\begin{lemme}\label{gpq} 
La transformée anti-Laplace de la fonction
\[
H\mapsto\Gamma_P^\R(H,X)
\]
définie par
\[
\gamma_P^\R(\Lambda,X)=\int_{\ga_P^\R} \ee^{\Lambda(H)}\,\Gamma_P^\R(H,X)\dd H
\]
\newindex{gammaPR(lambda,x)@$\gamma_P^\R(\Lambda,X)$}{gammapr}%
est égale\textup, pour $\Lambda$ régulier\textup, à
\[
\sum_{P\subset\Q} (-1)^{a_P-a_\Q}\ee^{\Lambda(X_\Q^\R)}\hat\epsilon_P^\Q(\Lambda)\epsilon_\Q^\R(\Lambda)
\]
où $X_\Q^\R$ est la projection de $X$ sur $\ga_\Q^\R$. Cette expression se prolonge en une fonction entière sur $\gao^*\otimes\CM$. La fonction
\[
X\mapsto\gamma_P^\R(X)\coloneqq \gamma_P^\R(0,X)= \int_{\ga_P^\R}\Gamma_P^\R(H,X)\dd H
\]
est un polynôme homogène de degré $n=a_P-a_\R=\dim\ga_P^\R$.
\end{lemme}

\begin{proof}
Par définition
\[
\Gamma_P^\R(H,X)=\sum_{P\subset\Q}(-1)^{a_\Q-a_\R}\,\tau_P^\Q(H)\htau_\Q^\R(H-X)\ptf
\]
On a donc, pour $-\Lambda$ régulier,
\[
\int_{\ga_P^\R} \ee^{\Lambda(H)}\,\Gamma_P^\R(H,X)\dd H= \sum_{P\subset\Q}(-1)^{a_\Q-a_\R}\ee^{\Lambda(X_\Q^\R)} \hat\epsilon_P^\Q(-\Lambda) \epsilon_\Q^\R(-\Lambda)
\]
soit encore
\[
\sum_{P\subset\Q}(-1)^{a_P-a_\Q}\ee^{\Lambda(X_\Q^\R)} \hat\epsilon_P^\Q(\Lambda) \epsilon_\Q^\R(\Lambda)
\]
où $X_\Q^\R$ est la projection de $X$ sur $\ga_\Q^\R$.
Comme d'après le lemme~\ref{Gammacar} on intègre une fonction à support compact, la transformée anti-Laplace se prolonge en une fonction entière de $\Lambda\in\gao^*\otimes\CM$. Sa valeur en $\Lambda=0$ est donnée par la somme des termes de degré 0 du développement en série de Laurent des fonctions
\[
t\mapsto(-1)^{a_P-a_\Q}\ee^{t\Lambda(X_\Q^\R)} \hat\epsilon_P^\Q(t\Lambda) \epsilon_\Q^\R(t\Lambda)\ptf
\]
On a donc, pour $\Lambda$ en dehors du lieu des zéros des $\hat\epsilon_P^\Q(\Lambda)\,\epsilon_\Q(\Lambda)$:
\[
\int_{\ga_P^\R}\Gamma_P^\R(H,X)\dd H=\frac{1}{n!}\sum_{P\subset\Q}(-1)^{a_P-a_\Q}\,\Lambda(X_\Q^\R)^{n}
\,\hat\epsilon_P^\Q(\Lambda)\,\epsilon_\Q^\R(\Lambda)
\]
où $n=a_P-a_\R$. Cette expression, indépendante de $\Lambda$, est en tant que fonction de $X$ un polynôme homogène de degré~$n$.
\end{proof}

\begin{lemme}\label{intphis} Supposons $\M$ standard et
$\Lambda$ régulier et considérons $s\in\weyl(\ga_\M)$. On note $P$ le sous-groupe parabolique de sous-groupe de Levi $\M$ et tel que $s(P)$ soit standard. Alors\textup, on a
\[
\int_{\ga_\M^\G} \ee^{\Lambda(H)}\phi_{\M,s}(H- X)\dd H= (-1)^{a(s)}\ee^{\Lambda(X_\M^\G)}\epsilon_{P}^{\G}(\Lambda)\ptf
\]
\end{lemme}

\begin{proof}
Tout d'abord on observe que si $H$ appartient au support de la fonction caractéristique $\phi_{\M,s}^\chm$
avec $\chm$ régulier, on a
\[
\langle\chm,H\rangle =\sum_{\alpha\in\Delta(\M,s)} \alpha(\chm)\vpi_\alpha(H)\le0\ptf
\]
Donc le support de $\phi_{\M,s}=\phi_{\M,s}^\chm$ est contenu dans le cône défini par
\[
\vpi(H)\le0\qquad\forall\vpi\in\hDelta_\M^\G\ptf
\]
On en déduit que l'intégrale converge pour $\Lambda\in\ga^*$ régulier et on a
\[
\int_{\ga_\M^\G} \ee^{\Lambda(H)}\phi_{\M,s}(H- X)\dd H=
\ee^{\Lambda(X_\M^\G)}\int_{\ga_\M^\G} \ee^{\Lambda(H)}\phi_{\M,s}(H)\dd H\ptf
\]
Mais, si $sP$ est standard pour $P$ de Levi $\M$, on voit que pour $\Lambda$ régulier
\[
\int_{\ga_\M^\G} \ee^{\Lambda(H)}\phi_{\M,s}(H)\dd H =\lvert\epsilon_P^\G(\Lambda)\rvert
=(-1)^{a(s)}\epsilon_P^\G(\Lambda)\ptf\qedhere
\]
\end{proof}

\begin{lemme}\label{LaplaceGamma}
Considérons une famille orthogonale $\XX$. Supposons $\M$ standard et posons
\[
\gamma_\M(\Lambda,\XX)=\int_{\ga_\M^\G} \ee^{\Lambda(H)}\,\Gamma_\M(H,\XX)\dd H\ptf
\]
La fonction
$\Lambda\mapsto \gamma_\M(\Lambda,\XX)$ est une fonction entière de $\Lambda\in\gao^*\otimes\CM$. Pour $\Lambda$ en dehors des murs on a
\begin{equation}
\gamma_\M(\Lambda,\XX)= \sum_{P\in\Parab(\M)}\epsilon_{P}^{\G}(\Lambda)\ee^{\Lambda(X_P^\G)}\ptf\label{eq1.1e}
\end{equation}
De plus,
\[ 
\gamma_\M(\XX)\coloneqq \gamma_\M(0,\XX)
\]
est un polynôme en $\mathcali X$, donné par la formule suivante, qui est indépendante du choix de $\LL$ en dehors des murs:
\begin{equation}
\gamma_\M(\XX)=\frac{1}{n!}\sum_{P\in\Parab(\M)}\Lambda(X_P^\G)^{n}\,\epsilon_P^\G(\Lambda)\ptf\label{eq1.2e}
\end{equation}
On a la décomposition
\begin{equation}
\gamma_\M^\Q(\Lambda,\XX)=\sum_{P\in\Parab^\Q(\M)}\gamma_P^\Q(\Lambda,X_P)\ptf\label{eq1.3e}
\end{equation}
\end{lemme}

\begin{proof}
Comme, d'après le corollaire~\ref{compact}, on intègre une fonction à support compact, l'intégrale définit une fonction entière.
Par ailleurs, d'après la proposition~\ref{preconv}
\[
\Gamma_\M(H,\XX)=\sum_{s\in\weyl(\ga_{\M})}(-1)^{a(s)}\phi_{\M,s}(H-X_s)\ptf
\]
On peut supposer $\Lambda$ régulier dans $\gao^*$; il suffit alors d'invoquer le lemme~\ref{intphis} pour obtenir la formule~\eqref{eq1.1e}. Comme cette expression se prolonge en une fonction entière, la formule est encore vraie en dehors des pôles des termes du membre de droite, \cad pour $\Lambda$ en dehors des murs. Sa valeur en $\Lambda=0$ est donnée par la somme des termes de degré 0 du développement en série de Laurent des fonctions:
\[
t\mapsto \sum_{P\in\Parab(\M)}\epsilon_{P}^{\G}(t\Lambda)\ee^{t\Lambda(X_P^\G)}\ptf
\]
La formule~\eqref{eq1.2e} en résulte immédiatement. La décomposition~\eqref{eq1.3e} résulte de l'équation~\eqref{eq1.2c} du lemme~\ref{decomp} en remarquant que les fonctions
\[
H\mapsto \delta_\M^P(H)\Gamma_P^\R(H,X_P)
\]
sont négligeables, et donc de transformée de Fourier nulle, sauf si $P$ est de Levi $\M$.
\end{proof}

\section[\GM-familles]{\mathversion{bold}\GM-familles}\label{gmfam}

Soit $\M$ un sous-groupe de Levi. On appelle \GM\hyph famille la donnée d'une famille de fonctions\footnote{Nous utiliserons la notation $c(\Lambda,P)$ plutôt que $c_P(\Lambda)$ utilisé par Arthur, pour éviter la confusion avec\vadjust{\vspace*{-1.4pt}} les
\[
c_P^\Q(\Lambda)\Qquad{et}c_\M^\Q(\Lambda)\vadjust{\vspace*{-1.4pt}}
\]
introduits plus bas et qui ont une toute autre signification même pour $\Q=\G$. La notation $c_\M^\Q(\Lambda)$ est celle d'Arthur et Langlands. Lorsque $\Q=\G$ Arthur utilise $c'_P$ au lieu de notre $c_P^\G$ et par ailleurs la notation $c_P^\Q$ désigne chez Arthur une notion que nous n'introduisons pas ici. Langlands observe que les notations sont mauvaises. Il serait bien de trouver des notations non ambigües et simples. Une solution correcte, mais très lourde, serait de remplacer $c_P^\Q$ par $\gamma_P^\Q(\bullet,c)$ où $c$ désigne la \GM\hyph famille et de même $\gamma_\M^\Q(\bullet,c)$ pour $c_\M^\Q$.} à valeurs dans un espace vectoriel topologique
\[
\Lambda\mapsto c(\Lambda,P)\qquad\text{sur $\ima\ga_P^*$}
\]
indexées par les sous-groupes paraboliques $P\in\FF(\M)$ (\cad contenant $\M$), qui sont lisses et sujettes aux conditions suivantes: si $P\subset\Q$ alors
\[
c(\Lambda,P)=c(\Lambda,\Q) \qquad\text{pour $\Lambda\in\ima\ga_{\Q}^*$}\ptf
\]
On prolongera $c(\Lambda,P)$ en une fonction sur $\ima\gao^*$ en la supposant constante sur les fibres de la projection
\[
\ima\gao^*\to\ima\ga_P^*\ptf
\]
Il en résulte que si $P$ et $\Q$ sont des sous-groupes paraboliques dans $\Parab(\M)$
adjacents qui correspondent à des chambres
séparées par le mur $\ga_\R$ où $\R$ est le sous-groupe parabolique
engendré par $P$ et $\Q$, alors
\[
c(\Lambda,P)=c(\Lambda,\Q)=c(\Lambda,\R) \qquad\text{pour $\Lambda\in\ima\ga_{\R}^*$}\ptf
\]
Pour définir une \GM\hyph famille il suffit donc de se donner les $c(\Lambda,P)$
pour $P\in\Parab(\M)$, \cad pour les $P$ admettant $\M$ comme sous-groupe de Levi,
satisfaisant la condition
\[
c(\Lambda,P)=c(\Lambda,\Q) \qquad\text{pour $\Lambda\in\ima\ga_\R^*$}
\]
pour des sous-groupes paraboliques $P$ et $\Q$ adjacents.

On rappelle que l'on a noté $\gHM$ l'espace vectoriel des familles $\M$\hyph orthogonales et $\pi_P$ la projection de $\gHM$ sur $\ga_P$ associée à chaque sous-groupe parabolique $P\in\FF(\M)$. C'est un espace de dimension finie. Son dual $\gHM^*$ est muni d'injections
\[
\iota_P\colon\ga_P^*\to\gHM^*
\]
transposées des projections
\[
\pi_P\colon \gHM\to\ga_P\ptf
\]
On dispose dans chaque $\ga_P^*/\ga_\G^*$ de la base $\hDelta_P$ et on rappelle que si $P\subset\Q$ on a $\hDelta_\Q\subset\hDelta_P$. On dispose donc dans $\gHM^*/\ga_\G^*$ d'une base naturelle $\base$ formée par l'union des images de ces bases:
\[
\base_P\coloneqq \iota_P(\hDelta_P)
\]
pour $P\in\FF(M)$. La base $\base$ est l'ensemble des $e_\Q=\iota_\Q(\vpi_\Q)$ où $\Q$ parcourt l'ensemble des sous-groupes paraboliques maximaux propres de $\G$ et on a noté $\vpi_\Q$ l'unique élément de $\hDelta_\Q^\G$\footnote{On remarquera que $\vpi_\Q$ et $\vpi_{\widebar\Q}=-\vpi_\Q$ correspondent à deux sous-groupes paraboliques maximaux opposés définissant deux éléments $e_\Q$ et $e_{\widebar\Q}$ distincts de $\base$.}. Pour $\XX\in\gHM$ et $\LL\in\ima\ga_P^*$ on a
\[
\iota_P(\Lambda)(\XX)=\Lambda\bigl(\pi_P(\XX)\bigr)=\Lambda(X_P)\ptf
\]
On obtiendra une \GM\hyph familles en considérant une fonction $\ext$ lisse sur $\ighMst$ et en posant:
\[
c(\Lambda,P)=(\ext\circ\iota_P)(\Lambda)\ptf
\]
Réciproquement on a la proposition suivante:

\begin{proposition}\label{exist}
Étant donné une \GM\hyph famille $c(\Lambda,P)$ il existe une fonction $\ext$ lisse sur $\ighMst$ telle que
\[
c(\Lambda,P)=(\ext\circ\iota_P)(\Lambda)\ptf
\]
\end{proposition}

\begin{proof}
Notons $\chi$ une fonction sur $\RM$ lisse à support compact, telle que ${\chi(0)=1}$. On a, pour chaque $P\in\FF(\M)$, une partition de la base canonique $\base$ de $\gHM^*$:
\[
\base=\base_P\cup\base^P
\]
où
\[
\base_P=\{e_\Q\mid \Q\supset P\}\Qquad{et}\base^P=\{e_\Q\mid \Q\nsupset P\}
\]
qui induit une décomposition en somme directe et pour $\lambda\in \ighMst$ on note
\[
\lambda=\lambda_P\oplus\lambda^P
\]
la décomposition associée. On définit $\chi_P(\lambda^P)$ en posant
\[
\chi_P(\lambda^P)=\prod_{\vpi\in\base^P} \chi (x_\vpi) \qquad\text{si $\lambda^P=\sum_{\vpi\in\base^P}ix_\vpi\,\vpi$}\ptf
\]
On peut alors définir la fonction
\[
\ext(\lambda)= \sum_{\Q\in\FF(\M)}(-1)^{a_\Q-a_\M} c(\lambda_\Q,\Q)\chi_\Q(\lambda^\Q)\ptf
\]
C'est une fonction lisse sur $\ighMst$. Nous devons calculer $\ext\circ\iota_P$ pour $P\in\Parab(\M)$. Il suffit de le faire lorsque $P$ est anti-standard \cad correspondant à l'opposé de la chambre positive dans $\ga_\M^\G$. D'après le lemme~\ref{soma}, on doit donc calculer
\[
\sum_{s\in\weyl(\ga_{\M})}\,\sum_{\Q\in\FF_s(\M)}(-1)^{a_\Q-a_\M} c(\iota_P(\LL)_\Q,\Q)\chi_\Q(\iota_P(\LL)^\Q)\ptf
\]
On observe maintenant que, si on note $F_P$ la facette associée à $P$ dans $\ga_\M$, on a
\[
\widebar F_P\cap\widebar F_{\Q_s}=\widebar F_P\cap\widebar F_{\Q}
\]
pour $\Q\in\FF_s(\M)$ et donc
\[
\iota_P(\ga^*_\M)\cap\iota_{\Q_s}(\ga^*_\M)=\iota_P(\ga^*_\M)\cap\iota_{\Q}(\ga^*_{\Q})\ptf
\]
Il en résulte que, pour $\LL\in\ima\ga_\M^*$ et $\Q\in\FF_s(\M)$, la projection $\iota_P(\LL)_\Q$ sur
\[
\iota_P(\ga^*_\M)\cap\iota_{\Q}(\ga^*_\Q)
\]
ne dépend que de $s$. De plus, si $s$ est l'élément du groupe de Weyl associé à l'opposé de la chambre positive dans $\ga_\M^\G$, on a
\[
\Q_s=\Q^s=P
\]
et
\[
c(\iota_P(\LL)_P,P)\chi_P(\iota_P(\LL)^P)= c(\LL,P)\chi_P(0)= c(\LL,P)\ptf
\]
Il résulte alors de ces observations et de la seconde assertion de lemme~\ref{soma} que la fonction $\ext$ a les propriétés désirées, \cad que
\[
(\ext\circ\iota_P)(\Lambda)=c(\Lambda,P)\ptf\qedhere
\]
\end{proof}

Dans la suite de cette section on se limitera aux \GM\hyph familles à valeurs dans un espace vectoriel de dimension finie.

On dira qu'une mesure de Radon $\mm$ sur un espace vectoriel $V$ est à décroissance rapide si son produit avec n'importe quelle fonction polynôme $p$ sur $V$ fournit une mesure bornée:
\[
\int_V \lvert p(x)\rvert\dd \lvert \mm\rvert(x)<+\infty\ptf
\]
En particulier $\mm$ est bornée et admet une transformée de Fourier qui est une fonction lisse. Supposons maintenant que $\ext$ est la transformée de Fourier d'une mesure de Radon $\mm$ sur $\gHM$, à décroissance rapide:
\[
\ext(\lambda)=\int_{\gHM} \ee^{\lambda(\XX)}d\mm(\XX)\ptf
\]
La \GM\hyph famille définie par $\ext$ peut alors s'écrire
\[
c(\Lambda,P) =\int_{\gH_\M} \,\ee^{\iota_P(\Lambda)(\XX)}\dd \mm(\XX) =\int_{\gH_\M} \,\ee^{\Lambda(X_P)}\dd \mm(\XX)
\]
où $X_P=\pi_P(\XX)$.

\begin{corollaire}\label{transfour}
Si les composantes $c(P,\LL)$ d'une \GM\hyph famille sont des fonctions dans l'espace de Schwartz sur $\ima\ga_\M^*$\textup, il existe une fonction $\vf$ dans l'espace de Schwartz sur $\gHM$ fournissant la \GM\hyph famille par transformation de Fourier\textup:
\[
c(\Lambda,P)=\int_{\gHM}\ee^{\iota_P(\Lambda)(\XX)}\vf(\XX)\dd \XX\ptf
\]
\end{corollaire}

\begin{proof}
La construction donnée dans la proposition~\ref{exist} fournit une fonction $\ext$ dans l'espace de Schwartz de $\ima\gHM^*$. Il suffit de prendre pour $\vf$ sa transformée de Fourier.
\end{proof}

Étant donné une \GM\hyph famille on introduit, pour $\Lambda$ en dehors des murs,
\[
c_P^\R(\Lambda)=\sum_{P\subset\Q\subset\R}(-1)^{a_P-a_\Q}\,\hat\epsilon_P^\Q(\Lambda)\epsilon_\Q^\R(\Lambda)c(\Lambda,\Q)
\]
et
\[
c_\M^\Q(\Lambda)=\sum_{P\in\Parab^\Q(\M)}\epsilon_{P}^{\Q}(\Lambda)c(\Lambda,P)\ptf
\]

\begin{lemme}\label{gmorth}
Soit $\XX=\{X_P\}$ une famille orthogonale.
Alors
\[
c(\Lambda,P)=\ee^{\Lambda(X_P)}\qquad\text{pour $P\in\Parab(\M)$}
\]
est une \GM\hyph famille. Dans ce cas on a
\begin{equation}
c_P^\Q(\Lambda)= \int_{\gaPQ} \ee^{\Lambda(H+X_\Q)}\,\Gamma_P^\Q(H,X_P)\dd H =\ee^{\Lambda(X_\Q)}\gamma_P^\Q(\Lambda,X_P)\label{eq1.1f}
\end{equation}
et
\begin{equation}
c_\M^\Q(\Lambda)=\int_{\ga_\M^\Q} \ee^{\Lambda(H+X_\Q)}\,\Gamma_\M^\Q(H,\XX)\dd H =\ee^{\Lambda(X_\Q)}\gamma_\M^\Q(\Lambda,\XX)
\ptf\label{eq1.2f}
\end{equation}
Plus généralement\textup, si
\[
c(\Lambda,P)=(\ext\circ\iota_P)(\Lambda)\quad\text{pour $P\in\Parab(\M)$}
\]
où $\ext$ est la transformée de Fourier d'une mesure $\mm$ à décroissance rapide on a
\begin{equation}
c_P^\Q(\Lambda)=
\int_{\gHM} \int_{\gaPQ} \ee^{\Lambda(H+X_\Q)}\,\Gamma_P^\Q(H,X_P)\dd H\dd \mm(\XX)
\label{eq1.3f}
\end{equation}
soit encore
\begin{equation}
c_P^\Q(\Lambda) =\int_{\gHM}\ee^{\Lambda(X_\Q)}\gamma_P^\Q(\Lambda,X_P)\dd\mm(\XX)
\tag{\ref{eq1.3f}$'$}\label{eq1.3f'}
\end{equation}
où $X_P=\pi_P(\XX)$
et
\begin{equation}
c_\M^\Q(\Lambda)=\int_{\gHM}\int_{\ga_\M^\Q} \ee^{\Lambda(H+X_\Q)}\,\Gamma_\M^\Q(H,\XX)
\dd H\dd \mm(\XX)
\label{eq1.4f}
\end{equation}
 soit encore
\begin{equation}
c_\M^\Q(\Lambda)=\int_{\gHM}\ee^{\Lambda(X_\Q)}\gamma_\M^\Q(\Lambda,\XX)\dd\mm(\XX)\ptf
\tag{\ref{eq1.4f}$'$}\label{eq1.4f'}
\end{equation}
\end{lemme}

\begin{proof} 
D'après le lemme~\ref{gpq} on sait que
\[
\ee^{\Lambda(X_\R)}\gamma_P^\R(\Lambda,X)=\int_{\ga_P^\R} \ee^{\Lambda(H+X_\R)}\,\Gamma_P^\R(H,X)\dd H
\]
est égale, pour $\Lambda$ régulier, à
\[
\sum_{\{\Q\mid P\subset\Q\}}(-1)^{a_P-a_\Q}\ee^{\Lambda(X_\Q)} \hat\epsilon_P^\Q(\Lambda)\epsilon_\Q^\R(\Lambda)
\]
qui est la définition de $c_P^\R(\Lambda)$ si
$c(\Lambda,P)=\ee^{\Lambda(X_P)}$, d'où l'assertion~\eqref{eq1.1f}. L'assertion~\eqref{eq1.3f} résulte, au moins formellement, de~\eqref{eq1.1f} et la convergence
résulte de ce que
\[
\lvert \gamma_P^\R(\Lambda,X)\rvert
\] est majoré par un polynôme en $X$.
De même~\eqref{eq1.2f} et~\eqref{eq1.4f} résultent du lemme~\ref{LaplaceGamma}.
\end{proof}

\begin{lemme}\label{lisse}
Soit $\{c(\Lambda,P)\}$ une \GM\hyph famille. Les fonctions $c_P^\R(\Lambda)$ et\linebreak $c_\M^\Q(\Lambda)$ définies
pour $\Lambda$ en dehors des murs se prolongent en fonctions lisses partout et on~a
\begin{equation}
c_\M^\Q(\Lambda)=\sum_{P\in\Parab^\Q(\M)}c_P^\Q(\Lambda)\ptf
\tag{$*$}\label{eq1.*}
\end{equation}
\end{lemme}

\begin{proof}
Le problème est local en $\Lambda$; il nous est donc loisible de supposer les $c(\Lambda,P)$ à support compact. Dans ce cas, d'après le corollaire~\ref{transfour} on peut les supposer de la forme
\[
c(\Lambda,P)=(\ext\circ\iota_P)(\Lambda)\qquad\text{pour $P\in\Parab(\M)$}
\]
avec $\ext$ dans l'espace de Schwartz. On a donc que pour $\Lambda\in\ima\ga_P^*$
\[
c(\Lambda,P)=\int_{\gHM} \ee^{\Lambda(X_P)}\vf(\XX)\dd \XX
\]
avec $\vf=\widehat\ext$ à décroissance rapide. Mais, le lemme~\ref{gmorth} montre que pour $\Lambda\in\ima\ga_P^*$
\[
c_P^\Q(\Lambda) =\int_{\gHM}\ee^{\Lambda(X_\Q)}\gamma_P^\Q(\Lambda,X_P)\vf(\XX)\dd \XX\ptf
\]
La lissité de $c_P^\Q$ résulte du lemme~\ref{gpq} et de ce que
\[
\gamma_P^\Q(\LL,X)
\]
est une fonction lisse en $\LL$, majorée par un polynôme en $X$, alors que $\widehat\ext$ est à décroissance rapide. De même, compte tenu du lemme~\ref{gmorth}, on a pour $\Lambda\in\ima\ga_\M^*$
\[
c_\M^\Q(\Lambda)=\int_{\gHM}\ee^{\Lambda(X_\Q)}\gamma_\M^\Q(\Lambda,\XX)\vf(\XX)\dd \XX\ptf
\]
La décomposition~\eqref{eq1.*} ainsi que la lissité de $c_\M^\Q$ résultent du lemme~\ref{LaplaceGamma}.
\end{proof}

\begin{lemme}\label{decompGM}
Supposons que $c$ et $d$ soient deux \GM\hyph familles\textup, la première étant à valeurs scalaires\textup, et considérons la \GM\hyph famille produit
\[
e(\Lambda,P)=c(\Lambda,P)d(\Lambda,P)\ptf
\]
Si on pose
\[
e_\M^\R(\Lambda)= \sum_{P\in\Parab^\R(\M)} \epsilon_{P}^{\R}(\Lambda)e(\Lambda,P)
\]
on a
\[
e_\M^\R(\Lambda)=\sum_{Q\in\FF^\R(\M)}c_\M^\Q(\Lambda)d_\Q^\R(\Lambda)\ptf
\]
\end{lemme}

\begin{proof}
Ici encore le problème est local en $\Lambda$; il nous est donc loisible de supposer les familles $c$ et $d$ à support compact; on peut alors, d'après le corollaire~\ref{transfour}, supposer l'existence de fonction $\vf$ et $\psi$ dans l'espace de Schwartz sur $\gHM$ telles que
\[
c(\Lambda,P)d(\Lambda,P) =\int_{\gHM}\int_{\gHM} \ee^{\Lambda(\XX+\YY)}\varphi(\XX)\psi(\YY)\dd \XX\dd \YY\ptf
\]
Dans ce cas $e_\M^\R(\Lambda)$ est égal à l'intégrale triple
\[
\int_{H\in\ga_P^\R}\int_{\XX\in\gHM}\int_{\YY\in\gHM} \ee^{\Lambda(H+X_\R+Y_\R)} \,\Gamma_\M^\R(H, \XX+\YY)\varphi(\XX)\psi(\YY)\dd \XX\dd \YY\dd H
\]
Le lemme résulte alors du lemme~\ref{convol}.
\end{proof}

\chapter{Espaces tordus}\label{esptor}

\section{Sorites}\label{sorites}

La notion d'espace tordu peut se définir dans diverses catégories.
Les définitions ci-dessous s'entendent soit dans la catégorie des ensembles, soit
dans la catégorie des espaces localement compacts soit encore dans la catégorie
des variétés algébriques.

Rappelons qu'un espace tordu
est la donnée d'un couple $(\G,\tG)$ où
$\G$ est un groupe et $\tG$ un $\G$\hyph torseur
(\ie un $\G$\hyph espace principal homogène) à gauche,
muni d'une application $\G$\hyph équivariante dans le groupe $\Aut(\G)$ des automorphismes de $\G$:
\[
\AdtG\colon \tG\to\Aut(\G)\ptf
\]
L'équivariance signifie que pour tout $x\in\G$ et tout $\delta\in\tG$ on a
\[
\AdtG(x\delta)=\AdG(x)\circ\AdtG(\delta)
\]
où $\AdG(x)$ est l'automorphisme intérieur défini par $x$.
L'application $\AdtG$ n'est pas injective en général: ses fibres
sont des torseurs sous le centre $Z_\G$ de $\G$.
On notera $\Int(\G)$ le groupe des automorphismes intérieurs
et $\Out(G)$ le groupe des automorphismes extérieurs.
La suite exacte
\[
1\to\Int(\G)\to\Aut(\G)\to\Out(\G)\to1
\]
montre que la classe d'isomorphisme de $\tG$ est déterminée
par l'unique élément image de $\tG$ dans $\Out(\G)$.\,\footnote{Le lecteur prendra garde à ce que le groupe des automorphismes
dépend fortement de la catégorie où on se place.}
On définit une action à droite de $\G$ sur $\tG$ en posant
\[
\delta\,x=\theta(x)\,\delta\qquad\text{avec $\theta=\AdtG(\delta)$}\ptf
\]
On dispose alors sur $\tG$ d'une structure de $\G$\hyph torseur à droite et à gauche
donc en particulier d'une action par conjugaison
de $\G$ sur $\tG$ et de la notion de classe de $\G$\hyph conjugaison dans $\tG$.
On note $Z_G({\tG})$ le centralisateur de $\tG$ dans $\G$.
Il est facile de voir que
\[
Z_G({\tG})=(Z_\G)^\theta
\]
le sous-groupe des $\theta$\hyph invariants dans le centre de $\G$.

On peut regarder un espace tordu $(\G,\tG)$ comme les composantes d'indice~$0$ et~$1$:
\[
\G=\GG_0\Qquad{et} \tG=\GG_1
\]
d'un groupe gradué par $\ZM$
\[
\GG=\coprod_{n\in\ZM}\GG_n\ptf
\]
Tous les $\GG_n$ sont des $\G$\hyph espaces tordus et
on dispose en particulier du $\G$\hyph espace tordu inverse $\tG\moins$:
\[
\tG\moins=\GG_{-1}\ptf
\]
L'espace tordu $(\G,\tG\moins)$ peut être défini au moyen d'une
bijection: $\tG\to\tG\moins$ notée $\delta\mapsto\delta\moins$
vérifiant, pour $x$ et $y$ dans $\G$:
\[
(x\delta y)\moins=y\moins\delta\moins x\moins
\Qquad{et}
\tAd(\delta\moins)=\tAd(\delta)\moins\ptf
\]
On dispose alors d'une application $\tG\times\tG\moins\to\G$ qui sera notée comme un produit
(c'est en effet le produit dans $\GG$):
\[
(\tau,\delta\moins)\mapsto\tau\delta\moins
\]
avec la propriété suivante:
si $\tau=u\delta$ pour $u\in\G$
on a $\tau\delta\moins=u$ et plus généralement
pour tout $x\in\G$ on a
\[
\tau\,x\,\delta\moins=u\,\theta(x)
\]
où $\theta=\tAd(\delta)$.
La donnée de $\theta=\AdtG(\delta)$ fournit un isomorphisme d'espace tordus
\[
\tG\to\G\rtimes\theta\subset \G\rtimes\Aut(\G)
\]
défini par
\[
x\delta\mapsto x\rtimes\theta\Quad{pour}x\in\G\ptf
\]
Toutefois, cet isomorphisme n'est pas canonique; il dépend du choix de $\delta$.

Dans le cadre des espaces localement compacts, un espace tordu $\tG$ est un espace tordu ensembliste où $\G$ est un groupe
localement compact et où les morphismes considérés dans la définition de la structure sont continus. Une mesure
$\G$\hyph invariante à droite ou à gauche sur $\tG$ sera appelée une mesure de Haar.
La donnée d'une mesure de Haar à gauche $\mu$ sur $\G$
permet de définir une mesure de Haar à gauche $\tilde\mu$ sur $\tG$ en posant
pour $f\in\ctc(G)$:
\[
\tilde\mu(f)=\int_\G f(x\delta)\dd\mu(x)\ptf
\]
Un espace tordu $\tG$ localement compact sera dit unimodulaire si pour tout $\delta\in\tG$
l'automorphisme $\theta=\AdtG(\delta)$
est de module $1$. Ceci implique que
$\G$ est unimodulaire.

On dispose également de la notion d'espace tordu dans la
catégorie des variétés algébriques.
Considérons un espace tordu algébrique $(\G,\tG)$ où $\G$ est
groupe linéaire algébrique connexe défini sur un corps $F$. S'il est non vide,
l'ensemble $\tG(F)$ est un espace tordu sous $\G(F)$ au sens
ensembliste, et si $A$ est une $F$\hyph algèbre localement compacte alors $\tG(A)$ est un
espace tordu localement compact.

\section{Exemples}

Un des exemples, important pour les applications, est le suivant.
Soit $V$ un espace vectoriel de dimension finie sur un corps $F$
et soit $V^*$ son dual. Considérons le groupe $\G=\GL(V)$; on dispose de la
représentation contragrédiente de $\G$ dans $V^*$ définie pour $x\in\G$ par
\[
x\mapsto x^\vee={}^t x\moins\ptf
\]
On notera
\[
\tG=\Isom(V,V^*)
\]
l'espace des isomorphismes $V\to V^*$
ou, si on préfère, l'espace des formes bilinéaires non dégénérées sur $V\times V$.
C'est un $\G$\hyph torseur à droite et à gauche en posant pour $x\in\G$ et $\delta\in\tG$:
\[
x\delta=x^\vee\circ\delta\Quad{et} \delta x=\delta\circ x\ptf
\]
Pour tout $\delta\in\tG$ on définit un automorphisme $\theta=\AdtG(\delta)$
de $\G$ en posant pour $x\in\G$
\[
\theta(x)=(\delta\,\circ x\,\circ\delta\moins)^\vee
\]
et ceci munit $\tG$ d'une structure de $\G$\hyph espace tordu.
Si on munit $V$ d'une base on dispose alors d'un isomorphisme de groupes
\[
\iota\colon \GL(V)\to \GL(n,F)
\]
et de l'application
\[
\delto\in \Isom(V,V^*)
\]
qui envoie la base de $V$ sur la base duale dans $V^*$. Dans ce cas
l'automorphisme de $\GL(n,F)$ associé à $\theto=\tAd(\delto)$
est l'inverse de la transposée pour les matrices:
\[
\iota\circ\theto(x) = {}^t \iota(x)\moins\ptf
\]
En d'autres termes le choix d'une base dans $V$ fournit un isomorphisme
\[
\tG\simeq \GL(n)\rtimes\ve
\]
où $\ve(m)={}^tm\moins$ pour $m\in\GL(n)$.

Dans le \textenglish{Morn\-ing Sem\-i\-nar} on suppose dans les deux premiers exposés que
\[
\tG=\G\rtimes\theta\subset\G\rtimes\langle\theta\rangle
\]
où $\theta$ est un automorphisme d'ordre fini; dans les exposés suivants (3 à~15) on considère,  comme dans \cite{ALB} et les autres articles d'Arthur sur le cas tordu, le cas un peu plus général où  $\tG$ (noté $\G$ chez Arthur) est une composante connexe d'un groupe réductif non connexe (noté $\G'$ dans \cite{MS} et $\G^+$ chez Arthur et dont la composante neutre est notée $\G$ dans \cite{MS} et $\G^0$ par Arthur).

\section{Représentations tordues}\label{reptor}

Soit $\omega$ un caractère de $\G$ et $\tG$ un $\G$\hyph espace tordu.
Soit $V$ un espace vectoriel.
On appelle représentation tordue de $\G$ dans $V$
pour le couple $(\tG,\omega)$, ou simplement représentation de $(\tG,\omega)$,
la donnée pour tout $\delta\in\tG$ d'un endomorphisme inversible
\[
\tpi(\delta,\omega)\in \GL(V)
\]
et d'une représentation $\Gpi$ de $\G$ dans $V$:
\[
\Gpi\colon G\to \GL(V)
\]
vérifiant pour $x,y\in\G$ et $\delta\in\tG$
\[
\tpi(x\,\delta\,y,\omega)=\Gpi(x)\,\tpi(\delta,\omega)\,(\Gpi\otimes\omega)(y)\ptf
\]
En particulier
\[
\tpi(\delta\,x,\omega)=\tpi(\delta,\omega)\,(\Gpi\otimes\omega)(x)
=\tpi(\theta(x)\,\delta,\omega)=\Gpi\bigl(\theta(x)\bigr)\,\tpi(\delta,\omega)
\]
et donc $\tpi(\delta,\omega)$ entrelace
$\pi\otimes\omega$ et $\pi\circ\theta$ si $\theta=\AdtG(\delta)$.
La donnée de $\tpi$ détermine $\Gpi$; on dira que $\Gpi$ est la restriction
de $\tpi$ à $\G$ et on écrira $\pi=\tpi\rest G$.

Réciproquement $\tpi$
est déterminé par la donnée de la représentation $\Gpi$ et,
pour un $\delta\in\tG$, d'un opérateur $A$ qui entrelace
$\pi\otimes\omega$ et $\pi\circ\theta$ avec $\theta=\AdtG(\delta)$:
\[
A\,(\pi\otimes\omega)(x)=(\pi\circ\theta)(x)\,A\ptf
\]
On reconstruit $\tpi$ en posant
\[
\tpi(x\delta,\omega)=\pi(x)A\qquad\text{pour $x\in\G$}\ptf
\]

Si $V$ est un espace de Hilbert
on dira que $\tpi$ est unitaire si $\tpi$ prend ses valeurs dans
le groupe unitaire de $V$.
Si $\tpi$ est unitaire et si $\pi$ est irréductible le lemme de Schur montre que
$\pi$ détermine $$A=\tpi(\delta,\omega)$$ à un scalaire non nul près, indépendant de $\delta$.

On dira que deux représentations tordues
$(\tpi,V)$ et $(\tpi',V')$
sont équivalentes s'il existe un opérateur d'entrelacement inversible
\[
I\colon V\to V'
\]
tel que, pour tout $\delta\in\tG$ on ait
\[
I\,\tpi(\delta,\omega)=\tpi'(\delta,\omega)\,I\ptf
\]
Fixons $\delta\in\tG$ et posons $A=\tpi(\delta,\omega)$ et $A'=\tpi'(\delta,\omega)$. On doit avoir pour tout $x\in\G$:
\[
I\,\tpi(x\delta,\omega)=\tpi'(x\delta,\omega)\,I
\]
soit encore
\[
I\,\pi(x)A=\pi'(x)A'\,I
\]
et en particulier $I\,A=A'\,I$ avec $A$ inversible et donc
\[
I\,\pi(x)=\pi'(x)\,I\ptf
\]
C'est dire que $\pi$ et $\pi'$ sont équivalentes.
Mais la réciproque est fausse puisque, même si $\pi$ est unitaire irréductible,
la classe de $\pi$ ne détermine $A$ qu'à un scalaire non nul près.

Supposons que $(\tpi,V)$ est une représentation unitaire et que $(\pi,V)$ est
une somme directe hilbertienne (finie ou dénombrable)
de représentations irréductibles.
On dira que $\pi$ est quasi simple (relativement à $\tG$)
s'il existe une représentation irréductible $(\sigma,W)$
et un entier $\ell = \ell (\pi)$ positif ou nul
tels que
\[
(\pi,V)=\widehat{\bigoplus_{r\in\ZM/\ell(\pi)}}(\sigma_r,W_r)
\]
où $W_r=W$ en tant qu'espace vectoriel pour tout $r$, mais est muni de la représentation $\sigma_r$ définie par
\[
\sigma_r=(\sigma\circ\theta^r)\otimes \omega^{-r}
\]
où $\theta=\AdtG(\delta)$ avec $\delta$ fixé dans $\tG$.
De plus les $(\sigma_r,W_r)$ sont deux à deux inéquivalentes
pour des $r$ non congrus modulo $\ell$ alors que
\[
\sigma_r\simeq\sigma_{r+\ell}\ptf
\]

\begin{lemme}
Supposons que $(\tpi,V)$ est une représentation unitaire et que $(\pi,V)$ est
une somme directe hilbertienne \textup(finie ou dénombrable\textup)
de représentations irréductibles.
Alors $\pi$ est quasi simple si et seulement si $\tpi$ est irréductible.
\end{lemme}

\begin{proof}
Par hypothèse, $\pi$ est somme
finie ou dénombrable de représentations irréductibles.
On peut décomposer $(\pi,V)$ en somme de composants isotypiques.
Soit $(\sigma_0,W_0)$ une des composantes isotypiques.
C'est un multiple d'une représentation irréductible $(\sigma,W)$ de $\G$.
Posons $$A=\tpi(\delta,\omega)$$ et notons $W_r$ l'espace $A^rW_0$.
On a, pour $w\in W_r$,
\[
\pi(x)w=A^r\pi(x)A^{-r}w=\sigma_r (x) w
\]
avec
\[
\sigma_r=(\sigma_0\circ\theta^r)\otimes \omega^{-r}\ptf
\]
Donc $(\pi,V)$ contient tous les $(\sigma_r,W_r)$.
En particulier, si $\pi$ est quasi simple $\tpi$ est irréductible.
Examinons la réciproque. Deux cas sont alors possibles:
\begin{asparaenum}
\item Il y a un nombre fini de composantes isotypiques.
Il existe donc un plus petit entier $\ell\ge1$ tel
que
\[
\sigma_0\simeq\sigma_\ell\ptf
\]
Dans ce cas l'opérateur $A^\ell$ peut s'écrire comme
une somme directe finie
\[
A^\ell=\bigoplus_{r\in\ZM/\ell(\pi)} B_r
\]
où $B_r$ est la restriction de $A^\ell$ au sous-espace
isotypique $W_r$.
Mais, tout projecteur spectral non trivial de $B_0$ permet de construire
un sous-espace $\tG$\hyph invariant non trivial de $V$ et comme $(\tpi,V)$ est irréductible
$B_0$ est nécessairement scalaire. Il en est donc de même de $A^\ell$.
On en déduit qu'un sous espace
$\G$ invariant $W$ dans $W_0$ engendre un sous-espace $\tG$\hyph invariant
qui n'est l'espace $V$ tout entier que si $W=W_0$. On en déduit que
$\sigma_0$ est irréductible et que donc $\pi$ est quasi simple.
\item Il y a un nombre infini de composantes isotypiques. Soit $W$ un sous espace
irréductible dans $W_0$. L'adhérence de la somme directe
des $A^rW$ avec $r\in\ZM$ est un sous-espace $\tG$\hyph invariant et
c'est donc l'espace $V$ tout entier par irréductibilité de~$\tpi$.\qedhere
\end{asparaenum}
\end{proof}

Nous supposerons dans la suite de cette discussion que $\tG$ est un
espace tordu localement compact unimodulaire,
muni d'une mesure de Haar et que toutes les représentations unitaires
considérées sont continues. Soit
$V$ l'espace vectoriel des fonctions de carré intégrable sur $\G$. On dispose
d'une représentation unitaire naturelle, appelée représentation régulière, de $\tG\times\tG$ dans $V$
en posant pour $\varphi\in V$, $x\in\G$ et $(\delta,\tau)\in\tG\times\tG$:
\[
(\reg(\delta,\tau)\varphi)(x)=\varphi(\delta\moins\,x\,\tau)\ptf
\]
C'est une variante de cette représentation qui est au cœur
de la théorie de la formule des traces tordue.

Soit $(\tpi,V)$ une représentation unitaire et soit $(\pi,V)$ sa restriction à $\G$.
Considérons une fonction $g\in\ctyc(\G)$.
Il est classique de considérer l'opérateur défini par l'intégrale (qui a un sens
pour la topologie forte sur l'espace de opérateurs)
\[
\pi(g)=\int_\G g(x)\pi(x)\dd x\ptf
\]
De même on posera pour $f\in\ctyc(\tG)$
\[
\tpi(f)=\int_{\tG} f(y)\tpi(y)\dd y\coloneqq \int_\G f(x\delta)\tpi(x\delta)\dd x\ptf
\]

\begin{lemme}\label{tracenulle}
Soit $\tG$ un espace tordu localement compact unimodulaire
muni d'une mesure de Haar.
Soit $(\tpi,V)$ une représentation unitaire irréductible.
On suppose que $\pi(g)$ est un opérateur à trace pour
toute fonction $g\in\ctyc(\G)$.
Alors\textup, l'opérateur
$\tpi(f)$ sera à trace pour toute fonction
$f\in\ctyc(\tG)$ et\textup, si $\ell(\pi)\ne1$\textup, on aura
\[
\tr\tpi(f)=0\ptf
\]
\end{lemme}

\begin{proof}
On suppose que $R=\pi(g)$ est un opérateur à trace pour
toute fonction $g\in\ctyc(\G)$.
En particulier $\pi$ est somme discrète d'irréductibles.
Comme $(\tpi,V)$ est irréductible $(\pi,V)$ est quasi simple.
On a $f\in\ctyc(\tG)$ et on pose
$g(x)=f(x\delta)$
et $A=\tpi(\delta,\omega)$.
L'opérateur $A=\tpi(\delta,\omega)$ est un opérateur unitaire
et donc
\[
\tpi(f)=RA
\]
est un opérateur à trace.
On rappelle que
\[
V=\widehat{\bigoplus_{r\in\ZM/\ell(\pi)} }W_r
\]
et on observe que
\[
RA(W_r)\subset W_{r+1}
\]
et donc si $\ell\ne1$
\[
\tr(RA)=0\ptf\qedhere
\]
\end{proof}

\section{Multiplicités des représentations tordues}\label{multi}

Considérons une représentation tordue $(\treg,H)$
somme directe de représentations de représentations irréductibles
$(\tpi,V_\pi)$
avec multiplicité $\mtG(\tpi)$, \cad que, comme $\tG$\hyph module,
\[
H=\bigoplus_{\tpi} M(\tpi)\otimes V_\pi
\]
où $M(\tpi)$ est un espace de dimension $\mtG(\tpi)$ sur lequel
$\tG$ agit trivialement.
Soit $\tpi$ une représentation de $\tG$ dont la restriction $\Gpi$ à $G$
reste irréductible. On suppose que $\tpi$ intervient dans $(\treg,H)$ avec la multiplicité
$\mtG(\tpi)$. Si on note $\mG(\Gpi)$ la multiplicité de $\Gpi$ dans $(\reg,H)$,
la restriction à $\G$ de $(\treg,H)$,
on a
\[
\mG(\pi)=\sum_{\tpi\rest G\simeq\pi}\mtG(\tpi)
\]
et en particulier
\[
\mtG(\tpi)\le \mG(\Gpi)\ptf
\]
De fait on a
\[
H=\bigoplus_\pi M(\pi)\otimes V_\pi
\qquad\text{avec $M(\pi)=\bigoplus_{\tpi'\rest G\simeq\pi} M(\tpi')$}\ptf
\]

Cette notion naïve de multiplicité n'est pas la bonne lorsqu'on
souhaite exploiter l'indépendance linéaire des traces et,
pour ce faire, il sera nécessaire de regrouper les contributions des
diverses $\tpi$ ayant même restriction $\pi$ à $\G$.
En effet, soient $\tpi$ et $\tpi'$ deux représentations
irréductibles dans un même espace $V$ et qui ont la
même restriction $\pi$ à $\G$ avec $\pi$ irréductible.
Alors les opérateurs $A=\tpi(\delta,\omega)$ et $A'=\tpi'(\delta,\omega)$ sont proportionnels.
Donc, si les opérateurs $\tpi(f)$ sont à trace, les formes linéaires
\[
f\mapsto\tr\tpi(f)\Quad{et} f\mapsto\tr\tpi'(f)
\]
sont proportionnelles. Notons
\[
\lambda(\tpi',\tpi)\in\CM^\times
\]
le scalaire tel que, pour tout $\delta\in\tG$ on ait
\[
\tpi'(\delta,\omega)=\lambda(\tpi',\tpi)\tpi(\delta,\omega)\ptf
\]
Si l'ensemble des $\tpi$ ayant même restriction $\pi$
et intervenant avec une multiplicté non nulle
est fini et si $\tpi(f)$ est un opérateur à trace pour
toute $f\in\ctyc(\tG)$, l'objet qu'il convient de considérer est la somme
\[
\sum_{\tpi'\rest G=\pi}\mtG(\tpi')\,\tr\tpi'(f)= m(\pi,\tpi)\tr\tpi(f)
\]
où l'on a posé
\[
m(\pi,\tpi)=\sum_{\tpi'\rest G=\pi}\lambda(\tpi',\tpi)\,\mtG(\tpi') \ptf
\]
L'ensemble des $\tpi$ ayant même restriction $\pi$ forment
un torseur sous $\CM^\times$ de même que les
nombres $\tr\tpi(f)$ (que l'on pourrait appeler la trace tordue).
L'ensemble des nombres $m(\pi,\tpi)$ peut
être vu comme un torseur
sous $\CM^\times$ que l'on appelera la multiplicité tordue.
Le produit des deux torseurs
\[
m(\pi,\tpi)\tr\tpi(f)
\]
est un nombre indépendant du choix du point base $\tpi$: il ne dépend que de $\pi$.

\section{Espaces tordus réductifs}\label{tordalg}

On suppose désormais que $(\G,\tG)$
\newindex{G@$\protect\tG$}{torG}%
est un espace tordu algébrique où $\G$ est un
groupe linéaire algébrique connexe défini sur un corps
de nombres $F$ et on suppose que $\tG(F)$ est non vide.
On peut alors définir l'espace tordu $\tGadef$ des points adéliques de $\tG$.
Tout élément $y\in\tGadef$ est de la forme
\[
y=x\delta
\]
avec $x\in\Gadef $ et $\delta\in\tG(F)$.
L'automorphisme induit par
\[
\theta=\AdtG(\delta)
\]
sur $\ga_\G$ sera encore noté $\theta$;
il est indépendant du choix de $\delta\in\tG(F)$. On notera
$\ga_{\tG}$ le sous-espace vectoriel des points fixes
sous $\theta$ dans $\ga_\G$:
\[
\ga_{\tG}=(\ga_\G)^\theta
\]
et on pose
\[
a_{\tG}\coloneqq \dim\ga_{\tG}\ptf
\]
Supposons $\G$ réductif.
On notera $\gAtG$ le sous-groupe des points fixes
sous $\theta$ dans $\gA_G$ et $\HG$ induit un isomorphisme
\[
\gAtG\to\ga_{\tG}\ptf
\]
Pour toutes les applications envisagées à ce jour il est loisible de
supposer que l'automorphisme induit par $\theta$ sur $\ga_\G$ est semi-simple
voire même d'ordre fini. Toutefois, Kottwitz et Shelstad font dans \cite{KS}*{\parag6.1}
une hypothèse moins restrictive: ils supposent simplement
que l'application entre invariants et coinvariants
\[
(\ga_\G)^\theta\to (\ga_\G)_\theta
\]
est un isomorphisme. Si on note $\ga_\G^{1-\theta}$ le noyau de
la surjection sur les coinvariants, la condition peut se reformuler ainsi:
on a une décomposition en somme directe
\[
\ga_\G=\ga_{\tG}\oplus\ga_\G^{1-\theta}\ptf
\]
Ceci est encore équivalent à demander que le jacobien
\[
\jtG=\lvert\det(\theta-1\rest \ga_\G/\ga_{\tG})\rvert
\]
\newnot{j(G)@$\jtG$}{jtGG}
soit non nul. Nous le supposerons désormais.

On a choisi un sous-groupe parabolique minimal $\PO$ dans $\G$
et un sous-groupe de Levi $\M_0\subset\PO$ définis sur $F$.
\newnot{P0@$\PO$}{po}%
\newnot{M0@$\MO$}{mo}%
\newnot{delta0@$\delto$}{deltonot}%
\newnot{theta0@$\theto$}{thetonot}%
Compte tenu de la conjugaison sur $F$ des sous-groupes paraboliques
minimaux et de leurs sous-groupes de Levi, il est possible de
choisir $\delto\in\tG(F)$ de sorte que
\[
\theto=\AdtG(\delto)
\]
préserve $\PO$ et $\MO$. On suppose désormais $\delto$ choisi ainsi. Il est uniquement déterminé modulo~$\MO(F)$.

\section{\'Eléments semi-simples ou elliptiques}\label{ellip}

On dit, suivant \cite{KS}*{Section~1.1, p.~13}, qu'un élément $\delta\in\tG$ est quasi semi-simple
si $\AdtG(\delta)$ induit un automorphisme semi-simple de $\G_{der}$.
Cela revient à demander que $\AdtG(\delta)$ préserve une paire de Borel $(B,T)$,
où $B$ est un sous-groupe de Borel et $\T$ un tore maximal dans $B$, définis sur la clôture algébrique.
On dispose pour les éléments d'un tel espace tordu de la décomposition de Jordan:
tout $\delta\in\tG(F)$ s'écrit de manière unique
\[
\delta=s_\delta n_\delta=n_\delta s_\delta
\]
avec $s_\delta$ quasi semi-simple dans $\tG(F)$ et $n_\delta$ unipotent dans $\G(F)$.
On notera $\G^\delta$ le centralisateur de $\delta\in\tG$ dans $\G$ et
$\G_\delta$ la composante neutre de $\G^\delta$:
\[
\G_\delta=(\G^\delta)^0\ptf
\]
On appelle centralisateur stable, noté $\centd$, le sous-groupe engendré par
$\G_\delta$ et $Z_G({\tG})$ le centralisateur de $\tG$ dans $\G$.
Dans le cas non tordu (\ie $\G=\tG$) on a $G_\delta=\centd$ lorsque $\delta$ est semi-simple.
Mais en général le groupe $\centd$ est  non connexe.

On dit qu'un élément $\delta\in\tG(F)$ est elliptique
s'il est quasi semi-simple et si de plus
le tore déployé maximal du centre du centralisateur stable $\centd$
ou, ce qui est équivalent, du centralisateur connexe $G_\delta$
est égal au tore déployé maximal du centralisateur $Z_G({\tG})$ de $\tG$ dans $\G$.
Une condition équivalente est que
\[
\ga_{\G_\delta}=\ga_{\tG}\ptf
\]

\section{Sous-espaces paraboliques}

On dit que $\tP\subset\tG$ est un
\newindex{P@$\protect\tP$}{torP}%
sous-espace parabolique si $\tP$ est le normalisateur dans $\tG$ d'un sous-groupe parabolique
$P$ de $\G$ et s'il est non vide.
Un sous-groupe parabolique étant son propre normalisateur dans $\G$ il en
résulte que $\tP$ est un $P$\hyph espace tordu. Le radical unipotent $\N$ de $P$ est invariant
par $$\theta=\AdtG(\delta)$$ pour tout $\delta\in\tP$ et si $\M$ est un sous-groupe de Levi
de $P$ on peut choisir $\delta$ de sorte que $\tM=\M\delta$ soit un $\M$\hyph espace tordu.
On dit que $\tM$ est un sous-ensemble de Levi de $\tP$ et
on a la décomposition de Levi tordue:
\[
\tP=\tM\N\ptf
\]
On définit l'espace $$\gatP=\ga_{\tM}$$
comme le sous-espace des vecteurs dans $\ga_P=\ga_\M$ fixes
sous l'automorphisme $\theta$ induit
par l'un quelconque des éléments $\delta\in\tM$.
Si $\tP\subset\tQ$
sont deux sous-ensembles paraboliques
on a $\ga_{\tQ}\subset\gatP$ et un supplémentaire canonique
\[
\gatP=\ga_{\tQ}\oplus\ga_{\tP}^{\tQ}\ptf
\]

Rappelons que l'on a choisi $\delto\in\tG(F)$ tel que
\[
\theto=\AdtG(\delto)
\]
préserve $\PO$.
Le sous-ensemble $\tPO=\PO.\delto$ est donc un sous-ensemble parabolique minimal.
Lorsque $P$ est standard, dire que son normalisateur dans $\tG$ est non vide
équivaut à dire que $P$ est $\theto$\hyph stable.
Soient $\M$ le sous-groupe de Levi de $P$ contenant $\M_0$
et $N$ le radical unipotent de $P$. On observe que $\M$ et $\N$ sont $\theto$\hyph invariants
et $\tM=\M\delto$ est un sous-ensemble de Levi de $\tP$.
Soient $\tP\subset\tQ$ deux sous-ensembles paraboliques standard (\cad contenant $\tPO$).
Puisque l'automorphisme $\theto$ préserve $P$ et $\Q$,
il induit une permutation de l'ensemble fini $\Delta_P^\Q$ et donc induit un endomorphisme
d'ordre fini $\ell$ dans $\gaPQ$.
Une racine $\alpha\in\Delta_P^\Q$ définit, par restriction, une forme linéaire $\talpha$
sur $\ga_{\tP}^{\tQ}$ qui ne dépend que de l'orbite de $\alpha$ sous $\theto$.
On observera d'ailleurs que la forme linéaire $\talpha$ sur $\ga_{\tP}^{\tQ}$
est aussi la restriction à cet espace de la moyenne
\[
\frac{1}{\ell}\sum_{r=0}^{\ell-1}\theta_0^r(\alpha)
\]
et on pourra identifier $\talpha$ à cette moyenne sur l'orbite.
On note
\[
\Delta_{\tP}^{\tQ}
\]
\newindex{DeltaPQ@$\Delta_{\protect\tP}^{\protect\tQ}$}{tordeltapq}%
l'ensemble de ces orbites (ou des formes linéaires associées).
On définit de même $\tvpi$ pour $\vpi\in\hDelta_P^\Q$ et on note
\[
\hDelta_{\tP}^{\tQ}
\]
l'ensemble de ces orbites.
Les lemmes suivants sont la clef de l'extension au cas tordu de
la combinatoire:

\begin{lemme}\label{bijpartor}
L'application
\[
\tP\mapsto\Delta_{\tPO}^{\tP}
\]
est une bijection entre l'ensemble des sous-ensembles paraboliques standard de $\tG$
et l'ensemble des parties de $\Delta_{\tPO}^{\tG}$.
\end{lemme}

\begin{proof}
Il suffit d'observer que $\Delta_{\tPO}^{\tP}$ est un ensemble d'orbites
sous $\theto$ dans $\Delta_{\PO}^\G$ et d'invoquer le lemme~\ref{bijpar}.
\end{proof}

\begin{lemme}\label{robtuses}
L'ensemble $\Delta_{\tP}^{\tQ}$
est une base obtuse et $\hDelta_{\tP}^{\tQ}$ une base aigüe
du dual de $\ga_{\tP}^{\tQ}$.
\end{lemme}

\begin{proof}
On rappelle que l'on peut représenter
les éléments de
\[
\Delta_{\tP}^{\tQ}\Qquad{et}\hDelta_{\tP}^{\tQ}
\]
par des moyennes sur les orbites correspondantes
ce qui permet le calcul des produits scalaires.
Les assertions résultent alors du lemme~\ref{obtu}.
\end{proof}

\begin{lemme}\label{deter}
Le sous-espace $\gatP$ détermine $\ga_P$ et
le centralisateur de $\gatP$ et de $\gaP$ dans $\weyl$ coïncident.
\end{lemme}

\begin{proof}
Il suffit d'observer que le point $X_P\in\gaP$ défini par
\[
X_P=\sum_{\vpi\in\hDelta_P} \vpi^\vee
\]
appartient à l'intersection de $\gatP$ et de
la chambre définie par $P$ dans $\ga_P$.
Son centralisateur dans $\weyl$, est
le groupe de Weyl $\weyl^\M$
où $\M$ est le sous-groupe de Levi de $P$.
\end{proof}

\section{Chambres et facettes: cas tordu}\label{chambresfacettes}

Le quotient du normalisateur de $\MO$ dans $\tG$ par $\MO$
est l'ensemble de Weyl de $\tG$ et sera noté $\weyl^{\tG}$
ou simplement $\tweyl$. On a
\[
\tweyl=\weyl\rtimes\theto\ptf
\]
Soit $\tM$ et $\tM'$ deux sous-ensembles de Levi.
Pour alléger un peu les notations on écrira parfois $\tga$ pour $\ga_{\tM}$
et $\tga'$ pour $\ga_{\tM'}$.
On note
$\weyl(\tga,\tga')$ l'ensemble des restrictions à $\tga$
des applications induites par des $s\in\weyl$ qui induisent un isomorphisme
\[
\tga\to\tga'\ptf
\]

\begin{lemme}\label{thetastab}
Soient $\tM$ et $\tM'$ deux sous-ensembles de Levi standard.
L'ensemble $\weyl(\tga,\tga')$ est le sous-ensemble des points fixes sous
$\theto$ dans $\weyl(\ga,\ga')$ \cad que $\weyl(\tga,\tga')$ est
en bijection avec l'ensemble des
$s\in\weyl^\G$ tels que
\begin{enumerate}[(i)]
\item $s(\ga)=\ga'$
\item $s\alpha>0$ pour toute $\alpha\in\Delta_{\PO}^{\M}$
\item $\theto(s)=s$.
\end{enumerate}
La condition~\textup{(ii)}
est équivalente à la condition
\begin{enumerate}[(i$'$)]\addtocounter{enumi}{1}
\item $s\moins\alpha>0$ pour toute $\alpha\in\Delta_{\PO}^{\M'}$.
\end{enumerate}
\end{lemme}

\begin{proof}
D'après le lemme~\ref{deter}, un élément
$s\in\weyl(\tga,\tga')$ est la restriction d'un unique élément, encore noté $s$, dans $\weyl(\ga,\ga')$.
D'après le lemme~\ref{weylga}
il est représenté par un unique élément dans $\weyl$ encore noté $s$ de longueur minimale
dans sa classe modulo $\weyl^{\M}$,
ce qui équivaut à demander que $s\alpha>0$ pour toute
$\alpha\in\Delta_{\PO}^{M}$ ou, ce qui est équivalent, que $s\moins\alpha>0$ pour toute $\alpha\in\Delta_{\PO}^{\M'}$.
Mais $\theto(s)$ a les mêmes propriétés et donc $s=\theto(s)$.
\end{proof}

On notera $\weyl(\ga_{\tM})$ l'union disjointe des $\weyl(\ga_{\tM},\ga_{\tM'})$
lorsque $\tM'$ parcourt l'ensemble des sous-ensembles de Levi standard.
On dispose dans $\ga_{\tM}$ des chambres de Weyl complémentaires
des hyperplans définis par les orbites des racines.
Soit $s\in\weyl(\ga_{\tM})$, on note $\Delta(\tM,s)$ l'ensemble
des projections $\talpha$ sur $\ga_{\tM}$ des $\alpha\in
s\moins(\Delta_{\PO})$ qui sont non nulles.
On définit une chambre $C_{\tM}(s)$ dans $\ga_{\tM}$ par les inégalités
\[
\talpha(H)>0\Quad{pour}\talpha\in\Delta(\tM,s)\ptf
\]
La chambre $C_{\tM}(s)$ est l'ensemble des points fixes
sous $\theto$ dans $C_\M(s)$.

\begin{lemme}\label{weylbijc}
Il y a une bijection naturelle entre les trois ensembles suivants
\begin{enumerate}[(i)]
\item L'ensemble $C_{\tM}$ des chambres de Weyl dans $\ga_{\tM}$
\item L'ensemble $\Parab(\tM)$ des sous-ensembles paraboliques $\tP$ admettant $\tM$ comme sous-ensemble de Levi
\item L'ensemble $\weyl(\ga_{\tM})$.
\end{enumerate}
\end{lemme}

\begin{proof}
Considérons $s\in\weyl(\ga_{\tM})$ alors
$C_{\tM}(s)$ est la projection sur $\ga_{\tM}$ de la chambre $C_\M(s)$ dans
$\ga_\M$ ou, si on préfère, l'intersection de $C_\M(s)$ et $\ga_{\tM}$.
Une telle chambre dans $\ga_\M$ a une intersection non triviale avec
$\ga_{\tM}$ si et seulement si elle est $\theto$\hyph invariante.
On conclut en invoquant le lemme~\ref{thetastab}.
\end{proof}

\section{Combinatoire: extension au cas tordu}

Soit $\tM$ un sous-ensemble de Levi, $\tP$, $\tQ$ et $\tR$
des sous-ensembles paraboliques.
On définit des fonctions caractéristiques $$\tau_{\tP}^{\tQ}\Quad{et} \htautPQ$$ de cônes
dans $\ga_{\tP}^{\tQ}$,
ainsi que des fonctions
\[
\phi_{\tP}^{\tQ,\tR},\qquad \Gamma_{\tP}^{\tQ}\Qquad{et}\Gamma_{\tM}^{\tQ}
\]
en remplaçant les bases
$\Delta_P^\Q$ et $\hDelta_P^\Q$ dans $\gaPQ$ par $\Delta_{\tP}^{\tQ}$ et $\hDelta_{\tP}^{\tQ}$ dans $\ga_{\tP}^{\tQ}$.
Une famille $\M$\hyph orthogonale $\XX$ définit une famille $\tM$\hyph orthogonale
en définissant $X_{\tP}$ pour $\tP\in\FF(\tM)$
comme la projection de $X_P$ sur $\ga_{\tP}$.
La notion de \tGM-famille est aussi définie de manière naturelle.
Une \GM-famille étant donnée, on lui associe une \tGM-famille comme suit:
soit $\tP\in\FF(\tM)$ on définit $c(\Lambda,\tP)$ comme la restriction
aux $\Lambda\in\ima\ga_{\tP}^*$ de $c(\Lambda,P)$.

\begin{lemme}\label{tgmfour}
Considérons une \GM-famille définie par transformée de Fourier
\[
c(\Lambda,P)=\int_{\gHM}\ee^{\iota_P(\Lambda)(\XX)}
\varphi(\XX)\dd\XX\ptf
\]
La \tGM-famille associée est définie par:
\[
c(\Lambda,\tP)=\int_{\gHM}\ee^{\iota_{\tP}(\Lambda)(\XX)}
\varphi(\XX)\dd\XX
\]
pour $\Lambda\in\ima\ga_{\tP}^*$.
\end{lemme}

\begin{lemme}\label{recomb}
L'analogue des assertions~\ref{adja} à~\ref{decompGM} sont encore valables
pour les fonctions
pour $\tau_{\tP}^{\tQ}$, $\htautPQ$, $\Gamma_{\tP}^{\tQ}$ et $\Gamma_{\tM}^{\tQ}$
et les \tGM-familles.
\end{lemme}

\begin{proof}
Pour traiter le cas tordu il
convient de remplacer partout les racines simples par leurs orbites sous $\theto$ et
les espaces vectoriels par leur sous-espaces de points fixes sous $\theto$.
À ceci près, et compte tenu des lemmes~\ref{bijpartor} à~\ref{weylbijc},
les preuves s'étendent \emph{verbatim} au cas tordu.
\end{proof}

En particulier, le lemme~\ref{recomb} fournit les énoncés suivants:

\begin{proposition}\label{envconvt}
Considérons une famille orthogonale $\XX$ et posons
\[
\Gamma_{\tM}(H,\XX)=\sum_{\tP\in\FF(\tM)}(-1)^{a_{\tP}-a_{\tG}}\,
\htau_{\tP}(H-X_P)\ptf
\]
Supposons que la famille orthogonale $\XX$ est régulière.
Alors\textup, la fonction
\[
H\mapsto \Gamma_{\tM}(H,\XX)
\]
est la fonction caractéristique de l'ensemble des $H$ dont la projection sur
$\ga_{\tM}^{\tG}$ appartient à l'enveloppe convexe des projections des $X_s$
avec $s\in\weyl(\ga_{\tM})$ ou\textup, ce qui est équivalent\textup, l'enveloppe convexe des projections des $X_P$
pour $\tP\in\Parab(\tM)$.
\end{proposition}

\begin{proposition}\label{GammaHXt}
On a
\begin{gather}
\tau_{\tP}^{\tR}(H)=\sum_{\tP\subset\tQ\subset\tR}\Gamma_{\tP}^{\tQ}(H,X)\tau_{\tQ}^{\tR}(H-X)\label{eq2.1} \\
\htau_{\tP}^{\tR}(H-X)=\sum_{\tP\subset\tQ\subset\tR}(-1)^{a_{\tQ}-a_{\tR}}\,
\htau_{\tP}^{\tQ}(H)\Gamma_{\tQ}^{\tR}(H,X)\ptf\label{eq2.2}\\
\Gamma_{\tP}^{\tR}(H,X+Y)=\sum_{\tP\subset\tQ\subset\tR}\Gamma_{\tP}^{\tQ}(H,X)
\Gamma_{\tQ}^{\tR}(H-X,Y)\label{eq2.3}
\end{gather}
et
\begin{equation}
\Gamma_{\tM}^{\tR}(H,\XX+\YY)=\sum_{\tQ\in\FF^{\tR}(\tM)}
\Gamma_{\tM}^{\tQ}(H,\XX)\Gamma_{\tQ}^{\tR}(H-X_{\tQ}, Y_{\tQ})
\label{eq2.4}
\end{equation}
Enfin, si $e$ est produit de deux \tGM-familles, on a
\begin{equation}
e_{\tM}^{\tR}(\Lambda)=\sum_{\tQ\in\FF^{\tR}(\tM)}c_{\tM}^{\tQ}(\Lambda)d_{\tQ}^{\tR}(\Lambda)
\ptf\label{eq2.5}
\end{equation}
\end{proposition}

\begin{proof}
Les équations~\eqref{eq2.1}, \eqref{eq2.2} et~\eqref{eq2.3} sont les variantes tordues des
équations~\eqref{eq1.1b}, \eqref{eq1.2b} et~\eqref{eq1.3b} du lemme~\ref{GammaHX} respectivement.
L'équation~\eqref{eq2.4} est la variante tordue du lemme~\ref{convol}.
Enfin~\eqref{eq2.5} est la variante tordue du lemme~\ref{decompGM}.
\end{proof}


\section{Volumes de convexes et polynômes}\label{volconvpol}

On dispose de l'ensemble des racines réduites $\Rac$. Plus généralement, soit $\M$ un sous-groupe de Levi
semi-standard, on notera $\Rac_\M$ l'ensemble des racines réduites défini par la projection
sur $\ga_\M^*$ des racines de $\G$ dans $\gao$. On prendra garde que ce n'est pas en général un système de racines.
Pour $P\in\Parab(\M)$ on notera $\Rac_P$ le sous-ensemble des racines de $\Rac_\M$ positives sur la chambre
associée à $P$.

On appellera famille $M$\hyph radicielle la donnée de nombres r\'eels $z_\beta$ pour chaque $\beta\in\Rac_M$ et on pose 
\[
X_P=\sum_{\beta\in\Rac_P} z_\beta\beta^\vee\ptf
\]
 Si $P$ et $Q$ sont adjacents
le mur étant défini par $\gamma$ on a
\[
X_P-X_Q=\sum_{\beta\in\Rac_P} z_\beta\beta^\vee-\sum_{\beta\in\Rac_\Q} z_\beta\beta^\vee
=(z_\gamma+z_{-\gamma})\gamma^\vee
\]
car $\Rac_P\cap\Rac_\Q$ est le complémentaire de $\gamma$  dans $\Rac_P$ et de $-\gamma$ dans $\Rac_Q$.
La famille des $X_P$ définie à partir de la collection des $z_\beta$ est donc une famille $M$\hyph orthogonale.
Notons $\gZ_M$ l'espace vectoriel des familles de scalaires  $\zz=\{z_\beta\}$ pour $\beta\in\Rac_M$.
L'application qui à $\zz$ associe la famille des 
\[
X_P=\sum_{\beta\in\Rac_P} z_\beta\beta^\vee
\]
est une application linéaire 
\[
j\colon \gZ_M\to\gHM
\]
dont l'image sera notée $\gY$.

On rappelle que l'on a introduit en \ref{LaplaceGamma} des polyn\^omes $\gamma_\M$.
Le lemme suivant est une variante des lemmes~7.1 et~7.2 de \cite{AeisII}.

\begin{lemme} \label{volpol} 
Soit 
$\XX$ la famille orthogonale associée à une famille radicielle $\zz$.
Le polynôme $$\gamma_\M(\XX)=\gamma_\M\circ j(\zz)$$
peut s'écrire sous la forme
\[
\gamma_\M\circ j(\zz)=\sum_F  c_F\prod_{\beta\in F} z_\beta
\] 
où $F$ parcourt l'ensemble des bases de $\ga_\M^G$ formées de racines réduites pour $M$
et où
\[
c_F=\vol(F)
\]
est le volume du parallélépipède engendré par $F$.
\end{lemme}

\begin{proof} 
Par définition
\[
\gamma_\M\circ j(\zz)=\lim_{\LL\to0}\sum_{P\in\Parab(M)}
\ee^{\LL(X_P)}\,\epsilon_P^\G(\Lambda)
\qquad\text{avec $X_P=\sum_{\beta\in\Rac_P} z_\beta\beta^\vee$}\ptf
\]
Donc
\[
\frac{\partial}{\partial z_\beta}\gamma_\M\circ j(\zz)
=\lim_{\LL\to0}\sum_{P\in\Parab(M)}
\frac{\partial}{\partial z_\beta}
\ee^{\LL(X_P)}\,\epsilon_P^\G(\Lambda)
\]
est égal à
\[
\lim_{\LL\to0}\sum_{\{P\in\Parab(M)\mid \beta\in\Rac_P\}}\LL(\beta^\vee)
\ee^{\LL(X_P)}\,\epsilon_P^\G(\Lambda)\ptf
\]
On rappelle que
\[
\epsilon_P^\G(\Lambda)=\vol(\check\Delta_P^\G) \prod_{\alpha\in\Delta_P^\G} \Lambda (\alpha^\vee)\moins\ptf
\]
 Si nous supposons que $\LL=\LL_0+t\beta$ avec $\LL_0$ générique dans $\ga_L$
 l'orthogonal de $\beta^\vee$ et $t\ne0$ alors 
\[
\frac{\partial}{\partial z_\beta}\gamma_\M\circ j(\zz)
=\lim_{\LL_0\to0}\lim_{t\to0}\sum_{\{P\in\Parab(M)\mid\beta\in\Rac_P\}}t\beta(\beta^\vee)
\ee^{\LL(X_P)}\,\epsilon_P^\G(\LL_0+t\beta)\ptf
\]
Maintenant, si $\beta\in\Delta_P\subset\Rac_P$, on a
\[
\epsilon_P^\G(\Lambda_0+t\beta)=\frac{\lvert\beta^\vee\rvert}{t\beta(\beta^\vee)}\epsilon_\Q^\G(\Lambda_0)
\]
où $Q$ est le sous-groupe parabolique tel que 
\[
\Delta_P=\Delta_\Q\cup{\beta}
\]
et où $\lvert \beta^\vee\rvert$ est la longueur de $\beta^\vee$. 
On a utilisé que
\[
\vol(\check\Delta_P^G)=\lvert \beta^\vee\wedge\alpha^\vee_1\wedge\dots\wedge\alpha^\vee_r\rvert
=\lvert \beta^\vee\rvert\,\lvert\widebar\alpha^\vee_1\wedge\dots\wedge\widebar\alpha^\vee_r\rvert
=\lvert\beta^\vee\rvert\vol(\check\Delta_Q^G)
\]
où $\widebar\alpha$ est la projection de $\alpha$ sur l'orthogonal de $\beta$.
Par contre $\epsilon_P^\G(\Lambda_0+t\beta)$ a une limite finie lorsque $t\to0$
si  $\beta\notin\Delta_P$.
Donc, si on note $L$ le sous-groupe de Levi  tel que $\ga_L$ soit
l'orthogonal de $\beta^\vee$ alors
\[
\frac{\partial}{\partial z_\beta}\gamma_\M\circ j(\zz)=\lvert\beta^\vee\rvert\,
\lim_{\LL_0\to0}\sum_{Q\in\Parab(L)}
\ee^{\LL_0(X_\Q)}\,\epsilon_\Q^\G(\LL_0)=\lvert\beta^\vee\rvert\,\gamma_L(\LL_0)\ptf
\]
 Maintenant les $X_Q$ sont indépendants de la variable $z_\beta$.
Donc $\gamma_\M\circ j(\zz)$ est somme monômes  de degré${}\le1$ en chaque variable $z_\beta$. Plus précisément
on a
\[
\gamma_\M\circ j(\zz)=z_\beta\lvert\beta^\vee\rvert\,\gamma_L(\LL_0)+\text{des termes ne contenant pas $z_\beta$}\ptf
\]
Chaque racine $\gamma\in\Rac_P-\{\beta\}$ se projette en un multiple d'une racine réduite $\delta$ pour $Q$ et on a donc
\[
X_Q=\sum_{\gamma\in\Rac_P-\{\beta\}} z_\gamma\widebar\gamma^\vee=
\sum_{\delta\in\Rac_Q} y_\delta\delta^\vee\qquad\text{avec $y_\delta=\sum_{\gamma\to\delta} n_\gamma z_\gamma$}
\]
et donc, si $\delta_1,\cdots,\delta_r$ est une base de $\ga_L$ formée d'éléments de $\Rac_L$ on a
\[
y_{\delta_1}\dotsm y_{\delta_r}\lvert\delta_1^\vee\wedge\dots\wedge\delta_r^\vee\rvert
=\sum z_{\gamma_1}\dotsm z_{\gamma_r}\lvert \widebar\gamma_1^\vee\wedge\cdots\wedge\widebar\gamma_r^\vee\rvert
\]
la somme portant sur les familles $\gamma_1,\dots,\gamma_r$ d'éléments de $\Rac_M$
se projetant sur des multiples de $\delta_1,\dots,\delta_r$. Elle sont donc telles que
\[
\{\beta,\gamma_1,\cdots,\gamma_r\}
\]
est une base de $\ga_\M$.
On voit alors, par récurrence sur le nombre de racines, que
$\gamma_\M\circ j(\zz)$ est la somme des monômes 
\[
c_F \prod_{\beta\in F}z_{\beta}
\]
où $F$ est une base de $\ga_\M^G$ formé de coracines réduites pour $M$
et où 
\[
c_F=\vol(F)=\lvert\beta^\vee_0\wedge\dots\wedge\beta^\vee_r\rvert
\]
 est le volume du parallélépipède engendré par les $\beta_i^\vee \in F$.
\end{proof}

Considérons une famille radicielle $\{z_\beta\}$ et soit $\XX$ la famille orthogonale associée.
Soit $\tL$ un sous-ensemble de Levi de $\tG$. Pour tout $\tP$ de Levi $\tM$ on définit
$X_{\tP}$ comme la projection de $X_P$ sur $\ga_{\tP}=\ga_{\tL}$. Les $X_{\tP}$ définissent
une famille $\tM$\hyph orthogonale.

\begin{lemme} \label{volpone}
\[
\gamma_{\tL}^{\tG}\circ j(\zz)=\sum_F  c_{F_{\tL}}\prod_{\beta\in F} z_\beta
\] 
où la somme porte sur les familles $F$ de racines $\gamma$ telles que l'ensemble $F_{\tL}$ de
leur projections $\widebar\gamma$ sur $\ga_{\tL}^{\tG}$ soit une base de cet espace et où
\[
c_{F_{\tL}}=\lvert\widebar\gamma_1^\vee\wedge\dots\wedge\widebar\gamma_r^\vee\rvert\ptf
\]
\end{lemme}

\begin{proof} 
D'après la variante tordue du lemme~\ref{volpol} on sait que
\[
\gamma_{\tL}^{\tG}\circ j(\zz)=\sum_B  c_B\prod_{\delta\in B} y_\delta
\]
où $B$ parcourt les bases de l'ensemble $\Rac_{\tL}$. Comme dans le lemme précédent on observe que
chaque racine $\gamma\in\Rac_P$ se projette en un multiple d'une racine réduite $\delta$ pour $\tP$ et on a donc
\[
X_{\tP}=\sum_{\gamma\in\Rac_P} z_\gamma\widebar\gamma^\vee=
\sum_{\delta\in\Rac_{\tP}} y_\delta\delta^\vee\Quad{avec}y_\delta=\sum_{\gamma\to\delta} n_\gamma z_\gamma
\]
et donc, si $\delta_1,\dots,\delta_r$ est une base $B$ de $\ga_{\tP}$ formée d'éléments de $\Rac_{\tP}$ on a
\[
y_{\delta_1}\dots y_{\delta_r}\lvert\delta_1^\vee\wedge\dots\wedge\delta_r^\vee\rvert
=\sum_{\{F\mid F_{\tL}=B\}} z_{\gamma_1}\dotsm z_{\gamma_r}\lvert\widebar\gamma_1^\vee\wedge\dots\wedge\widebar\gamma_r^\vee\rvert
\]
la somme portant sur les familles $\gamma_1,\dots,\gamma_r$ d'éléments de $\Rac_P$
se projetant sur des multiples de $\delta_1,\dots,\delta_r$.
\end{proof}


On dira qu'une \GM\hyph famille $c$ est radicielle si elle est obtenue de la manière suivante:
\[
c(\Lambda,P)=(\exta\circ\iota_P)(\Lambda)
\]
et
\[
\exta=\extb\circ j^*
\]
où 
\[
j^*\colon\gH_M^*\to\gZ_M^*
\]
est l'application linéaire duale de l'application $j$ définie plus haut.

Le corollaire suivant reproduit et étend au cas tordu 
le résultat d'Arthur (\cite{AeisII}*{Corollary~7.3} pour les \GM\hyph familles radicielles.
C'est un cas particulier de résultats de Finis et Lapid \cite{FL}.

\begin{corollaire}\label{gmrad}
Soit $\{c(\Lambda,P)\}$ une \GM\hyph famille de la forme
\[
c(\Lambda,P)=(\extb\circ j^*\circ\iota_P)(\Lambda)\ptf
\]
On a
\[
c_{\tL}^{\tG}(0)=(D_{\tL}^{\tG}\extb)(0)
\] 
où $D_{\tL}^{\tG}$ est  l'opérateur différentiel déduit du polynôme $\gamma_{\tL}^{\tG}\circ j$ par transformation de Fourier.
C'est une combinaison linéaire de monômes différentiels
produits de dérivées partielles par rapport aux variables $z_\beta$
où chaque variables intervient au plus une fois\textup:
\[
D_{\tL}^{\tG}=\sum_{F}\, c_{F_{\tL}}D_F
\qquad\text{avec 
$D_F=\prod_{\beta^\vee\in F} \partial_{z_\beta}$}
\]
le produit des dérivations portant sur les  coracines dans $F$.
\end{corollaire}

\begin{proof} 
Comme ci-dessus il suffit de traiter le cas où $\exta=\extb\circ j^*$ avec $\extb$ à support compact.
Sa transformée de Fourier $\mm$ est une mesure à décroissance rapide de support contenu dans $\gY$.
D'après l'équation~\eqref{eq1.4f'} du lemme~\ref{gmorth}  on a
\[
c_{\tL}^{\tG}(0)=\int_{\gHM}\gamma_{\tL}^{\tG}(\XX)\dd\mm(\XX)\ptf
\]
Si on note $\extc$ de composé de $\extb$ avec l'injection 
\[
\gY^*\to\gZ^*
\]
on a donc
\[
c_{\tL}^{\tG}(0)
=\int_{\gY}\gamma_{\tL}^{\tG}(\yy)\widehat\extc(\yy)\dd\yy=\int_{\gZ}\gamma_{\tL}^{\tG}\circ j(\zz)\widehat\extb(\zz)\dd\yy
\ptf
\]
L'assertion résulte alors immédiatement du lemme~\ref{volpone} par transformation de Fourier.
\end{proof}

\section[Les fonctions $\sQR$ et $\tsQR$]{\mathversion{bold}Les fonctions $\sQR$ et $\tsQR$}

Soit $\Q$ un sous-groupe parabolique standard. On notera
$\Qp$ le sous-groupe parabolique standard dont le sous-groupe
de Levi admet pour racines simples les éléments des orbites
sous $\theto$ des racines dans $\Delta_{\PO}^\Q$.
On notera $\Qm$ le sous-groupe parabolique standard dont le sous-groupe
de Levi admet pour racines simples les $\alpha\in\Delta_{\PO}^\Q$
telles que l'orbite de $\alpha$ sous $\theto$ soit toute entière
contenue dans $\Delta_{\PO}^\Q$. Les sous-groupes paraboliques
$\Qp$ et $\Qm$ sont stables sous $\theto$ et on note
$\tQp$ (resp.~$\tQm$) les sous-ensembles paraboliques
associés.

\begin{lemme}\label{qplus}
Soient $\Q$ et $\R$ deux sous-groupes paraboliques standard.
Il existe un sous-ensemble parabolique \textup(standard\textup) $\tP$ avec
\[
\Q\subset P\subset\R
\]
\newindex{Q+@$\Qp$}{Qplus}%
\newindex{R-@$\Rm$}{Rmoins}%
si et seulement si $\Qp\subset\Rm$. Dans ce cas on a
\[
\Q\subset\Qp\subset P\subset\Rm\subset\R\ptf
\]
En d'autres termes $\tQp$ est le plus petit sous-ensemble parabolique
$\tP$ avec $\Q\subset P\subset\R$
et $\tRm$ est le plus grand.
\end{lemme}

\begin{proof}
D'après le lemme~\ref{bijpartor}
les sous-ensembles paraboliques $\tP$ avec $\Q\subset P\subset\R$
sont en bijection avec les sous-ensembles $\Delta_{\PO}^P$ vérifiant
\[
\Delta_{\PO}^\Q\subset\Delta_{\PO}^P\subset\Delta_{\PO}^\R
\]
et formés d'orbites sous $\theto$. Le lemme est alors conséquence des définitions de $\Qp$ et $\Rm$.
\end{proof}

\begin{lemme}\label{qplusa}
Supposons $\Qp\subset\Rm$.
Le sous-espace $\tgaQR$
des $\theto$\hyph invariants dans $\ga_\Q^\R$
est égal au sous-espace $\ga_{\tQp}^{\tRm}$.
\end{lemme}

\begin{proof}
Tout d'abord il est clair que $\tgaQR$ contient $\ga_{\tQp}^{\tRm}$.
Réciproquement,
considérons $H\in\tgaQR$. Pour $\alpha\in\Delta_{\PO}^\Q$ on a
\[
\alpha(H)=\talpha(H)=0
\]
et ces $\talpha$ forment l'ensemble $\Delta_{\tPO}^{\tQp}$.
On en déduit que $\tgaQR\subset\ga_{\tQp}^{\tG}$.
De même pour $\vpi\in\hDelta_\R^\G$ on a
\[
\vpi(H)=\tvpi(H)=0
\]
d'où on déduit que $\tgaQR\subset\ga_{\tPO}^{\tRm}$.
\end{proof}

On considère deux sous-groupes paraboliques $\Q\subset\R$.
On définit $\sQR$ comme la fonction caractéristique de l'ensemble des $H$
\newindex{sigmaQR@$\sQR$}{sqr}%
tels que
\begin{enumerate}[(i)]
\item $\alpha(H)>0$ pour tout $\alpha\in\Delta_\Q^\R$
\item $\alpha(H)\le0$ pour tout $\alpha\in\Delta_\Q^\G-\Delta_\Q^\R$
\item $\vpi(H)>0$ pour tout $\vpi\in\hDelta_\R$
\end{enumerate}
On observera que, d'après le lemme~\ref{cone}, $\sQR(H)=1$ implique que
$\vpi(H)>0$ pour tout $\vpi\in\hDelta_\Q$
et il est immédiat de vérifier que, pour $P$ et $\Q$ fixés,
\[
\sum_{\{\R\mid P\subset\R\}}\sQR =\tau_\Q^P\htau_P
\]
(voir le lemme~\ref{repart} pour une preuve dans un cas plus général).

Nous aurons aussi besoin de la variante tordue de cette fonction caractéristique de cône.
On suppose $\Qp\subset\Rm$ et soit $\tP$ un sous-ensemble
parabolique tel que
\[
\Q\subset P\subset\R\ptf
\]
On définit ${}_{\tP}\sQR$ comme la fonction caractéristique de l'ensemble des $H$
tels que
\begin{enumerate}[(i)]
\item $\alpha(H)>0$ pour tout $\alpha\in\Delta_\Q^\R$
\item $\alpha(H)\le0$ pour tout $\alpha\in\Delta_\Q-\Delta_\Q^\R$
\item $\tvpi(H)>0$ pour tout $\tvpi\in\hDelta_{\tP}$
\end{enumerate}

\begin{lemme}\label{indepb}
La fonction caractéristique ${}_{\tP}\sQR$ est indépendante de $\tP$.
\end{lemme}

\begin{proof}
Considérons $H$ tel que
$\alpha(H)>0$ pour tout $\alpha\in\Delta_\Q^\R$ et $\tvpi(H)>0$ pour tout $\tvpi\in\hDelta_{\tP}$.
En particulier $\widebar\alpha(H)>0$ pour tout $\widebar\alpha\in\Delta_{\Qp}^P$.
Ceci implique $\talpha(H)>0$ pour tout $\talpha\in\Delta_{\tQp}^{\tP}$.
On a donc
\[
\tau_{\tQp}^{\tP}(H)\htau_{\tP}(H)=1
\]
et il résulte alors de la variante tordue du lemme~\ref{cone} que
\[
\htau_{\tQp}(H)=1
\]
\cad que
$\tvpi(H)>0$ pour tout $\tvpi\in\hDelta_{\tQp}$.
On a donc ${}_{\tP}\sQR={}_{\tQp}\sQR$.
\end{proof}

Si $\Qp\subset\Rm$ et compte tenu du lemme~\ref{indepb} il est loisible de poser
\[
\tsQR={}_{\tP}\sQR
\]
\newindex{sigmaQR@$\protect\tsQR$}{tsqr}%
où $\tP$ est l'un quelconque des sous-ensembles paraboliques avec
$\Q\subset P\subset\R$.
Par convention $\tsQR=0$ si $\Qp\nsubset\Rm$.

\begin{lemme}\label{nullos}
Si $\Q=\R$ alors $\tsQR=0$ sauf si $\Q=\G$ auquel cas $\widetilde\sigma_\G^\G=1$.
\end{lemme}

\begin{proof}
Considérons $\tP$ avec $\Q=P=\R$.
Les hypothèses impliquent
\[
\alpha(H)\le0\qquad\forall\,\alpha\in\Delta_P
\]
alors que
\[
\tvpi(H)>0\qquad\forall\,\tvpi\in\hDelta_{\tP}
\]
ce qui, d'après le lemme~\ref{precone},
impose $\Delta_P=\varnothing$ et donc $P=\G$.
\end{proof}

\begin{lemme}\label{repart}
Pour tout $Q\subset P$ on a l'\'egalit\'e:
$$\sum_{\R\supset P}{\tsQR}=\tau_\Q^P\htautP$$
\end{lemme}

\begin{proof}
On observe tout d'abord que les supports des diverses
$\tsQR$ sont disjoints lorsque $\R$ varie.
La fonction $\tau_\Q^P\htautP$ est la fonction caractéristique des $H$ tels que
\[
\alpha(H)>0\quad\forall\alpha\in\Delta_\Q^P\Qquad{et}\tvpi(H)>0\quad\forall\tvpi\in\hDelta_{\tP}\ptf
\]
Fixons $H$ dans le support de $\tau_\Q^P\htautP$.
Soit $\R$ le sous-groupe parabolique avec $\R\supset\Q$ et tel que
\[
\Delta_\Q^\R=\{\alpha\in\Delta_\Q\mid\alpha(H)>0\}
\]
alors $\R\supset P$ et par définition ${}_{\tP}\sQR(H)=\tsQR(H)=1$.
\end{proof}

\begin{lemme}\label{croiss}
Considérons $H\in \gao$ de la forme $H=H_1+H_2$ avec $H_1\in\gao$
et
\begin{equation}
H_2\in\ga_\R^\G\tag{i}\label{eq2.i}
\end{equation}
ou bien
\begin{equation}
H_2\in\ga_{\tRm}^{\tG}\ptf\tag{ii}\label{eq2.ii}
\end{equation}
En particulier\textup, dans le second cas\textup, on a $H_2=\theto(H_2)$.
Supposons que
\[
\tsQR(H)=1\ptf
\]
Il existe une constante $c$ telle que
\[
\lVert H_2\rVert \le c\lVert H_1\rVert \ptf
\]
\end{lemme}

\begin{proof}
Supposons $H_2$ de la forme~\eqref{eq2.i}.
On peut écrire
\[
H_2=\sum_{\alpha\in\Delta_\Q^\G}a_\alpha\vpi_\alpha
\qquad\text{avec $a_\alpha=\alpha(H_2)$}\ptf
\]
Par hypothèse $a_\alpha=0$ pour $\alpha\in\Delta_\Q^\R$.
Par ailleurs pour $\alpha\in\Delta_\Q^\G-\Delta_\Q^\R$ on a $\alpha(H)\le0$ et donc
pour un tel $\alpha$, on a
\begin{equation}
a_\alpha=\alpha(H_2)=\alpha(H)-\alpha(H_1)\le -\alpha(H_1)
\le c_1\lVert H_1\rVert \label{eq2.1a}
\end{equation}
pour une constante $c_1>0$.
On en déduit que
\begin{equation}
a_\alpha\le c_1\lVert H_1\rVert\qquad\text{pour $\alpha\in\Delta_\Q^\G$}\ptf
\tag{\ref{eq2.1a}$'$}\label{eq2.1a'}
\end{equation}
Par hypothèse
\[
\tvpi(H)=\tvpi(H_2)+\tvpi(H_1)>0
\]
pour $\tvpi\in\hDelta_{\tP}$ et donc il existe $c_2>0$ telle que
\begin{equation}
\tvpi(H_2)>-c_2\lVert H_1\rVert\ptf\label{eq2.2a}
\end{equation}
Maintenant, si $\tvpi_\alpha\in\hDelta_{\tRm}$ est la moyenne sur l'orbite de $\vpi_\alpha\in\hDelta_\R$
on a
\[
\tvpi_\alpha(H_2)=a_\alpha\langle \vpi_\alpha,\tvpi_\alpha\rangle
+\sum_{\substack{\beta\in\Delta_\Q^\G\\ \alpha\ne\beta}}a_\beta\langle\vpi_\beta,\tvpi_\alpha\rangle\ptf
\]
On remarque que les produits scalaires $\langle\vpi_\beta,\tvpi_\alpha\rangle$ sont tous positifs ou nuls.
Compte tenu de~\eqref{eq2.1a'} et~\eqref{eq2.2a} on a
\[
- c_2\lVert H_1\rVert <
\tvpi_\alpha(H_2)\le rc_1\lVert H_1\rVert +
a_\alpha\langle \vpi_\alpha,\tvpi_\alpha\rangle
\]
où $r$ est le rang de $\G$. Comme $\langle \vpi_\alpha,\tvpi_\alpha\rangle$ est strictement positif on obtient
\begin{equation}
a_\alpha\ge -c_3\lVert H_1\rVert \label{eq2.3a}
\end{equation}
Les assertions du lemme se déduisent immédiatement de~\eqref{eq2.1a'} et~\eqref{eq2.3a} pour le cas~\eqref{eq2.i}.
Dans le cas~\eqref{eq2.ii} on peut encore écrire
\[
H_2=\sum_{\alpha\in\Delta_\Q^\G}a_\alpha\vpi_\alpha
\qquad\text{avec $a_\alpha=\alpha(H_2)$}
\]
mais cette fois on a $a_\alpha=0$ seulement pour $\alpha\in\Delta_\Q^{\Rm}$.
Par ailleurs, comme ci-dessus, on a~\eqref{eq2.1a}
pour $\alpha\in\Delta_\Q^\G-\Delta_\Q^\R$ et donc on a
\begin{equation}
a_\alpha\le c_1\lVert H_1\rVert \qquad\text{pour $\alpha\in\Delta_\Q^\G-\Delta_\Q^\R$}\ptf
\tag{\ref{eq2.1a}$''$}\label{eq2.1a''} \end{equation}
Comme $H_2$ est supposé $\theto$\hyph invariant, le nombre $a_\alpha$ est constant sur l'orbite
de $\alpha$; l'inégalité~\eqref{eq2.1a''} est donc encore vraie pour tout $\alpha\in\Delta_\Q^\G-\Delta_\Q^{\Rm}$.
On en déduit l'inégalité~\eqref{eq2.1a'} et on conclut
comme dans le cas~\eqref{eq2.i}.
\end{proof}

\section{Quelques inégalités géométriques}

Dans toute cette section $\Q$ est un sous-groupe parabolique standard de $\G$.
On utilisera le symbole
\[
\ll
\]
\newindex{<@$\ll$}{infeerieur}%
qui signifie qu'il existe $c>0$ tel que le membre de gauche
soit inférieur à $c$ fois celui de droite.

\begin{lemme}\label{bigron}
Soient $P$ un sous-groupe parabolique $\theto$\hyph stable contenant $\Q$ et
\[
s=\so\rtimes\theto\in\weyl^{\tP}=\weyl^P\rtimes\theto\ptf
\]
On a deux possibilités\textup:
\begin{asparaenum}
\item il existe une constante $c>0$ telle que
\[
\langle X,X\rangle -\langle X,s X\rangle \ge c\langle X,X\rangle
\]
ou\textup, ce qui est équivalent\textup,
\[
\lVert (1-s)X\rVert \gg \lVert X\rVert
\]
pour tout $X\in\widebar C_0\cap\ga_\Q^P$
où $\widebar C_0$ est l'adhérence de la chambre de Weyl positive\textup,
\item il existe un sous-ensemble parabolique propre $\tP_1\varsubsetneq\tP$
de sous-ensemble de Levi $\tM_1$ avec $\Q\subset P_1\subset P$ et
\[
s=\so\rtimes\theto\in\weyl^{\tM_1}=\weyl^{\M_1}\rtimes\theto\ptf
\]
\end{asparaenum}
\end{lemme}

\begin{proof}
Tout d'abord on observe que $s$ étant une isométrie, on a
\[
2(\langle X,X\rangle-\langle X,s X\rangle)=\langle X-s X,X-s X\rangle=\lVert (1-s)X\rVert^2\ge0\ptf
\]
Supposons désormais que $\Q\ne P$, sinon le lemme est trivial.
Sur le compact, intersection de $\widebar C_0\cap\ga_\Q^P$
et de la sphère de rayon~1 (intersection qui est non vide puisque $\Q\ne P$), il existe $Y$ où la fonction
\[
X\mapsto\langle X,X-s X\rangle
\]
atteint son minimum $c\ge0$.
Si $c>0$ le lemme est démontré.
Supposons $c=0$; ceci équivaut à $Y=s Y$.
On a alors
\begin{equation}
\langle Y,Y-s Y\rangle =\langle Y,Y-\theto Y\rangle +\langle Y,\theto Y-s Y\rangle =0\ptf
\label{eq2.1b}
\end{equation}
Mais, $Y$ est dans l'adhérence de la chambre de Weyl positive
$C_0$ et il en est de même de $\theto Y$ puisque $\theto$
préserve la chambre positive.
D'après le lemme~\ref{wT},
\[
\theto Y- \so(\theto Y)
\]
est combinaison à coefficients positifs ou nuls
de racines positives. On en déduit que
\begin{equation}
\langle Y,\theto Y- s Y\rangle \ge0\ptf\label{eq2.2b}
\end{equation}
Enfin on a
\begin{equation}
2\langle Y,Y-\theto Y\rangle =\langle Y-\theto Y,Y-\theto Y\rangle\ge 0\label{eq2.3b}
\end{equation}
La conjonction de~\eqref{eq2.1b}, \eqref{eq2.2b} et~\eqref{eq2.3b} implique
\[
\langle Y-\theto Y,Y-\theto Y\rangle =0
\]
et donc
\[
Y=\theto Y\Qquad{et}Y=\so(Y)\ptf
\]
Le sous-groupe parabolique standard $P_1$ dont le sous-groupe de Levi $\M_1$
a pour racines simples les
\[
\alpha\in\Delta_{\PO}^P\qquad\text{telles que $\alpha(Y)=0$}
\]
est donc un sous-groupe parabolique qui est $\theto$\hyph stable puisque $Y=\theto Y$
et qui contient $\Q$ puisque $Y\in\ga_\Q^P$.
Enfin on a
\[
\so\in\weyl^{\M_1}
\]
puisque $Y=\so(Y)$. Il reste à observer que $P_1$ est strictement plus petit que
$P$ puisque $Y\ne0$.
\end{proof}

\begin{lemme}\label{rebigron}
Soient $s\in\weyl^{\tG}=\weyl\rtimes\theto$ et $R$ un sous-groupe parabolique contenant $Q$.
Supposons qu'il existe un unique sous-espace parabolique $\tP$ vérifiant les deux conditions suivantes\textup:
\[
\Q\subset P\subset\R\Qquad{et} s\in\weyl^{\tP}\ptf
\]
Considérons $H\in\ga_\Q^\G$ avec $\tsQR(H)=1$ alors
\[
\lVert (1-s)H\rVert \gg \lVert H\rVert \ptf
\]
\end{lemme}

\begin{proof}
On a une décomposition orthogonale:
\[
H=H_0+H_1+H_2
\]
avec $H_0\in\ga_\Q^P$, $H_1\in\ga_{\tP}^{\tG}$ et $H_2\in\gb$ où
$\gb$ est l'orthogonal de $\ga_{\tP}^{\tG}$ dans $\ga_P^\G$. On observe que $(1-s)$
envoie $\ga_\Q^P$ dans $\ga_{\PO}^P$,
préserve le sous-espace $\ga_P^\G$ et agit comme $(1-\theto)$ sur ce sous-espace;
en particulier il
s'annule sur $\ga_{\tP}^{\tG}$; de plus
\begin{equation}
\lVert (1-s)H_2\rVert \gg \lVert H_2\rVert\label{eq2.1c}
\end{equation}
par injectivité de $(1-s)$ sur $\gb$. D'après l'équation~\eqref{eq2.ii} du lemme~\ref{croiss} on sait que si $\tsQR(H)=1$ alors
\begin{equation}
\lVert H_1\rVert\ll \lVert H_0+H_2\rVert\label{eq2.2c}
\end{equation}
et, sous la même hypothèse le lemme~\ref{bigron} montre que
\begin{equation}
\lVert (1-s)H_0\rVert \gg \lVert H_0\rVert \label{eq2.3c}
\end{equation}
et donc, compte tenu de~\eqref{eq2.1c}
\begin{equation}
\lVert (1-s)H\rVert =\lVert (1-s)(H_0+H_2)\rVert \gg \lVert H_0\rVert + \lVert H_2\rVert \label{eq2.4c}
\end{equation}
et au total on obtient que si $\tsQR(H)=1$ alors
\begin{equation}
\lVert (1-s)H\rVert \gg \lVert H\rVert \ptf\qedhere\label{eq2.5c}
\end{equation}
\end{proof}

\section{Une application omniprésente}\label{omni}

Ce paragraphe introduit une application $q$ qui sera présente fréquemment dans la quatrième partie.
Les lemmes ci-dessous sont des variantes des lemmes~\ref{bigron} et~\ref{rebigron}
ci-dessus.

Soit $\Q$ un sous-groupe parabolique standard.
Pour $X\in\gao$ on note $X_\Q$ sa projection sur $\ga_\Q$.
On considère l'application linéaire
\[
q\colon \ga_0^G\to\ga_Q^G
\]
définie par
\[
X\mapsto\bigl((1-\theto) X\bigr)_Q
\]
et on note $\gk$ son noyau. Les sous-groupes
\[
\Qdo=\Q\cap\theta_0\moins\Q\Qquad{et}\theto(\Qdo)=\theto(Q)\cap\Q
\]
sont aussi des sous-groupes paraboliques standard.
Puisque
\[
\ga_0^{\Qdo}\subset \ga_0^Q\Quad{et} \theto(\ga_0^{\Qdo})\subset \ga_0^Q
\]
on a
\[
\ga_0^{\Qdo}\subset \gk\ptf
\]
D'où l'égalité $q(X_{\Qdo})=q(X)$ pour tout $X$.

\begin{lemme} \label{recroissa}
On a une majoration
\[
\lVert X\rVert \ll \lVert q(X)\rVert
\]
pour tout $X\in \ga_{\Qdo}^{\Qp}$ vérifiant $\tau_\Q^{\Qp}(X)\phi_{\Qdo}^Q(X)=1$.
\textup(La fonction $\phi_{\Qdo}^\Q$ a été introduite dans le lemme~\ref{partition}\textup).
\end{lemme}

\begin{proof}
Par hypothèse l'élément $X$ appartient au cône $\mathcal{C}$ engendré par
les $\varpi^\vee$, pour $\varpi\in \hDelta_Q^{\Qp}$, et les
$-\widebar\alpha^\vee$ pour pour $\alpha\in\Delta_{\Qdo}^\Q$
où $\widebar\alpha$ est la projection de $\alpha\in\Delta_0^\Q-\Delta_0^{\Qdo}$ sur le dual de $\ga_{\Qdo}^\Q$.
Il suffit de prouver que le cône engendré par les
images de ces éléments par l'application $q$ est un \og vrai\fg
cône, c'est-à-dire ne contient pas de sous-espace non nul.
Il suffit encore de prouver que $\widebar{\mathcal{C}}\cap \gk=\{0\}$. Soit
donc $X$ dans cette intersection. On pose
$Y=X^Q$ et $Z=X_Q$. On a
\[
0=q(X)=Z-(\theto Z)_Q-(\theto Y)_Q\ptf
\]
Par produit scalaire avec $Z$, on obtient
\begin{equation}
\langle Z,Z-(\theto Z)_Q\rangle=\langle Z,(\theto Y)_Q\rangle\ptf\tag{$*$}\label{eq2.*a}
\end{equation}
L'élément $(\theto Y)_Q$ appartient au cône engendré par les $- \bigl(\theto(\widebar\alpha^\vee)\bigr)_Q$ pour $\alpha\in
\Delta_0^Q-\Delta_0^{\Qdo}$.
Puisque $\theto(\Qdo)\subset Q$, on a
\[
(\theto\widebar\alpha^\vee)_{\Q}=(\theto\alpha^\vee)_{\Q}\ptf
\]
D'autre part, $\theto$ envoie injectivement $\Delta_0^Q-\Delta_0^{\Qdo}$ sur un sous-ensemble de
$\Delta_0^{\Qp}-\Delta_0^Q$: en effet si $\alpha$ et $\theto\alpha$ appartiennent à $\Delta_0^\Q$
alors $\alpha$ appartient à $\Delta_0^{\Qdo}$.
Donc, d'une part $(\theto Y)_Q$ appartient au cône engendré par les $- \beta^\vee$,
pour $\beta\in \Delta_\Q^{\Qp}$, d'autre part cet
élément n'est nul que si $Y=0$. L'élément $Z$ appartient au cône engendré par les $\varpi^\vee$, pour $\varpi\in \hat
{\Delta}_Q^{\Qp}$. Il en résulte que le membre de droite de~\eqref{eq2.*a}
est négatif ou nul. Or celui de gauche est positif ou nul: en effet, comme $Z=Z_\Q$ on a
\[
\langle Z,Z-(\theto Z)_Q\rangle =\langle Z,Z-\theto Z\rangle =\frac{1}{2}\lVert (1-\theto)Z\rVert^2\ptf
\]
Les deux membres de l'équation~\eqref{eq2.*a} sont donc nuls. La nullité
de celui de gauche entraîne que
\[
Z=(\theto Z )_Q=\theto Z\ptf
\]
Le lemme~\ref{qplusa} entraîne alors $Z=0$. Donc $ 0=q(X)=-(\theto Y)_Q$. On a déjà dit que cela impliquait $Y=0$.
\end{proof}

\begin{corollaire}\label{recroissb}
Soient $\Q\subset\R$ deux sous-groupes paraboliques standard.
On suppose que $\Qp=\Rm$. On a une majoration
\[
\lVert X\rVert \ll \lVert q(X)\rVert
\]
pour tout $X\in \ga_{\Qdo}^\G$ vérifiant $\tsQR(X)\phi_{\Qdo}^Q(X)=1$.
\end{corollaire}

\begin{proof}
On note $\gb$ le supplémentaire orthogonal de $\ga_{\tRm}^{\tG}$ dans $\ga_{\Rm}^\G$. On décompose $X$ en une somme de vecteurs orthogonaux:
\[
X=X_0+X_1+X_2\qquad\text{avec $X_0\in\ga_{\Qdo}^{\Rm}$, $X_1\in\gb$ et $X_2\in\ga_{\tRm}^{\tG}$}\ptf
\]
Tout d'abord l'équation~\eqref{eq2.ii} du lemme~\ref{croiss} montre que
\[
\lVert X_2\rVert \ll \lVert X_0+ X_1\rVert\ptf
\]
Comme $\Qp=\Rm$ on a $X_0\in\ga_{\Qdo}^{\Qp}$.
La condition $\tsQR(X)=1$ implique $\tau_\Q^\R(X)=1$ et en particulier
$\tau_\Q^{\Qp}(X_0)=1$. Alors, d'après le lemme~\ref{recroissa}, on a
\[
\lVert X_0\rVert \ll \lVert q(X_0)\rVert\ptf
\]
Comme $\ga_{\tRm}^{\tG}$ est le noyau de la restriction de $q$ à $\ga_{\Rm}^\G$,
l'application $q$ est injective sur $\gb$ et on a donc
\[
\lVert X_1\rVert \ll \lVert q(X_1)\rVert\ptf
\]
Compte tenu de ces trois inégalités et en observant que $q(X_0)$ et $q(X_1)$ sont orthogonaux on obtient
\[
\lVert X\rVert \ll \lVert q(X_0)\rVert +\lVert q(X_1)\rVert \ll \lVert q(X_0+X_1)\rVert =\lVert q(X)\rVert\ptf\qedhere
\]
\end{proof}

\begin{lemme}\label{ex10.3.2}
Soit $\T\in\gao$ tel que $T=\theto T$ et soit $P'$ avec $\Qdo\subset P'\subset\Q$.
On a une majoration
\[
\lVert (H-T_{\Qdo})\rVert \ll \lVert q(H)\rVert +\lVert (H-T)_{P'}^Q\rVert
\]
pour tous $T$\textup, $H$ tels que $\tsQR(H-T)\phi_{\Qdo}^{P'}(H-T)\tau_{P'}^Q(H-T)=1$.
\end{lemme}

\begin{proof}
On a $q(T_{\Qdo})=0$ en vertu de l'égalité $T=\theto T$.
En remplaçant $H$ par $H+T_{\Qdo}$, on
est ramené au cas $T=0$. On applique le lemme~\ref{recroissb} à
$H_Q$ et $H^{P'}$. On en déduit
\[
\Vert H_Q+H^{P'}\Vert <<\Vert q(H_Q+H^{P'})\Vert \ptf
\]
On a aussi
\[
\lVert q(H_Q+H^{P'})\rVert \ll \lVert q(H)\rVert +\lVert q(H^Q_{P'})\rVert <<\lVert q(H)\rVert +
\lVert H^Q_{P'}\rVert \ptf
\]
Enfin
\[
\lVert H\rVert \ll \lVert H_Q+H^{P'}\rVert +\lVert H^Q_{P'}\rVert\ptf
\]
Le lemme en résulte.
\end{proof}

\chapter{Théorie de la réduction}\label{ch3}

\section[Les fonctions $\HP$]{\mathversion{bold}Les fonctions $\HP$}

Rappelons que l'on a choisi un sous-groupe parabolique minimal $\PO$ de $\G$ sur $F$
et un sous-groupe de Levi $\MO$. Pour chaque place $v$
de $F$ on choisit un sous-groupe parabolique $\POO$ 
minimal sur $F_v$, de sous-groupe de Levi $\MOO$. 
On les choisit de sorte que $\MOO\subset\MO$ et $\POO\subset\PO$.
On note $\AOO$ la composante déployée d'un tore maximal
de $\MOO$ et on rappelle que
\[
\Norm_\G(\MOO)=\Norm_\G(\AOO)\ptf
\]
Lorsque $v$ est une place finie
on associe à $\AOO$ un appartement $\mathcal{A}$ de
l'immeuble $\mathfrak{B}$ de $\G_v$.
On dit que $\K_v$ est un sous-groupe compact spécial de $\G(F_v)$
s'il est le stabilisateur d'un point spécial $s\in\mathcal{A}$
(\cf \cite{T}*{1.9}).

\begin{lemme}\label{speci}
Soit $v$ une place finie.
Un sous-groupe spécial $\K_v$ vérifie les propriétés suivantes
\begin{enumerate}[(i)]
\item$\G(F_v)=\POO(F_v)\K_v$
\item Si $\M$ est un sous-groupe de Levi de $P$ contenant $\MOO$
alors $\K_v\cap\M(F_v)$ est un sous-groupe compact maximal spécial de $\M(F_v)$.
\item $\Norm_{\G(F_v)}(\M)\subset\M(F_v)\K_v$.
\end{enumerate}
\end{lemme}

\begin{proof}
Le point~(i) est la décomposition d'Iwasawa \cite{T}*{3.3.2}.
Le point~(ii) résulte de ce que $\mathcal{A}$ s'identifie
à un appartement de l'immeuble $\mathfrak{B}_\M$
de façon compatible à un plongement $\mathfrak{B}_\M\subset\mathfrak{B}$
et $s$ est \emph{a fortiori} spécial pour $\M$. Le groupe $\Norm_{\G(F_v)}(\MOO)$
opère sur $\mathcal{A}$ par transformations affines.
Comme $\MOO(F_v)$ opère par translation, le groupe
$\weyl=\Norm_{\G(F_v)}(\MOO)/\MOO(F_v)$
opère naturellement sur l'espace vectoriel $\Vect(\mathcal{A})$
associé à l'espace affine $\mathcal{A}$. On rappelle que d'après \cite{T}*{1.9}
$\K_v\cap \Norm_{\G(F_v)}(\MOO)$ s'envoie surjectivement sur le groupe de Weyl $\weyl$.
Ceci établit~(iii) pour $\MOO$. Le cas général résulte de ce que étant donné
$n\in \Norm_{\G(F_v)}(\M)$ il existe $m\in\M(F_v)$ avec $m\,n\in \Norm_{\G(F_v)}(\MOO)$.
\end{proof}

Nous dirons que $\K$ est un bon sous-groupe compact dans
$\Gadef $ s'il est de la forme
\[
\K=\prod_v\K_v
\]
où $\K_v$ est un sous-groupe compact maximal de $\G(F_v)$ pour toute
place $v$ et de plus, aux places finies, $\K_v$ est sous-groupe maximal spécial
de $\G(F_v)$; enfin, aux places réelles, on suppose que l'involution de Cartan relative
à $\K_v$ laisse stable $\MOO$.

\begin{lemme} \label{bon}
Soit $\KG$ un bon sous-groupe compact maximal.
\begin{enumerate}[(i)]
\item On a la décomposition d'Iwasawa\textup: $\Gadef =\PO(\adef)\,\KG$.
\item Si $\M$ est un sous-groupe de Levi de $P$ contenant $\MO$
alors
\[
\K_\M=\KG\cap\M(\adef)
\]
est un bon sous-groupe compact maximal de $\M(\adef)$.
\item On a
\[
\Norm_{\Gadef }(\M)\subset\M(\adef)\KG\ptf
\]
\end{enumerate}
\end{lemme}

\begin{proof}
C'est une conséquence facile du lemme~\ref{speci}.
\end{proof}

Nous choisirons désormais
un bon sous-groupe compact $\KG$ (le plus souvent noté simplement $\K$) de $\Gadef$.
On dispose de la fonction
\[
\HP\colon P(\adef)\to\ga_P
\]
que l'on prolonge, au moyen de la décomposition d'Iwasawa, en une fonction
\[
\HP\colon \Gadef \to\gao\ptf
\]
en posant
\[
\HP(pk)=\HP(p)
\]
et on observe que $\HP(\xix)=\HP(x)$ pour tout $\xi\in P(F)$.
La fonction $\HP$ dépend du choix de $\K$.

Lorsque $P=\PO$ nous noterons $\HO$ la fonction $\HPO$.
\newindex{H0@$\HO$}{hzero}%
Nous prolongerons $\HO$ en une fonction $\tHO$ sur $\tGadef$ en posant
\[
\tHO(x\delto)=\HO(x)\ptf
\]

\section{Hauteurs}

La construction que nous donnons ici est essentiellement celle proposée
dans la section~I.2.2 du livre \cite{MW} auquel nous renvoyons pour des preuves détaillées.

Soit $V$ un $F$-espace vectoriel de dimension finie muni d'une base $\{e_i\}_{1\le i\le n}$. On pose
\[
\lVert a\rVert_0=\prod_v \sup_i \lvert a_i \rvert_v\qquad\text{pour $a=\sum a_i e_i\in V\otimes\adef$}\ptf
\]
Cette fonction vérifie:
\begin{equation}
\lVert\lambda\,a\rVert_0=\lvert \lambda\rvert\ldot \lVert a\rVert_0\qquad\text{pour $\lambda\in\adef^\times$}\ptf
\label{eq3.1}
\end{equation}
On pourra remarquer que $\lVert a\rVert_0\ge \lvert a_i\rvert$ si $a_i\in\adef^\times$ et en particulier
\begin{equation}
\lVert \xi\rVert_0\ge1\qquad\text{pour $\xi\in V-\{0\}$}\ptf\label{eq3.2}
\end{equation}
On appellera hauteur sur $V$ une fonction
\[
\lVert \ \rVert\colon V\otimes\adef\to\RM
\]
vérifiant la propriété~\eqref{eq3.1} et équivalente à $\lVert\ \rVert_0$ \cad que
\[
c_1\lVert a\rVert_0\le \lVert a\rVert\le c_2\lVert a\rVert_0
\]
pour des constantes $c_i$ strictement positives.
En particulier, les hauteurs construites à partir de deux bases distinctes sont équivalentes.

De manière analogue, pour $g\in\End(V\otimes\adef)$ de matrice $g_{ij}$
on pose
\[
\lVert g\rVert_0=\prod_v \sup_{i,j} \lvert g_{ij} \rvert_v\ptf
\]
Il existe une constante $c_3$ telle que pour $g\in \End(V\otimes\adef)$ et $v\in V\otimes\adef$ on ait
\begin{equation}
\lVert gv\rVert_0\le c_3\lVert g\rVert_0\ldot \lVert v\rVert_0\label{eq3.3}
\end{equation}
et une constante $c_4$ telle que
\begin{equation}
\lVert g_1g_2\rVert_0\le c_4\lVert g_1\rVert_0\ldot \lVert g_2\rVert_0\ptf\label{eq3.4}
\end{equation}
Une hauteur est bornée sur les compacts. Donc, si $\K$ est un sous-groupe compact du groupe $\GL(V,\adef)$ on peut, par intégration sur $\K$, construire une hauteur bi-invariante sur $V$:
\[
\lVert kv\rVert=\lVert vk\rVert =\lVert v\rVert\qquad\forall k\in\K\ptf
\]
On définit une hauteur sur $\GL(V\otimes\adef)$ en
considérant une hauteur sur l'espace vectoriel
\[
\End(V)\oplus\End(V)
\]
et en posant pour $g\in\GL(V\otimes\adef)$:
\[
\lvert g\rvert=\lVert(g,{}^t g\moins)\rVert \ptf
\]
En particulier, on a\,\footnote{Dans \cite{MW} les auteurs utilisent $\SL(V)$ plutôt que $\GL(V)$
et n'utilisent pas la composition avec la diagonale $\lvert g\rvert =\lVert(g,g\moins)\rVert$. C'est la raison de leur inégalité~(iii) p.~20, qui pour nous est simplement $\lvert g\rvert =\lvert g\moins\rvert$.}
\[
\lvert g\rvert=\lvert g\moins\rvert\ptf
\]
Compte tenu de~\eqref{eq3.3} et~\eqref{eq3.4} on voit qu'il existe
une constante $c'_3$ telle que
\begin{equation}
\lVert g\,v\rVert \le c'_3\lvert g\rvert \ldot \lVert v\rVert\label{eq3.5}
\end{equation}
et une constante $c'_4$ telle que
\begin{equation}
\lvert g_1 g_2\rvert \le c'_4 \lvert g_1\rvert\ldot \lvert g_2\rvert \ptf\label{eq3.6}
\end{equation}
Enfin, il existe $c_0$ telle que
\begin{equation}
\lvert g\rvert\ge c_0\ptf \label{eq3.7}
\end{equation}

Supposons donnée une représentation linéaire fidèle de $\G$ dans $V$
\[
\rho\colon \G\to\GL(V)
\]
cela permet de définir une hauteur pour les éléments de $\Gadef $ en
posant
\[
\lvert x\rvert =\lvert \rho(x)\rvert =\bigl\lVert\bigl(\rho(x),{}^t\rho(x\moins)\bigr)\bigr\rVert \ptf
\]
Une telle hauteur vérifie encore~\eqref{eq3.5}, \eqref{eq3.6} et~\eqref{eq3.7}.

\begin{lemme}\label{proj}
Il existe des constantes $c_5$ et $N$ telles que
l'ensemble des $\xi\in\G(F)$ tels que $\lvert \xi\rvert \le A$
est un ensemble fini de cardinal majoré par $c_5\ldot A^N$.
\end{lemme}

\begin{proof}
L'application
\[
g\mapsto(g,{}^tg\moins)
\]
composée avec la projection sur l'espace projectif associé à l'espace vectoriel
\[
\End(V)\oplus\End(V)
\]
a des fibres de cardinal au plus~$2$. L'assertion résulte alors
des propriétés classiques des hauteurs sur un espace projectif.
\end{proof}

\begin{lemme}[\cf \cite{MW}*{assertion~(v), p.~20}]\label{isahauteur}
Il existe une constante $c_6$ telle que pour tout $x\in\Gadef $
\[
\lVert \HO(x)\rVert \le c_6\bigl(1+\bigl\lvert\log\lvert x\rvert\bigr\rvert\bigr)\ptf
\]
\end{lemme}

\begin{proof}
Il suffit de le prouver pour $x=p=mn\in\PO(\adef)$. Considérons
dans $\GL(V)$ un sous-groupe des matrices triangulaires supérieures
par blocs $P=\M\N$ où $M$ est diagonal par blocs et $N$ unipotent.
Alors $p=mn\in\M\N$ peut s'écrire
$p=m+X$
avec $X$ dans l'algèbre de Lie de $\N$, identifiée à un sous-espace
vectoriel de $\End(V)$, et
donc on a $\lVert m\rVert \le \lVert p\rVert$. Maintenant
soit $\rho$ la représentation rationnelle utilisée pour définir
la hauteur sur $G$. Quitte à changer de base et donc à utiliser une hauteur équivalente,
on peut supposer que
\[
\rho(\PO)\subset P\qquad\rho(\MO)\subset\M\qquad\rho(\NO)\subset\N
\]
où $P$ est, comme ci-dessus, un sous-groupe des matrices triangulaires supérieures
par blocs; alors pour $p=mn\in\PO(\adef)$
on a $\lVert \rho(p)\rVert \ge \lVert \rho(m)\rVert$. L'assertion en résulte facilement.
\end{proof}

\section[Calcul de $\HO(wn)$]{\mathversion{bold}Calcul de $\HO(wn)$}\label{hown}

\begin{lemme} \label{nilneg}
Soit $w\in\G(F)$ un représentant d'un élément $s$ du groupe de Weyl
de $\G$.
Il existe une constante $c$ telle que pour tout $n\in N_0(\adef)$
et tout $\vpi\in\hDelta_{\PO}$ on ait
\begin{equation}
\vpi\bigl(\HO(wn)\bigr)\le c\ptf\tag{i}\label{eq3.i}
\end{equation}
Supposons que $s=s_\alpha$ est une symétrie
par rapport à une racine simple $\alpha$ on a
\begin{equation}
\HO(w_\alpha\,n)=t_\alpha(n)\,\alpha^\vee\in\RM\alpha^\vee\tag{ii}\label{eq3.ii}
\end{equation}
avec $t_\alpha(n)\le c$.
Le nombre $t_\alpha(n)$ est indépendant du choix du représentant $w_\alpha\in\G(F)$.
\end{lemme}

\begin{proof}
Soit $\TO$ un tore maximal dans $\MO$. Choisissons un ordre sur les
racines de $\TO$ (sur la clôture algébrique)
compatible avec l'ordre sur les racines de $\AO$ (\cad les racines
relatives) déjà choisi.
Un caractère rationnel de $\MO$
\[
\lambda\in X_F(\MO)
\]
définit un caractère de $\TO$. Supposons le dominant.
Soit $V_\lambda$ la représentation rationnelle irréductible de $\G$ de
poids dominant $\lambda$ et de vecteur de plus haut poids $e_{\lambda}$.
On observe que $\lambda$ est encore un poids dominant
de la restriction à $\MO$ de la représentation $V_\lambda$.
Il en résulte que $\MO$
agit sur le vecteur $e_{\lambda}$ par le caractère $\lambda$.
On a donc, pour $m\in\MO(\adef)$ et $n\in\NO(\adef)$
\[
m\,e_{\lambda}=m^\lambda\,e_{\lambda}\Qquad{et} n\,e_{\lambda}=e_{\lambda}\ptf
\]
On choisit une hauteur $\K$-invariante sur cet espace,
normalisée de sorte que $\lVert e_{\lambda}\rVert=1$.
Alors si $x=n\,m\,k$ avec $k\in\K$ on a
\[
\lVert x\moins e_{\lambda}\rVert =\lVert k\moins\,m\moins\,n\moins e_{\lambda}\rVert =\lVert m\moins\,e_{\lambda}\rVert
=\lvert m^{-\lambda}\rvert \ldot \lVert e_\lambda\rVert\ptf
\]
Mais, comme
\[
\lvert m^{\lambda}\rvert=\exp\Bigl(\lambda\bigl(\HO(x)\bigr)\Bigr)
\]
on a
\[
\lambda\bigl(\HO(x)\bigr)=-\log\lVert x\moins e_{\lambda}\rVert
\]
et en particulier si on pose $w\moins\,e_\lambda= e_{s\moins(\lambda)}$ on aura
\[
-\log\lVert e_{s\moins(\lambda)}\rVert=-\log \lVert w\moins\,e_\lambda\rVert=\lambda\bigl(\HO(w)\bigr)
\]
et
\[
\lambda\bigl(\HO(wn)\bigr)=-\log\lVert n\moins e_{s\moins(\lambda)}\lVert\ptf
\]
Comme $n$ est dans le radical unipotent
\[
n\moins e_{s\moins(\lambda)}= e_{s\moins(\lambda)}+ v
\]
où $v$ est une somme de vecteurs de poids supérieurs à $s\moins(\lambda)$.
Donc il existe une constante $c_1$ telle que
\[
\lVert n\moins e_{s\moins(\lambda)}\rVert \ge c_1\lVert e_{s\moins(\lambda)}\rVert
\]
soit encore
\[
\lambda\bigl(\HO(wn)\bigr)\le c_2\ptf
\]
D'après Borel et Tits (\cf \cite{BTI}*{\parag 12 p.~141, commentaires suivant la proposition~12.13})
pour tout $\vpi\in\hDelta_{\PO}$ il existe un entier $d$ tel que $\lambda=d\vpi$
soit le poids dominant d'une représentation rationnelle. Ceci établit l'assertion~\eqref{eq3.i}.
Maintenant soit $s_\alpha$ une symétrie relativement à une racine simple
et soit $\beta\ne\alpha$ une autre racine simple. Si $\lambda=d_\beta\vpi_\beta$ on a
$s_\alpha(\lambda)=\lambda$ et donc, par unicité à un scalaire près du vecteur
de poids dominant on voit que $w_\alpha\moins e_{\lambda}$
est proportionnel à $e_{\lambda}$ ce qui implique
\[
n\moins w_\alpha\moins e_{\lambda}= c_\alpha \,n\,e_{\lambda}=
c_\alpha e_{\lambda}\qquad\text{avec $c_\alpha\in F^\times$}\ptf
\]
Il en résulte que
\[
\vpi_\beta\bigl(\HO(w_\alpha n)\bigr)=\vpi_\beta\bigl(\HO(w_\alpha)\bigr)=0
\]
d'où on déduit que
\[
\HO(w_\alpha n)\in\RM\alpha^\vee\ptf
\]
Enfin, l'inégalité $t_\alpha(n)\le c$ résulte de~\eqref{eq3.i}.
\end{proof}

\begin{lemme} \label{wn}
Soient $w_s\in\G(F)$ un représentant d'un élément $s$ du groupe de Weyl et
$n\in\N(\adef)$.
\begin{enumerate}[(i)]
\item Il existe des réels $h_\beta(s,n)\ge - c$\textup, où $c$ est la constante du lemme~\ref{nilneg}\textup, tels que
\begin{equation}
s\moins\HO(w_s\,n)=\sum_{\beta\in\Rac(s)}h_\beta(s,n)\,\beta^\vee\label{eq3.1a}
\end{equation}
soit encore
\begin{equation}
\HO(w_s\,n)=\sum_{\gamma\in\Rac(s\moins)}k_\gamma(s,n)\,\gamma^\vee\label{eq3.2a}
\end{equation}
avec des réels $k_\gamma(s,n)\le c$.
\item La fonction
\[
s\mapsto Y_s(n,T)=s\moins\bigl(T-\HO(w_s n)\bigr),\qquad s\in\weyl
\]
 est une famille orthogonale, régulière si $\dPO(T)>c$.
\item Plus généralement\textup, la fonction
\[
s\mapsto Y_s(x,\T)=s\moins\bigl(T-\HO(w_s x)\bigr)
\]
est une famille orthogonale
\newnot{Ys(x,T)@$Y_s(x,\T)$}{ysxt}%
et, si $x=mnk$ est une décomposition d'Iwasawa, on a
\begin{equation}
Y_s(x,\T)+\HO(x)+\sum_{\beta\in\Rac(s)}h_\beta(s,n,\T)\,\beta^\vee=0\tag{\ref{eq3.1a}$'$}\label{eq3.1a'}
\end{equation}
avec $h_\beta(s,n,\T)\ge0$ si $\dPO(\T)\ge c$.
\end{enumerate}
\end{lemme}

\begin{proof} 
Supposons $s=s_\alpha t$ avec $l(s)=l(t)+1$
et $s_\alpha$ la symétrie définie par rapport à la racine simple $\alpha$.
On écrit
\[
w_t n=m\,n'\, k
\] 
d'où
\[
\HO(w_t n)=\HO(m)
\] 
et $w_s\equiv w_\alpha\,w_t$ modulo $\M(F)$; donc
\[
\HO(w_s n)=\HO(w_\alpha m\,n'\,k)=s_\alpha\bigl(\HO(m)\bigr)+\HO(w_\alpha n')
\]
soit encore
\[
s\moins\HO(w_s\,n)=t\moins\HO(w_t\,n)+s\moins\HO(w_\alpha\,n')\ptf
\]
Il résulte alors de l'équation~\eqref{eq3.ii} du lemme~\ref{nilneg} que
\[
t\moins\HO(w_t\,n)-s\moins\HO(w_s\,n)=t_\alpha(n')\gamma^\vee
\]
avec $\gamma=t\moins\alpha^\vee$.
Les hypothèses du lemme~\ref{ws} sont donc vérifiées pour
\[
s\mapsto s\moins\HO(w_s\,n)\ptf
\]
Les formules~\eqref{eq3.1a} et~\eqref{eq3.2a} en résultent et
la construction par récurrence,
suivant~\ref{ws}, des coefficients
$h_\beta(s,n)$ et $k_\gamma(s,n)$ fournit les majorations souhaitées.
Compte tenu de~\ref{wT} l'assertion~(ii) en résulte immédiatement.
Pour établir~(iii) on observe que si
$x=mnk$ est une décomposition d'Iwasawa alors
\[
Y_s(x,\T)=Y_s(n,\T)-\HO(x)\ptf\qedhere
\]
\end{proof}

\begin{lemme}\label{Yorth}
Soit $w_s\in\G(F)$ un représentant d'un élément $s$ du groupe de Weyl.
\begin{enumerate}[(i)]
\item Le vecteur $\HO(w_s )\in\gao$ est indépendant du choix de $w_s$ et de $\PO$.
\item Il existe un point $\TK\in\gao^\G$ tel que
\[
\HO(w_s )=\TK-s\TK\Qquad{et} \HO(w_s\moins)=\TK-s\moins\TK\ptf
\]
\item
L'élément $\TK$ est égal à
\[
\TK=\sum_{\alpha\in\Delta_0}t_\alpha(1) \vpi_\alpha^\vee\ptf
\]
où $\vpi_\alpha^\vee\in\gao^\G$ est l'élément de la base duale
correspondant à $\alpha\in\Delta_0$ et $t_\alpha(1)$ le réel introduit à l'équation~\eqref{eq3.ii} du lemme~\ref{nilneg} pour $n=1$.
\end{enumerate}
\end{lemme}

\begin{proof} On rappelle que l'on a fixé $\K$ et  $\MO$.
D'après le lemme~\ref{bon}(iii), on peut écrire $w_s=m_sk_s$
avec $m_s\in\MO(\adef)$, bien défini modulo $\MO(F)$, et $k_s\in\K$
et donc
\[
\HO(w_s )=\HO(m_s)=\HH_{\MO}(m_s)
\]
est indépendant du choix de $w_s$ et de $\PO$. Comme $w_t=m_tk_t$, on voit que
\[
\HO(w_s w_t)=\HO(w_s m_t)=\HO(w_s)+s\HO(m_t)
\]
et donc
\[
\HO(w_s w_t)=\HO(w_s)+s\HO(w_t)\ptf
\]
Cette relation montre que
\[
s\mapsto \HO(w_s )
\]
est un 1-cocycle de $\weyl$ à valeurs dans $\gao$. Comme
la multiplication par l'entier $\lvert \weyl\rvert$ est inversible la cohomologie est nulle.
C'est donc un cobord et on obtient l'existence d'un $\TK\in\gao$ tel que
\[
\HO(w_s )=\TK-s\TK\ptf
\]
Pour achever la preuve de (ii) on remarque que comme
$w_s\moins$ et $w_{s\moins}$ sont congrus modulo $\M(F)$ on a
\[
\HO(w_s\moins)=\HO(\,w_{s\moins}) \Qquad{et donc} \HO(w_s\moins)=\TK-s\moins\TK\ptf
\]
Pour établir (iii) écrivons
\[
\TK=\sum_{\alpha\in\Delta_0} c_\alpha\vpi_\alpha^\vee\ptf
\]
On observe que
\[
s_\alpha\vpi_\beta=\vpi_\beta\quad\text{si $\beta\ne\alpha$}
\Qquad{et} s_\alpha\vpi_\alpha=\vpi_\alpha-\alpha^\vee
\]
 et donc
\[
 \HO(w_\alpha)=\TK-s_\alpha(\TK)=c_\alpha(\vpi_\alpha^\vee-s_\alpha\vpi_\alpha^\vee)
=c_\alpha\alpha^\vee\ptf
\]
Pour conclure on observe que d'autre part, d'après l'équation~\eqref{eq3.ii} du lemme~\ref{nilneg},
\[
 \HO(w_\alpha)= t_\alpha(1)\alpha^\vee\ptf
\]
\end{proof}

On posera
\[
Y_s(\T)=s\moins T+\TK-s\moins\TK\ptf
\]
\newnot{Ys(T)@$Y_s(T)$}{Yst}%
On a donc, avec les notations du lemme~\ref{wn} et compte tenu du lemme~\ref{Yorth}:
\[
Y_s(\T)=Y_s(1,\T)\ptf
\]
La fonction $s\mapsto Y_s(\T)$ définit une famille orthogonale.
On observera que
\[
Y_s(\TK)=\TK
\]
et est donc indépendant de $s$.
Soit $\M$ le sous-groupe de Levi d'un sous-groupe parabolique standard $P$;
à tout $S\in\Parab(\M)$
on associe, suivant les conventions de la section~\ref{famorth}, le vecteur
\[
Y_S(T)\in\ga_\M
\]
\newnot{YS(T)@$Y_S(T)$}{YSt}%
qui est la projection de $Y_s(\T)$ sur $\ga_\M$
lorsque $s\in\weyl$ est tel que $sS$ est standard.

\section[Espaces $\XP$, $\XPG$ et $\YP$]{\mathversion{bold}Espaces $\XP$, $\XPG$ et $\YP$}\label{espacesXP}

L'étude des formes automorphes et de la formule des traces amène à considérer divers espaces
homogènes pour lesquels nous allons fixer les notations. L'espace le plus important est
\[
\XG=\gA_G\G(F)\bs\Gadef \ptf
\]
\newindex{XG@$\XG$}{espacesXG}%
\newindex{XPG@$\XPG$}{espacesXPG}%
\newindex{YP@$\YP$}{espacesYP}%
C'est sur cet espace que vivent les formes automorphes pour $\G$, qui sont de plus
invariantes sous $\gA_\G$; ce sont les seules que nous considérerons ici.
Plus généralement si $P$ est un sous-groupe parabolique
de sous-groupe de Levi $\M$ et de radical unipotent $\N_P$
on aura besoin de considérer l'espace\,\footnote{On rappelle que, par définition, $\gA_P=\gA_{\M}$.}
\[
\XP=\gA_P P(F)\N_P(\adef)\bs\Gadef
\]
ainsi que sa variante
\[
\XPG=\gA_\G P(F)N_P(\adef)\bs\Gadef
\]
de sorte que
\[
\XP=\gA_P\bs\XPG\ptf
\]
On observera que l'on dispose d'une injection naturelle
\[
\XM\to\XP
\]
et d'une bijection
\[
\XM/\K_\M\simeq\XP/\KG\ptf
\]
Nous aurons besoin de considérer également les espaces
\[
\YP=\gA_\G P(F)\bs\Gadef \ptf
\]
On observera que $\YG=\XG$. Enfin on a une surjection
\[
\YP\to\XPG
\]
dont les fibres sont isomorphes à $\N(F)\bs\N(\adef)$.

\section{Ensembles de Siegel}\label{esphom}

Désormais, le sous-groupe compact
$\KG$ est simplement noté $\K$. Soit $t\in\RM$;
on appelle ensemble de Siegel un sous-ensemble de $\Gadef $ de la forme
\[
\gStOm=\Omega \AO(t)\K
\]
\newindex{StOmega@$\gStOm$}{espaces de Siegel}%
où $\Omega$ est un sous-ensemble compact de $\PO(\adef)$ et
\[
\AO(t)=\{\exp(H)\mid \text{$H\in\gao^\G$ et $\alpha(H)>t$, $\forall\alpha\in\Delta_{\PO}^\G$}\}\ptf
\]
On souligne que par construction on a $\HG(x)=0$ pour $x\in\gStOm$.

\begin{lemme}\label{siegcomp}
Il existe un compact $\Omega'\subset\Gadef $ tel que
tout $x\in\gStOm$ peut s'écrire la forme
\[
x=ac\qquad\text{avec $a\in\AO(t)$ et $c\in\Omega'$}\ptf
\]
\end{lemme}

\begin{proof}
On peut choisir $\Omega$ sous la forme $\Omega_1\times\Omega_2$
avec
\[
\Omega_1\subset\N_0(\adef)\Qquad{et} \Omega_2\subset\M_0(\adef)\ptf
\]
On a alors
\[
\gStOm=\Omega_1 \AO(t)\Omega_3\qquad\text{où $\Omega_3=\Omega_2 \K$}\ptf
\]
Comme $a\moins\omega a$ reste dans un compact lorsque $a\in\AO(t)$ et $\omega\in\Omega_2$
on a
\[
\gStOm\subset \AO(t)\Omega'
\]
où $\Omega'$ est un compact.
\end{proof}

Nous aurons besoin du théorème suivant pour lequel on renvoie au livre de Borel \cite{B}:

\begin{theoreme} \label{siegel}
Pour $t\in\RM$ donné
l'ensemble des $\gamma\in\G(F)$ tels que
\[
\gamma\ldot \gStOm\cap\gStOm\ne\varnothing
\]
est fini. De plus, pour $t$ assez petit et $\Omega$ assez gros\textup, on a
\[
\Gadef =\gA_\G\G(F)\ldot\gStOm\ptf
\]
En d'autres termes\textup, l'application naturelle
\[
\gStOm\to\XG
\]
est à fibres finies de cardinal borné. De plus\textup, pour $t$
assez petit et $\Omega$ assez gros\textup, l'application est surjective.
\end{theoreme}

Il en résulte que l'espace homogène
\[
\XG=\gA_G\G(F)\bs\Gadef
\]
est de volume fini, et que de plus il est compact si $\G_{\mathrm{der}}$, le groupe dérivé, est anisotrope.

\begin{proposition}\label{Psiegel}
Soit $\Q$ un sous-goupe parabolique.
Pour $t$ donné assez petit\textup, il existe des compacts
$\Omega_1\subset\N_\Q(\adef)$ et $\Omega_2\in\Gadef$
tels que tout $x\in\Gadef $ on ait
\[
x=\eta\,n\,a\,\omega\qquad\textup{avec $\eta\in\Q(F)$, $n\in\Omega_1$, $a\in\gA_0$, $\omega\in\Omega_2$}
\]
et
\[
\alpha\bigl(\HO(a)\bigr)>t\qquad\forall\alpha\in\Delta_{\PO}^\Q\ptf
\]
En particulier il existe une constante $c$ telle que pour $x\in\Gadef $
il existe $x_0\in\Gadef $ et $\eta\in\Q(F)$ avec $x=\eta\,x_0$ et
\[
\log\lvert x_0\rvert\le c\,(1+\lVert \HO(x_0)\rVert)\ptf
\]
\end{proposition}

\begin{proof}
Ceci résulte de la décomposition d'Iwasawa
\[
\Gadef =\Q(\adef)\K=\N_Q(\adef)\M_\Q(\adef)\K
\]
ainsi que du lemme~\ref{siegcomp} et du théorème~\ref{siegel} appliqués au
sous-groupe de Levi $\M_\Q$ de $\Q$.
\end{proof}

\begin{lemme}[\cf \cite{Fr}*{Theorem~1(1)}]\label{franke}
Il existe $c\in \RM$ tel que\textup, pour tout $x\in \Sieg^G$ et tout $\gamma\in \G(F)$\textup, on ait
\begin{equation}
\varpi_{\alpha}\bigl(\HO(\gamma x)-\HO(x)\bigr)\leq c \label{eq3.1b}
\end{equation}
pour tout $\alpha\in \Delta_0$.
\end{lemme}

\begin{proof}
Par décomposition de Bruhat et d'Iwasawa on a
\[
\HO(\gamma x)-\HO(x)=\HO(w_s n a)-\HO(a)=(s-1)\HO(a)+\HO(w n)
\]
avec $w$ représentant $s\in\weyl^\G$, $n\in\NO(\adef)$ et $a$ vérifiant
\[
\alpha\bigl(\HO(a)\bigr)>t\qquad\forall\alpha\in\Delta_{\PO}^\G\ptf
\]
Maintenant le lemme~\ref{wn} montre que
\[
\varpi_{\alpha}\bigl(\HO(wn)\bigr)\le c_1
\]
et le lemme~\ref{wT} montre que
\[
\varpi_{\alpha}\bigl((s-1)\HO(a)\bigr)\le c_2\ptf\qedhere
\]
\end{proof}

\begin{lemme}\label{xycomp}
Soient $\gamma\in\G(F)$\textup, $x\in \Sieg^G$ et $y\in \Sieg^G$. Si $x^{-1}\gamma y\in \Omega$ où $\Omega$ est un compact\textup, alors $\HO(x)-\HO(y)$ appartient à un compact.
\end{lemme}

\begin{proof}
Si $x^{-1}\gamma y\in \Omega$, on a $\gamma y\in x\Omega$ et donc
$\HO(\gamma y)- \HO(x)$ reste dans un compact.
Compte tenu du lemme~\ref{franke} on en
déduit une majoration
\[
\varpi_{\alpha}\bigl(\HO(x)-\HO(y)\bigr)\leq c_{2}
\]
pour tout $\alpha$.
La situation est symétrique en $x$ et $y$. Donc $\HO(x)-\HO(y)$ appartient à un
compact.
\end{proof}

\begin{lemme} [\cite{MW}*{1.2.2(vii)}] \label{ratmaj}
Il existe $c>0$ tel que\textup, pour tout $x\in\Sieg$ \textup(domaine de Siegel pour $G$\textup) et tout~$\xi\in G(F)$\textup:
\[
\lvert x\rvert\le c\lvert \xix\rvert\ptf
\]
\end{lemme}

\begin{proof}
Comme, d'après le lemme~\ref{siegcomp},
$x=a\omega$ avec $\omega$ dans un compact et $a\in\gA_0$,
il suffit de traiter le cas $x=a\in\gA_0$.
Le lemme résulte alors de la remarque suivante:
si $\xi$ est une matrice rationnelle dans $\GL(V)$
et $a$ une matrice diagonale dans $\GL(V\otimes\RM)$ on a
\[
\lVert \xi a\rVert_0=\prod_v\sup_{i,j}\lvert\xi_{ij}a_{jj}\rvert_v\ge
\sup_{i,j}\prod_v\lvert \xi_{ij}a_{jj}\rvert_v=\sup_{j}\lvert a_{jj}\rvert_{\RM}=
\lVert a\rVert_0\ptf\qedhere
\]
\end{proof}

\section[Une partition de $\XG$]{\mathversion{bold}Une partition de $\XG$}\label{unepartXG}

Pour établir la partition~\ref{FPQ} (qui est l'analogue pour $\XP$ de
la partition~\ref{partition} de l'espace vectoriel $\gao$)
nous aurons besoin du lemme suivant:

\begin{lemme}[\cf \cite{MS}*{Lecture 3, Erratum Lemma~3.2.3}]\label{unique}
Fixons $\Tsieg\in\gao$ et soit $P$ un sous-groupe parabolique standard.
Pour $T\in\gao$ assez régulier \textup(de façon précise si $\dPO(T)>C_1$
où $C_1$ est une constante dépendant de $\Tsieg$\textup)\textup,
l'ensemble des $\xi\in\G(F)$ tels que\textup, pour $x\in\Gadef$ donné on ait
\begin{equation}
 \alpha(\HO(\xix)-\Tsieg)>0\quad\ \alpha\in\Delta_{\PO}^P\Qquad{et}
\alpha(\xixT)>0\ \forall\alpha\in\Delta_{\PO}^\G-\Delta_{\PO}^P
\label{eq3.1c} 
\end{equation}
est soit vide soit forme une seule classe modulo $P(F)$.
\end{lemme}

\begin{proof} 
Pour $\T$ donné, il suffit de considérer le cas particulier où
\[
\alpha(\HO(x)-\Tsieg)>0\quad\forall\alpha\in\Delta_{\PO}^P\Qquad{et}
\alpha(\xT)>0\quad\forall\alpha\in\Delta_{\PO}^\G-\Delta_{\PO}^P\ptf
\]
Compte tenu de la décomposition de Bruhat
il reste à montrer que si $T$ est assez régulier et si $\xi=w_s$
représente un élément $s$ dans le groupe de Weyl de $\G(F)$, alors
les hypothèses sont satisfaites seulement si $s$ appartient au
groupe de Weyl du sous-groupe de Levi $\M$ de $P$.
On observe que si $x=mnk$ on a, d'après le lemme~\ref{wn},
\begin{equation}
s\moins\HO(w_s\,x)
=\HO(x)+\sum_{\beta\in\Rac(s)}h_\beta(s,n)\, \beta^\vee\label{eq3.2c}
\end{equation}
avec des scalaires $h_\beta(s,n)\ge- c$.
On notera $\demisomO$ la demi-somme des racines dans $\PO$.
Considérons
\[
\lambda=\demisomO-s\moins\demisomO=\sum_{\beta\in\Rac(s)} \beta\ptf
\]
On observe que pour tout $\beta\in\Rac(s)$
alors $\gamma=-s(\beta)$ est positive et donc
\begin{equation}
\lambda(\beta^\vee)=\demisomO(\beta^\vee+\gamma^\vee)>0\ptf \label{eq3.3c}
\end{equation}
De plus $s\lambda$ est l'opposé d'une somme de racines positives:
\begin{equation}
s\lambda=-\sum_{\gamma\in\Rac(s\moins)}\gamma\ptf \label{eq3.4c}
\end{equation}
Compte tenu de l'équation~\eqref{eq3.2c} et de l'inégalité~\eqref{eq3.3c}
il existe une constante $C$ telle que
\[
s\lambda(\HO(w_sx)-\Tsieg) +C\ge \lambda\bigl(\HO(x)\bigr)\ptf
\]
Supposons $T-\Tsieg$ régulier, ce qui est loisible; alors l'hypothèse~\eqref{eq3.1c} implique que
\[
\alpha(\HO(w_sx)-\Tsieg)>0\qquad\forall\alpha\in\Delta_{\PO}^\G
\]
et donc, d'après~\eqref{eq3.4c},
\[
s\lambda(\HO(w_sx)-\Tsieg)\le0
\]
ce qui implique
\begin{equation}
\lambda\bigl(\HO(x)\bigr) \le C\ptf\label{eq3.5c}
\end{equation}
On peut écrire
\[
\HO(x)=\sum_{\alpha\in\Delta_{\PO}^\G}d_\alpha\vpi_\alpha^\vee+\HG
\]
avec $\HG\in\ga_\G$ et~\eqref{eq3.5c} fournit l'inégalité
\begin{equation}
\sum_{\beta\in\Rac(s)} \,\sum_{\alpha\in\Delta_{\PO}^\G} c^\alpha_\beta\,d_\alpha\le C\label{eq3.6c} 
\end{equation}
où les $c^\alpha_\beta$ sont des entiers naturels définis par
\[
\beta=\sum_{\alpha\in\Delta_{\PO}^\G} c^\alpha_\beta \alpha\ptf
\]
Par hypothèse
\[
d_\alpha=\alpha\bigl(\HO(x)\bigr)>\alpha(\Tsieg)\qquad\forall\alpha\in\Delta_{\PO}^P
\]
et
\[
d_\alpha=\alpha\bigl(\HO(x)\bigr)>\alpha(T)\qquad\forall\alpha\in\Delta_{\PO}^\G-\Delta_{\PO}^P\ptf
\]
L'inégalité~\eqref{eq3.6c} n'est possible, si $T$ est assez régulier, que si les $\beta\in\Rac(s)$
ne font intervenir que les $\alpha\in\Delta_{\PO}^P$ auquel cas
$s$ appartient au groupe de Weyl de $\M$, le sous-groupe de Levi de $P$,
ainsi qu'il résulte du lemme~\ref{swp}.
\end{proof}

On introduit, pour $\PO\subset\Q$, l'ensemble
$\gSP^\Q(\Tsieg,\T)$
des
\[
x=nac\in\Gadef
\] a
vec $ n\in\N_\Q(\adef)$, $a=\ee^H$ pour
$H\in\gao$, $c\in C_\Q$ où $C_\Q$ est un compact (qui sera choisi assez gros
dans la proposition~\ref{FPQ}),
et vérifiant
\[
 \alpha(\HO(x)-\Tsieg)>0\quad\forall\alpha\in\Delta_{\PO}^\Q
\Quad{et}
\vpi(\xT)\le0\quad\forall\vpi\in\hDelta_{\PO}^\Q
\]
soit encore, d'après le lemme~\ref{Gammacar}, tel que
\[
\Gamma_{\PO}^\Q(\HO(x)-\Tsieg,T-\Tsieg)=1\ptf
\]
On note $\FPQ(\bullet,\T)$
\newindex{FP0Q@$\FPQ$}{deffpq}%
la fonction caractéristique de l'ensemble
\[
\Q(F)\gSP^\Q(\Tsieg,\T)\ptf
\]

\begin{lemme}\label{gSPQ}
Soit $x\in\Gadef$ tel que $\FPQ(x,T)=1$. Il existe $\gamma\in\Q(F)$ tel que
\[
\gamma\,x=nac
\]
avec $c$ dans un compact\textup,
$n\in\N_\Q(\adef)$ et
la projection de $\HO(a)$ dans $\gao^\Q$ est bornée.
\end{lemme}

\begin{proof}
Par définition de $\FPQ$ on peut choisir $\gamma\in\Q(F)$ tel que
$\gamma\,x$ appartienne à $\gSP^\Q(\Tsieg,\T)$.
Le lemme résulte alors de ce que, d'après le lemme~\ref{convcomp},
la projection de
\[
\HO(\gamma x)=\HO(a)+\HO(c)
\]
sur $\gao^\Q$ est bornée.
\end{proof}

\begin{proposition}[\cite{MS}*{Proposition~3.2.1}] \label{FPQ}
Si le compact $C_\Q$ implicite dans la définition
de $\gSP^\Q(\Tsieg,\T)$ est assez gros alors\textup,
pour $-\Tsieg$ et $T$ assez réguliers\textup,
\cad $\dPO(\T)\ge c$ et $\dPO(-\Tsieg)\ge c'$
où $c$ et $c'$ sont des contantes dépendant de $\G$\textup,
on a
\[
\sum_{\{\Q\mid\PO\subset\Q\subset P\}}
\sum_{\xi\in\Q(F)\bs P(F)}\FPQ(\xix,T)\tau_\Q^P(\xixT)=1\ptf
\]
\end{proposition}

\begin{proof}
D'après la proposition~\ref{Psiegel}, si $-\Tsieg$ est assez régulier, il existe un compact $C$
tel que pour tout $x\in\Gadef $
il existe
\[
\xi\in P(F)
\]
tel que:
\[
\xix=nac
\] 
avec $ n\in\NO(\adef)$, $a=\ee^H$ pour
$H\in\gao$ et $c\in C$, vérifiant
\[
\tau_{\PO}^P(\HO(\xix)-\Tsieg)=1\ptf
\]
Maintenant l'équation~\eqref{eq2.1} du lemme~\ref{GammaHX} montre que
\[
\tau_{\PO}^P(\HO(\xix)-\Tsieg)=\sum_{\PO\subset\Q\subset P}
\Gamma_{\PO}^\Q(\HO(\xix)-\Tsieg,T-\Tsieg)\tau_\Q^P(\xixT)\ptf
\]
On observe que si, pour un certain $\Q$, notre $\xix=nac$ vérifie de plus
\[
\Gamma_{\PO}^\Q(\HO(\xix)-\Tsieg,T-\Tsieg)\tau_\Q^P(\xixT)=1
\]
alors, si $C_\Q\supset C$ on a
\[
\xix\in\gSP^\Q(\Tsieg,\T)
\]
et donc
\[
\FPQ(\xix,T)\tau_\Q^P(\xixT)=1
\]
Il reste à observer que d'après le lemme~\ref{recone} les hypothèses du lemme~\ref{unique} sont satisfaites et donc un tel $\xi$
est uniquement déterminé
modulo $\Q(F)$ si $T$ est assez régulier.
\end{proof}

\begin{lemme}\label{FPQR}
Pour $P$ fixé et sous les hypothèses de~\ref{FPQ}\textup, on a
\[
\sum_{\Q\subset P\subset \R} \, \sum_{\xi\in \Q(F)\bs P(F)}
\FPQ(\xix,T) \sQR(\xixT)=\htau_P(\xT)\ptf
\]
\end{lemme}

\begin{proof}
Ceci résulte de la décomposition
\[
 \tau_\Q^P\htau_P=\sum_{\{\R\mid P\subset\R\}}\sQR
\]
(\cf lemme~\ref{repart}) et de la proposition~\ref{FPQ}.
\end{proof}

\begin{lemme}\label{newF}
Supposons $T$ assez régulier \textup(comme à la proposition~\ref{FPQ} ci-dessus\textup). Soit $x\in\gSP^\Q(\Tsieg,\T)$ avec
\[
{\tsQR}(\xT)=1\ptf
\] Alors
\[
\alpha(\HO(x)-\Tsieg)>0\qquad\forall\alpha\in\Delta_{\PO}^\R\ptf
\]
\end{lemme}

\begin{proof}
Par hypothèse
\[
\alpha(\xT)>0\qquad\forall\alpha\in\Delta_\Q^\R\ptf
\]
Comme on suppose d'autre part $x\in\gSP^\Q(\Tsieg,\T)$ on a
\[
\vpi(\xT)\le0\qquad\forall\vpi\in\hDelta_{\PO}^\Q
\]
il résulte du lemme~\ref{recone} que
\begin{equation}
\alpha(\xT)>0\qquad\forall\alpha\in\Delta_{\PO}^\R-\Delta_{\PO}^\Q\tag{i}\label{eq3.ia}
\end{equation}
mais comme $\alpha(T-\Tsieg)>0$ pour tout $\alpha\in\Delta_{\PO}^\G$
on a aussi
\begin{equation}
\alpha(\HO(x)-\Tsieg)>0\qquad\forall\alpha\in\Delta_{\PO}^\R\ptf\qedhere\tag{ii}\label{eq3.iia} 
\end{equation}
\end{proof}

\begin{lemme}\label{aloin}
Soient $\tP$ un sous-ensemble parabolique,
$\Q$ et $\R$ deux sous-groupes paraboliques standard tels que
$\Q\subset P\subset\R$. Considérons
$\delta\in\tP(F)$\textup, $n, n'\in\NO(\adef)$ et $a\in \gAO$.
Soit $\Omega$ un sous-ensemble compact de $\tGadef$.
Supposons que
\begin{equation}
a\moins\,n\,\delta\,n'\,a\in\Omega\tag{i}\label{eq3.ib}
\end{equation}
et que $a$ satisfasse aux inégalités
\begin{equation}
\alpha(\HO(a)-T)>0\qquad\forall\alpha\in\Delta_{\PO}^\R-\Delta_{\PO}^\Q\tag{ii}\label{eq3.iib}
\end{equation}
et
\begin{equation}
\alpha(\HO(a)-\Tsieg)>0\qquad\forall\alpha\in\Delta_{\PO}^\R\ptf\tag{iii}\label{eq3.iiib}
\end{equation}
Alors\textup, il existe une constante $c(\Omega)$ telle que
si $\dPO(T)>c(\Omega)$ alors\textup, avec les notations du lemme~\ref{qplus}\textup, on a
\[
\delta\in\tQp(F)\ptf
\]
\end{lemme}

\begin{proof}
La décomposition de Bruhat permet d'écrire
\[
\delta=\nu\eta\,w_s\,\nu'
\]
avec $\nu,\nu'\in N_0(F)$, $\eta\in \MO(F)$ et où $w_s$ représente un élément
\[
s=\so\rtimes\theto
\]
de l'ensemble de Weyl de $\tM$. On a donc
\[
\tHO(a\moins\,n\delta\,n' a)=\tHO(a\moins \,w_{s}\, a\,n'')
=\HO(a\moins)+s\HO(a)+\tHO( \,w_{s} n'')\ptf
\]
Posons
\[
A_s=s\moins\bigl(\HO(a)-\tHO(w_{s}n'')\bigr)\ptf
\]
L'hypothèse~\eqref{eq3.ib} implique que
\[
A_1-A_s=\HO(a)-s\moins\HO(a)+s\moins\tHO(w_{s} n'')
\]
appartient au compact $s\moins\tHO(\Omega)$.
Fixons $T_1\in\gao$ vérifiant $\dPO(T_1)\ge c$, où
$c$ est la constante du lemme~\ref{nilneg}.
On introduit
\[
B_s=s\moins\bigl(T_1-\tHO(w_{s} n'')\bigr)\Qquad{et}C_s=s\moins(\HO(a)-T_1)\ptf
\]
On observe que $A_s=B_s+C_s$.
Soit
\[
\lambda=\sum_{\vpi\in{\hDelta}_{\PO}^P}\vpi\ptf
\]
C'est une forme linéaire sur $\gao^P$ qui est $\theto$-invariante
et strictement positive sur toutes les racines positives pour $\M$.
Compte tenu de la $\theto$-invariance on a
\[
\lambda(s\moins X)=\lambda(\so\moins X)\ptf
\]
On sait par le lemme~\ref{wn} que $\lambda(B_1-B_s)$ est positive. La condition de compacité sur $A_1-A_s$
impose que
$\lambda(C_1-C_s)$ est borné supérieurement par une constante $c_0(\Omega)$.
Si $\tP\ne\tQp$ et si $w_{s}\ne\tQp(F)$ la décomposition réduite de $\so$ fait intervenir
une racine simple $\alpha\in\Delta_{\PO}^P-\Delta_{\PO}^{\Qp}$.
Mais le lemme~\ref{wT} montre que $\lambda(C_1-C_s)$ est la somme de
\[
\lambda(\beta^\vee)\alpha(\HO(a)-T_1)
\]
où $\beta$ est une racine positive et
d'autres termes qui sont minorés d'après~\eqref{eq3.iiib}.
On en déduit que, d'après l'hypothèse~\eqref{eq3.iib} et pour une certaine constante $c_1(\Omega)$
\[
\alpha(T-T_1)<\alpha(\HO(a)-T_1)\le c_1(\Omega)
\]
ce qui est impossible
si par ailleurs $\dPO(T)>c(\Omega)$,
pour une constante $c(\Omega)$ bien choisie.
\end{proof}

\begin{corollaire}[\cf \cite{MS}*{Lemma~4.1.3}]\label{deltaQ}
Soit $\Omega$ un compact de $\tGadef$.
Supposons que
\[
\FPQ(x,T){\tsQR}(\xT)\ne0
\]
et
\[
x\moins\, n\,\delta\,n'\,x\in\Omega
\]
avec $\delta\in\tP(F)$\textup, $n\in\NO(\adef)$ et $x\in\Gadef $. Si $\dPO(T)>c(\Omega)$ ceci
implique $\delta\in\tQp(F)$.
\end{corollaire}

\begin{proof}
Quitte à changer $x$ en $\xix$ avec $\xi\in\Q(F)$ on peut supposer
vérifiées les inégalités~\eqref{eq3.ia} et~\eqref{eq3.iia}
de la preuve du lemme~\ref{newF}. On écrit
$x=n_1\,m\,a\,k$ avec $m\in\MO(\adef)$ et $\HO(m)=0$.
Quitte à changer $n$ et $n'$ on peut supposer $n_1=1$.
Maintenant, modulo conjugaison de $\delta$ par $\gamma\in\PO(F)$
on peut supposer que $m$ appartient à un ensemble compact.
Donc $a$ satisfait les conditions~\eqref{eq3.ib}, \eqref{eq3.iib} et~\eqref{eq3.iiib} de~\ref{aloin}.
\end{proof}

\section{Lemmes de finitude}\label{sec3.7}

\begin{lemme} \label{finitude}
Pour $x$ fixé\textup, il existe des constantes $C$\textup, $N$ et $A$
telles que
l'ensemble des $\xi\in\G(F)$ vérifiant\textup, pour $X\in\gao$\textup,
\[
\htau_P^\G(\HO(\xix)-X )\ne0
\]
est un ensemble fini de classes modulo $P(F)$
dont les représentants peuvent être choisis de sorte que
\[
\lvert \xi\rvert\le
C \lvert x\rvert^{N+1}\ee^{A \lVert X\rVert}\ptf
\]
En particulier cet ensemble peut être choisi indépendant de $x$\textup,
lorsque $x$ reste dans un compact.
\end{lemme}

\begin{proof}
Si cet ensemble est non vide on peut, quitte à changer $\xi$ et $x$, supposer que $x$ vérifie
\[
\htau_P^\G(\HO(x)-X)\ne0
\]
\cad
\[
 \vpi(\HO(x)-X)>0\qquad\forall\vpi\in\hDelta_P^\G\ptf
\]
De plus, compte tenu de la proposition~\ref{Psiegel}, on peut supposer que
\[
 \alpha(\HO(x)-\Tsieg)>0\qquad\forall\alpha\in\Delta_{\PO}^P\ptf
\]
Ceci implique (d'après le lemme~\ref{cone})
\[
\vpi(\HO(x)-\Tsieg)>0\qquad\forall\vpi\in\hDelta_{\PO}^\G\ptf
\]
Maintenant on peut écrire $\xi=n_1w\,\eta\,n_2$
et donc, si $s$ est l'image de $w$ dans le groupe de Weyl, on a
\[
\HO(\xix)=\HO(w\,n_2\,x)=s\HO(x)+\HO(w\,n)
\]
pour un certain $n$ (dépendant de $x$)
et, d'après le lemme~\ref{wn},
\[
\HO(w_s\,n)=\sum_{\gamma\in\Rac(s\moins)}k_\gamma(s,n)\,\gamma^\vee
\]
avec des réels $k_\gamma(s,n)\le c$.
Donc il existe une constante $c'$ telle que
\begin{equation}
 \vpi\bigl(\HO(\xix)\bigr)\le\vpi\bigl(s\HO(x)\bigr)+c'\qquad\forall\vpi\in\hDelta_{\PO}^\G\ptf
\label{eq3.1d} 
\end{equation}
D'après la proposition~\ref{Psiegel} on peut choisir $\xi$,
modulo $P(F)$ à gauche, de sorte que
\begin{equation}
\alpha(\Tsieg)\le\alpha\bigl(\HO(\xix)\bigr)\qquad\forall\alpha\in\Delta_{\PO}^P\ptf\label{eq3.2d} 
\end{equation}
Par ailleurs, si
\[
\htau_P^\G(\HO(\xix)-X)\ne0
\]
on a
\begin{equation}
\vpi(X)\le\vpi\bigl(\HO(\xix)\bigr)\qquad\forall \vpi\in\hDelta_P^\G\ptf\label{eq3.3d} 
\end{equation}
En combinant~\eqref{eq3.1d}, \eqref{eq3.2d} et~\eqref{eq3.3d} ainsi que le lemme~\ref{isahauteur} on obtient que
\[
\lVert\HO(\xix)\rVert\le c''(1+\log\lvert x\rvert+\lVert X\rVert)\ptf
\]
En invoquant la deuxième assertion de la proposition~\ref{Psiegel}
on voit qu'avec notre choix de $\xix$ on a
\[
\log\lvert \xix\rvert\le c'''(1+\log\lvert x\rvert+\lVert X\rVert)\ptf
\] et donc
il existe des constantes $C_1$, $N$ et $A$ telles que
\[
\lvert \xix\rvert\le
C_1\lvert x\rvert^N\ee^{A \lVert X\rVert}
\]
et donc
\[
\lvert \xi\rvert\le C_2\lvert\xix\rvert\ldot \lvert x\moins\rvert\le
C\lvert x\rvert^{N+1}\ee^{A \lVert X\rVert}
\]
ce qui, d'après le lemme~\ref{proj}, impose à $\xi$ d'être dans un ensemble fini.
\end{proof}

\begin{lemme} \label{finitudebis}
Soit $\tP$ un sous-ensemble parabolique.
Pour $x$ fixé\textup, il existe des constantes $C'$\textup, $N'$ et $A'$
telles que
l'ensemble des $\xi\in\G(F)$ vérifiant
\[
\htau_{\tP}^{\tG}(\xixT )\ne0
\]
est un ensemble fini de classes modulo $P(F)$
dont les représentants peuvent être choisis de sorte que
\[
\lvert \xi\rvert \le
C' \lvert X\rvert^{N'+1}\ee^{A' \lVert \T\rVert}\ptf
\]
En particulier cet ensemble peut être choisi indépendant de $x$\textup,
lorsque $x$ reste dans un compact.
\end{lemme}

\begin{proof}
La preuve est identique à celle du lemme~\ref{finitude} à ceci près qu'au lieu de
\eqref{eq3.3d} on a
\begin{equation}
\tvpi(T)\le\tvpi\bigl(\HO(\xix)\bigr)\qquad\forall
\tvpi\in\hDelta_{\tP}^{\tG}\ptf\tag{\ref{eq3.3d}$'$}\label{eq3.3d'} 
\end{equation}
La conclusion est identique (avec des constantes différentes).
\end{proof}

On dira que $\delta\in\tG(F)$ est primitif si sa classe de conjugaison
ne rencontre aucun $\tP(F)$ lorsque $\tP$ parcourt l'ensemble des
sous-ensembles paraboliques propres, \cad $\tP\ne\tG$.
On notera $\tG(F)_{\mathrm{prim}}$ l'ensemble des éléments primitifs.

\begin{lemme}\label{geofini}
Soit $\Omega$ un compact de $\tGadef$ et $\Sieg$ un ensemble de Siegel.
L'ensemble des $\delta\in\tG(F)_{\mathrm{prim}}$ tels qu'il existe
$x\in\Sieg$ avec
\[
x\moins\, \delta\,x\in\Omega
\]
est fini.
\end{lemme}

\begin{proof}
La décomposition de Bruhat permet d'écrire
\[
\delta=\eta\,w_s\,\xi\,\eta'
\]
avec $\xi\in\MO(F)$ et $w_s$ représente un élément de l'ensemble de Weyl
$\weyl\rtimes\theto$.
Grâce au lemme~\ref{siegcomp}
et au théorème~\ref{siegel} on a
\[
a\moins \delta\, a = a\moins \eta\,w_s\,\xi\,\eta'\, a \in\Omega'
\]
pour un $a\in \AO(t)$ et pour un compact $\Omega'$, soit encore
\[
n\, a'\,w_s\, \xi\,n'\in\Omega'
\]
avec $n\in\NO(\adef)$, $n'=a\moins\eta' a\in\NO(\adef)$
et $a'=a\moins w_s aw_s\moins$ et donc
\begin{equation}
\HO(n\, a'\,w_s\, \xi\,n')=\HO(a')+\HO(w_s n')\tag{$*$}\label{eq3.*}
\end{equation}
appartient à un compact.
On observe que
\[
\HO(a')=(s-1)\HO(a)\ptf
\]
Puisque $a$ est dans un domaine de Siegel,
$X=\HO(a)+S$ est dans la chambre positive pour un certain $S\in\gao$.
C'est dire que
\[
X=\sum a_\alpha\vpi^\vee_\alpha
\]
avec $a_\alpha>0$.
Maintenant~\eqref{eq3.*} montre que
\[
\langle X,(s-1)X\rangle +\langle X,\HO(w_sn')\rangle
\] est borné.
D'après le lemme~\ref{wn},
on a
\[
\langle\tvpi^\vee,\HO(w_sn')\rangle\le c
\]
pour tout $n'$ et tout $\tvpi$.
On a donc
\[
C_1\le\langle X,(s-1)X\rangle+\langle X,\HO(w_sn')\rangle \le\langle X,(s-1)X\rangle +C_2\ptf
\]
Il existe donc une constante $C$ telle que
\[
\langle X,(1-s)X\rangle\le C
\]
Mais, si $\langle X,(1-s)X\rangle$ reste borné
alors que $\lVert X\rVert$ tend vers l'infini, il résulte du lemme~\ref{bigron} qu'il
existe un sous-ensemble parabolique standard $\tP$
strictement plus petit que
$\tG$ avec
\[
w_s\xi\in\tP(F)\ptf
\]
Donc $\delta$ appartient à
$\tP(F)$, ce qui contredit la primitivité de $\delta$.
On en déduit que $X=\HO(a)+S$ doit rester borné
ce qui impose à $a$ de rester dans un compact
et donc $\delta$ appartient à un ensemble fini.
\end{proof}

\begin{lemme}\label{ssfini}
Soit $\Omega$ un compact de $\tGadef$ et $\tP$ un sous-ensemble parabolique.
L'ensemble des $\delta\in\tM(F)$ qui sont quasi semi-simples et
tels qu'il existe $x\in\Gadef$ et $n\in\Nadef$ avec
\[
x\moins\delta\,n\,x\in\Omega
\]
appartiennent à un ensemble fini de classes de $\M(F)$\hyph conjugaison.
\end{lemme}

\begin{proof}
Pour tout élément quasi semi-simple $\delta\in\tM(F)$
il existe un sous-ensemble parabolique
\[
\tP_1=\tM_1\N_1\subset\tP
\]
et $\delta_1\in\tM_1(F)$ tel que $\delta_1$ soit un conjugué
de $\delta$, et soit un élément primitif pour $\tM_1$. On a donc
\[
\delta=\gamma\moins\delta_1\gamma
\]
pour un $\gamma\in\M(F)$
et
\[
x_1\moins\delta_1\,n'\,x_1\in\Omega
\]
avec $x_1=\gamma x$. Mais $x_1=m_1n_1k_1$ avec $m_1n_1\in\M_1(\adef)\N_1(\adef)$
et donc
\[
m_1\moins\delta_1\,m_1n''\in\K\Omega\K
\]
d'où on déduit que
\[
m_1\moins\delta_1\,m_1\in\Omega'
\]
où $\Omega'$ est un compact dans $\tM_1(\adef)$.
Compte tenu du théorème~\ref{siegel}
on voit qu'il existe $m_2\in\M_1(\adef)$
appartenant à un domaine de Siegel pour $\M_1$ et $\xi\in\tM_1(F)$
conjugué de $\delta_1$ tels que l'on ait
\[
m\moins\xi\,m\in\Omega'\ptf
\]
On invoque alors le lemme~\ref{geofini} et la finitude du nombre de classes
de conjugaison de sous-ensembles paraboliques.
\end{proof}

\part{Théorie spectrale, troncatures et noyaux}

\chapter{L'opérateur de troncature}\label{ch4}

Nous utiliserons les mesures
de Tamagawa sur les groupes adéliques unipotents. En particulier, si $\N$ est
un groupe unipotent,
le quotient $\N(F)\bs\N(\adef)$ est de volume~1.

\section{Définition et une propriété d'annulation}\label{defpropannul}

Considérons une fonction $\varphi$ dans
\[
L^1_{\mathrm{loc}}(\Q(F)\bs\Gadef) \ptf
\]
Soit $P$ un sous-groupe parabolique.
Le terme constant de $\vf$ le long de $P$ sera noté
$\Pi_P\vf$ ou $\vf_P$
\[
\Pi_P\vf(x)=\varphi_P(x)=\int_{\N_P(F)\bs\N_P(\adef)}\varphi(n\,x)\dd n
\]
où $N_P$ est le radical unipotent de $P$.

On 
définit pour $T\in\gao$ un opérateur de troncature par
\[
\tronc^{T,\Q}\varphi(x)=\sum_{\PO\subset P\subset\Q}
 (-1)^{\ga_P-\ga_\Q}\sum_{\xi\in P(F)\bs\Q(F)}
\htau_P^\Q(\xixT)\,\varphi_P(\xix)\ptf
\]
\newindex{LambdatQ@$\Lambda^{T,Q}$}{deftronc}%
On observe que, d'après le lemme~\ref{finitude}, les sommes en $\xi$ portent sur des ensembles finis.
Dans le cas où $Q=G$ on écrira le plus souvent $\tronc^\T$ pour $\tronc^{T,\G}$.
Une propriété importante de l'opérateur de troncature 
est donnée par le lemme suivant:

\begin{lemme}[\cite{ATFII}*{Lemma~1.1}]\label{troncnul} 
Supposons $\dPO(T)\ge c$ où $c$ est la constante du lemme~\ref{nilneg}.
Soit $\Q =M_\Q N_\Q$ un sous-groupe parabolique\textup,
\[
(\Pi_\Q\circ\tronc^\T \varphi)(x)\coloneqq\int_{N_\Q (F)\bs N_\Q (\adef)} (\tronc^\T \varphi)(n x)\dd n\ne0
\]
implique
\begin{equation}
\vpi(\xT)\le 0\qquad\forall\vpi\in\hDelta_\Q^\G\tag{i}\label{eq4.i}
\end{equation}
ou\textup, ce qui est équivalent\textup,
\begin{equation}
\phi_\Q^\G(\xT)=1\ptf\tag{ii}\label{eq4.ii}
\end{equation}
En particulier\textup, pour $T$ assez régulier et $\Q \ne\G$
\begin{equation}
\htau_{\Q }(\xT)\int_{N_\Q (F)\bs N_\Q (\adef)} (\tronc^\T \varphi)(n x)\dd n=0\ptf\tag{iii}\label{eq4.iii}
\end{equation}
\end{lemme}

\begin{proof} 
Considérons un sous-groupe parabolique standard $P=\M\N$ et posons
\[
A_P=\int_{\N_\Q (F)\bs\N_\Q (\adef)} \sum_{\xi\in P(F)\bs G(F)} \psi(\xi n\,x)\dd n 
\]
où
\[
\psi(x)= \htau_P(\xT) \varphi_P(x) \ptf
\]
On rappelle que d'après le lemme~\ref{weylgab}
\[
\weyl(\gao,P) = \{s\in \weyl \mid \text{$s^{-1} \alpha >0$ pour $\alpha\in \Delta_{\PO}^P$}\}
\]
est l'ensemble des représentants de longueur minimale pour les classes du quotient
\[
\weyl^\M\bs\weyl^\G\ptf
\]
Si $\Rac^\M \subset (\gao^G)^*$ est l'ensemble des racines réduites de $M$, 
on pourra remarquer que pour 
$s\in \weyl(\gao,P)$ on a l'équivalence
$\alpha>0 \Leftrightarrow s^{-1}\alpha>0$ 
pour les $\alpha\in \Rac^\M$. On a la décomposition de Bruhat
\[
G=\coprod_{s\in \weyl(\gao,P)} P\, w_s\, N_0 
\]
où l'on a choisi un représentant $w_s$ pour chaque $s$; de plus
\[
P\bs P\, w_s \, N_0 \cong N_s \bs N_0
\]
 où l'on a posé $N_s = w_s^{-1}N_0 w_s \cap N_0$.
On a donc
\[
A_P=\sum_{s\in\weyl(\gao,P)}A_{s,P}
\]
avec
\[
A_{s,P}=\int_{N_\Q (F)\bs N_\Q (\adef)}\ \sum_{\nu\in N_s(F)\bs N_0(F)} \psi(w_s\nu n x)\dd n\ptf
\]
Fixons $s$ et $w_s$.
Soit $N_0^1 =N_0 \cap M_\Q $. Alors
\[
N_\Q (F)\bs N_\Q (\adef) = N_0(F)\bs N_0^1(F) N_\Q (\adef)
\]
et $A_{s,P}$ se récrit :
\[
A_{s,P}=\int_{N_s(F) \bs N_0^1(F) N_\Q (\adef)} \psi(w_s\,n_1 x)\dd n_1
\]
où la mesure $\dd n_1$ contient la mesure discrète sur $N_0^1(F)$.
 écrivons ceci comme une intégrale itérée 
\[
A_{s,P}=\int_{\N^s}\biggl(\int_{N_s^*}\psi(w_s\,n_s^* n^sx)\dd n_s^*\biggr)\dd n^s
\]
où $n^s$ décrit
\[
N^s = w_s^{-1} N_0(\adef)w_s\cap N_0^1 (F) N_\Q (\adef) \bs N_0^1(F)N_\Q (\adef)
\]
et où $n_s^*$ décrit
\[
N_s^*= N_s(F) \bs w_s^{-1}N_0(\adef) w_s \cap N_0^1(F) N_\Q (\adef)\ptf
\]
La décomposition radicielle montre que
\[
N_0^1(F)N_\Q (\adef) \cap w_s^{-1}N_0(\adef)w_s 
\]
est égal à
\[
 (N_0^1(F) \cap w_s^{-1}N_0(F)w_s)(N_\Q (\adef) \cap w_s^{-1}N_0(\adef)w_s)\;,
\]
le premier facteur étant contenu dans $N_s(F)$. On remarque que
 $N_s(F)$ est en fait le produit de $N_\Q (F) \cap w_s^{-1}N_0(F)w_s$ 
et de ce facteur, et l'on peut donc récrire l'intégrale sur $N_s^*$ comme une intégrale sur
\[
 N_* = N_\Q (F) \cap w_s^{-1} N_0(F) w_s\bs N_\Q (\adef) \cap w_s^{-1} N_0(\adef)w_s
\]
\cad
\[
A_{s,P}=\int_{\N^s}\biggl(\int_{N_*}\psi(w_s\, n_* n^sx)\dd n_*\biggr)d\dd ^s\ptf
\]
On a remarqué que $w_s N_0 w_s^{-1} \cap M = N_0 \cap M$. Le sous-groupe 
\[
P'_s = w_s \Q w_s^{-1} \cap M
\]
de $M$ contient $\NO\cap M$; c'est donc un sous-groupe parabolique standard de $M$, 
de radical unipotent
\[
N'_s=w_sN_\Q w_s^{-1} \cap M\ptf
\]
En particulier $N'_s \subset N_0$, et la décomposition $N_0 =N(M\cap N_0)$ 
implique que
\[
N_0 \cap w_s N_\Q w_s^{-1}= N''_sN'_s
\]
où
\[
N''_s=N\cap w_sN_\Q w_s^{-1}\ptf
\]
Le changement de variable $w_sn_*w_s^{-1}=n''_s n'_s$ donne alors une intégrale sur le produit
\[
\bigl(N''_s(F) \bs N''_s(\adef)\bigr) \times \bigl(N'_s(F)\bs N'_s(\adef)\bigr)\;,
\]
et on obtient
\[
A_{s,P}=\iiint\psi (n''_s n'_sw_s\, n^sx)\dd n''_s\dd n'_s\dd n^s\ptf
\]
Comme $N''_s\subset\N$,
l'intégrale sur le quotient 
$N''_s(F) \bs N''_s(\adef)$
peut être omise compte tenu de l'invariance à gauche de $\psi$ par $N(\adef)$, les mesures étant normalisées
de sorte que ce quotient soit de volume~1. 
On a donc obtenu pour $A_{s,P}$ l'expression suivante : 
\[
\int_{N^s}\htau_P(\HO(w_sn^sx)-T)\biggl(\int_{N'_s(F)\bs N'_s(\adef)}
\varphi_P( n'_s w_s n^s x) \dd n'_s\biggr)\dd n^s
\]
en utilisant que $\HO$ est invariant à gauche par $n'_s \in N_0(\adef)$. 
Par ailleurs le sous-groupe $P'_s$ est l'intersection avec $M$ 
d'un unique sous-groupe parabolique standard $R=P'_sN$ de $G$; en particulier $\R\subset P$.
Son radical unipotent 
est le sous-groupe $N_\R= N'_sN$.
On a donc, en désignant par $\varphi_\R$ le terme constant de $\varphi$ le long de $\R$:
\[
A_{s,P}=\int_{N^s} \varphi_\R(w_s n^sx) \htau_P(\HO(w_sn^sx)-T)\dd n^s\ptf
\]
Posons 
\begin{align*}
\Sigma^1 &= \{\alpha\in \Delta_{\PO} \mid s^{-1} \alpha>0, s^{-1} \alpha\rest \ga_{\Q }= 0\}\\
\Sigma_1 &= \{\alpha\in \Delta_{\PO} \mid s^{-1} \alpha>0,s ^{-1} \alpha\rest \ga_{\Q } \ne 0\}
\end{align*}
et
\[
\Sigma=\Sigma_1\cup\Sigma^1=\{\alpha\in \Delta_{\PO} \mid s^{-1} \alpha>0\}
\ptf
\]
On rappelle que $R=P'_sN$ et donc le sous-groupe de Levi $\M_\R$ de $\R$ est contenu dans 
\[
P'_s = w_s \Q w_s^{-1} \cap M
\]
et donc
\[
\M_\R=\M\cap w_s\M_\Q w_s\moins
\]
d'où on déduit que les racines simples dans $\M_\R$ sont les racines simples de $P$ qui
s'annulent sur $s(\ga_\Q)$. Elles vérifient $s\moins\alpha>0$ puisque $s\in\weyl(\gao,P)$
et donc $\Delta_{\PO}^P\subset\Sigma$. On en déduit que

\[
\Delta_{\PO}^{\R}=\Delta_{\PO}^P\cap\Sigma^1 \subset\Sigma^1\ptf
\]
Notons $S$ le sous-groupe parabolique tel que $\Delta_{\PO}^S=\Delta_{\PO}^\R\cup\Sigma_1$.
Il résulte des remarques qui précèdent que
\[
\Delta_{\PO}^P\subset\Delta_{\PO}^S\subset\Sigma\ptf
\]
Les $P$ contenant $\R$ et tels que
\[
R\cap\M=P'_s
\]
sont en bijection avec les sous ensembles de $\Sigma_1$
ce qui est équivalent à demander que $\R\subset P\subset S$.
On veut calculer
\[
(\tronc^\T \varphi)_\Q(x)=\int_{N_\Q (F)\bs N_\Q (\adef)} (\tronc^\T \varphi)(n x)\dd n
=\sum_{\PO\subset P\subset\G} (-1)^{\ga_P-\ga_\G} \sum_{s\in\weyl(\gao,P)}A_{s,P}\ptf
\]
La dernière expression peut se récrire comme une somme de termes associés aux couples $(s,\R)$ où $\R$ est un sous-groupe parabolique admettant un sous-groupe de Levi vérifiant:
\[
\M_\R\subset w_s\M_\Q w_s\moins\ptf
\]
On obtient
\[
(\tronc^\T \varphi)_\Q(x)=\sum_{(s,\R)}\sum_{\{P\mid \R\subset P\subset S\}} (-1)^{a_P-a_\G}A_{s,P}\ptf
\]
On observe que l'ensemble $N^s$ ne dépend pas de $P$ et donc
\[
(\tronc^\T \varphi)_\Q(x)=\sum_{(s,\R)}\int_{\N^s}\varphi_\R(w_s\,n^sx)
\sum_{\{P\mid \R\subset P\subset S\}} (-1)^{a_P-a_\G} \htau_P(\HO(w_s\,n^sx)-T)\dd n^s
\]
soit encore
\[
(\tronc^\T \varphi)_\Q(x)=\sum_{(s,\R)}(-1)^{a_S-a_\G}\int_{\N^s}\varphi_\R(w_s\,n^sx) \phi_\R^{S,\G}(\HO(w_s\,n^sx)-T)\dd n^s\ptf
\]
où (avec les notations du lemme~\ref{prepar})
\[
\phi_\R^{S,\G}(H) = \sum_{\{P\mid \R\subset P\subset S\}} (-1)^{a_P-a_S} \,\htau_P(H)\ptf
\]
Pour conclure la preuve de~\eqref{eq4.i}
on invoque le lemme~\ref{phinul} ci-dessous. L'équivalence de~\eqref{eq4.i} et~\eqref{eq4.ii} n'est autre que le lemme~\ref{prepar}.
L'assertion~\eqref{eq4.iii} est alors immédiate.
\end{proof}

\begin{lemme}\label{phinul}
Soient $\Q$\textup, $\R$ et $S$ comme ci-dessus. Alors\textup, si $\dPO(T)\ge c$\textup,
\[
\phi_\R^{S,\G}(\HO(w_s\,n^sx)-T)\ne0
\]
implique
\[
\vpi(\xT)\le0\qquad\forall\vpi\in\hDelta_\Q^\G\ptf
\]
\end{lemme}

\begin{proof}
D'après le lemme~\ref{prepar}, $\phi_\R^{S,\G}$
est la fonction caractéristique de l'ensemble des $H$ tels que $ \vpi(H)>0$ pour 
$\vpi\in \hDelta_S^\G$ et $ \vpi(H) \le 0$ 
pour $\vpi\in \hDelta_R^\G-\hDelta_S^\G$.
On suppose
\[
\phi_\R^{S,\G}(\HO(w_sn^sx)-T)\ne0
\]
\cad que
\[
\HO(w_sn^sx) -T = \sum_\alpha t_\alpha\, \alpha^\vee
\]
avec
\[
t_\alpha>0 \quad \alpha\in \Delta_{\PO}^\G-\Delta_{\PO}^S\Qquad{et}
t_\alpha \le 0 \quad \alpha\in \Delta_{\PO}^S-\Delta_{\PO}^\R=\Sigma_1\ptf
\]
Pour $\vpi \in \hDelta_{\Q }$ on a
\[
\vpi(s^{-1}\bigl(\HO(w_sn^sx) -T)\bigr) = \sum_{\alpha\in \Delta_{\PO}} t_\alpha\, \vpi(s^{-1}\alpha^\vee)\le0\ptf
\]
En effet, chaque terme est $\le0$:
\begin{itemize}
\item  si $\alpha\notin\Sigma$ alors $\alpha\notin\Delta_{\PO}^S\subset\Sigma$
et donc $s^{-1}\alpha^\vee<0$ et $t_\alpha>0$; 
\item si $\alpha\in \Sigma_1$ alors $s^{-1}\alpha^\vee>0$ et $t_\alpha\le 0$; 
\item si $\alpha\in \Sigma^1$  alors $s\moins\alpha\in\Delta_{\PO}^\Q$ et donc $\vpi(s^{-1}\alpha^\vee)=0$. 
\end{itemize}
 Mais par ailleurs,
\[
s^{-1} \HO(w_sn^sx) = \HO(x) +s^{-1} \HO(w_sn)
\]
pour un $n\in N_0(\adef)$ et donc
\[
s^{-1} (\HO(w_sn^sx) -T) =( \xT) +Y_1(n,\T)-Y_s(n,\T)
\]
avec
\[
Y_s(n,\T)= s^{-1} \bigl(T - \HO(w_sn)\bigr)\ptf
\]
D'après le lemme~\ref{wn}, si $\dPO(T)\ge c$, la famille des $Y_s(n,\T)$ est une famille orthogonale régulière. 
Dans ce cas
\[
\vpi \bigl( Y_1(n,\T)-Y_s(n,\T)\bigr)\ge 0
\]
pour tout $\vpi$ et donc $ \vpi(\xT)\le 0$ pour $\vpi \in{\hDelta}_{\Q}$.
\end{proof}

Le lemme~\ref{troncnul} a pour conséquence immédiate le

\begin{corollaire}\label{invol}
Supposons $\dPO(T)\ge c$. L'opérateur $\tronc^\T$ est un idempotent\textup:
\[
\tronc^\T(\tronc^\T\varphi) = \tronc^\T\varphi\ptf
\]
\end{corollaire}

\section{Un raffinement}

Soit $\Q$ un sous-groupe parabolique standard. Pour $X\in \ga_\Q^G$ et $\T\in\ga_0^\G$, on définit un élément
\[
T[X]\in \ga_0^\Q
\]
comme suit: c'est l'unique élément de $\ga_0^\Q$ tel que, pour tout 
$\varpi\in \hDelta_0^\G-\hDelta_\Q^\G$, on ait l'égalité 
\[
\varpi(T[X])=\varpi(T-X)\ptf
\]
En d'autres termes, si

\[
T-X=\sum_{\alpha\in\Delta_0^\G}x_{\alpha}\alpha^\vee\qquad\text{alors 
$T[X]=\sum_{\alpha\in\Delta_0^\Q}x_{\alpha}\alpha^\vee$}\ptf
\]

\begin{lemme}\label{raff}
On a
\begin{equation}
T[X]=T^\Q-\sum_{\alpha\in \Delta_0^\G-\Delta_0^\Q}x_{\alpha}(\alpha^\vee)^\Q\ptf\label{eq4.1}
\end{equation}
De plus, si
\begin{equation}
\phi_\Q^\G(X-T)=1\label{eq4.2}
\end{equation}
alors $T[X]$ est \og plus régulier\fg 
que $\T^\Q$:
 pour $\beta\in\Delta_0^\Q$ on a 
\[
\beta(\T[X])\ge\beta(\T^\Q)\ptf
\]
\end{lemme}

\begin{proof}
On observe que $(\alpha^\vee)^\Q=\alpha^\vee$ pour $\alpha\in\Delta_0^\Q$ et que
\[
T^\Q=(T-X)^\Q=\sum_{\alpha\in \Delta_0^\G}x_{\alpha}(\alpha^\vee)^\Q
\]
et~\eqref{eq4.1} s'en déduit.
Maintenant la condition~\eqref{eq4.2} implique que
les $x_{\alpha}$ sont positifs ou nuls pour $\alpha\notin 
\Delta_0^\Q$ et comme $\beta(X)=0$ pour $\beta\in\Delta_0^\Q$ on a donc,
compte tenu de~\eqref{eq4.1}:
\[
\beta(\T[X])=\beta(T)-\sum_{\alpha\notin 
\Delta_0^\Q}x_{\alpha}\beta(\alpha^\vee)\ge\beta(\T)=\beta(\T^\Q)
\]
puisque $\beta(\alpha^\vee)\le0$.
\end{proof}

Rappelons que l'opérateur $\tronc^{T,Q}$ ne dépend en fait que de $T^Q$. 
Le lemme~\ref{troncnul} se raffine en celui qui suit. 

\begin{lemme} \label{LemmeB}
 Il existe $c'>0$ tel que\textup, pour tout $c>0$ et tout $\T\in\ga_0^\G$ vérifiant $\dPO(T)\geq c'(c+1)$\textup, la propriété 
suivante soit vérifiée. Soit $x\in G({\adef})$ et notons
$X$ la projection de $\HO(x)$ sur $\ga_\Q^\G$.
On suppose que
\[
\lVert X-T_\Q\rVert \leq c\ptf
\]
Alors\textup, pour toute fonction $\vf$ sur $G(F) \backslash G({\adef})$\textup, on a l'égalité 
\[
(\tronc^{T}\vf)_\Q(x)=(-1)^{a_\Q-a_\G} \phi_\Q^{\G}(X-T)\tronc^{T[X],\Q}(\vf_\Q)(x)\ptf
\]
\end{lemme}

\begin{proof} 
D'après la preuve des lemmes~\ref{troncnul} et~\ref{phinul}, on sait que
\[
(\tronc^\T \varphi)_\Q(x)=\sum_{(s,\R)}(-1)^{a_S-a_\G}\int_{\N^s}\varphi_\R(w_s\,n^sx) 
\,\phi_\R^{S,\G}(\HO(w_s\,n^sx)-T)\dd n^s
\]
et pour que
\[
\phi_\R^{S,\G}(\HO(w_s\,n^sx)-T)
\]soit non nul on doit avoir
\[
\vpi\bigl(s^{-1}(\HO(w_sn^sx) -T)\bigr) = \sum_{\alpha\in \Delta_{\PO}} t_\alpha\, \vpi(s^{-1}\alpha^\vee)\le0
\]
pour $\vpi \in \hDelta_{\Q }$. On a vu (\cf lemme~\ref{phinul}) que
\[
s^{-1} (\HO(w_sn^sx) -T) =( \xT) +Y_1(n,\T)-Y_s(n,\T)\ptf
\]
Par hypothèse 
\[
\vpi( \xT)=\vpi(X-\T_\Q)
\]
reste borné
et donc $\vpi\bigl(Y_1(n,\T)-Y_s(n,\T)\bigr)$ doit être majoré. Pour $\dPO(\T)$ assez grand, de façon précise si
$\dPO(T)\geq c'(c+1)$ avec $\lVert X-T_\Q\rVert \leq c$,
ceci impose $s\in\weyl^\Q$. On a alors 
\[
\Sigma^1\subset\Delta_{\PO}^\Q\Qquad{et donc}\R\subset\Q
\]
et
\[
\Delta_{\PO}^S=\Delta_{\PO}^\R\cup\Sigma_1
\qquad\text{avec $\Sigma_1=\Delta_{\PO}^\G-\Delta_{\PO}^\Q$}\ptf
\]
 On en déduit que $S$ ne dépend que de $\R$ et $\Q$:
 de fait, on a
\[
\hDelta_S=\hDelta_\R-\hDelta_\Q\ptf
\]
On rappelle que $N^s$
est un quotient de $N_\Q(\adef)\bigl(\NO(F)\cap\Q(F)\bigr)$.
On obtient alors
\begin{equation}
(\tronc^\T \varphi)_\Q(x)=\sum_{\{\R\mid\PO\subset\R\subset\Q\}}
(-1)^{a_\R-a_\G}
\sum_{\xi\in\R(F)\bs\Q(F)}\varphi_\R(\xix) 
\,\phi_\R^{S,\G}(\HO(\xix)-T)\ptf\label{eq4.1a}
\end{equation}
D'après le lemme~\ref{prepar}, la fonction
$\phi_\R^{S,\G}$ est la fonction caractéristique des $H\in\gao$ tels que 
\[
\vpi(H)\le0\qquad\text{pour tous les $\vpi\in\hDelta_\R^\G-\hDelta_S^\G=\hDelta_\Q^\G$}
\]
 et 

\[
\vpi(H)>0\qquad\text{pour tous les $\vpi\in\hDelta_S^\G=\hDelta_\R^\G-\hDelta_\Q^\G$}\ptf
\]
On a donc
\[
\phi_\R^{S,\G}(H)=\phi_\Q^{\G}(H)\htau_\R^\Q\circ\QpQ(H)
\]
où $\QpQ(H)$ est l'élément de $\ga_{\PO}^\Q$
défini par les équations
\[
\vpi\bigl(\QpQ(H)\bigr)=\vpi(H)\qquad\text{pour tous les $\vpi\in\hDelta_{\PO}^\G-\hDelta_\Q^\G$}\ptf
\]
Donc 
\[
\phi_\R^{S,\G}(\HO(\xix)-T)=\phi_\Q^{\G}(\xixT)\htau_\R^\Q\circ{\QpQ}(\HO(\xix)-T)\ptf
\]
On observera que
\[
\HO(\xix)_\Q=\HO(x)_\Q=X
\]
et que, pour tous les $\vpi\in\hDelta_{\PO}^\G-\hDelta_\Q^\G$ on a,
par définition de $\T[X]$,
\[
\vpi(\xixT)=\vpi(\HO(x)_\Q+\HO(\xix)^\Q-\T)=\vpi(\HO(\xix)^\Q-\T[X])
\]
\cad que
\[
\QpQ(\HO(\xix)-T)=\HO(\xix)^\Q-T[X]\ptf
\]
On a donc aussi
\begin{equation}
\phi_\R^{S,\G}(\HO(\xix)-T)=\phi_\Q^{\G}(X-T)\htau_\R^\Q(\HO(\xix)-T[X])\ptf\label{eq4.2a}
\end{equation}
Le lemme résulte immédiatement de~\eqref{eq4.1a} et~\eqref{eq4.2a}.
\end{proof}

\section{Troncature et décroissance}

Fixons $\Q\subset\R$ deux sous-groupes paraboliques.
Pour $P=M_P\N_P$ avec $\Q \subset P\subset \R$ on a
$\N_\R\subset \N_P\subset N_\Q$. 
Posons $\Sigma_P=\Delta_{\PO}^\R -\Delta_{\PO}^P$.
Le radical unipotent $N_P$ de $P$ est le produit de $N_\R$ et de $U=M_\R\cap\N_P$.
 Soit $P_\alpha$ le groupe associé à $\alpha\in\Sigma_P$
 \cad tel que
\[
\Sigma_{P_\alpha}=\{\alpha\}
\]
et soit $U_\alpha$ l'intersection de son radical unipotent avec $\M_\R$. 
Alors, 
\[
N_P= N_\R \prod_{\alpha\in \Sigma_P} U_\alpha\ptf
\]
Soit $\psi$ une fonction continue sur $N_\Q (F)\bs N_\Q(\adef)$. Considérons
\[
\psi_P(n_1)=\Pi_P \psi(n_1)\coloneqq \int_{N_P(F)\bs \N_P(\adef)} \psi(nn_1)\dd n\ptf
\]
La fonction $\Pi_P \psi$
est encore invariante à gauche par $N_\Q (F)$. 
Tous les groupes unipotents considérés sont contenus dans $N_\Q $ et 
invariants par celui-ci. On peut donc introduire
\[
\Theta \psi = \sum_{\Q \subset P\subset \R}(-1)^{a_P-a_\R}\Pi_P \psi= \sum_{\Q \subset P\subset \R}(-1)^{a_P-a_\R} \psi_P
\ptf
\]

On dispose de l'action à droite, notée $\reg(X)$, des opérateurs $X$ de l'algèbre enveloppante $\goth U(\gn_\Q)$,
sur les fonctions lisses sur $\N_\Q(F)\bs\N_\Q(\adef)$.

\begin{lemme}\label{unipcroiss} 
Soit $\mathcal{O}$ un sous-groupe ouvert compact du groupe
des adèles finis de $\N_\Q$. Pour tout entier $r\ge 0$ 
il existe des opérateurs différentiels $X_{Q,R}$\textup, définis par des éléments
de l'algèbre enveloppante de $N_\Q (F\otimes\RM)$\textup, de la forme suivante\textup:
\[
X_{Q,R}=\prod_{\alpha\in\Delta_\Q^\R} \Biggl( \sum_{j=1}^{n_\alpha} Y_{\alpha,j}^r\Biggr)\qquad\text{avec
$Y_{\alpha,j}\in\gu_\alpha=\Lie U_\alpha$}
\]
et tels que
\[
\lVert \Theta \psi \rVert_\infty \le  \lVert \reg(X_{Q,R})\psi\rVert_\infty
\]
pour toute fonction $\psi$ lisse sur $N_\Q(F) \bs \N_\Q(\adef)/\mathcal{O}$.
\end{lemme}

\begin{proof} 
On peut supposer $\psi$ invariante à gauche par $N_\R (\adef)$.
On a alors
\[
\Theta \psi=\sum_P (-1)^{\lvert\Sigma_P\rvert} \prod_{\alpha\in\Sigma_P}\Pi_{P_\alpha}\psi
=\sum_P (-1)^{\lvert\Sigma_P\rvert} \prod_{\alpha\in\Sigma_P}\Pi_{U_\alpha}\psi
= \prod_{\alpha\in\Delta_\Q^\R} \bigl(1 - \Pi_{U\alpha}\bigr) \psi
\]
où l'on a noté $\Pi_{\U_\alpha}$ l'intégrale sur $U_\alpha(F)\bs U_\alpha(\adef)$. 
Comme l'intégration sur
 un sous-groupe unipotent (agissant à gauche)
 commute avec l'action à droite d'un opérateur différentiel,
on peut traiter séparément chaque facteur $(1 - \Pi_{U\alpha})$. 
On doit donc estimer
\[
 \psi(n_1) - \int_{U_\alpha(F)\bs U_\alpha(\adef)} \psi(un_1)\dd u\ptf
\]
On considère une suite de composition de~$U_\alpha$
\[
\{1\} \subset V_1 \subset V_2 \subset \dots \subset V_{n_\alpha} = U_\alpha
\]
dont les quotients sont isomorphes au groupe additif. On remarque que
\[
1-\Pi_{V_j} = 1 -\Pi_{V_{j-1}} + (1- \Pi_{V_{j-1}\bs V_j}) \Pi_{V_{j-1}}
\]
et donc
\[
1-\Pi_{U_\alpha} = \sum_{j=1}^{n_\alpha} (1- \Pi_{V_{j-1}\bs V_j}) \Pi_{V_{j-1}}\ptf
\]
En observant que, si $V$ est l'un quelconque de $V_j$, on a
\[
\lVert\Pi_V \psi \rVert_\infty \le \lVert\psi\rVert_\infty
\]
on est ramené à traiter le cas d'un seul facteur $1-\Pi_V$
 ou $V$ est de dimension~1. Le lemme résulte alors
 de ce que, si
$\psi$ est une fonction sur $\RM/\ZM$ ayant pour développement
de Fourier
\[
\psi(x)=\sum_{n\in\ZM}\ a_n\ee^{2 \pi \ima nx}
\]
on a pour $r>0$
\[
\lVert(1-\Pi_V)\psi\rVert_\infty= \lVert\psi -a_0\rVert_\infty \le \sum_{n\ne 0}\lvert a_n\rvert 
\le \Bigl(\sum n^{2r} \lvert a_n\rvert^2\Bigr)^{1/2} \Biggl(\sum_{n\ne 0}\frac{1}{n^{2r}}\Biggr)^{1/2}
\]
et donc
\[
\lVert\psi -a_0\rVert_\infty \le c_r \biggl\lVert\frac{\partial^r \psi}{\partial n^r} 
\biggr\rVert_\infty\ptf\qedhere
\]
\end{proof}

Soit $\varphi$ une fonction sur $\XG=\gA_G G(F) \bs G(\adef)$. On dit que $\varphi$ est 
à croissance lente si, pour un certain $N$ et pour tout
$x$ dans un domaine de Siegel $\Sieg$, on a
\[
\lvert \varphi(x)\rvert \ll \lvert x\rvert^N\ptf
\]
On dira que $\varphi$ est à croissance uniformément lente, si $\varphi$ est $\K$\hyph finie 
(à droite) et à 
croissance lente ainsi que toutes ses dérivées (pour l'action à droite des opérateurs de l'algèbre enveloppante)
pour le même exposant $N$.
On dit que $\varphi$ est à décroissance rapide si, pour tout $N$,
\[
\lvert \varphi(x)\rvert \le c_N\lvert x\rvert^{-N}
\]
pour $x\in\Sieg$.

\begin{proposition}\label{rapdec}
Soit $\vf$ une fonction lisse sur $\XG$
 à croissance uniformément lente.
Soit $\Sieg$ un domaine de Siegel dans $\G(\adef)$.
Pour tout couple d'entiers positifs $A$ et
$B$\textup, il existe un ensemble fini $X_1,\dots,X_r$ d'éléments de $\goth U(\frakg)$ 
 tels que\textup:
\[
\lvert \tronc^\T \varphi(x)\rvert\,\lvert x\rvert^{B} \ll \sum_i \sup_{y\in \Gadef}
 \lvert \reg(X_i) \varphi(y)\rvert\, \lvert y\rvert^{-A} 
\]
pour tout $x\in\Sieg$. En d'autres termes\textup:
si $\varphi$ est à croissance uniformément lente sur $\XG$\textup, alors
$\tronc^\T\varphi$ est à décroissance rapide.
\end{proposition}

\begin{proof}
Commençons par insérer dans l'expression de $\tronc^\T$ l'identité du lemme \ref{FPQR}
\[
\sum_{\Q\subset P\subset \R} \, \sum_{\xi\in \Q(F)\bs P(F)} 
\FPQ(\xix,T) \sQR(\xixT)=\htau_P(\xT)
\]
pour chaque $P$. Donc
\[
\tronc^\T \varphi(x)=
\sum_{\{\Q,\R\mid \Q \subset\R\}}A_{Q,R}(x)
\]
avec
\[
A_{Q,R}(x)= \sum_{\xi\in \Q(F) \bs G(F)}
 \FPQ(\xix,T)\sQR(\xixT)\sum_{\{P\mid \Q \subset P\subset\R\}}(-1)^{a_P-a_\G}\varphi_P(\xix)\ptf
\]
Fixons $\Q$ et $\R$ et supposons tout d'abord que $\Q=\R$.
On a alors $\sQR=0$ sauf dans le cas
$\Q=\R=\G$. Dans ce cas on simplement
\[
A_{G,G}(x)= F_{\PO}^\G(x,T)\varphi(x)\ptf
\]
Cette fonction est à décroissance rapide car
 $F_{\PO}^\G$ est à support compact et 
 l'inégalité est vérifiée en prenant pour $X$ un opérateur
 de degré zéro \cad une constante non nulle (dépendant de $\vf$ et de $T$).
Supposons désormais $\Q\ne\R$. Il nous suffit de majorer
\[
\sum_{\{P\mid \Q\subset P\subset\R\}} (-1)^{a_P-a_\G} \varphi_P(\xix)\;,
\]
sous la condition
\begin{equation}
\FPQ(\xix,T) \sQR (\xixT)=1\ptf\label{eq4.1b}
\end{equation}
Rappelons que pour $x$ fixé il y a au plus un $\xi$ modulo $\Q(F)$
pour lequel un tel terme est non nul (\cf lemme~\ref{unique}).
Compte tenu du lemme~\ref{ratmaj} il suffit de montrer que pour tout $A$ et tout $B$
on peut trouver des opérateurs différentiels $X_i$ tels que
\begin{equation}
\lvert \xix\rvert^B\Biggl\lvert\sum_{\{P\mid \Q\subset P\subset\R\}}(-1)^{a_P-a_\G}
\varphi_P(\xix)\Biggr\rvert \ll \sum_i\sup_y \lvert \reg(X_i)\varphi(y)\rvert \, \lvert y\rvert^{-A}\ptf\label{eq4.2b}
\end{equation}
On s'intéresse aux $x\in\Sieg$; on a donc
\[
\HG(\xix)=\HG(x)=0
\]
et on remarque 
que l'on est libre de multiplier $\xix$ par un élément de $\Q(F)$ à gauche.
Mais $\xix$ vérifie \eqref{eq4.1b}, en particulier
$\FPQ(\xix,T)=1$. Il résulte alors du lemme~\ref{gSPQ} que l'on peut supposer 
\[
\xix=n_1ac \qquad\text{avec $n_1\in\N_\Q(\adef)$, $a\in \gA_0^\G$ et $c\in\Omega$}
\]
où $\Omega$ est un compact et où
la projection de $H=\HO(a)$ dans $\gao^\Q$ est bornée.
La translation par un compact à droite ne change rien aux estimées.
Il suffit donc d'étudier le cas où
\[
\xix=n_1a \qquad\text{avec $a=\ee^H$ et $H\in\ga_\Q^\G$}\ptf
\]
De plus on doit avoir
\[
\sQR(H-T) =1\ptf
\]
On décompose $H$ sous la forme 
\[
H=H_1+H_2
\]
où $H_1 \in \ga_{\Q }^\R$ et $H_2\in \ga_\R^\G$.
D'après l'équation~\eqref{eq2.i} du lemme~\ref{croiss} on a
\[
\lVert H_2\rVert \le c\lVert H_1-T\rVert \le \lVert H_1\rVert +\lVert T\rVert \ptf
\]
 Quitte à changer les constantes, il nous suffit de montrer que
\begin{equation}
\ee^{B'\lVert H_1\rVert} \Bigl\lvert\sum (-1)^{a_P} \varphi_P(n_1a)\Bigr\rvert
\label{eq4.3b}
\end{equation}
est dominé par le membre de droite de~\eqref{eq4.2b}. 
En appliquant le lemme~\ref{unipcroiss} aux fonctions de la forme
\[
\psi(n_1) = \varphi(n_1 a)
\]
on obtient
\[
\sup_{n_1}\Bigl\lvert\sum (-1)^{a_P} \varphi_P(n_1\,a)\Bigr\rvert \le 
\sup_{n_1}\lvert \reg(\Ad(a^{-1})X_{Q,R})\varphi(n_1a)\rvert\ptf
\]
Les opérateurs $X_{Q,R}$ sont de la forme
\[
X_{Q,R}=\sum_k Z_k\qquad\text{avec $Z_k=\prod_{\alpha\in\Delta_\Q^\R} Y_{\alpha,j_k}^r$}
\]
et les $Y_{\alpha,j_k}$ se transforment sous $\Ad(a)$ par une racine $\beta$
dont la décomposition en racines simples fait intervenir
 $\alpha$ avec un coefficient${}\ge1$.
Donc, pour tout $a=\ee^{H_1}$ avec $H_1\in\ga_\Q^\R$ on a
\[
\Ad(a)Z_k= \ee^{\lambda_k(H_1)}Z_k
\]
où $\lambda_k$ est une combinaison linéaire à coefficients entiers${}\ge r$ des racines
dans $\Delta_\Q^\R$. En particulier 
il existe une constante $a>0$ telle que
\[
\ee^{\lambda_k(H_1)}\gg \ee^{ra\lVert H_1\rVert}
\]
puisque $\sQR(H-\T)=1$ (la constante implicite dépendant de $\T$).
Il en résulte que, pour $r$ assez grand, 
\eqref{eq4.3b} est dominé par
\[
\ee^{-A\lVert H\rVert} \sup_i \lvert\reg(X_{Q,R})\varphi(n_1a)\rvert 
\ll \sup_{i,y} \lvert\reg(X_{Q,R})\varphi(y)\rvert\, \lvert y\rvert^{-A}\ptf\qedhere
\]
\end{proof}

La proposition~\ref{rapdec} admet la variante suivante.
 
\begin{proposition}\label{propoC}
Soit $\vf$ une fonction lisse sur $\XG$ à croissance uniformément lente.
 Pour tout couple 
d'entiers positifs $A$ et $B$\textup, il existe un ensemble fini d'éléments de $\goth U(\frakg)$\textup:  $X_1,\dots, X_
{r}$ tels que
\[
\lvert \tronc^T\vf(x)-F_{\PO}^G(x,T)\vf(x)\rvert\, \lvert x\rvert ^B \ll
 \ee^{-A\dPO(T)}\sum_{i}\sup_{y\in G({\adef})}\lvert 
\rho(X_{i})\vf(y)\rvert\,\lvert y\rvert^{-A}
\]
 pour tout $x\in \Sieg^G$.
\end{proposition}

\begin{proof} 
Il suffit de reprendre la preuve de la proposition~\ref{rapdec} ci-dessus, en observant que
pour contrôler les termes 
\[
A_{Q,R}(x)= \sum_{\xi\in \Q(F) \bs G(F)}
 \FPQ(\xix,T)\sQR(\xixT)\sum_{\{P\mid \Q \subset P\subset\R\}}(-1)^{a_P-a_\G}\varphi_P(\xix)
\]
 avec $\Q\ne\R$, on considère des $H=\HO(\xix)$ vérifiant
\[
\sQR(H-T) =1
\]
 et donc $\alpha(H)>\alpha(\T)$ pour $\alpha\in\Delta_\Q^\R$.
 Comme cet ensemble est non vide on a
\[
\lVert H\rVert \gg (\lVert H\rVert+\lvert\alpha(\T)\rvert)\ge\bigl(\lVert H\rVert+\dPO(\T)\bigr)\ptf\qedhere
\]
 \end{proof}
 
 Cette proposition admet elle-même une variante immédiate quand on remplace $\tronc^T$ par $\tronc^{T,Q}$. 

\section[$\tronc^\T$ comme projecteur]{\mathversion{bold}$\tronc^\T$ comme projecteur}

\begin{proposition}\label{autoad}
Soient $\varphi$ et $\psi$ des fonctions à croissance uniformément lente 
sur $\XG/\K_f$. Alors
\[
\langle \tronc^\T\varphi,\psi\rangle = \langle \varphi,\tronc^\T \psi\rangle
\]
les produits scalaires étant absolument convergents.
\end{proposition}

\begin{proof}
Si $\varphi$ est à croissance uniformément lente
 alors d'après la proposition~\ref{rapdec} la fonction
$\tronc^\T\varphi$ est à décroissance rapide et donc le produit scalaire 
$\langle \tronc^\T\varphi,\psi\rangle$
 est absolument convergent et dépend continûment de $\varphi$ pour la semi-norme:
\[
\lVert\varphi\rVert= \sum_i \Sup_{x\in \Gadef} \lvert\reg(X_i) \varphi(x)\rvert\, \lvert x\rvert^{-N}\ptf
\]
Maintenant si $\varphi$ est à support compact sur $\XG$ on a
\[
 \int_{\XG} \varphi(x)\Biggl(\sum_{P(F)\bs G(F)} 
\htau_P(\xixT)\,\widebar{\psi_P(\xix)}\Biggr) \dd x
\]
soit encore
\[
{}= \int_{\gA_G P(F)\bs G(\adef)} \htau_P (\xT) \, 
 \varphi(x)\,\widebar{\psi_P(x)}\dd x
\]
 et aussi
\[
{}= \int_{\gA_G P(F)N(\adef)\bs G(\adef)} 
\htau_P(\xT)\varphi_P(x) \widebar{\psi_P(x)} \dd x
\]
vu l'invariance de $\HO$ sous $N_0(\adef)$; l'expression finale est symétrique.
L'assertion est donc démontrée pour les fonctions à support compact;
le cas général en résulte par continuité, les fonctions lisses et à support compact étant denses.
\end{proof}

\begin{corollaire} 
Supposons $T$ assez régulier \textup(comme au lemme~\ref{troncnul}\textup).
L'opérateur $\tronc^\T$ s'étend en un projecteur autoadjoint sur $L^2(\XG)$.
\end{corollaire}

\begin{proof}
On vérifie à l'aide des propositions~\ref{rapdec} et~\ref{autoad} 
et du corollaire~\ref{invol}
que
\[
\langle(1-\tronc^\T)\varphi, \tronc^\T\varphi\rangle =0
\]
pour $\varphi$ lisse et à support compact sur $\XG$.
Il en résulte que
\[
\langle \varphi,\varphi\rangle = \langle \tronc^\T\varphi,\tronc^\T \varphi\rangle + 
\langle (1-\tronc^\T)\varphi, (1-\tronc^\T)\varphi\rangle
\]
ce qui implique que $\tronc^\T$ se prolonge en un opérateur continu involutif autoadjoint
dans $L^2$ puisque de telles fonctions sont évidemment denses.
\end{proof}

\chapter{Formes automorphes et produits scalaires}\label{ch5}

\section[Formes automorphes sur $\XP$]{\mathversion{bold}Formes automorphes sur $\XP$}

On rappelle qu'une forme automorphe sur
\[
\XG=\gA_G\G(F)\bs\Gadef 
\]
est une fonction lisse, $\K$\hyph finie et $\goth z(\frakg)$\hyph finie\,\footnote{$\goth z(\frakg)$ désigne le centre de l'algèbre enveloppante de $\frakg$.}
et qui est à croissance lente. On sait qu'une telle fonction est
alors automatiquement à croissance uniformément lente.
Plus généralement,
soit $P$ un sous-groupe parabolique de sous-groupe de Levi $\M$
et de radical unipotent $\N_P$ et considérons
\[
\XP=\gA_P P(F)\N_P(\adef)\bs\Gadef\ptf
\]
On appellera forme automorphe sur $\XP$ une fonction $\Phi$
sur $\XP$ qui est $\K$\hyph finie à droite et telle que pour tout
$k\in\K$ la fonction sur $\XM$ définie par
\[
m\mapsto\Phi(m\,k)
\]
soit automorphe.
On notera
\[
\autom(\XP)
\]
l'espace de formes automorphes sur $\XP$.
On trouvera dans \cite{MW}*{p.~37} une autre définition de la notion
de forme automorphe sur $\XP$ qui est démontrée être
équivalente à celle-ci.
Soit $\sigma$ une représentation automorphe de $\M$.
On notera
\[
\autom(\XP,\sigma)
\]
l'espace des formes automorphes sur $\XP$
telles que pour tout $x\in\Gadef$ la fonction
\[
m\mapsto\Phi(mx)\qquad\text{pour $m\in\M(\adef)$}
\]
soit une forme automorphe de l'espace
isotypique de $\sigma$ dans $L^2_{\mathrm{disc}}(\XM)$.
On dira, par abus de langage, que $\Phi$ est cuspidale sur $\XP$ si
\[
m\mapsto\Phi(mx)
\]
est cuspidale.
L'espace $\autom_{\mathrm{cusp}}(\XP)$, des formes cuspidales
(resp.~$\autom_{\mathrm{disc}}(\XP)$, des formes de carré intégrables),
est muni d'une structure d'espace
préhilbertien par le produit scalaire
\[
\langle \Phi,\Psi\rangle_P=\int_{\XP}\Phi(x)\widebar{\Psi(x)}\dd x=
\int_\K\int_{\XM}
\Phi(mk)\widebar{\Psi(mk)}\dd m\dd k\ptf
\]

\section{Opérateurs d'entrelacement et séries d'Eisenstein}

Soit $\Phi$
une fonction lisse sur $\XP$.
Pour
\[
\lambda\in\ga_P^*\otimes\CM
\]
on pose
\[
\Phi(x,\lambda)=
\Phi(x)\ee^{\langle \lambda+\demisom_P,\HP(x)\rangle}
\]
où $\HP(x)$ est la projection de $\HO(x)$ sur $\ga_P$
et $\demisom_P$ est la demi-somme des racines dans $\N_P$.
Considérons un sous-groupe parabolique $\Q$ associé à $P$
et
\[
s\in\weyl(\ga_P,\ga_\Q)\ptf
\]
On notera $\mathbf{s}$ l'opérateur défini par
\[
\mathbf{s}\Phi(x)=\Phi(w_s\moins\,x)
\]
\newnot{s@$\mathbf{s}$}{bfs}%
et on pose
\[
 \N_{s,P,\Q}=\N_\Q\cap w_s\N_P w_s\moins\bs\N_\Q\ptf
\]
Supposons $\Phi$ cuspidale sur $\XP$.
Pour $\lambda$ assez régulier dans la chambre associée à $P$ dans $\ga_P^*\otimes\CM$,
il existe une fonction $\Psi$ sur
\[
\XQ=\gA_\Q\Q(F)\N_\Q(\adef)\bs\Gadef
\]
telle que
\[
\Psi(x,s\lambda)=\int_{\N_{s,P,\Q}(\adef)}\mathbf{s}\Phi(n\,x,\lambda)\dd n
\]
l'intégrale étant absolument convergente.
On définit l'opérateur d'entrelacement
$\Mint_{\Q\rest P}(s,\lambda)$  par
\[
\qquad \Mint_{\Q\rest P}(s,\lambda)\Phi=\Psi\ptf
\]
\newnot{MQP(s,lambda)@$\Mint_{Q\string\rest{}P}(s,\string\lambda)$}{mintqp}%
En d'autres termes
\[
\Mint_{\Q\rest P}(s,\lambda)\Phi(x) =\ee^{-\langle s\lambda+\demisom_\Q,\HQ(x)\rangle}
\int_{\N_{s,P,\Q}(\adef)}\Phi(w_s\moins n\,x)
\ee^{\langle \lambda+\demisom_P,\HP(w_s\moins n\,x)\rangle}\dd n\ptf
\]
On posera
\[
\Mint_{P\rest \Q}(\lambda)=\Mint_{P\rest \Q}(1,\lambda)\ptf
\]

\begin{lemme}\label{Mintsimple}
\[
\Mint_{s(P)\rest P}(s,\lambda) =\ee^{\langle\lambda+\demisom_P,\TK-s\moins\TK\rangle}\mathbf{s}\ptf
\]
\end{lemme}

\begin{proof}
Lorsque
\[
\Q=w_sP\,w_s\moins=s(P)
\]
le groupe $\N_{s,P,\Q}$
est trivial et on a simplement
\[
\Mint_{s(P)\rest P}(s,\lambda)\Phi=\ee^{\langle \lambda+\demisom_P,\HP(w_s\moins)\rangle}\mathbf{s}\Phi
\ptf
\]
On conclut en observant que, compte tenu du lemme~\ref{Yorth}, on a pour tout $P$ semi-standard
\[
\langle\lambda+\demisom_P,\HP(w_s\moins)\rangle=\langle\lambda+\demisom_P,\HO(w_s\moins)\rangle
\]
et
\[
\HO(w_s\moins)=\TK-s\moins\TK\ptf\qedhere
\]
\end{proof}

Lorsque $P$ et $\Q$ sont standard, et $P$ fixé,
la donnée du couple $(s,\lambda)$ suffit à déterminer l'opérateur
$\Mint_{\Q\rest P}(s,\lambda)$
qui sera parfois noté simplement $\Mint(s,\lambda)$.
\newnot{M(s,lambda)@$\string\Mint(s,\lambda)$}{mintsl}%
Ce sera en particulier le cas dans 
\ref{eqfonct}(ii) ci-dessous.

Soit $P\subset\Q$ une paire de sous-groupes paraboliques standard
et soit $\Phi$ une forme automorphe cuspidale sur $\XP$.
On définit, pour
$\lambda\in\ga_P^*\otimes\CM$ assez régulier, une série d'Eisenstein sur $\Q$ en posant\,\footnote{Nous utilisons la notation $E^\Q$, suivant en cela Arthur et \cite{MW}
plutôt que $E_\Q$ utilisé par Langlands \cite{MS}*{Lecture~15} pour éviter les confusions avec
le terme constant le long de $\Q$ souvent noté ainsi.}
\[
E^\Q(x,\Phi,\lambda)=\sum_{\gamma\in P(F)\bs\Q(F)}
\Phi(\gamma\,x,\Lambda)\ptf
\]
La série d'Eisenstein sera notée simplement $E(x,\Phi,\lambda)$
lorsque $\Q=\G$.

Soit $\R$ un sous-groupe parabolique standard de $\G$;
on rappelle que l'on note $\Pi_\R E$, ou parfois $E_\R$, le terme constant de $E$ le long de $\R$
(\cf section~\ref{defpropannul}):
\[
\Pi_\R E(x,\Phi,\lambda)=
\int_{\N_\R(F)\bs\N_\R(\adef)}E(nx,\Phi,\lambda)\dd n\ptf
\]
On rappelle que l'on a introduit dans le lemme~\ref{weylgab} le sous-ensemble
$\weyl(\ga_P,\R)$ du groupe de Weyl.
On peut alors énoncer le théorème fondamental de Langlands.

\begin{theoreme}\label{eqfonct}
\begin{asparaenum}[(i)]
\item L'opérateur d'entrelacement $\Mint_{\Q\rest P}(s,\lambda)$ possède un prolongement méromorphe à tout $\ga_P^*\otimes\CM$.
Si $P$\textup, $\Q$ et $\R$ sont trois sous-groupes paraboliques \textup(semi-standard\textup)\textup,
et $s$ et $t$ sont deux éléments du groupe de Weyl tels que
\[
s\in\weyl(\ga_P,\ga_\Q)\Qquad{et}t\in\weyl(\ga_\Q,\ga_\R)
\]
alors on a l'équation fonctionnelle:
\begin{equation}
\Mint_{\R\rest \Q}(t,s\lambda)\Mint_{\Q\rest P}(s,\lambda)=\Mint_{\R\rest P}(ts,\lambda)\ptf\label{eq5.1}
\end{equation}
De plus\textup,
\begin{equation}
\Mint_{\Q\rest P}(s,-\widebar\lambda)^*=\Mint_{\Q\rest P}(s,\lambda)\moins\ptf\label{eq5.2}
\end{equation}
En particulier\textup, pour $\lambda$ imaginaire pur
l'opérateur d'entrelacement $\Mint_{P\rest \Q}(s,\lambda)$ est une isométrie.
\item
La série d'Eisenstein $E(x,\Phi,\lambda)$ converge absolument si $\RE(\lambda)>\demisomP$.
Elle admet un prolongement méromorphe à tout $\ga_P^*\otimes\CM$
et définit ainsi une forme automorphe qui est à croissance uniformément lente
lorsque le paramètre $\lambda$ reste dans un compact du domaine d'holomorphie.
Elle satisfait les équations fonctionnelles
\begin{equation}
E(x,\Mint(s,\lambda)\Phi,s\lambda)=E(x,\Phi,\lambda)
\qquad\text{pour $s\in\weyl(\ga_P)$}\ptf\label{eq5.3}
\end{equation}
Supposons $\Phi$ cuspidale sur $\XP\ptf$
On a\,\footnote{On observera que, compte tenu des équations fonctionnelles~\ref{eqfonct}\eqref{eq5.3}, le choix dans~\ref{eqfonct}\eqref{eq5.4}
du représentant $s$, dans la classe modulo $\weyl^\R$ définie par
un élément de $\weyl(\ga_P,\R)$, est indifférent.}
\begin{equation}
\Pi_\R E(x,\Phi,\lambda)=\sum_{s\in\weyl(\ga_P,\R)}E^\R(x,\Mint(s,\lambda)\Phi,s\lambda)\ptf\label{eq5.4}
\end{equation}
Dans le cas particulier où $P$ et $\R$ sont associés on a simplement
\begin{equation}
\Pi_\R E(x,\Phi,\lambda)=\sum_{s\in\weyl(\ga_P,\ga_\R)}(\Mint(s,\lambda)\Phi)(x,s\lambda)\ptf\label{eq5.5}
\end{equation}
\end{asparaenum}
\end{theoreme}

\begin{proof} 
Pour le prolongement analytique et les équations fonctionnelles on renvoie
à \cite{Eisen} (voir aussi \cite{MW}*{IV.1.10}).
Pour le calcul du terme constant on pourra consulter \cite{Bould} ou \cite{MW}*{II.1.7}.
\end{proof}

\section[La \GM-famille spectrale]{\mathversion{bold}La \GM-famille spectrale}\label{gmspec}

Soit $\M$ un sous-groupe de Levi et soient $P$ et $\Q$ dans $\Parab(\M)$.

\begin{lemme}\label{adjmur}
Supposons que $P$ et $\Q$ correspondent à des chambres
adjacentes dans $\ga_\M$.
Si $\tronc^\vee$ appartient au mur séparant les deux chambres on a
\[
\Mint_{\Q\rest P}(\lambda+\Lambda)=\Mint_{\Q\rest P}(\lambda)\ptf
\]
\end{lemme}

\begin{proof}
On rappelle que, par définition
\[
\Mint_{\Q\rest P}(\lambda)=\Mint_{\Q\rest P}(1,\lambda)
\]
et donc, dans le domaine de convergence,
\[
\Mint_{\Q\rest P}(\lambda+\Lambda)\Phi(x) =\ee^{-\langle \lambda+\Lambda+\demisom_\Q,\HQ(x)\rangle}
\int_{\N_{1,P,\Q}(\adef)}\Phi( n\,x)
\ee^{\langle \lambda+\Lambda+\demisom_P,\HP( n\,x)\rangle}\dd n\ptf
\]
Il suffit alors d'observer que si $\tronc^\vee$ appartient au mur séparant les deux chambres on a
\[
\langle\Lambda,\HQ(x)\rangle=\langle \Lambda,\HP (nx)\rangle
\]
pour tout $n\in\NO(\adef)$.
\end{proof}

\begin{corollaire}\label{GMspeca} 
Pour $P$ dans $\Parab(\M)$
et $\lambda$ donnés\textup, la famille de fonctions à valeurs opérateurs
\[
\MM(P,\lambda;\Lambda,\Q)= \Mint_{\Q\rest P}(\lambda)\moins\Mint_{\Q\rest P}(\lambda+\Lambda)
\]
indexée par $\Q\in\Parab(\M)$ est une \GM-famille.
\end{corollaire}

\begin{proof}
L'équation fonctionnelle~\ref{eqfonct} montre que
\[
\Mint_{\R\rest P}(\lambda+\Lambda)=\Mint_{\R\rest \Q}(\lambda+\Lambda)\Mint_{\Q\rest P}(\lambda+\Lambda)\ptf
\]
Donc,
\begin{equation}
\MM(P,\lambda;\Lambda,\R)=\Mint_{\Q\rest P}(\lambda)\moins\Mint_{\R\rest \Q}(\lambda)\moins\Mint_{\R\rest \Q}(\lambda+\Lambda)
\Mint_{\Q\rest P}(\lambda+\Lambda)\ptf\label{eq5.1a}
\end{equation}

Si nous supposons maintenant que $\Q$ et $\R$ sont associés à des chambres
adjacentes et si $\Lambda$ appartient au mur séparant les deux chambres on
sait d'après le lemme~\ref{adjmur} que
\[
\Mint_{\R\rest \Q}(\lambda+\Lambda)=\Mint_{\R\rest \Q}(\lambda)\ptf
\]
Dans ce cas~\eqref{eq5.1a} se récrit
\[
\MM(P,\lambda;\Lambda,\R)=\Mint_{\Q\rest P}(\lambda)\moins\Mint_{\Q\rest P}(\lambda+\Lambda)=\MM(P,\lambda;\Lambda,\Q)\ptf \qedhere
\]
\end{proof}


Rappelons que l'on a introduit (à la section~\ref{hown}) la famille orthogonale
\[
Y_s(\T)=s\moins T+\TK-s\moins\TK
\]
et que, si $\M$ est le sous-groupe de Levi d'un sous-groupe parabolique standard $P$,
on associe à tout $S\in\Parab(\M)$, le vecteur
\[
Y_S(T)\in\ga_\M
\]
qui est la projection de $Y_s(\T)$ sur $\ga_\M$ lorsque $s\in\weyl$ est tel que $sS$ est standard. Enfin on écrira parfois
\[
Y_s\quad\text{pour $Y_s(0)$}\Qquad{ainsi que}Y_S\quad\text{pour $Y_S(0)$}
\]
\newnot{Ys@$Y_s$}{Yszero}%
\newnot{YS@$Y_S$}{YSzero}%
et on observera que $Y_S(\TK)=\TK$ et est donc indépendant de $S$.

\begin{proposition}\label{GMspec}
La famille de fonctions méromorphes de $\lambda$ et $\Lambda$
définie pour
\[
S\in\Parab(\M)
\]
par
\begin{equation}
\MM(P,\T,\lambda;\Lambda,S)=\ee^{\langle\Lambda,Y_S(T)\rangle}\,\,
\Mint_{S\rest P}(\lambda)\moins\Mint_{S\rest P}(\lambda+\Lambda)\label{eq5.1b}
\end{equation}
est une \GM-famille. Considérons $\Q\in\FF(\M)$.
Les fonctions méromorphes
\begin{equation}
\MM_{\M}^\Q(P,\lambda;\Lambda)=
\sum_{S\in\Parab^{\Q}(\M)}\epsilon_S^{\Q}(\Lambda)\MM(P,\lambda;\Lambda,S)\label{eq5.2b}
\end{equation}
et
\begin{equation}
\MM_{\M}^\Q(P,\T,\lambda;\Lambda)=
\sum_{S\in\Parab^{\Q}(\M)}\epsilon_S^{\Q}(\Lambda)\MM(P,\T,\lambda;\Lambda,S)\ptf
\label{eq5.3b}
\end{equation}
sont lisses pour les valeurs imaginaires pures des paramètres $\lambda$ et $\Lambda$.
De plus\textup:
\begin{equation}
\MM_{\M}^\Q(P,\TK,\lambda;0)=\MM_{\M}^\Q(P,\lambda;0)\ptf\label{eq5.4b}
\end{equation}
\end{proposition}

\begin{proof}
Le fait que la formule~\eqref{eq5.1b} définisse une \GM-famille résulte des observations suivantes.
Tout d'abord la famille de fonctions
\[
\cbf(\T;\Lambda,S)= \ee^{\langle \Lambda,Y_S(T)\rangle}
\]
définit une \GM-famille d'après le lemme~\ref{gmorth} car
$Y_S(\T)$ est une famille $\M$\hyph orthogonale. Maintenant la famille de fonctions
\[
\MM(P,\lambda;\Lambda,S)=\Mint_{S\rest P}(\lambda)\moins\Mint_{S\rest P}(\lambda+\Lambda)
\]
définit aussi une \GM-famille d'après le corollaire~\ref{GMspeca}. Donc
\[
\MM(P,\T,\lambda;\Lambda,S)=\cbf(\T;\Lambda,S)\MM(P,\lambda;\Lambda,S)
\]
est le produit de deux \GM-familles: c'est une \GM-famille.
La lissité de~\eqref{eq5.2b} et~\eqref{eq5.3b} pour les valeurs imaginaires pures des paramètres
résulte alors du lemme \ref{lisse}.  L'assertion~\eqref{eq5.4b} provient de ce que, pour tout $s$, on a
$Y_s(\TK)=\TK$ et donc
\[
\MM_{\M}^\Q(P,\TK,\lambda;\Lambda)=\ee^{\langle \Lambda,\TK\rangle}
\MM_{\M}^\Q(P,\lambda;\Lambda)\ptf\qedhere
\]
\end{proof}

\begin{lemme}\label{sta}
Soient $P$\textup, $\Q$ et $\R$ trois sous-groupes
paraboliques. Soient $s$ et $t$ deux éléments du groupe de Weyl
avec $s(\ga_P)=\ga_\R$ et $t(\ga_\Q)=\ga_\R$ \cad
\[
s\in\weyl(\ga_P,\ga_\R),\qquad t\in\weyl(\ga_\Q,\ga_\R)\ptf
\]
Posons
\[
\us=t\moins s,\qquad S=t\moins\R\Qquad{et}\Lambda=\us\lambda-\mu\ptf
\]
Alors,
\begin{equation}
\Mint_{\R\rest \Q}(t,\mu)\moins\Mint_{\R\rest P}(s,\lambda)
=\ee^{\langle\Lambda,Y_t\rangle}
\Mint_{S\rest \Q}(\mu)\moins
\Mint_{S\rest \Q}(\mu+\Lambda)\Mint_{\Q\rest P}(\us,\lambda)\ptf\label{eq5.1c}
\end{equation}
De plus,
\begin{equation}
\Mint_{\Q\rest P}(\us,\lambda)=
\ee^{\langle\lambda+\demisom_{P},Y_{\us}\rangle}
\Mint_{\Q\rest \us P}(\mu+\Lambda)\mathbf{\us}\ptf\label{eq5.2c}
\end{equation}
En particulier, si $P$, $\Q$ et $\R$ sont standard et si $s=t$ on a
$P=\Q$ et
\begin{equation}
\Mint(s,\mu)\moins\Mint(s,\lambda)=
\ee^{\langle\Lambda,Y_s\rangle}\Mint_{S\rest \Q}(\mu)\moins\Mint_{S\rest \Q}(\lambda)\ptf\label{eq5.3c}
\end{equation}
\end{lemme}

\begin{proof}
L'équation fonctionnelle pour les opérateurs d'entrelacement~\ref{eqfonct}(\ref{eq5.1})\linebreak montre que
\[
\Mint_{\R\rest \Q}(t,\mu)=\Mint_{\R\rest S}(t,\mu)\Mint_{S\rest \Q}(1,\mu)\ptf
\]
Mais d'après le lemme~\ref{Mintsimple} on sait que
\[
\Mint_{\R\rest S}(t,\mu)=\Mint_{tS\rest S}(t,\mu)
=\ee^{\langle\mu+\demisom_S,\TK-t\moins\TK\rangle}\mathbf{t}=\ee^{\langle\mu+\demisom_S,Y_t\rangle}\mathbf{t}\ptf
\]
De même, en posant $\us=t\moins s$, on a
\[
\Mint_{\R\rest P}(s,\lambda)=
\Mint_{\R\rest S}(t,\us\lambda)\Mint_{S\rest P}(\us,\lambda)
\]
soit encore
\[
\Mint_{\R\rest P}(s,\lambda)=\ee^{\langle\us\lambda+\demisom_S,Y_t\rangle}\mathbf{t}
\Mint_{S\rest P}(\us,\lambda)\ptf
\]
On a donc
\[
\Mint_{\R\rest \Q}(t,\mu)\moins\Mint_{\R\rest P}(s,\lambda)=\ee^{\langle\Lambda,Y_t\rangle}
\Mint_{S\rest \Q}(\mu)\moins\Mint_{S\rest P}(\us,\lambda)\ptf
\]
Par ailleurs en posant
\[
\Lambda=\us\lambda-\mu
\]
l'équation fonctionnelle fournit
\[
\Mint_{S\rest P}(\us,\lambda)=
\Mint_{S\rest \Q}(\mu+\Lambda)\Mint_{\Q\rest P}(\us,\lambda)
\]
et on voit alors que
\[
\Mint_{\R\rest \Q}(t,\mu)\moins\Mint_{\R\rest P}(s,\lambda)
\]
est égal à
\[
\ee^{\langle\Lambda,Y_t\rangle}\Mint_{S\rest \Q}(\mu)\moins
\Mint_{S\rest \Q}(\mu+\Lambda)\Mint_{\Q\rest P}(\us,\lambda)\ptf
\]
Ceci établit la première assertion. Maintenant
\[
\Mint_{\Q\rest P}(\us,\lambda)=
\Mint_{\Q\rest \us P}(\mu+\Lambda)\Mint_{\us P\rest P}(\us,\lambda)\ptf
\]
Mais, d'après le lemme~\ref{Mintsimple}
\[
\Mint_{\us P\rest P}(\us,\lambda)=\Mint_{\us P\rest P}(\us,\lambda)=
\ee^{\langle\lambda+\demisom_{P},Y_{\us}\rangle}\mathbf{\us}
\]
on a donc
\[
\Mint_{\Q\rest P}(\us,\lambda)=
\ee^{\langle\lambda+\demisom_{P},Y_{\us}\rangle}
\Mint_{\Q\rest \us P}(\mu+\Lambda)\mathbf{\us}\ptf\qedhere
\]
\end{proof}

Soient $P$ et $\Q$ deux sous groupes paraboliques standard.
On introduit la fonction méromorphe en
$\lambda$ et $\mu$, à valeurs opérateurs
\[
\gmomega_{\Q\rest P}^\T(\lambda,\mu)=\sum_\R\sum_{s\in\weyl(\ga_P,\ga_\R)}\sum_{t\in\weyl(\ga_\Q,\ga_\R)}
\ee^{\langle s\lambda-t\mu,T\rangle}\epsilon_\R^\G(s\lambda-t\mu)
\Mint(t,\mu)\moins\Mint(s,\lambda)
\]
où $\R$ parcourt l'ensemble des sous-groupes paraboliques standard associés à $P$.
On observera que cette expression n'est non nulle que si $P$ et $\Q$ sont associés.

\begin{lemme}\label{OmegaT}
Notons $\M$ le sous-groupe de Levi de $P$. On a\textup:
\[
\gmomega_{\Q\rest P}^\T(\lambda,\mu)=\sum_{\us\in\weyl(\ga_P,\ga_\Q)}\MM_{\M}^\G(\Q,\T,\mu;\Lambda)
\Mint_{\Q\rest P}(\us,\lambda)\ptf
\]
Cette fonction
est lisse pour les valeurs imaginaires pures des paramètres $\lambda$ et $\mu$.
\end{lemme}

\begin{proof}
Compte tenu du lemme~\ref{sta} on voit que
$\gmomega_{\Q\rest P}^\T(\lambda,\mu)$
est égal à
\[
\sum_{\us\in\weyl(\ga_P,\ga_\Q)}
\sum_{S\in\Parab(\M)}
\ee^{\langle\Lambda,Y_S(\T)\rangle}\epsilon_S^G(\Lambda)\Mint_{S\rest \Q}(\mu)\moins
\Mint_{S\rest \Q}(\mu+\Lambda)\Mint_{\Q\rest P}(\us,\lambda)\ptf
\]
On invoque alors la proposition~\ref{GMspec}.
\end{proof}

\section{Séries d'Eisenstein et troncature}

Le calcul du produit scalaire de deux séries d'Eisenstein tronquées
\[
\int_{\XG}\tronc^{T}E(x,\Phi,\lambda)\widebar{\tronc^{T}E(x,\Psi,-\widebar\mu)}\dd x
\]
est un résultat classique, dû à Langlands \cite{Eisen}, que nous rappelons au théorème~\ref{prodscal}. Compte tenu de l'autoadjonction (proposition~\ref{autoad}) et de l'involutivité (corollaire~\ref{invol}) de $\tronc^{T}$ on a
\[
\int_{\XG}\tronc^{T}E(x,\Phi,\lambda)\widebar{\tronc^{T}E(x,\Psi,-\widebar\mu)}\dd x=
\int_{\XG}\tronc^{T}E(x,\Phi,\lambda)\widebar{E(x,\Psi,-\widebar\mu)}\dd x\ptf
\]
En utilisant la deuxième expression nous allons donner, dans le cas cuspidal, une preuve\,\footnote{C'est la preuve donnée dans \cite{MS}*{Lecture 13}.} du théorème~\ref{prodscal} beaucoup plus simple que celle donnée par Arthur dans \cite{ATFII}*{p.~113-119}).
Pour le passage du cas cuspidal au cas général on renvoie à la littérature.

\begin{proposition}\label{tronceis}
Soit $\Phi$ cuspidale sur
\[
\XQ=\gA_\Q\Q(F)\N_\Q(\adef)\bs\Gadef\ptf
\]
Pour $\lambda$ dans le domaine de convergence de la série d'Eisenstein\textup,
la série d'Eisenstein tronquée
\[
\tronc^{T}E(x,\Phi,\lambda)
\]
est donnée par l'expression
\[
\sum_S
\sum_{s\in\weyl(\ga_\Q,\ga_S)}\,\sum_{\xi\in S(F)\bs\G(F)}(-1)^{a(s)}
\phi_{\M,s}\bigl( s\moins (\xixT)\bigr) (\Mint(s,\lambda)\Phi)(\xix,s\lambda)
\]
où $S$ parcourt l'ensemble des sous-groupes paraboliques standard
associés à $\Q$.
\end{proposition}

\begin{proof}
Par définition de $\tronc^\T$
on a
\[
\tronc^{T}E(x,\Phi,\lambda)=\sum_{\{P\mid\PO\subset P\}}
(-1)^{\ga_P-\ga_\G}\sum_{\xi\in P(F)\bs\G(F)}
\htau_P^\G(\xixT )\Pi_P E(\xix,\Phi,\lambda)
\]
qui par l'équation~\eqref{eq5.3} du théorème~\ref{eqfonct} est égal à
\[
\sum_{\{P\mid\PO\subset P\}} (-1)^{\ga_P-\ga_\G}
\sum_{\xi\in P(F)\bs\G(F)}
\htau_P^\G(\xixT )\sum_
{s\in\weyl(\ga_\Q,P)}E^P(\xix,\Mint(s,\lambda)\Phi,s\lambda)
\]
soit encore
\[
\sum_{\{P\mid\PO\subset P\}}(-1)^{\ga_P-\ga_\G}
\sum_{s\in\weyl(\ga_\Q,P)}\sum_{\xi\in P(F)\bs\G(F)}
\htau_P^\G( \xixT)E^P(\xix,\Mint(s,\lambda)\Phi,s\lambda)\ptf
\]
Si on note $S$ le sous-groupe parabolique standard
avec
\[
\ga_S=s(\ga_\Q)
\]
et puisque nous sommes dans le domaine de convergence,
l'expression
\[
\sum_{\xi\in P(F)\bs\G(F)}\htau_P^\G( \xixT)E^P(\xix,\Mint(s,\lambda)\Phi,s\lambda)
\]
est égale à
\[
\sum_{\xi\in S(F)\bs\G(F)}
\htau_P^\G( \xixT)(\Mint(s,\lambda)\Phi)(\xix,s\lambda)\ptf
\]
Par ailleurs, d'après le lemme~\ref{preparb},
\[
\sum_{\{P\mid s\moins(P)\in\FF_s(\M)\}}(-1)^{a_P-a_\G}\,\htau_P(H)
=(-1)^{a(s)}\phi_{\M,s} (s\moins H)\ptf
\]
Donc, compte tenu du lemme~\ref{somb}, on obtient que $\tronc^{T}E(x,\Phi,\lambda)$ est égal à
\[
\sum_S\sum_{s\in\weyl(\ga_\Q,\ga_S)}\,\sum_{\xi\in S(F)\bs\G(F)}
(-1)^{a(s)} \phi_{\M,s}\bigl( s\moins (\xixT)\bigr)(\Mint(s,\lambda)\Phi)(\xix,s\lambda)
\]
où $S$ parcourt l'ensemble des sous-groupes paraboliques standard associés à $\Q$.
\end{proof}

\begin{theoreme}\label{prodscal}
On suppose donnés $\T\in\gao$ ainsi que
\[
\lambda\in\ga_\Q^*\otimes\CM\Quad{et}\mu\in\ga_\R^*\otimes\CM
\]
qui coïncident sur $\ga_\G$.
\begin{enumerate}[(i)]
\item
Lorsque $\Phi$ et $\Psi$ sont cuspidales sur $\Q$ et $\R$ alors
on a l'égalité de fonctions méromorphes\textup:
\[
\int_{\XG}\tronc^{T}E(x,\Phi,\lambda)
\widebar{E(x,\Psi,-\widebar\mu)}\dd x=\langle\gmomega_{\R\rest \Q}^\T(\lambda,\mu)\Phi,\Psi\rangle\ptf\label{eq5.1d}
\]
En particulier\textup, ce produit scalaire est nul si $\Q$ et $\R$ ne sont pas associés.\,\footnote{On observera que $\T$ n'intervient que via sa projection $\T^\G$ sur $\gao^\G$.}
\item
Dans le cas général \textup(où $\Phi$ et $\Psi$ ne sont plus nécessairement cuspidales\textup)
pour $\lambda$ et $\mu$ dans des compacts fixés de $\ima\ga_\Q^*$ et $\ima\ga_\R^*$\textup,
il existe $A>0$ tel que
\begin{equation}
\biggl\lvert \int_{\XGa}\tronc^{T}E(x,\Phi,\lambda)\widebar{E(x,\Psi,-\widebar\mu)}\dd x-
\langle \gmomega_{\R\rest \Q}^\T(\lambda,\mu)\Phi,\Psi\rangle 
\biggr\rvert \ll
\ee^{-A\,\dPO(T)}\ptf\label{eq5.2d}
\end{equation}
\end{enumerate}
\end{theoreme}

\begin{proof}
Considérons $\lambda$ dans le domaine de convergence de la série d'Eisenstein
$E(x,\Phi,\lambda)$ et soit $\mu$ une valeur non singulière pour $E(x,\Psi,-\widebar\mu)$.
D'après la proposition~\ref{tronceis} l'intégrale est égale à
\[
\sum_S\sum_{s\in\weyl(\ga_\Q,\ga_S)}\int_{\XSG}
(-1)^{a(s)} \phi_{\M,s}(s\moins Z(x,T) )\, A(x,s)
\]
avec
\[
\XSG= \gA_\G S(F)\N_S(\adef)\bs\Gadef
\]
et
\[
A(x,s)=(\Mint(s,\lambda)\Phi)(x,s\lambda)
\,\Pi_S \widebar{E(x,\Psi,-\widebar\mu)}\dd x\ptf
\]
Il suffit alors d'invoquer l'équation~\eqref{eq5.5} du théorème~\ref{eqfonct}
pour obtenir
\[
A(x,s)=\sum_{t\in\weyl(\ga_\R,\ga_S)}
\ee^{\langle s\lambda-t\mu+2\demisom_S,H_S(x)\rangle}
(\Mint(s,\lambda)\Phi)(x)
\widebar{(\Mint(t,-\widebar\mu)\Psi)(x)}\ptf
\]
Maintenant on remarque que la fonction
\[
x\mapsto(\Mint(s,\lambda)\Phi)(x)
\widebar{(\Mint(t,-\widebar\mu)\Psi)(x)}
\]
est invariante par $\gA_S S(F)\N_S(\adef)$
et on pose
\[
\langle \Mint(s,\lambda)\Phi,\Mint(t,-\widebar\mu)\Psi\rangle=\int_{\XS}
(\Mint(s,\lambda)\Phi)(x)
\widebar{(\Mint(t,-\widebar\mu)\Psi)(x)}\dd x
\]
avec
\[
\XS=\gA_S S(F)\N_S(\adef)\bs\Gadef\ptf
\]
En tenant compte de ce que
\[
\Mint(t,-\widebar\mu)^*=\Mint(t,\mu)\moins
\]
on obtient
\[
\langle \Mint(s,\lambda)\Phi,\Mint(t,-\widebar\mu)\Psi\rangle =\langle\Mint(t,\mu)\moins\Mint(s,\lambda)\Phi,\Psi\rangle
\ptf
\]
Par ailleurs
\[
x\mapsto \HH_S(x)
\]
est invariante sous le noyau de $S(\adef)\to\ga_S$
et on a une fibration
\[
\ga_S^\G\to\XSG\to \XS
\]
où
\[
\XSG=\gA_\G S(F)\N_S(\adef)\bs\Gadef\ptf
\]
On en déduit que, au moins formellement, le produit scalaire
est égal à la somme sur $s$, $t$, et $S$ de
\begin{equation}
\langle \Mint(t,\mu)\moins\Mint(s,\lambda)\Phi,\Psi\rangle
\int_{\ga_S^\G}\ee^{\langle s\lambda-t\mu,H\rangle}(-1)^{a(s)}\phi_{\M,s}(s\moins(H-T) )\dd H\ptf
\label{eq5.3d}
\end{equation}
La convergence est assurée d'après le lemme~\ref{intphis} si, pour $\mu$ fixé,
$\lambda$ est assez régulier et l'expression~\eqref{eq5.3d} ci-dessus est alors égale à
\[
\langle \Mint(t,\mu)\moins\Mint(s,\lambda)\Phi,\Psi\rangle
\ee^{\langle s\lambda-t\mu,T_S^\G\rangle}\epsilon_S^G(s\lambda-t\mu)
\]
où $T_S^\G$ est la projection de $\T$ sur $\ga_S^\G$.
De plus, comme
les formes linéaires $s\lambda$ et $t\mu$
sont triviales sur $\ga^S\oplus\ga_\G$ on a
\[
\langle s\lambda-t\mu, T_S^\G\rangle= \langle s\lambda-t\mu,T\rangle\ptf
\]
On a ainsi établi~\eqref{eq5.1d} pour $\lambda$ assez régulier.
Maintenant les deux membres de l'équation~\eqref{eq5.1d}sont des fonctions méromorphes en $\lambda$ et $\mu$:
c'est clair pour
\[
\langle\gmomega_{\R\rest \Q}^\T(\lambda,\mu)\Phi,\Psi\rangle
\]
compte tenu du théorème~\ref{eqfonct}; par ailleurs
pour $\lambda$ et $\mu$ dans des compacts de l'ouvert où les séries d'Eisenstein
sont holomorphes ces séries définissent des fonctions à croissance uniformément lente
et donc l'intégrale définissant le produit scalaire
est uniformément convergente et définit une fonction holomorphe sur cet ouvert.
L'égalité~\eqref{eq5.1d} est donc encore vraie pour tout $\lambda$ et $\mu$
en tant qu'égalité entre fonctions méromorphes.
Le passage du cas cuspidal au cas général est dû à Arthur \cite{AeisIII}*{Corollaire~9.2}.
Nous renvoyons à l'article \cite{AeisIII} pour une preuve.
\end{proof}

Soit $a=\ee^{H_\G}$ avec
$H_\G\in\ga_\G$.
On notera $\XGa$ le sous-ensemble
du quotient
\[
\G(F)\bs\Gadef
\]
formé des $x$ tels que
\[
\HG(x)=\HG(a)=H_\G\ptf
\]
On observera que l'application naturelle
$\XGa\to\XG$
est une bijection.

\begin{theoreme}\label{prodscalvar} 
Soient\textup, comme ci-dessus
$\lambda\in\ga_\Q^*\otimes\CM$ et $\mu\in\ga_\R^*\otimes\CM$\textup;
mais on ne suppose plus nécessairement
qu'ils coïncident sur $\ga_\G$.
\begin{enumerate}[(i)]
\item
Lorsque $\Phi$ et $\Psi$ sont cuspidales sur $\Q$ et $\R$ alors\textup, si $a=\ee^{H_\G}$
on a l'égalité de fonctions méromorphes\textup:
\[
\int_{\XGa}\tronc^{T}E(x,\Phi,\lambda)\widebar{E(x,\Psi,-\widebar\mu)}\dd x=
\langle \gmomega_{\R\rest \Q}^{H_\G+\T^\G}(\lambda,\mu)\Phi,\Psi\rangle
\]
où
$T^\G$ est la projection de $T$ sur $\gao^\G$.
\item
Dans le cas général\textup, pour $\lambda$ et $\mu$ dans des compacts fixés de $\ima\ga_\Q^*$ et $\ima\ga_\R^*$\textup,
il existe $A>0$ tel que
\[
\biggl\lvert \int_{\XGa}\tronc^{T}E(x,\Phi,\lambda)\widebar{E(x,\Psi,-\widebar\mu)}\dd x-
\langle\gmomega_{\R\rest \Q}^{H_\G+\T^\G}(\lambda,\mu)\Phi,\Psi\rangle
\biggr\rvert \ll \ee^{-A\,\dPO(T)}\ptf
\]
\end{enumerate}
\end{theoreme}

\begin{proof} 
La preuve est une variante de la preuve du théorème~\ref{prodscal}. La seule différence
est dans le résultat de l'intégration sur $\ga_S^G$.
On a utilisé que cette intégration permet un changement de variable
\[
H_S^\G\mapsto H_S^\G+T_S^\G
\]
où $T_S^\G$ est la projection de $T$ sur
$\ga_S^\G$. On peut encore ici remplacer $\T_S^\G$ par $\T^\G$ car
les formes linéaires $s\lambda$ et $t\mu$ sont triviales sur $\ga^S$
mais, comme on ne suppose plus que $\lambda$ et $\mu$ coïncident sur
$\ga_\G$ on ne peut pas remplacer $\T^\G$ par $\T$. Enfin, la présence de $H_\G$
provient de ce qu'on intègre sur l'espace des $x$ avec $\HG(x)=\HG(a)=H_\G$.
\end{proof}

\chapter{Le noyau intégral}\label{ch6}

\section{Les opérateurs en question}\label{noyint}

L'espace tordu $\tGadef$ agit sur 
l'espace homogène
\[
\XG=\gA_G\G(F)\bs\Gadef
\]
via l'action naturelle de $\tG(F)\times\tGadef$ sur $\XG$
suivant les conventions de la section~\ref{reptor}. Rappelons en la définition.
Considérons un point $\dot x\in\XG$ et $y\in\tGadef$.
Choisissons un représentant $x\in\Gadef$ de $\dot x$
et un élément $\delta\in\tG(F)$. Alors,
\[
\delta\moins\,x\,y
\]
définit un élément de $\Gadef$. Il est immédiat de voir
que la classe dans $\XG$ de cet élément est indépendante des choix de $x$ et de $\delta$;
nous la noterons
\[
\dot x*y\ptf
\]
La représentation régulière gauche $\reg$ de $\Gadef $ dans $L^2(\XG)$ 
admet un prolongement naturel en une représentation $\treg$
de $\tGadef$ définie par
\[
\treg(y)\vf(\dot x)=\vf(\dot x*y)
\]
pour $\vf\in L^2(\XG)$ et $y\in\tGadef$.
Nous utiliserons,
comme dans \cite{KS}, un objet un peu plus général:
on considère de plus un caractère unitaire $\omega$ de $\Gadef $
trivial sur $\gA_\G\G(F)$ et l'opérateur $\treg(y,\omega)$ défini par
\[
(\treg(y,\omega)\vf)(\dot x)=(\omega\vf)(\dot x*y)=\omega(\delta\moins x\,y)\,\vf(\delta\moins x\,y)\ptf
\]
Si on pose $y=g\delto$ avec $g\in\Gadef$ et où $\delto$ est l'élément choisi à la section~\ref{tordalg} dans $\tG(F)$
préservant $\PO$,
on a
\[
\treg(y,\omega)=\treg(g\delto,\omega)=A(\omega)\circ B(\theto)\circ\reg(g)
\]
où $\theto=\Ad(\delto)$ et où les opérateurs $A(\omega)$, $B(\theto)$ et $\reg(g)$ sont définis par
\[
(A(\omega)\vf)(\dot x)=\omega(\dot x)\,\vf(\dot x),\quad
(B(\theto)\vf)(\dot x)=\vf\bigl(\theta_0\moins(\dot x)\bigr)\Quad{et}
(\reg(g)\vf)(\dot x)=\vf(\dot x*g)\ptf
\]
On a ainsi défini une représentation
unitaire de $(\tGadef,\omega)$ au sens de la section~\ref{reptor}.
En effet, pour $x,z\in\Gadef $ et $y\in\tGadef$ on a
\[
\treg(x\,y\,z,\omega)=\reg(x)\treg(y,\omega)(\reg\otimes\omega)(z)\ptf
\]
Par intégration contre une fonction $\fff\in\ctyc\bigl(\tGadef\bigr)$
on définit l'opérateur
\[
\treg(\fff,\omega)=\int_{\tGadef} \fff(y)\treg(y,\omega)\dd y\ptf
\]
\newnot{rho(f,omega)@$\protect\treg(\fff,\omega)$}{tregfom}%
En d'autres termes on a
\[
(\treg(\fff,\omega)\vf)(\dot x)=\int_{\tG(\adef)}\,\fff(y)(\omega\vf)(\dot x*y)\dd y\ptf
\]
Si on pose $y=g\delto$ et
$h(g)=\fff(g\delto)$
on aura
\[
\treg(\fff,\omega)=A(\omega)B(\theto)\reg(h)\ptf
\]
On aurait pu utiliser, comme on le fera parfois dans un cadre plus général (voir ci-dessous), la notation
\[
\reg(\delta,y,\omega)\quad\text{au lieu de $\treg(y,\omega)$} \Qquad{et}
\reg(\delta,\fff,\omega)\quad\text{au lieu de $\treg(\fff,\omega)$}\ptf
\]
Mais comme ici l'opérateur $\reg(\delta,\fff,\omega)$ est indépendant de $\delta\in\tG(F)$ cette notation
est inutilement lourde.

Plus généralement,
soit $P$ un sous-groupe parabolique
et soit $\delta\in\tG(F)$.
Notons $\Q$ le sous-groupe parabolique
obtenu par conjugaison par $\delta$:
$\Q=\delta P\delta\moins=\theta(P)$
où $\theta=\tAd(\delta)$.
Considérons un point $\dot x\in\XQ$ et $y\in\tGadef$.
Choisissons un représentant $x\in\Gadef$ de $\dot x$. Alors,
\[
\delta\moins\,x\,y
\]
définit un élément de $\Gadef $ dont la classe dans $\XP$
est indépendante du choix de $x$.
On considère de plus un caractère unitaire $\omega$ de $\Gadef $
trivial sur $\gA_\G\G(F)$. Un élément $y\in\tGadef$ définit alors un opérateur
noté $\reg(\delta,y,\omega)$ entre l'espace des fonctions sur $\XP$
et l'espace des fonctions sur $\XQ$:
\[
\reg(\delta,y,\omega)\colon\Phi\mapsto \Psi
\]
avec
\[
\Psi(x)=\Phi(\delta\moins\,x\,y)\omega(\delta\moins x\,y)
\]
(\cf section~\ref{reptor}).
Plus généralement, considérons une représentation automorphe $\sigma$ de $\M$ et
soit $\Phi\in\autom(\XP,\sigma)$.
Pour $\mu\in\ga_P^*\otimes\CM$ on définit un opérateur
\[
 \reg_{P,\sigma,\mu}(\delta,y,\omega)
\]
entre l'espace des fonctions de carré intégrable engendré par
$\autom(\XP,\sigma)$ et celui engendré par
$\autom(\XQ,\tau)$ où
$\tau=\sigma\circ\theta\moins$ en posant
\[
\Psi=\reg_{P,\sigma,\mu}(\delta,y,\omega)\Phi
\]
avec
\[
\Psi(x)=\ee^{-\langle \theta(\mu+\demisom_P),\HQ(x)\rangle}
(\omega\Phi)(\delta\moins\,x\,y)\ee^{\langle\mu+\demisom_P,\HP(\delta\moins\,x\,y)\rangle}
\ptf
\]
On rappelle que
$\HP(x)$ est l'image dans $\ga_P$ de $p\in P(\adef)$:
$\HP(x)=\HP(p)$
si $x=p\,k$ est une décomposition d'Iwasawa de $x$.
Enfin, $\demisomP$ (resp.~$\demisomQ$)
est la demi-somme des racines dans $\N_P$ (resp.~$\N_\Q$).
Par intégration contre une fonction $\fff\in\ctyc\bigl(\tGadef\bigr)$
on définit l'opérateur
\[
\reg_{P,\sigma,\mu}(\delta,\fff,\omega)=\int_{\tGadef} \fff(y)\reg_{P,\sigma,\mu}(\delta,y,\omega)\dd y\ptf
\]
\newnot{rhoPsigmamu(delta,f,omega)@$\reg_{P,\sigma,\mu}(\delta,\fff,\omega)$}{regmusigma}%
On posera
\[
\treg_{P,\sigma,\mu}(\fff,\omega)\coloneqq\reg_{P,\sigma,\mu}(\delto,\fff,\omega)\ptf
\]
\newnot{rhoPsigmamu(f,omega)@$\treg_{P,\sigma,\mu}(\fff,\omega)$}{tregmusigma}%

\begin{lemme}\label{asu}
Supposons que $\delta=w_\us\delto$ avec $\us\in\weyl$.
Si $P$ est standard\textup, on a
\[
\reg_{P,\sigma,\mu}(\delta,\fff,\omega)=
\ee^{\langle\theto(\mu+\demisom_P),\HO(w_\us\moins)\rangle}
\mathbf{\us}\treg_{P,\sigma,\mu}(\fff,\omega)
\ptf
\]
\end{lemme}

\begin{proof}
On observe que, si on pose $\Q=\theta(P)$
et $\Q_0=\theto(P)$,
\[
(\reg_{P,\sigma,\mu}(\delta,\fff,\omega)\Phi)(x)
\]
est égal à
\[
\ee^{\langle \theto(\mu+\demisom_P),\HH_{\Q_0}(w_\us\moins x)\rangle
-\langle\theta(\mu+\demisom_P),\HQ(x)\rangle}
(\mathbf{\us}\treg_{P,\sigma,\mu}(\fff,\omega)\Phi)(x)\ptf
\]
Maintenant
\[
\langle\theto(\mu+\demisom_P),\HH_{\Q_0}(w_\us\moins x)\rangle
-\langle\theta(\mu+\demisom_P),\HQ(x)\rangle=\langle\theto(\mu+\demisom_P),\HH_{\Q_0}(w_\us\moins )\rangle
\]
et on conclut en observant que $\Q_0$ est standard.
\end{proof}

\section{Le noyau de la formule des traces}

L'opérateur $\treg(\fff,\omega)$ est représenté par un noyau intégral sur $\XG$:
\[
\treg(\fff,\omega)\vf(x)=\int_{\XG}\,K_{\tG}(\fff,\omega;x,y)\vf(y)\dd y
\]
avec
\[
K_{\tG}(\fff,\omega;x,y)=\sum_{\delta\in\tG(F)}\omega(y)\fff^1(x\moins\delta\,y)
\]
où  
\[
\fff^1(x)=\int_{z\in\gA_\G}\fff(z\,x)\dd z\ptf
\]
\newnot{f1@$\fff^1$}{fffun}%
Comme $\fff$ est à support compact la fonction
\[
y\mapsto K_{\tG}(\fff,\omega;x,y)
\]
est à support compact sur $\XG$ pour $x$ fixé.
Le noyau $K_{\tG}(\fff,\omega;x,y)$ sera noté $K_{\tG}(\fff;x,y)$ si $\omega=1$ voire même
simplement $K_{\tG}(x,y)$ 
\newnot{KG(x,y)@$K_{\protect\tG}(x,y)$}{ktgxy}%
si aucune confusion n'est à craindre,
la fonction $\fff$ et le caractère $\omega$ étant fixés.

\begin{lemme}\label{premaj}
Il existe des constantes $c(f)$ et
$N$ telles que pour tout $x$ et tout~$y$
\[
\lvert K_{\tG}(x,y)\rvert \le c(f)\lvert x\rvert^N\,\lvert y\rvert^N\ptf
\]
\end{lemme}

\begin{proof}
On observe que, pour une certaine constante $c$ on a
\[
\lvert \xi\rvert \le c \lvert x\rvert\,\lvert y\rvert \,\lvert x\moins \xi\,y\rvert
\]
et donc, pour une certaine constante $c'$ dépendant du support de $\fff$,
\[
\lvert \xi\rvert \le c' \lvert x\rvert\,\lvert y\rvert
\]
si $x\moins \xi\,y$ appartient au support de $\fff$ (qui est compact).
L'assertion résulte alors du lemme~\ref{proj}.
\end{proof}

\section{Factorisation de Dixmier-Malliavin}

\begin{theoreme}\label{dixmal}
Toute fonction $\fff\in\ctyc\bigl(\tGadef\bigr)$
est une somme finie de produits de convolution\textup:
\[
\fff=\sum_i g_i*h_i^*
\]
avec $g_i\in\ctyc\bigl(\tGadef\bigr)$ et $h_i\in\ctyc(\Gadef)$
qui peuvent de plus être choisies $\K$\hyph finies à gauche\,\footnote{Les fonctions peuvent probablement être choisies $\K$\hyph finies à droite et à gauche.
Mais cela nécessiterait un raffinement du théorème de Dixmier-Malliavin que nous n'avons pas tenté d'établir.
On obtient aussi une factorisation en fonctions $\K$\hyph finies à droite et à gauche si on utilise, comme Arthur,
la technique de Duflo-Labesse;
toutefois cela impose de se contenter d'un ordre fini de différentiabilité, ce qui est suffisant pour les
applications en vue, mais alourdit légèrement les énoncés.
Nous ne l'emploierons pas.}
si $\fff$ est $\K$\hyph finie à droite et à gauche.
\end{theoreme}

\begin{proof}
Ceci résulte du théorème de factorisation de Dixmier-Malliavin \cite{DiM}.
\end{proof}

On dira qu'un noyau $K(x,y)$ est $\autom$\hyph admissible si
pour $x$ fixé dans $\XG$ et pour tout $k\in\K$
la fonction
\[
y\mapsto K(xk,y)
\]
est lisse à support compact sur $\XG$
et si de plus l'espace vectoriel engendré par ces fonctions
est de dimension finie lorsque $k$ parcourt $\K$.

\begin{lemme} \label{factor}
Si $\fff$ est $\K$\hyph finie à droite et à gauche le noyau
$K_{\tG}(\fff,\omega;x,y)$ est somme finie de produits
$K_iH_i^*$ où les $K_i$ et les $H_i$ sont
des noyaux $\autom$\hyph admissibles.
\end{lemme}

\begin{proof}
Il résulte du théorème~\ref{dixmal} que le noyau $K_{\tG}(\fff,\omega;x,y)$ peut s'écrire
\[
K_{\tG}(\fff,\omega;x,y)=\sum_i \int_{\XG} K_{\tG}(g_i;x,z)K_\G^*(h_i,\omega;y,z)\dd z
\]
 où les $g_i$ et les $h_i$ sont $\K$\hyph finies à gauche.
\end{proof}

\section{Propriétés du noyau tronqué}

On utilisera l'opérateur de troncature sur un noyau
$K(x,y)$ en le faisant agir sur la première ou la seconde variable. On pose:
\[
\tronc_1^{T}K(x,y)=\tronc^{T}\phi(x)\qquad\text{pour $\phi(x)=K(x,y)$}
\]
et
\[
\tronc_2^{T}K(x,y)=\tronc^{T}\psi(y)\qquad\text{pour $\psi(y)=K(x,y)$}\ptf
\]
On rappelle que pour $\fff\in\ctyc(\Gadef)$ on a posé (\cf section~\ref{noyint})
\[
K_{\tG}(x,y)=\sum_{\delta\in\tG(F)}\fff^1(x\moins\delta y)\ptf
\]
\begin{lemme}\label{admiss}
Le noyau
\[
H(x,y)=\tronc_1^\T K_{\tG}(x,y)
\]
est $\autom$\hyph admissible si $\fff\in\ctyc(\Gadef)$ est $\K$\hyph finie à gauche.
\end{lemme}

\begin{proof}
L'opérateur de troncature appliqué à la première variable du noyau $K_{\tG}(x,y)$:
\[
\tronc_1^{T}K_{\tG}(x,y)
\]
fait intervenir une somme finie de termes indexés par des sous-groupes paraboliques
comportant chacun une intégration sur un compact (pour le calcul du terme constant
le long de $P$) et une somme en $\xi$ qui est finie pour $x$ fixé:
\[
 (-1)^{\ga_P-\ga_\Q} \sum_{\xi\in P(F)\bs\Q(F)}
\htau_P^\Q(\xixT )\int_{\N_P(F)\bs\N_P(\adef)} K_{\tG}(n\,\xi \,x,y)\dd n\ptf
\]
La fonction
\[
y\mapsto \tronc_1^\T K_{\tG}(x,y)
\]
est donc à support compact. On obtient une fonction $\K$\hyph finie en $x$ si $\fff$ est $\K$\hyph finie à gauche car
l'opérateur de troncature commute à l'action de $\K$.
\end{proof}

\begin{lemme}\label{rapdeK}
Il existe $N$ et
pour tout $M$ des constantes $c_M(f)$ et $c'_M(f)$ avec
pour $x$ dans un domaine de Siegel et pour tout $y$\textup:
\[
\lvert \tronc_1^{T} K_{\tG}(x,y)\rvert\le c_M(f)\,\lvert x\rvert ^{-M}\,\lvert y\rvert^N\ptf
\]
Enfin,\textup pour $x$ et $y$ dans un domaine de Siegel
\[
\lvert \tronc_1^{T} \Lambda_2^{T} K_{\tG}(x,y)\rvert \le c'_M(f)\,\lvert x\rvert^{-M}\,\lvert y\rvert^{-M}\ptf
\]
\end{lemme}

\begin{proof}
Les assertions résultent du lemme~\ref{premaj} et de la proposition~\ref{rapdec}.
\end{proof}
\chapter{Décomposition spectrale}\label{dspK}

\section{Sorites}

Soient $(X,\dd x)$, $(Y,\dd y )$ et $(\Lambda,\dd\lambda)$ trois espaces
localement compacts dénombrables à l'infini, munis de mesures de Radon.
On note $\langle \phi,\psi\rangle_\bullet$ le produit scalaire de deux fonctions dans $L^2(\bullet)$.
On suppose donnée une fonction continue sur $Y\times\Lambda$:
\[
E^Y(y,\lambda)\in \mathcal{C}(Y\times\Lambda)
\]
telle que, si on pose
\[
\widehat\phi(\lambda)=\int_Y \phi(y)E^Y(y,\lambda)\dd y
\]
pour $\phi$ continue et à support compact sur $Y$
on ait
\begin{equation}
\langle\phi,\psi\rangle_Y=
\int_\Lambda \widehat\phi(\lambda)\widebar{\widehat\psi(\lambda)}\dd\lambda
=\langle\widehat\phi,\widehat\psi\rangle_\Lambda
\ptf\tag{$*$}\label{eq7.*}
\end{equation}
On dit alors que $(Y,\dd y)$ est muni d'une décomposition spectrale supportée par
$(\Lambda,d\lambda)$. Considérons
un noyau intégral $K(x,y)$ (\cad une fonction sur $X\times Y$)
représentant un opérateur entre $L^2(Y)$ et $L^2(X)$.
On dira que $K$ est admissible si, pour tout $x\in X$ la fonction
\[
y\mapsto K(x,y)
\]
est continue à support compact. On considère un noyau $H$ de la forme
$K_1K_2^*$:
\[
H(x,y)=\int_Y K_1(x,z)K^*_2(z,y)\dd z\coloneqq\int_Y K_1(x,z)\widebar{K_2(y,z)}\dd z
\]
avec $K_1$
sur $X\times Y$ et $K_2$ sur $Y\times Y$ admissibles.
On pose
\[
\widehat{K}_i(z,\lambda)=\int_Y K_i(z,y)\widebar{E^Y(y,\lambda)}\dd y
\]
et
\[
H(x,y;\lambda)= \widehat{K}_1(x,\lambda) \widebar{\widehat{K}_2(y,\lambda)}\ptf
\]

\begin{proposition}\label{absdspK}
Si $H=K_1K_2^*$ on a
\begin{equation}
H(x,y)=\int_\Lambda H(x,y;\lambda)\dd\lambda
=\int_\Lambda \widehat{K}_1(x,\lambda)
\widebar{\widehat{K}_2(y,\lambda)}\dd\lambda\ptf\label{eq7.1}
\end{equation}
En particulier\textup, si $X=Y$ et
$H=KK^*$ avec $K$ admissible\textup, on a
\begin{equation}
H(x,x)=\int_\Lambda \lvert\widehat{K}(x,\lambda) \rvert^2\dd\lambda\ptf
\label{eq7.2}
\end{equation}
Plus généralement\textup, si $H=K_1K_2^*$ on a
\begin{equation}
\lvert H(x,y)\rvert\le K_1K_1^*(x,x)^{1/2}K_2K_2^*(y,y)^{1/2}
\ptf\label{eq7.3}
\end{equation}
\end{proposition}

\begin{proof} On pose
\[
\phi(z)=K_1(x,z)\Qquad{et}\psi(z)=K_2(y,z)
\]
et comme les $K_i$ sont admissibles les fonctions $\phi$ et $\psi$ sont
continues et à support compact.
La décomposition spectrale~\eqref{eq7.*} fournit~\eqref{eq7.1} et~\eqref{eq7.2}.
Pour obtenir~\eqref{eq7.3} on observe que
\[
\lvert H(x,y)\rvert \le\int_\Lambda \lvert H(x,y;\lambda)\rvert \dd\lambda=\int_\Lambda \lvert \widehat{K}_1(x,\lambda)
\widebar{\widehat{K}_2(y,\lambda)}\rvert \dd\lambda
\]
puis on invoque l'inégalité de Schwarz.
\[
\int_\Lambda \lvert\widehat{K}_1(x,\lambda)
\widebar{\widehat{K}_2(y,\lambda)}\rvert \dd\lambda
\le K_1K_1^*(x,x)^{1/2}K_2K_2^*(y,y)^{1/2}
\ptf\qedhere
\]
\end{proof}

\section{Le cas automorphe}\label{casautom}

On appelle donnée cuspidale pour $\G$ un couple
$(\M,\sigma)$ où $\M$ est un sous-groupe de Levi standard dans $\G$
et $\sigma$ une représentation cuspidale pour $\M$ triviale sur $\ga_\M$.
On dit que deux données cuspidales $\chi=(\M,\sigma)$ et $\chi'=(\M',\sigma')$ sont équivalentes s'il existe $g\in\G(F)$ tel que
\[
g\M g\moins=\M'\Qquad{et}\sigma'\circ\Ad(g)\simeq\sigma\ptf
\]
Il est bien connu que l'on peut décomposer $L^2(\XG)$
suivant les classes d'équivalence de données cuspidales:
on construit des sous-espaces $L^2_\chi(\XG)$
de $L^2(\XG)$, attachés à chaque donnée cuspidale $\chi=(\M,\sigma)$
au moyen des pseudo-séries d'Eisenstein (aussi appelées \og séries theta\fg).
Les données cuspidales inéquivalentes donnent naissance à des sous-espaces orthogonaux.
On renvoie le lecteur à \cite{Eisen} et \cite{MW} pour des définitions et des preuves détaillées.
En particulier la \og décomposition suivant les données cuspidales\fg est
le titre (et l'unique objet) du chapitre~II de \cite{MW}.

On notera $\Pi(\chi)$ le projecteur sur le sous-espace $L^2_\chi(\XG)$.
L'espace total est engendré par ces sous-espaces:
\[
L^2(\XG)=\widehat{\bigoplus_{\chi\in\XX}} L^2_\chi(\XG)
\]
où $\XX$ désigne l'ensemble des classes de données cuspidales.
C'est la décomposition suivant les données cuspidales. On prendra garde que
le projecteur $\Pi(\chi)$ commute à la représentation de $\Gadef$
mais pas, en général, à la représentation tordue.
La décomposition suivant les données cuspidales
fournit une décomposition du noyau $K_{\tG}$:
\[
K_{\tG}(x,y)=\sum_{\chi\in\XX} K_{\tG,\chi}(x,y)
\]
où $K_{\tG,\chi}$ est le noyau de l'opérateur
\[
\treg(\fff,\omega)\circ\Pi(\chi)\colon L^2_\chi(\XG)\to L^2(\XG)
\ptf
\]
La décomposition suivant les données cuspidales est aussi appelée
décomposition spectrale grossière.
Nous allons maintenant donner la décomposition spectrale fine.

On note $\Levi^\G$ l'ensemble des sous-groupes de Levi de $G$
contenant le sous-groupe de Levi minimal fixé $M_0$.
Soit $M\in\Levi^\G$ un sous-groupe de Levi de $G$; on notera $\weyl^\G(M)$ le quotient de
l'ensemble des $s\in \weyl^\G$ tels que $s(M)=M$ par $\weyl^M$, le groupe
de Weyl de $M$. On observera que $\weyl^G(M)$ est un groupe.
On note (\cf section~\ref{chambresfaces})
\[
w^\G(\M)= \lvert\weyl^\G(\M)\rvert
\]
son ordre.

Soit $P$ un sous-groupe parabolique
admettant $\M$ comme sous-groupe de Levi.
Pour toute représentation $\sigma$ automorphe de $\M$
choisissons une base orthonormale $\base^P(\sigma)$
dans l'espace préhilbertien $\autom(\XP,\sigma)$.
Comme une représentation automorphe est admissible tout vecteur $\Phi\in\autom(\XP,\sigma)$
est combinaison linéaire finie d'éléments de $\base^P(\sigma)$
autrement dit $\base^P(\sigma)$ est une base de l'espace vectoriel $\autom(\XP,\sigma)$.
Soit $\phi$ une fonction continue à support compact sur $\XG$.
Pour $\Psi\in\base^P(\sigma)$ et
\[
\mu\in(i\ga_P^\G)^*
\]
on définit
\[
\widehat\phi(\Psi,\mu)=\int_{\XG}\phi(x)\,\widebar{E(x,\Psi,\mu)}\dd x\ptf
\]
Soient $\phi$ et $\psi$ deux fonctions
continues et à support compact sur $\XG$. On note
\[
\langle \phi,\psi\rangle _\G\coloneqq \int_{\XG} \phi(x)\widebar{\psi(x)}\dd x
\]
leur produit scalaire.

\begin{theoreme}\label{dspth}
Le produit scalaire admet la décomposition spectrale suivante\textup:
\[
\langle\phi,\psi\rangle_\G=\sum_\chi \sum_{\M\in\Levi^\G/\weyl^\G}
\frac{1}{w^\G(\M)}\sum_{\sigma\in\Pi_{\mathrm{disc}}(M)_\chi}
\int_{\ima(\ga_\M^\G)^*}
\sum_{\Psi\in\base^P(\sigma)}
\widehat\phi(\Psi,\mu)\,\widebar{\widehat\psi(\Psi,\mu)}\dd\mu
\]
où $\dd\mu$ est la mesure de Haar sur le groupe
$\ima(\ga_\M^\G)^*$ duale, au sens de la transformation de Fourier\textup,
de celle utilisée sur $\ga_\M^\G$\,\footnote{Un facteur $(1/2\ima\pi)^{a_\M^\G}$
intervient lorsqu'on qu'on munit $\ga^*$ de la mesure duale de celle sur $\ga$
au sens des espaces vectoriels au lieu de la mesure duale sur $\ima\ga^*$
au sens de la transformation de Fourier.}.
La somme sur $\sigma$ porte sur l'ensemble $\Pi_{\mathrm{disc}}(M)_\chi$
des classes de représentations automorphes de $\M$
intervenant discrètement dans $L^2_\chi(\XM)$.
\end{theoreme}

\begin{proof}
On renvoie le lecteur à \cite{Eisen} ou \cite{MW} pour une preuve.
\end{proof}

Soit $K(x,y)$ un noyau intégral sur $\XG$. On rappelle que le noyau $K$ est dit $\autom$\hyph admissible si
pour $x$ fixé dans $\XG$ et pour tout $k\in\K$ la fonction
\[
y\mapsto K(xk,y)
\]
est lisse à support compact sur $\XG$
et si de plus l'espace vectoriel engendré par ces fonctions
est de dimension finie lorsque $k$ parcourt $\K$.

\begin{proposition}\label{dspker}
On considère une sous-groupe parabolique $P$.
Soit $H(x,y)$ un noyau intégral de la forme $K_1K_2^*$ où les $K_i$ sont des noyaux
$\autom$\hyph admissibles sur $\X_P$.
Supposons de plus qu'il existe un automorphisme $\theta$ de $\G$ tel que\textup:
si $S$ est un
sous-groupe parabolique de $P$
et $\sigma$ une représentation automorphe discrète de $\M$\textup,
le sous-groupe de Levi de $S$\textup,
on ait\textup, pour $\mu\in{i(\ga_\M^\G)^*}$\textup,
des opérateurs bornés de rang fini\textup:
\[
A_{1,\sigma,\mu}\in\Hom\Bigl(\autom(\X_S,\sigma),\autom\bigl(\X_{\theta(S)},\theta(\sigma)\bigr)\Bigr)
\]
et
\[
A_{2,\sigma,\mu}\in\Hom\bigl(\autom(\X_S,\sigma),\autom(\X_S,\sigma)\bigr)
\]
vérifiant
\[
\int_{\X_P} K_1(x,y)\,E^P(y,\Psi,\mu)\dd y= E^Q(x,A_{1,\sigma,\mu}\Psi,\theta\mu)
\]
où $\Q=\theta(P)$ et
\[
\int_{\X_P} K_2(x,y)\,E^P(y,\Psi,\mu)\dd y=E^P(x,A_{2,\sigma,\mu}\Psi,\mu)\ptf
\]
Alors le noyau $H_\chi(x,y)$ admet la décomposition spectrale
suivante\textup:
\begin{equation}
H_\chi(x,y)=\sum_{\M\in\Levi^P/\weyl^P}
\frac{1}{w^P(\M)}\sum_{\sigma\in\Pi_{\mathrm{disc}}(M)_\chi}
\int_{\ima(\ga_\M^\G)^*} H_{\sigma}(x,y;\mu)
\dd\mu \label{eq7.1a}
\end{equation}
avec
\[
H_{\sigma}(x,y;\mu)=\sum_{\Psi\in\base^P(\sigma)}
E^\Q(x,B_{\sigma,\mu}\Psi,\theta\mu)
\widebar{E^P(y,\Psi,\mu)}
\]
où
\[
B_{\sigma,\mu}=A_{1,\sigma,\mu}A_{2,\sigma,\mu}^*
\]
et la somme en $\Psi$ porte sur un ensemble fini.
Enfin si on pose
\[
h_\chi(x,y)=\sum_{\M\in\Levi^P/\weyl^P}
\frac{1}{w^P(\M)}\sum_{\sigma\in\Pi_{\mathrm{disc}}(M)_\chi}
\int_{\ima(\ga_\M^\G)^*}\lvert H_{\sigma}(x,y;\mu)\rvert\dd\mu 
\]
on a la majoration
\begin{equation}
\sum_\chi\lvert H_\chi(x,y)\rvert\le\sum_\chi
h_\chi(x,y)\le K_1K_1^*(x,x)^{1/2}K_2K_2^*(y,y)^{1/2}\ptf \label{eq7.2a}
\end{equation}
\end{proposition}

\begin{proof} 
Il résulte des généralités sur la décomposition spectrale des
noyaux produits (proposition~\ref{absdspK}\eqref{eq7.1}) et de la forme explicite de la décomposition
spectrale automorphe (théorème~\ref{dspth}) que $H(x,y)$ est donné par
\[
\sum_{\M\in\Levi^P/\weyl^P}
\frac{1}{w^P(\M)}\sum_{\sigma\in\Pi_{\mathrm{disc}}(M)_\chi}
\int_{\ima(\ga_\M^\G)^*}\sum_{\Psi\in\base^P(\sigma)}
H_{\sigma}(x,y;\Psi,\mu)\dd \mu
\]
où
\[
H_{\sigma}(x,y;\Psi,\mu)=
E^\Q(x,A_{1,\sigma,\mu}\Psi,\theta\mu)
\widebar{E^P(y,A_{2,\sigma,\mu}\Psi,\mu)}\ptf
\]
On observe que comme par hypothèse les opérateurs $A_{i,\sigma,\mu}$ sont supposés de rang fini,
l'opérateur produit
\[
B_{\sigma,\mu}=A_{1,\sigma,\mu}A_{2,\sigma,\mu}^*
\]
est tel que l'ensemble des $\Psi\in\base^P(\sigma)$
pour lesquels $B_{\sigma,\mu}\Psi\ne0$ est fini.
On peut donc définir
un noyau $H_{\sigma}(x,y;\mu)$ comme fonction lisse de trois variables
en posant:
\[
H_{\sigma}(x,y;\mu)=\sum_{\Psi\in\base^P(\sigma)}
E^\Q(x,B_{\sigma,\mu}\Psi,\theta\mu)
\widebar{E^P(y,\Psi,\mu)}
\]
puisque, pour chaque $\sigma$,
les sommes en $\Psi$ sont finies.
Soit $F$ une partie finie de $\base^P(\sigma)$. Si $F$ est assez grand
le sous-espace vectoriel engendré par $F$
contient les images des $A_{i,\sigma,\mu}$ et de $B_{\sigma,\mu}$.
Un calcul élémentaire d'algèbre linéaire montre alors que pour un tel $F$ on a
\[
\sum_{\Psi\in F}
E^\Q(x,B_{\sigma,\mu}\Psi,\theta\mu)
\widebar{E^P(y,\Psi,\mu)}=\sum_{\Psi\in F}
E^\Q(x,A_{1,\sigma,\mu}\Psi,\theta\mu)
\widebar{E^P(y,A_{2,\sigma,\mu}\Psi,\mu)}\ptf
\]
Donc pour tout $F$ assez grand on a
\[
H_{\sigma}(x,y;\mu)=\sum_{\Psi\in F}H_{\sigma}(x,y;\Psi,\mu)
\]
et par passage à la limite on en déduit que
\[
H_{\sigma}(x,y;\mu)=\sum_{\Psi\in\base^P(\sigma)}
H_{\sigma}(x,y;\Psi,\mu)\ptf
\]
Ceci établit~\eqref{eq7.1a}. L'assertion~\eqref{eq7.2a} résulte de l'équation~\eqref{eq7.3} de la proposition~\ref{absdspK}.
\end{proof}

\section{Estimée d'un noyau}\label{estimm}

On reprend les notations de la section~\ref{noyint} et de la proposition~\ref{dspker}
mais on écrit $\pQ$ pour $P$. On pose $\Q=\theta(\pQ)$.
Soient $S$ un sous-groupe parabolique standard
de $\pQ$, $\M$ son Levi standard et $\sigma$ une représentation automorphe pour $\M$.
Si $\Psi$ appartient à l'espace des fonctions de type $\sigma$
sur $\X_S$ \cad si
\[
\Psi\in\autom(\X_S,\sigma)
\]
on définit une fonction $\Phi$ sur $\X_{\theta(S)}$:
en posant
\[
\treg_{S,\sigma,\mu}(\fff,\omega)\Psi=\Phi
\]
On a
\[
\treg(\fff,\omega)E^{\pQ}(x,\Psi,\mu)=E^\Q(x,\treg_{S,\sigma,\mu}(\fff,\omega)\Psi,\theto\mu)\ptf
\]
On  considère le noyau
\[
K_{\Q,\delta}(x,y)=
\int_{\N_\Q(F)\bs \N_\Q(\adef)} \omega(x)
\sum_{\eta\in \Q(F)}
\fff^1 (x\moins\,n_\Q\moins\,\eta\moins\,\delta\,y)\dd n_\Q\ptf
\]
La fonction $K_{\Q,\delta}(x,y)$ est le noyau intégral représentant l'opérateur
défini par
\[
\psi\mapsto \phi=\reg(\delta,\fff,\omega)\psi \qquad\text{où 
$\phi(x)=\int_{\tGadef}(\omega\psi)(\delta\moins x\,y)\fff(y)\dd y$}
\]
entre $L^2(\XpQG)$ et $L^2(\XQG)$.
La décomposition de $K_{\Q,\delta}(x,y)$ suivant les données cuspidales s'écrit
\[
K_{\Q,\delta}(x,y)=\sum_\chi K_{\Q,\delta,\chi}(x,y)
\]
où $K_{\Q,\delta,\chi}(x,y)$ est le noyau de l'opérateur
entre
$L^2_\chi(\XpQG)$ et $L^2(\XQG)$
défini par
$\reg(\delta,\fff,\omega)\circ\Pi(\chi)$.

\begin{proposition} \label{majo} Consid\'erons une fonction $\fff\in\ctyc(\Gadef)$ qui soit $\K$\hyph finie
 à droite et à gauche.
\begin{enumerate}[(i)]
\item Le noyau $ K_{\Q,\delta,\chi}(x,y)$ admet le développement  spectral
suivant\textup:
\[
\sum_{\M\in\Levi^{\pQ}/\weyl^{\pQ}}
\frac{1}{w^{\pQ}(\M)}\sum_{\sigma\in\Pi_{\mathrm{disc}}(M)_\chi}
\int_{\ima(\ga_\M^\G)^*} K_{\Q,\pQ,\sigma}(x,y;\mu)\dd\mu
\]
avec
\[
K_{\Q,\pQ,\sigma}(x,y;\mu)=
\sum_{\Psi\in\base^{\pQ}(\sigma)}
E^\Q(x,\treg_{S,\sigma,\mu}(\fff,\omega)\Psi,\theta\mu)
\widebar{E^{\pQ}(y,\Psi,\mu)}\ptf
\]
\item
Il existe $N$ tel que
pour tout $M$ il existe $c$ tel que
\[
\sum_\chi\sum_{\M\in\Levi^{\pQ}/\weyl^{\pQ}}
\frac{1}{w^{\pQ}(\M)}\sum_{\sigma\in\Pi_{\mathrm{disc}}(M)_\chi}
\int_{i(\ga_\M^\G)^*} \lvert \tronc_1^{T,\Q }K_{\Q,\pQ,\sigma}(x,y;\mu)
\rvert\dd\mu
\]
est majoré par
\[
c\,\lvert x\rvert^{-M}\,\lvert y\rvert^{N}\ptf
\]
\end{enumerate}
\end{proposition}

\begin{proof}
L'assertion~(i) est une conséquence immédiate de la proposition~\ref{dspker}. On en déduit que
\[
\tronc_1^{T,\Q }K_{\Q,\delta,\chi}(x,y)
\]
est égal à 
\[
\sum_{\M\in\Levi^{\pQ}/\weyl^{\pQ}}
\frac{1}{w^{\pQ}(\M)}\sum_{\sigma\in\Pi_{\mathrm{disc}}(M)_\chi}
\int_{\ima(\ga_\M^\G)^*} \tronc_1^{T,\Q }K_{\Q,\pQ,\sigma}(x,y;\mu)\dd\mu\ptf
\]
D'après le théorème~\ref{dixmal} le noyau
$ \tronc_1^{T,\Q }K_{\Q,\pQ,\sigma}$
est égal à une somme finie de produits de noyaux:
\[
 \tronc_1^{T,\Q }K_{\Q,\pQ,\sigma}(x,y,\mu)=\sum_i \int_{\XG} \tronc_1^{\T,\Q} K_{\Q,\pQ,\sigma}(g_i;x,z,\mu)
K_{\pQ,\pQ,\sigma}^*(h_i,\omega;y,z,\mu)\dd z\ptf
\]
Suivant l'équation~\eqref{eq7.3} de la proposition~\ref{absdspK},
l'inégalité de Schwartz montre
que l'expression
\[
\sum_\chi\sum_{\M\in\Levi^{\pQ}/\weyl^{\pQ}}
\frac{1}{w^{\pQ}(\M)}\sum_{\sigma\in\Pi_{\mathrm{disc}}(M)_\chi}
\int_{i(\ga_\M^\G)^*} \lvert \tronc_1^{T,\Q }K_{\Q,\pQ,\sigma}(x,y;\mu)
\rvert\dd\mu
\]
est majorée par une somme finie de termes
du type
\[
\lvert \tronc_1^{\T,\Q} \tronc_2^{\T,\Q} K_\Q(g*g^*;x,x)\rvert^{1/2}
\lvert K_{\pQ}(h*h^*;y,y)\rvert^{1/2}
\]
enfin on invoque les lemmes~\ref{premaj} et~\ref{rapdeK} pour établir
qu'il existe $N$ tel que
pour tout $M$ il existe $c$ tel que
ceci soit majorée par
\[
c\,\lvert x\rvert^{-M}\,\lvert y\rvert^{N}
\]
pour tout $x$ et tout $y$ dans un domaine de Siegel de $\M_\Q$.
On en déduit~(ii).
\end{proof}

\begin{corollaire}\label{majod} Consid\'erons une fonction $\fff\in\ctyc(\Gadef)$ qui soit $\K$\hyph finie
 à droite et à gauche. Alors,
\[
\sum_\chi\lvert \tronc_1^{\T,\Q} K_{\Q,\delta,\chi}(x,y)\rvert \le c\,\lvert x\rvert^{-M}\,\lvert y\rvert^{N}
\]
pour tout $x$ et tout $y$ dans un domaine de Siegel de $\M_\Q$.
\end{corollaire}

\begin{proof}
L'assertion est une conséquence immédiate de la proposition~\ref{majo}(ii). On aurait pu aussi
invoquer directement l'équation~\eqref{eq7.2a} de la proposition~\ref{dspker} puis, comme ci-dessus faire appel aux lemmes~\ref{premaj} et~\ref{rapdeK}.
\end{proof}

\part{La formule des traces grossière}

\chapter{Formule des traces: état zéro}\label{ch8}

\section{La problématique}\label{problemat}

Soit $\fff\in\ctyc(\Gadef)$. On a posé 
\[
\fff^1(\delta)=\int_{z\in\gA_\G}\fff(z\,\delta)\dd z\qquad\text{pour $\delta\in\tGadef$}\ptf
\] 
On introduit également:
\[
\tff(\delta)=\jtG\int_{\gAtG}\fff(z\,\delta)\dd z\qquad\text{avec $\jtG=\lvert\det(\theto-1\rest \ga_\G/\ga_{\tG})\rvert$}\ptf
\]
\newnot{f@$\tff$}{ffo}%
On a l'identité suivante: pour tout $\delta\in\tGadef$ 
\[
\fff^1(\delta)=\int_{\gAtG\bs\gAG}\tff(a\moins \delta\,a)\dd a\ptf
\]

On rappelle (\cf section~\ref{noyint}) que l'on a défini pour $y\in\tGadef$
l'opérateur $\treg(y,\omega)$ par
\[
(\treg(y,\omega)\vf)(\dot x)=(\omega\vf)(\delta\moins x\,y)
\]
pour $\vf$ dans $L^2(\XG)$
avec $\delta$ quelconque dans $\tG(F)$ 
et que pour $\fff\in\ctyc(\Gadef)$ l'opérateur 
\[
\treg(\fff,\omega)=\int_{\tGadef} \fff(y)\treg(y,\omega)\dd y
\]
est représenté par le noyau 
\[
\K_{\tG}(x,y)=\sum_{\delta\in\tG(F)}\fff^1(x\moins\delta y)\qquad\text{pour $x,y\in\Gadef$}\ptf
\]
Le noyau est une fonction lisse sur $\XG\times\XG$.

Considérons $\delta\in\tG(F)$ quasi semi-simple.
Suivant les conventions de la section~\ref{ellip} on note
$\G^\delta$ le centralisateur de $\delta\in\tG(F)$, $\G_\delta$ sa composante neutre
et on appelle centralisateur stable le groupe $\centd=\G_\delta.Z_{\tG}$.
On introduit
\[
\orb_\delta(\fff,\omega)
=\int_{\centd(\adef)\bs\Gadef}\omega(x)\fff(x\moins\delta\,x)\dd \dot x
\]
\newnot{Odelta(f,omega)@$\orb_\delta(\fff,\omega)$}{orbdelta}%
si $\centd(\adef)$ est dans le noyau de $\omega$
et $\orb_\delta(\fff,\omega)=0$
sinon. C'est
l'intégrale orbitale tordue par $\omega$.
Nous aurons aussi besoin du nombre $a^\G(\delta)$ est défini par
\[
a^\G(\delta)
=\iota(\delta)\moins\vol\bigl(\gA_{\tG}\centd(F)\bs \centd(\adef)\bigr)
\]
\newnot{aG(delta)@$a^G(\delta)$}{ahatG}%
où $\iota(\delta)$ est l'ordre du quotient $\G^\delta(F)/\centd(F)$.

\begin{remarque}\label{idgd}
Pour éviter de manipuler des groupes non connexes, ce qui est en général le cas pour $\centd$,
le lecteur pourra à sa guise remplacer systématiquement $\centd$ par $G_\delta$
dans la définition des intégrales orbitales $\orb_\delta(\fff,\omega)$ et des coefficients $a^\G(\delta)$
ainsi que dans les expressions analogues intervenant dans la section~\ref{explicitss}. 
C'est d'ailleurs le point de vue adopté par Arthur dans \citelist{\cite{Ageom} \cite{AinvI} \cite{AinvII}} par exemple. 
Toutefois, l'introduction de $\centd$ semble indispensable dans l'étude de la stabilisation
dans le cas tordu (\cf \citelist{\cite{KS} \cite{Ltw}}).
\end{remarque}

On notera
$\tGamma$ un ensemble de représentants des classes de conjugaison
dans $\tG(F)$.
On a défini à la section~\ref{ellip} la notion d'élément elliptique.
On notera $\tG(F)_{\mathrm{ell}}$ l'ensemble des éléments
elliptiques dans $\tG(F)$ et $\tGamma_{\mathrm{ell}}$ un ensemble
de représentants des classes de $\G(F)$\hyph conjugaison dans $\tG(F)_{\mathrm{ell}}$.
On pose
\[
k_{\mathrm{ell}}(x)=\sum_{\delta\in\tG(F)_{\mathrm{ell}}}\omega(x)
\fff^1(x\moins \delta\,x)\ptf
\]

\begin{proposition} \label{kell}
L'intégrale de $k_{\mathrm{ell}}(x)$ est convergente. On pose
\begin{equation}
J_{\mathrm{ell}}^{\tG}(\fff,\omega)=\int_{\XG} k_{\mathrm{ell}}(x)\dd x\ptf\label{eq8.1}
\end{equation}
C'est une distribution invariante et
\begin{equation}
J_{\mathrm{ell}}^{\tG}(\fff,\omega)=
\sum_{\delta\in\tGamma_{\mathrm{ell}}}a^\G(\delta)\orb_\delta(\tff,\omega)\ptf\label{eq8.2}
\end{equation}
\end{proposition}

\begin{proof}
Un calcul élémentaire fournit l'égalité~\eqref{eq8.2} au moins formellement.
Grâce au lemme~\ref{ssfini}, on voit que dans l'expression~\eqref{eq8.2}
la somme porte en fait
sur un ensemble fini (dépendant du support de $\fff$).
On en déduit la convergence et l'égalité.
On renvoie au théorème~\ref{geoconv} pour une preuve plus détaillée dans un cas plus général.
\end{proof}

La représentation de $\Gadef$ dans $L^2(\XG)$ comporte en général un spectre discret 
et un spectre continu:
\[
L^2(\XG)=L^2_{\mathrm{disc}}(\XG)\oplus L^2_{\mathrm{cont}}(\XG)\ptf
\]
Nous noterons
\[
\Pi(\tG,\omega)
\]
l'ensemble des représentations irréductibles $\pi$ de $\Gadef$
admettant un prolongement tordu $\tpi$ et
\[
\Pi_{\mathrm{disc}}(\tG,\omega)
\]
le sous-ensemble des classes 
qui apparaissent dans $L^2_{\mathrm{disc}}(\XG)$. 
On renvoie le lecteur à la discussion de la notion de multiplicité tordue $m(\pi,\tpi)$
dans la section~\ref{multi} et on rappelle que le nombre
\[
m(\pi,\tpi)\tr\tpi(\fff)
\]
est indépendant du choix du prolongement $\tpi$. 

\begin{proposition} \label{tracedisc}
L'opérateur $\treg(\fff,\omega)$ est à trace dans le spectre discret.
On note
\[
J_{\G,\mathrm{disc}}^{\tG}(\fff,\omega)=
\tr\bigl(\treg(\fff,\omega)\rest L^2_{\mathrm{disc}}(\XG)\bigr)
\]
cette trace\,\footnote{On verra que d'autres termes \og discrets\fg indexés par les
classes de sous-groupes de Levi
interviennent dans la formule des traces.}.
On a alors
\[
J_{\G,\mathrm{disc}}^{\tG}(\fff,\omega)=
\sum_{\pi\in\Pi_{\mathrm{disc}}(\tG,\omega)}m(\pi,\tpi)\,\tr\tpi(\fff)\ptf
\]
\end{proposition}

\begin{proof}
On sait grâce à W.~Müller \cite{Mu} que pour $h\in\ctyc\bigl(\Gadef\bigr)$
l'opérateur $\reg(h)$ est à trace dans le spectre discret. Comme 
$\treg(\fff,\omega)$ est produit de $\reg(h)$ avec $h(x)=\fff(x\delto)$ et d'opérateurs
unitaires:
\[
\treg(\fff,\omega)=A(\omega)B(\theto)\reg(h)
\]
on en déduit que l'opérateur $\treg(\fff,\omega)$ est encore à trace dans le spectre discret.
Comme observé au lemme~\ref{tracenulle}, seules les représentations $\tpi$
dont la restriction $\pi$ à $\G(\adef)$ sont irréductibles peuvent donner une contribution non nulle
à la trace de $\treg(\fff,\omega)$ dans le spectre discret. On obtient ainsi la formule souhaitée.
\end{proof}

La formule des traces de Selberg, dans le cas compact, est
l'égalité entre l'intégrale du noyau $K_{\tG}(x,y)$ sur la diagonale
et, d'autre part,
la trace de cet opérateur développée suivant la décomposition spectrale
de $L^2(\XG)$, qui est une somme discrète avec multiplicité finie
de représentations irréductibles.

\begin{proposition}
Lorsque $\XG$ est compact\textup, \cad lorsque $\G_{\mathrm{der}}$ est anisotrope sur $F$,
on a
\[
\sum_{\delta\in\tGamma}a^\G(\delta)\orb_\delta(\tff,\omega)=
\sum_{\pi\in\Pi(\tG,\omega)}m(\pi,\tpi)\,\tr\tpi(\fff)\ptf
\]
\end{proposition}

\begin{proof} 
La compacité de $\XG$ implique que
\[
\tr\treg(\fff,\omega)=\int_{\XG}\,K_{\tG}(x,x)\dd x\ptf
\]
Lorsque $\G_{\mathrm{der}}$ est anisotrope toutes les classes de conjugaison sont elliptiques,
et la proposition~\ref{kell} s'écrit simplement
\[
\int_{\XG}\,K_{\tG}(x,x)\dd x=
\sum_{\delta\in\tGamma}a^\G(\delta)\orb_\delta(\tff,\omega)\ptf
\]
Comme $\XG$ est compact le spectre continu est nul et la proposition~\ref{tracedisc} s'écrit:
\[
\tr\treg(\fff,\omega)=
\sum_{\pi\in\Pi(\tG,\omega)}m(\pi,\tpi)\,\tr\tpi(\fff)\ptf\qedhere
\]
\end{proof}

Lorsque $\G_{\mathrm{der}}$ n'est pas anisotrope l'existence de classes non elliptiques a pour conséquence que,
pour le noyau tout entier, l'intégrale sur la diagonale est (en général) divergente.
Comme $\XG$ n'est plus compact l'opérateur n'est pas (en général)
à trace à cause de l'existence du spectre continu. 
La formule des traces est l'égalité du développement géométrique 
et du développement spectral pour une \og trace renormalisée\fg de l'opérateur. 
La renormalisation se fait en soustrayant les contributions divergentes,
au moyen de troncatures qui dépendent d'un paramètre $T\in\gao$ et
du choix de $\PO$, $\MO$ et $\K$.
Cette dépendance implique que la distribution obtenue
n'est pas invariante par conjugaison. 
Nous noterons $k^\T$ la restriction à la diagonale du noyau tronqué.
On disposera de deux expressions pour $k^\T$, notées $k^\T_{\mathrm{geom}}$ et $k^\T_{\mathrm{spec}}$,
dont l'égalité est appelée identité fondamentale (\cf proposition~\ref{fonfon}).
On démontrera que l'intégrale de $k^\T$ est convergente 
et on exhibera une expression asymptotique 
$J^\T$, de l'intégrale sur $\XG$ de $k^\T$, qui est polynomiale en $T$
(\cf théorème~\ref{poly}). On rappelle que l'on a introduit dans le lemme~\ref{Yorth} un élément $\TK$.
La trace renormalisée $J(\fff,\omega)$
sera, par définition, la valeur en $T=\TK$ de $J^\T$.
On obtiendra ainsi une expression indépendante du choix de $\PO$ (pour
$\MO$ et $\K$ fixés).
L'expression $k^\T_{\mathrm{geom}}$ se prête bien au développement suivant les classes de conjugaison,
appelé développement géométrique,
ainsi qu'à la preuve du caractère polynomial en $T$ de 
l'expression asymptotique $J^\T$.
Pour le développement suivant les classes de représentations
automorphes, appelé développement
spectral, on utilisera l'expression $k^\T_{\mathrm{spec}}$.
Lorsque $\XG$ est compact on a les égalités 
\[
J(\fff,\omega)=\sum_{\delta\in\tGamma_{\mathrm{ell}}}a^\G(\delta)\orb_\delta(\tff,\omega)=
\sum_{\pi\in\Pi_{\mathrm{disc}}(\tG,\omega)}m(\pi,\tpi)\,\tr\tpi(\fff)
\]
mais elles cessent d'être vraies en général; c'est toutefois vrai pour certaines fonctions $\fff$.
On dit alors que l'on dispose d'une formule des traces simple. De nombreux cas ont été étudiés
dans la littérature, mais nous n'en dirons rien de plus ici.

\section{L'identité fondamentale}

On va utiliser diverses troncatures. Nous aurons besoin, pour les définir
d'un analogue du noyau $K$ pour chaque sous-espace parabolique.
Soit $\tP$ un sous-espace parabolique de sous-espace de Levi $\tM$
contenant $\tM_0$ et de radical unipotent $N$. On pose
\[
K_{\tP}(x,y)=\int_{\N(F)\bs\N(\adef)}
\sum_{\delta\in\tP(F)}\omega(y)
\fff^1(x\moins\delta\, n\,y)\dd n\ptf
\]
\newindex{KP@$K_{\protect\tP}$}{noyautP}%
C'est le noyau de la représentation naturelle de $(\tGadef,\omega)$
dans $L^2(\XPG)$ où, avec les notations de la section~\ref{espacesXP}, on a posé
\[
\XPG=\gA_\G\,P(F)N(\adef)\bs\Gadef \ptf
\]
On utilisera l'opérateur de troncature sur ces noyaux en le faisant agir sur la première
variable \cad que par définition, si $\Q\subset P$:
\[
\tronc_1^{T,\Q}K_{\tP}(x,y)=\tronc^{T,\Q}\phi(x)\qquad\text{pour $\phi(x)=K_{\tP}(x,y)$}\ptf
\]

\begin{lemme}\label{idenfon}
Soit $\tP$ un sous-ensemble parabolique et soit $\phi$ une fonction sur $P(F)\bs\Gadef)$. On a
\[
\sum_{\{\Q,\R\mid \Q\subset P\subset\R\}}\,\sum_{\xi\in\Q(F)\bs P(F)}
\tsQR(\xixT)
\,\tronc^{T,\Q}\phi(\xix)=\htau_{\tP}(\xT)\,\phi_P(x)\ptf
\]
\end{lemme}

\begin{proof} 
On commence par observer que, par définition de l'opérateur de troncature,
\[
\sum_{\Q\subset P}\,\sum_{\xi\in\Q(F)\bs P(F)}
\tau_\Q^P(\xixT )\,\tronc^{T,\Q}\phi(\xix)
\]
est égal à
\[
\sum_{S\subset\Q\subset P}\,\sum_{\xi\in S(F)\bs P(F)}
(-1)^{\ga_S-\ga_\Q}\,\tau_\Q^P(\xixT )\,\htau_S^\Q(\xixT )\,\phi_S(\xix)
\]
où les sommes en $\xi$
portent sur des ensembles finis (grâce aux lemmes~\ref{cone} et~\ref{finitude}).
En effectuant d'abord la somme en $\Q$ sur les sous-groupes paraboliques 
tels que $S\subset\Q\subset P$ et compte tenu de la proposition~\ref{totoun} on voit que
seul le terme avec $S=\Q=P$ subsiste et l'expression se réduit à $\phi_P(x)$. On a donc
\[
\sum_{\Q\subset P}\,\sum_{\xi\in\Q(F)\bs P(F)}
\tau_\Q^P(\xixT )\,\tronc^{T,\Q}\phi(\xix)=\phi_P(x)\ptf
\]
Pour conclure il reste à observer que d'après~\ref{repart} on a
\[
\sum_{\{\R\mid P\subset\R\}}\tsQR(\xT)
\,\tronc^{T,\Q}\phi(x)=
\htau_{\tP}(\xT)\,\tau_\Q^P(\xT)\,\tronc^{T,\Q}\phi(x)\ptf\qedhere
\]
\end{proof}

On pose
\[
k^\T_{\tP,\mathrm{geom}}(x)=\htau_{\tP}(\xT)\,K_{\tP}(x,y)
\]
et 
\[
k^\T_{\tP,\mathrm{spec}}(x)=\sum_{\{\Q,\R\mid \Q\subset P\subset\R\}}
\,\sum_{\xi\in\Q(F)\bs P(F)}
\tsQR(\xixT )
\,\tronc^{T,\Q}_1K_{\tP}(\xix,y)\ptf
\]

Nous sommes maintenant en mesure d'énoncer l'identité fondamentale
appelée \textenglish{Ba\-sic Iden\-ti\-ty} dans (\cite{MS}*{Lectures~1, 2,~9}).
On pose 
\[
k^\T_{\mathrm{geom}}(x)=\sum_{\tP\supset\tPO}(-1)^{a_{\tP}-a_{\tG}}\,
\sum_{\xi\in P(F)\bs\G(F)}\htautP(\xixT )\,K_{\tP}(\xix,\xix)
\]
et
\[
k^\T_{\mathrm{spec}}(x)=\sum_{\tP\supset\tPO}(-1)^{a_{\tP}-a_{\tG}}
\,\sum_{\Q\subset P\subset\R}\,\sum_{\xi\in\Q(F)\bs\G(F)}
\tsQR(\xixT )
\,\tronc_1^{T,\Q}K_{\tP}(\xix,\xix)\ptf
\]

\begin{proposition}\label{fonfon}
Les fonctions $k^\T_{\tP,\mathrm{geom}}$\textup, $k^\T_{\mathrm{geom}}$\textup, $k^\T_{\tP,\mathrm{spec}}$ et $k^\T_{\mathrm{spec}}$
ne dépendent que de la projection
de $\T$ sur le sous-espace des $\theto$\hyph invariants dans $\gao^\G$
et on a les identités
\[
 k^\T_{\tP,\mathrm{geom}}=k^\T_{\tP,\mathrm{spec}}\Qquad{et}
k^\T_{\mathrm{geom}} = k^\T_{\mathrm{spec}}\ptf
\]
\end{proposition}

\begin{proof}
Par définition, on a
\[
k^\T_{\mathrm{geom}}(x)=\sum_{\tP\supset\tPO}(-1)^{a_{\tP}-a_{\tG}}\,
\sum_{\xi\in P(F)\bs\G(F)}k^\T_{\tP,\mathrm{geom}}(\xix)
\]
et
\[
k^\T_{\mathrm{spec}}(x)=\sum_{\tP\supset\tPO}(-1)^{a_{\tP}-a_{\tG}}\,
\sum_{\xi\in P(F)\bs\G(F)}k^\T_{\tP,\mathrm{spec}}(\xix)
\]
On observe tout d'abord que les sommes en $\xi$
portent sur des ensembles finis, d'après le lemme~\ref{finitudebis}, et 
les expressions sont donc trivialement convergentes.
Maintenant on a
\[
k^\T_{\tP,\mathrm{geom}}(x)=k^\T_{\tP,\mathrm{spec}}(x)
\]
d'après le lemme~\ref{idenfon} 
appliqué à $\phi(x)=K_{\tP}(x,y)$
en observant que dans ce cas on a $\phi_P=\phi$.
\end{proof}

\chapter{Développement géométrique}\label{ch9}

\section{Convergence: côté géométrique}\label{devgeo}

On dira que deux éléments dans $\tG(F)$ sont ss\hyph conjugués si leurs parties quasi semi-simples
sont conjuguées. On notera $\gO$ l'ensemble des
classes de ss\hyph conjugaison.
On peut décomposer $k^\T_{\mathrm{geom}}(x)$ suivant les classes
de ss\hyph conjugaison:
\[
k^\T_{\mathrm{geom}}(x)=\sum_{\go\in\gO} k^\T_\go(x)
\]
où $k^\T_\go$ ne comporte que les contributions
d'éléments $\delta$ dont la partie quasi semi-simple appartient à une même classe de conjugaison:
\[
k^\T_\go(x)=\sum_{\tP\supset\tPO}(-1)^{a_{\tP}-a_{\tG}}\,
\sum_{\xi\in P(F)\bs\G(F)}\htautP(\xixT)\,K_{\tP,\go}(\xix,\xix)\ptf
\]
avec
\[
K_{\tP,\go}(x,x)=\int_{\N(F)\bs\N(\adef)}
\sum_{\delta\in\go\cap\tP(F)}\omega(x)
\fff^1(x\moins\delta\, n\,x)\dd n\ptf
\]

On considère deux sous-groupes paraboliques standard
$\Q\subset\R$. On a discuté dans le lemme~\ref{qplus}
les propriétés des sous-groupes paraboliques $\Qp$ et $\Rm$.

\begin{proposition}\label{geoconvb}
On considère l'espace quotient
\[
\YQ=\gA_\G\Q(F)\bs\Gadef \ptf
\]
Supposons $T$ assez régulier
\cad $\dPO(T)\ge c$ où $c$ est une constante  dépendant  du support de $\fff$. L'intégrale
\[
\int_{\YQ}\FPQ(x,T){\tsQR}(\xT)\Biggl\lvert
\sum_{\{\tP\mid \tQp\subset\tP\subset\tRm\}}(-1)^{a_{\tP}-a_{\tG}}\,
K_{\tP,\go}(x,x)\Biggr\rvert \dd x
\]
est convergente.
\end{proposition}

\begin{proof}
On rappelle que
\[
K_{\tP,\go}(x,x)=\int_{\N(F)\bs\N(\adef)}
\sum_{\delta\in\go\cap\tP(F)}\omega(x)
\fff^1(x\moins\delta\, n\,x)\dd n\ptf
\]
D'après le corollaire~\ref{deltaQ}, les $\delta\in P(F)\cap\go$ qui donnent une contribution
non nulle appartiennent à $\tQp(F)\cap\go$, pourvu que $T$ soit assez r\'egulier.
Il suffit donc d'établir la convergence de l'intégrale
\[
\int_{\YQ}\FPQ(x,T){\tsQR}(\xT)
\Biggl\lvert\sum_{\{\tP\mid \tQp\subset\tP\subset\tRm\}}(-1)^{a_{\tP}-a_{\tG}}\,\Phi_{\tP,\go}(x)\,\Biggr\rvert\dd x
\]
avec
\[
\Phi_{\tP,\go}(x)=\int_{\N(F)\bs\N(\adef)}
\sum_{\delta\in\tQp(F)\cap\go}\,\fff^1(x\moins\,\delta\,n\,x) \dd n
\]
soit encore
\[
\Phi_{\tP,\go}(x)=\sum_{\eta\in\tMQp(F)\cap\go}\Phi_{\tP,\eta,\go}(x)
\]
avec
\[
\Phi_{\tP,\eta,\go}(x)=
\int_{\N(F)\bs\N(\adef)}\sum_{\nu\in \NQp(F)}
\,\fff^1(x\moins\,\eta\nu\,n\,x) \dd n
\]
où $\tMQp$ est le sous-ensemble de Levi de $\tQp$ contenant $\tM_0$ et
$\NQp$ son radical unipotent. On observe que dans cette expression seule l'intégrale
sur $\N(F)\bs\N(\adef)$ dépend de $\tP$.
Posons
\[
\Xi_\Q^\R(x)=\tsQR(\xT)\sum_{\eta\in\tMQp(F)\cap\go}
\Biggl\lvert\sum_{\{\tP\mid \Q\subset P\subset\R\}}(-1)^{a_{\tP}-a_{\tG}}\,
\Phi_{\tP,\eta,\go}(x)\Biggr\rvert\ptf
\]
Nous devons montrer que l'intégrale
\[
\int_{\YQ}\FPQ(x,T)\,\Xi_\Q^\R(x)\dd x
\]
est convergente. Nous allons tout d'abord d'estimer l'intégrale
\begin{equation}
\Theta_\Q^\R(n,x)=\int_{\gA_\G\bs\gA_\Q}\Xi_\Q^\R(nax)\,\delta_\Q(a)\moins\dd a
\tag{$*$}\label{eq9.*}
\end{equation}
où $\delta_\Q$ est le module pour $\Q$,
de façon uniforme lorsque $x$ reste dans un compact fixe.
On observe que puisque $\fff$ est à support compact
la somme sur $\eta$ ne porte que sur un ensemble fini
(dépendant \emph{a priori} de $x$ et $a$).
L'homomorphisme de $\ga_\Q^\G$ dans $\gao$
\[
H\mapsto \theta(H)-H
\]
a pour noyau le sous-espace des $\theta$\hyph invariants,
qui d'après le lemme~\ref{qplusa}
est l'espace
\[
\tga_\Q^\G=\ga_{\tQp}^{\tG}
\]
et on notera $\gb_\Q^\G$,
le supplémentaire orthogonal de ce sous-espace
dans $\ga_\Q^\G$; on a donc une décomposition
en sous-espaces deux à deux orthogonaux:
\[
\ga_\Q^\G=\tga_\Q^\G\oplus\gb_\Q^\G\ptf
\]
On observe que si $\eta\in\tMQp(F)$ et $a=\ee^H$
avec $H\in\ga_\Q^\G$ on a
\[
a\moins\,n_1\eta\,n a=n_1'a\moins\theta(a)\eta n'
\]
et si cette expression reste dans un compact il en est nécessairement
de même pour $a\moins\theta(a)$. Pour le voir on décompose $a$
en $a_0a_1a_2$ avec
\[
a_0\in\tga_\Q^\G=\ga_{\tQp}^{\tG},\qquad
a_1\in\gb_\Q^{\Qp}\Qquad{et}a_2\in\gb_{\Qp}^\G
\]
et on observe que
\[
a\moins\theta(a)=a_1\moins\theta(a_1)\ldot a_2\moins\theto(a_2)\ptf
\]
Maintenant si la décomposition de Bruhat de $\eta$ s'écrit $\eta=\xi w_s \xi' \mu$
on voit tout d'abord que
\[
(s-1)\HO(a_1)+(\theto-1)\HO(a_2)+\HO(w_s n'')
\]
reste borné et comme
\[
(s-1)\HO(a_1)+\HO(w_s n'')\in\gao^{\Qp}\Qquad{et}(1-\theto)\HO(a_2)\in\ga_{\Qp}^\G
\]
on en déduit que que
\[
(s-1)\HO(a_1)+\HO(w_s n'')\Qquad{et}(\theto-1)\HO(a_2)
\]
restent aussi bornés.
On remarque que l'application linéaire $(\theta-1)$ est injective sur $\gb_\Q^\G$.
En particulier $a_2$ reste dans un compact.
On rappelle que par hypothèse on a
\[
\tau_\Q^{\Qp}(\HO(a_1)-\T_1)=\tau_\Q^{\Qp}(\HO(a)-\T)=1
\]
où $\T_1$ est la projection de $\T$ sur $\ga_\Q^{\Qp}$.
Posons $X=\HO(a_1)-\T_1$. On a donc
\[
X=\sum_{\alpha\in\Delta_\Q^{\Qp}} a_\alpha \vpi_\alpha
\qquad\text{avec $a_\alpha>0$}\ptf
\]
Maintenant (comme dans le lemme~\ref{geofini}) on observe que le lemme~\ref{wn} implique que
\[
\langle X,\HO(w_sn'')\rangle\le c
\]
d'où on déduit
qu'il existe une constante $C$ telle que
\[
\langle X,(1-s)X\rangle \le C\ptf
\]
Il résulte alors du lemme~\ref{bigron} que $X$ et donc $a_1$ restent dans un compact.
Pour de tels $a$ l'ensemble
des $\eta$ intervenant est contenu dans un compact.
En particulier l'intégrale en $b=\ee^H$
avec $H\in\gb_\Q^\G$ porte sur un compact.
Il nous reste à estimer l'intégrale en $a_0=\ee^H$ avec
\[
H\in\tga_\Q^\G=\ga_{\tQp}^{\tRm}\oplus\ga_{\tRm}^{\tG}\ptf
\]
Maintenant, compte tenu de l'équation~\eqref{eq2.ii} du lemme~\ref{croiss} ,
il suffit de considérer les
\[
H\in\ga_{\tQp}^{\tRm}
\]
qui vérifient de plus
\[
\alpha(H)>\alpha(T)-C\qquad\forall\alpha\in\Delta_\Q^\R
\]
où $C$ dépend du support de $\fff$.
Nous sommes ainsi essentiellement ramenés à la situation traitée par
Arthur dans \cite{ATFI}. Rappelons en les étapes.
Notons $\gnQ$ l'algèbre de Lie de $\NQp$ et
$\widehat\gnQ$ son dual. L'application exponentielle définit une bijection
entre $\gnQ(F)$ et $\NQp(F)$. Soit $\psi$ un caractère
non trivial de $\adef/F$. Posons
\[
\hff(x,Y,\delta,n)=\int_{\gnQ(\adef)}\psi(\langle X,Y\rangle)
\fff^1(x\moins\,\delta\,\ee^X\,n\,x)\dd X\ptf
\]
On note $\gnQ_\bot$ l'orthogonal
de $\gnQ(F)$ \cad l'ensemble des $Y\in \widehat\gnQ(\adef)$ tels que
\[
\psi(\langle X,Y\rangle)= 1\qquad\forall X\in\gnQ(F)\ptf
\]
La formule de Poisson montre que
\[
\sum_{X\in \gnQ(F)} \fff^1(x\moins\,\eta\,\ee^X\,n\,z\,x)=\sum_{Y\in\gnQ_\bot}\hff(x,Y,\eta,n)
\]
et
\[
\int_{\N(F)\bs\N(\adef)}\sum_{X\in\gnQ(F)}\fff^1(x\moins\,\eta\,\ee^X\,n\,x)\dd n
=\sum_{Y\in\gnQ_\bot(P)}\hff(x,Y,\eta,1)
\]
où cette fois la somme porte sur le sous-ensemble $\gnQ_\bot(P)$
des éléments de $\gnQ_\bot$ qui sont triviaux sur $\gn(\adef)$.
On fait maintenant intervenir la somme alternée sur
les sous-ensembles paraboliques $\tP$ et compte tenu du lemme~\ref{binome} on obtient que
\[
\sum_{\{\tP\mid \Q\subset P\subset\R\}}(-1)^{a_{\tP}-a_{\tG}}\,
\Phi_{\tP,\eta,\go}(x)=\sum_{Y\in\gnQ_\bot(\Q,\R)}\hff(x,Y,\eta,1)
\]
où $\gnQ_\bot(\Q,\R)$ est le sous-ensemble des $Y\in\gnQ_\bot$ ayant la propriété que $Y\in\gnQ_\bot(P)$
pour un seul sous-ensemble parabolique $\tP$
tel que
\[
\Q\subset P\subset\R\ptf
\]
Comme
\[
\gnQ_\bot(P)\subset\gnQ_\bot(\Rm)
\]
on a donc
\[
\gnQ_\bot(\Q,\R)=\gnQ_\bot(\Rm)-\bigcup_{\substack{\tP\subset\tRm\\ P\ne\Rm}}\gnQ_\bot(P)
\]
et
\[
\Xi_\Q^\R(x)=\tsQR(\xT)\sum_{\eta\in\tMQp(F)\cap\go}
\Biggl\lvert\sum_{Y\in\gnQ_\bot(\Q,\R)}\hff(x,Y,\eta,1)\Biggr\rvert\ptf
\]
On rappelle que l'on souhaite estimer l'intégrale~\eqref{eq9.*} définissant $\Theta_\Q^\R$:
\[
\Theta_\Q^\R(n,x)=\int_{\gA_\G\bs\gA_\Q}\Xi_\Q^\R(nax)\,\delta_\Q(a)\moins\dd a
\ptf
\]
Considérons donc,
avec les notations du lemme~\ref{qplusa},
les $a=\ee^H$ pour des $H\in\tgaQR$ qui vérifient de plus
\[
\alpha(H)>\alpha(T)\qquad\forall\alpha\in\Delta_\Q^\R\ptf
\]
On va voir que, si l'espace $\tgaQR$ n'est pas réduit à zéro,
de tels $a$ agissent par dilatation non triviale
sur au moins une des coordonnées de chaque
élément de $\gnQ_\bot(\Q,\R)$.
Pour cela on décompose $\gnQ_\bot$ suivant les caractères
de l'action coadjointe de $\tgaQR=\ga_{\tQp}^{\tRm}$:
\[
\gnQ_\bot=\bigoplus_\Qbeta \gnQ_\bot(\Qbeta)\ptf
\]
L'ensemble de ces caractères peut s'identifier avec l'ensemble des restrictions à $\ga_{\tQp}^{\tRm}$ des
racines de $\gao$ dans $\gnQ$.
Un élément $Y\in\gnQ_\bot$ peut donc s'écrire
\[
Y=\sum Y_\Qbeta
\]
et on a $Y\in\gnQ_\bot(\Q,\R)$ seulement si pour tout $\alpha\in\Delta_{\Qp}^{\Rm}$
il existe $\Qbeta$ avec
\[
Y_\Qbeta\ne0\Qquad{et}\langle \Qbeta,\tvpi_\alpha\rangle\ne0\ptf
\]
En effet, s'il existait $\alpha$ tel que pour tout $\Qbeta$ avec $Y_\Qbeta\ne0$ on ait
$\langle\Qbeta,\tvpi_\alpha\rangle=0$, alors $Y$ appartiendrait à
$\gnQ_\bot(\Rm)$ et à $\gnQ_\bot(P)$ où $P$ est le sous-groupe parabolique
de $\Rm$ qui admet pour sous-groupe de Levi le centralisateur de $\tvpi_\alpha$, ce qui est exclu.
Par ailleurs les $\tvpi_\alpha$ avec
$\alpha\in\Delta_{\Qp}^{\Rm}$ forment une base du dual de $\ga_{\tQp}^{\tRm}$ et
\[
\Ad(a)Y=\sum_\Qbeta \ee^{\Qbeta(H)}Y_\Qbeta=\sum_\Qbeta\prod_{\talpha\in\Delta_{\tQp}^{\tRm}}\ee^{h_{\talpha}\langle \Qbeta,\tvpi_\alpha\rangle}\,Y_\Qbeta\ptf
\]
On observe de plus que
\[
\delta_\Q(a)\moins\hff(ax,Y,\delta,1)=\delta_\Q(a)\moins\delta_{\Qp}(a)\,\hff(x,\Ad(a)Y,\delta,1)\ptf
\]
Pour $a=\ee^H$ avec $H\in\ga_{\Qp}$ on a $\delta_\Q(a)=\delta_{\Qp}(a)$.
Maintenant $\hff$ est une fonction lisse à décroissance rapide en $Y$ comme transformée de
Fourier d'une fonction lisse à support compact.
Il en résulte que l'intégrale définissant $\Theta_\Q^\R$
est absolument convergente, uniformément lorsque $x$
reste dans un compact. En utilisant la décomposition
d'Iwasawa on voit que l'intégrale
\[
\int_{\YQ}\FPQ(x,T)\,\Xi_\Q^\R(x)\dd x
\]
est égale à
\[
\int_\K\int_{\gA_\Q\M_\Q(F)\bs\M_Q(\adef)}
\int_{\N_\Q(F)\bs\N_Q(\adef)}\FPQ(m,T)\Theta_\Q^\R(n,mk)\dd n\dd m\dd k
\ptf
\]
Il nous reste à observer que, compte tenu du lemme~\ref{Gammacar} et du théorème~\ref{siegel},
l'intégrale en $m$ porte sur un compact et que donc la somme sur $\eta$
dans la définition de $\Xi_\Q^\R$ ne porte que sur un ensemble fini
qui peut être choisi indépendant de $m$, $n$ et $a$.
\end{proof}

\begin{theoreme}\label{geoconv} 
Supposons $T$ assez régulier
\cad $\dPO(T)\ge c$ où $c$ est une constante ne dépendant  que du support de $\fff$.
L'expression
\[
\sum_{\go\in\gO}\int_{\XG} \lvert k^\T_\go(x)\rvert\dd x
\]
est convergente. Seul un ensemble fini de $\go$ fournit une
contribution non nulle \textup(cet ensemble dépend du support de $\fff$\textup).
\end{theoreme}

\begin{proof}
La finitude résulte de ce que, d'après le lemme~\ref{ssfini}, les fonctions
\[
x\mapsto K_{\tP,\go}(x,x)
\]
ne sont non identiquement nulles que pour un ensemble fini de $\go$.
Il suffit donc de considérer une classe $\go$ et de prouver la convergence de
\[
\int_{\XG} \lvert k^\T_\go(x)\rvert\dd x\ptf
\]
On utilise la partition du lemme~\ref{FPQ}
\[
\sum_{\{\Q\mid \PO\subset\Q\subset P\}}\,
\sum_{\xi\in\Q(F)\bs P(F)}\FPQ(\xix,T)\tau_\Q^P(\xixT)=1
\]
et on obtient
\[
k^\T_\go(x)=\sum_{\{\tP,\Q\mid \PO\subset\Q\subset P\}}\,
\sum_{\xi\in\Q(F)\bs\G(F)}(-1)^{a_{\tP}-a_{\tG}}\,k^\T_{\tP,\Q,\go}(\xix)
\]
avec
\[
k^\T_{\tP,\Q,\go}(x)=\FPQ(x,T)\tau_\Q^P(\xT)
\htautP(\xT)\,K_{\tP,\go}(x,x)\ptf
\]
On rappelle que
\[
\YQ=\gA_\G\Q(F)\bs\Gadef \ptf
\]
On va montrer que pour tout $\Q$ l'intégrale
\[
\int_{\YQ}\Biggl\lvert
\sum_{\{\tP\mid \Q\subset P\}}(-1)^{a_{\tP}-a_{\tG}}\,k^\T_{\tP,\Q,\go}(x)\Biggr\rvert
\dd x
\]
est convergente.
On rappelle que d'après le lemme~\ref{repart} on a
\[
\sum_{\R\supset P} \tsQR=\tau_\Q^P\htautP\ptf
\]
La convergence souhaitée est donc conséquence de la convergence des intégrales
\[
\int_{\YQ}\FPQ(x,T)\tsQR(\xT)\Biggl\lvert
\sum_{\{\tP\mid \Q\subset P\subset\R\}}(-1)^{a_{\tP}-a_{\tG}}\,
K_{\tP,\go}(x,x)\Biggr\rvert \dd x
\]
qui a été établie dans la proposition~\ref{geoconvb}.
\end{proof}

\section{Dévelopement géométrique grossier}

Soit $\go$ une classe de ss\hyph conjugaison.
On rappelle que l'on a prouvé la convergence de l'intégrale sur $\XG$ de
\[
k^\T_\go(x)=\sum_{\tP\supset\tPO}(-1)^{a_{\tP}-a_{\tG}}\,
\sum_{\xi\in P(F)\bs\G(F)}\htautP(\xixT)\,K_{\tP,\go}(\xix,\xix)\ptf
\]
avec
\[
K_{\tP,\go}(x,x)=\int_{\N(F)\bs\N(\adef)}
\sum_{\delta\in\go\cap\tP(F)}\omega(x)
\fff^1(x\moins\delta\, n\,x)\dd n\ptf
\]
On aura besoin d'une variante de l'expression
$k^\T_\go(x)$ de même intégrale sur $\XG$: on pose
\[
j^\T_\go(x)=\sum_{\tP\supset\tPO}(-1)^{a_{\tP}-a_{\tG}}\,\,
\sum_{\xi\in P(F)\bs\G(F)}\htautP(\xixT)\,j_{\tP,\go}(\xix)
\]
avec
\[
j_{\tP,\go}(x)=\sum_{\delta\in\go\cap\tP(F)}\int_{\N(\delta,F)\bs\N(\delta,\adef)}\omega(x)
\fff^1(x\moins\delta\, n\,x)\dd n\ptf
\]
où $\N(\delta)=\N(\delta_s)$ est le sous-groupe de $N$ qui centralise la partie semi-simple
$\delta_s$ de $\delta$.

\begin{lemme}\label{surjnil}
On pose $\theta=\AdtG(\delta)$. L'application
de
\[
\N\times\N(\delta)\to\N
\]
définie par
\[
n\times n'\mapsto n\moins\,n'\,\theta(n)
\]
est surjective. L'image réciproque
du sous-ensemble $n\moins\,\N(\delta)\,\theta(n)$
est
\[
\N(\delta)\ldot n\times N(\delta)\ptf
\]
\end{lemme}

\begin{proof}
Soit $\theta=\theta_s\theta_u$ la décomposition de Jordan de $\theta$.
On peut munir le groupe nilpotent $\N$ d'une filtration par des sous-groupes normaux
stables sous $\theta$:
\[
\N=\N_0\supset\N_1\dotsm\supset\N_r=\{1\}
\]
tels que $n_1\moins\,n\moins\,n_1\,n\in\N_{i+1} $ pour $n\in\N$ et $n_1\in\N_i$
et $\theta_u(n)n\moins \in\N_{i+1} $ pour $n\in\N_i$.
Posons
\[
S_k=\N_k\ldot \N(\delta)=\N(\delta).\N_k
\ptf
\]
Nous allons montrer, par récurrence descendante
sur $k$, que l'application
de
\[
\N_k\times\N(\delta)\to S_k
\]
définie par
\[
n\times n'\mapsto n\moins\,n'\,\theta(n)
\]
est surjective avec $(\N(\delta)\cap N_k)\ldot n\times\N(\delta)$ comme image réciproque dans $S_k\times\N(\delta)$ au dessus de
$n\moins\,\N(\delta)\,\theta(n)$.
Le lemme est le cas particulier $k=0$. L'assertion est claire pour $k=r$. Supposons la vraie pour $k+1$;
il en résulte que l'ensemble des
$n\moins n'\,\theta(n)$
avec $n\in\N_k$ et $n'\in\N(\delta)$ est aussi égal à l'ensemble des
$n\moins n''\,\theta(n)$ où cette fois on prend $n''\in S_{k+1}$.
Mais $S_{k+1}$ est normal dans $S_{k}$ avec pour quotient
un groupe abélien
muni d'une structure d'espace vectoriel
et où $\theta$ agit par un endomorphisme
semi-simple qui n'admet pas la valeur propre~1.
L'assertion en résulte.
\end{proof}

\begin{lemme}\label{intsurjnil}
Soit $\phi$ une fonction sur $\tP(\adef)$.
\[
\int_{\N(F)\bs\N(\adef)}
\int_{\N(\delta,F)\bs\N(\delta,\adef)}\sum_{\delta\in\go\cap\tP(F)}
\phi( n\moins\delta\,n_1\,n)\dd n_1\dd n
\]
est égal à
\[
\int_{\N(F)\bs\N(\adef)} \sum_{\delta\in\go\cap\tP(F)}\phi(\delta\,n)\dd n\ptf
\]
\end{lemme}

\begin{proof}
C'est une conséquence facile du lemme~\ref{surjnil}.
\end{proof}

\begin{proposition}
\[
\int_{\XG} k^\T_\go(x)\dd x=\int_{\XG} j^\T_\go(x)\dd x\ptf
\]
\end{proposition}

\begin{proof}
La convergence de l'intégrale dans le membre de droite s'établit comme au théorème~\ref{geoconv}.
Maintenant on a l'égalité
\[
K_{\tP,\go}(x,x)=\int_{\N(F)\bs\N(\adef)} j_{\tP,\go}(nx)\dd n
\]
qui est la conséquence du lemme~\ref{intsurjnil}. La proposition en résulte.
\end{proof}

\begin{theoreme}\label{geoconvc} 
Supposons $T$ assez régulier.
L'expression
\[
\sum_{\go\in\gO}\int_{\XG} \lvert j^\T_\go(x)\rvert \dd x
\]
est convergente.
\end{theoreme}

\begin{proof}
Les arguments sont une reprise, presque mot à mot, de ceux de la preuve du théorème~\ref{geoconv}.
Nous les laissons en exercice pour le lecteur.
(Voir aussi \cite{ATFI}).
\end{proof}

Ce sont les intégrales
\[
\int_{\XG} j^\T_\go(x)\dd x
\]
qui donnent naissance au développement géométrique fin et
qui permettent en particulier d'obtenir des expressions explicites pour les contributions
des classes quasi semi-simples (\cf section~\ref{explicitss}).

\section{Termes quasi semi-simples}\label{explicitss}

Soit $\delta$ un élément quasi semi-simple dans $\tG(F)$.
Soit $\centd$ le centralisateur stable de $\delta$ (\cf section~\ref{ellip} et
la remarque~\ref{idgd}).
On note $\gc$ la classe de conjugaison de $\delta$.
Considérons un tore déployé maximal $S_\delta$ dans
le centre de $\centd$ (ou, ce qui revient au même, dans la composante neutre $\G_\delta$ du centralisateur de $\delta$).
Le centralisateur de $S_\delta$ est un sous-groupe de Levi
$\M_\delta$ et on pose $\tM_\delta=\M_\delta.\delta$.
On observe que $\centd\subset\M_\delta$ et que
$\delta$ est elliptique dans $\tM_\delta$.
Le centralisateur $\G^\delta$ normalise $\M_\delta$.
Ceci fournit (dans les notations du lemme~\ref{thetastab})
une application
\[
\G^\delta(F) \to \weyl(\ga_{\tM_\delta},\ga_{\tM_\delta})
\]
de noyau $\M_\delta\cap\G^\delta(F)$.
A conjugaison près on peut supposer
que $\tM_{\delta}$ est un sous-ensemble de Levi standard.
Soit $\tP=\tM\,\N$ un sous-ensemble parabolique standard et
supposons que $\tM$ contient un conjugué $\tM_1$ de $\tM_\delta$.
Considérons
\[
\delta_1\in\gc\cap\tM_1(F)\ptf
\]
On observe qu'alors $\cent_{\delta_1}\subset\M_1\subset\M$ et donc
le sous-groupe $N(\delta_1)$ est trivial. A conjugaison près dans $\M$
on peut supposer $\M_1$ standard.
Dans ces conditions il existe
\[
s\in\weyl(\ga_{\tM_\delta},\ga_{\tM_1})
\]
de représentant $w_s$ tel que $w_s\delta w_s\moins=\delta_1$.
Comme dans le lemme~\ref{weylgab}, mais dans le cas tordu,
on note
\[
\weyl(\ga_{\tM_\delta},\tP)
\]
l'union (disjointe) des quotients
\[
\weyl^P(\ga_{\tM_1},\ga_{\tM_1})\bs\weyl(\ga_{\tM_\delta},\ga_{\tM_1})
\]
où $\tM_1$ parcourt les sous-ensembles de Levi standard de $\tM$
à conjugaison près par $\M$. On introduit
\[
j_{\tP,\gc}(x)=\iota(\delta)\moins\sum_{s\in\weyl(\ga_{\tM_\delta},\tP)}
\,\,\sum_{\eta\in \centdi\bs P(F)}
\omega(x)\fff^1(x\moins\eta \moins\,\delti\,\eta\,x)
\]
où $\delti=w_s\delta\, w_s\moins$
et
\[
\iota(\delta)=\#\centd(F)\bs\G^\delta(F)\ptf
\]
On pose
\[
j^\T_{\gc}(x)=\sum_{\tP\supset \tPO}(-1)^{a_{\tP}-a_{\tG}}\,
\sum_{\xi\in P(F)\bs\G(F)}\htautP(\xixT)\,j_{\tP,\gc}(\xix)\ptf
\]
On définit de manière analogue $k^\T_\gc$.
On va donner une expression pour
\[
\int_{\XG} k^\T_\gc(x)\dd x=\int_{\XG} j^\T_\gc(x)\dd x
\]
au moyen d'intégrales orbitales pondérées.
Pour cela nous aurons besoin d'introduire les objets suivants.
On a considéré au lemme~\ref{wn} la famille orthogonale $\YY(x,\T)$ définie par les:
\[
Y_s(x,\T)=s\moins\bigl(T-\HO(w_s x)\bigr)\qquad\text{pour $s\in\weyl$}\ptf
\]
On pose (\cf proposition~\ref{envconvt})
\[
v_{\tM_\delta}^\T(x)=\int_{\tga_{\M_\delta}^\G}\Gamma_{\tM_\delta}^{\tG}\bigl(H,\YY(x,\T)\bigr)\dd H\ptf
\]
Si $\T$ est assez régulier
\[
H\mapsto \Gamma_{\tM_\delta}^{\tG}\bigl(H,\YY(x,\T)\bigr)
\]
est la fonction caractéristique
de l'enveloppe convexe des projections sur $\tga_{\M_\delta}^\G$ des
$Y_s(x,\T)$ pour $s\in\weyl(\tga_{\M_\delta})$.
Dans ce cas
\[
v_{\tM_\delta}^\T(x)
\]
est le volume de cette enveloppe convexe.
C'est un polynôme en $\T$, ne dépendant que de la projection de $\T$
sur le sous espace des vecteurs $\theto$\hyph invariants, de degré
$(a_{\tM_\delta}-a_{\tG})$. On rappelle enfin que dans la section~\ref{problemat} on a
associé à $\fff$ une autre fonction $\tff$.

\begin{proposition}
L'intégrale
\[
\int_{\XG} k^\T_{\gc}(x)\dd x=\int_{\XG} j^\T_{\gc}(x)\dd x
\]
est égale à
\[
\frac{\vol\bigl(\gA_{\centd}\centd(F)\bs\centd(\adef)\bigr)}{\iota(\delta)}\int_{\centd(\adef)\bs\Gadef}\omega(x)
v_{\tM_\delta}^\T(x)\tff(x\moins\delta\,x)\dd x
\]
si $\centd(\adef)$ est dans le noyau de $\omega$ et zéro sinon.
\end{proposition}

\begin{proof}
En effet
\[
j^\T_{\gc}(x)=\iota(\delta)\moins\omega(x)\sum_{\xi\in \centd(F)\bs\G(F)}
e_{\tM_\delta}(\xix,T)\fff^1(x\moins\xi\moins\delta\,\xix)
\]
où
\[
e_{\tM_\delta}(x,T)=
\sum_{\tP\supset\tPO}(-1)^{a_{\tP}-a_{\tG}}\,\sum_{s\in\weyl(\ga_{\tM_\delta},\tP)}\htautP(\HO(w_s x)-T)
\]
soit encore, avec les notations de la proposition~\ref{preconv} (ou plutôt de sa variante tordue),
\[
e_{\tM_\delta}(x,T)=\sum_{s\in\weyl(\ga_{\tM_\delta})}\,
\sum_{s\moins(\tP)\in\FF_s(\tM_\delta)}(-1)^{a_{\tP}-a_{\tG}}\,\htautP(\HO(w_s x)-T)
\]
et donc
\[
e_{\tM_\delta}(x,T)=\sum_{s\in\weyl(\ga_{\tM_\delta})}\,
\sum_{\tQ\in\FF_s(\tM_\delta)}(-1)^{a_{\tQ}-a_{\tG}}\,\htau_{\tQ}\bigl(s\moins(\HO(w_s x)-T)\bigr)\ptf
\]
On obtient la formule
\begin{equation}
\int_{\XG} j^\T_{\gc}(x)\dd x=\int_{\gAtG\centd(F)\bs\Gadef}\iota(\delta)\moins\omega(x)
e_{\tM_\delta}(x,T)\tff(x\moins\delta\,x)\dd x\ptf\tag{$*$}\label{eq9.*a}
\end{equation}
On observe que, si $a=\exp H$ avec $H\in\gao$
alors, avec les notations du lemme~\ref{wn},
\[
\htau_{\tQ}\bigl(s\moins(\HO(w_s ax)-T)\bigr)=\htau_{\tQ}\bigl(H-Y_s(x,\T)\bigr)
\]
et on déduit de la proposition~\ref{envconvt} que
\[
 e_{\tM_\delta}(a \,x,\T)=\Gamma_{\tM_\delta}^{\tG}\bigl(H,\YY(x,\T)\bigr)\ptf
\]
En observant que
\[
\ga_{\centd}\simeq\tga_{\M_\delta}\coloneqq \ga_{\tM_\delta}
\]
la formule~\eqref{eq9.*a} se récrit
\[
\int_{\XG} j^\T_{\gc}(x)\dd x=\int_{\gA_{\centd}\centd(F)\bs\Gadef}\iota(\delta)\moins\omega(x)
v_{\tM_\delta}^\T(x)\tff(x\moins\delta\,x)\dd x\ptf
\]
On observe enfin que cette expression est nulle si $\centd(\adef)$ n'est pas dans le noyau de $\omega$.
Sinon on obtient la formule souhaitée en décomposant l'intégrale par passage au quotient par $\centd(\adef)$.
\end{proof}

\begin{remarque}
Le lecteur observera que nous avons muni les espaces vectoriels isomorphes
$\ga_{\tM_\delta}$ et $\ga_{\centd}$ de la même mesure de Haar. Toutefois,
si on souhaite utiliser les mesures de Tamagawa, les mesures de Haar canoniques
qui leur sont associées sur ces espaces vectoriels, seront en général différentes.
Cette remarque joue un rôle dans l'étude de la stabilisation (voir par exemple \cite{Ltw}).
\end{remarque}

\section{Développement géométrique fin}

Le cas des $\delta$ non quasi semi-simples suppose le traitement préalable
des contributions unipotentes; le cas général s'en déduit
pas descente au centralisateur. Ceci a fait l'objet de deux articles d'Arthur:
\cite{Aunip} et \cite{Ageom} (qui eux-mêmes reposent sur \cite{ALB}
publié ultérieurement)
où il étudie les termes géométriques, y compris
dans le cas tordu en s'appuyant sur \cite{MS}.

Dans ces articles, l'espace tordu $\tG$ (resp.~$\G$ dans la notation d'Arthur)
est une composante d'un groupe réductif non connexe $\G^+$
de composante neutre $\G$ (resp.~$\G^0$); ceci revient à demander que
l'automorphisme $\theto$ ait une puissance $\theta_0^\ell$ qui est un automorphisme intérieur
représenté par un élément rationnel:
\[
\theta_0^\ell=\Ad_\G(x)\qquad\text{avec $x\in G(F)$}
\]
et de plus Arthur ne considère pas de caractère $\omega$ non trivial.

Tout ceci est légèrement restrictif par rapport
à notre cadre, mais cela est sans conséquence sérieuse sur les preuves.
Le développement géométrique fin est donné dans \cite{Ageom}.
Nous n'avons rien à ajouter à ces résultats (sauf l'introduction d'un caractère $\omega$)
et nous renvoyons le lecteur à ces articles.
\chapter{Développement spectral grossier}\label{ch10}

Dans ce chapitre, ainsi d'ailleurs que dans toute la suite de ce livre, 
les fonctions $\fff$ sont suppos\'ees $\K$\hyph finies (\`a droite et \`a gauche).

\section{Convergence: côté spectral}

Rappelons
que le membre de droite de l'identité fondamentale~\ref{fonfon} est
\[
k^\T_{\mathrm{spec}}(x) = \sum_{\tP}(-1)^{a_{\tP}-a_{\tG}}\sum_{\Q \subset P\subset \R }\,
\sum_{\xi\in\Q(F) \bs\G(F)}\tsQR(\xixT) \cdot \tronc_1^{T,\Q}K_{\tP} (\xix,\xix)\ptf
\]
où
\[
K_{\tP}(x,y) = \int_{\N(F)\bs \Nadef} \sum_{\delta\in \tP(F)} \omega(y) \fff^1 (x\moins\,\delta\,n\, y)\dd n
\]
soit encore
\[
K_{\tP}(x,y) = \int_{ \Nadef} \sum_{\delta\in \tM(F)} \omega(y) \fff^1 (x\moins\,n\moins\delta\, y)\dd n\ptf
\]
Posons
\[
k_{\mathrm{spec}}^\T(\Q,\R,x)=\tsQR(\xT)\sum_{\{\tP\mid \tQp \subset P\subset \tRm\}}(-1)^{a_{\tP}-a_{\tG}}\,\tronc_1^{T,\Q}K_{\tP} (x,x)\ptf
\]
On a donc
\[
k^\T_{\mathrm{spec}}(x) = \sum_{\Q \subset \R }\sum_{\xi\in\Q(F) \bs\G(F)}k_{\mathrm{spec}}^\T(\Q,\R,\xix)\ptf
\]
On pose
\[
\tve(\Q,\R)=(-1)^{a_{\tRm}-a_{\tG}}
\]
\newnot{epsilon(Q,R)@$\string\tve(Q,R)$}{tveqr}%
si $\Qp\subset\Rm$,
et $0$ sinon. Notons $\tG(\Q,\R)$
l'ensemble des $\delta$ qui appartiennent à $\tP(F)$ pour un seul $\tP$
avec $\tQp\subset\tP\subset\tRm$ \cad $\delta\in\tRm(F)$ mais $\delta\notin\tP(F)$ si
\[
\tQp\subset\tP\varsubsetneq\tRm\ptf
\]
Posons
\[
K_{\Q,\R}(x,y)=\int_{\N_\Q(F)\bs \N_\Q(\adef)}
\sum_{\delta\in \tG(\Q,\R)} \omega(y) \fff^1 (x\moins\,\,n_\Q\moins\,\delta\, y)\dd n_\Q\ptf
\]

\begin{lemme}\label{prenul}
Avec ces notations on a
\[
k^\T_{\mathrm{spec}}(x) = \sum_{\Q \subset \R }\tve(\Q,\R)\sum_{\xi\in\Q(F) \bs\G(F)}
\tsQR(\xixT)\tronc_1^{T,\Q}K_{\Q,\R}(x,\xix)\ptf
\]
\end{lemme}

\begin{proof}
On observe que
\[
\tronc_1^{T,\Q}K_{\tP} (x,y)=\tronc_1^{T,\Q}\Pi_\Q K_{\tP} (x,y)
\]
où
\[
\Pi_\Q K_{\tP} (x,y)=\int_{\N_\Q(F)\bs \N_\Q(\adef)} K_{\tP} (nx,y)\dd n_\Q
\]
et comme $\N\subset\N_\Q$ on a
\[
\Pi_\Q K_{\tP} (x,y)=\int_{\N_\Q(F)\bs \N_\Q(\adef)}
\sum_{\delta\in \tP(F)} \omega(y) \fff^1 (x\moins \,n_\Q\moins\,\delta\, y)\dd n_\Q\ptf
\]
Mais
\[
\sum_{\{\tP\mid \Q \subset P\subset \R\}}(-1)^{a_{\tP}-a_{\tG}}\Pi_\Q K_{\tP} (x,y)
\]
est égal à
\[
\int_{\N_\Q(F)\bs \N_\Q(\adef)}
\sum_{\delta\in\tRm(F)} \omega(y) \sum_{\{\tP\mid\delta\in\tP(F),\Q \subset P\subset \R\} }
(-1)^{a_{\tP}-a_{\tG}}\,\fff^1 (x\moins\,\,n_\Q\moins\,\delta\, y)\dd n_\Q\ptf
\]
La somme alternée des termes contenant $\delta$ est nulle sauf
si $\delta$ appartient à $\tP(F)$ pour un seul $\tP$
avec $\Q\subset P\subset\R$ \cad si $\delta\in\tG(\Q,\R)$.
On a donc
\[
\sum_{\{\tP\mid \Q \subset P\subset \R\}}\
(-1)^{a_{\tP}-a_{\tG}}\,\Pi_\Q K_{\tP} (x,y)=\tve(\Q,\R)K_{\Q,\R}(x,y)
\]
ce qui fournit
\[
k_{\mathrm{spec}}^\T(\Q,\R,x)=\tve(\Q,\R)\,\tsQR(\xT)\tronc_1^{T,\Q}K_{\Q,\R}(x,x)\ptf\qedhere
\]
\end{proof}

Posons, pour $\delta\in\tG(F)$,
\[
K_{\Q,\delta}(x,y)=
\int_{\N_\Q(F)\bs \N_\Q(\adef)} \omega(x)
\sum_{\eta\in \Q(F)}
\fff^1 (x\moins\,\,n_\Q\moins\,\eta\,\delta\,y)\dd n_\Q
\]
soit encore
\[
K_{\Q,\delta}(x,y)=
\int_{\N_\Q(\adef)} \omega(x)
\sum_{\mu\in \M_\Q(F)}
\fff^1 (x\moins\,\,n_\Q\moins\,\mu\,\delta\,y)\dd n_\Q\ptf
\]

\begin{lemme}\label{invar}
\[
K_{\Q,\delta}(n_1 x,n_2y)=K_{\Q,\delta}(x,y)=K_{\Q,\delta}(\xix,y)
\]
si $n_1\in\N_\Q(F)$\textup, $n_2\in \delta\moins \N_\Q(F)\delta$ et $\xi\in\Q(F)$.
\end{lemme}

\begin{proof} 
Cela résulte immédiatement de la définition de $K_{\Q,\delta}$.
\end{proof}

\begin{lemme}\label{hoyx}
Si $K_{\Q,\delta}(x,y)\ne0$ alors il existe un compact $C\subset\gao$
et $\eta\in\Q(F)$ tels que
\[
\HO\bigl(\eta\,\theta(y)\bigr)-\HO(x)\in C
\]
si $\theta$ est l'automorphisme de $\G$ défini par $\delta$.
\end{lemme}

\begin{proof}
Puisque $\fff$ est à support compact on a
\[
x\moins n_\Q\,\eta\,\delta\,y\in C_1
\]
où $C_1$ est un compact. La décomposition d'Iwasawa
$x=nak$ de $x$ montre que
\[
a\moins n\moins \eta\,\delta\,y\in C_2
\]
où $C_2$ est encore compact. On a donc
\[
\tHO(a\moins n\moins\eta\,\delta\,y)=\tHO( \eta\,\theta(y)\delta)+\HO(a\moins)
\]
borné et $\HO(a\moins)=-\HO(x)$.
\end{proof}

\begin{lemme}\label{xifini}
La somme sur $\xi$ dans
\[
\sum_{\xi\in\Q_\delta(F)\bs\Q(F)}\tronc_1^{T,\Q }K_{\Q,\delta}(x, \xi\,y)
\]
porte sur un ensemble fini dont le cardinal est
majoré par $c \lvert x\rvert^A\lvert y\rvert^B$.
\end{lemme}

\begin{proof}
En effet on déduit du lemme~\ref{hoyx} qu'il existe un compact $C_{\Q,\delta}\subset\ga_\Q$,
dépendant du support de $\fff$, tel que les $\xi$ qui interviennent vérifient
\[
\HQ\bigl(\theta(\xi\,y)\bigr)-\HQ(x)\in C_{\Q,\delta}
\]
où $\theta$ est l'automorphisme induit par $\delta$
et $\HQ(g)$ la projection de $\HO(g)$ sur $\ga_\Q$. En particulier, on a
\[
\htau_\Q\bigl(\HO\bigl(\theta(\xi\,y)\bigr)-\HO(x)-T_1\bigr)=1
\]
pour un certain $T_1$ (défini par le compact $C_{\Q,\delta}$). Comme $\theta(\Q_\delta)\subset\Q$
il résulte alors du lemme~\ref{finitude} qu'on peut choisir $\xi'\in\Q(F)\theta(\xi)$ de sorte que
\[
\lvert \xi'\rvert \le c_1\lvert x\rvert^{A_1}\lvert y\rvert^{B_1}
\]
et on observe que l'application qui à $\xi\in\Q_\delta(F)\bs\Q(F)$ associe
la classe $\Q(F)\theto(\xi)$ est injective. On invoque enfin le lemme~\ref{proj}.
\end{proof}

\begin{lemme}\label{maja}
Considérons $\delta\in\Q(F)\bs\tG(\Q,\R)$\textup, $k,k_1\in\K$\textup, $n\in\N_\Q(\adef)$\textup,
$m$ et $m_1\in\M_\Q(\adef)$ avec
\[
\HQ(m)=\HQ(m_1)=0
\]
et $a\in\gA_\Q^\G$.
Supposons que\textup, pour $\xi\in\Q(F)$\textup, on ait
\[
K_{\Q,\delta}(m_1ak_1,\xi namk)\ne0\Qquad{et} \tsQR\bigl(\HO(a)\bigr)=1\ptf
\]
Alors il existe une constante $c>0$ \textup(ne dépendant que de $\fff$\textup) telle que
\[
\lVert \HO(a)\rVert \le c(1+\lVert \HO(m)\rVert)\ptf
\]
\end{lemme}

\begin{proof}
Rappelons que l'on a choisi dans la section~\ref{tordalg}
un élément $\delto\in\tMO(F)$
et on a noté $\theto$ l'automorphisme associé.
On remarque que, comme on peut modifier $\delta$ à gauche par un élément de $Q(F)$,
on peut supposer $\delta\xi$ choisi de la forme
\[
\delta\xi=w_s\eta=w_{\so}\delto\eta
\]
avec $\eta\in\NO(F)$ et
où $w_{\so}$ représente un élément
$\so$ du groupe de Weyl de $\M_{\Rm}$ tel que
\[
\so\moins\alpha>0\qquad\text{pour toute $\alpha\in\Delta_{\PO}^\Q$}\ptf
\]
Nous supposerons donc désormais que
\[
\delta=w_s=w_{\so}\delto\Qquad{et} \xi=\eta\in\NO(F)\ptf
\]
Notons $\theta$ l'automorphisme de $\G$ défini par $\delta$.
On observe que $\theta$ et $\theto$ ont la même restriction à $\ga_{\Rm}$.
On a supposé
\[
K_{\Q,\delta}(m_1ak_1,\eta namk)\ne0
\]
et donc
\[
k_1\moins\,m_1\moins\,a\moins\,n_\Q\moins\mu\,\delta\,\eta\,n\,\,a\,m\,k\in\Support(\fff)
\]
ce qui implique que
\[
m_1\moins\,a\moins\,n_\Q\moins\mu\,\delta\,\eta\,n\,a\,m\in\Omega
\]
où $\Omega$ est un compact.
On décompose $H=\HO(a)$ en
\[
H=H_1\oplus H_2\oplus H_3
\]
au moyen de la décomposition en somme directe
\[
\ga_\Q^\G=\ga_\Q^{\Rm}\oplus \gb_{\Rm}^\G\oplus \ga_{\tRm}^{\tG}
\]
où $\gb_{\Rm}^\G$ est l'orthogonal de $\ga_{\tRm}^{\tG}$, le sous-espace 
des $\theto$\hyph invariants, dans $\ga_{\Rm}^\G$.
Posons $a_i=\ee^{H_i}$. On a
\[
m_1\moins\,a_3\moins\,a_2\moins\,a_1\moins\,n_\Q\moins\,\mu\,\delta\,\eta\,n\,a_1a_2 a_3\,m=
m_1\moins\,n_1\moins\mu_1\,a_1\moins\theta(a_1)\,b\,\delta\,n'\,m
\]
où
\[
b=\ee^{B}\qquad\text{avec $B=(\theto-1)H_2$}\ptf
\]
Comme par hypothèse $\HQ(m)=\HQ(m_1)=0$ on a
\[
\HH_{\Rm}(m_1)=\HH_{\Rm}\bigl(\theto(m)\bigr)=0
\]
ce qui implique
\[
\HH_{\Rm}\bigl(m_1\moins\,a\moins\,n_\Q\moins\,\mu\,\theta(\eta\,n\,a\,m)\bigr)=\HH_{\Rm}(b)=B
\]
et on en déduit que $B$ reste borné. L'homomorphisme:
\[
H_2\mapsto B=(\theto-1)H_2
\]
est injectif et donc $H_2$ reste aussi borné.
Maintenant, compte tenu de l'équation~\eqref{eq2.ii} du lemme~\ref{croiss}, on contrôle $H_3$
au moyen de $H_1+H_2$ et donc :
\[
\lVert H_3\rVert \ll (1+\lVert H_1\rVert)\ptf
\]
Pour conclure il reste à montrer qu'il existe $c'_0>0$ avec
\begin{equation}
\lVert H_1\rVert \ll ( 1+\lVert \HO(m)\rVert)\ptf\tag{$*$}\label{eq10.*}
\end{equation}
Comme $\tsQR\bigl(\HO(a)\bigr)=1$ on a
\[
\alpha(H_1)>0\qquad\forall\alpha\in\Delta_\Q^{\Rm}\ptf
\]
Comme $H_1$ appartient à $\ga_\Q^{\Rm}$ on a
\[
\alpha(H_1)=0
\qquad\forall\alpha\in\Delta_{\PO}^\Q
\]
et donc
\begin{equation}
\alpha(H_1)\ge0
\qquad\forall\alpha\in\Delta_{\PO}^{\Rm}\ptf\tag{$**$}\label{eq10.**}
\end{equation}
Comme $\theto$ préserve la chambre positive, les inégalités~\eqref{eq10.**} sont aussi vérifiées par $\theto(H_1)$.
Donc, pour $\alpha\in\Delta_{\PO}^\Q$ on a
\begin{equation}
\alpha\bigl(H_1-\theta(H_1)\bigr)=-\alpha\bigl(\theta(H_1)\bigr)=-\so\moins\alpha\bigl(\theto(H_1)\bigr)\le0
\label{eq10.1}
\end{equation}
pour notre choix de $\so$. On observe que
comme
\[
m_1\moins\,a\moins\,n_\Q\moins\mu\,\delta\eta\,\,n\,a\,m\in\Omega
\]
où $\Omega$ est un compact alors
\[
\HQ\bigl(m_1\moins\,a\moins\,n_\Q\moins\,\mu\,\theta(\eta\,n\,a\,m)\bigr)=\HQ\bigl(a\moins\theta(\eta\,n\,a\,m)\bigr)
\]
est borné. Il en résulte que la projection de
\begin{equation}
H_1-\theta(H_1)-\theta\bigl(\HO(m)\bigr)-\HO(w_sn')\label{eq10.2}
\end{equation}
sur $\ga_\Q^{\Rm}$ est bornée.
En combinant~\eqref{eq10.1} et~\eqref{eq10.2} on obtient que, compte tenu du lemme~\ref{nilneg}, on a
\begin{equation}
H_1-\theta(H_1)=X-Y\label{eq10.3}
\end{equation}
avec
\begin{equation}
\lVert X\rVert \le c_1( 1+\lVert \HO(m)\rVert)\Qquad
{et}\vpi(Y)\ge0\qquad\forall\vpi\in\hDelta_{\PO}^{\Rm}\ptf \tag{\ref{eq10.3}$'$}\label{eq10.3'}
\end{equation}
Par ailleurs
\[
H_1-\theta(H_1)=\bigl(H_1-\theto(H_1)\bigr)+(1-\so)\theto(H_1)
\]
On a vu que $H_1$ et $\theto(H_1)$ sont
dans l'adhérence de la chambre de Weyl positive; il résulte du lemme~\ref{ws} que
$(1-\so)\theto(H_1)$ est combinaison linéaire à coefficients
positifs ou nuls de coracines positives. Donc on a aussi
\begin{equation}
H_1-\theto(H_1)=X-Y_1\label{eq10.4}
\end{equation}
avec
\begin{equation}
\vpi(Y_1)\ge0\qquad\forall\vpi\in\hDelta_{\PO}^{\Rm}\ptf \tag{\ref{eq10.4}$'$}\label{eq10.4'}
\end{equation}
Si $X_0$ et $Y_0$ sont les projections de $X$ et $Y_1$ sur le sous-espace des $\theto$\hyph invariants
dans $\gao^{\Rm}$,
il résulte de~\eqref{eq10.4} que l'on a
\[
0=X_0-Y_0
\]
et donc, compte tenu de~\eqref{eq10.3'} on a
\begin{equation}
\lVert Y_0\rVert = \lVert X_0\rVert \ll ( 1+\lVert \HO(m)\rVert)\ptf\tag{\ref{eq10.3}$''$}\label{eq10.3''}
\end{equation}
Mais
\[
\lVert Y_0\rVert =\lVert Y_1\rVert \cos(Y_1,Y_0)
\]
et
le cosinus de l'angle entre $Y_1$ et $Y_0$ est minoré par une constante $c_1>0$:
\[
\cos(Y_1,Y_0)\ge c_1\ptf
\]
En effet, sinon on pourrait trouver $Z\in\gao^{\Rm}$ non nul
satisfaisant~\eqref{eq10.4'} et dont la projection sur les invariants serait nulle. Mais alors, si $\ell$
est l'ordre de $\theto$, on a
\[
Z+\sum_{r=1}^{r=\ell-1}\theta_0^r(Z)=0
\]
et donc $- Z$ satisfait aussi~\eqref{eq10.4'} ce qui impose $\vpi(Z)=0$ pour tout $\vpi\in\hDelta_{\PO}^{\Rm}$
et donc $Z=0$
ce qui est absurde. Donc
\[
\lVert Y_1\rVert \ll \lvert Y_0\rVert
\]
et l'inégalité~\eqref{eq10.3''} ci-dessus implique l'inégalité
\begin{equation}
\lVert Y_1\rVert \ll ( 1+\lVert \HO(m)\rVert)\ptf\label{eq10.5}
\end{equation}
Donc~\eqref{eq10.4} et~\eqref{eq10.5} montrent que
\[
\lVert H_1-\theto(H_1)\rVert \ll ( 1+\lVert \HO(m)\rVert)\ptf
\]
Comme $\theto$ est une isométrie
on a des inégalités analogues pour $\lVert \theta_0^rH_1-\theta_0^{r+1}H_1\rVert$
et donc aussi pour
$\lVert H_1-\theta_0^{r}H_1\rVert$ (avec une autre constante). On en déduit que
si $H_0$ est la moyenne sur la $\theto$\hyph orbite de $H_1$ on a encore une inégalité similaire:
\begin{equation}
\lVert H_1-H_0\rVert \ll ( 1+\lVert \HO(m)\rVert)\ptf\label{eq10.6}
\end{equation}
Mais
\[
H_0-\so H_0=(H_0-H_1)
+\bigl(H_1-\theta(H_1)\bigr)+\theta(H_1-H_0)
\]
et donc~\eqref{eq10.3} et~\eqref{eq10.6} impliquent que
\[
H_0-\so H_0=X_2-Y_2
\]
avec
\begin{equation}
\lVert X_2\rVert \ll ( 1+\lVert \HO(m)\rVert)\Qquad
{et}\vpi(Y_2)\ge0\quad\forall\vpi\in\hDelta_{\PO}^{\Rm}\ptf \tag{\ref{eq10.6}$'$}\label{eq10.6'}
\end{equation}
Comme
$\alpha(H_1)\ge0$ pour tout $\alpha\in\Delta_{\PO}^{\Rm}$,
il en est de même de $\alpha(H_0)$ et
la projection de $(H_0-\so H_0)$ sur $\gao^{\Rm}$
est une combinaison à coefficients positifs de racines positives;
on a donc pour tout $\vpi\in\hDelta_{\PO}^{\Rm}$
\[
\vpi(H_0-\so H_0)\ge 0
\]
et on en déduit que
\begin{equation}
\lVert H_0-\so H_0\rVert \ll ( 1+\lVert \HO(m)\rVert)\ptf\label{eq10.7}
\end{equation}
On peut écrire $H_1$ sous la forme
\[
H_1=\sum_{\alpha\in\Delta_{\PO}^{\Rm}-\Delta_{\PO}^\Q}\alpha(H_1)\vpi_\alpha^\vee\ptf
\]
Soit $\alpha\in\Delta_{\PO}^{\Rm}$; trois cas se présentent:
\par\noindent
Cas (i): L'orbite de $\alpha$ sous $\theto$ rencontre $\Delta_{\PO}^\Q$. Soit
$\alpha'\in\Delta_{\PO}^\Q$ un élément de cette orbite. On a alors $\alpha'(H_1)=0$ et donc
\[
\alpha(H_0)=\alpha'(H_0)=\alpha'(H_0-H_1)
\]
et on déduit de~\eqref{eq10.6} l'inégalité
\[
\lvert \alpha(H_0)\rvert \ll ( 1+\lVert \HO(m)\rVert)\ptf
\]
\par\noindent Cas (ii):
L'orbite de $\alpha$ sous $\theto$ ne rencontre pas $\Delta_{\PO}^\Q$
et pour au moins un $\alpha'$ dans l'orbite on a
\[
\vpi_{\alpha'}^\vee\ne\so\vpi_{\alpha'}^\vee\ptf
\]
On observe que, d'après le lemme~\ref{ws}, pour tout $X$ dans l'adhérence de la chambre positive, on a
\[
X-\so X=\sum_{\gamma\in\Delta_{\PO}^{\Rm}} c_\gamma(X,\so)\gamma^\vee
\]
avec
\[
c_\gamma(X,\so)\ge0\ptf
\]
et donc $c_\gamma(X,\so)=0$ pour tout $\gamma$ équivaut à $X=\so X$. Comme nous supposons
$\vpi_{\alpha'}\ne \so\vpi_{\alpha'}$, il existe
$\gamma$ tel que
\[
 c_\gamma(\vpi_{\alpha'},\so)= c_2>0
\]
et on a donc
\[
\vpi_\gamma(H_0-\so H_0)=\sum_{\beta} \beta(H_0)
c_\gamma(\vpi_\beta,\so)\ge c_2\,\alpha'(H_0)=c_2\,\alpha(H_0)
\]
et on a donc encore
\[
\lvert \alpha(H_0)\rvert \ll ( 1+\lVert \HO(m)\rVert)\ptf
\]
\par\noindent Cas (iii):
La dernière possibilité serait que l'orbite de $\alpha$ sous $\theto$ ne rencontre pas $\Delta_{\PO}^\Q$
et que pour tout $\alpha'$ dans l'orbite on ait
\[
\vpi_{\alpha'}^\vee=\,\so\vpi_{\alpha'}^\vee\ptf
\]
Mais dans ce cas l'ensemble des racines $\beta\in\Delta_{\PO}^{\Rm}$ qui sont orthogonales à tous
les $\vpi_{\alpha'}$ est l'ensemble $\Delta_{\PO}^P$
des racines simples du sous-groupe de Levi $\M$
d'un sous-groupe parabolique $P$ qui est
$\theto$\hyph stable, qui est tel que
\[
\Q\subset P\varsubsetneq\Rm\;,
\]
et dont le groupe de Weyl $\weyl^\M$ contient ${\so}$.
On aurait donc
\[
\delta=w_{\so}\delto\eta\in\tP(F)\varsubsetneq\tRm(F)\ptf
\]
Ceci est impossible puisque par hypothèse $\delta$ appartient
à $\tG(\Q,\R)$.
\par\noindent
Il résulte de cette discussion que
\[
\lvert \alpha(H_0)\rvert \ll ( 1+\lVert \HO(m)\rVert)
\]
pour tout $\alpha$ et donc
\begin{equation}
\lVert H_0\rVert \ll ( 1+\lVert \HO(m)\rVert)\ptf\label{eq10.8}
\end{equation}
Compte tenu de~\eqref{eq10.6} et~\eqref{eq10.8} on obtient, comme espéré, l'inégalité~\eqref{eq10.*}:
\[
\lVert H_1\rVert \ll ( 1+\lVert \HO(m)\rVert)\ptf\qedhere
\]
\end{proof}

Rappelons que l'on a posé
\[
\YQ=\gA_\G\Q(F)\bs\Gadef\;,
\]
et introduisons la sous-groupe parabolique
\[
\Q_\delta=\Q\cap\delta\moins\Q\delta\ptf
\]
Avec ces notations on a la proposition suivante:

\begin{proposition}\label{specconva}
Pour tout couple de sous-groupes paraboliques standard $\Q\subset\R$
et tout $\delta\in\tG(\Q,\R)$ l'intégrale
\[
 \int_{\mathbf{Y}_{\Q_\delta}} \tsQR (\xT)
\bigl\lvert\tronc_1^{T,\Q } K_{\Q,\delta}(x,x)\bigr\rvert \dd x
\]
est convergente si $\dPO(T)\ge c$ où $c$ est une constante ne dépendant que du support de~$\fff$
\end{proposition}

\begin{proof}
Considérons $x=nmak$ avec $k\in\K$,
$n\in\Omega\subset\N_\Q(\adef)$ où $\Omega$ est un compact,
$m\in\goth{S}_\Q\subset\M_\Q(\adef)$ où
$\goth{S}_\Q$ est un ensemble de Siegel pour $\M_\Q$
(en particulier $\HQ(m)=0$) et $a\in\gA_\Q^\G$.
On doit estimer
\begin{multline*}\
\smash[b]{\int_\Omega\int_{\gA_\Q^\G}\int_{\goth{S}_\Q}}
\ee^{-2\demisomQ\left(\HO(a)\right)}
\tsQR (\aT)\\
{}\times\sum_{\xi\in\Q_\delta(F)\bs\Q(F)}
\bigl\lvert\tronc_1^{T,\Q }K_{\Q,\delta}(mak,\xi nmak)\bigr\rvert\dd n\dd a\dd m
\end{multline*}
\cad
\[
\int_\Omega\int_{\gA_\Q^\G}\int_{\goth{S}_\Q}
\tsQR (\aT) \sum_{\xi\in\Q_\delta(F)\bs\Q(F)}
\bigl\lvert\tronc_1^{T,\Q }K_{\Q,\delta}(mk,\xi a^{1-\delta} n'mk)\bigr\rvert\dd n\dd a\dd m 
\]
avec $n'=a\moins n a$.
D'après le lemme~\ref{premaj}, l'opérateur de troncature
fournit un noyau
\[
(m_1,m_2)\mapsto
\tronc_1^{T,\Q } K_{\Q,\delta}(m_1k,\xi a^{1-\delta}n'm_2k)
\]
à décroissance rapide en $m_1$ et à croissance lente en $m_2$
sur le domaine de Siegel $\goth{S}_\Q$
de $\M_\Q$ et, par restriction à la diagonale, on obtient
une fonction $\phi$ à décroissance rapide en $m=m_1=m_2$ sur $\goth{S}_\Q$
en ce sens que, pour tout $N$,
elle est majorée par
\[
c_N(\phi)\,\lvert m\rvert^{-N}
\]
ce qui, au vu des lemmes~\ref{xifini} et~\ref{maja}, permet de compenser la croissance éventuelle
due à la somme sur $\xi$, et de contrôler l'intégrale sur $a$.
\end{proof}

Notons $\tweyl(\Q,\R)$ un ensemble de représentants de $\tG(\Q,\R)$ modulo
$\Q(F)$ à droite et à gauche:
\[
\tweyl(\Q,\R)\simeq\Q(F)\bs \tG(\Q,\R)/\Q(F)\ptf
\]
C'est un ensemble fini. On a
\[
 K_{\Q,\R}(x,y)=\sum _{\delta\in \tweyl(\Q,\R)}\,\sum_{\xi\in\Q_\delta(F)\bs\Q(F)}
K_{\Q,\delta}(x,\xi\,y)\ptf
\]

\begin{proposition}\label{specconvb} 
Si $\dPO(T)\ge c$\textup,
où $c$ est la constante de la proposition~\ref{specconva}, on a
\[
\int_{\XG}\lvert k_{\mathrm{spec}}^\T(x) \rvert \dd x < \infty\ptf
\]
\end{proposition}

\begin{proof}
On observe que
\[
\int_{\XG}\lvert k_{\mathrm{spec}}^\T(x)\rvert \dd x \le\sum_{\Q \subset \R } \tve(\Q,\R)
\int_{\YQ}\tsQR (\xT)
\lvert\tronc_1^{T,\Q } K_{\Q,\R}(x,x)\rvert \dd x
\]
et que
\[
 K_{\Q,\R}(x,x)=\sum _{\delta\in \tweyl(\Q,\R)}\,\sum_{\xi\in\Q_\delta(F)\bs\Q(F)}
K_{\Q,\delta}(x,\xix) \ptf
\]
L'assertion est alors une conséquence immédiate du lemme~\ref{prenul} et
de la proposition~\ref{specconva}.
\end{proof}

On aurait pu déduire cette proposition la conjonction de
l'identité fondamentale~\ref{fonfon} et du théorème~\ref{geoconv}.
Mais la preuve donnée ci-dessus va pouvoir se raffiner
pour établir la convergence du développement spectral grossier,
\cad du développement suivant les données
cuspidales.

\section{Annulations supplémentaires}\label{annulsup}

On considère comme ci-dessus $\Q\subset\R$ et
on suppose qu'il existe un sous-ensemble parabolique $\tP$
avec
\[
\Q\subset P\subset\R\ptf
\]
On
note $\tQp$ le plus petit (resp.~$\tRm$ le plus grand)
sous-ensemble parabolique avec cette propriété, \cad
que, dans les notations du lemme~\ref{qplus}, on a
\[
\Q\subset\Qp\subset P\subset\Rm\subset\R\ptf
\]
On choisira des représentants de l'ensemble fini de doubles classes
\[
\tweyl(\Q,\R)\simeq\Q(F)\bs \tG(\Q,\R)/\Q(F)
\]
de la forme $\delta=w_s$ où $w_s$ représente un élément
\[
s=\so\rtimes\theto
\]
appartenant à l'ensemble de Weyl $\weyl^{\tM_{\Rm}}$ de $\tM_{\Rm}$. En choisissant $\so$ de longueur
minimale dans sa double classe il résulte des lemmes~\ref{echangebis} et~\ref{weylg} que
\[
s\alpha>0\Quad{et} s\moins\alpha>0\qquad\forall\alpha\in\Delta_{\PO}^\Q\ptf
\]
Donc, plus généralement, $s\alpha>0$ pour toute racine positive pour $\M_\Q$.
On en déduit que,
\[
\Ms=Q_\delta\cap\M_\Q=\M_\Q\cap w_s\moins\Q\, w_s
\]
est un sous-groupe
parabolique standard de $\M_\Q$. En effet, les racines dans le radical unipotent
de $\Ms$ sont les restrictions à un sous-espace de $\gao$ de racines $\alpha$ pour $\M_\Q$
telles que la restriction de
$s\alpha$ à $\ga_\Q^\G$
soit une racine pour le radical unipotent de $\Q$.
Comme $\Q$ est standard de telles racines $s\alpha$ sont positives
et donc nécessairement $\alpha$ est positive.
On note $S$ le sous-groupe parabolique
standard de $\G$ tel que
\[
\Ms=S\cap\M_\Q\ptf
\]
On note $\N_S$ son radical unipotent. On observera que
\[
\ga_S\supset\ga_\Q\Quad{et}s(\ga_S)\supset\ga_\Q
\]
et que donc
\[
s(\gao^S)\perp\ga_\Q\ptf\qedhere
\]

\begin{lemme}\label{binul}
Supposons $\delta\in \tweyl(\Q,\R)$ et considérons l'expression
\[
\tsQR(\xT)\int_{\N_S^\Q(F)\bs\N_S^\Q(\adef)}
\tronc_1^{T,\Q } K_{\Q,\delta}(n_S\,x,n_\Q\,x)\dd n_S
\]
où $N_S^\Q=\N_S\cap\M_\Q$.
Alors\textup, si $\T$ avec $\dPO(T)\ge c\bigl(1+N(\fff)\bigr)$ où $c$ est une constante ne dépendant que de $\G$
et $N(\fff)$ dépend
du support de $\fff$\textup, l'intégrale est nulle pour tout $x\in\Gadef$
et tout $n_\Q\in\N_\Q(\adef)$
sauf peut-être si $\delta\equiv\delto$ comme double classe modulo $\Q(F)$.
\end{lemme}

\begin{proof}
Si l'intégrale double est non nulle alors il résulte du lemme~\ref{troncnul} que
\begin{equation}
\vpi(\xT)\le0\qquad\forall\vpi\in\hDelta_S^\Q\ptf\label{eq10.1a}
\end{equation}
De plus, le lemme~\ref{hoyx} montre que la projection de
\[
\tHO(w_s\,n_\Q\,x)-\HO(x)
\]
sur $\ga_\Q^\G$ reste bornée; plus précisément, reste dans une boule
dont le rayon dépend du support de $\fff$.
On a, pour un certain $n\in\NO(\adef)$,
\[
\tHO(w_s\,n_\Q\,x)-\HO(x)=\HO(w_{\so}\,n)+s\HO(x)-\HO(x)
\]
et donc, pour toute forme linéaire $\lambda$ sur $\ga_\Q^\G$, l'expression
\[
\lambda\bigl(\HO(w_{\so}\,n)+s\HO(x)-\HO(x)\bigr)
\]
reste bornée.
Maintenant, d'après le lemme~\ref{nilneg}, il existe une constante $c\ge0$ telle que pour tout $n\in\NO(\adef)$
\[
\vpi\bigl(\HO(w_{\so}\,n)\bigr)\le c
\]
pour tout $\vpi\in\hDelta_{\PO}^\G$.
Choisissons pour $\lambda$ la somme des $\vpi$ dans $\hDelta_{\Qp}^{\Rm}$.
C'est une forme linéaire $\theto$\hyph invariante et positive sur la chambre obtuse dans $\ga_{\Qp}^{\Rm}$. On a
\begin{equation}
\lambda\bigl(\HO(x)-s\HO(x)\bigr)\le C(\fff)\label{eq10.2a}
\end{equation}
pour une autre constante $C(\fff)$.
Posons
\[
X=\xT
\]
et decomposons $X$ sous la forme
\[
X=X_0^S+X_S^\Q+X_\Q^R+X_R
\]
où les $X_A^B$ sont les projections de $X$ sur les sous-espaces $\ga_A^B$.
En particulier
\begin{equation}
X_S^\Q=\sum_{\alpha\in\Delta_S^\Q} c_\alpha\alpha^\vee\label{eq10.3a}
\end{equation}
avec les $c_\alpha\le 0$ d'après~\eqref{eq10.1a} et
\begin{equation}
X_Q^\R=\sum_{\vpi\in\hDelta_\Q^R} c_\vpi\vpi^\vee\label{eq10.4a}
\end{equation}
avec les $c_\vpi>0$ si l'on suppose $\tsQR(X)\ne0$.
Maintenant~\eqref{eq10.2a} se récrit
\[
\lambda\bigl((\xT)-s(\xT)+(T-s T)\bigr)\le C
\]
soit encore
\begin{equation}
\lambda\bigl((X-s X)+(T-s T)\bigr)\le C\label{eq10.5a}
\end{equation}
On observe que, pour notre choix de $\lambda$,
\[
\lambda(X)=\lambda(X_\Q^R)
\]
mais il résulte du lemme~\ref{wT} que, en posant $s'_0=\thetomoins(\so)$
\[
X_\Q^R-s'_0 X_\Q^R
\]
est une combinaison linéaire de racines positives avec pour coefficients
des $\beta(X_\Q^R)$, où $\beta$ est une racine simple,
qui sont des réels positifs d'après~\eqref{eq10.4a}
et donc
\[
\lambda(X_\Q^R-s X_\Q^R)\ge0
\]
et~\eqref{eq10.5a} implique
\[
\lambda(T-s T)\le C+\lambda\bigl(s (X_0^S+X_S^\Q+X_R)\bigr)\ptf
\]
Supposons l'intégrale de l'énoncé non nulle pour $\delta\in \tweyl(\Q,\R)$ avec $\delta\nequiv\delto$ comme double classe modulo $\Q(F)$.
En particulier $\delta=w_s$ avec $s=\so\rtimes\theto$ et $\so\ne1$; donc
il existe une racine $\alpha>0$ dans le système de racines de $\Rm$
avec $s'_0\alpha<0$ et il résulte du lemme~\ref{wT} que
\[
\lambda(T-s T)=\lambda(T-s'_0 T)
\]
est arbitrairement grand
pour $T$ assez régulier. Pour montrer que ceci est impossible
il suffit de montrer que
\[
\lambda\bigl(s (X_0^S+X_S^\Q+X_\R)\bigr)\le0\ptf
\]
Comme $s'_0 X_R\in\ga_R$ on a $\lambda(s X_\R)=0$.
On rappelle que
\[
\ga_S\supset\ga_\Q\Qquad{et}s(\ga_S)\supset\ga_\Q
\]
et donc
\[
s(\gao^S)\perp\ga_\Q
\]
d'où on déduit que
$\lambda(s X_0^S)=0$.
On a
\[
\lambda(s X_S^\Q)=
\sum_{\alpha\in\Delta_S^\Q} c_\alpha\lambda(s \alpha^\vee)
\]
avec $c_\alpha\le0$ d'après~\eqref{eq10.3a}.
Mais, pour $\alpha\in\Delta_S^\Q$ on a
$\alpha^\vee=\beta^\vee+\gamma^\vee$ où
$\beta$ est la racine de $\Delta_{\PO}^\Q-\Delta_{\PO}^S$ qui se projette sur $\alpha\in\Delta_S^\Q$
et $\gamma^\vee\in\gao^S$.
Mais
\[
\Delta_{\PO}^\Q\subset\Delta_{\PO}^{\Qp}
\]
et donc $s\beta$ est une racine positive par choix des représentants dans $\tweyl(\Q,\R)$.
Il reste à observer que
\[
\lambda(s \alpha^\vee)=\lambda(s \beta^\vee)\ge0
\]
puisque $s\gamma^\vee$ est orthogonal à $\ga_\Q$.
\end{proof}

\begin{proposition}\label{ranul}
Supposons $\delta\in \tweyl(\Q,\R)$ et $\delta\nequiv\delto$ comme double classe modulo $\Q(F)$.
Alors\textup, pour $T$ assez régulier \textup(comme au lemme~\ref{binul}\textup)\textup, l'intégrale
\[
\int_{\gA_\G\Q_\delta(F) \bs \Gadef }\tsQR (\xT) \tronc_1^{T,\Q } K_{\Q,\delta}(x,x)\dd x 
\]
est nulle.
\end{proposition}

\begin{proof}
L'expression
\[
\int_{\gA_\G\Q_\delta(F) \bs \Gadef }\tsQR (\xT)
\tronc_1^{T,\Q } K_{\Q,\delta}(x,x)\dd x 
\]
s'écrit encore
\[
\int_{\gA_\G\Q_s(F)\N_s(\adef) \bs \Gadef } k_{\mathrm{spec}}^\T(\Q,\delta,x)\dd x
\]
avec
\[
k_{\mathrm{spec}}^\T(\Q,\delta,x)= \int_{\Q_\delta(F)\bs\Q_s(F)\N_s(\adef)}\tsQR (\xT)
\tronc_1^{T,\Q } K_{\Q,\delta}(mx,mx)\dd m \ptf
\]
Mais $k_{\mathrm{spec}}^\T(\Q,\delta,x)$
est égal au produit de $\tsQR (\xT)$ et de
l'intégrale double
\[
\int_{\N_s^\Q(F)\bs\N_s^\Q(\adef)}\int_{\N_\Q(\adef)\cap\Q_\delta(\adef)\bs\N_\Q(\adef)}
\tronc_1^{T,\Q } K_{\Q,\delta}(n_\Q\,n_s\,x,n_\Q\,n_s\,x)\dd n_\Q\dd n_s 
\]
où $N_s^\Q=\N_s\cap\M_\Q$. Compte tenu des invariances de $K_{\Q,\delta}$ observées
dans le lemme~\ref{invar},
l'intégrale double s'écrit encore
\[
\int_{\N_s^\Q(F)\bs\N_s^\Q(\adef)}\int_{\N_\Q(\adef)\cap\Q_\delta(\adef)\bs\N_\Q(\adef)}
\tronc_1^{T,\Q } K_{\Q,\delta}(n_s\,x,n_\Q\,x)\dd n_\Q\dd n_s \ptf
\]
On invoque alors le lemme~\ref{binul}.
\end{proof}

Par abus de notation
nous écrirons
\[
\delto\in\tweyl(\Q,\R)
\]
pour exprimer que
la double classe modulo $\Q(F)$ définie par $\delto$ appartient
à l'ensemble $\tweyl(\Q,\R)$.
On rappelle que, par définition de $ \tweyl(\Q,\R)$, on ne peut avoir
\[
\delto\in\tweyl(\Q,\R)\Qquad{et} \delto\in\tP(F)
\]
que pour un seul sous-ensemble parabolique $\tP$ avec $\Q\subset P\subset\R$.
Comme
\[
\delto\in\tQp(F)\subset\tRm(F)
\]
on voit que
\[
\Qp=\Rm \Qquad{équivaut à}\delto\in\tweyl(\Q,\R)\ptf
\]
Dans ce cas, si $\tP$ est le seul
sous-ensemble parabolique qui vérifie $\Q\subset P\subset\R$ on a
\[
\tve(\Q,\R)\coloneqq (-1)^{a_{\tRm}-a_{\tG}}=(-1)^{a_{\tP}-a_{\tG}}\ptf
\]
Nous poserons
\[
 \tvedQR=\sum_{\{\tP\mid \Q\subset P\subset\R\}}(-1)^{a_{\tP}-a_{\tG}}\ptf
\]
Ce nombre est nul sauf si un seul sous-ensemble parabolique $\tP$ vérifie
$\pQ\subset P\subset\R$ auquel cas
\[
 \tvedQR=(-1)^{a_{\tP}-a_{\tG}}\ptf
\]
On a donc
\[
\tvedQR=
\newnot{eta(Q,R)@$\tvedQR$}{tvedqr}%
\begin{cases}
\tve(\Q,\R)& \text{si $\Qp=\Rm$}\\
0& \text{sinon.}
\end{cases}
\]
Considérons maintenant
\[
\YQdo=\gA_\G\Qdo(F)\bs\Gadef
\qquad\text{avec $\Qdo=\Q_{\delto}=\Q\cap\delta_0\moins\Q\delto$}\ptf
\]
On a observé que $\Qdo$ est un sous-groupe parabolique standard.

\begin{corollaire}\label{explicite}
Si $T$ est assez régulier \textup(comme au lemme~\ref{binul}\textup)\textup, l'intégrale
\[
\int_{\XG}k_{\mathrm{spec}}^\T(x)\dd x 
\]
est égale à la somme
\[
\sum_{\{Q,R\mid  \PO\subset Q\subset R\}} \tvedQR
\int_{\YQdo} \tsQR (\xT)
\tronc_1^{T,\Q } K_{\Q,\delto}(x,x)\dd x
\ptf
\]
Dans le cas non tordu on a simplement
\[
\int_{\XG}k_{\mathrm{spec}}^\T(x)\dd x =\int_{\XG}\tronc^\T_1K_\G(x,x)\dd x
\]
si $T$ est assez régulier.
\end{corollaire}

\begin{proof} On rappelle que
\[
\int_{\XG}k_{\mathrm{spec}}^\T(x)\dd x = \sum_{\Q \subset \R } \tve(\Q,\R)
\int_{\YQ}\tsQR (\xT)
\tronc_1^{T,\Q } K_{\Q,\R}(x,x)\dd x
\]
et que
\[
 K_{\Q,\R}(x,x)=\sum _{\delta\in \tweyl(\Q,\R)}\,\sum_{\xi\in\Q_\delta(F)\bs\Q(F)}
K_{\Q,\delta}(x,\xix) \ptf
\]
La première assertion résulte alors immédiatement de la proposition~\ref{ranul}.
La seconde assertion résulte de ce que, dans le cas non tordu, la condition
$\delto\in\tweyl(\Q,\R)$ implique $\Q=\R$ mais le lemme~\ref{nullos} montre que
$\sigma_\Q^\Q=0$ sauf si $\Q=\G$.
\end{proof}

\section[Contrôle du développement en $\chi$]{\mathversion{bold}Contrôle du développement en $\chi$}

\begin{proposition}[\cite{ATFII}*{Lemma 2.3}]\label{arthurnul}
Soit $\Q=\N_\Q\M_\Q$ un sous-groupe parabolique
de $\G$. Soit $\tP_i$ des sous-ensembles paraboliques de $\tG$ avec $\Q\subset P_i$.
Soit $\chi$ une donnée cuspidale.
Supposons données une famille finie d'éléments $x_i$ et $y_i$ dans $\Gadef$ et des constantes $c_i$
telles que
\[
\sum c_i \int_{\N_\Q(F)\bs\N_\Q(\adef)} K_{\tP_i}(nmx_i,y_i)\dd n=0
\]
pour tout $m\in \M_\Q(\adef)$ tel que $\HQ(m)=0$\textup, alors
\[
\sum c_i \int_{\N_\Q(F)\bs\N_\Q(\adef)} K_{\tP_i,\chi}(nx_i,y_i)\dd n=0\ptf
\]
\end{proposition}

\begin{proof}
Posons
\[
\phi(m)=\int_{\N_\Q(F)\bs\N_\Q(\adef)}\sum c_i K_{\tP_i}(nmx_i,y_i)\dd n
\]
et
\begin{equation}
A(\psi)=\int_{\M(F)\bs\M(\adef)^1}\phi_\chi(m)\psi(m)\dd m\label{eq10.1b}
\end{equation}
pour $\psi\in L^2(\XM)$.
La fonction $\phi$ est le terme constant suivant $\N_\Q$ de
\[
\sum c_i K_{\tP_i}(mx_i,y_i)
\]
ce qui
annule les éventuelles contributions des données cuspidales attachées à un sous-groupe
parabolique contenant $\Q$ strictement. On invoque la décomposition spectrale de $\phi$.
L'orthogonalité des contributions relatives à des données cuspidales
inéquivalentes montre que
$A(\psi)=0$
si $\psi$ est de type $\chi'\ne\chi$. Mais par ailleurs, si $\psi$ est de type $\chi$
on a
\[
A(\psi)=\int_{\goth{S}}\phi(m)\psi(m)\dd m=0
\]
par hypothèse.
Il en résulte que l'intégrale \eqref{eq10.1b} est nulle pour toute $\psi$ ce qui implique la nullité
de la fonction continue
\[
m\mapsto \phi_\chi(m)\ptf\qedhere
\]
\end{proof}

\begin{lemme}\label{hoyxchi}
Soit $\Q=\N_\Q\M_\Q$ un sous-groupe parabolique
de $\G$.
Soit $\chi$ une donnée cuspidale pour $\Q$.
Supposons que
\[
K_{\Q,\delta,\chi}(x,y)\ne0\ptf
\]
Alors il existe un compact $C\subset\gao$,
$m\in\M_\Q(\adef)$ avec $\HQ(m)=0$
et $\eta\in\Q(F)$ tels que
\[
\HO\bigl(\eta\,\theta(y)\bigr)-\HO(mx)\in C\ptf
\]
\end{lemme}
\begin{proof} En reprenant la preuve de la proposition~\ref{arthurnul} avec
$$\phi(m)=K_{\Q,\delta}(m\,x,y)$$ on voit que
$K_{\Q,\delta,\chi}(x,y)\ne0$
implique l'existence d'un $m\in\M_\Q(\adef)$ avec $\HQ(m)=0$
tel que $K_{\Q,\delta}(m\,x,y)\ne0$. On conclut en invoquant le lemme~\ref{hoyx}.
\end{proof}

\begin{corollaire}\label{binulb}
Pour tout $\chi$
\[
\int_{\N_s^\Q(F)\bs\N_s^\Q(\adef)}
\tronc_1^{T,\Q } K_{\Q,\delta,\chi}(n_s\,x,n_\Q\,x)\dd n_s =0
\]
sauf peut-être si $\delta\equiv\delto$ comme double classe modulo $\Q(F)$.
\end{corollaire}

\begin{proof} 
On reprend la preuve du lemme~\ref{binul} en invoquant
Le lemme~\ref{hoyxchi} au lieu du lemme~\ref{hoyx}.
\end{proof}

Nous allons maintenant énoncer un raffinement de la proposition~\ref{specconvb}
et du corollaire~\ref{explicite}. On pose
\[
k^\T_\chi(x)=\sum_{\tP\supset\tPO}(-1)^{a_{\tP}-a_{\tG}}\sum_{\Q\subset P\subset\R}\sum_{\xi\in\Q(F)\bs\G(F)}
\tsQR(\xixT )
\,\tronc_1^{T,\Q}K_{\tP,\chi}(\xix,\xix)\ptf
\]

\begin{theoreme}\label{chiconv} Si $T$ est assez régulier \textup(comme au lemme~\ref{binul}\textup)\textup, on a
\[
\int_{\XG}k_{\chi}^\T(x)\dd x
=\sum_{\{Q,R\mid  \PO\subset Q\subset R\}} \tvedQR\,A_\Q^\R(\T,\fff,\omega,\chi)
\]
avec
\[
A_\Q^\R(\T,\fff,\omega,\chi)=
\int_{\YQdo} \tsQR (\xT)
\tronc_1^{T,\Q } K_{\Q,\delto,\chi}(x,x)\dd x
\ptf
\]
Dans le cas non tordu on a simplement
\[
\int_{\XG}k_{\chi}^\T(x)\dd x =\int_{\XG}\tronc^\T_1K_\chi(x,x)\dd x
\]
si $T$ est assez régulier.
La fonction $\fff$ \'etant suppos\'ee $\K$\hyph finie (\`a droite et \`a gauche)
on a:
\[
\sum_{\chi\in\XX}\int_{\XG}\bigl\lvert k^\T_\chi(x) \bigr\rvert \dd x < \infty\ptf
\]
\end{theoreme}

\begin{proof}
On a
\[
k^\T_{\chi}(x) = \sum_{\Q \subset \R }\,\sum_{\xi\in\Q(F) \bs\G(F)}k_{\chi}^\T(\Q,\R,\xix)
\]
avec
\[
k_{\chi}^\T(\Q,\R,x)=\tsQR(\xT)
\tronc_1^{T,\Q}\sum_{\{\tP\mid \Q \subset P\subset \R\}}
(-1)^{a_{\tP}-a_{\tG}}\,K^\Q_{\tP,\chi} (x,x)\ptf
\]
mais compte tenu de la proposition~\ref{arthurnul} on voit en reprenant
les arguments du lemme~\ref{prenul} que
\[
k_{\chi}^\T(\Q,\R,x)= \tve(\Q,\R)\,\,\tsQR(\xT)
\tronc_1^{T,\Q}K_{\Q,\R,\chi}(x,x)\dd x\ptf
\]
L'analogue du résultat d'annulation~\ref{binul}, mais pour $\chi$ fixé,
s'obtient lui aussi grâce à la proposition~\ref{arthurnul}.
Reprenons alors la démonstration de la proposition~\ref{specconvb}. Il suffit de
démontrer l'analogue de la proposition~\ref{specconva}, à savoir la finitude de
\[
\sum_\chi \int_{\mathbf{Y}_{\Q_\delta}} \tsQR (\xT)
\bigl\lvert\tronc_1^{T,\Q } K_{\Q,\delta,\chi}(x,x)\bigr\rvert \dd x\ptf
\]
Il convient d'abord d'avoir pour
\[
\sum_\chi \bigl\lvert\tronc_1^{T,\Q }K_{\Q,\delta,\chi}\bigr\rvert
\]
des estimées similaires à celles obtenues pour $\lvert \tronc_1^{T,\Q }K_{\Q,\delta}\rvert$;
en particulier, l'analogue du lemme~\ref{premaj}
résulte du corollaire~\ref{majod}.
Grâce à la proposition~\ref{arthurnul} on voit que
les propriétés~\ref{xifini} et~\ref{maja} restent vraies pour $K_{\Q,\delta,\chi}$
puisque son support est contrôlé par celui de $K_{\Q,\delta}$.
La convergence de
\[
\sum_\chi \int_{\XG}\lvert k^\T_\chi(x) \rvert \dd x < \infty
\]
en résulte. 
\end{proof}

\chapter{Formule des traces: propriétés formelles}\label{ch11}

On a
\[
k^\T_{\mathrm{geom}}(x)=\sum_{\go\in\gO} k^\T_\go(x)\Qquad{et}
k_{\mathrm{spec}}^\T(x)=\sum_{\chi\in\XX}k^\T_\chi(x)\ptf
\]
On sait que
\[
k^\T_{\mathrm{geom}}=k_{\mathrm{spec}}^\T
\]
et on note $k^\T$ la valeur commune de ces deux fonctions.
Dans ce qui suit l'indice $\bullet$ peut représenter une classe de conjugaison quasi semi-simple
$\go$ ou encore une donnée cuspidale $\chi$ ou enfin être vide.
Nous allons rencontrer le noyau de la formule des traces pour divers espaces tordus
et diverses fonctions.
Pour tenir compte de cette dépendance nous écrirons
\[
k^{T,\tG}_\bullet(\fff,\omega;x)
\]
au lieu de $k^{T}_\bullet(x)$.
Rappelons enfin que la convergence des intégrales
\[
\int_{\XG} k^{T,\tG}_\bullet(\fff,\omega;x)\dd x
\]
a été l'objet des théorèmes~\ref{geoconv} et~\ref{chiconv}.

\section{Le polynôme asymptotique}

On a introduit et calculé au lemme~\ref{gpq} une fonction $\gamma_\Q(X)$; nous utiliserons ici son avatar tordu:
soit $\tQ$ un sous-ensemble parabolique, on pose
\[
\gamma_{\tQ}(X)=
\int_{\gA_{\tG}\bs\gA_{\tQ}}\Gamma_{\tQ}(\HO(a),X)\dd a=
\int_{\ga_{\tG}\bs\ga_{\tQ}}\Gamma_{\tQ}(H,X)\dd H\ptf
\]
C'est un polynôme en $X$ homogène de degré $a_{\tQ}-a_{\tG}$.

On notera $\fff_{\tQ}$ une fonction dans $\ctyc\bigl(\tM_\Q(\adef)\bigr)$ telle que
pour tout $m\in\tM_\Q(\adef)$
de la forme $m=m_0\delto $ avec $m_0\in\M_\Q(\adef)$ et
$\HQ(m_0)\in\ga_\G$
on ait
\[
\int_{\gA_\Q}\fff_{\tQ}( z m)\dd z=
\int_\K\int_{\N_\Q(\adef)}\int_{\gA_\G\gA_{\tQ}\bs\gA_\Q}\fff^1(k\moins\,
a\moins\,m\,a\,n\,k)\dd a
\dd n\dd k \ptf
\]
Il est facile de voir que de telles fonctions $\fff_{\tQ}$ existent.

\begin{theoreme}\label{poly}
Il existe une fonction polynôme
\[
T \mapsto J^\T_\bullet(\fff,\omega)
\]
dont le degré est inférieur
ou égal à
\[
a_{\tPO}-a_{\tG}=\dim \gatPO^\G
\]
telle que
si $\dPO(T)\ge c(\fff)$ on ait
\[
J^{T,\tG}_\bullet(\fff,\omega)=\int_{\XG} k^{T,\tG}_\bullet(\fff,\omega;x)\dd x
\]
et
\[
J^{T+X,\tG}_\bullet(\fff,\omega)=\sum_{\tQ} \gamma_{\tQ}(X) J^{T,\tQ}_\bullet(\fff_{\tQ},\omega)
\ptf
\]
La constante $c(\fff)$ ne dépend que
du support de $\fff$.
\end{theoreme}

\begin{proof} D'après les
théorèmes~\ref{geoconv} et~\ref{chiconv} les intégrales sont convergentes.
Rappelons que
\[
k^{T,\tG}_\bullet(\fff,\omega;x)=\sum_{\tP\supset\tPO}(-1)^{a_{\tP}-a_{\tG}}\,
\sum_{\xi\in P(F)\bs\G(F)}\htautP(\xixT)\,K_{\tP,\bullet}(\xix,\xix)\ptf
\]
C'est la définition lorsque $\bullet$ est soit vide soit $\bullet=\go$ une classe de conjugaison quasi semi-simple.
Dans le cas $\bullet=\chi$ il convient d'observer, en utilisant les propositions~\ref{fonfon} et~\ref{arthurnul}, que l'identité fondamentale~\ref{fonfon} est encore valable pour $k^\T_\chi$.
Pour alléger la notation on omettra dans le reste de la preuve l'indice $\bullet$.
Posons pour tout ensemble parabolique standard $\tQ$
\[
 k^{T,\tG}_{\tQ}(\fff,\omega;x)= \sum_{\{\tP\mid \tPO\subset\tP\subset\tQ\}}
\,\sum_{\xi\in P(F)\bs\Q(F)}(-1)^{a_{\tP}-a_{\tQ}}\,\htautPQ(\xixT )\,K_{\tP}(\xix,\xix)\ptf
\]
Compte tenu de la proposition~\ref{GammaHXt}, on a
\[
\htautP(H-X)=\sum_{\tP\subset\tQ\subset\tR}(-1)^{a_{\tQ}-a_{\tR}}
\Gamma_{\tQ}(H,X)\, \htautPQ(H)\ptf
\]
On observera que la fonction $H\mapsto\Gamma_{\tQ}(H,X)$ ne dépend que de la projection
de $H$ sur $\ga_{\tQ}$.
On a alors
\[
k^{T+X,\tG}(\fff,\omega;x)=\sum_{\tQ}\sum_{\eta\in\Q(F)\bs\G(F)}
\Gamma_{\tQ}(\xT,X) \, k^{T,\tG}_{\tQ}(\fff,\omega;x)
\]
et donc, pour $T$ et $X$ assez réguliers, on a
\[
\int_{\XG} k^{T+X,\tG}(\fff,\omega;x)\dd x
=\sum_{\tQ}\int_{\YQ}\Gamma_{\tQ}(\xT,X) \, k^{T,\tG}_{\tQ}(\fff,\omega;x)\dd x
\]
avec
\[
\YQ=\gA_\G\Q(F)\bs\Gadef
\]
les intégrales étant absolument convergentes. On obtient
\[
\int_{\XG} k^{T+X,\tG}(\fff,\omega;x)\dd x= \sum_{\tQ} \gamma_{\tQ}(X)\int_{\gA_{\tQ}\gA_\G\Q(F)\bs\Gadef}
k^{T,\tG}_{\tQ}(\fff,\omega;x)\dd x\ptf
\]
Compte tenu de la définition de $\fff_{\tQ}$ et de la décomposition d'Iwasawa on voit que
\[
\int_{\gA_{\tQ}\gA_\G \Q(F)\bs\G(\adef)} k^{T,\tG}_{\tQ}(\fff,\omega;x)\dd x= \int_{\XMQ} k^{T,\tQ}(m,\fff_{\tQ})\dd m
\]
et donc en posant, pour $T$ et $X$ assez réguliers,
\[
J^{T,\tG}(\fff,\omega)=\int_{\XG} k^{T,\tG}(\fff,\omega;x)\dd x
\]
on a
\[
J^{T+X,\tG}(\fff,\omega)=\sum_{\tQ} \gamma_{\tQ}(X) J^{T,\tQ}(\fff_{\tQ},\omega)\ptf
\]
Il reste à observer que
\[
X\mapsto \gamma_{\tQ}(X)
\]
est une fonction polynôme sur $\gao$ de degré $a_{\tQ}-a_{\tG}$.
\end{proof}

\section{Action de la conjugaison}

Soit $y\in\Gadef$ et posons
\[
\fff^y(x)=\fff(y\,x\,y\moins)
\]
et soit $\fff_{\tQ,y}\in\ctyc(\tM_\Q(\adef))$ telle que,
pour tout $m\in\tM_\Q(\adef)$, l'intégrale
\[
\int_{\gA_Q}\fff_{\tQ,y}(zm)dz
\]
soit égale à
\[
\frac{\jtG}{\jtQ}\boldsymbol{\delta}_{\tQ}(m)^{1/2}
\int_\K\int_{\N_\Q(\adef)}\int_{\gA_\Q}\fff^1(k\moins z m\,n\,k)\,u_{\tQ}(k,y)\dd z \dd n\dd k
\]
avec
\[
u_{\tQ}(k,y)= \int_{\ga_{\tG}\bs\ga_{\tQ}}\Gamma_{\tQ}(H,-\HO(k\,y))\dd H\ptf
\]
De plus, si $\fff$ est $\K$\hyph invariante \ie si $\fff(k\,x\,k\moins)=\fff(x)$ pour tout
$k\in\K$, alors
\[
 \fff_{\tQ,y}(x)=u_{\tQ}(y)\,\fff_{\tQ}(x)
\qquad\text{où $u_{\tQ}(y)=\int_\K u_{\tQ}(k,y)\dd k$}\ptf
\]

\begin{proposition}\label{conjug}
On a
\[
J^{T,\tG}_\bullet(\fff^y,\omega)=\sum_{\tQ} J^{T,\tQ}_\bullet(\fff_{\tQ,y},\omega)
\]
la somme portant sur les sous-ensembles paraboliques standard.
\end{proposition}

\begin{proof}
On utilisera les notations de la preuve du théorème~\ref{poly} et on pose
\[
 k^{T,\tG}_{\tQ}(\fff,\omega;x,X)=\Gamma_{\tQ}(\xT,X) \,k^{T,\tG}_{\tQ}(\fff,\omega;x)\ptf
\]
On observe que si $x=n\,m\,k$ est une décomposition d'Iwasawa on a
\[
\HO(x\,y)=\HO(x)+\HO(k\,y)
\]
d'où on déduit que,
\[
k^{T,\tG}(\fff^y,\omega;x)=\sum_{\tQ}\sum_{\eta\in\Q(F)\bs\G(F)}
k^{T,\tG}_{\tQ}(\eta\,x,\fff,-\HO(k\,y))\ptf
\]
On pose
\[
u_{\tQ}(x,y)=u_{\tQ}(k,y)
\]
si $x=n\,m\,k$, ce qui fournit
\[
\int_{\XG} k^{T,\tG}(\fff^y,\omega;x)\dd x= \sum_{\tQ} \int_{\YQ}
u_{\tQ}(x,y)\, k^{T,\tG}_{\tQ}(\fff,\omega;x)\dd x\ptf
\]
Il reste à observer que
\[
\int_{\YQ}
u_{\tQ}(x,y)\,
k^{T,\tG}_{\tQ}(\fff,\omega;x)\dd x= \int_{\XMQ} k^{T,\tQ}(m,\fff_{\tQ,y})\dd m
\]
\end{proof}

\section{La formule des traces grossière}

\begin{proposition}\label{indepPO}
Les sous-groupes $\MO$ et $\K$ étant fixés\textup,
la valeur du polynôme $J^{\T,\tG}_\bullet(\fff,\omega)$ évalué en $T=\TK$
est indépendante du choix de $\PO$.
\end{proposition}

\begin{proof}
Soit $\PO'$ un autre sous-groupe parabolique
minimal de sous-groupe de Levi $\MO$. Il existe $s\in\weyl$ représenté par $w_s$ tel que
\[
\PO'=w_s\moins\PO\, w_s\ptf
\]
Si
$x=n'm'k'$ est une décomposition d'Iwasawa relative à $\PO'$ on pose
\[
\HO'(x)=\HO(m')\ptf
\]
mais on a
\[
\HO(w_s\,x)=\HO(w_s\,n'\,m'\,w_s\moins\,w_s\,k')=s\HO(m')+\HO(w_s)
\]
et donc
\[
\HO'(x)=s\moins(\HO(w_s\,x)-\HO(w_s))\ptf
\]
Mais d'après le lemme~\ref{Yorth} on a
\[
\HO(w_s)=\TK-s\TK\ptf
\]
et donc
\[
\HO'(x)-\TK=s\moins(\HO(w_s\,x)-\TK)
\]
ce qui implique par exemple que si $\tP'=w_s\moins\tP\,w_s$ on a
\[
\htau_{\tP}(\HO(w_s\,x)-\TK)=\htau_{\tP'}(\HO'(x)-\TK)\ptf
\]
On conclut en observant que
\[
K_{\tP}(w_s\,x,\,w_s\,y)=K_{\tP'}(x,y)\ptf\qedhere
\]
\end{proof}

La valeur de $J^{\T,\tG}_\bullet(\fff,\omega)$ en $T=\TK$ sera notée $J^{\tG}_\bullet(\fff,\omega)$.

\begin{theoreme}\label{gross}  Soit $\fff\in\ctyc(\Gadef)$ qui est $\K$\hyph finie.
La forme grossière de la formule des traces est l'identité\textup:
\[
\sum_\go J^{\tG}_\go(\fff,\omega)=\sum_\chi J^{\tG}_\chi(\fff,\omega)\ptf
\]
La somme sur $\go$ ne comporte qu'un nombre fini de termes non nuls
\textup(dépendant du support de $\fff$\textup). Les divers termes sont indépendants du choix de $\PO$
lorsque $\MO$ et $\K$ sont fixés.
\end{theoreme}

\begin{proof}
On rappelle que, compte tenu de l'identité fondamentale~\ref{fonfon}:
\[
k^\T_{\mathrm{geom}}(x)=k_{\mathrm{spec}}^\T(x)
\]
on a
\[
\sum_{\go\in\gO} k^\T_\go(x)=\sum_{\chi\in\XX}k^\T_\chi(x)\ptf
\]
Pour $\T$ assez régulier on sait, d'après les théorèmes~\ref{geoconv} (ou~\ref{geoconvc} si on préfère)
et~\ref{chiconv},
que
\[
\sum_{\go\in\gO}\int_{\XG} \lvert k^\T_\go(x)\rvert\dd x< \infty
\Qquad{et}
\sum_{\chi\in\XX}\int_{\XG}\lvert k^\T_\chi(x) \rvert \dd x < \infty\ptf
\]
On a donc pour $\T$ assez régulier
\[
\sum_{\go\in\gO}\int_{\XG} k^\T_\go(x)\dd x=
\sum_{\chi\in\XX}\int_{\XG}k^\T_\chi(x) \dd x
\]
ce qui fournit l'identité de polynômes en $\T$:
\[
\sum_\go J^{\T,\tG}_\go(\fff,\omega)=\sum_\chi J^{\T,\tG}_\chi(\fff,\omega)\ptf
\]
Son évaluation en $T=\TK$ fournit l'identité cherchée.
La finitude de la somme sur $\go$ résulte du théorème~\ref{geoconv}. L'indépendance du choix de $\PO$
résulte de la proposition~\ref{indepPO}.
\end{proof}

\part{Forme explicite des termes spectraux}

\chapter[Introduction d'une fonction $B$]{\mathversion{bold}Introduction d'une fonction $B$}\label{IntroB}

{On rappelle que dans toute la suite de ce livre les fonctions $\fff$ sont suppos\'ees $\K$\hyph finies (\`a droite et \`a gauche).}

\section{La formule de départ}\label{formdep}

Soit $\Q$ un sous-groupe parabolique de $\G$. On rappelle que l'on a posé
\[
\pQ=\theta_0^{-1}(Q),\qquad\Qdo=Q\cap \pQ\Qquad{et}
\YQdo=\gA_\G\Qdo(F)\bs\Gadef\ptf
\]
Si $S\subset\pQ$ est un sous-groupe parabolique
on note $n^{\pQ}(S)$ le nombre de chambres dans $a_S^{\pQ}$.
Soit maintenant $\chi$ une donnée cuspidale. On reprend les notations de la section~\ref{estimm}; en particulier
on note
\[
\base^{\pQ}_{\chi}(\sigma)
\]
une base
formée de vecteurs $\K$\hyph finis de type $\chi$
dans la composante isotypique $\autom(\X_{\pQ},\sigma)$.

\begin{proposition}\label{debut}
Le polynôme $J^\T_{\chi}(\fff,\omega)$\textup, introduit au théorème~\ref{poly}\textup,
admet la décomposition spectrale suivante\textup:
{\multlinegap0pt
\begin{multline*}
J^\T_{\chi}(\fff,\omega)\\
\shoveleft{=\sum_{\{Q,R\mid \PO\subset Q\subset R\}}\tvedQR \int_{\YQdo}\tsQR(\HQ(x)-T)}\\
{}\times\smash[b]{\Biggl(\sum_{\{S\mid \PO\subset S\subset\pQ\}}\mspace{-5mu}
\frac{1}{n^{\pQ}(S)}\sum_{\sigma\in \Pi_{\mathrm{disc}}(M_S)}\sum_{\WPhi\in \base^{\pQ}_{\chi}(\sigma)}
\int_{\ima(\ga_S^\G)^*}}\mspace{-20mu}\tronc^{T,Q}E^Q(x,\treg_{S,\sigma,\mu}\bigl((f,\omega)\WPhi,\theto\mu\bigr)\\
{}\times\widebar{E^{\pQ}(x,\WPhi,\mu)}\dd \mu\Biggr)\dd x
\end{multline*}}%
pourvu que $\dPO(T)\geq c(f)$.\,\footnote{On observera que l'on n'affirme pas la convergence absolue de l'intégrale multiple
mais simplement la convergence des sommations itérées dans l'ordre indiqué.
Un meilleur contrôle de la convergence est l'objet de la section~\ref{iter}.}
\end{proposition}

\begin{proof}
D'après le théorème~\ref{poly} le polynôme $J^\T_{\chi}(\fff,\omega)$ admet l'expression suivante pour $\T$ assez régulier:
\[
J^\T_{\chi}(\fff,\omega)=\int_{\XG} k^\T_{\chi}(\fff,\omega;x)\dd x
\]
soit encore, suivant le théorème~\ref{chiconv},
\[
J^\T_{\chi}(\fff,\omega)= \sum_{\{Q,R\mid \PO\subset Q\subset R\}}
\tvedQR\,A_\Q^\R(\T,\fff,\omega,\chi)
\]
où
\[
A_\Q^\R(\T,\fff,\omega,\chi)=\int_{\YQdo} \tsQR (\xT)
\tronc_1^{T,\Q } K_{\Q,\delto,\chi}(x,x)\dd x\ptf
\]
Ici $K_{\Q,\delto,\chi}$ est la restriction du noyau $K_{\Q,\delto}$
à $L^2_{\chi}(\XpQG)$.
La décomposition spectrale (\cf proposition~\ref{majo}) fournit pour le noyau
$K_{\Q,\delto,\chi}$ une expression de la forme suivante:
\[
K_{\Q,\delto,\chi}(x,y)=\sum_{\M\in\Levi^{\pQ}/\weyl^{\pQ}}
\frac{1}{w^{\pQ}(\M)}\sum_{\sigma\in\Pi_{\mathrm{disc}}(M)_{\chi}}
\int_{i(\ga_\M^\G)^*} K_{\Q,\pQ,\sigma}(x,y;\mu)\dd \mu
\]
avec
\[
K_{\Q,\pQ,\sigma}(x,y;\mu)=\sum_{\WPhi\in\base^{\pQ}_{\chi}(\sigma)}
E^\Q(x,\treg_{S,\sigma,\mu}(\fff,\omega)\WPhi,\theto\mu)\widebar{E^{\pQ}(y,\WPhi,-\widebar\mu)}\ptf\qedhere
\]
\end{proof}

On remarquera que puisque
l'ensemble $\base^{\pQ}_{\chi}(\sigma)$ est une base de vecteurs $\K$\hyph finis de type $\chi$
dans la composante isotypique $\autom(\X_{\pQ},\sigma)$ et que $\fff$ est supposée $\K$\hyph finie, la somme
sur $\Phi$
est une somme finie: en effet il n'y a qu'un nombre fini de $\WPhi$ pour lesquels
\[
\treg_{S,\sigma,\mu}(f,\omega)\WPhi\ne 0\ptf
\]
De plus, les résultats de
Langlands sur la décomposition spectrale montrent de plus qu'il n'y a qu'un nombre fini de $\sigma$ pour lesquels
$\base^{\pQ}_{\chi}(\sigma)$ est non vide.

On aura besoin d'une variante de cette proposition faisant intervenir les multiplicateurs d'Arthur, imitant en cela
les sections~3 et~4 de \cite{AeisI} dont on rappelle brièvement le contenu.
On considère le groupe de Lie
$\G_\infty=\G(F\otimes\RM)$.
Considérons
\[
\cgh=\ima\cgh_\K\oplus\cgh_0
\]
où $\cgh_\K$ est
une sous-algèbre de Cartan du sous-groupe compact maximal
$\K_\infty$
et $\cgh_0$ l'algèbre de Lie d'un tore déployé maximal de $\G_\infty$.
En particulier $\cgh_\CM$ est une sous-algèbre de Cartan pour $\frakg_\CM$.
On notera $\weyl_\CM$ le groupe de Weyl complexe de $\G_\infty$
et $w_\CM$ son ordre.
On dispose de la théorie des multiplicateurs d'Arthur ce qui permet de
construire des fonctions $\fff_X\in\ctyc\bigl(\Gadef\bigr)$
pour $X\in\cgh$ vérifiant
\[
\pi(\fff_X)=\ee^{\nu_\pi(X)}\pi(\fff)=\frac{1}{w_\CM}
\sum_{s\in\weyl_\CM}\ee^{\langle\nu_\pi, s\moins X\rangle}\pi(\fff)
\]
pour toute représentation admissible irréductible $\pi$ et où
$\nu_\pi$ est le caractère infinitésimal de $\pi_\infty$.
On considère ici $\nu_\pi$ soit comme une forme linéaire $\weyl_\CM$\hyph invariante
sur $\cgh$ soit comme un élément de $\cgh^*\otimes\CM$. L'extension au cas tordu est immédiate.
En particulier on a
\[
\treg_{S,\sigma,\mu}(\fff_X,\omega)=\frac{1}{w_\CM}
\sum_{s\in\weyl_\CM}\ee^{\langle\nu_\sigma+\mu, s\moins X\rangle}\treg_{S,\sigma,\mu}(\fff,\omega)
\ptf
\]

\begin{corollaire}\label{lemmeA}
Pour tout $X\in \cgh$\textup, si $\dPO(\T)\ge c(\fff)(1+\lVert X\rVert)$ on a
\[
p^\T_{\chi}(X)=J^\T_{\chi}(\fff_X,\omega)=\sum_{\{Q,R\mid \PO\subset Q\subset R\}} \tvedQR\,A_\Q^\R(\T,\fff_X,\omega,\chi)
\ptf
\]
et
{\multlinegap0pt\begin{multline*}
J^\T_{\chi}(\fff_X,\omega)\\
\shoveleft{=\sum_{\{Q,R\mid \PO\subset Q\subset R\}}\tvedQR
\int_{\YQdo}\tsQR(\HQ(x)-T)}\\
\smash[b]{{}\times\Biggl(\sum_{\{S\mid \PO\subset S\subset\pQ\}}\mspace{-5mu}
\frac{1}{n^{\pQ}(S)}\sum_{\sigma\in \Pi_{\mathrm{disc}}(M_S)}\sum_{\WPhi\in \base^{\pQ}_{\chi}(\sigma)}
\int_{\ima(\ga_S^\G)^*}}\mspace{-10mu}\tronc^{T,Q}E^Q(x,\treg_{S,\sigma,\mu}(f_X,\omega)
\WPhi,\theto\mu)\\
{}\times\widebar{E^{\pQ}(x,\WPhi,\mu)}\dd \mu\Biggr)\dd x\ptf
\end{multline*}}%
\end{corollaire}

\begin{proof} 
Ceci résulte de la proposition~\ref{debut} compte tenu de la dépendance en $X$
du support de $\fff_X$. On renvoie le lecteur à \cite{AeisI}*{Proposition~3.1} pour un énoncé précis
de cette dépendance.
\end{proof}

\section{Estimations}\label{estimations}

Dans cette section on établit des raffinements des estimées~\ref{premaj} et~\ref{majo}.

\begin{lemme}\label{PREMAJbis}
Soit $h $ une fonction à support compact sur $G(\adef)$\textup, à valeurs${}\geq0$. Alors il existe $c>0$ tel
que
\[
\sum_{\gamma\in G(F)}\int_{\gA_\G}h(x^{-1}z\gamma y)\dd z\le c\delta_{\PO}(x)^{1/2}\delta_{\PO}(y)^{1/2}
\]
pour tous $x,y\in \Sieg^G$.
\end{lemme}

\begin{proof}
Soit $h $ une fonction à support compact sur $G(\adef)$, à valeurs${}\geq0$, et soit $x,y\in \Sieg^G$. On veut
évaluer
\[
\sum_{\gamma\in G(F)}\int_{\gA_{G}}h(x^{-1}z\gamma y)\dd z\ptf
\]
Il suffit d'évaluer le nombre de $\gamma$ tel que $x^{-1}\gamma y\in \Omega$, où $\Omega$ est l'intersection du support de
$h$ avec $G(\adef)^1$. Cet ensemble $\Omega$ est
compact. D'après le lemme~\ref{xycomp}
\[
\HO(x)-\HO(y)
\]
appartient à un
compact.
Quitte à agrandir $\Omega$, on peut
donc supposer $x=y$. Fixons un élément régulier $T_1$ et utilisons la partition~\ref{FPQ}: il existe un
unique parabolique standard $R$ tel que
\begin{equation}
F_{\PO}^R(x,T_1)
\tau_{R}(\HO(x)-T_1)=1\ptf\tag{$*$}\label{eq12.*}
\end{equation}

Si $x^{-1}\gamma x\in \Omega$, on a $\gamma x\in x\Omega$ et
quitte à agrandir encore $\Omega$, on peut supposer $x \in \gA_0(t)$. On a donc
\[
\tau_{R}(\HO(\gamma x)-T_{2})=1
\]
pour un $T_{2}\in T_1+\HO
(\Omega)$. En prenant $T_1$ assez grand, le lemme~\ref{unique} implique que $\gamma\in R(F)$.
On a déjà supposé $x\in\gA_0(t)$, on peut écrire $x=\ee^H$ avec $H\in\gao$. La condition~\eqref{eq12.*}
entraîne que $H^\R$ reste dans un compact. La condition $x\moins\gamma x\in\Omega$ entraîne
donc
\[
\ee^{-H_\R}\gamma \ee^{H_\R}\in\Omega'
\]
où $\Omega'$ est un compact plus gros. En notant $M_{R}$ le Levi
standard de $R$, on est ramené à évaluer le nombre
de
\[
(\delta,\eta)\in M_{R}(F)\times N_{R}(F)
\]
tels que $x^{-1}\delta\eta x\in \Omega$. Puisque $\ee^{H_\R}$ commute à $\delta$, cela
entraîne que $\delta$ reste dans un compact indépendant
de $\ee^{H_\R}$. Ces $\delta$ sont en nombre fini et on est ramené à évaluer le nombre de $\eta\in N_{R}(F)$ tels que $x^{-1}\eta x\in
C$, où $C$ est un sous-ensemble compact de $N_{R}
(\adef)$. Par l'exponentielle, on descend à l'algèbre de Lie $\gn_{R}$ et on doit évaluer le nombre de $X\in
\gn_{R}(F)$ tels que $\ad(x)^{-1}(X)\in \mathfrak{C}$,
où $\mathfrak{C}$ est un sous-ensemble compact de $\gn_{R}(\adef)$. On peut majorer la fonction
caractéristique de $\mathfrak{C}$ par une fonction $\psi\in C^{\infty}
\bigl(\gn_{R}(\adef)\bigr)$ à valeurs positives ou nulles. Notre nombre d'éléments est majoré par
\[
\sum_{X\in \gn_{R}(F)}\psi\bigl(\ad(x)^{-1}(X)\bigr)\ptf
\]
On utilise la formule de Poisson, en identifiant le dual de $\gn_{R}$ à l'algèbre opposée $\gn_{\bar{R}}$. La
transformée de Fourier de $\psi\circ \ad(x)^{-1}$ est $\delta_{R}(x)\hat{\psi}\circ \ad(x)^{-1}$. La somme ci-dessus est égale à
\[
\delta_{R}(x)\sum_{X\in \gn_{\bar{R}}(F)}\widehat{\psi}\bigl(\ad(x)^{-1}(X)\bigr)\ptf
\]
Puisque $x\in \gA_0(t)$, $\ad(x)^{-1}$ dilate $\gn_{\bar{R}}$ et la dernière série est bornée
indépendamment de $x$. On obtient une majoration par $\delta_{R}(x)$.
Puisque $F_{\PO}^R(x,T_1)=1$, ce terme est lui-même essentiellement borné \,\footnote{\og essentiellement
borné\fg signifie pour nous qu'il existe $c$ tel que le terme soit
majoré par $c\delta_{\PO}(x)$} par $\delta_{\PO}(x)$.
Enfin, puisque $xy^{-1}$ reste dans un compact, ce dernier terme est essentiellement borné
par $\delta_{\PO}(x)^{1/2}\delta_{\PO}(y)^{1/2}$.
\end{proof}

On rappelle que ${(\ga_S^\G)^*}$ est naturellement un sous-espace de $\cgh^*$
(avec $\cgh$ comme dans la section~\ref{formdep}).
Notons $\cgh^{S,*}$ son orthogonal. À la représentation $\sigma$ est associé un paramètre
$\lambda(\sigma)\in(\cgh_\CM^S)^*$.

\begin{lemme}\label{W1.2.4}
Soit $\varphi$ une fonction de Paley-Wiener sur $(\ga_{S,\CM}^G)^*$. Alors il existe une
fonction $\phi$ sur $\cgh_\CM^*$ vérifiant
les conditions suivantes\textup:
\begin{enumerate}[(i)]
\item $\phi$ est de Paley-Wiener\textup;
\item $\phi$ est invariante par $W_{\CM}$\textup;
\item pour tout $\mu\in \ima(\ga_S^\G)^*$\textup, $\phi(\lambda(\sigma)+\mu)\ge\lvert \varphi(\mu)\rvert^2$.
\end{enumerate}
\end{lemme}

\begin{proof}
On choisit une fonction de Paley-Wiener $\varphi^S$ sur $(\cgh_\CM^S)^*$ telle que
\[
\varphi^S\bigl(\lambda(\sigma)\bigr)=1\ptf
\]
On définit une fonction $\phi_1$ sur
$\cgh_\CM^*$ par
\[
\phi_1(\nu)=\varphi(\nu_S)\varphi^S(\nu^S)
\]
pour $\nu\in \cgh_\CM^*$, où
$\nu_S$ et $\nu^S$ sont les projections orthogonales de
$\nu$ sur $(\ga^G_{S,\CM})^*$ et $(\cgh_\CM^S)^*$. Décomposons $\lambda(\sigma)$ en ses
parties réelles et imaginaires: $\lambda(\sigma)=X(\sigma)
+\ima Y(\sigma)$. Notons $W'\subset W_{\CM}$ le fixateur de $X(\sigma)$ dans $W_{\CM}$. Pour $w\notin W'$, on a
\[
w^{-1}\bigl(X(\sigma)\bigr)^S\ne X(\sigma)
\]
(sinon, par
comparaison des normes, on a
\[
w^{-1}\bigl(X(\sigma)\bigr)=w^{-1}\bigl(X(\sigma)\bigr)^S=X(\sigma)
\]
et $w\in W'$). On peut donc choisir
$\beta_{w}\in \cgh^S$ tel que $\beta_{w}\Bigl(w^{-1}\bigl(X(\sigma)\bigr)\Bigr)
\ne\beta_{w}\bigl(X(\sigma)\bigr)$. Définissons $\phi_{2}$ par
\[
\phi_{2}(\nu)=\phi_1(\nu)\prod_{w\notin W'}\Bigl(\beta_{w}\bigl(w^{-1}(\nu)\bigr)-\beta_{w}\bigl(\lambda(\sigma)\bigr)\Bigr),
\]
puis $\phi_{3}$ par
\[
\phi_{3}(\nu)=\phi_{2}(\nu)\widebar{\phi_{2}\bigl(-\bar{\nu}+2X(\sigma)\bigr)},
\]
enfin $\phi$ par
\[
\phi(\nu)=\sum_{w\in W_{\CM}}\phi_{3}\bigl(w(\nu)\bigr)\ptf
\]
Cette fonction vérifie évidemment les deux premières conditions de l'énoncé. Vérifions la dernière, soit donc
$\nu=\lambda(\sigma)+\mu$, avec $\mu\in \ima(\ga_S^\G)^*$. Pour
$w\notin W'$ le terme $\phi\bigl(w(\nu)\bigr)$ contient le facteur $\beta_{w}(\nu)-\beta_{w}\bigl(\lambda(\sigma)\bigr)$. Puisque $\beta_{w}$ est
orthogonal à ${(\ga_S^\G)^*}$, on a $\beta_{w}(\nu)=
\beta_{w}\bigl(\lambda(\sigma)\bigr)$ et le facteur précédent est nul. Donc
\[
\phi(\nu)=\sum_{w\in W'}\phi_{3}\bigl(w(\nu)\bigr)\ptf
\]
Pour $w\in W'$, la partie réelle de $w(\nu)$ est $X(\sigma)$. Donc
\[
-\widebar{w(\nu)}+2X(\sigma)=w(\nu)
\]
et
\[
\phi_{3}\bigl(w(\nu)\bigr)=\phi_{2}\bigl(w(\nu)\bigr)\widebar{\phi_{2}\bigl(w(\nu)\bigr)}\geq0\ptf
\]
On peut abandonner les $w\ne 1$: $\phi(\nu)\geq \phi_{3}(\nu)=\lvert \phi_{2}(\nu)\rvert ^2$.
Pour $w\notin W'$, la partie réelle de
\[
\beta_{w}\bigl(w^{-1}(\nu)\bigr)-\beta_{w}\bigl(\lambda(\sigma)\bigr)
\]
est
$\beta_{w}\Bigl(w^{-1}\bigl(X(\sigma)\bigr)\Bigr)-\beta_{w}\bigl(X(\sigma)\bigr)$ qui n'est pas nulle.
Donc
\[
\bigl\lvert \beta_{w}\bigl(w^{-1}(\nu)\bigr)-\beta_{w}\bigl(\lambda(\sigma)\bigr)\bigr\rvert 
\]
est minoré par un nombre strictement positif.
On en déduit une minoration
\[
\lvert \phi_{2}(w)\rvert \geq c \lvert \phi_1(\nu)\rvert =c\lvert \varphi(\mu)\rvert
\]
pour un certain $c>0$. D'où
\[
\phi(\nu)\geq c^2\lvert \varphi(\mu)\rvert ^2\ptf\qedhere
\]
\end{proof}

On fixe $\Q$ et $\R$ avec $\tvedQR\ne0$. Cette condition équivaut à ce qu'il existe un et un seul parabolique, que
l'on note $P$, avec $Q\subset P\subset R$ et
$\theto(P)=P$. On fixe $S$ et $\sigma$. Pour $\WPsi\in \base^{\pQ}_{\chi}(\sigma)$, on peut écrire
\[
\treg_{S,\sigma,\mu}(f,\omega)
\WPsi=\sum_{\WPhi\in \mathcal{B}^Q(\theto\sigma)_{\chi}}\hat{f}_{\WPhi,\WPsi}(\mu)\WPhi\;,
\]
où la somme est finie et $\hat{f}_{\WPhi,\WPsi}$ est une fonction de Paley-Wiener sur $(\ga_{S,\CM}^G)^*$.
L'expression souhaitée pour $J^\T_{\chi}(\fff,\omega)$ est donc combinaison linéaire d'intégrales itérées
\begin{equation}
 \int_{\YQdo}\tsQR(\HQ(x)-T)\biggl(\int_{\ima(\ga_S^\G)^*}
\tronc^{T,Q}E^Q(x,\WPhi,\theto\mu)\widebar{E^{\pQ}(x,\WPsi,\mu)}
\varphi(\mu)\dd \mu\biggr)\dd x\label{eq12.1}
\end{equation}
où $\varphi$ est une fonction de Paley-Wiener sur $(\ga_{S,\CM}^G)^*$
et dont nous devons montrer la convergence.

On note $L$, $L'$, $L_0$ les Levi standard et $N_Q$, $N_{\pQ}$, $N_{\Qdo}$ les radicaux unipotents de $Q$, $\pQ$ et $Q_
0$. On fixe un ensemble de Siegel $\Sieg^L$ pour le
quotient $L(F)\backslash L(\adef)^1$. On choisit un sous-ensemble compact
$\Omega_{N_Q}\subset N_Q(\adef)$ tel que $N_Q(\adef)=N_Q(F)\Omega_{N_Q}$
et on pose
\[
\Sieg^Q=\Omega_{N_Q}\gA^G_Q\Sieg^LK\ptf
\]
C'est un ensemble de Siegel pour le quotient $Q(F)\backslash G(\adef)^1$. On
introduit de même des ensembles $\Sieg^{L'}$, $\Sieg^{\pQ}$, $\Sieg^{L_0}$, $\Sieg^{\Qdo}$ et
un ensemble $\Sieg^G$.

\begin{proposition}\label{eqqvf}
Soient $\WPhi$\textup, $\WPsi$ et $\varphi$ comme ci-dessus. Alors il existe $c>0$ tel que
\[
\int_{\ima(\ga_S^\G)^*}\lvert E^Q(x,\WPhi,\theto\mu)
\widebar{E^{\pQ}(y,\WPsi,\mu)}\varphi(\mu)\rvert \dd \mu\leq c\delta_{\PO}(x)^{1/2}\delta_{\PO}(y)^{1/2}
\]
pour tous $x\in \Sieg^Q$ et $y\in \Sieg^{\pQ}$.
\end{proposition}

\begin{proof}
Posons
\[
J(x,y,\WPhi,\WPsi,\varphi)=\int_{\ima (\ga_S^\G)^*}\lvert E^Q(x,\WPhi,\theto\mu)
\widebar{E^{\pQ}(y,\WPsi,\mu)}\varphi(\mu)\rvert \dd \mu\ptf
\]
L'élément $\WPhi$ est $K$\hyph fini. La fonction
\[
x\mapsto \delta_Q(x)^{-1/2}E^Q(x,\WPhi,\theto\mu)
\]
est invariante à gauche par
$N_Q(\adef)$ et sa valeur absolue est invariante à
gauche par $\gA_Q$. Les termes relatifs à $\pQ$ vérifient des propriétés similaires. Il en résulte l'existence d'un
nombre fini de couples $(\WPhi_{i},\WPsi_{i})$ tels que pour
$n\ee^Hxk\in \Sieg^Q$, avec $u\in N_Q(\adef)$, $H\in \ga_Q$, $x\in \Sieg^L$, $k\in K$, et pour un
élément similaire $u'\ee^{H'}yk'\in \Sieg^{\pQ}$, on
ait une majoration
\[
J(n\ee^Hxk,u'\ee^{H'}yk',\WPhi,\WPsi,\varphi)\leq \delta_Q(\ee^H)^{1/2}
\delta_{\pQ}(\ee^{H'})^{1/2}\sum_{i}J(x,y,\WPhi_{i},\WPsi_{i},\varphi)\ptf
\]
On peut donc se limiter à majorer $J(x,y,\WPhi,\WPsi,\vf)$ pour $x\in \Sieg^L$ et $y\in \Sieg^{L'}$. D'après un
théorème de Dixmier et Malliavin \cite{DiM}, on peut décomposer $\varphi$
en somme finie de produits de deux fonctions de Paley-Wiener. Cela nous ramène au cas où $\varphi$ est produit de deux
telles fonctions $\varphi_1$ et $\varphi_{2}$. Par l'inégalité
de Schwartz, on voit que
\[
 J(x,y,\WPhi,\WPsi,\varphi)^2
\]
est majoré par
\[
\int_{\ima(\ga_S^\G)^*}\lvert E^Q(x,\WPhi,\theto\mu)\varphi_1(\mu)\rvert ^2\dd \mu\int_{\ima
(\ga_S^\G)^*}\lvert E^{\pQ}(y,\WPsi,\mu)\varphi_{2}(\mu)
\rvert ^2\dd \mu\ptf
\]
Les deux facteurs sont du même type, à quelques changements inessentiels près. Cela nous ramène à prouver que, pour
$\WPsi$ et $\varphi$ fixés, il existe $c$ tel que
\[
 \int_{\ima(\ga_S^\G)^*}\lvert E^{\pQ}(y,\WPsi,\mu)\varphi(\mu)\rvert^2 \dd \mu \leq c\delta_{\PO}(y)
\]
pour tout $y\in \Sieg^{L'}$. Il existe $\WPhi^{L'}$ appartenant à une induite convenable pour le groupe $L'$ tel que
\[
E^{\pQ}(y,\WPsi,\mu)=E^{L'}(y,\WPhi^{L'},\mu)
\]
pour tout $y\in \Sieg^{L'}$.
On peut donc aussi bien supposer ici $\pQ=G$. On supprime alors les primes et on pose
\[
J(y,\WPhi,\varphi^2)= \int_{\ima(\ga_S^\G)^*}\lvert E^{G}(y,\WPhi,\mu)\varphi(\mu)\rvert^2 \dd \mu\ptf
\]
De nouveau, on peut supposer que $\varphi=\varphi_1\varphi_{2}$ et on utilise l'inégalité de Schwartz:
\[
J(y,\WPhi,\varphi^2)^2\leq J(y,\WPhi,\varphi_1^4)J(y,\WPhi,\varphi_{2}^4)\ptf
\]
Cela nous ramène au cas où $\varphi$ est un carré, disons $\varphi=\varphi_0^2$. On applique le lemme~\ref{W1.2.4} à
$\varphi_0$, soit $\phi$ la fonction qui s'en déduit. Alors
\begin{equation}
J(y,\WPhi,\varphi_0^4)\leq \int_{\ima(\ga_S^\G)^*}E^G(y,\WPhi,\mu)\phi(\lambda(\sigma)+\mu)
\widebar{E^G(y,\WPhi,\mu)\phi(\lambda(\sigma)+\mu)}d\mu\ptf\label{eq12.2}
\end{equation}
On peut inclure $\WPhi$ dans un sous-espace $V$ de dimension finie de l'induite qui est une somme de composantes isotypiques
pour l'action de $K$. Soit $V'$ le supplémentaire de $V$
invariant par $K$. D'après \cite{CD}*{théorème~3} (qui est une conséquence d'un théorème d'Arthur), il existe une fonction
$h\in \ctyc\bigr(G(\adef)\bigr)$ telle que $\reg_{S,\sigma,\mu}(h)$
(il s'agit de l'action non tordue usuelle) annule $V'$ et agisse sur $V$ par multiplication par
$\phi(\lambda(\sigma)+\mu)$. Considérons alors le noyau $K_{G}((h^*\star h)
^1;y,y)$. Son expression spectrale est somme de termes tous positifs ou nuls et l'un d'eux est le membre de droite de la relation
\eqref{eq12.2}. On en déduit l'inégalité
\[
J(y,\WPhi,\varphi_0^4)\leq K_{G}((h^*\star h)^1;y,y)\ptf
\]
Mais
\[
K_{G}((h^*\star h)1;y,y)=\sum_{\gamma\in G(F)}\int_{\gA_{G}}h^*\star h(y^{-1}z\gamma y)\dd z\ptf
\]
Il reste à appliquer le lemme~\ref{PREMAJbis}
pour obtenir la majoration cherchée.
\end{proof}

\begin{proposition}\label{psipsivf}
Soient $\WPhi$\textup, $\WPsi$ et $\varphi$ comme ci-dessus. Alors l'intégrale
\[
\int_{\ima(\ga_S^\G)^*} \tronc^{T,Q}E^Q(x,\WPhi,\theto\mu)\widebar{E^{\pQ}(y,\WPsi,\mu)}\varphi(\mu) \dd \mu
\]
est absolument convergente. Il existe un réel $D$ et\textup, quel que soit le réel $r$\textup, il existe $c>0$ tel que cette intégrale soit
majorée par
\[
c\lvert y\rvert^D\ee^{D\lVert H\rVert }\lvert x^L\rvert ^{-r}
\]
pour tout $y$ et tout $x=\ee^Hx^L$, avec $H\in \ga_Q$ et $x^L\in \Sieg^L$.\,\footnote{Ici, $T$ est considéré comme fixé; on ne se demande pas comment le $c$ ci-dessus dépend de~$T$.}
\end{proposition}

\begin{proof}
L'opérateur $\tronc^{T,Q}$ est une combinaison d'intégrales sur des compacts et de sommes finies, affectées de signes.
Considérons l'opérateur (idiot) où on supprime les signes,
notons-le $\tronc^{T,Q}_{+}$. Il est clair que
\[
\lvert \tronc^{T,Q}(h)\rvert\leq \tronc^{T,Q}_{+}(\lvert h\rvert )
\]
pour toute fonction
$h$ sur $Q(F)\backslash G(\adef)$. Alors
l'expression
\[
\int_{\ima(\ga_S^\G)^*} \lvert \tronc^{T,Q}E^Q(x,\WPhi,\theto\mu)\widebar{E^{\pQ}(y,\WPsi,\mu)}\varphi(\mu)\rvert \dd \mu
\]
est majorée par l'image par cet opérateur $\tronc^{T,Q}_{+}$ de la fonction
\[
x\mapsto J(x,y,\WPhi,\WPsi,\varphi)\ptf
\]
Or cette image
est définie par une intégrale convergente, grâce à la
proposition~\ref{eqqvf} et parce que, comme on vient de le dire, $\tronc^{T,Q}_{+}$ ne fait intervenir que des intégrales sur des
compacts et des sommes finies.
Considérons maintenant la même expression sans les valeurs absolues:
\[
\int_{\ima(\ga_S^\G)^*} \tronc^{T,Q}E^Q(x,\WPhi,\theto\mu)\widebar{E^{\pQ}(y,\WPsi,\mu)}\varphi(\mu) \dd \mu\ptf
\]
Elle est obtenue en appliquant $\tronc^{T,Q}$ (portant sur la variable $x$) à l'expression
\[
\int_{\ima(\ga_S^\G)^*} E^Q(x,\WPhi,\theto\mu)\widebar{E^{\pQ}(y,\WPsi,\mu)}\varphi(\mu) \dd \mu\ptf
\]
Celle-ci est à croissance modérée en les deux variables d'après la proposition~\ref{eqqvf}. Ses dérivées en $x$ sont des
expressions similaires, la fonction est donc uniformément à
croissance modérée. La majoration de l'énoncé se déduit alors de la proposition~\ref{rapdec}.
\end{proof}

\section{Convergence d'une intégrale itérée}\label{iter}

Dans la suite
le terme $T$ est un élément régulier de $\ga^G_0$, fixe par $\theto$\,\footnote{Il faut garder en mémoire que l'on va \emph{in fine} évaluer les polynômes
en $\T=\TK$ qui n'est pas nécessairement $\theto$\hyph invariant.
Mais ce n'est pas une difficulté étant donné que les polynômes
à évaluer ne dépendent que de la projection de $\T$
sur les invariants}.
On le limite à un domaine défini par des inégalités
\[
c_1<\alpha(T)\leq c_{2} \dPO(T)
\]
pour tout $\alpha\in \Delta_0$, où $c_1$ et $c_{2}$ sont des réels strictement
positifs arbitraires mais fixés (avec $\Delta_0=\Delta_{\PO}$). Dans un tel domaine, les fonctions $\dPO(T)$, $\lVert
T\rVert$ et $\alpha(T)$ pour $\alpha\in \Delta_0$ sont équivalentes.

Considérons l'intégrale itérée
\[
\int_{\YQdo}\tsQR(\HQ(y)-T)\biggl\lvert \int_{\ima(\ga_S^\G)^*}\tronc^{T,Q}E^Q(y,\WPhi,\theto\mu)
\widebar{E^{\pQ}(y,\WPsi,\mu)}\varphi(\mu)\dd \mu\biggr\rvert\dd y
\]
où $\vf$ est une fonction de Paley-Wiener.
L'intégration sur $ \YQdo$ se décompose en une intégration sur le produit
\[
\bigr(N_{Q_0}(F)\backslash N_{\Qdo}(\adef)\bigr)\times (L_0(F)
\backslash L_0(\adef)^1)\times \ga^G_{\Qdo}\times K\ptf
\]
Par ailleurs, la fonction que l'on intègre est invariante
à gauche par $N_Q(\adef)\cap N_{\pQ}(\adef)$.
Posons $\Q_0^L=\Qdo\cap L$, $\Q_0^{L'}=\Qdo\cap L'$. L'application naturelle
\[
N_{\Qdo}(F)\bigl(N_Q(\adef)\cap N_{\pQ}(\adef)\bigr)\backslash N_{\Qdo}(\adef)
\to \bigl(N_{\Q_0^L}(F)\backslash N_{\Q_0^L}(\adef)\bigr)\times \bigl(N_{\Q_0^{L'}}(F)\backslash
N_{\Q_0^{L'}}(\adef)\bigr)
\]
est un isomorphisme. On peut aussi bien intégrer sur
\[
\bigl(N_{\Q_0^L}(F)\backslash N_{\Q_0^L}(\adef)\bigr)\times \bigl(N_{\Q_0^{L'}}(F)\backslash N_{\Q_0^{L'}}(\adef)\bigr)\times
(L_0(F)\backslash L_0(\adef)^1)\times \ga^G_{\Qdo}\times K\ptf
\]
On remplace $y$ par $nn'x\ee^Hk$ (avec $x$ appartenant à
$L_0(F)\backslash L_0(\adef)^1$). La mesure $\dd y$ se transforme en
$\delta_{\Qdo}(\ee^H)^{-1}\dd n\dd n'\dd x\dd H\dd k$. L'intégrale sur $K$ est inoffensive, on l'oublie. Posons
\begin{equation}
I(nn'x\ee^H)= \int_{\ima(\ga_S^\G)^*}\tronc^{T,Q}E^Q(nx\ee^H,\WPhi,\theto\mu)\widebar{E^{\pQ}(n'x\ee^H,\WPsi,\mu)}\varphi(\mu)
\dd \mu\ptf\label{eq12.1a}
\end{equation}

\begin{lemme}\label{essentiel}
Il existe un sous-ensemble compact $\omega\subset\ga_Q^G$ tel que\textup, si 
l'int\'egrale $I(nn'x\ee^H)$ est non nulle\textup,
alors \textup(dans les notations de la section~\ref{omni}\textup) $q(H)\in \omega$.
\end{lemme}

\addtocounter{equation}{1}
\begin{proof}
En effet, on a
{\multlinegap0pt
\begin{multline}\label{eq12.2a}
I(nn'x\ee^H)=\int_{\ima(\ga_S^{\pQ})^*}\int_{\ima\ga_{\pQ}^{G,*}}\tronc^{T,Q}
E^Q\bigr(nx\ee^H,\WPhi,\theto(\mu_{\pQ}+\mu^{\pQ})\bigr)\\
{}\times\widebar{E^{\pQ}(n'x\ee^H,\WPsi,\mu_{\pQ}+\mu^{\pQ})}\varphi(\mu_{\pQ}+\mu^{\pQ})\dd \mu_{\pQ}\dd \mu^{\pQ}\ptf
\end{multline}}%
On a aussi
\begin{multline*}
\tronc^{T,Q}E^Q\bigl(nx\ee^H,\WPhi,\theto(\mu_{\pQ}+\mu^{\pQ})\bigr)\\
=\delta_Q(\ee^{H_Q})^{1/2}\ee^{\langle H_Q,\theto(\mu_{\pQ})\rangle}
\tronc^{T,Q}E^Q\bigl(nx\ee^{H^Q},\WPhi,\theto(\mu^{\pQ})\bigr)
\end{multline*}
et
\[
\widebar{E^{\pQ}(n'x\ee^H,\WPsi,\mu_{\pQ}+\mu^{\pQ})} =\delta_{\pQ}(\ee^{H_{\pQ}})^{1/2}\ee^{-\langle H_{\pQ},\mu_{\pQ}\rangle}
\widebar{E^{\pQ}(n'x\ee^{H^{\pQ}},\WPsi,\mu^{\pQ})} \ptf
\]
L'expression pour $I(nn'x\ee^H)$ contient donc la sous-intégrale
\begin{multline}\label{eq12.3a}
\smash[b]{\int_{\ima\ga_{\pQ}^{G,*}}}\ee^{\langle H_Q,\theto(\mu_{\pQ})\rangle -\langle H_{\pQ},\mu_{\pQ}\rangle}\varphi(\mu_{\pQ}+\mu^{\pQ})\dd \mu_{\pQ}\\
=\int_{\ima\ga_{\pQ}^{G,*}}\ee^{\langle\theta_0^{-1}(H_Q)-H_{\pQ},\mu_{\pQ}\rangle}\varphi(\mu_{\pQ}+\mu^{\pQ})\dd \mu_{\pQ}\ptf
\end{multline}
C'est la transformée de Fourier partielle de la fonction $\varphi$, évaluée au point
\[
\theta_0^{-1}(H_Q)-H_{\pQ}\ptf
\]
Puisque $\varphi$ est de Paley-Wiener, cette transformée de Fourier est
à support compact. Plus précisément, il existe un sous-ensemble compact $\omega'\subset \ga_{\pQ}^G$ tel que, quel
que soit $\mu^{\pQ}$, le support de cette transformée de
Fourier soit inclus dans $\omega'$. Donc
\[
\theta_0^{-1}(H_Q)-H_{\pQ}\in \omega'
\]
ce qui équivaut à $q(H)\in \theto(\omega')$.
\end{proof}

\begin{proposition}
\label{tsqrpsi}
L'intégrale
\[
\int_{\YQdo}\tsQR(\HQ(y)-T)\biggl\lvert \int_{\ima(\ga_S^\G)^*}\tronc^{T,Q}E^Q(y,\WPhi,\theto\mu)
\widebar{E^{\pQ}(y,\WPsi,\mu)}\varphi(\mu)\dd \mu\biggr\rvert\dd y
\]
est convergente.
\end{proposition}

\addtocounter{equation}{3}
\begin{proof}
On doit prouver que le produit de $I(nn'x\ee^H)$ avec
\[
\tsQR(H -T)\delta_{\Qdo}(\ee^H)^{-1}
\]
est absolument intégrable.
On peut découper le domaine d'intégration en $H$ grâce à la partition du lemme~\ref{partition} appliquée au cas $P=\Qdo$ et $R=Q$. C'est-à-dire
que l'on peut fixer un sous-groupe parabolique $P'$
avec $\Qdo\subset P'\subset Q$ et imposer que
\[
\phi_{\Qdo}^{P'}(H-T)\tau_{P'}^Q(H-T)=1\ptf
\]
Compte tenu du lemme~\ref{ex10.3.2}, on en déduit la majoration
\begin{equation}
\lVert (H-T_{\Qdo})\rVert \ll 1 +\lVert (H-T)_{P'}^Q\rVert\label{eq12.4a}
\end{equation}
pour tous $T$, $H$ tels que
\[
\tsQR(H-T)\phi_{\Qdo}^{P'}(H-T)\tau_{P'}^Q(H-T)=1\Qquad{et}q(H)\in \omega\ptf
\]
Au lieu d'intégrer en
\[
x\in L_0(F)\backslash L_0(\adef)^1
\]
on peut intégrer sur $x\in\Sieg^{L_0}$. Pour tous
$n$, $x$ et $H$, on peut choisir $\gamma\in L(F)$ tel que
\[
y'=\gamma nx\ee^{H^Q}\in \Sieg^L\ptf
\]
D'après la proposition~\ref{psipsivf}, il existe $D$ tel que, pour tout $r$, on ait une majoration
\begin{equation}
\delta_{\Qdo}(\ee^H)^{-1} \lvert I(nn'x\ee^H)\rvert \ll \lVert x\ee^H\rVert ^D\ee^{-r\lVert \HO(y')
\rVert}\ptf\label{eq12.5a}
\end{equation}
Montrons que l'on a la relation
\begin{equation}
\lVert H_{P'}^Q\rVert +\lVert \HO(x)\rVert \leq 1+\lVert \HO(y')\rVert \ptf\label{eq12.6a}
\end{equation}
D'après le lemme~\ref{franke}, il existe $c\in \RM$ tel que
\[
\varpi_{\alpha}\bigl(\HO(\gamma\moins y')-\HO(y')\bigr)\leq c 
\]
pour tout $\alpha\in \Delta_0^\Q$ soit encore
\[
\varpi_{\alpha}\bigl(H^Q+\HO(x)-\HO(y')\bigr)\leq c \ptf
\]
C'est dire qu'il existe $H'$ tel que $H'-\HO(y')$ soit borné et tel que
\[
H'=H^Q+\HO(x)+X
\]
où $X$ est combinaison linéaire d'éléments à coefficients positifs ou
nuls de coracines $\alpha^\vee$ pour $\alpha\in\Delta_0^\Q$. Écrivons $X$ comme une somme:
\[
X=X_1+X_{2}+X_{3}
\]
où les $X_i$ sont combinaison linéaire d'éléments à coefficients positifs ou
nuls de coracines $\alpha^\vee$ pour $\alpha\in\Sigma_i$ avec
\[
\Sigma_1=\Delta_0^Q-\Delta_0^{P'},\qquad \Sigma_2=\Delta_0^{P'}-\Delta_0^{\Qdo}
\Qquad{et}\Sigma_3=\Delta_0^{\Qdo}\ptf
\]
Parce que
\[
\tau_{P'}^Q(H-T)=1
\]
$H_{P'}^Q$ est dans le cône engendré par les $\varpi^\vee$ pour $\varpi\in\hDelta_{P'}^Q$, lequel est contenu
dans celui engendré par les $\alpha^\vee_{P'}$, pour $
\alpha\in \Delta_0^Q-\Delta_0^{P'}$. L'élément $X_{1,P'}$ appartient à ce dernier cône. On en déduit
\[
\lVert H_{P'}^Q\rVert +\lVert X_{1,P'}\rVert \ll \lVert H_{P'}^Q+X_{1,P'}\rVert =\lVert H'_{P'}\rVert
\ll \lVert H'\rVert ,
\]
puis, parce que $X_1\mapsto X_{1,P'}$ est injective sur le cône auquel appartient $X_1$,
\begin{equation}
 \lVert H_{P'}^Q\rVert +\lVert X_1\rVert \ll \lVert H'\rVert \ptf
\label{eq12.7a}
\end{equation}
On a
\[
H_{\Qdo}^{P'}+X_{1,\Qdo}^{P'}+X_{2,\Qdo}=H^{\prime P'}_{\Qdo}\ptf
\]
D'où
\[
\lVert X_{2}\rVert \ll \lVert X_{2,\Qdo}\rVert \ll \lVert H_{\Qdo}^{P'}\rVert +\lVert H'\rVert +\lVert
X_1\rVert \ptf
\]
La relation~\eqref{eq12.4a} entraîne
\[
\lVert H_{\Qdo}^{P'}\rVert \ll 1+\lVert H_{P'}^Q\rVert 
\]
(la constante implicite dépend de $T$). En utilisant~\eqref{eq12.7a}, la relation ci-dessus devient
\begin{equation}
 \lVert X_{2}\rVert \ll 1+\lVert H'\rVert \ptf\label{eq12.8a}
\end{equation}
Parce que $x$ appartient à $\Sieg^{L_0}$, $\HO(x)$ appartient, à une translation fixe près, au cône
engendré par les $(\varpi^\vee_{\alpha})^{\Qdo}$ pour $\alpha\in
\Delta_0^{\Qdo}$, \emph{a fortiori} à celui engendré par les $\alpha^\vee$. On en déduit
\[
\lVert \HO(x)\rVert \ll 1+\rVert \HO(x)+X_{3}\rVert \ptf
\]
D'autre part, on a
\[
\HO(x)+X_{3}=(H')^{\Qdo}-X_1^{\Qdo}-X_{2}^{\Qdo}\ptf
\]
D'où
\[
\lVert \HO(x)\rVert \ll 1+\lVert H'\rVert +\lVert X_1\rVert +\lVert X_{2}\rVert \ptf
\]
Grâce à~\eqref{eq12.7a} et~\eqref{eq12.8a}, on obtient encore
\[
\lVert \HO(x)\rVert \ll 1+\lVert H'\rVert \ptf
\]
Cette relation, jointe à~\eqref{eq12.7a} et au fait que $H'-\HO(y')$ est borné, entraîne~\eqref{eq12.6a}.
En utilisant~\eqref{eq12.3a} et~\eqref{eq12.6a}, et en se rappelant que le $r$ de la relation~\eqref{eq12.5a} est quelconque, cette dernière relation entraîne
\[
\delta_{\Qdo}(\ee^H)^{-1}\lvert I(nn'x\ee^H)\rvert
\ll \ee^{-r(\lVert H_{P'}^Q\rVert +\lVert \HO(x)\rVert)}
\]
pour tout $r$. L'intégrale de l'énoncé, limitée comme on l'a dit au domaine défini par
\[
\phi_{\Qdo}^{P'}(H-T)\tau_{P'}^Q(H-T)=1
\]
est alors bornée par l'intégrale de l'expression de
droite ci-dessus sur le domaine suivant: $x$ parcourt $\Sieg^{L_0}$, $H$ parcourt un sous-ensemble de $\ga
^G_{\Qdo}$ sur lequel
\[
\lVert H\rVert \ll 1+\lVert H_{P'}^Q\rVert
\]
et $n$ et $n'$ restent dans des compacts. Pour $r$ assez grand, cette intégrale est finie.
\end{proof}
\setcounter{equation}{0}

\section[Transformation de l'opérateur $\tronc^{T,Q}$]{\mathversion{bold}Transformation de l'opérateur $\tronc^{T,Q}$}

On veut calculer l'expression:
\[
\int_{\YQdo}\tsQR(\HQ(y)-T)\int_{\ima(\ga_S^\G)^*}
\tronc^{T,Q}E^Q(y,\treg_{S,\sigma,\mu}(f,\omega)\WPhi,\theto\mu)
\widebar{E^{\pQ}(y,\WPhi,\mu)}\dd \mu\dd y\ptf
\]
On décompose l'intégrale sur $\YQdo$ comme dans la preuve précédente. On peut commencer par intégrer sur
\[
\bigl(N_{\Q_0^L}(F)\backslash N_{\Q_0^L}(\adef)\bigr)\times \bigl(N_{\Q_0^{L'}}(F)\backslash N_{\Q_0^{L'}}(\adef)\bigr)\ptf
\]
Cette
intégrale étant à support compact, on peut la permuter
avec l'intégrale sur $\ima(\ga_S^\G)^*$. On obtient pour composée de ces deux intégrales l'expression
\begin{equation}
\int_{\ima(\ga_S^\G)^*}(\tronc^{T,Q}E^Q)_{\Qdo}\bigl(x\ee^Hk,\WPhi,\theto(\mu)\bigr)
\widebar{E^{\pQ}_{\Qdo}(x\ee^Hk,\WPsi,\mu)}\varphi(\mu)\dd \mu\label{eq12.1b}
\end{equation}
les indices $\Qdo$ signifiant que l'on prend les termes constants. Le lemme~\ref{phinul} implique que ceci est nul si
$\phi_{Q_0}^Q(H-T)\ne 1$. Dans la preuve précédente, on avait
découpé le domaine d'intégration en $H$ selon des paraboliques $P'$. On voit que maintenant, seul le domaine correspondant
à $P'=Q$ donne une contribution non nulle.

Remarquons
que les diverses relations que l'on a établies dans la preuve précédente s'appliquent aussi bien à l'intégrale ci-dessus. La
relation~\eqref{eq12.4a} implique que, pour les $H$ qui vérifiant
\[
\tsQR(H-T)=1
\]
l'intégrale est nulle hors d'un domaine $\lVert H-T_{\Qdo}\rVert \leq c$, où $c$ est indépendant de $T$.

\emph{On fixe un réel $\eta$ avec $0<\eta<1$ que l'on précisera dans la proposition~\ref{W1.7} et
sera à l'œuvre dans la section~\ref{retourdep} \textup(on le supposera alors assez voisin de $0$\textup).}

Pour $Z\in \ga_0^G$, on note $\kappa^Z$ la fonction caractéristique du sous-ensemble des $X\in \ga_0^G$
tels que $\lVert X\rVert \leq \lVert Z\rVert$.
Remarquons que, quitte à agrandir le $c'$ ci-dessus, les relations
\[
\dPO(T)\geq c'(c+1)\Qquad{et}\lVert H-T_{\Qdo}\rVert \leq c
\]
entraînent
\[
\lVert H^Q-T^Q_{Q_0}\rVert \leq \lVert \eta T\rVert 
\]
autrement dit $\kappa^{\eta T}(H^Q-T^Q_{\Qdo})=1$.
En utilisant le lemme~\ref{LemmeB}, on obtient que, pourvu que
$\dPO(T)\geq c'(c+1)$, l'expression~\eqref{eq12.1b} multipliée par $\tsQR(H-T)$, vaut
\begin{multline*} 
\kappa^{\eta T}(H^Q-T^Q_{\Qdo})\tsQR(H-T)\phi_{\Qdo}^Q(H-T)\\
{}\times\int_{\ima(\ga_S^\G)^*}\tronc^{T
[H^Q],\Qdo}E^Q_{\Qdo}\bigl(x\ee^Hk,\WPhi,\theto(\mu)\bigr)
\widebar{E^{\pQ}_{\Qdo}(x\ee^Hk,\WPsi,\mu)}\varphi(\mu)\dd \mu
\end{multline*}
si $\lVert {H-T_{\Qdo}}\rVert \leq c$, et $0$ sinon. Mais la preuve de la relation~\eqref{eq12.4a} de la proposition~\ref{tsqrpsi} s'applique aussi bien à l'expression
ci-dessus: cette expression est nulle si $\lVert H-T_{Q_0}\rVert >c$. Donc l'expression~\eqref{eq12.1b} multipliée par $\tsQR(H-T)$ est égale à l'expression ci-dessus pour tout~$H$.

Il est utile de préciser le nombre $c$, qui dépend de $\varphi$. Pour $X\in \cgh$, notons $\ee^{\langle X,\bullet\rangle}$ la fonction
\[
\mu\mapsto \ee^{\langle X,\mu\rangle}\qquad\text{sur $\ima(\ga_S^\G)^*$}\ptf
\]
Que se passe-t-il quand on remplace $\varphi$ par $\varphi \ee^{\langle X,\bullet\rangle}$?
En examinant les preuves, on voit que le nombre $c$ est
essentiellement borné par le sup des normes des éléments du
compact $\omega$ du lemme~\ref{essentiel}.
Ce dernier est lui-même essentiellement le support d'une transformée de Fourier partielle de $\varphi$. Quand on remplace $\varphi$ par son produit avec
$\ee^{\langle X,\bullet\rangle}$, le nouvel $\omega$ est essentiellement un translaté du $\omega$ initial par une projection de l'élément $X$. Le
sup des normes de ses éléments est donc essentiellement
borné par $1+\lVert X\rVert $. Il en est donc de même de la constante $c$. On a obtenu la proposition ci-dessous.

\begin{proposition}\label{1.6}
Il existe $c(\varphi)>0$ tel que\textup:
\begin{enumerate}[(i)]
\item pour $\dPO(T)\geq c(\varphi)$\textup, on a l'égalité entre
\[
\int_{\YQdo}\tsQR(\HQ(y)-T)\int_{\ima(\ga_S^\G)^*}
\tronc^{T,Q}E^Q(y,\WPhi,\theto\mu)\widebar{E^{\pQ}(y,\WPsi,\mu)}\varphi(\mu)\dd \mu\dd y
\]
et
\begin{multline*}
\int_{\ga_{\Qdo}^G}\int_{L_0(F)\backslash L_0(\adef)^1}\int_{\K}
\kappa^{\eta T}(H^Q-T^Q_{\Qdo})\tsQR(H-T)\phi_{\Qdo}^Q(H-T)\delta_{\Qdo}(\ee^H)^{-1}\\
{}\times\int_{\ima(\ga_S^\G)^*}\tronc^{T[H^Q],\Qdo}E^Q_{\Qdo}(x\ee^Hk,\WPhi,\theto\mu)
\widebar{E^{\pQ}_{\Qdo}(x\ee^Hk,\WPsi,\mu)}\varphi(\mu)\dd \mu\dd k\dd x\dd H;
\end{multline*}
\item l'intégrale intérieure en $\mu$ du membre de droite ci-dessus est nulle pour tous $x,k$ si
\[
\lVert H-T_{\Qdo}\rVert > c(\varphi)\;;
\]
\item pour $\varphi$ fixée\textup, on a une majoration $c(\varphi \ee^{\langle X,\bullet\rangle})\ll (1+\lVert X\rVert)$ pour tout $X\in \cgh$.
\end{enumerate}
\end{proposition}

\section{De nouvelles majorations}

\begin{proposition}\label{W1.7}
Pour $H\in \ga_{\Qdo}^G$\textup, considérons
\begin{multline*}
\smash[b]{\int_{L_0(F)\backslash L_0(\adef)^1}\int_\K\int_{\ima(\ga_S^\G)^*}}
\lvert \tronc^{T[H^Q],\Qdo}E^Q_{\Qdo}(x\ee^Hk,\WPhi,\theto\mu)\\
{}\times\widebar{E^{\pQ}_{\Qdo}(x\ee^Hk,\WPsi,\mu)}\varphi(\mu)\rvert \dd \mu\dd k\dd x\ptf
\end{multline*}
\begin{enumerate}[(i)]
\item On suppose $\phi_{\Qdo}^Q(H-T)=1$. L'expression ci-dessus est convergente.
\item Il existe $\eta_0$ avec $0<\eta_0<1$ tel que si $\eta$ vérifie $0<\eta<\eta_0$\textup, la propriété suivante soit vérifiée. Il
existe $c>0$ telle que l'expression ci-dessus soit majorée
par $c\delta_{\Qdo}(\ee^H) \dPO(T)^{\dim(\ga_0^{\Qdo})}$ pour tout $T$ et tout $H$ vérifiant
\[
\tsQR(H-T)\phi_{\Qdo}^Q(H-T)\kappa^{\eta T}(H^Q-T^Q_{\Qdo})=1\ptf
\]
\end{enumerate}
\end{proposition}

\begin{proof} 
L'intégrale sur $K$ est inessentielle, on l'oublie. On veut majorer l'intégrale intérieure. Comme dans la preuve de la
proposition~\ref{eqqvf}, on se ramène à majorer deux types d'intégrales:
\begin{equation}
\int_{\ima(\ga_S^\G)^*}
\lvert \tronc^{T[H^Q],\Qdo}E^Q_{\Qdo}(x\ee^H,\WPhi,\theto\mu)\varphi(\mu)\rvert ^2\dd \mu\label{eq12.1c}
\end{equation}
et
\begin{equation}
 \int_{\ima(\ga_S^\G)^*} \lvert E^{\pQ}_{\Qdo}(x\ee^H,\WPsi,\mu)\varphi(\mu)\rvert^2\dd \mu\ptf\label{eq12.2c}
\end{equation}
Considérons la seconde, que l'on peut écrire
\[
\int_{\ima(\ga_S^\G)^*} E^{\pQ}_{\Qdo}(x\ee^H,\WPsi,\mu)\widebar{E^{\pQ}_{\Qdo}(x\ee^H,\WPsi,\mu)}\lvert \varphi(\mu)
\rvert^2\dd \mu\ptf
\]
Sous cette forme, on voit qu'elle se déduit de l'intégrale
\[
\int_{\ima(\ga_S^\G)^*} E^{\pQ}(y,\WPsi,\mu)\widebar{E^{\pQ}(y',\WPsi,\mu)}\lvert \varphi(\mu)\rvert^2\dd \mu
\]
en prenant les termes constants en chacune des variables $y,y'$, puis en posant $y=y'=\ee^Hx$ (prendre des termes constants
consiste à intégrer sur des compacts, cette opération
commute à l'intégrale sur $\ima(\ga_S^\G)^*$). Fixons une fonction $h^{\pQ}$ sur $\pQ(F)\backslash G(\adef)^1$, à
valeurs positives et telle que
\[
h^{\pQ}(y)\ll \delta_{\PO}(y)^{1/2}\ll h^{\pQ}(y)
\]
pour tout $y\in \Sieg^{\pQ}$: par exemple la fonction
\[
h^{\pQ}(y)=\sum_{\gamma\in \pQ(F)}\delta_{\PO}(\gamma y)^{1/2}\mathbf{1}_{\Sieg^{\pQ}}(\gamma y)
\]
où $\mathbf{1}_{\Sieg^{\pQ}}$ est la fonction caractéristique de $\Sieg^{\pQ}$.
La proposition~\ref{eqqvf} nous dit que la dernière intégrale ci-dessus est essentiellement bornée par
$h^{\pQ}(y)h^{\pQ}(y')$.
Donc \eqref{eq12.2c} est essentiellement bornée par
\[
h^{\pQ}_{\Qdo}(x\ee^H)h^{\pQ}_{\Qdo}(x\ee^H)\ptf
\]
La fonction $h^{\pQ}$ est à
croissance modérée, donc $h^{\pQ}_{\Qdo}$ aussi. Pour $H$
fixé, l'expression~\eqref{eq12.2c} est donc essentiellement bornée par $\lvert x\rvert^D$ pour un entier $D$ assez grand.
Considérons l'expression~\eqref{eq12.1c}. Elle se déduit de même de
\begin{equation}
 \int_{\ima(\ga_S^\G)^*} E^Q(y,\WPhi,\mu)\widebar{E^Q(y',\WPhi,\mu)}\lvert \varphi(\mu)\rvert^2\dd \mu
\label{eq12.3c}
\end{equation}
en prenant en chaque variable les termes constants puis en appliquant l'opérateur $\tronc^{T[H^Q],\Qdo}$, enfin en égalant
$y=y'=x\ee^H$. Quand on prend les termes constants, on
obtient comme ci-dessus une fonction essentiellement bornée par
\[
h^Q_{\Qdo}(y)h^Q_{\Qdo}(y')
\]
où $h^Q$ est l'analogue
de $h^{\pQ}$. Mais une majoration analogue vaut pour les dérivées en $y$ et $y'$
de notre fonction: en effet, par les procédés que l'on a déjà employés, de telles dérivées se
majorent par des combinaisons linéaires d'intégrales
similaires. Donc notre fonction est à croissance uniformément modérée en les deux variables $y$ et $y'$. Quand on applique
ensuite les opérateurs $\tronc^{T[H^Q],\Qdo}$, on
obtient une fonction à décroissance rapide en les deux variables grâce à la proposition~\ref{rapdec} (l'hypothèse $\phi_{\Qdo}^Q(H-T)=1$
assure que $T[H^Q]$ est régulier). Donc, pour tout $r$,
l'expression~\eqref{eq12.1c} est essentiellement bornée par $\lvert x\rvert^{-r}$ pour $x\in \Sieg^{L_0}$. Il en résulte que
l'intégrale intérieure de l'expression de la
proposition est à décroissance rapide en $x$. La première assertion de la proposition s'ensuit.
Pour la seconde assertion, on a besoin d'un ingrédient supplémentaire.
Montrons qu'il existe $D$ tel que l'on ait une majoration
\begin{equation}
h^Q_{\Qdo}(x\ee^H)\ll \delta_{\Qdo}(\ee^H)^{1/2}\lvert x\rvert^D\label{eq12.4c}
\end{equation}
pour
tout $H$ tel que $\tau^Q_{\Qdo}(H)=1$ et pour tout $x\in\Sieg^{L_0}$.

On ne perd rien à supposer le temps de la preuve que $Q=G$. On peut aussi supposer que $x=\ee^{\HO(x)}$. Posons
$Y=H+\HO(x)$. Considérons l'ensemble des
paraboliques standard $P'$ tels que $\Qdo\subset P'$ et $\alpha(Y)>0$ pour toute racine $\alpha\in \Sigma(N_{P'})$, où on
désigne ainsi l'ensemble des racines intervenant dans le
radical unipotent $N_{P'}$ de $P'$. Remarquons que pour deux paraboliques standard $P'_1$ et $P'_{2}$, et en notant
$P'_{3}=P'_1\cap P'_{2}$, on a
\[
\Sigma(N_{P'_{3}})=\Sigma(N_{P'_1})\cup \Sigma(N_{P'_{2}})\ptf
\]
Notre ensemble de paraboliques est donc stable par intersection, il possède en conséquence
un plus petit élément que l'on note $P'_0$. Si $P'_
0\ne \Qdo$, soit $\alpha\in \Delta^{P'_0}_0-\Delta_0^{\Qdo}$. Le parabolique $P'_{\alpha}$ tel que
$\Delta_0^{P'_{\alpha}}=\Delta_0-\{\alpha\}$ contient $\Qdo$ mais pas
$P'_0$. Il existe donc $\beta\in \Sigma(N_{P'_{\alpha}})$ tel que $\beta(Y)\leq 0$.
On fixe un tel $\beta$ que l'on décompose dans la
base $\Delta_0$. Le coefficient de $\alpha$ est
strictement positif. Puisque $\tau_{\Qdo}(H)=1$, on a donc $0<\alpha(H)\leq \beta(H)$.
D'où $0<\alpha(H)\leq -\beta\bigl(\HO(x)\bigr)$. Cela étant vrai pour tout
$\alpha\in \Delta^{P'_0}_0-\Delta_0^{\Qdo}$, on obtient une majoration
\[
\lVert H^{P'_0}\rVert \ll \lVert \HO(x)\rVert\ptf
\]
\emph{A fortiori}
\[
\lVert Y^{P'_0}\rVert \ll \lVert \HO(x)\rVert\ptf
\]
On a supposé $P'_0\ne \Qdo$ mais cette majoration reste vraie si $P'_0=\Qdo$, auquel cas
$Y^{P'_0}=\HH_0(x)$. Fixons $v\in \weyl^{P'_0}$ tel que $\alpha(vY)\geq0$ pour tout $\alpha\in \Delta_0^{P'_0}$.
Pour $\alpha\in \Delta_0-\Delta_0^{P'_0}$, on a
$\alpha(vY)=(v^{-1}\alpha)(Y)$. Puisque $v\in \weyl^{P'_0}$, $v^{-1}\alpha$ appartient à
$\Sigma(N_{P'_0})$, donc $(v^{-1}\alpha)(Y)>0$. On a donc $\alpha(vY)\geq0$ pour tout $\alpha\in \Delta_0$.

On a
\[
h^G_{\Qdo}(\ee^{H+\HO(x)})=\int_{N_{\Qdo}(F)\backslash N_{\Qdo}(\adef)}h^G(n\ee^{Y})\dd n\ptf
\]
Pour tout $n$ dans un ensemble de représentants de $N_{\Qdo}(F)\backslash N_{\Qdo}(\adef)$, fixons
$\gamma\in G(F)$ tel que $\gamma n\ee^{Y}\in \Sieg^G$. Calculons
$\HO(\gamma n\ee^Y)$. En utilisant la décomposition de Bruhat-Tits, on peut supposer que $\gamma=\nu'w\nu$ où $w$
normalise le Levi minimal et $\nu',\nu\in N_{\PO}(F)$. Alors
\[
\HO(\gamma n\ee^{Y})=\HO(\ee^{wY}wn')=wY+\HO(wn')
\]
où $n'=\ee^{-Y}\nu n \ee^{Y}$. D'après le lemme~\ref{nilneg},
on a une majoration $\varpi\bigl(\HO(wn')\bigr)\leq c$ pour
tout $\varpi\in \hDelta_0$, où $c$ est une certaine constante. On a
\[
wY=vY+(wv^{-1}-1)(vY)\ptf
\]
Puisque $vY$ est dans la
chambre positive fermée, on a $\varpi\bigl((wv^{-1}-1)(vY)\bigr)\leq0$
pour tout $\varpi\in \hDelta_0$. Donc $\varpi(\HO(\gamma n\ee^Y)-vY)\leq c$ pour tout $\varpi$. Il en résulte que
\[
h^G(n\ee^Y)\ll \delta_{\PO}(\gamma n\ee^Y)^{1/2}\ll \delta_{\PO} (\ee^{vY})^{1/2}\ptf
\]
On a
\[
\delta_{\PO}(\ee^{vY})=\delta_{P'_0}(\ee^{(vY)_{P'_0}})\delta_{\PO}^{P'_0}(\ee^{(vY)^{P'_0}})\ptf
\]
On a $(vY)^{P'_0}=v(Y^{P'_0})$ et, par construction, $(vY)_{P'_0}=H_{P'_0}$, donc
\[
\delta_{P'_0}(\ee^{(vY)_{P'_0}})=\delta_{P'_0}(\ee^{H_{P'_0}})
=\delta_{\Qdo}(\ee^{H})\delta_{\Qdo}^{P'_0}(\ee^{H^{P'_0}})^{-1}\ptf
\]
On obtient
\[
h^G(n\ee^Y)\ll \delta_{\Qdo}(\ee^{H})^{1/2}\delta_{\Qdo}^{P'_0}(\ee^{H^{P'_0}})^{-1/2}
\delta_{\PO}^{P'_0}(\ee^{v(Y^{P'_0})})\ptf
\]
On a montré que $\lVert H^{P'_0}\rVert$ et $\lVert Y^{P'_0}\rVert$ étaient essentiellement bornés par
$\lVert \HO(x)\rVert$. Il en résulte que le produit
des deux derniers termes ci-dessus est borné par $\lvert x\rvert ^D$ pour $D$ assez grand.
L'intégration en $u$ se faisant sur un compact, \eqref{eq12.4c} en résulte. Montrons que
pourvu que $\eta$ soit assez petit, l'hypothèse
\begin{equation}
\tsQR(H-T)\phi_{\Qdo}^Q(H-T)\kappa^{\eta T}(H^Q-T^Q_{\Qdo})=1\label{eq12.5c}
\end{equation}
implique
\[
\tau_{\Qdo}^P(H)=1\ptf
\]
En effet, pour
\[
\alpha\in \Delta_{\Qdo}^P-\Delta_{\Qdo}^Q
\]
l'hypothèse $\tsQR(H-T)$ implique
$\alpha(H_Q)>\alpha(T_Q)$. L'hypothèse $\phi_{\Qdo}^Q(H-T)=1$
implique que $H^Q-T^Q_{\Qdo}$ est combinaison linéaire à coefficients négatifs ou nuls de $\check{\beta}$ pour
$\beta\in \Delta_{\Qdo}^Q$. On a $\alpha(\check{\beta})\leq0$,
donc $\alpha(H^Q)\geq \alpha(T^Q)$ et finalement $\alpha(H)>\alpha(T)>0$. Pour $\alpha\in \Delta_{\Qdo}^Q$, il existe une
constante absolue $c>0$ telle que l'hypothèse
$\kappa^{\eta T}(H^Q-T^Q_{\Qdo})=1$ implique
\[
\lvert \alpha(H-T_{\Qdo})\rvert < c\eta \alpha(T)
\]
(rappelons que $T$ reste dans un cône fixé, \cf section~\ref{iter}). D'où
\[
\alpha(H)> \alpha(T_{\Qdo})-c\eta \alpha(T)\geq (1-c\eta)\alpha(T)\ptf
\]
Il suffit que $c\eta<1$ pour que cela entraîne $\alpha(H) \gg \alpha(T)>0$.
On suppose désormais
\[
\tsQR(H-T)\phi_{\Qdo}^Q(H-T)\kappa^{\eta T}(H^Q-T^Q_{\Qdo})=1\ptf
\]
On suppose
aussi $\eta$ tel que la conclusion de~\eqref{eq12.5c} soit vérifiée. Pour
simplifier, notons $\tronc$ l'opérateur $\tronc^{T[H^Q],\Qdo}$ et $C$ l'opérateur qui multiplie une fonction sur
\[
\YQdo\simeq \Qdo(F)\backslash G(\adef)^1\qquad\text{par la fonction $x\mapsto F_{\PO}^{\Qdo}(x,T[H^Q])$}\ptf
\]
L'intégrale intérieure de l'expression de l'énoncé se majore par la somme de deux intégrales
analogues où on remplace $\tronc$ soit par $\tronc-C$,
soit par $C$. Notons
\[
I_{\tronc-C}(x,H)\Qquad{et}I_{C}(x,H)
\]
ces deux intégrales. Commençons par majorer la première. De
nouveau, on doit majorer~\eqref{eq12.2c} et l'analogue, disons~(\ref{eq12.1c}$'$), de~\eqref{eq12.1c}
où $\tronc$ est remplacé par $\tronc-C$. Il résulte de~\eqref{eq12.4c} (appliqué à $\pQ$: l'hypothèse $\tau_{\Qdo}^{\pQ}(H)=1$ est
satisfaite) que~\eqref{eq12.2c} est essentiellement majoré par $\delta_{\Qdo}(\ee^H)\lvert x\rvert ^D$ pour un $D$ assez grand. On a une majoration analogue pour la fonction déduite de~\eqref{eq12.3c}
par passage aux termes constants. Comme on l'a expliqué,
on l'a même pour ses dérivées, avec des constantes implicites dépendant de la dérivation mais un $D$ uniforme. En
appliquant la proposition~\ref{propoC}, on voit que~(\ref{eq12.1c}$'$) est essentiellement majoré par
\[
\delta_{\Qdo}(\ee^H)\ee^{-r\mathbf{d}_{\PO\cap L_0}^{L_0}(T[H^Q])}\lvert x\rvert ^{-r}
\]
pour n'importe quel $r$.
On a déjà observé (juste avant le lemme~\ref{LemmeB})
que $T[H^Q]$ était \og plus régulier\fg que $T$, donc
\[
\dPO(T)\ll \mathbf{d}_{\PO\cap L_0}^{L_0}(T[H^Q])\ptf
\]
Il en résulte une majoration
\[
I_{\tronc-C}(x,H)\ll \delta_{\Qdo}(\ee^H)\ee^{-r\dPO(T)}\lvert x\rvert^{-r}
\]
pour tout $r$, puis
\begin{equation}
\int_{L_0(F)\backslash L_0(\adef)^1}I_{\tronc-C}(x,H)\dd x \ll \delta_{\Qdo}(\ee^H)\ee^{-r\dPO (T)}\ptf\label{eq12.6c}
\end{equation}
Majorons maintenant $I_{C}(x,H)$. L'opérateur $\tronc$ n'intervient plus.
Le procédé de la preuve de la proposition~\ref{eqqvf}
nous conduit à majorer~\eqref{eq12.2c} et une intégrale analogue où
$Q$ remplace $\pQ$, mais sous les restrictions
\[
\tsQR(H-T)\phi_{\Qdo}^Q(H-T)\kappa^{\eta T}(H^Q-T^Q_{\Qdo})=1
\]
et $F_{\PO}^{\Qdo}(x,T[H^Q])=1$. On peut aussi
supposer $x\in \Sieg^{L_0}$. Montrons que:
\begin{itemize}
\item il existe $c\in \RM$ tel que
\begin{equation}
\alpha\bigl(H+\HO(x)\bigr)\geq c\qquad\text{pour tout $\alpha\in \Delta_0^{\Qdo}$}\;;\label{eq12.7c}
\end{equation}
\item si $\eta$ est assez petit, il existe $c'>0$ tel que
\begin{equation}
\alpha\bigl(H+\HO(x)\bigr)\geq c'\dPO(T) \qquad\text{pour tout $\alpha\in \Delta_0^P-\Delta_0^{\Qdo}$}\ptf\label{eq12.8c}
\end{equation}
\end{itemize}
Pour $\alpha\in \Delta_0^{\Qdo}$, c'est clair puisque $x\in \Sieg^{L_0}$. Soit $\alpha\in \Delta_0^P -\Delta_0^{Q_0}$. L'hypothèse $F_{\PO}^{\Qdo}(x,T[H^Q])=1$ entraîne que $\HO(x)-T[H^Q]$
est combinaison linéaire à coefficients négatifs ou nuls de $\check{\beta}$ pour
$\beta\in \Delta^{\Qdo}_0$. On a $\alpha(\check{\beta})\geq0$,
donc
\[
\alpha\bigl(\HO(x)\bigr)\geq \alpha(T[H^Q])
\]
et on est ramené à considérer $\alpha(H+T[H^Q])$. Ecrivons
\[
T_{\Qdo}^Q-H^Q=\sum_{\beta\in \Delta_0^Q-\Delta_0^{\Qdo}}x_{\beta}\check{\beta}_{\Qdo}
\]
avec des $x_{\beta}\geq0$ (c'est l'hypohèse $\phi_{\Qdo}^Q(H-T)=1$).
D'après le lemme~\ref{raff} on a
\[
T[H^Q]=T^{\Qdo}-\sum_{\beta\in \Delta_0^Q-\Delta_0^{\Qdo}}x_{\beta}\check{\beta}^{\Qdo}\ptf
\]
Il en résulte que
\[
H+T[H^Q]=H_Q+T^Q-\sum_{\beta\in \Delta_0^Q-\Delta_0^{\Qdo}}x_{\beta}\check{\beta}=
H_Q-T_Q+T-\sum_{\beta\in \Delta_0^Q-\Delta_0^{\Qdo}}x_{\beta}\check{\beta}\ptf
\]
Si $\alpha\in \Delta_0^P-\Delta_0^Q$, les $\alpha(\check{\beta})$ sont négatifs ou nuls et
$\alpha(H_Q-T_Q)>0$ par l'hypothèse $\tsQR(H-T)=1$. Donc
\[
\alpha(H+T[H^Q])\geq \alpha(T)\geq \dPO(T)\ptf
\]
Supposons enfin $\alpha\in \Delta_0^Q-\Delta_0^{\Qdo}$. Alors
\[
 \alpha(H+T[H^Q])=\alpha(T)-\sum_{\beta\in \Delta_0^Q-\Delta_0^{\Qdo}}x_{\beta}\alpha(\check{\beta})
\]
et il existe une constante absolue $c_1>0$ telle que
\[
 \alpha(H+T[H^Q])\geq \dPO(T)-c_1\sup_{\beta\in \Delta_0^Q-\Delta_0^{\Qdo}}x_{\beta} \ptf
\]
Il existe une constante absolue $c_{2}>0$ telle que la condition $\kappa^{\eta T}(H^Q-T_{\Qdo}^Q)=1$ implique
\[
\sup_{\beta\in \Delta_0^Q-\Delta_0^{\Qdo}}x_{\beta}\leq c_{2}\eta \mathbf{d}_{P_0}(T)\ptf
\]
Donc
\[
 \alpha(H+T[H^Q])\geq (1-c_1c_{2}\eta)\dPO(T)\ptf
\]
Si $c_1c_{2}\eta<1$, la conclusion de~\eqref{eq12.7c} est vérifiée.
En conséquence de~\eqref{eq12.8c}, les éléments
\[
H^Q+\HO(x)^Q\Quad{et}H^{\pQ}+\HO(x)^{\pQ}
\]
restent dans des domaines de
Siegel relatifs à $Q$ et $\pQ$ (éventuellement plus gros
que ceux que l'on a fixés, mais peu importe). On a alors une majoration
\[
h^{\pQ}(nx\ee^H)\ll \delta_{\PO}(\ee^{H+\HO(x)})^
{1/2}
\]
pour tout $u\in N_{\Qdo}(\adef)$ d'où
\[
h^{\pQ}_{\Qdo}(x\ee^H)\ll \delta_{\PO}(\ee^{H+\HO(x)})^{1/2}\ptf
\]
Donc~\eqref{eq12.2c} est essentiellement majoré par
$\delta_{\PO}(\ee^{H+\mathbf{H}_0(x)})$. Il en est de même de l'analogue de~\eqref{eq12.2c} relatif
au parabolique $Q$ et donc aussi de $I_{C}(x,H)$. Alors
\[
\int_{L_0(F)\backslash L_0(\adef)^1}I_{C}(x,H)\dd x\,\ll \delta_{\Qdo}(\ee^H)
\int_{L_0(F)\backslash L_0(\adef)^1}F_{\PO}^{\Qdo}(x,T[H^Q])\delta_{\PO}(x)\dd x\ptf
\]
On majore la dernière intégrale en se limitant à un domaine de Siegel et en écrivant $x=vak$, avec $v$ dans un compact de
$(\PO\cap L_0)(\adef)${,} $a\in \gA^{L_0}(t)$ {et} $k\in K\cap L_0(\adef)\ptf$
Comme on sait, la décomposition des mesures introduit un $\delta_{\PO}(a)^{-1}$.
L'intégrale est donc essentiellement bornée par la mesure
du sous-ensemble des $a\in \gA^{L_0}(t)$ tels que $F_{\PO}^{\Qdo}(a,T[H^Q])=1$.
En écrivant $a=\ee^Y$, l'élément $Y$ reste dans l'intérieur d'un polyèdre de $\ga_0^{L_0}$
dont les côtés ont une longueur essentiellement bornée par $\lVert T[H^Q]\rVert $, ou encore par
\[
\lVert T\rVert +\lVert H^Q\rVert \ptf
\]
La condition $\kappa^{\eta T}(H^Q-T^Q_{\Qdo})=1$ entraîne une majoration
\[
\lVert H^Q\rVert \ll \lVert T\rVert\ptf
\]
Donc le volume du polyèdre est borné par
$\lVert T\rVert ^{\dim(\ga_0^{\Qdo})}$. On obtient
\[
\int_{L_0(F)\backslash L_0(\adef)^1}I_{C}(x,H)\dd x \ll 
\delta_{\Qdo}(\ee^H)\lVert T \rVert^{\dim(\ga_0^{\Qdo})}\ptf
\]
Jointe à~\eqref{eq12.6c}, cette majoration entraîne celle de l'énoncé.
\end{proof}

\section{Retour à la formule de départ}\label{retourdep}

\begin{quotation}\itshape
\hspace*{\parindent}Dorénavant\textup, on suppose $0<\eta<\eta_0$\textup, où $\eta_0$ vérifie les conditions de la proposition~\ref{W1.7}.
\end{quotation}

La proposition~\ref{tsqrpsi} entraîne la convergence dans l'ordre indiqué des doubles intégrales figurant dans l'expression
$J_{\chi}^T(f)$ de la section~\ref{formdep}. La proposition~\ref{1.6} entraîne l'égalité
\begin{multline*} 
J^\T_{\chi}(\fff,\omega)\\
\shoveleft{=\smash[t]{\sum_{\{Q,R\mid \PO\subset Q\subset R\}}
\tvedQR\sum_{\{S\mid \PO\subset S\subset\pQ\}}
\frac{1}{n^{\pQ}(S)}}}\\
{}\times\sum_{\sigma\in \Pi_{\mathrm{disc}}(M_S)}
\sum_{\WPsi\in \mathcal{B}^{\pQ}(\sigma)_{\chi}}
\int_{\ga_{\Q_0}^G}\int_{L_0(F)\backslash L_0(\adef)^{1}}\int_\K\kappa^{\eta T}(H^Q-T^Q_{\Qdo})
\tsQR(H-T)\\
{}\times\phi_{\Qdo}^Q(H-T)\delta_{\Qdo}(\ee^H)^{-1}
\smash[b]{\int_{\ima(\ga_S^\G)^*}}\tronc^{T[H^Q],\Qdo}E^Q_{\Qdo}(x\ee^Hk,\treg_{S,\sigma,\mu}(f,\omega)\WPsi,
\theto\mu)\\
{}\times\widebar{E^{\pQ}_{\Qdo}(x\ee^Hk,\WPsi,\mu)}\dd \mu\dd k\dd x
\dd H\ptf
\end{multline*}
Cela est vrai sous l'hypothèse $\dPO(T)\geq c'(c(f)+1)$ où $c'>0$ est une constante absolue et $c(f)>0$ dépend de
$f$. Quand on remplace $f$ par $f_{X}$, pour $X\in
\cgh$, les fonctions $\varphi$ qui interviennent dans le calcul sont changées en des combinaisons linéaires de fonctions
$\varphi \ee^{\langle X,\bullet\rangle}$, avec $s\in W_{\CM}$. D'après la proposition~\ref{1.6}, on a donc une majoration
\[
c(f_{X})\ll 1+\lVert X\rVert\ptf
\]
Précisément, $J^T_{\chi}(f_{X})$ est l'expression déduite de celle ci-dessus en glissant dans la dernière intégrale le terme
\[
\lvert W_{\CM}\rvert ^{-1}\sum_{s\in W_{\CM}}\ee^{\langle sX,\theto(\lambda(\sigma)+\mu)\rangle}\ptf
\]
On écrit $\lambda(\sigma)=X(\sigma)+iY(\sigma)$. En suivant Arthur, on considère l'ensemble des quadruplets $(Q,R,
\sigma,s)$ qui interviennent dans l'expression de $J^T_{\chi}(f_{X})
$. Il y a une application
\[
(Q,R,\sigma,s)\mapsto s^{-1}\theto\bigl(X(\sigma)\bigr)
\]
définie sur cet ensemble, à valeurs dans $\cgh$. On note $\boldsymbol{\Gamma}$ l'ensemble des fibres de cette
application. Pour une fibre $\Gamma$, on note $X_{\Gamma}$
son image. On peut alors écrire
\[
J_{\chi}^T\fff_X,\omega)=\sum_{\Gamma\in \boldsymbol{\Gamma}}\ee^{\langle X_{\Gamma},X\rangle}\psi^T_{\Gamma}(X,\fff,\omega),
\]
où $\psi^T_{\Gamma}(X,\fff,\omega)$ est la sous-somme de $J_{\chi}^T\fff_X,\omega)$
limitée aux $(Q,R,\sigma,s)\in \Gamma$ et multipliée par
$\ee^{-\langle X_{\Gamma},X\rangle}$ (autrement dit, on a sorti de $
\psi^T_{\Gamma}(X,\fff,\omega)$ la partie non unitaire des exponentielles).

\section{De nouveaux polynômes}

\begin{lemme}\label{W1.9.1}
Soit $\Gamma\in \boldsymbol{\Gamma}$. Pour tout opérateur différentiel $D$ sur $\cgh$ à
coefficients constants, on a une majoration
\[
\lvert D\psi_{\Gamma}^T(X,\fff,\omega)\rvert \ll (1+\lVert X\rVert )^{\dim(\ga^G_0)}\dPO(T)^{\dim(\ga_0^G)}
\]
pour tout $T$ et tout $X\in \cgh$.
\end{lemme}

\begin{proof} 
Le terme $\psi_{\Gamma}^T(X,\fff,\omega)$ est somme finie de termes
\begin{multline*}
\int_{\ga_{\Qdo}^G}\int_{L_0(F)\backslash L_0(\adef)^{1}}\int_\K\kappa^{\eta T}(H^Q-T^Q_{\Qdo})
\tsQR(H-T)\phi_{\Qdo}^Q(H-T)\delta_{\Qdo}(\ee^H)^{-1}\\
{}\times\smash[b]{\int_{\ima(\ga_S^\G)^*}}\tronc^{T[H^Q],\Qdo}E^Q_{\Qdo}(x\ee^Hk,\WPhi,\theto\mu)\widebar{E^{\pQ}_{\Qdo}(x\ee^Hk,\WPsi,
\mu)}\varphi(\mu)\\
{}\times\ee^{\langle sX,\ima Y(\sigma)+\mu\rangle} \dd \mu\dd k\dd x
\dd H\ptf
\end{multline*}
Appliquer l'opérateur $D$ ne fait que remplacer $\varphi$ par une autre fonction de Paley-Wiener. La proposition~\ref{W1.7} nous
permet de majorer essentiellement la triple intégrale intérieure
par
\[
\delta_{\Qdo}(\ee^H)\lVert T\rVert ^{\dim(\ga_0^{\Qdo})}\ptf
\]
D'autre part, d'après la proposition~\ref{1.6}, cette
triple intégrale est nulle sauf si $H$ vérifie une majoration
\[
\lVert H-T_{\Qdo}\rVert \ll 1+\lVert X\rVert \ptf
\]
Donc $D\psi_{\Gamma}^T(X,\fff,\omega)$ est essentiellement borné par le
produit de $\lVert T\rVert ^{\dim(\ga_0^{Q_0})}$ et de la mesure du sous-ensemble des $H\in \ga^G_{\Qdo}$ vérifiant cette majoration, laquelle est
essentiellement bornée par
\[
(1+\lVert X\rVert )^{\dim(\ga^G_{\Qdo})}\ptf
\]
Le lemme en résulte.
\end{proof}

\begin{proposition}\label{W1.9.2}
Pour tout $\Gamma\in \boldsymbol{\Gamma}$\textup, il existe une unique fonction $p_{\Gamma}^T(X,\fff,\omega)$ qui est
lisse en $X$ et polynomiale en $T$ de degré au plus $\dim(\ga_0^G)$ et qui vérifie les conditions suivantes\textup:
\begin{enumerate}[(i)]
\item il existe $c>$ tel que
\[
J_{\chi}^T\fff_X,\omega)=\sum_{\Gamma\in \boldsymbol{\Gamma}}\ee^{\langle X_{\Gamma},X\rangle}p_{\Gamma}^T(X,\fff,\omega)
\]
 si
$\dPO(T)\geq c(1+\lVert X\rVert )$\textup;
\item pour tout opérateur différentiel $D$ à coefficients constants sur $\cgh$\textup, il existe $R>0$ et $c_1,c_{2}>0$ tel que
\[
\bigl\lvert D\bigl(\psi_{\Gamma}^T(X,\fff,\omega)-p_{\Gamma}^T(X,\fff,\omega)\bigr)\bigr\rvert \leq c_1\ee^{-R\dPO(T)}
\]
si $\dPO(T)\geq c_{2}(1+\lVert X\rVert )$\textup;
\item pour tout opérateur différentiel $D$ à coefficients constants sur $\cgh$\textup, il existe $R\in \NM$ et $c'>0$ tels
que
\[
\lvert Dp_{\Gamma}^T(X,\fff,\omega)\rvert \leq c'(1+\lVert X\rVert)^R\dPO(T)^{\dim(\ga_0^G)}
\]
pour tous $X,T$.
\end{enumerate}
\end{proposition}

\begin{proof} 
Grâce au lemme~\ref{lemmeA} et à la proposition~\ref{1.6}, on est presque dans la situation de la proposition~5.1 de \cite{AeisI}. La seule différence est que Arthur dispose d'une majoration
\[
\lvert D\psi_{\Gamma}^T(X,\fff,\omega)\rvert \ll \dPO(T)^{\dim(\ga_0^G)}
\]
alors que le lemme ci-dessus nous fournit seulement
\[
\lvert D\psi_{\Gamma}^T(X,\fff,\omega)\rvert \ll (1+\lVert X\rVert )^{\dim(\ga_0^G)}
\dPO(T)^{\dim(\ga_0^G)}\ptf
\]
Un examen de la preuve d'Arthur montre que celle-ci s'applique encore,
avec bien sûr une conclusion plus faible concernant l'assertion~(iii). 
\end{proof}

\section{Permutation de deux intégrales}\label{sec12.8}

Soient $Q,R,S,\sigma,\WPhi,\WPsi$ comme dans la section~\ref{estimations}. Pour une fonction $\varphi$ sur $\ima(\ga_S^\G)^*$, posons au moins
formellement
\[
A^T(\varphi)=\int_{\ga_{\Qdo}^G}A^T(\varphi,H)\dd H
\]
où
{\multlinegap0pt
\begin{multline*}
A^T(\varphi,H)= \int_{L_0(F)\backslash L_0(\adef)^{1}}
\int_\K\kappa^{\eta T}(H^Q-T^Q_{\Qdo})\tsQR(H-T)\phi_{\Qdo}^Q(H-T)\delta_{\Qdo}(\ee^H)^{-1}\\
{}\times\int_{\ima(\ga_S^\G)^*} \tronc^{T[H^Q],\Qdo}E^Q_{\Qdo}(x\ee^Hk,\WPhi,\theto\mu)
\widebar{E^{\pQ}_{\Qdo}(x\ee^Hk,\WPsi,\mu)}\varphi(\mu) \dd \mu\dd k\dd x\ptf
\end{multline*}}%
Fixons $\varphi$ de Paley-Wiener sur $\ima(\ga_S^\G)^*$. Considérons une fonction
$B\in C_{c}^{\infty}(\ima\cgh^{G,*})$ et sa transformée de Fourier inverse $\widehat{B}$ sur
$\cgh^G$. On peut restreindre $B$ en une fonction sur $\ima(\ga_S^\G)^*$.

\begin{lemme}\label{W1.10}
Les expressions
\[
A^T(\varphi \ee^{\langle X,\bullet\rangle},H)\quad\text{pour $X\in \cgh$}, \Qquad{et} A^T(\varphi B,H)
\]
sont absolument
convergentes. Les expressions
\[
A^T(\varphi \ee^{\langle X,\bullet\rangle})\Qquad{et} A^T(\varphi B)
\]
sont absolument convergentes. L'intégrale
\[
\int_{\cgh^G}A^T(\varphi \ee^{\langle X,\bullet\rangle})\widehat{B}(X)\dd X
\]
est absolument convergente et est égale à $A^T(\varphi B)$.
\end{lemme}

\begin{proof}
Résumons ce que l'on a déjà prouvé. La proposition~\ref{W1.7} (que l'on peut aussi bien appliquer à
$\varphi \ee^{\langle X,\bullet\rangle}$ qui a même valeur absolue que $\varphi$) montre que
\begin{enumerate}
\item $A^T(\varphi \ee^{\langle X,\bullet\rangle},H)$ est absolument convergente et $\lvert A^T(\varphi \ee^{\langle X,\bullet\rangle},H)\rvert $ est borné
indépendamment de $X$ et $H$, $T$ étant fixé; plus précisément, l'intégrale
obtenue en remplaçant dans $A^T(\varphi \ee^{\langle X,\bullet\rangle},H)$ toutes les fonctions par leurs valeurs absolues est bornée.
\end{enumerate}
\noindent Par ailleurs,
d'après la proposition~\ref{1.6}:
\begin{enumerate} \addtocounter{enumi}{1}
\item il existe $c>0$ tel que
\[
A^T(\varphi \ee^{\langle X,\bullet\rangle},H)=0\qquad\text{sauf si $\lVert H-T_{\Qdo}\rVert \leq c(1+\lVert X\rVert )$}\ptf
\]
\end{enumerate}
L'intégrale $A^T(\varphi B,H)$ est aussi absolument convergente puisque $\varphi B$ est essentiellement bornée par $\varphi$.
Utilisons la formule d'inversion de Fourier:
\[
B(\mu)=\int_{\cgh^G}\widehat{B}(X)\ee^{\langle X,\mu\rangle}\dd X\ptf
\]
On obtient formellement
\[
A^T(\varphi B,H)=\int_{\cgh^G}\widehat{B}(X)A^T(\varphi \ee^{\langle X,\bullet\rangle},H)\dd X\ptf
\]
Ce calcul est justifié par~(1): quand on remplace toutes les fonctions par leurs valeurs absolues, l'expression ci-dessus reste
essentiellement bornée par
\[
\int_{\cgh^G}\bigl\lvert \widehat{B}(X) \bigr\rvert \dd X
\]
qui est convergente. En utilisant~(2), on obtient
\[
A^T(\varphi B,H)=\int_{X\in \cgh^G; \lVert H-T_{\Qdo}\rVert \leq c(1+\lVert X\rVert)}
\widehat{B}(X)A^T(\varphi \ee^{\langle X,\bullet\rangle},H)\dd X\ptf
\]
Mais l'intégrale
\[
\int_{\ga_{\Qdo}^G}\int_{X\in \cgh^G; \lVert H-T_{\Qdo}\rVert \leq c(1+\lVert X\rVert)}\lvert
\widehat{B}(X)A^T(\varphi \ee^{\langle X,\bullet\rangle},H)\rvert \dd X\dd H
\]
est convergente. En effet, on peut oublier le terme $A^T(\varphi \ee^{\langle X,\bullet\rangle},H)$ d'après~(1). L'intégrale en $H$ est
essentiellement bornée par
\[
(1+\lVert X\rVert )^{\dim(\ga^G_{\Qdo})}
\]
et l'intégrale restante en $X$ est convergente puisque $\widehat{B}$ est de Schwartz. Cela prouve la convergence
de $A^T(\varphi B)$. Cela prouve aussi que l'on peut
intervertir les intégrales:
\[
A^T(\varphi B)=\int_{\cgh^G}\widehat{B}(X)\int_{H\in \ga^G_{\Qdo}; \lVert H-T_{\Qdo}\rVert\leq
c(1+\lVert X\rVert)}A^T(\varphi \ee^{\langle X,\bullet\rangle},H)\dd H\dd X\ptf
\]
Toujours d'après~(2), on peut aussi bien supprimer la condition
\[
\lVert H-T_{\Qdo}\rVert\leq c(1+\lVert X\rVert)\ptf
\]
L'intégrale intérieure devient $A^T(\varphi \ee^{\langle X,\bullet\rangle})$ et
on obtient que $A^T(\varphi B)$ est donnée par l'intégrale de l'énoncé, laquelle est absolument convergente.
\end{proof}

\section[Un polynôme associé à la fonction $B$]{\mathversion{bold}Un polynôme associé à la fonction $B$}

Soit $B$ une fonction $C^{\infty}$ à support compact sur $\ima\cgh^*$, que l'on suppose invariante par $W_{\CM}$.
Pour des données $Q,R,S,\sigma$ intervenant dans
l'expression ci-dessous, on définit $B_{\sigma}$ sur $ \ima(\ga_S^\G)^*$ par $B_{\sigma}(\mu)=B(\ima Y(\sigma)+\mu)$. On
pose
{\multlinegap0pt
\begin{multline*} 
J_{\chi}^T(B,\fff,\omega)\\
\shoveleft{\quad=\sum_{\{Q,R\mid \PO\subset Q\subset R\}}\tvedQR
\sum_{\{S\mid \PO\subset S\subset\pQ\}}\frac{1}{n^{\pQ}(S)}}\\
{}\times\sum_{\sigma\in \Pi_{\mathrm{disc}}(M_S)}\sum_{\WPsi\in \mathcal{B}^{\pQ}(\sigma)_{\chi}}
\int_{\ga_{\Qdo}^G}\int_{L_0(F)\backslash L_0(\adef)^{1}}\int_\K\kappa^{\eta T}(H^Q-T^Q_{\Qdo})\tsQR(H-T)\\
\shoveright{{}\times\phi_{\Qdo}^Q(H-T)\delta_{\Qdo}(\ee^H)^{-1}
\int_{\ima(\ga_S^\G)^*}\tronc^{T[H^Q],\Qdo}E^Q_{\Qdo}(x\ee^Hk,\treg_{S,\sigma,\mu}(f,\omega)\WPsi,\theto\mu)}\\
{}\times\widebar{E^{\pQ}_{\Qdo}(x\ee^Hk,\WPsi,\mu)}B_{\sigma}(\mu)\dd \mu\dd k\dd x\dd H\ptf
\end{multline*}}%
Cette expression est combinaison linéaire finie d'expressions $A^T(\varphi B')$ du paragraphe précédent, où $\varphi$ est
de Paley-Wiener et $B'$ est $C^{\infty}$ et à support
compact. Donc les intégrales sont convergentes dans l'ordre indiqué. Pour $\epsilon>0$, définissons $B^{\epsilon}$ par
$B^{\epsilon}(\nu)=B(\epsilon\nu)$ pour tout $\nu\in\ima\cgh^
{G,*}$.

\begin{theoreme}\label{thWA}
\begin{enumerate}[(i)]
\item Pour tout $B$ comme ci-dessus\textup, il existe un unique polynôme $p^\T_{\chi}(B,\fff,\omega)$ en $T$ de degré au plus
$\dim(\ga_0^G)$ tel que
\[
\lim_{\dPO(T)\to \infty} \bigl(J_{\chi}^T(B,\fff,\omega)-p^\T_{\chi}(B,\fff,\omega)\bigr)=0\ptf
\]
\item Supposons $B(0)=1$. Alors il existe $c>0$ tel que
\[
J^\T_{\chi}(\fff,\omega)=\lim_{\epsilon\to 0} p^\T_{\chi}(B^{\epsilon},\fff,\omega)
\]
si $\dPO(T)\geq c$.
\end{enumerate}
\end{theoreme}

\begin{proof} 
C'est celle d'Arthur que nous ne reproduisons que pour nous rassurer. Pour $X\in \cgh^G$, on a écrit
\[
J^T_{\chi}(\fff_X,\omega)=\sum_{\Gamma\in \boldsymbol{\Gamma}}\ee^{\langle X_{\Gamma},X\rangle}\psi^T_{\Gamma}(X,\fff,\omega)\ptf
\]
Chaque $\psi^T_{\Gamma}(X,\fff,\omega)$ peut s'écrire comme une somme finie
\[
\sum_{(Q,R,\sigma,s)\in \Gamma}\sum_{\WPhi,\WPsi,\varphi}\ee^{\langle \ima Y(\sigma),sX\rangle}
A^T(Q,R,\sigma,s,\WPhi,\WPsi;\varphi \ee^{\langle sX,\bullet\rangle})
\]
où $(\WPhi,\WPsi,\varphi)$ décrit un ensemble fini indépendant de $X$ et où les
\[
A^T(Q,R,\sigma,s,\WPhi,\WPsi;\varphi \ee^{\langle sX,\bullet\rangle})
\]
sont les termes du paragraphe précédent (on a
simplement précisé leur notation). Par définition, $J^T_{\chi}(f)$ est l'expression ci-dessus pour $X=0$. Pour obtenir
\[
J^T_{\chi}(B,\fff,\omega)
\]
on doit glisser les fonctions $B_{\sigma}$ dans les
intégrales intérieures des termes
\[
A^T(Q,R,\sigma,s,\WPhi,\WPsi;\varphi \ee^{\langle sX,\bullet\rangle})\ptf
\]
Autrement dit
\[
J^T_{\chi}(B,\fff,\omega)=\sum_{\Gamma\in \boldsymbol{\Gamma}}\psi^T_{\Gamma}(B,\fff,\omega),
\]
où
\[
\psi^T_{\Gamma}(B,\fff,\omega)=\sum_{(Q,R,\sigma,s)}\sum_{\WPhi,\WPsi,\varphi}A^T(Q,R,\sigma,s,\WPhi,\WPsi;\varphi B_{\sigma})\ptf
\]
La fonction $B_{\sigma}$ est la restriction à $\ima(\ga_S^\G)^*$ de la translatée de $B$ par $\ima Y(\sigma)$. La
transformée de Fourier inverse de cette translatée est la fonction $X
\mapsto \ee^{\langle X,\ima Y(\sigma)\rangle}\widehat{B}(X)$. En appliquant le lemme~\ref{W1.10}, on obtient
\[
\psi^T_{\Gamma}(B,\fff,\omega)= \sum_{(Q,R,\sigma,s)}\int_{\cgh^G}\widehat{B}(X)\ee^{\langle \ima Y(\sigma),X\rangle}\sum_{\WPhi,\WPsi,\varphi}A^T(Q,R,\sigma,s,\WPhi,\WPsi;\varphi_{X})\dd X,
\]
cette expression étant absolument convergente.On remplace $X$ par $sX$ ce qui ne change pas $\widehat{B}(X)$ puisque $B$ est
invariante par $W_{\CM}$. On obtient
\[
\psi^T_{\Gamma}(B,\fff,\omega)=\int_{\cgh^G}\widehat{B}(X)\psi^T_{\Gamma}(X,\fff,\omega)\dd X\ptf
\]
On pose
\[
p^\T_{\chi}(B,\fff,\omega)=\sum_{\Gamma\in \boldsymbol{\Gamma}}\int_{\cgh^G}\widehat{B}(X)p^\T_{\Gamma}(X,\fff,\omega)\dd X\ptf
\]
Les intégrales sont convergentes puisque $p^\T_{\Gamma}(X,\fff,\omega)$ est à croissance modérée en $X$
(proposition~\ref{W1.9.2}(iii)) et $\widehat{B}$ est à décroissance rapide. C'est un polynôme
en $T$ de degré au plus $dim(\ga_0^G)$ puisqu'il en est de même de $p^\T_{\Gamma}(X,\fff,\omega)$. Alors
\[
J^T_{\chi}(B,\fff,\omega)-p^\T(B,\fff,\omega)=\sum_{\Gamma\in \boldsymbol{\Gamma}}\int_{\cgh^G}
\widehat{B}(X)(\psi^T_{\Gamma}(X,\fff,\omega)-p^\T_{\Gamma}(X,\fff,\omega)\dd X\ptf
\]
On découpe chaque intégrale en deux: l'une sur un domaine $(1+\lVert X\rVert )\leq c\dPO(T)$, l'autre sur le
complémentaire, où $c$ est une constante convenable.
Dans la première, on a, d'après la proposition~\ref{W1.9.2}(ii)
\[
\lvert \psi^T_{\Gamma}(X,\fff,\omega)-p^\T_{\Gamma}(X,\fff,\omega)\rvert \ll \ee^{-R \dPO(T)}
\]
pour un certain $R>0$ et l'intégrale vérifie la même majoration puisque $\widehat{B}$ est intégrable. Dans la deuxième, on a
\[
\lvert \psi^T_{\Gamma}(X,\fff,\omega)\rvert +\lvert p^\T_{\Gamma}(X,\fff,\omega)\rvert\ll (1+\lVert X\rVert )^D\dPO(T)^{\dim(\ga_0^G)}
\]
pour un certain entier $D$ d'après le lemme~\ref{W1.9.1} et la proposition~\ref{W1.9.2}(iii). L'intégrale est essentiellement bornée par
\[
\dPO(T)^{\dim(\ga_0^G)}\int_{X\in \cgh^G(c,T)}\lvert \widehat{B}(X)\rvert (1+\lVert X\rVert )^D\dd X
\]
où
\[
 \cgh^G(c,T)=\{X\in\cgh^\G\mid (1+\lVert X\rVert )\geq c\dPO(T)\}\ptf
\]
Puisque $\widehat{B}$ est de Schwartz, cette expression est essentiellement majorée par
\[
\dPO(T)^{-r}
\]
pour tout réel $r$. Cela prouve la première assertion de l'énoncé.
D'après la proposition~\ref{W1.9.2}(i), on a
\[
J^\T_{\chi}(\fff,\omega)=\sum_{\Gamma\in \boldsymbol{\Gamma}}p^\T_{\Gamma}(0,\fff,\omega)
\]
pourvu que $\dPO(T)$ soit assez grand. Par inversion de Fourier et d'après l'hypothèse $B(0)=1$,
l'intégrale de $\widehat{B^{\epsilon}}$ vaut $1$. Donc
\[
J^\T_{\chi}(\fff,\omega)=\sum_{\Gamma\in \boldsymbol{\Gamma}}\int_{\cgh^G}\widehat{B^{\epsilon}}(X)p^\T_{\Gamma}(0,\fff,\omega)\dd X\ptf
\]
On a aussi
\[
p^\T(B^{\epsilon},\fff,\omega)=\sum_{\Gamma\in \boldsymbol{\Gamma}}
\int_{\cgh^G}\widehat{B^{\epsilon}}(X)p^\T_{\Gamma}(X,\fff,\omega)\dd X\ptf
\]
Le calcul de
\[
\widehat{B^{\epsilon}}(X)=\epsilon^{-\dim(\cgh^G)}\widehat{B}(\epsilon^{-1}X)
\]
puis un changement de variables
entraîne l'égalité
\[
J^\T_{\chi}(\fff,\omega)-p^\T(B^{\epsilon},\fff,\omega)=\sum_{\Gamma\in \boldsymbol{\Gamma}}\int_{\cgh^G}
\widehat{B}(X)(p_{\Gamma}^T(0,\fff,\omega)-p^\T_{\Gamma}(\epsilon X,\fff,\omega)\dd X\ptf
\]
Le théorème des accroissements finis et la proposition~\ref{W1.9.2}(iii) entraînent une majoration
\[
\lvert p_{\Gamma}^T(0,\fff,\omega)-p^\T_{\Gamma}(\epsilon X,\fff,\omega)\rvert \ll \epsilon\lVert X\rVert
(1+\lVert X\rVert )^D\mathbf{d}_{\PO}(T)^{\dim(\ga_0^G)}\ptf
\]
Toujours parce que $\widehat{B}$ est de Schwartz, la différence ci-dessus est alors essentiellement bornée par $\epsilon$, d'où la
seconde assertion de l'énoncé.
\end{proof}
\chapter[Calcul de $A^T( B)$]{\mathversion{bold}Calcul de $A^T( B)$}\label{atb}

\section{Une majoration uniforme}

Soit $\phi$ une forme automorphe sur $G(F)\backslash G(\adef)$ (non nécessairement de carré intégrable). On sait définir
les termes constants cuspidaux de $\phi$: ce sont les
\og composantes cuspidales\fg $\phi_{P,\mathrm{cusp}}$ des termes constants $\phi_P $ de $\phi$ pour les différents paraboliques standard
$P$ de $G$, \cf \cite{MW}*{I.3.5}. Un tel terme $\phi_{P,\mathrm{cusp}}$
s'écrit sous la forme
\begin{equation}
\phi_{P,\mathrm{cusp}}(x)=\sum_{i=1}^{n_P }\ee^{\langle\rho_P +\lambda_{P,i},\HH_P (x)\rangle}
\sum_{j=1}^{n_{P,i}}p_{P,i,j}\bigl(\HH_P (x)\bigr)\phi_{P,i,j}(x),
\label{eq13.1}
\end{equation}
où les $\lambda_{P,i}$ sont des éléments de $\ga_{P,\CM}^*$, $\rho_P $ est la demi-somme usuelle des
racines intervenant dans $N_P $, les $p_{P,i,j}$ sont des
polynômes sur $\ga_P $ et les $\phi_{P,i,j}$ sont des formes automorphes cuspidales sur
\[
\gA_P P(F)\bs \Gadef\ptf
\]
Pour tout $P$, fixons un sous-
ensemble compact $\Gamma_P \subset \ga_{P,\CM}^*$, deux entiers naturels $\mathbf{n}_P $ et $d_P $ et un
sous-espace de dimension finie $V_P $ de l'espace des
formes automorphes cuspidales sur $P(F)\gA_P \backslash G(\adef)$. Notons
\[
A\bigl((V_P ,d_P ,\Gamma_P ,\mathbf{n}_P)_{P; \PO\subset P\subset G}\bigr)
\]
l'ensemble des
formes automorphes $\phi$ telles que, pour tout $P$, $\phi_{P,\mathrm{cusp}}$ puisse s'écrire sous la forme~\eqref{eq13.1}, avec
$n_P \leq \mathbf{n}_P $, des $\lambda_{P,i}\in \Gamma_P $, des $p_{P,i,j}$ de
degré inférieur ou égal à $d_P $ et des $\phi_{P,i,j}\in V_P $ (la condition $n_P \leq \mathbf{n}_P $ empêche cet ensemble
d'être stable par addition). On munit cet ensemble d'une
\og norme\fg que l'on note $\lVert \cdot\rVert_{\mathrm{cusp}}$ de la façon suivante. Pour tout $P$, notons
\[
\Pol_{\leq d_P }(\ga_P)
\]
l'espace des polynômes de degré${}\leq d_P $
sur $\ga_P $. Munissons
\[
\Pol_{\leq d_P }(\ga_P)\otimes V_P
\]
d'une norme $\lVert \cdot\rVert$
(il s'agit d'un espace de dimension finie, toutes les normes sont
équivalentes). Dans l'expression~\eqref{eq13.1}, on peut supposer les $\lambda_{i}$ tous distincts. L'élément
\[
\sum_{j=1}^{n_{P,i}}p_{P,i,j}\otimes \phi_{P,i,j}\in \Pol_{\leq d_P }(\ga_P)\otimes V_P
\]
est alors bien déterminé. On pose
\[
\lVert \phi_{P,\mathrm{cusp}}\rVert_{\mathrm{cusp}}=\sum_{i=1}^{n_P }\Biggl\lVert
\sum_{j=1}^{n_{P,i}}p_{P,i,j}\otimes \phi_{P,i,j}\Biggr\rVert
\]
puis,
\[
\lVert\phi\rVert_{\mathrm{cusp}}=\sum_P \lVert \phi_{P,\mathrm{cusp}}\rVert_{\mathrm{cusp}}
\]
pour
\[
\phi\in A\bigl((V_P ,d_P ,\Gamma_P ,\mathbf{n}_P)_{P; \PO\subset P\subset G}\bigr)\ptf
\]

\begin{lemme} \label{2.1.1}
\begin{enumerate}[(i)]
\item Il existe un réel $D$ et\textup, pour tout $X\in \mathfrak{U}(\frakg)$\textup, il existe $c>0$ tel que\textup, pour tout $\phi
\in A\bigl((V_P ,d_P ,\Gamma_P ,\mathbf{n}_P)_{P; \PO\subset
P\subset G}\bigr)$\textup, on ait la majoration
\[
\lvert X\phi(x)\rvert\leq c\lVert \phi\rVert _{\mathrm{cusp}}\lvert x\rvert^D
\]
pour tout $x\in G(\adef)$.
\item Pour tout $\lambda\in \ga_0$\textup, il existe $c>0$ tel que\textup, pour tout
\[
\phi\in A\bigl((V_P ,d_P ,\Gamma_P ,\mathbf{n}_P)_{P; \PO\subset P\subset G}\bigr)
\]
et tout $x\in \gA_{G}
\Sieg^G$, on ait la majoration
\[
\lvert \phi(x)\rvert\leq c\lVert \phi\rVert _{\mathrm{cusp}}\sum_P \sum_{i=1}^{n_P }
\ee^{\langle \rho_P +\lambda^P+\RE(\lambda_{P,i}),\HO(x)\rangle}\bigl(1+\HH_P (x)\bigr)^{d_P }\ptf
\]
\end{enumerate}
\end{lemme}

\addtocounter{equation}{1}
\begin{proof}
La conjonction de \cite{MW}*{lemmes~I.4.4(b) et~I.4.3} assure l'existence d'un réel $D$ et d'un $c_1>0$ tels que
\[
\lvert \phi(x)\rvert\leq c_1\lVert \phi\rVert_{\mathrm{cusp}}\lvert x\rvert^D
\]
pour tout $x$ et tout
\[
\phi\in A\bigl((V_P ,d_P ,\Gamma_P ,\mathbf{n}_P)_{P; \PO\subset P\subset G}\bigr)\ptf
\]
Le (a) du même lemme~I.4.4 assure l'existence de $h\in \ctyc\bigl(G(\adef)\bigr)$,
d'une constante $c_{2}>0$ et, pour tout
\[
\phi\in A\bigl((V_P ,d_P ,\Gamma_P ,\mathbf{n}_P)_{P; \PO\subset P\subset G}\bigr)
\]
d'une forme automorphe $\phi'$ telle que
\begin{gather}
\delta(h)\phi'=\phi\label{eq13.2}\\
\sup_{x\in G(\adef)} \lvert \phi'(x)\rvert\, \lvert x\rvert^{-D} \leq c_{2} \sup_{x\in G(\adef)}
\lvert \phi(x) \rvert \,\lvert x\rvert^{-D}\ptf\label{eq13.3}
\end{gather}
Pour $X\in \mathfrak{U}(\frakg)$, on a $X\phi=\delta(Xh)\phi'$. Puisque $h$ est à support compact, cela entraine
l'existence de $c_{3}>0$ dépendant de $X$ mais pas de $\phi$ tel
que
\[
\sup_{ x\in G(\adef)} \lvert X\phi(x)\rvert \,\lvert x\rvert^{-D}\leq c_{3}\sup_{ x\in G(\adef)}
\lvert \phi'(x)\rvert \,\lvert x\rvert^{-D}\ptf
\]
Alors
\[
\sup_{x\in G(\adef)}\lvert X\phi(x)\rvert \,\lvert x\rvert^{-D}\leq c_{2 }c_{3}
\sup_{x\in G(\adef)}\lvert \phi(x)\rvert \,\lvert x\rvert^{-D}\leq c_1c_{2}c_{3}\lVert\phi\rVert_{\mathrm{cusp}},
\]
ce qui prouve~(i).
Le lemme~I.4.1 de \cite{MW} énonce la majoration~(ii) de façon plus imprécise: le terme $c\lVert \phi\rVert_{\mathrm{cusp}}$ y est
remplacé par un réel $c(\phi)>0$ qui n'est pas précisé.
Il suffit de reprendre la démonstration pour voir que l'on peut prendre
\[
c(\phi)=c\lVert \phi\rVert_{\mathrm{cusp}}
\] avec un $c$
indépendant de $\phi$. Le seul point un peu subtil est de
montrer que dans la majoration des fonctions $s\psi_{Q,\xi}$ de la page~50 de cette référence, on peut remplacer la constante
non précisée par $c\lVert \phi\rVert_{\mathrm{cusp}}$ où
$c$ est indépendant de $\phi$. Mais cela résulte du~(i) que l'on vient de prouver et du corollaire~I.2.11 de \cite{MW}. \end{proof}

\setcounter{equation}{0}
\section{Majoration des termes constants}\label{W2.2}

On fixe, dans les sections suivantes jusqu'en~\ref{Ws2.8} inclusivement,
des paraboliques standard $Q\subset R$, avec $\tvedQR=1$, un
parabolique standard $S\subset \pQ$ et une représentation $
\sigma\in \Pi_{\mathrm{disc}}(M_S)$.

La représentation $\sigma$ intervient dans le spectre discret de $M_S(F)\backslash M_S(\adef)^1$ (où $M_S$ est le
Levi standard de $S$). D'après Langlands, les éléments
de celui-ci sont des résidus de séries d'Eisenstein issus d'une représentation cuspidale. Rappelons plus précisément les
propriétés dont nous avons besoin. Considérons:
\begin{itemize}
\item un parabolique standard $S_{\mathrm{cusp}}\subset S$;
\item une représentation automorphe cuspidale $\sigma_{\mathrm{cusp}}$ de $M_{S_{\mathrm{cusp}}}(\adef)$;
\item un opérateur différentiel $D$ à coefficients polynomiaux sur $\ga_{S_{\mathrm{cusp}},\CM}^{S,*}$;
\item un point $\nu_0\in \ga_{S_{\mathrm{cusp}}}^{S,*}$.
\end{itemize}
Pour $\WPhi_{\mathrm{cusp}}\in \mathcal{B}^{M_S}(\sigma_{\mathrm{cusp}})$ et $\nu\in \ga_{S_{\mathrm{cusp}},\CM}^{S,*}$, formons la série
d'Eisenstein
\[
E^{M_S}(y,\WPhi_{\mathrm{cusp}},\nu)
\]
(la variable $y$ appartient à $M_S(\adef)$). On applique l'opérateur $D$. On suppose que
\begin{equation}
\text{\emph{la fonction $DE^{M_S}(y,\WPhi_{\mathrm{cusp}},\nu)$ est holomorphe en $\nu=\nu_0$.}}\label{eq13.1a}
\end{equation}
On note
\[
D_{\nu=\nu_0}E^{M_S}(y,\WPhi_{\mathrm{cusp}},\nu)
\]
sa valeur en $\nu=\nu_0$. La première assertion est que l'espace de
$\sigma$ est engendré par de telles fonctions
\[
D_{\nu=\nu_0}E^{M_S}(y,\WPhi_{\mathrm{cusp}},\nu)
\]
pour des données $S_{\mathrm{cusp}}$, $\sigma_{\mathrm{cusp}}$, $D$, $\nu_0$ et $\WPhi_{\mathrm{cusp}}$ vérifiant les
conditions précédentes. Ces données vérifient en fait
des conditions supplémentaires. En tout cas, l'assertion précédente s'induit à $Q$ ou $\pQ$. En choisissant convenablement la
base $\base^{\pQ}(\sigma)$, on peut supposer que pour
tout élément $\WPsi$ de cette base, il existe des données $S_{\mathrm{cusp}}$, $\sigma_{\mathrm{cusp}}$, $D$, $\nu_0$ et
$\WPhi_{\mathrm{cusp}}\in \mathcal{B}^{\pQ}(\sigma_{\mathrm{cusp}})$ de sorte que
\[
E^{\pQ}(y,\WPsi,\mu)=D_{\nu=\nu_0}E^{\pQ}(y,\WPhi_{\mathrm{cusp}},\nu+\mu)
\]
pour tout $\mu\in (\ga_{S,\CM}^G)^*$. Prendre un terme constant et prendre un résidu sont des opérations qui
commutent. En appliquant l'équation~\eqref{eq5.4} du théorème~\ref{eqfonct}, on obtient
\begin{equation}
E^{\pQ}_{\Qdo}(y,\WPsi,\mu)=D_{\nu=\nu_0}\sum_{s\in \weyl^{\pQ}(\ga_{S_{\mathrm{cusp}}},\Qdo)}
E_{\Qdo}\bigr(y,\mathbf{M}(s,\nu+\mu)\WPhi_{\mathrm{cusp}},s(\nu+\mu)\bigr)\ptf\label{eq13.2a}
\end{equation}
Rappelons que, nos sous-groupes paraboliques étant standard, on désigne par
\[
\weyl^{\pQ}(\ga_{S_{\mathrm{cusp}}},\Qdo)
\]
l'ensemble des restrictions à $\ga_{S_{\mathrm{cusp}}}$ d'éléments
$s\in \weyl^{\pQ}$ tels que $s(\ga_{S_{\mathrm{cusp}}})\supset \ga_{\Qdo}$ et $s$ est de longueur minimale dans
sa classe $\weyl^{\Qdo} s$ (ce qui se traduit par $s(S_{\mathrm{cusp}})
\cap L_0$ est standard dans $L_0$) (\cf lemme~\ref{weylgab}).

Calculons les termes constants cuspidaux de cette forme automorphe, c'est-à-dire ses termes constants relatifs aux
sous-groupes paraboliques standard de $\Qdo$ associés à
$S_{\mathrm{cusp}}$ dans $\pQ$. Pour un tel sous-groupe parabolique $\Spcusp$, on a
{\multlinegap0pt
\begin{multline*}
E^{\pQ}_{\Spcusp}(y, \WPsi,\mu)\\
=D_{\nu=\nu_0}\mspace{-10mu}\sum_{s\in \weyl^{\pQ}(\ga_{S_{\mathrm{cusp}}},\Qdo)}
\sum_{s'\in \weyl^{\Qdo}
(\ga_{s(S_{\mathrm{cusp}})},\ga_{\Spcusp})} \mspace{-10mu}\Mint
(s's,\nu+\mu)\WPhi_{\mathrm{cusp}}\bigl(y,s's(\nu+\mu)\bigr)\ptf
\end{multline*}
Les exposants cuspidaux des différents termes sont les $s's(\nu_0+\mu)$. Notons $\mathcal{W}(\nu_0)$ le stabilisateur de
$\nu_0$ dans $\weyl^S(M_{S_{\mathrm{cusp}}})$. Regroupons les différents $s'$, $s$ selon la classe $s's\mathcal{W}(\nu_0)$.
On obtient que $E^{\pQ}_{\Spcusp}(y, \WPsi,\mu)$ est la valeur en $\nu=\nu_0$ de
\begin{equation}
\sum_{w\in \weyl^{\pQ}(\ga_{S_{\mathrm{cusp}}},\ga_{\Spcusp})/\mathcal{W}(\nu_0)}\mspace{-10mu}
D\Biggl(\sum_{\{s,s'\mid s's\in w\mathcal{W}(\nu_0)\}} \mspace{-10mu}\Mint(s's,\nu+\mu)\WPhi_{\mathrm{cusp}}\bigl(y,s's(\nu+\mu)\bigr)\Biggr)\ptf\label{eq13.3a}
\end{equation}
Remarquons que pour $s'_1,s_1,s'_{2},s_{2}$ tels que $s'_1s_1\notin s'_{2}s_{2}\mathcal{W}(\nu_0)$, on a
\[
s'_1s_1(\nu_0+\mu)\ne s'_{2}s_{2}(\nu_0+\mu)
\]
pour un point
$\mu\in \ima(\ga_S^\G)^*$ en position générale. Au moins en un tel point $\mu$,
l'holomorphie de l'expression ci-dessus en $\nu=\nu_0$ entraîne que chaque composante l'est
aussi. D'où
\begin{multline}\label{eq13.4a}
E^{\pQ}_{\Spcusp}(y, \WPsi,\mu)\\
=\mspace{-10mu}
\sum_{w\in \weyl^{\pQ}(\ga_{S_{\mathrm{cusp}}},\ga_{\Spcusp})/\mathcal{W}(\nu_0)}\mspace{-10mu}
D_{\nu=\nu_0}\Biggl(\sum_{s,s'; s's\in w\mathcal{W}(\nu_0)}\mspace{-10mu}\Mint(s's,\nu+\mu)\WPhi_{\mathrm{cusp}}\bigl(y,s's(\nu+\mu)\bigr)\Biggr)\ptf
\end{multline}}%

On lève l'hypothèse que $\mu$ est en position générale en considérant cette égalité comme une égalité de fonctions
méromorphes en $\mu$. Pour
\[
w\in \weyl^{\pQ}(\ga_{S_{\mathrm{cusp}}},\ga_{\Spcusp})
\]
notons $\pQ_{w}$ le plus petit sous-groupe parabolique standard de $\pQ$ tel que
$\ga_{\pQ_{w}}\subset w(\ga_S)$. On sait que:
\begin{enumerate}\setcounter{enumi}{4}
\item\label{eq13.5a}\itshape pour $\mu\in \ima(\ga_S^\G)^*$\textup, les parties réelles des exposants cuspidaux de
\[
E^{\pQ}_{\Spcusp}(y, \WPsi,\mu)
\]
sont de la forme $w\nu_0$ pour des $w\in \weyl^{\pQ}(\ga_{S_{\mathrm{cusp}}},\ga_{\Spcusp})$ tels que
\[
\htau_{\Spcusp}^{\pQ_{w}}(-w\nu_0)=1\ptf
\]
\end{enumerate}
(\Cf \cite{Eisen}*{lemme~7.5 et théorème~7.1}, repris dans \cite{MW}*{corollaire~V.3.16 et proposition~VI.1.6(c)}.
L'exposé le plus clair, mais sans démonstration, est sans doute \cite{Fr}*{\parag 5.2}).

Les termes de l'expression~\eqref{eq13.4a} indexés par des $w$ ne vérifiant pas la conclusion de~\eqref{eq13.5a} sont donc nuls.
Décomposons $E^{\pQ}_{\Qdo}(y,\WPsi,\mu)$ en
\[
E^{\pQ}_{\Qdo,\mathrm{unit}}(y,\WPsi,\mu)+E^{\pQ}_{\Qdo,+}(y,\WPsi,\mu)\ptf
\]
Le premier terme est la sous-somme de~\eqref{eq13.2a} indexée par les $s$ tels que
$s(\ga_0^{S})\subset \ga_0^{\Qdo}$, le second est la sous-somme restante. Par restriction à
$\ga_S$, le premier ensemble de sommation s'identifie à
$\weyl^{\pQ}(\ga_S,\Qdo)$. Pour un élément $s$ de cet ensemble, on dispose de l'opérateur
\[
\Mint(s,\mu)\colon L^2_{\mathrm{disc}}\bigl(S(F)\gA_{G}\backslash G(\adef)\bigr)\to L^2_{\mathrm{disc}}\bigl(S_{s}(F)\gA_{G}\backslash G(\adef)\bigr),
\]
où $S_{s}$ est le parabolique standard de Levi $s(M_S)$. On a (\cf \cite{MW}*{p.~268, égalité~(5)})\,\footnote{Waldspurger dixit:
\og Je cite \cite{MW} parce que je le connais mieux, mais le résultat est bien sûr dû à Langlands.\fg}:
\begin{enumerate}\setcounter{enumi}{5}
\item\label{eq13.6a}\itshape la fonction $DE_{\Qdo}\bigl(y,\Mint(s,\nu+\mu)\WPhi_{\mathrm{cusp}},s(\nu+\mu)\bigr)$ est holomorphe en
$\nu=\nu_0$ et sa valeur en $\nu_0$ est $E_{\Qdo}(y,\Mint(s,\mu)\WPhi,\mu)$.
\end{enumerate}
D'où\addtocounter{equation}{2}
\begin{equation}
E^{\pQ}_{\Qdo,\mathrm{unit}}(y,\WPsi,\mu)=\sum_{s\in \weyl^{\pQ}(\ga_S,\Qdo)} E_{\Qdo}(y,\Mint(s,\mu)\WPhi,\mu)\ptf
\label{eq13.7a}
\end{equation}
Cela entraîne que $ E^{\pQ}_{\Qdo,\mathrm{unit}}(y,\WPsi,\mu)$ est holomorphe en $\mu$. La différence
\[
E^{\pQ}_{\Qdo,+}(y,\WPsi,\mu)= E^{\pQ}_{\Qdo}(y,\WPsi,\mu)- E^{\pQ}_{\Qdo,\mathrm{unit}}(y,\WPsi,\mu)
\]
l'est donc aussi.

\begin{proposition} \label{2.2.1}
Soient $\Omega$ un sous-ensemble compact de $\ima(\ga_S^\G)^*$ et $T_1$ un élément de
$\ga_0^{\Qdo}$.
\begin{asparaenum}[(i)]
\item Pour tout $X\in \mathfrak{U}(\mathfrak{l}_0)$\textup, il existe un entier $N\in \NM$ et un réel $c>0$ tels que
\[
\lvert XE^{\pQ}_{\Qdo}(x\ee^Hk,\WPhi,\mu)\rvert \leq c\delta_{\PO}(\ee^H)^{1/2}(1+\lVert H\rVert)^N\lvert x\rvert ^N
\]
pour tous $k\in K$\textup, $\mu\in \Omega$\textup, $x\in \Sieg^{L_0}$\textup, $H\in \ga_{\Qdo}^G$ tel que
$\tau_{\Qdo}^{\pQ}(H)=1$.
\item Pour tout $X\in \mathfrak{U}(\mathfrak{l}_0)$\textup, il existe un entier $N\in \NM$ et un réel $c>0$ tels que
\[
\lvert XE^{\pQ}_{\Qdo}(x\ee^Hk,\WPhi,\mu)\rvert \leq
c\delta_{\PO}(x\ee^H)^{1/2}(1+\lVert H\rVert)^N(1+\lVert\HO(x)\rVert)^N
\]
pour tous $k \in  K$\textup, $\mu \in  \Omega$\textup, $x \in  \Sieg^{L_0}$\textup, $H \in  \ga_{\Qdo}^G$ tels que
${\tau_{\PO}^{\pQ}(H + \HO(x) + T_1) = 1}$.
\item Pour tout $X\in \mathfrak{U}(\mathfrak{l}_0)$\textup, il existe un entier $N\in \NM$ et un réel $c>0$ tels que
\[
\lvert XE^{\pQ}_{\Qdo,\mathrm{unit}}(x\ee^Hk,\WPhi,\mu)\rvert \leq c\delta_{\PO}(x\ee^H)^{1/2}(1+\lVert H\rVert)^N(1+\lVert\HO(x)\rVert)^N
\]
pour tous $k\in K$\textup, $\mu\in \Omega$\textup, $x\in \Sieg^{L_0}$ et $H\in \ga_{\Qdo}^G$.
\item Il existe un réel $R>0$ et\textup, pour tout $X\in \mathfrak{U}(\mathfrak{l}_0)$\textup, il existe un entier $N\in \NM$ et un réel $c>0$ tels que
\begin{multline*}
\lvert XE^{\pQ}_{\Qdo,+}(x\ee^Hk,\WPhi,\mu)\rvert\\ \leq c\delta_{\PO}(x\ee^H)^{1/2}
(1+\lVert H\rVert)^N(1+\lVert \HO(x)\rVert) ^N
\sup_{\alpha\in \Delta_0^{\pQ}-\Delta_0^{\Qdo}}\ee^{-R\alpha\left(H+\HO(x)\right)}
\end{multline*}
pour tous $k \in  K$\textup, $\mu \in  \Omega$\textup, $x \in  \Sieg^{L_0}$\textup, $H \in  \ga_{\Qdo}^G$ tels que
${\tau_{\PO}^{\pQ}(H + \HO(x) + T_1)  = 1}$.
\end{asparaenum}
\end{proposition}

\addtocounter{equation}{7}
\begin{proof}
La fonction
\[
XE^{\pQ}_{\Qdo}(y,\WPhi,\mu)
\]
est combinaison linéaire de fonctions $E^{\pQ}_{\Qdo}(y,\WPsi,\mu)$, avec
pour coefficients des fonctions $C^{\infty}$ de $\mu$. Puisque
$\mu$ reste dans un compact, des majorations pour ces dernières fonctions entraînent les mêmes majorations pour
$XE^{\pQ}_{\Qdo}(y,\WPhi,\mu)$. On peut donc se limiter au cas $X=1$.
Comme toujours, le $k$ ne compte guère, on l'oublie. Notons $\Xi$ l'image de l'application
\[
\begin{aligned}
\weyl^{\pQ}(\ga_{S_{\mathrm{cusp}}},\Qdo)&\to \ga_{\Qdo}\\ 
s&\mapsto (s\nu_0)_{\Qdo}\\
\end{aligned}\ptf
\]
Soit $\xi\in \Xi$. Notons $\weyl^{\pQ}(\ga_{S_{\mathrm{cusp}}},\Qdo)_{\xi}$ la fibre au-dessus de $\xi$ et considérons la fonction
\begin{equation}
D\Biggl(\sum_{s\in \weyl^{\pQ}(\ga_{S_{\mathrm{cusp}}},\Qdo)_{\xi}}
E_{\Qdo}\bigl(y,\Mint(s,\nu+\mu)\WPhi_{\mathrm{cusp}},s(\nu+\mu)\bigr)\Biggr)\ptf
\label{eq13.8a}
\end{equation}
Son terme constant cuspidal relatif à un parabolique $\Spcusp$ est la sous-somme de~\eqref{eq13.3a} où on ne garde que les
$s\in \weyl^{\pQ}(\ga_{S_{\mathrm{cusp}}},\Qdo)_{\xi}$. Mais pour $s$,
$s'$ et $w$ comme dans la relation~\eqref{eq13.3a}, cette condition sur $s$ se lit sur $w$: elle équivaut à $(w\nu_0)_{\Qdo}=\xi$.
Donc le terme constant cuspidal de la fonction ci-dessus est la sous-somme de~\eqref{eq13.3a} indexée par les $w$ vérifiant
$(w\nu_0)_{\Qdo}=\xi$. On a déjà dit qu'une telle somme était holomorphe
en $\nu=\nu_0$. On l'avait dit pour $\mu$ en position générale. Mais en reprenant l'argument,
on voit que c'est vrai pour tout $\mu$: si $w$ et $w'$ sont dans deux sous-sommes
distinctes, c'est-à-dire si $(w\nu_0)_{\Qdo}\ne (w'\nu_0)_{Q_0}$,
on a $w(\nu_0+\mu)\ne w'(\nu_0+\mu)$ pour tout $\mu$. Il résulte alors de \cite{MW}*{lemme I.4.10} que l'expression~\eqref{eq13.8a} est
elle-même holomorphe en $\nu=\nu_0$. Notons
$\phi_{H,\mu,\xi}(x)$ la valeur de~\eqref{eq13.8a}en $\nu=\nu_0$ et $y=x\ee^H$. On a
\[
E^{\pQ}_{\Qdo}(x\ee^H,\WPhi,\mu)=\sum_{\xi\in \Xi}\phi_{H,\mu,\xi}(x)\ptf
\]
Pour majorer le membre de gauche, on peut fixer $\xi\in \Xi$ et majorer $\phi_{H,\mu,\xi}(x)$.
Pour $\xi$ fixé, considérons cette
dernière fonction comme une forme automorphe en $x$
dépendant de paramètres $H$ et $\mu$. Son terme constant cuspidal relatif à $\Spcusp$
est la sous-somme de~\eqref{eq13.4a} indexée par les $w$ tels que $(w\nu_0)_{\Qdo}=\xi$. On peut
aussi imposer que $w$ vérifie~\eqref{eq13.5a}, sinon le terme correspondant est nul.
Fixons une fonction $B\in \ctyc(\ima(\ga_S^G)^*)$ qui vaut $1$ sur $\Omega$, notons $\Omega'$ le
support de $B$. Le calcul que l'on vient de faire des termes constants cuspidaux montre que, quand $\mu$
reste dans $\Omega'$, la fonction $\phi_{H,\mu,\xi}$ reste dans un ensemble
\[
A\bigl((V_{P'},d_{P'},\Gamma_{P'},\mathbf{n}_{P'})_{P'; \PO\cap L'\subset P'\subset L'}\bigr)
\]
comme dans le paragraphe précédent. Les polynômes en $\HH_{\Spcusp}(x)$ qui apparaissent ont
des coefficients qui sont eux-mêmes fonctions de $H$ et $\mu$.
Ce sont des produits de $\delta_{\Qdo}(\ee^H)^{1/2}$, de
termes exponentiels $\ee^{\langle s(\nu_0+\mu),H\rangle}$, de polynômes
en $H$ de degrés bornés et de fonctions méromorphes de $\mu$.
On ne peut pas affirmer que ces dernières fonctions sont
holomorphes: c'est seulement le terme constant tout entier
que l'on sait holomorphe et cela n'entraîne pas que chacune de ses composantes le soit. Appelons fonction affine réelle sur
un espace vectoriel réel la somme d'une forme linéaire réelle et d'une constante réelle.
Comme toujours dans la théorie des séries d'Eisenstein, on peut, pour tout $\mu_0$, fixer
une famille finie $(\alpha_{b})_{b=1,\dots,d}$ de fonctions affines
réelles non nulles sur $(\ga_{\Qdo}^G)^*$ de sorte que le produit de $\prod_{b}\alpha_{b}(i\mu)$ avec chacun des
coefficients de nos polynômes en $\HH_{\Spcusp}(x)$
soit holomorphe en $\mu_0$. Puisque $\Omega'$ est compact,
on peut fixer une telle famille qui vaut pour chaque point de $\Omega'$. On définit
\[
\psi'_{H,\mu,\xi}=B(\mu)\psi_{H,\mu,\xi}\prod_{b}\alpha_{b}(\ima\mu)\ptf
\]
Le\hspace*{-.2pt} lemme~\ref{2.2.2}\hspace*{-.2pt} (démontré\hspace*{-.2pt} plus\hspace*{-.2pt} bas)\hspace*{-.2pt}
nous\hspace*{-.2pt} dit\hspace*{-.2pt} que, \hspace*{-.2pt}pour\hspace*{-.2pt} majorer\hspace*{-.2pt} la\hspace*{-.2pt} fonction \hspace*{-.2pt}$\phi'_{H,\mu,\xi}(x)$ pour tout $\mu\in \ima(\ga_S^\G)^*$,
\emph{a fortiori} pour majorer $\phi_{H,\mu,\xi}(x)$ pour tout $\mu\in \Omega$, il suffit de
majorer $D'\psi'_{H,\mu,\xi}(x)$ pour un nombre fini d'opérateurs différentiels $D'$ (portant sur la variable $\mu$).
Fixons un tel $D'$. La fonction $D'\psi'_{H,\mu,\xi}$ reste dans un
ensemble
\[
A\bigl((V_{P'},d_{P'},\Gamma_{P'},\mathbf{n}_{P'})_{P'; \PO\cap L'\subset P'\subset L'}\bigr)
\]
(peut-être plus gros que le précédent car une dérivation augmente les degrés des
polynômes). On applique le lemme~\ref{2.1.1}(ii), pour un $\lambda$ que l'on précisera plus tard.
Les $P$ de ce lemme sont ici les $\Spcusp$. Les $\lambda_{P,i}$ sont les $s's(\nu_0+\mu)$ où
$s'$ et $s$ interviennent dans la sous-somme de~\eqref{eq13.4a} correspondant à $\xi$. Leurs parties réelles sont les $w\nu_0$ pour
\[
w\in \weyl^{\pQ}(\ga_{S_{\mathrm{cusp}}},\ga_{\Spcusp})/ \mathcal{W}(\nu_0)
\]
tels que $(w\nu_0)_{\Qdo}=\xi$ et $w$ vérifie~\eqref{eq13.5a}. En notant $\mathcal{W}_{\Spcusp,\xi}$ cet ensemble de $w$,
on obtient une majoration
{\multlinegap0pt
\begin{multline*}
\sup_{\mu\in \Omega}\lvert \psi_{H,\mu,\xi}(x)\rvert \leq c \sup_{\mu\in \ima(\ga_S^G)^*}
\lVert D'\psi'_{H,\mu,\xi}\rVert_{\mathrm{cusp}}\\
{}\times \sum_{\Spcusp}\delta_{\Spcusp}(x)^{1/2}\mspace{-10mu}
\sum_{w\in \mathcal{W}_{\Spcusp,\xi}}\mspace{-10mu}\ee^{\langle\lambda^{\Spcusp}+w\nu_0,\HO(x)\rangle}\bigl(1+\HO(x)\bigr)^N
\end{multline*}}%
où $N$ est un entier assez grand et $c$ une constante absolue. On majore aisément
$\lVert D'\psi'_{H,\mu,\xi}\rVert_{\mathrm{cusp}}$ grâce à la description que l'on a faite ci-dessus
des coefficients des polynômes en $\HH_{\Spcusp}(x)$. En multipliant par la fonction $\prod_{b}\alpha_{b}(\ima\mu)$, on a
supprimé les pôles en $\mu$. Donc ces coefficients sont des
produits de $\delta_{\Qdo}(\ee^H)^{1/2}$, de $\ee^{\langle s(\nu_0+\mu),H\rangle}$, de polynômes de degrés bornés en $H$ et de
fonctions $C^{\infty}$ de $\mu$. Puisque $\mu$ reste dans le
compact $\Omega'$ et que les $(s\nu_0)_{\Qdo}$ sont égaux à $\xi$, on obtient une majoration
\[
\sup\{\lVert D'\psi'_{H,\mu,\xi}\rVert_{\mathrm{cusp}}\mid \mu\in \ima(\ga_S^\G)^*\}\leq
c'\delta_{\Qdo}(\ee^H)^{1/2}\ee^{\langle H,\xi\rangle}(1+\lVert H\rVert)^N
\]
(quitte à accroître le $N$ précédent) avec une constante absolue $c'$. D'où
{\multlinegap0pt\begin{multline}\label{eq13.9a}
\sup_{\mu\in \Omega}\lvert \psi_{H,\mu,\xi}(x)\rvert
\leq cc'\delta_{\PO}(x\ee^H)^{1/2}\sum_{\Spcusp}
\delta_{\Spcusp}(x)^{1/2}\\
{}\times\sum_{w\in \mathcal{W}_{\Spcusp,\xi}}\mspace{-10mu} \ee^{\langle\lambda^{\Spcusp}+w\nu_0,\HO(x)\rangle+\langle\xi,H\rangle}
(1+\lVert \HO(x)\rVert)^N(1+\lVert H\rVert)^N\ptf
\end{multline}}%
Pour obtenir le~(i) de l'énoncé, on prend $\lambda=0$. Quitte à agrandir $N$, on peut essentiellement majorer
l'expression ci-dessus par
\[
\delta_{\Qdo}(\ee^H)^{1/2}\ee^{\langle \xi,H\rangle}(1+\lVert H\rVert)^N\lvert x\rvert^N\ptf
\]
D'après \eqref{eq13.5a}, $\xi_{\Qdo}$ est une combinaison linéaire à coefficients négatifs ou nuls d'éléments de
$\Delta_{\Qdo}^{\pQ}$. Pour $\tau_{\Qdo}^{\pQ}(H)=1$, on a donc
\[
\ee^{\langle\xi,H\rangle}\leq1
\]
et la majoration ci-dessus est celle du~(i).
Pour obtenir le~(ii) de l'énoncé, on doit prouver que l'on peut choisir $\lambda$ tel que, pour tous $\Spcusp$,
$w$ intervenant dans \eqref{eq13.9a} et pour $H$, $x$ vérifiant les hypothèses de~(ii),
\[
\langle \lambda^{\Spcusp}+w\nu_0,\HO(x)\rangle +\langle \xi,H\rangle
\]
reste borné supérieurement. On a $\Spcusp\subset \Qdo$ et $(w\nu_0)_{\Qdo}=\xi$. Le terme précédent est donc égal à
\[
\langle \lambda^{\Spcusp}+w\nu_0,\HO(x)+H\rangle\ptf
\]
D'après \eqref{eq13.5a}, $w\nu_0$ est une somme à coefficients négatifs ou
nuls de projections sur $\ga_{\Spcusp}$ d'éléments de $\Delta_0^{\pQ}$.
En prenant $\lambda$ assez négatif, on peut assurer que
\[
\lambda^{\Spcusp}+w\nu_0
\]
est
une somme à coefficients négatifs ou nuls d'éléments de
$\Delta_0^{\pQ}$. Puisque
\[
\tau_{\PO}^{\pQ}(H+\HO(x)+T_1)=1
\]
le terme $\langle\lambda^{\Spcusp}+w\nu_0,\HO(x)+H\rangle$ est bien borné supérieurement. Cela prouve~(ii).
Montrons que
\begin{equation}
E_{\Qdo,\mathrm{unit}}^{\pQ}(x\ee^H,\WPhi,\mu)=\psi_{H,\mu,0}(x)\ptf\label{eq13.10a}
\end{equation}
Si $s$ vérifie $s(\ga_0^S)\subset \ga_0^{\Qdo}$, on a certainement $(s\nu_0)_{\Qdo}=0$,
c'est-à-dire
\[
s\in \weyl^{\pQ}(\ga_{S_{\mathrm{cusp}}},\Qdo)_0\ptf
\]
Considérons la forme automorphe
\[
\psi'_{H,\mu,0}(x)=\psi_{H,\mu,0}(x)-E_{\Qdo,\mathrm{unit}}^{\pQ}(x\ee^H,\WPhi,\mu)\ptf
\]
On a vu que les exposants cuspidaux de $\psi_{H,\mu,0}$ relatifs à un
parabolique $\Spcusp$ vérifiaient la propriété~\eqref{eq13.5a}. Mais ceux de
\[
E_{\Qdo,\mathrm{unit}}^{\pQ}(x\ee^H,\WPhi,\mu)
\]
la vérifient aussi. En effet, on peut appliquer à chaque composante
\[
E_{\Qdo}(y,\Mint(s,\mu)\WPhi,\mu)
\]
(\cf~\eqref{eq13.7a}) l'analogue de \eqref{eq13.5a} où l'on remplace $\pQ$ par $\Qdo$.
Cet analogue nous dit que ses exposants cuspidaux sont de la forme $s's\nu_0$ pour des
\[
s'\in \weyl^{\Qdo}(\ga_{s(S_{\mathrm{cusp}})},\ga_{\Spcusp})
\]
tels que
\[
\htau^{Q_{0,s',s}}_{\Spcusp}(-s's\nu_0)=1
\]
où $Q_{0,s',s}$ est le plus petit sous-groupe parabolique
standard de $\Qdo$ tel que
\[
\ga_{\pQ_{0,s',s}}\subset s'(\ga_{s(S)})\ptf
\]
En posant $w=s's$, on vérifie que la condition $s(\ga_0^S)\subset \ga_0^{Q_0}$
entraîne l'égalité
\[
Q_{0,s',s}=\Q'_{w}
\]
d'où l'assertion. Donc les exposants cuspidaux de $\psi'_{H,\mu,0}$
vérifient aussi~\eqref{eq13.5a}. Par construction, ils sont aussi de la forme
$s's\nu_0$ avec
\[
(s\nu_0)_{\Qdo}=0,\qquad s(\ga_0^S)\nsubset \ga_0^{\Qdo}\Qquad{et}
s'\in \weyl^{\Qdo}(\ga_{s(S_{\mathrm{cusp}})},\ga_{\Spcusp})\ptf
\]
En posant $w=s's$, on a encore $(w\nu_0)_{\Qdo}=0$. La relation
\[
\htau^{\pQ{w}}_{\Spcusp}(w\nu_0)=1
\]
entraîne alors que
\[
\Delta_{\Spcusp}^{\Q'_{w}}\subset \Delta_{\Spcusp}^{\Qdo}
\]
autrement dit $\Q'_{w}\subset \Qdo$.
Par définition de $\Q'_{w}$, on a alors $\ga_{\Qdo}
\subset w(\ga_S)$, ce qui équivaut à $w(\ga_0^S)\subset \ga_0^{\Qdo}$, ou encore à
$s(\gao^S)\subset \gao^{\Qdo}$. C'est une
contradiction, sauf si l'ensemble des exposants cuspidaux de $\psi'_{H,\mu,0}$ est vide.
Donc cet ensemble est vide et cela implique comme on le sait que $\psi'_{H,\mu,0}=0$. D'où~\eqref{eq13.10a}.
Pour démontrer le~(iii) de l'énoncé, on reprend la preuve ci-dessus dans le cas $\xi=0$. Cette hypothèse fait disparaître le
$H$ dans l'expression
\[
\langle\lambda^{\Spcusp}+w\nu_0,\HO(x)+H\rangle\ptf
\]
On peut donc la borner supérieurement sous la seule hypothèse que $x\in \Sieg^{L_0}$. D'où~(iii).
D'après~\eqref{eq13.10a}, on a
\[
E_{\Qdo,+}^{\pQ}(x\ee^H,\WPhi,\mu)=\sum_{\xi\in \Xi; \xi\ne 0}\psi_{H,\mu,0}(x)\ptf
\]
Pour démontrer~(iv), il suffit
de prouver que l'on peut choisir $\lambda$ tel que, pour tout
$\xi\ne 0$, tout $\Spcusp$ et tout $w\in \mathcal{W}_{\Spcusp,\xi}$ l'on ait une majoration
\[
\langle\lambda^{\Spcusp}+w\nu_0,\HO(x)+H\rangle \leq c-R\inf_{ \alpha\in \Delta_0^{\pQ}-\Delta_0^{Q_0}}
\alpha\bigl(H+\HO(x)\bigr)
\]
pour des constantes $c>0$, $R>0$ convenables. On peut supposer que
\[
\lambda^{\Spcusp}+w\nu_0
\]
est une combinaison linéaire à coefficients négatifs ou nuls d'éléments de
$\Delta_0$. L'hypothèse $\xi\ne 0$, c'est-à-dire $(w\nu_0)_{\Qdo}\ne 0$, entraîne qu'il y a au moins un
$\alpha\in \Delta_0^{\pQ}-\Delta_0^{\Qdo}$ dont le coefficient est non
nul. L'assertion s'ensuit. 
\end{proof}

\begin{lemme}\label{2.2.2}
Soient $m\geq1$ un entier et $(\beta_{b})_{b=1,\ldots,d}$ une famille de fonctions affines réelles sur
$\RM^m$\textup, non nulles. Pour tout opérateur différentiel $D$
à coefficients polynomiaux sur $\RM^m$\textup, il existe une famille finie $D'_1,\dots,D'_{k}$ de tels opérateurs tels que pour
toute fonction $h$ lisse sur $\RM^m$\textup, on ait la
majoration
\[
\sup_{ y\in \RM^m}\lvert Dh(y)\rvert \leq \sum_{l=1}^{k}\sup_{ y\in \RM^m}\lvert D'_{l}H(y)\rvert
\]
où $H$ est la fonction
\[
H(y)=h(y)\prod_{b=1}^{d}\beta_{b}(y)\ptf
\]
\end{lemme}

\begin{proof} 
Ce lemme est élémentaire et certainement bien connu, on donne une preuve pour la simple commodité
du lecteur. Par récurrence sur $d$, on se ramène au
cas où $d=1$. Par changement de variables, on peut supposer que l'unique forme affine est la première coordonnée
$y\mapsto y_1$. On peut aussi se limiter aux opérateurs $D$ de
la forme $D_1D_{2}$ où $D_1$ ne dépend que de la variable $y_1$ et $D_{2}$ ne dépend que des variables
$y_{2},\dots,y_{m}$. Supposons le lemme résolu pour $m=1$. Il associe
à $D_1$ des opérateurs $D'_{1,1},\dots,D'_{1,k}$. Alors les opérateurs $D'_{1,1}D_{2},\dots,D'_{1,k}D_{2}$
valent pour notre opérateur $D_1D_{2}$. On peut donc supposer $m=1$ et que
l'unique forme affine est l'identité $y\mapsto y$. On peut aussi supposer que $D=\delta^{i}y^j$, où $\delta$ est
l'opérateur $\dd/\dd y$. Si $j\geq1$, on prend $k=1$ et $D'_1=\delta^{i}y^{j-1}$.
Supposons $j=0$. Il existe un opérateur $D[i]$ obtenu par les règles de dérivation usuelles tel que, pour $h$ de
Schwartz et $H(y)=yh(y)$, on ait $\delta^{i}h(y)=y^{-i-1}D[i]H(y)$.
La fonction $D[i]H$ s'annule à l'ordre au moins $i+1$ en $0$. Par une application successive du théorème des accroissements
finis, on trouve pour tout $y$ des points $y_1,\dots,y_{i}$ tels que
\[
\lvert D[i]H(y)\rvert \leq\lvert y\delta D[i]H(y_{1})\rvert \le \lvert yy_1\delta^2  D[i]H(y_{2})\rvert
\leq\dots\leq\lvert yy_{1}\dotsm y_{i}\delta^{i+1}D[i]H(y_{i+1})\rvert.
\]
et chaque point appartient au segment joignant $0$ à son prédécesseur. En particulier
\[
\lvert y_{i+1}\rvert \leq\dots\leq\lvert y_{1}\rvert \leq\lvert y\rvert \ptf
\]
Donc
\[
\lvert D[i]H(y)\rvert \leq \lvert y\rvert ^{i+1}\sup_{y'\in \RM}\lvert\delta^{i+1}D[i]H(y')\rvert \;,
\]
puis
\[
\lvert \delta^{i}h(y)\rvert \leq \sup_{y'\in \RM}\lvert\delta^{i+1}D[i]H(y')\rvert \ptf
\]
En posant $k=1$ et $D'_1=\delta^{i+1}D[i]$, la conclusion du lemme est vérifiée.
\end{proof}

\section{Simplification du terme constant}

On fixe deux éléments $\WPhi\in \mathcal{B}^Q(\theto\sigma)_{\chi}$, $\WPsi\in \base^{\pQ}_\chi(\sigma)$.
On fixe une fonction $B\in C_{c}^{\infty}(\ima(\ga_S^\G)^*)$. On définit comme
au lemme~\ref{W1.10} un terme $A^T(B,H)$ pour $H\in \ga_{\Qdo}$, puis un terme $A^T(B)$.
En remplaçant dans ces définitions
les fonctions $E_{\Qdo}^Q$ et $E_{\Qdo}^{\pQ}$ par les
fonctions $E_{\Qdo,\mathrm{unit}}^Q$ et $E_{\Qdo,\mathrm{unit}}^{\pQ}$ du paragraphe précédent, on obtient de nouveaux termes
$A^T_{\mathrm{unit}}(B,H)$ et $A^T_{\mathrm{unit}}(B)$.

\begin{lemme}\label{2.3}
\begin{enumerate}[(i)]
\item Les intégrales définissant les expressions $A^T(B,H)$\textup, $A^T(B)$\textup,  $A^T_{\mathrm{unit}}(B,H)$ et $A^T_{\mathrm{unit}}(B)$ sont absolument
convergentes.
\item Pour tout réel $r$\textup, il existe $c>0$ tel que
\[
\lvert A^T(B)-A^T_{\mathrm{unit}}(B)\rvert \leq c\dPO(T)^{-r}\ptf
\]
\end{enumerate}
\end{lemme}

\begin{proof} 
On a déjà démontré au lemme~\ref{W1.10} les assertions du (i) concernant $A^T(B,H)$ et $A^T(B)$.
On l'a démontré pour une fonction $\varphi B$, où $\varphi$ était de Paley-Wiener,
mais ce n'est pas une restriction: toute fonction $\cty$ à support compact est produit d'une telle fonction et d'une fonction
de Paley-Wiener. Donnons une nouvelle démonstration qui
s'applique aussi bien aux termes $A^T_{\mathrm{unit}}(B,H)$ et $A^T_{\mathrm{unit}}(B)$.
Soit $H$ tel que
\[
\kappa^{\eta T}(H^Q-T_{\Qdo}^Q)\tsQR(H-T)\phi_{\Qdo}^Q(H-T)=1\ptf
\]
D'après l'équation~\eqref{eq12.5c} de la proposition~\ref{W1.7}, cela
implique $\tau_{\Qdo}^P(H)=1$. On peut appliquer la proposition~\ref{2.2.1}(i) à
\[
E_{\Qdo}^{\pQ}(x\ee^Hk,\WPsi,\mu)
\]
pour $\mu$ dans le support de $B$.
On peut aussi appliquer la proposition similaire à
\[
E_{\Qdo}^Q(x\ee^Hk,\WPhi,\tmu)
\]
où on a posé $\tmu=\theto(\mu)$
On applique ensuite à cette dernière fonction la proposition~\ref{rapdec}.
Il en résulte pour tout réel $r$ une majoration
\[
\lvert \tronc^{T[H^Q],\Qdo}E_{\Qdo}^Q(x\ee^Hk,\WPhi,\tmu)\widebar{E_{\Qdo}^{\pQ}(x\ee^Hk,\WPsi,\mu)}B(\mu)\rvert
\leq c\lvert B(\mu)\rvert\,\lvert x\rvert^{-r},
\]
où $c$ dépend de $r$, $T$ et $H$, mais pas de $x$, $k$, $\mu$.
L'expression ci-dessus est intégrable en ces trois dernières variables, ce qui prouve l'assertion~(i) pour $A^T(B,H)$.
Comme dans la preuve de la proposition~\ref{W1.7}, on décompose $A^T(B,H)$ en
\[
A^T_{C}(B,H)+A^T_{\tronc-C}(B,H)
\]
où
{\multlinegap0pt\begin{multline*}
A^T_{C}(B,H)=\int_{L_0(F)\backslash L_0(\adef)^1}\int_\K\kappa^{\eta T}(H^Q-T_{\Qdo}^Q)\tsQR(H-T)
\phi_{\Qdo}^Q(H-T)\delta_{\Qdo}(\ee^H)^{-1}\\
{}\times\int_{\ima(\ga_S^\G)^*}F_{\PO}^{\Qdo}(x,T[H^Q])E_{\Qdo}^Q(x\ee^Hk,\WPhi,\tmu)
\widebar{E_{\Qdo}^{\pQ}(x\ee^Hk,\WPsi,\mu)}B(\mu) \dd \mu\dd k\dd x
\end{multline*}}%
et
{\multlinegap0pt\begin{multline*}
A^T_{\tronc-C}(B,H)\\
\shoveleft{=\int_{L_0(F)\backslash L_0(\adef)^1}\int_\K \kappa^{\eta T}(H^Q-T_{\Qdo}^Q)
\tsQR(H-T)\phi_{\Qdo}^Q(H-T)\delta_{\Qdo}(\ee^H)^{-1}}\\
{}\times\smash[b]{\int_{\ima(\ga_S^\G)^*}}
\bigl(\tronc^{T[H^Q],\Qdo}(E_{\Qdo}^Q)(x\ee^Hk,\WPhi,\tmu)-F_{\PO}^{\Qdo}(x,T[H^Q])
E_{\Qdo}^Q(x\ee^Hk,\WPhi,\tmu)\bigr)\\
{}\times\widebar{E_{\Qdo}^{\pQ}(x\ee^Hk,\WPsi,\mu)}B(\mu) \dd\mu\dd k\dd x\ptf
\end{multline*}}%

Considérons la première intégrale. On a
\begin{multline*} 
F_{\PO}^{\Qdo}(x,T[H^Q])E_{\Qdo}^Q(x\ee^Hk,\WPhi,\tmu)\widebar{E_{\Qdo}^{\pQ}(x\ee^Hk,\WPsi,\mu)}\\
=\delta_Q(\ee^{H_Q})^{1/2}\delta_{\pQ}(\ee^{H_{\pQ}})^{1/2}\ee^{\langle \theta_0^{-1}H_Q-H_{\pQ},\mu_{\pQ}\rangle}
F_{\PO}^{\Qdo}(x,T[H^Q])\\
{}\times E_{\Qdo}^Q\bigl(x\ee^{H^Q}k,\WPhi,\theto(\mu^{\pQ})\bigr)
\widebar{E_{\Qdo}^{\pQ}(x\ee^{H^{\pQ}}k,\WPsi,\mu^{\pQ})}\ptf
\end{multline*}
On peut écrire l'intégrale intérieure sous la forme
\begin{multline*}
\delta_Q(\ee^{H_Q})^{1/2}\delta_{\pQ}(\ee^{H_{\pQ}})^{1/2}\\
{}\times\int_{\ima(\ga_S^{\pQ})^*}F_{\PO}^{\Qdo}(x,T[H^Q]) E_{\Qdo}^Q\bigl(x\ee^{H^Q}k,\WPhi,\theto(\mu^{\pQ})\bigr)
\widebar{E_{\Qdo}^{\pQ}(x\ee^{H^{\pQ}}k,\WPsi,\mu^{\pQ})}\\
{}\times\smash[t]{\int_{\ima(\ga_{\pQ}^G)^*}}\ee^{\langle\theta_0^{-1}H_Q-H_{\pQ},\mu_{\pQ}\rangle} B(\mu^{\pQ}+\mu_{\pQ})\dd\mu_{\pQ}\dd\mu^{\pQ}\ptf
\end{multline*}
L'intégrale intérieure de cette expression est la transformée de Fourier partielle de $B$ en la variable $\mu_{\pQ}$. On peut la
majorer par
\[
\widehat{B}_1(\theta_0^{-1}H_Q-H_{\pQ})B_{2}(\mu^{\pQ})
\]
où $\widehat{B}_1$ est une fonction de Schwartz sur $a_{\pQ}^G$ et $B_{2}\in C_{c}^{\infty}(\ima(\ga_S^{\pQ})^*)$, toutes
deux à valeurs positives ou nulles. D'où
\begin{multline*}
\lvert A^T_{C}(B,H)\rvert\\
 \leq \kappa^{\eta T}(H^Q-T_{\Qdo}^Q)\tsQR(H-T)\phi_{\Qdo}^Q(H-T)
\delta_{\Qdo}(\ee^{H^Q+H^{\pQ}})^{-1/2}\widehat{B}_1(\theta_0^{-1}H_Q-H_{\pQ})\\
{}\times\int_{L_0(F)\backslash L_0(\adef)^1}\int_\K\int_{i(\ga_S^{\pQ})^*}
F_{\PO}^{\Qdo}(x,T[H^Q])B_{2}(\mu^{\pQ})\qquad\qquad\\
{}\times\bigl\lvert E_{\Qdo}^Q\bigl(x\ee^{H^Q}k,\WPhi,\theto(\mu^{\pQ})\bigr)
\widebar{E_{\Qdo}^{\pQ}(x\ee^{H^{\pQ}}k,\WPsi,\mu^{\pQ})}\bigr\rvert \dd\mu^{\pQ}\dd k\dd x\ptf
\end{multline*}
D'après l'équation~\eqref{eq12.7c} de la proposition~\ref{W1.7} il existe $T_1$ (ne dépendant d'aucune variable) tel que
\[
\tau_{\PO}^P(H+\HO(x)+T_1)=1
\]
pour
tous $x\in \Sieg^{\Qdo}$ et $H$ tels que les fonctions
ci-dessus soient non nulles. On peut appliquer la proposition~\ref{2.2.1}(ii) aux deux séries d'Eisenstein.
On obtient que l'intégrale intérieure est essentiellement majorée par
\[
\delta_{\Qdo}(\ee^{H^Q+H^{\pQ}})^{1/2}(1+\lVert H\rVert)^D\int_{\Sieg^{\Qdo}}\delta_{\PO}(x)
F_{\PO}^{\Qdo}(x,T[H^Q])(1+\lVert\HO(x)\rVert)^D\dd x
\]
pour un entier $D$ assez grand. Comme dans la preuve de la proposition~\ref{W1.7},
cette dernière intégrale est essentiellement majorée par $\mathbf{d}_{\PO}(T)^D$, quitte à accroître $D$. D'où
\begin{multline*}
\lvert A^T_{C}(B,H)\rvert \leq \dPO(T)^D\kappa^{\eta T}(H^Q-T_{\Qdo}^Q)\tsQR(H-T)
\phi_{\Qdo}^Q(H-T)\\
{}\times\widehat{B}_1(\theta_0^{-1}H_Q-H_{\pQ})(1+\lVert H\rVert)^D\ptf
\end{multline*}
D'après le lemme~\ref{ex10.3.2} appliqué à $P'=Q$, on a sur le support de cette fonction une majoration
\[
\lVert H-T_{\Qdo}\rVert \ll  \lVert q(H)\rVert \ptf
\]
Puisque $\widehat{B}_1$ est de Schwartz, on a pour tout réel $r$ une majoration
\[
\widehat{B}_1(\theta_0^{-1}H_Q-H_{\pQ})\ll (1+\lVert \theta_0^{-1}H_Q-H_{\pQ}\rVert)^{-r}\ll 
(1+\lVert q(H)\rVert)^{-r}\ptf
\]
D'où une majoration
\begin{equation}
\lvert A^T_{C}(B,H)\rvert \ll  \dPO(T)^D (1+\lVert H-T_{\Qdo}\rVert)^{-r}\ptf\label{eq13.1b}
\end{equation}
Cette expression est intégrable en $H$. Considérons maintenant l'intégrale
\[
A^T_{\tronc-C}(B,H)\ptf
\]
La première partie du raisonnement ci-dessus s'applique. D'où
{\multlinegap0pt\begin{multline*}
\lvert A^T_{\tronc-C}(B,H)\rvert\\
 \leq \kappa^{\eta T}(H^Q-T_{\Qdo}^Q)\tsQR(H-T) \phi_{\Qdo}^Q(H-T)
\delta_{Q_0}(\ee^{H^Q+H^{\pQ}})^{-1/2}\widehat{B}_1(\theta_0^{-1}H_Q-H_{\pQ})\\
{}\times\smash[b]{\int_{L_0(F)\backslash L_0(\adef)^1}\int_\K\int_{\ima(\ga_S^{\pQ})^*}}
\bigl\lvert
\tronc^{T[H^Q],\Qdo}E_{\Qdo}^Q\bigl(x\ee^{H^Q}k,\WPhi,\theto(\mu^{\pQ})\bigr)\qquad\\
\shoveright{{}-F_{\PO}^{\Qdo}(x,T[H^Q]) E_{\Qdo}^Q\bigl(x\ee^{H^Q}k,\WPhi,\theto(\mu^{\pQ})\bigr) \bigr\rvert}\\
{}\times\bigl\lvert E_{\Qdo}^{\pQ}(x\ee^{H^{\pQ}}k,\WPsi,\mu^{\pQ})\bigr\rvert B_{2}(\mu^{\pQ}) \dd\mu^{\pQ}\dd k\dd x\ptf
\end{multline*}}%
On applique la proposition~\ref{2.2.1} puis la proposition~\ref{propoC}.
On obtient pour un certain entier $D$ et pour tout réel $r$ une majoration
\begin{multline*}
\bigl\lvert \tronc^{T[H^Q],\Qdo}E_{\Qdo}^Q\bigl(x\ee^{H^Q}k,\WPhi,\theto(\mu^{\pQ})\bigr)-F_{\PO}^{\Qdo}(x,T[H^Q])
E_{\Qdo}^Q\bigl(x\ee^{H^Q}k,\WPhi,\theto(\mu^{\pQ})\bigr)\bigr\rvert\\
\shoveright{{}\times\lvert E_{\Qdo}^{\pQ}(x\ee^{H^{\pQ}}k,\WPsi,\mu^{\pQ})\rvert B_{2}(\mu^{\pQ})}\\
\ll 
\ee^{-r\mathbf{d}_{\PO\cap L_0}^{L_0}(T[H^Q])}\delta_{\Qdo}(\ee^{H^Q+H^{\pQ}})^{1/2}(1+\lVert H\rVert)^D\lvert x\rvert^{-r}B_{2}(\mu^{\pQ})
\end{multline*}
valable pour tout $x\in\Sieg^{\Qdo}$, $k$ $\mu^{\pQ}$ et $H$ tel que
\[
\kappa^{\eta T}(H^Q-T_{\Qdo}^Q)\tsQR(H-T)\phi_{\Qdo}^Q(H-T)=1\ptf
\]
Cette expression est intégrable en $x$, $k$ et $\mu^{\pQ}$. D'autre part, $T[H^Q]$
est \og plus régulier\fg que $T$, donc
\[
\dPO(T)\ll \mathbf{d}_{P_0\cap L_0}^{L_0}(T[H^Q])\ptf
\]
D'où
\begin{multline*}
\lvert A^T_{\tronc-C}(B,H)\rvert\ll  \ee^{-r\dPO(T)}\kappa^{\eta T}(H^Q-T_{\Qdo}^Q)\\
{}\times\tsQR(H-T)\phi_{Q_0}^Q(H-T) \widehat{B}_1(\theta_0^{-1}H_Q-H_{\pQ})(1+\lVert H\rVert)^D\ptf
\end{multline*}
On peut reprendre la preuve de~\eqref{eq13.1b} et on obtient cette fois pour tout $r$ une majoration
\begin{equation}
\lvert A^T_{\tronc-C}(B,H)\rvert \ll  \ee^{-r\dPO(T)} (1+\lVert H-T_{\Qdo}\rVert)^{-r}\ptf
\label{eq13.2b}
\end{equation}
Ceci est encore intégrable en $H$. Cela prouve d'abord que l'intégrale $A^T(B)$ est absolument convergente.
On obtient de plus une majoration
\begin{equation}
\lvert A^T(B)-A^T_{C}(B)\rvert=\lvert A^T_{\tronc-C}(B)\rvert \ll \ee^{-r\dPO(T)},
\label{eq13.3b}
\end{equation}
où
\[
A^T_{C}(B)=\int_{\ga_{\Qdo}^G}A^T_{C}(B,H)\dd H
\]
et
\[
A^T_{\tronc-C}(B)=\int_{\ga_{\Qdo}^G}A^T_{\tronc-C}(B,H)\dd H\ptf
\]
La même démonstration s'applique aux termes $A^T_{\mathrm{unit}}(B,H)$ et $A^T_{\mathrm{unit}}(B)$.
On obtient aussi pour ces termes une relation analogue à~\eqref{eq13.3b}.
Il résulte de~\eqref{eq13.3b} et de la relation analogue que, pour prouver le~(ii) de l'énoncé, il suffit de majorer
\[
A^T_{C}(B)-A^T_{C,\mathrm{unit}}(B)\ptf
\]
On a
{\multlinegap0pt
\begin{multline*}
A^T_{C}(B,H)-A^T_{C,\mathrm{unit}}(B,H)\\
=\int_{L_0(F)\backslash L_0(\adef)^1}\int_\K\kappa^{\eta T}(H^Q-T_{\Qdo}^Q)
\tsQR(H-T)\phi_{\Qdo}^Q(H-T)\delta_{\Qdo}(\ee^H)^{-1}\qquad\qquad\\
{}\times\int_{\ima(\ga_S^\G)^*}F_{\PO}^{\Qdo}(x,T[H^Q])E_{\Qdo,+}^Q(x\ee^Hk,\WPhi,\tmu)
\widebar{E_{\Qdo}^{\pQ}(x\ee^Hk,\WPsi,\mu)}B(\mu) \dd \mu\dd k\dd x\ptf
\end{multline*}}%
On utilise la proposition~\ref{2.2.1}(iv) pour majorer cette expression. Il apparaît un terme
\[
\sup_{\alpha\in \Delta_0^{\pQ}-\Delta_0^{\Qdo}}\ee^{-R\alpha\left(H+\HO(x)\right)}\ptf
\]
Or, d'après l'équation~\eqref{eq12.7c} de la proposition~\ref{W1.7}, on a une minoration
\[
\inf_{\alpha\in \Delta_0^{\pQ}-\Delta_0^{\Qdo}}\alpha\bigl(H+\HO(x)\bigr)\geq c\dPO(T)
\]
pour un $c>0$ convenable, pourvu que $H$ et $x\in \Sieg^{L_0}$ appartiennent aux supports de nos fonctions.
Le terme ci-dessus est donc majoré par
\begin{equation}
e^{-Rc\dPO(T)}\ptf\label{eq13.4b}
\end{equation}
Le calcul se poursuit comme précédemment et on obtient une relation analogue à~\eqref{eq13.1b}, où se glisse ce terme
\eqref{eq13.4b} et donc
\[
\lvert A^T_{C}(B,H)-A^T_{C,\mathrm{unit}}(B,H)\rvert \ll \ee^{-Rc\dPO(T)}\dPO(T)^D(1-\lVert H-T_{\Qdo}\rVert)^{-r}
\]
où $r$ est quelconque. L'intégrale en $H$ de cette expression est évidemment essentiellement bornée par
$\dPO(T)^{-r}$ pour tout réel $r$.
\end{proof}

\section{Simplification du produit scalaire}

Pour $H\in \ga_{\Qdo}^G$ et $\mu,\nu\in \ima(\ga_S^\G)^*$ et $\tmu=\theto\mu$ posons
\begin{multline*}
\omega^{T,\Qdo}(H,\tmu,\nu)
=\sum_{\Sp; \PO\subset \Sp\subset \Qdo}\sum_{s\in
\weyl^Q(\theto(\ga_S),\ga_{\Sp})}\sum_{t\in \weyl^{\pQ}(\ga_S,\ga_{\Sp})}
\ee^{\langle s\tmu-t\nu,H+T[H^Q]\rangle}\\
{}\times \epsilon_{\Sp}^{\Qdo}(s\tmu-t\nu)
\langle\Mint(t,\nu)^{-1}\Mint(s,\tmu)\WPhi,\WPsi\rangle\ptf
\end{multline*}
Le produit scalaire est celui de deux éléments de $\mathcal{A}(\mathbf{X}_S,\sigma)$.
Il résulte de la proposition~\ref{GMspec} que $\omega^{T,\Qdo}
(H,\tmu,\nu)$ est holomorphe en $\mu$ et $\nu$. On note
\[
\omega^{T,\Qdo}(H,\mu)=\omega^{T,\Qdo}(H,\theto\mu,\mu)\ptf
\]
Posons
\[
A^T_{\mathrm{pure}}(B,H)=\kappa^{\eta T}(H^Q-T_{\Qdo}^Q)
\tsQR(H-T)\phi_{\Qdo}^Q(H-T)\int_{\ima(\ga_S^\G)^*}
\omega^{T,\Qdo}(H,\mu)B(\mu)\dd \mu
\]
et
\[
A^T_{\mathrm{pure}}(B)=\int_{\ga_{\Qdo}^G}A^T_{\mathrm{pure}}(B,H)\dd H\ptf
\]

\begin{proposition}\label{W2.4}
Les deux intégrales ci-dessus sont absolument convergentes. Pour tout réel $r$\textup, il existe $c>0$ tel que
\[
\lvert A^T(B)-A^T_{\mathrm{pure}}(B)\rvert \leq c\dPO(T)^{-r}\ptf
\]
\end{proposition}

\begin{proof} 
L'expression $A^T_{\mathrm{pure}}(B,H)$ est l'intégrale d'une fonction $C^{\infty}$ à support compact,
elle est donc absolument convergente. Remarquons que les opérateurs d'entrelacement
$\Mint(t,\mu)$ et $\Mint(s,\tmu)$ ne dépendent en fait que de $\mu^{\pQ}$. D'autre part
\[
s\tmu-t\mu=\bigl(s\theto(\mu^{\pQ})-t(\mu^{\pQ})\bigr)+\theto(\mu_{\pQ})-\mu_{\pQ}
\]
et $T[H^Q]$
appartient à $\ga_0^{\Qdo}$. D'où l'égalité
\[
\omega^{T,\Qdo}(H,\mu)=\ee^{\langle\theto(\mu_{\pQ})-\mu_{\pQ},H\rangle}\omega^{T,\Qdo}(H,\mu^{\pQ})\ptf
\]
Comme dans la preuve du paragraphe précédent, on a
\[
A^T_{\mathrm{pure}}(B,H)
=
\int_{\ima(\ga_S^{\pQ})^*}\omega^{T,\Qdo}(H,\mu^{\pQ})\int_{\ima(\ga_{\pQ}^G)^*}\!\!\!
\ee^{\langle\theto(\mu_{\pQ})-\mu_{\pQ},H\rangle}B(\mu^{\pQ}+\mu_{\pQ})\dd \mu_{\pQ}\dd \mu^{\pQ}
\]
d'où une majoration
\begin{multline*}
\lvert A^T_{\mathrm{pure}}(B,H)\rvert \ll \kappa^{\eta T}(H^Q-T_{\Qdo}^Q)\tsQR(H-T)\phi_{\Qdo}^Q(H-T)\\
{}\times\widehat{B}_1(\theta_0^{-1}(H_Q)-H_{\pQ})
\int_{\ima(\ga_S^{\pQ})^*}\lvert \omega^{T,\Qdo}(H,\mu^{\pQ})\rvert B_{2}(\mu^{\pQ})\dd \mu^{\pQ}\ptf
\end{multline*}
Montrons qu'il existe un entier $D$ tel que l'on ait une majoration
\begin{equation}
\lvert \omega^{T,\Qdo}(H,\mu^{\pQ})\rvert B_{2}(\mu^{\pQ})\ll (1+\lVert H\rVert)^D\dPO(T)^D\ptf\label{eq13.1c}
\end{equation}
Il suffit de majorer $\omega^{T,\Qdo}(H,\tmu,\nu)$
pour $\mu$ et $\nu$ parcourant un sous-ensemble compact $\Omega$ de
$\ima(\ga_S^\G)^*$. Si chaque terme de la somme définissant
$\omega^{T,\Qdo}(H,\tmu,\nu)$ était holomorphe en $\mu^{\pQ}$, la majoration serait évidente. Il y a des pôles
dus aux fonctions $\epsilon_{\Sp}^{\Qdo}(s\tmu-t\nu)$.
On peut raisonner comme dans la preuve de la proposition~\ref{2.2.1}:
le lemme~\ref{2.2.2} nous ramène à majorer un nombre fini de
fonctions
\[
Dp(\tmu,\nu)\omega^{T,\Qdo}(H,\tmu,\nu)
\]
où
$D$ est un opérateur différentiel en $\tmu$, $\nu$ et $p(\tmu,\nu)$
est un produit de fonctions affines réelles tel qu'après
multiplication par $p$, plus aucune de nos fonctions n'ait de pôles pour $\mu,\nu\in \Omega$.
La majoration de
\[
Dp(\tmu,\nu)\omega^{T,\Qdo}(H,\tmu,\nu)
\]
est évidente et~\eqref{eq13.1c} s'ensuit.
De~\eqref{eq13.1c} résulte une majoration
\begin{multline*}
\lvert A^T_{\mathrm{pure}}(B,H)\rvert \ll \kappa^{\eta T}(H^Q-T_{\Qdo}^Q)\tsQR(H-T)\phi_{\Qdo}^Q(H-T)\\
{}\times\widehat{B}_1(\theta_0^{-1}(H_Q)-H_{\pQ}) (1+\lVert H\rVert)^D\dPO(T)^D\ptf
\end{multline*}
En faisant appel au lemme~\ref{ex10.3.2}, appliqué au cas $P'=Q$, on voit que cette dernière expression est intégrable en $H$.
D'où la première assertion de l'énoncé.
En reprenant la preuve ci-dessus et celle du lemme~\ref{2.3}, on voit que
\begin{multline*}
\lvert A^T_{\mathrm{pure}}(B,H)-A^T_{\mathrm{unit}}(B,H)\rvert \\
\ll 
\kappa^{\eta T}(H^Q-T_{\Qdo}^Q)\tsQR(H-T)\phi_{\Qdo}^Q(H-T)
\widehat{B}_1(\theta_0^{-1}(H_Q)-H_{\pQ})\\
{}\times\int_{\ima(\ga_S^{\pQ})^*}\lvert X(T,H,\mu^{\pQ})\rvert B_{2}(\mu^{\pQ})\dd \mu^{\pQ}
\end{multline*}
où $X(T,H,\mu^{\pQ})$ est égal à
{\multlinegap0pt\begin{multline*}
\omega^{T,\Qdo}(H,\mu^{\pQ})-\delta_{\Qdo}(\ee^H)^{-1}\\
{}\times\int_{L_0(F)\backslash L_0(\adef)^1}
\int_\K 
\tronc^{T[H^Q],\Qdo}
E^Q_{\Qdo,\mathrm{unit}}\bigl(x\ee^Hk,\WPhi,\theto(\mu^{\pQ})\bigr)\\
{}\times
\widebar{E^{\pQ}_{\Qdo,\mathrm{unit}}(x\ee^Hk,\WPsi,\mu^{\pQ})}\dd k\dd x\ptf
\end{multline*}}%
Pour $s\in \weyl^Q(\theto(\ga_S),\Qdo)$, $ t\in \weyl^{\pQ}(\ga_S,\Qdo)$ et
$\mu,\nu\in\ima(\ga_S^G)^*$, posons $\tmu=\theto\mu$ et
\begin{multline*}
\omega^{T,\Qdo}(H,\tmu,\nu,s,t)\\
=\sum_{\{\Sp\mid \PO\subset \Sp\subset \Qdo\}}
 \sum_{s'\in \weyl^{\Qdo}(\ga_{s\theto(S)},\ga_{\Sp})}\sum_{t'\in \weyl^{\Qdo}(t(\ga_S),\ga_{\Sp})}
\ee^{\langle s's\tmu-t't\nu,H+T[H^Q]\rangle} \\
{}\times\epsilon_{\Sp}^{\Qdo}(s\tmu-t\nu)\langle \Mint(t't,\nu)^{-1}\Mint(s's,\tmu)\WPhi,\WPsi\rangle\ptf
\end{multline*}
Ceci est encore holomorphe en $\tmu$, $\nu$. On note
\[
\omega^{T,\Qdo}(H,\mu,s,t)=\omega^{T,\Qdo}(H,\theto\mu,\mu,s,t)\ptf
\]
On a l'égalité
\[
\omega^{T,\Qdo}(H,\mu)=\sum_{s\in \weyl^Q(\theto(\ga_S),\Qdo)}
\sum_{ t\in \weyl^{\pQ}(\ga_S,\Qdo)}\omega^{T,\Qdo}(H,\mu,s,t)\ptf
\]
Pour $s$ et $t$ comme ci-dessus, posons
{\multlinegap0pt\begin{multline*}
Y(T,H,\mu^{\pQ},s,t)=\delta_{\Qdo}(\ee^H)^{-1}\\
{}\times\smash[b]{\int_{L_0(F)\backslash L_0(\adef)^1}\int_\K}
\tronc^{T[H^Q],\Qdo}E_{\Qdo}\bigl(x\ee^Hk,\Mint\bigl(s,\theto (\mu^{\pQ})\bigr)\WPhi,\theto(\mu^{\pQ})\bigr)\\ 
{}\times\widebar{E_{\Qdo}(x\ee^Hk,\Mint(t,\mu^{\pQ})\WPsi,\mu^{\pQ})}\dd k\dd x\ptf
\end{multline*}}%
Grâce à l'équation~\eqref{eq13.7a} de la section~\ref{W2.2}, on peut décomposer $X(T,H,\mu^{\pQ})$ en une somme
\[
\sum_{s\in \weyl^\Q(\theto(\ga_S),\Qdo), t\in \weyl^{\pQ}(\ga_S,\Qdo)}
\bigl(\omega^{T,\Qdo}(H,\mu^{\pQ},s,t)-Y(T,H,\mu^{\pQ},s,t)\bigr)\ptf
\]
Fixons $s$ et $t$. Remarquons que, dans l'expression $Y(T,H,\mu^{\pQ},s,t)$,
on peut sortir le $\ee^H$ des séries d'Eisenstein et on obtient
{\multlinegap0pt\begin{multline*}
Y(T,H,\mu^{\pQ},s,t)=\ee^{\langle s\theto (\mu^{\pQ})-t(\mu^{\pQ}),H\rangle}\\
\smash[b]{{}\times \int_{L_0(F)\backslash L_0(\adef)^1}\int_\K}
\tronc^{T[H^Q],\Qdo}E_{\Qdo}\bigl(xk,\Mint\bigl(s,\theto (\mu^{\pQ})\bigr)\WPhi,\theto(\mu^{\pQ})\bigr)\\
{}\times \widebar{E_{\Qdo}(xk,\Mint(t,\mu^{\pQ})\WPsi,\mu^{\pQ})}\dd k\dd x\ptf
\end{multline*}}
On sait (\cf théorème~\ref{prodscalvar}) qu'il existe un réel $R>0$ tel que, pour $\mu^{\pQ}$ dans le support de
$B_{2}$, on ait la majoration
\[
\lvert \omega^{T,\Qdo}(H,\mu^{\pQ},s,t)-Y(T,H,\mu^{\pQ},s,t)\rvert \ll 
\sup_{\alpha\in \Delta_0^{\Qdo}} \ee^{-R\alpha(T[H^Q])}\ptf
\]
Comme on l'a dit plusieurs fois, ce dernier terme se majore par $\ee^{-R\dPO(T)}$.
D'où une majoration
\begin{multline*}
\lvert A^T_{\mathrm{pure}}(B,H)-A^T_{\mathrm{unit}}(B,H)\rvert\\
 \ll 
\ee^{-R\dPO(T)}\kappa^{\eta T}(H^Q-T_{\Qdo}^Q)
\tsQR(H-T)\phi_{\Qdo}^Q(H-T)\widehat{B}_1(\theta_0^{-1}(H_Q)-H_{\pQ}).
\end{multline*}
On continue le calcul comme précédemment et on obtient
\[
\lvert A^T_{\mathrm{pure}}(B)-A^T_{\mathrm{unit}}(B)\rvert \ll \ee^{-R\dPO(T)}\ptf
\]
La dernière assertion de l'énoncé résulte maintenant du lemme~\ref{2.3}(ii). 
\end{proof}

\section[Décomposition de $A^T_{\mathrm{pure}}(B)$]{\mathversion{bold}Décomposition de $A^T_{\mathrm{pure}}(B)$}\label{Apure}

Le terme $A^T_{\mathrm{pure}}(B)$ est une intégrale en $H\in \ga_{\Qdo}^G$.
Changeons de variable en remplaçant $H$
par $T_{\Qdo}+H-Y$, où le nouvel $H$ appartient à
$\ga_Q^G$ et $Y\in \ga_{\Qdo}^Q$. La condition
\[
\kappa^{\eta T}(H^Q-T_{\Qdo}^Q)\tsQR(H-T)\phi_{\Qdo}^Q(H-T)=1
\]
devient
\[
\kappa^{\eta T}(Y)\tsQR(H)\phi_{\Qdo}^Q(-Y)=1\ptf
\]
Puisque $\phi_{\Qdo}^Q(-Y)=1$, on peut écrire
\[
Y=\sum_{\alpha\in \Delta_0^Q-\Delta_0^{\Qdo}}y_{\alpha}(\alpha^\vee)_{\Qdo}
\]
avec des $y_{\alpha}\geq0$. L'élément $H+T[H^Q]$ du paragraphe précédent devient
\[
T_{\Qdo}+H+Y+T[T_{\Qdo}^Q-Y]\ptf
\]
On a calculé $T[T_{\Qdo}^Q-Y]$ dans le lemme   \ref{raff}: ce terme vaut
\[
T^{Q_0}-\sum_{\alpha\in \Delta_0^Q-\Delta_0^{\Qdo}}y_{\alpha}(\alpha^\vee)^{\Qdo}\ptf
\]
Alors
\[
T_{\Qdo}+H+Y+T[T_{\Qdo}^Q+Y]=T+H-X
\]
où
\[
X=\sum_{\alpha\in \Delta_0^Q-\Delta_0^{\Qdo}}y_{\alpha}\alpha^\vee\ptf
\]
On peut encore changer de variables en remplaçant $Y$ par $X$. Ce $X$ appartient au
cône engendré par les éléments de
$\Delta_0^Q-\Delta_0^{\Qdo}$. On note $\mathcal{C}(Q,\Qdo)$ ce cône.
L'ancien $Y$ devient $X_{\Qdo}$. On a ainsi transformé notre intégrale sur $\ga_{\Qdo}^G$ en une intégrale sur
$\ga_Q^G\times \mathcal{C}(Q,\Qdo)$, l'ancienne fonction
\[
\kappa^{\eta T}(H^Q-T_{\Qdo}^Q)\tsQR(H-T)\phi_{\Qdo}^Q(H-T)
\]
devenant $\kappa^{\eta T}(X_{\Qdo})\tsQR(H)$ et le terme $H+T[Y^Q]$ devenant $T+H-X$.
On observera que la mesure en $X$ n'est plus la mesure euclidienne mais la transportée de la mesure euclidienne sur
$\ga_{\Qdo}^Q$ par la projection injective $X\mapsto X_{Q_0}$.
Revenons à la définition du terme $\omega^{T,\Qdo}(H,\tmu,\nu)$ du paragraphe précédent,
que l'on note maintenant $\omega^{T,\Qdo}(H,X,\tmu,\nu)$, avec nos nouvelles
variables $H$ et $X$. L'application
\[
\bigcup_{t\in \weyl^{\pQ}(\ga_S,\Qdo)}\weyl^{\Qdo}(t(\ga_S),\ga_{\Sp})\to\weyl^{\pQ}(\ga_S,\ga_{\Sp})
\]
qui, à $t'\in \weyl^{\Qdo}(t(\ga_S),\ga_{\Sp})$ associe $t't$, est bijective.
On remplace la variable d'origine $t$ par un tel produit $t't$. On remplace ensuite $s$ par $(t')^{-1}s$.
On obtient
{\multlinegap0pt\begin{multline*}
\omega^{T,\Qdo}(H,X,\tmu,\nu)\\
=\sum_{t\in \weyl^{\pQ}(\ga_S,\Qdo)}
\sum_{s\in \weyl^Q\left(\theto(\ga_S),t(\ga_S)\right)}
\sum_{\{\Sp\mid  \PO\subset \Sp\subset \Qdo\}}\,
\sum_{t'\in\weyl^{\Qdo}(t(\ga_S),\ga_{\Sp})}\ee^{\langle t's\tmu-t't\nu,T+H-X\rangle}\\
{}\times\epsilon_{\Sp}^{\Qdo}(t's\tmu-t't\nu)
\langle \Mint(t't,\nu)^{-1}\Mint(t's,\tmu)\WPhi,\WPsi\rangle\ptf
\end{multline*}}%
Le parabolique $tS$ n'a pas de raison d'être standard mais il existe un unique parabolique standard
$\tdS\subset \Qdo$ 
\newindex{St@$\tdS$}{inftdS}%
qui a même Levi $\tdM$ que $tS$. On peut remplacer ci-dessus $\Sp$ par $\tdS$.
La double somme en $\Sp$ et $t'$ se transforme en une somme sur $\mathcal{P}^{\Qdo}(\tdM)$ et on
obtient (\cf lemme~\ref{sta}\eqref{eq5.3c}):
\[
\omega^{T,\Qdo}(H,X,\tmu,\nu)=\sum_{t\in \weyl^{\pQ}(\ga_S,\Qdo)}\sum_{s\in \weyl^\Q\left(\theto(\ga_S),
t(\ga_S)\right)} \omega_{s,t}^{T,\Qdo}(H,X,\tmu,\nu),
\]
où
{\multlinegap0pt\begin{multline*}
\omega_{s,t}^{T,\Qdo}(H,X,\tmu,\nu)=\sum_{\Sp\in \mathcal{P}^{\Qdo}(\tdM)}\ee^{\langle s\tmu-t\nu,H+Y_{\Sp}(T-X)\rangle}\epsilon{\Qdo}_{\Sp}(s\tmu-t\nu)\\
{}\times\langle \Mint(t,\nu)^{-1}\Mint_{\Sp\rest \tdS}(t\nu)^{-1}\Mint_{\Sp\rest \tdS}(s\tmu)
\Mint(s,\tmu)\WPhi,\WPsi\rangle\ptf
\end{multline*}}

On doit intégrer
\[
\omega^{T,\Qdo}(H,X,\mu)=\omega^{T,\Qdo}(H,X,\theto(\mu),\mu)
\]
en $\mu$, $X$ et $H$.
Chaque $\omega_{s,t}^{T,\Qdo}(H,X,\tmu,\nu)$ est encore une fonction $C^{\infty}$ de $\tmu$,
$\nu$. On note $\omega_{s,t}^{T,\Qdo}(H,X,\mu)$ sa valeur en $\tmu=\theto(\mu)$, $\nu=\mu$. En reprenant les preuves du
paragraphe précédent, on voit que l'on peut sortir les sommes en
$s$ et $t$ des intégrales. On obtient
\begin{equation}
A_{\mathrm{pure}}^T(B)=\sum_{t\in \weyl^{\pQ}(\ga_S,\Qdo)}\sum_{s\in\weyl^Q\left(\theto(\ga_S),t(\ga_S)\right)}A_{s,t}^T(B),\label{eq13.1C}
\end{equation}
où
\[
A_{s,t}^T(B)= \int_{\ga_Q^G}\tsQR(H)\int_{\mathcal{C}(Q,\Qdo)}\kappa^{\eta T}(X_{\Qdo})
\int_{\ima(\ga_S^\G)^*}\omega_{s,t}^{T,\Qdo}(H,X,\mu)B(\mu)\dd \mu\dd X\dd H\ptf
\]
\newnot{AstT@$A_{s,t}^T$}{astt}%

On fixe $t\in \weyl^{\pQ}(\ga_S,\Qdo)$ et $s\in \weyl^Q\bigl(\theto(\ga_S),t(\ga_S)\bigr)$. On a défini
\[
\omega_{s,t}^{T,\Qdo}(H,X,\tmu,\nu)
\]
par une sommation sur
$\Sp\in \mathcal{P}^{\Qdo}(\tdM)$. Élargissons cette sommation en la faisant porter sur $\Sp\in \mathcal{P}^Q(\tdM)$.
Précisément, pour $H\in \ga_Q^G$, $\mu,\nu\in\ima(\ga_S^G)^*$ et
$\tmu=\theto\mu$, posons 
\newindex{omegastTQ@$\omega^{T,Q}_{s,t}$}{omegatqst}%
{\multlinegap0pt\begin{multline*}
\omega^{T,Q}_{s,t}(H,\tmu,\nu)=\smash[b]{\sum_{\Sp\in \mathcal{P}^Q(\tdM)}}\ee^{\langle s\tmu-t\nu,H+Y_{\Sp}(T)\rangle}
\epsilon^Q_{\Sp}(s\tmu-t\nu)\\
{}\times\langle \Mint(t,\nu)^{-1}\Mint_{\Sp\rest \tdS}(t\nu)^{-1}\Mint_{\Sp\rest \tdS}(s\tmu)\Mint(s,\tmu)\WPhi,\WPsi\rangle
\ptf
\end{multline*}}
C'est une fonction $C^{\infty}$ de $\mu$ et $\nu$. On note $\omega_{s,t}^{T,Q}(H,\mu)$
sa valeur en $\nu=\mu$.

\begin{proposition}\label{W2.5}
\begin{enumerate}[(i)]
\item L'intégrale itérée
\[
\mathbf{A}_{s,t}^T(B)=\int_{\ga_Q^G}\tsQR(H)\biggl(\int_{\ima(\ga_S^\G)^*}
\omega_{s,t}^{T,Q}(H,\mu)B(\mu)\dd \mu\biggr)\dd H
\]
\newnot{AstT@$\mathbf{A}_{s,t}^T$}{bfastt}%
est convergente dans l'ordre indiqué.
\item Pour tout réel $r$\textup, on a une majoration
\[
\lvert A^T_{s,t}(B)-\mathbf{A}_{s,t}^T(B)\rvert \ll \dPO(T)^{-r}\ptf
\]
\end{enumerate}
\end{proposition}

\begin{proof}
Les paragraphes~\ref{Ws2.6} et~\ref{Ws2.8} seront consacrés à la preuve de cette proposition.
\end{proof}

\section{Majoration de transformées de Fourier}\label{Ws2.6}

Considérons l'application linéaire
\[
s\theto-t\colon (\ga_S^\G)^*\to t(\ga_S^G)^{*}\ptf
\]
On note $\gb_S^*$ son noyau, $\gc_{\tdS}^*$ son image, $\gc_S^*$ l'orthogonal de
$\gb_{S}^*$ dans ${(\ga_S^\G)^*}$ et $\gb_{\tdS}^*$ l'orthogonal de $\gc_{\tdS}^*$ dans $t(\ga_S^G)^{*}$.
L'application ci-dessus se restreint en un
isomorphisme
\[
\iota\colon \gc_S^*\to \gc_{\tdS}^*\ptf
\]
Pour $\Lambda\in t(\ga_S^G)^{*}$, on note $\Lambda_{b}$ et $\Lambda_{c}$ ses projections sur
$\gb_{\tdS}^*$ et $\gc_{\tdS}^*$. On définit une application linéaire
\[
\begin{aligned}
\gb_S^*\times t(\ga_S^G)^{*}&\to{(\ga_S^\G)^*}\times\gb_{\tdS}^*\\ 
(\lambda,\Lambda)&\mapsto(\mu(\lambda,\Lambda),
\Lambda_{b})
\end{aligned}
\]
par $\mu(\lambda,\Lambda)=\lambda+\iota^{-1}(\Lambda_{c})$. C'est un isomorphisme. On fixe une fonction
$B_{b}\in C_{c}^{\infty}(\ima\gb_{\tdS}^*)$ telle que $B_{b}(0)=1$.
Rappelons que l'on note $\tdM$ et $L$ les Levi standard de $\tdS$ et $Q$. Pour
$\lambda\in\ima\gb_S^*$, $\Lambda\in\ima(\ga_{\tdM}^G)^*$ et $\Sp\in \mathcal{P}^Q(\tdM)$,
posons
{\multlinegap0pt
\begin{multline*}
\varphi(\lambda;\Lambda,\Sp)=B\bigl(\mu(\lambda,\Lambda)\bigr)B_{b}(\Lambda_{b})
\bigl\langle\Mint\bigl(t,\mu(\lambda,\Lambda)\bigr)^{-1}\Mint_{\Sp\rest \tdS}\bigl(t\mu(\lambda,\Lambda)\bigr)^{-1}\\
{}\times\Mint_{\Sp\rest \tdS}(t\mu(\lambda,\Lambda)+\Lambda)
\Mint(s,s^{-1}t\mu(\lambda,\Lambda)+s^{-1}\Lambda)\WPhi,\WPsi\bigr\rangle \ptf
\end{multline*}}%
Ces fonctions sont $C^{\infty}$ et à support compact en $\lambda$ et $\Lambda$.
Considérées comme des fonctions de $\Lambda$, dépendant d'un paramètre $\lambda$, elles
forment une $(L,\tdM)$\hyph famille. Pour un élément $Z\in \ga_0^G$, on définit une $(L,\tdM)$\hyph famille familière par
\[
d^Z(\Lambda,\Sp)=\ee^{\langle \Lambda,Y_{\Sp}(Z)\rangle}\ptf
\]
On pose $\varphi^Z(\lambda;\Lambda,\Sp)=\varphi(\lambda;\Lambda,\Sp)d^Z(\Lambda,\Sp)$.
En se limitant aux $\Sp\subset \Qdo$,
on obtient des $(L_0,\tdM)$\hyph familles. Conformément aux définitions de la section~\ref{gmfam}, on définit les fonctions
\[
\varphi^{Z,\Qdo}_{\tdM}(\lambda;\Lambda)=
\sum_{\Sp\in \mathcal{P}^{\Qdo}(\tdM)}\varphi^Z(\lambda;\Lambda,\Sp)\epsilon_{\Sp}^{\Qdo}(\Lambda)
\]
et
\[
\varphi^{Z,Q}_{\tdM}(\lambda;\Lambda)=\sum_{\Sp\in \mathcal{P}^Q(\tdM)}\varphi^Z(\lambda;\Lambda,\Sp)
\epsilon_{\Sp}^Q(\Lambda)
\]
qui sont encore $C^{\infty}$ et à supports compacts en $\lambda$ et $\Lambda$.

\begin{lemme}\label{W2.6}
Pour $\lambda\in\ima\gb_S^*$\textup, $\Lambda=\Lambda_{c}\in\ima\gc_{\tdS}^*$\textup,
$H\in \ga_Q^G$ et $X\in \mathcal{C}(Q,\Qdo)$\textup, on a les égalités
\begin{equation}
\omega^{T,\Qdo}\bigl(H,X,\mu(\lambda,\Lambda)\bigr)B\bigl(\mu(\lambda,\Lambda)\bigr)=
\varphi^{T+H-X,\Qdo}_{\tdM}(\lambda;\Lambda)\tag{i}\label{eq13.i}
\end{equation}
et
\begin{equation}
\omega^{T,Q}\bigl(H,\mu(\lambda,\Lambda)\bigr)B\bigl(\mu(\lambda,\Lambda)\bigr)=
\varphi^{T+H,Q}_{\tdM}(\lambda;\Lambda)\ptf\tag{ii}\label{eq13.ii}
\end{equation}
\end{lemme}

\begin{proof} 
Pour $\Lambda\in \ima t(\ga_S^G)^{*}$ en position générale, il résulte des définitions que
\[
\varphi^{T+H-X,\Qdo}_{\tdM}(\lambda;\Lambda)=
\omega^{T,\Qdo}\bigl(H,X,\mu(\lambda,\Lambda),\nu(\lambda,\Lambda)\bigr)
B\bigl(\mu(\lambda,\Lambda)\bigr)B_{b}(\Lambda_{b}),
\]
où
\[
\nu(\lambda,\Lambda)=\theta_0^{-1}s^{-1}t\mu(\lambda,\Lambda)+\theta_0^{-1}s^{-1}\Lambda\ptf
\]
Il résulte de la définition de $\mu(\lambda,\Lambda)$ que
\[
\nu(\lambda,\Lambda)=\mu(\lambda,\Lambda)+\theta_0^{-1}s^{-1}\Lambda_{b}\ptf
\]
Pour $\lambda$ et $\Lambda_{c}$ fixés (donc aussi $\mu(\lambda,\Lambda)$),
il reste à faire tendre $\Lambda_{b}$ vers $0$ pour
obtenir l'égalité~\eqref{eq13.i}. La preuve de~\eqref{eq13.ii} est similaire.
\end{proof}

En conséquence, la définition de $A_{s,t}^T(B)$ se récrit
{\multlinegap0pt
\begin{multline*}
A_{s,t}^T(B)
=\lvert \det(\iota)\rvert ^{-1}\int_{\ga_Q^G}\tsQR(H)\int_{\mathcal{C}(Q,\Qdo)}
\kappa^{\eta T}(X_{\Qdo})\\
{}\times\int_{\ima\gb_S^*}
\int_{\ima\gc_{\tdS}^*}\mspace{-5mu}\varphi^{T+H-X,\Qdo}_{\tdM}(\lambda;\Lambda)\dd \Lambda\dd \lambda\dd X\dd H\ptf
\end{multline*}}%
On introduit les transformées de Fourier inverses en $\Lambda$ des fonctions
\[
\varphi(\lambda;\Lambda,\Sp),\qquad \varphi^Z(\lambda;\Lambda,\Sp),
\qquad \varphi_{\tdM}^{Z}(\lambda;\Lambda)\Qquad{et}
\varphi_{\tdM}^{Z,\Qdo}(\lambda,\Lambda)
\]
que l'on note
\[
\hvf(\lambda;U,\Sp),\qquad \hvf^Z(\lambda;U,\Sp),\qquad \hat
{\varphi}_{\tdM}^{Z,Q}(\lambda;U)\Qquad{et}\hvf_{\tdM}^{Z,\Qdo}(\lambda;U)
\]
le paramètre $U$ appartenant à $t(\ga_S)^G$. Ce sont des transformées de Fourier
inverses de fonctions $C^{\infty}$ à support compact, donc
des fonctions de Schwartz en $U$. Par inversion de Fourier,
l'intégrale intérieure de l'expression ci-dessus est égale à
\[
\int_{\gb_{\tdS}}\hvf^{T+H-X,\Qdo}_{\tdM}(\lambda;U)\dd U,
\]
où $\gb_{\tdS}$ est l'annulateur de $\gc_{\tdS}^*$ dans $t(\ga_S)^{G}$. D'où le

\begin{corollaire}\label{W2.6(1)}
\begin{multline*} 
A_{s,t}^T(B)=\lvert \det(\iota)\rvert ^{-1}\int_{\ga_Q^G}\tsQR(H)\int_{\mathcal{C}(Q,\Qdo)} 
\kappa^{\eta T}(X_{\Qdo})\\
{}\times\int_{\ima\gb_S^*}
\int_{\gb_{\tdS}}\hvf^{T+H-X,\Qdo}_{\tdM}(\lambda;U)\dd U\dd \lambda\dd X\dd H\ptf
\end{multline*}
\end{corollaire}

\begin{lemme}\label{W2.7}
Fixons un réel $\rho>0$. Considérons les cinq expressions
\begin{gather}
\int_{\ga_Q^G}\tsQR(H)\int_{\mathcal{C}(Q,\Qdo)}
\int_{\ima\gb_S^*}\int_{\gb_{\tdS}}
\lvert \hvf^{T+H-X,\Qdo}_{\tdM}(\lambda;U)\rvert \dd U\dd \lambda\dd X\dd H;\label{eq13.1d}\\
\int_{\ga_Q^G}\tsQR(H)\int_{\mathcal{C}(Q,\Qdo)}\mspace{-10mu}\bigl(1- \kappa^{\eta T}(X_{\Qdo})\bigr)
\int_{\ima\gb_S^*}\int_{\gb_{\tdS}}\lvert\hvf^{T+H-X,\Qdo}_{\tdM}(\lambda;U)\rvert \dd U\dd \lambda\dd X\dd H;\label{eq13.2d}\\
\int_{\ga_Q^G}\tsQR(H)\int_{\mathcal{C}(Q,\Qdo)} \int_{\ima\gb_S^*}
\int_{\gb_{\tdS}}\bigl(1-\kappa^{\rho T}(U)\bigr)\lvert
\hvf^{T+H-X,\Qdo}_{\tdM}(\lambda;U)\rvert \dd U\dd \lambda\dd X\dd H;\label{eq13.3d}\\
\int_{\ga_Q^G}\tsQR(H) \int_{\ima\gb_S^*}\int_{\gb_{\tdS}}\lvert
\hvf^{T+H,Q}_{\tdM}(\lambda;U)\rvert \dd U\dd \lambda\dd H;\label{eq13.4d}\\
\int_{\ga_Q^G}\tsQR(H) \int_{\ima\gb_S^*}
\int_{\gb_{\tdS}}\bigl(1-\kappa^{\rho T}(U)\bigr)\lvert \hvf^{T+H,Q}_{\tdM}(\lambda;U)\rvert
\dd U\dd \lambda\dd H\ptf\label{eq13.5d}
\end{gather}
Alors
\begin{enumerate}[(i)]
\item Les cinq expressions sont convergentes.
\item Pour tout réel $r$\textup,~\eqref{eq13.2d} est essentiellement majorée par $\dPO(T)^{-r}$.
\item Il existe une constante absolue $\rho_0>0$ telle que\textup, si $\rho>\rho_0$\textup,
\eqref{eq13.3d} et~\eqref{eq13.5d} sont essentiellement majorées pour
tout réel $r$ par $\dPO(T)^{-r}$.
\end{enumerate}
\end{lemme}

\begin{proof} 
On traite d'abord l'expression~\eqref{eq13.4d}. Par définition, la $(L,\tdM)$\hyph famille
\[
\bigl(\varphi^Z(\lambda;\Lambda,\Sp)\bigr)_{\Sp\in\mathcal{P}^Q(\tdM)}
\]
est le produit des $(L,\tdM)$\hyph familles
\[
\bigl(\varphi(\lambda;\Lambda,\Sp)\bigr)_{\Sp\in\mathcal{P}^Q(\tdM)}
\]
et $\bigl(d^Z(\Lambda,\Sp)\bigr)_{\Sp\in \mathcal{P}^Q(\tdM)}$.
D'après le lemme~\ref{decompGM}, on a
\[
\varphi^{Z,Q}_{\tdM}(\lambda;\Lambda)=\sum_{P'\in \mathcal{F}^Q(\tdM)}
d^{Z,P'}_{\tdM}(\Lambda)\varphi^Q_{P'}(\lambda;\Lambda_{P'})\ptf
\]
Ou encore, grâce à l'équation~\eqref{eq1.2f} du lemme~\ref{gmorth}:
\[
\varphi^{Z,Q}_{\tdM}(\lambda;\Lambda)=\sum_{P'\in \mathcal{F}^Q(\tdM)}\varphi^Q_{P'}(\lambda;\Lambda_{P'})
\int_{\ga_{\tdM}^{P'}}\ee^{\langle\Lambda,U^{P'}+Y_{P'}(Z)\rangle}
\Gamma_{\tdM}^{P'}\bigl(U^{P'},\YY(Z)\bigr)\dd U^{P'},
\]
où
\[
\YY(Z)=\bigl(Y_{\Sp}(Z)\bigr)_{\Sp\in \mathcal{P}^Q(\tdM)}
\]
est la $(L,\tdM)$\hyph famille orthogonale associée à $Z$. La fonction
$\varphi^Q_{P'}(\lambda;\Lambda_{P'})$ est $C^{\infty}$ et à
support compact en $\lambda$ et $\Lambda_{P'}$. En introduisant sa transformée de Fourier inverse
$\hvf^Q_{P'}(\lambda;U_{P'})$, l'expression ci-dessus devient
\begin{align*}
\varphi^{Z,Q}_{\tdM}(\lambda;\Lambda)&=\sum_{P'\in \mathcal{F}^Q(\tdM)}\int_{\ga_{\tdM}^{G}}
\ee^{\langle\Lambda,U+Y_{P'}(Z)\rangle}\Gamma_{\tdM}^{P'}\bigl(U^{P'},\YY(Z)\bigr)
\hvf_{P'}^Q(\lambda;U_{P'})\dd U\\
&=\sum_{P'\in \mathcal{F}^Q(\tdM)}\int_{\ga_{\tdM}^{G}}\ee^{\langle \Lambda,U\rangle}\Gamma_{\tdM}^{P'}\bigl(U^{P'},\YY(Z)\bigr)\hat
{\varphi}_{P'}^Q\bigl(\lambda;U_{P'}-Y_{P'}(Z)\bigr)\dd U\ptf
\end{align*}
Cette formule détermine
\[
\hvf^{Z,Q}_{\tdM}(\lambda;U)=
\sum_{P'\in \mathcal{F}^Q(\tdM)}\Gamma_{\tdM}^{P'}\bigl(U^{P'},\YY(Z)\bigr)\hvf_{P'}^Q\bigl(\lambda;U_{P'}-Y_{P'}(Z)\bigr)\ptf
\]
L'expression~\eqref{eq13.4d} est donc majorée par
\[
\sum_{P'\in \mathcal{F}^Q(\tdM)}I^T_{\eqref{eq13.4d}}(P'),
\]
où
\begin{multline*}
I^T_{\eqref{eq13.4d}}(P')
=\smash[b]{\int_{\ga_Q^G}\tsQR(H)\int_{\ima\gb_S^*}\int_{\gb_{\tdS}}}
\Gamma_{\tdM}^{P'}\bigr(U^{P'},\YY(T)\bigr)\\
{}\times\lvert
\hvf_{P'}^Q(\lambda;U_{P'}-Y_{P'}(T)-H)\rvert \dd U\dd \lambda\dd H\ptf
\end{multline*}
On a simplifié $\YY(T+H)$ en $\YY(T)$ et $Y_{P'}(T+H)$ en $Y_{P'}(T)+H$, ainsi qu'il est loisible.
Fixons $P'\in \mathcal{F}^Q(\tdM)$. Posons $\mathfrak{d}=\gb_{\tdS}\cap \ga_{P'}^G$, notons $\mathfrak{e}$
l'orthogonal de cet espace dans $\ga_{P'}^G$ et
posons
\[
\gb_{\sharp} =\gb_{\tdS}\cap(t(\ga_S)^{P'}\oplus \mathfrak{e})\ptf
\]
On a l'égalité
\[
\gb_{\tdS}=\mathfrak{d}\oplus \gb_{\sharp}
\]
et la projection sur $t(\ga_S)^{P'}$ est injective sur $\gb_{\sharp}$. La fonction
$\hvf_{P'}^Q(\lambda;V)$ est à support compact en $\lambda$ et de Schwartz
en $V\in \ga_{P'}^G$. On peut fixer une fonction $C^{\infty}$ et à support compact $C$ sur $\ima\gb_S^*$ et
des fonctions de Schwartz $\phi_{d}$ sur $\mathfrak{d}$ et
$\phi_{e}$ sur $\mathfrak{e}$, toutes à valeurs positives ou nulles, de sorte que
\[
\lvert \hvf_{P'}^Q(\lambda;V)\rvert \leq C(\lambda)\phi_{d}(V_{d})\phi_{e}(V_{e}),
\]
où $V_{d}$ et $V_{e}$ sont les projections orthogonales sur $\mathfrak{d}$ et $\mathfrak{e}$.
On a alors la majoration
{\multlinegap0pt\begin{multline}\label{eq13.6d}
I^T_{\eqref{eq13.4d}}(P')\leq  \smash[b]{\int_{\ga_Q^G}\tsQR(H)\int_{\ima\gb_S^*}
\int_{\mathfrak{d}}\int_{\gb_{\sharp}}}C(\lambda)\Gamma_{\tdM}^{P'}\bigl(U^{P'},
\YY(T)\bigr)\phi_{d}(V-Y_{P',d}(T)-H_{d})\\
{}\times\phi_{e}(U_{P',e}-Y_{P',e}(T)-H_{e})\dd U\dd V\dd \lambda\dd H
\end{multline}}
où on a abrégé les notations, par exemple $Y_{P',e}(T)=\bigl(Y_{P'}(T)\bigr)_{e}$.
Par le changement de variable
\[
V\mapsto V+Y_{P',d}(T)+H_{d}
\]
et parce que $C$ et $\phi_{d}$
sont intégrables, on obtient
\begin{equation}
I^T_{\eqref{eq13.4d}}(P')\ll \int_{\ga_Q^G}\tsQR(H)\int_{\gb_{\sharp}}
\Gamma_{\tdM}^{P'}\bigl(U^{P'},\YY(T)\bigr)\phi_{e}(U_{P',e}-Y_{P',e}(T)-H_{e})\dd U\dd H,
\label{eq13.7d}
\end{equation}
ou encore
\begin{equation}
I^T_{\eqref{eq13.4d}}(P')\ll \int_{\ga_Q^G}\tsQR(H)\phi_{e}^T(-Y_{P',e}(T)-H_{e})\dd H,\label{eq13.8d}
\end{equation}
où, pour $V\in \mathfrak{e}$,
\[
\phi_{e}^T(V)=\int_{\gb_{\sharp}}\Gamma_{\tdM}^{P'}\bigl(U^{P'},\YY(T)\bigr)\phi_{e}(U_{P',e}+V)\dd U\ptf
\]
Puisque la projection $U\mapsto U^{P'}$ est injective sur $\gb_{\sharp}$, cette dernière intégrale est à support
compact. Donc $\phi_{e}^T$ est encore de Schwartz. Pour démontrer
que l'intégrale de droite de~\eqref{eq13.8d} est convergente, il suffit de prouver la majoration:
\begin{equation}
\lVert H\rVert \ll \lVert H_{e}\rVert \qquad\text{pour tout $H\in \ga_Q^G$ tel que $\tsQR(H)=1$}
\ptf\label{eq13.9d}
\end{equation}
On rappelle que l'on a noté $\gk$ le noyau de l'application $q$. L'espace $\gb_{\tdS}$
est l'annulateur dans $t(\ga_S)^G$ de l'image par $s\theto-t$ de $\ga_S^G
$. C'est donc le noyau de $\theta_0^{-1}s^{-1}-t^{-1}$ dans
$t(\ga_S)^G$. On a
\[
\theta_0^{-1}s^{-1}-t^{-1}=\theta_0^{-1}s^{-1}(1-w\theto)
\]
où $w=s\theto(t)^{-1}$,
donc $\gb_{\tdS}$ est le noyau de $1-w\theto$ dans
$t(\ga_S)^G$. Mais $w$ fixe $\ga_Q^G$. Donc
\[
(X-w\theto X)_Q=q(X)
\]
pour tout $X$ et le noyau de $1-w\theto$ est contenu dans $\gk$. A fortiori
\begin{equation}
\gb_{\tdS}\subset \gk\ptf\label{eq13.10d}
\end{equation}
Il en résulte une majoration
\[
\lVert q(H)\rVert \ll \lVert H_{e}\rVert \ptf
\]
On utilise le lemme~\ref{recroissb}: on a une majoration $\lVert H\rVert \ll \lVert q(H)\rVert $.
La majoration~\eqref{eq13.9d} en résulte.
Cela achève la preuve de la convergence de~\eqref{eq13.4d}.
Considérons l'expression~\eqref{eq13.5d}. On voit comme ci-dessus qu'elle est majorée par
\[
\sum_{P'\in \mathcal{F}^Q(\tdM)}I^T_{\eqref{eq13.5d}}(P'),
\]
avec
{\multlinegap0pt\begin{multline*}
I^T_{\eqref{eq13.5d}}(P')\ll \int_{\ga_Q^G}\tsQR(H)\int_{\mathfrak{d}}
\int_{\gb_{\sharp}}\bigl(1-\kappa^{\rho T}(U+V)\bigr)\Gamma_{\tdM}^{P'}\bigl(U^{P'},\YY(T)\bigr)\\
{}\times\phi_{d}(V-Y_{P',d}(T)-H_{d})\phi_{e}(U_{P',e}-Y_{P',e}(T)-H_{e})\dd U\dd V\dd H
\end{multline*}}%
(\cf~\eqref{eq13.6d} ci-dessus, où on s'est débarrassé de l'intégrale en $\lambda$).
Fixons $\rho'>0$ que l'on précisera plus tard. On
majore $I^T_{\eqref{eq13.5d}}(P') $ par $I^T_{\eqref{eq13.5d},\geq}(P')+I^T_{\eqref{eq13.5d},<}(P')$,
ces termes étant définis par l'intégrale précédente où l'on glisse la fonction $1-\kappa^{\rho' T}(H)$
dans le premier terme et la fonction $\kappa^{\rho' T}(H)$ dans le second.
Majorons d'abord $I^T_{\eqref{eq13.5d},\geq}(P')$. On a une majoration similaire à~\eqref{eq13.7d}:
{\multlinegap0pt\begin{multline}\label{eq13.11d}
I^T_{\eqref{eq13.5d},\geq}(P')\ll \int_{\ga_Q^G}\bigl(1-\kappa^{\rho' T}(H)\bigr)\tsQR(H)
\int_{\gb_{\sharp}}\Gamma_{\tdM}^{P'}\bigl(U^{P'},\YY(T)\bigr)\\
{}\times\phi_{e}(U_{P',e}-Y_{P',e}(T)-H_{e})\dd U\dd H
\end{multline}}
où la constante implicite ne dépend pas de $T$. Les constantes $c_1$, $c_{2}$, etc.,
que l'on va introduire sont absolues,
c'est-à-dire ne dépendent d'aucune des variables $T$, $H$,
etc. Elles ne dépendent pas non plus de $Q$, $R$, etc.,
simplement parce que ces dernières données ne parcourent que
des ensembles finis. Il existe $c_1>0$ tel que la condition
\[
\Gamma_{\tdM}^{P'}\bigl(U^{P'},\YY(T)\bigr)=1
\]
entraîne
\[
\lVert U^{P'}\rVert \leq c_1\lVert T\rVert \ptf
\]
Puisque l'application $U\mapsto U^{P'}$ est injective, il existe
$c_{2}>0$ tel que cette même condition entraîne $\lVert U\rVert \leq c_{2}\lVert T\rVert$, donc aussi
$\lVert U_{P',e}\rVert \leq c_{2}\lVert T\rVert$. Il existe
$c_{3}>0$ tel que cette relation entraîne
\[
\lVert U_{P',e}+Y_{P',e}(T)\rVert \leq c_{3}\lVert T\rVert \ptf
\]
D'après~\eqref{eq13.9d}, il existe $c_{4}>0$ tel que
$\lVert H_{e}\rVert \geq c_{4}\lVert H\rVert$. La condition $1-\kappa^{\rho' T}(H)=1$ entraîne
$\lVert H_{e}\rVert \geq c_{4}\rho'\lVert T\rVert$. On fixe $\rho'$ tel que
$c_{4}\rho'> c_{3}$. Alors les deux conditions ensemble entraînent l' inégalité
\[
\lVert U_{P',e}-Y_{P',e}(T)-H_{e}\rVert \geq \biggl(c_{4}-\frac{c_{3}}{\rho'}\biggr)\lVert H\rVert \ptf
\]
Puisque $\phi_{e}$ est de Schwartz, l'expression~\eqref{eq13.11d} est donc essentiellement majorée par
\[
\int_{\ga_Q^G}\bigl(1-\kappa^{\rho' T}(H)\bigr)\int_{\gb_{\sharp}}
\Gamma_{\tdM}^{P'}\bigl(U^{P'},\YY(T)\bigr)\lVert H\rVert ^{-r}\dd U\dd H
\]
pour n'importe quel réel $r$. L'intégrale en $U$ est essentiellement majorée par $\lVert T\rVert ^D$ pour un entier $D$
convenable. L'intégrale en $H$ est essentiellement
majorée par $\lVert T\rVert ^{-r}$. On en déduit une majoration
\begin{equation}
I^T_{\eqref{eq13.5d},\geq}(P')\ll \dPO(T)^{-r}\label{eq13.12d}
\end{equation}
pour tout réel $r$.
Traitons maintenant $I^T_{\eqref{eq13.5d},<}(P')$. Rappelons que
{\multlinegap0pt\begin{multline*}
I^T_{\eqref{eq13.5d},<}(P')=\int_{\ga_Q^G}\tsQR(H)\kappa^{\rho' T}(H)\int_{\mathfrak{d}}
\int_{\gb_{\sharp}}\bigl(1-\kappa^{\rho T}(U+V)\bigr)
\Gamma_{\tdM}^{P'}\bigl(U^{P'},\YY(T)\bigr)
\\
{}\times\phi_{d}(V-Y_{P',d}(T)-H_{d})\phi_{e}(U_{P',e}-Y_{P',e}(T)-H_{e})\dd U\dd V\dd H\ptf
\end{multline*}}
On a encore $\lVert U\rVert \leq c_{2}\lVert T\rVert$. Ajoutons la condition
\[
\bigl(1-\kappa^{\rho T}(U+V)\bigr)=1
\]
c'est-à-dire $\lVert U+V\rVert \geq \rho\lVert T\rVert $. Si $\rho>c_{2}$, les deux conditions entraînent que
\[
\lVert V\rVert \geq (\rho-c_{2})\lVert T\rVert
\]
c'est-à-dire
\[
\bigl(1-\kappa^{(\rho-c_{2})T}(V)\bigr)=1\ptf
\]
\emph{A fortiori}
$V\ne 0$ et l'espace $\mathfrak{d}$ n'est pas nul.
On a fixé $\rho'$ ci-dessus. Il existe $c_{5}>0$ tel que $\kappa^{\rho' T}(H)=1$ entraîne
\[
\lVert Y_{P',d}(T)-H_{d}\rVert \leq c_{5}\lVert T\rVert \ptf
\]
Supposons $\rho> c_{2}+c_{5}$. Alors nos conditions entraînent
\[
\lVert V-Y_{P',d}-H_{d}\rVert \geq \biggl(1-\frac{c_{5}}{\rho-c_{2}}\biggr)\lVert V\rVert\ptf
\]
Puisque $\phi_{d}$ est de Schwartz, $I^T_{\eqref{eq13.5d},<}(P')$ est essentiellement majorée par
{\multlinegap0pt\begin{multline*}
\smash[b]{\int_{\ga_Q^G}\kappa^{\rho' T}(H)\int_{\mathfrak{d}}\int_{\gb_{\sharp}}}
\Gamma_{\tdM}^{P'}\bigl(U^{P'},\YY(T)\bigr)
\lVert V\rVert ^{-r}
\bigl(1-\kappa^{(\rho-c_{2})T}(V)\bigr)\\
\phi_{e}(U_{P',e}-Y_{P',e}(T)-H_{e})\dd U\dd V\dd H
\end{multline*}}
pour n'importe quel réel $r$. On majore encore $\phi_{e}$ par une constante. Les intégrales en $H$ et $U$ sont
essentiellement bornées par $\lVert T\rVert ^D$ pour un entier
$D$ convenable. L'intégrale en $V$ est essentiellement bornée par $\lVert T\rVert ^{-r}$ pour n'importe quel réel $r$.
D'où une majoration
\[
I^T_{\eqref{eq13.5d},<}(P')\ll \dPO(T)^{-r}
\]
pour tout réel $r$. Jointe à~\eqref{eq13.12d}, elle prouve l'assertion de l'énoncé concernant l'expression~\eqref{eq13.5d}.
Considérons maintenant l'expression~\eqref{eq13.1d}. On voit comme précédemment qu'elle est essentiellement majorée par
\[
\sum_{P'\in \mathcal{F}^{\Qdo}(\tdM)}I^T_{(1)}(P'),
\]
où
{\multlinegap0pt\begin{multline*}
I^T_{\eqref{eq13.1d}}(P')= \int_{\ga_Q^G}\tsQR(H)\int_{\mathcal{C}(Q,\Qdo)} \int_{\ima\gb_S^*}
\int_{\gb_{\tdS}}\bigl\lvert \Gamma_{\tdM}^{P'}\bigl(U^{P'},\YY(T-X)\bigr)\bigr\rvert
\\
{}\times\lvert \hvf^{\Qdo}_{\tdM}(\lambda;U_{P'}-Y_{P'}(T-X)-H)\rvert \dd U\dd \lambda\dd X\dd H\ptf
\end{multline*}}
On fixe $P'\in \mathcal{F}^{\Qdo}(\tdM)$.
On pose les mêmes définitions que dans la partie de la preuve consacrée à
l'expression~\eqref{eq13.4d}. On obtient une majoration similaire à~\eqref{eq13.7d}:
{\multlinegap0pt\begin{multline*}
I^T_{\eqref{eq13.1d}}(P')\ll \smash[b]{\int_{\ga_Q^G}\tsQR(H) \int_{\mathcal{C}(Q,\Qdo)}\int_{\gb_{\sharp}}}\bigl\lvert
\Gamma_{\tdM}^{P'}\bigl(U^{P'},\YY(T-X)\bigr)\bigr\rvert \\
{}\times\phi_{e}(U_{P',e}-Y_{P',e}(T-X)-H_{e})\dd U\dd X\dd H\ptf
\end{multline*}}
Décomposons $\mathfrak{e}$ en $\mathfrak{e}_{\natural}\oplus\mathfrak{e}_{\flat}$,
où $\mathfrak{e}_{\flat}$ est l'orthogonal de $\gb_{\tdS}$ dans $\ga_{\Qdo}^G$ et
$\mathfrak{e}_{\natural}$ est son orthogonal dans $\mathfrak{e}$ (parce que $P'\subset \Qdo$, $\mathfrak{e}_{\flat}$
est bien contenu dans $\mathfrak{e}$).On peut fixer des fonctions de
Schwartz $\phi_{\natural}$ sur $\mathfrak{e}_{\natural}$ et $\phi_{\flat}$ sur $\mathfrak{e}_{\flat}$
de sorte que, pour $V\in \mathfrak{e}$,
$\phi_{e}(V)\leq \phi_{\natural}(V_{\natural})\phi_{\flat}(V_{\flat})$, avec le même genre de notations que précédemment.
Alors
{\multlinegap0pt\begin{multline*}
I^T_{\eqref{eq13.1d}}(P')\ll \smash[b]{\int_{\ga_Q^G}\tsQR(H)\int_{\mathcal{C}(Q,\Qdo)}\int_{\gb_{\sharp}}}
\bigl\lvert \Gamma_{\tdM}^{P'}\bigl(U^{P'},\YY(T-X)\bigr)\bigr\rvert \phi_{\natural}(U_{\natural}-Y_{P',\natural}(T-X)-H_{\natural})\\
{}\times\phi_{\flat}(-T_{\flat}-H_{\flat}+X_{\flat})\dd U\dd X\dd H\ptf
\end{multline*}}%
On a utilisé que la projection de $Y_{P',e}(T-X)$ dans $\mathfrak{e}_{\flat}$ n'était autre que $T_{\flat}-X_{\flat}$,
ce qui résulte de l'inclusion $P'\subset \Qdo$.
On majore $\phi_{\natural}$ par une constante. Puisque $U\mapsto U^{P'}$ est injective, l'intégrale en $U$ est
essentiellement bornée par
\[
\lVert T\rVert ^D+(1+\lVert X\rVert)^{D}
\]
pour un $D$ convenable. Puisque $\phi_{\flat}$ est de Schwartz, l'expression ci-dessus sera convergente si l'on
montre que pour tout $T$, tout $H\in \ga_Q^G$ tel que $\tsQR(H)=1$ et tout $X\in \mathcal{C}(Q,\Qdo)$, on a une
majoration
\begin{equation}
\lVert H\rVert+\lVert X\rVert \ll \lVert T_{\flat}+ H_{\flat}- X_{\flat}\rVert \ptf
\label{eq13.13d}
\end{equation}
L'espace $\gk$ contient à la fois $\ga_0^{\Qdo}$ et $\gb_{\tdS}$ d'après~\eqref{eq13.10d}. Il
contient donc le noyau
$\ga_0^{P'}\oplus \mathfrak{d}\oplus \mathfrak{e}_{\natural}$
de la projection sur $\mathfrak{e}_{\flat}$. On a donc
\[
\lVert q(V)\rVert \ll \lVert V_{\flat}\rVert
\]
pour tout $V$. On applique cette relation à
\[
V=T+H-X_{\Qdo}\ptf
\]
Puisque $q(T)=0$, il reste à appliquer le lemme~\ref{recroissb} pour en déduire~\eqref{eq13.13d}.
Cela prouve la convergence de l'expression~\eqref{eq13.1d}. Plus précisément, on déduit de~\eqref{eq13.13d} une majoration
\[
I^T_{\eqref{eq13.1d}}(P')\ll \int_{\ga_Q^G}\int_{\mathcal{C}(Q,\Qdo)}\lVert T\rVert ^{D}(1+\lVert X\rVert)^{-r}
(1+\lVert H\rVert)^{-r}\dd X\dd H
\]
pour tout réel $r$. Considérons maintenant l'expression~\eqref{eq13.2d}. La seule différence est que le terme
$1-\kappa^{\eta T}(X_{Q_0})$ se glisse dans les calculs. La majoration ci-dessus est
alors remplacée par
\[
I^T_{\eqref{eq13.1d}}(P')\ll \int_{\ga_Q^G}\int_{\mathcal{C}(Q,\Qdo)}\lVert T\rVert ^{A}
\bigl(1-\kappa^{\eta T}(X_{\Qdo})\bigr)(1+\lVert X\rVert)^{-r}(1+\lVert H\rVert)^{-r}\dd X\dd H\ptf
\]
Cette expression est convergente et majorée par $\lVert T\rVert ^{-r}$ pour tout réel $r$. Cela prouve le~(ii) de
l'énoncé.
Considérons maintenant l'expression~\eqref{eq13.3d}. Remarquons que, pour démontrer l'assertion la concernant,
on peut aussi bien y glisser
\[
\kappa^{\eta T}(X_{\Qdo})\ptf
\]
En effet, en notant (\ref{eq13.3d}$'$)
cette nouvelle expression, la différence entre~\eqref{eq13.3d} et (\ref{eq13.3d}$'$) vérifie grâce à~(ii) la majoration souhaitée.
Pour majorer (\ref{eq13.3d}$'$), on peut reprendre le raisonnement qui a servi à majorer~\eqref{eq13.5d}. Puisque $\lVert X_{\Qdo}\rVert$ reste
borné par $\eta \lVert T\rVert $ et que $X\mapsto X_{\Qdo}$
est injective sur le cône auquel appartient $X$, les coordonnées de $X$ qui s'introduisent dans les calculs ne
perturbent pas ce raisonnement. Il intervient une intégrale supplémentaire
en $X$ qui reste bornée par celle de $\kappa^{\eta T}(X_{\Qdo})$, donc par $\lVert T\rVert ^{D}$ pour un
$D$ convenable, ce qui ne change rien au résultat.
\end{proof}

\setcounter{equation}{0}
\section{Deux lemmes et fin de la preuve de la proposition~\ref{W2.5}}\label{Ws2.8}

On fixe un réel $\rho>\rho_0$, où $\rho_0$ est le réel du~(iii) du lemme précédent.
Posons
{\multlinegap0pt
\begin{multline*}
E^T_1=\lvert \det(\iota)\rvert ^{-1}\\
{}\times\int_{\ga_Q^G}\tsQR(H) \int_{\ima\gb_S^*}
\int_{\ima\gb_{\tdS}}\kappa^{\rho T}(U)
\int_{\mathcal{C}(Q,\Qdo)}\hat{\varphi}^{T+H-X,\Qdo}_{\tdM}(\lambda;U)\dd X\dd U\dd \lambda\dd H\ptf
\end{multline*}}%
En utilisant le corollaire~\ref{W2.6(1)} et les assertions du lemme précédent concernant les expressions~\eqref{eq13.1d}, \eqref{eq13.2d} et~\eqref{eq13.3d},
on voit que cette expression est absolument convergente et que l'on a une majoration
\begin{equation}
\lVert A^T_{s,t}-E^T_1\rVert \ll  \dPO(T)^{-r}\label{eq13.1e}
\end{equation}
pour tout réel $r$.

Notons $\Sigma^+_{\tdS}$ l'ensemble des racines de $\ga_{\tdM}$ qui sont positives pour le parabolique standard
$\tdS$. Pour tout $\Sp\in \mathcal{P}^{\Qdo}(\tdM)$, notons $a(\Sp)$ le nombre d'éléments de
\[
(-\Delta_{\Sp})\cap \Sigma^+_{\tdS}
\]
(ou encore de $(-\Delta_{\Sp}^{\Qdo})\cap \Sigma^+_{\tdS}$)
et $\mathcal{C}^{\Qdo}(\Sp)\subset t(\ga_S)^{\Qdo}
$ le cône formé des
\[
\Biggl(\sum_{\alpha\in \Delta_{\Sp}^{\Qdo}\cap \Sigma^+_{\tdS}}x_{\alpha}\alpha^\vee\Biggr)+
\Biggl(\sum_{\alpha\in(-\Delta_{\Sp}^{\Qdo})\cap \Sigma^+_{\tdS}}y_{\alpha}\alpha^\vee\Biggr)
\]
avec des $x_{\alpha}\geq0$ et des $y_{\alpha}>0$. En remplaçant les exposants $\Qdo$ par $Q$,
on définit de même le cône $\mathcal{C}^Q(\Sp)\subset t(\ga_S)^Q$ pour
tout $\Sp\in \mathcal{P}^Q(\tdM)$.

\begin{lemme} \label{W2.8.1}
Pour tous $Z\in \ga_0^G$\textup, $\lambda\in\gb_S^*$ et $U\in\gb_{\tdS}$\textup, on a
l'égalité
\[
\hvf^{Z,\Qdo}_{\tdM}(\lambda;U)=\sum_{\Sp\in \mathcal{P}^{\Qdo}(\tdM)}(-1)^{a(\Sp)}
\int_{\mathcal{C}^{\Qdo}(\Sp)}\hvf(\lambda;U-Y_{\Sp}(Z)+V,\Sp)\dd V\ptf
\]
\end{lemme}

\addtocounter{equation}{1}
\begin{proof} 
D'après le lemme~\ref{gmorth}, on a
\[
\varphi^{Z,\Qdo}_{\tdM}(\lambda;\Lambda)=\int_{\gH_{\tdM}(\Qdo)}\int_{t(\ga_S)^{\Qdo}}
\ee^{\langle \Lambda,U+Y_{\Qdo}\rangle}\Gamma_{\tdM}^{\Qdo}(U,\YY)
\hvf^{Z,\Qdo}(\lambda;\YY)\dd U\dd \YY,
\]
où
\[
\gH_{\tdM}(\Qdo)=\varprojlim_{P'\in \FF^{\Qdo}(\tdM)} \ga_{P'}^G
\]
et
$\hvf^{Z,\Qdo}(\lambda;\YY)$ est une fonction sur
cet espace qui \og globalise\fg les différentes fonctions $\hvf^Z(\lambda;Y,\Sp)$ pour $\Sp\in\mathcal{P}^{\Qdo}(\tdM)$.
La proposition~\ref{preconv} affirme que
\[
\Gamma_{\tdM}^{\Qdo}(U,\YY)=\sum_{\Sp\in \mathcal{P}^{\Qdo}(\tdM)} (-1)^{a(\Sp)}
\mathbf{1}_{\mathcal{C}^{\Qdo}(\Sp)}(Y^{\Qdo}_{\Sp}-U),
\]
où $\mathbf{1}_{\mathcal{C}^{\Qdo}(\Sp)}$ est la fonction caractéristique du cône $\mathcal{C}^{\Qdo}(\Sp)$.
Il en résulte l'égalité
{\multlinegap0pt\begin{multline*}
\varphi^{Z,\Qdo}_{\tdM}(\lambda;\Lambda)=\int_{t(\ga_S)^{\Qdo}}\ee^{\langle\Lambda,U+Y_{\Qdo}\rangle}
\sum_{\Sp\in \mathcal{P}^{\Qdo}(\tdM)}
(-1)^{a(\Sp)}\\ {}\times\int_{\ga_{\Sp}^G} \mathbf{1}_{\mathcal{C}^{\Qdo}(\Sp)}(Y^{\Qdo}-U)\hvf^Z(\lambda;Y,\Sp)\dd Y\dd U\ptf
\end{multline*}}%
D'où
{\multlinegap0pt\begin{multline*}
\hvf^{Z,\Qdo}_{\tdM}(\lambda;U)\\
=\sum_{\Sp\in \mathcal{P}^{\Qdo}(\tdM)} (-1)^{a(\Sp)}\int_{t(\ga_S)^{\Qdo}}
\mathbf{1}_{\mathcal{C}^{\Qdo}(\Sp)}(V-U^{\Qdo})\hvf^Z
(\lambda;V+U_{\Qdo},\Sp)\dd V\ptf
\end{multline*}}%
On vérifie que
\[
\hvf^Z(\lambda;Y,\Sp)=\hvf(\lambda;Y-Y_{\Sp}(Z),\Sp)\ptf
\]
Par le changement de variable
$V\mapsto U^{\Qdo}+V$, l'intégrale intérieure devient
\[
\int_{\mathcal{C}^{\Qdo}(\Sp)}\hvf(\lambda;U-Y_{\Sp}(Z)+V;\Sp)\dd V
\]
et la formule ci-dessus devient celle de l'énoncé.
\end{proof}

L'intégrale intérieure de l'expression $E_1^T$ devient
\[
\sum_{\Sp\in \mathcal{P}^{\Qdo}(\tdM)} (-1)^{a(\Sp)}\int_{\mathcal{C}(Q,\Qdo)}
\int_{\mathcal{C}^{\Qdo}(\Sp)}\hvf(\lambda;U-Y_{\Sp}(T-X)-H+V,\Sp)\dd V\dd X\;,
\]
cette expression étant évidemment absolument convergente. Pour $\Sp\in \mathcal{P}^{\Qdo}(\tdM)$, choisissons
$u\in \weyl^{\Qdo}$ tel que $u(\Sp)$ soit standard et notons $\widetilde{Y}_{\Sp}(X)$
la projection de $u^{-1}X$ sur $t(\ga_S)^G$. On a $Y_{\Sp}(T-X)=Y_{\Sp}(T)-\widetilde{Y}_{\Sp}(X)$.
Il résulte des définitions que l'application
\[
\begin{aligned}
\mathcal{C}(Q,\Qdo)\times \mathcal{C}^{\Qdo}(\Sp)&\to t(\ga_S)^Q\\ 
(X,V)&\mapsto\widetilde{Y}_{\Sp}(X)+V
\end{aligned}
\]
est injective et a pour image le cône $\mathcal{C}^Q(\Sp)$. Elle ne respecte pas les mesures euclidiennes.
Pour nous, les cônes $\mathcal{C}^{\Qdo}(\Sp)$ et $\mathcal{C}^Q(\Sp)$ sont bien munis de
ces mesures, mais $\mathcal{C}(Q,\Qdo)$ est muni de la mesure introduite à la proposition~\ref{W2.5}, pour laquelle la projection sur
$\ga_{\Qdo}^Q$ préserve les mesures, ce dernier espace étant
muni de la mesure euclidienne. On voit alors que l'application ci-dessus préserve les mesures.
L'expression précédente devient
\[
\sum_{\Sp\in \mathcal{P}^{\Qdo}(\tdM)} (-1)^{a(\Sp)}\int_{\mathcal{C}^Q(\Sp)}\hvf(\lambda;U-Y_{\Sp}(T)-H+V,\Sp)\dd V\ptf
\]
Il y a un lemme similaire au précédent où $\Qdo$ est remplacé par $Q$. La somme ci-dessus vaut donc
\[
\hvf^{T+H,Q}_{\tdM}(\lambda;U)-\sum_{\Sp\in \mathcal{P}^Q(\tdM)-\mathcal{P}^{\Qdo}(\tdM)} (-1)^{a(\Sp)}
\int_{\mathcal{C}^Q(\Sp)}\hvf(\lambda;U(S',H,T)+V,\Sp)\dd V
\]
avec
\[
U(S',H,T)=U-Y_{\Sp}(T)-H\ptf
\]
On en déduit l'égalité
\begin{equation}
E_1^T=E_{2}^T-E_{3}^T,\label{eq13.2e}
\end{equation}
où
\[
E_{2}^T=\lvert \det(\iota)\rvert ^{-1}\int_{\ga_Q^G}\tsQR(H) \int_{\ima\gb_S^*}
\int_{\gb_{\tdS}}\kappa^{\rho T}(U)\hvf^{T+H,Q}_{\tdM}
(\lambda;U)\dd U\dd \lambda\dd H,
\]
et
{\multlinegap0pt\begin{multline*}
E_{3}^T=\lvert \det(\iota)\rvert ^{-1}\int_{\ga_Q^G}\tsQR(H) \int_{\ima\gb_S^*}
\int_{\gb_{\tdS}}\kappa^{\rho T}(U) \sum_{\Sp\in \mathcal{P}^Q(\tdM)-\mathcal{P}^{\Qdo}(\tdM)} (-1)^{a(\Sp)}\\
{}\times\int_{\mathcal{C}^Q(\Sp)}
\hvf(\lambda;U-Y_{\Sp}(T)-H+V,\Sp)\dd V\dd U\dd \lambda\dd H\ptf
\end{multline*}}%
Cette décomposition est justifiée car l'expression $E^T_{2}$ est absolument convergente
(expression~\eqref{eq13.4d} du lemme~\ref{W2.7}), donc $E^T_{3}$ est convergente au moins dans l'ordre indiqué.
Pour $\Sp\in \mathcal{P}^Q(\tdM)-\mathcal{P}^{\Qdo}(\tdM)$, posons
{\multlinegap0pt\begin{multline*}
E_{\Sp}^T=\smash[b]{\int_{\ga_Q^G}\tsQR(H) \int_{\ima\gb_S^*}\int_{\ima\gb_{\tdS}}}\kappa^{\rho T}(U)\\
{}\times\int_{\mathcal{C}^Q(\Sp)}\lvert \hvf(\lambda;U-Y_{\Sp}
(T)-H+V,\Sp)\rvert \dd V\dd U\dd \lambda\dd H\ptf
\end{multline*}}%
On a évidemment
\begin{equation}
\lvert E_{3}^T\rvert\ll \sum_{\Sp\in \mathcal{P}^Q(\tdM)-\mathcal{P}^{\Qdo}(\tdM)}E_{\Sp}^T\ptf\label{eq13.3e}
\end{equation}

\begin{lemme} \label{W2.8.2}
Pour tout
\[
\Sp\in \mathcal{P}^Q(\tdM)-\mathcal{P}^{\Qdo}(\tdM)
\]
et tout réel $r$, on a une majoration
\[
E^T_{\Sp}\ll \dPO(T)^{-r}\ptf
\]
\end{lemme}

\addtocounter{equation}{3}
\begin{proof} 
On a noté $\gk$ le noyau de $q$. On note $\gk_t$ sa projection sur $t(\ga_S)^G$, ou
son intersection avec cet espace (cela revient au même
puisque $\ga_0^{\tdS}\subset \ga_0^{\Qdo}\subset \gk$). On note $\mathfrak{f}$ l'orthogonal de
$\gk_t$ dans $t(\ga_S)^G$ (qui est
aussi l'orthogonal de $\gk$ dans $\ga_0^G$). On fixe une fonction $C\in C_{c}^{\infty}(\ima\gb_S^*)$ et
des fonctions de Schwartz $\phi_{k}$ sur $\gk_t
$, $\phi_{f}$ sur $\mathfrak{f}$, à valeurs positives ou nulles, de sorte que
\[
\lvert \hvf(\lambda;V,\Sp)\rvert \leq C(\lambda)\phi_{k}(V_{k})\phi_{f}(V_{f}),
\]
avec des notations familières. Grâce au lemme~\ref{W2.7}\eqref{eq13.10d}, et puisque $C$ est intégrable, on obtient
{\multlinegap0pt\begin{multline*}
E_{\Sp}^T\ll \smash[b]{\int_{\ga_Q^G}\tsQR(H)\int_{\ima\gb_{\tdS}}\kappa^{\rho T}(U)\int_{\mathcal{C}^Q(\Sp)}}
\phi_{k}(U-Y_{\Sp,k}(T)-H_{k}+V_{k})\\
{}\times\phi_{f}(-Y_{\Sp,f}(T)-H_{f}+V_{f})\dd V\dd U\dd H\ptf
\end{multline*}}%
Notons $\mathcal{C}^Q(\Sp)_{\Qdo}$ la projection orthogonale de $\mathcal{C}^Q(\Sp)$ dans $\ga_{\Qdo}^Q$. Pour
\[
X\in \mathcal{C}^Q(\Sp)_{\Qdo}
\]
la fibre de cette projection est de la forme
$V(X)+\mathcal{C}(X)$ où $V(X)$ est un élément fixé de la fibre et $\mathcal{C}(X)$ est un sous-ensemble convexe de
$t(\ga_S)^{\Qdo}$. Si $X$ est en position générale, ce
sous-ensemble est d'intérieur non vide. On peut décomposer
l'intégrale sur $\mathcal{C}^Q(\Sp)$ en une intégrale sur $\mathcal{C}^Q(\Sp)_{\Qdo}$ d'intégrales sur les fibres. On peut majorer
brutalement ces dernières intégrales par l'intégrale sur $t(\ga_S)^{\Qdo}$ tout entier. En se rappelant que cet
espace est contenu dans $\gk_t$, on obtient
{\multlinegap0pt\begin{multline*}
E_{\Sp}^T\ll \int_{\ga_Q^G}\tsQR(H)\int_{\ima\gb_{\tdS}}\kappa^{\rho T}(U)
\int_{\mathcal{C}^Q(\Sp)_{\Qdo}}\phi_{f}(-Y_{\Sp,f}(T)-H_{f}+X_{f})\\
{}\times\int_{t(\ga_S)^{\Qdo}}\phi_{k}(U-Y_{\Sp,k}(T)-H_{k}+V(X)_{k}+V)\dd V\dd X\dd U\dd H\ptf
\end{multline*}}%
Puisque $\phi_{k}$ est de Schwartz, la dernière intégrale intérieure est majorée indépendamment de $T$, $H$, $X$ et $U$.
On intègre ensuite en $U$ la fonction $\kappa^{\rho T}(U)$.
Cette intégrale est essentiellement majorée par $\lVert T\rVert ^D$ pour un $D$ convenable. D'où
\[
E_{\Sp}^T\ll \lVert T\rVert ^D\int_{\ga_Q^G}\tsQR(H)\int_{\mathcal{C}^Q(\Sp)_{\Qdo}}
\phi_{f}(-Y_{\Sp,f}(T)-H_{f}+X_{f})\dd X\dd H\ptf
\]
On va montrer que pour $H$ tel que $\tsQR(H)=1$ et pour $X\in \mathcal{C}^Q(\Sp)_{\Qdo}$, on a une majoration
\begin{equation}
\lVert T\rVert +\lVert H\rVert +\lVert X\rVert \ll 1+\lVert -Y_{\Sp,f}(T)-H_{f}+X_{f}\rVert \ptf\label{eq13.4e}
\end{equation}
On en déduira une majoration
\[
E_{\Sp}^T\ll \lVert T\rVert ^{D-r}\int_{\ga_Q^G}\int_{\mathcal{C}^Q(\Sp)_{\Qdo}}(1+\lVert H\rVert)^{-r}
(1+\lVert X\rVert)^{-r}\dd X\dd H
\]
pour tout réel $r$. Ceci est essentiellement majoré par $\dPO(T)^{-r}$ pour tout $r$, ce qui démontrera le lemme.
Revenons à~\eqref{eq13.4e}. Le cône $\mathcal{C}^Q(\Sp)$ est engendré par des $\alpha^\vee$ pour
$\alpha\in\Sigma_{\tdS}^+\cap(\pm \Delta_{\Sp}^Q)$. Donc $\mathcal{C}^Q(\Sp)_{\Qdo}$ est contenu
dans le cône engendré par les $\alpha^\vee_{\Qdo}$ pour $\alpha\in \Delta_0^Q- \Delta_0^{\Qdo}$. Autrement dit, un
élément $X$ de ce cône vérifie $\phi_{\Qdo}^Q(-X)
=1$. On ne fait que renforcer~\eqref{eq13.4e} en remplaçant l'hypothèse
$X\in \mathcal{C}^Q(\Sp)_{\Qdo}$ par $X\in \ga_{\Qdo}^Q$
et $\phi_{\Qdo}^Q(-X)=1$. Puisque $\mathfrak{f}$ est
l'orthogonal de $\gk$, on a en tout cas
\[
\lVert q(Y_{\Sp}(T)+H-X)\rVert \ll \lVert -Y_{\Sp,f}(T)-H_{f}+X_{f}\rVert \ptf
\]
Soit $u\in \weyl^Q$ tel que $u(\Sp)$ soit standard. On a
\[
Y_{\Sp}(T)=(u^{-1}T+T_0-u^{-1}T_0)_{\tdS}=T'_{\tdS}-X(T)_{\tdS}
\]
où
\[
T'=T+T_0-u^{-1}T_0\Qquad{ et}X(T)=T-u^{-1}T\ptf
\]
Puisque $
\ga_0^{\tdS}\subset \ga_0^{\Qdo}\subset \gk$ et puisque $q(T)=0$, on a
\[
q\bigl(Y_{\Sp}(T)\bigr)=q(T'_{\tdS}-X(T)_{\tdS})=q(T'-X(T)_{\Qdo})=q(T_0-u^{-1}T_0-X(T)_{\Qdo})\ptf
\]
Remarquons que
\[
\phi_{\Qdo}^Q(-X(T)_{\Qdo})=1\ptf
\]
En appliquant le lemme~\ref{recroissb} à $H$ et $-X-X(T)_{\Qdo}$, on obtient
\[
\lVert H\rVert +\lVert X(T)_{\Qdo}+X\rVert \ll 1+\lVert q(Y_{\Sp}(T)+H-X)\rVert\;,
\]
le $1$ servant à se débarrasser du terme constant $q(T_0-u^{-1}T_0)$. Puisque $X(T)_{\Qdo}$ et $X$ sont dans un
même  \og vrai\fg cône, on peut remplacer
$\lVert X(T)_{\Qdo}+X\rVert $ par
\[
\lVert X(T)_{\Qdo}\rVert +\lVert X\rVert \ptf
\]
Pour obtenir~\eqref{eq13.4e}, il suffit maintenant de prouver une majoration
\[
\lVert T\rVert \ll \lVert X(T)_{\Qdo}\rVert \ptf
\]
D'après le lemme~\ref{wT}, $X(T)$ est combinaison linéaire des $\alpha^\vee$ pour
$\alpha\in \Rac(u)$, avec des coefficients supérieurs ou égaux à $\dPO(T)$. Or l'hypothèse que
$\Sp\nsubset \Qdo$ implique qu'il y a au moins un
$\alpha\in \Rac(u) $ tel que sa projection $\alpha^\vee_{\Qdo}$ ne soit pas nulle. Dans la base
\[
\{\alpha^\vee_{\Qdo}\mid \alpha\in \Delta_0^Q-\Delta_0^{\Qdo}\}
\]
de $\ga_{\Qdo}^Q$, au
moins l'un des coefficients de $X(T)_{\Qdo}$ est donc supérieur ou égal à $\dPO(T)$. Cela prouve la majoration
précédente et cela achève la démonstration du lemme~\ref{W2.8.2}.
\end{proof}

Nous pouvons maintenant achever la preuve de la proposition~\ref{W2.5}.
Posons
\[
E_{4}^T=\lvert \det(\iota)\rvert ^{-1}\int_{\ga_Q^G}\tsQR(H)\int_{\ima\gb_S^*}
\int_{\gb_{\tdS}}\hvf_{\tdM}^{T+H,Q}(\lambda;U)\dd U\dd 
\lambda\dd H\ptf
\]
Ceci est encore convergent et l'assertion du lemme~\ref{W2.7} concernant l'expression~\ref{W2.7}\eqref{eq13.5d} montre que
\begin{equation}
\lvert E_{2}^T-E_{4}^T\rvert \ll \dPO(T)^{-r}\label{eq13.5e}
\end{equation}
pour tout réel $r$. Par inversion de Fourier de l'intégrale intérieure, on a
\[
E_{4}^T=\lvert \det(\iota)\rvert ^{-1}\int_{\ga_Q^G}\tsQR(H)\int_{\ima\gb_S^*}
\int_{\ima\gc_{\tdS}^*}\varphi_{\tdM}^{T+H,Q}(\lambda;\Lambda)\dd \Lambda\dd \lambda\dd H\ptf
\]
Cette expression est convergente dans l'ordre indiqué.
On peut regrouper les deux intégrales intérieures (les fonctions sont à supports compacts en $\lambda$ et $\Lambda$). En
utilisant le lemme~\ref{W2.6}(ii) et par le changement de variable $(\lambda,\Lambda)\mapsto \mu(\lambda,\Lambda)$,
on obtient
\[
E_{4}^T=\int_{\ga_Q^G}\tsQR(H)\int_{\ima(\ga_S^\G)^*}\omega_{s,t}^{T,Q}(H,\mu)B(\mu)\dd \mu
\dd H=\mathbf{A}^T_{s,t}(B)\ptf
\]
Cela démontre que cette expression est convergente dans l'ordre indiqué.
En additionant~\eqref{eq13.1e}, \eqref{eq13.2e}, \eqref{eq13.3e}, \eqref{eq13.5e} et le lemme~\ref{W2.8.2}, on
obtient la majoration
\[
\lvert A_{s,t}^T-\mathbf{A}_{s,t}^T\rvert \ll \dPO(T)^{-r}
\]
pour tout $r$. Ceci achève la preuve de la proposition~\ref{W2.5}.
\qed

Rappelons les hypothèses sur $Q$ et $R$: on a $\PO\subset Q\subset R$ et $\tvedQR\ne 0$.
On a défini $A^T(B)$ dans la section~\ref{sec12.8}. Pour $t\in \mathbf{W}^{\pQ}(\ga_{S},Q_{0})$ et
$s\in \mathbf{W}^Q\bigl(\theta_{0}(\ga_{S}),t(\ga_{S})\bigr)$, on a défini $\mathbf{A}_{s,t}^T(B)$ à la section~\ref{Apure}.

\begin{corollaire}\label{Wbiscor}
Pour tout réel $r$\textup, on a la majoration
\[
\Biggl\lvert A^T(B)-\sum_{t\in W^{\pQ}(\ga_{S},Q_{0})}\sum_{s\in \mathbf{W}^Q\left(\theta_{0}(\ga_{S}),t(\ga_{S})\right)}
\mathbf{A}_{s,t}^T(B)\Biggr\lvert \ll  \mathbf{d}_{\PO}(T)^{-r}.
\]
\end{corollaire}

\begin{proof}
Cela résulte des propositions~\ref{W2.4}, \ref{W2.5} et de la relation~\ref{W2.5}(i).
\end{proof}

\section{Élargissement des sommations}\label{provi}

Rappelons (\cf lemme~\ref{weylgab}) que $\weyl^{\pQ}(\ga_S,\Qdo)$ est l'ensemble des restrictions à $\ga_S$
d'éléments $t\in \weyl^{\pQ}$ tels que $t(\ga_S)\supset
\ga_{\Qdo}$ et $t\in [\weyl^{\Qdo}\backslash \weyl^{\pQ}]$, où on note ainsi l'ensemble des $t\in \weyl^{\pQ}$ qui sont de
longueur minimale dans leur classe $\weyl^{\Qdo} t$.

\begin{lemme}\label{inclusi}
On a l'inclusion $\weyl^{\pQ}(\ga_S,\Qdo)\subset \weyl^\G(\ga_S,Q)$.
\end{lemme}

\begin{proof}
Il suffit de montrer que
\begin{equation}
[\weyl^{\Qdo}\backslash \weyl^{\pQ}]=[\weyl^Q\backslash \weyl^\G]\cap\weyl^{\pQ}\ptf\label{eq13.1f}
\end{equation}
D'après le lemme~\ref{weylg}, un élément $t\in [\weyl^{\Qdo}\backslash \weyl^{\pQ}]$
vérifie $t^{-1}\alpha>0$ pour $\alpha\in \Delta_0^{\Qdo}$.
Il vérifie aussi $t^{-1}\alpha>0$ pour $\alpha\in \Delta_0^Q-\Delta_0^{\Qdo}$
car une telle racine $\alpha$ intervient dans le radical unipotent de $\pQ$ et l'ensemble des racines intervenant dans ce
radical unipotent est conservé par $t$. Donc $t^{-1}\alpha>0$ pour tout $\alpha\in \Delta_0^Q$,
ce qui signifie, encore d'après le lemme~\ref{weylg},
que $t$ est de longueur minimale dans sa classe $\weyl^Qt$. Cela
prouve l'inclusion du membre de gauche de~\eqref{eq13.1f} dans celui de
droite. Inversement, soit $w\in [\weyl^Q\backslash \weyl^\G]\cap
\weyl^{\pQ}$.
Soit $t$ l'élément
de longueur minimale dans
la classe $\weyl^{\Qdo} w$. Par ce que l'on vient de prouver,
$t$ est aussi de longueur minimale dans $\weyl^\Q t$. Mais $\weyl^Qw=\weyl^Qt$, donc $w=t$.
\end{proof}

On va maintenant élargir les hypothèses sur $Q$ et $\R$:
on suppose seulement
\[
\PO\subset Q\subset R
\]
et on abandonne l'hypothèse $\tvedQR\ne 0$, mais on conserve les hypothèses sur $S$, $\sigma$ et $B$.
Pour $t\in \mathbf{W}^G(\ga_{S},Q)$, on pose
\[
\tvedQRt=\sum_{\tP} (-1)^{a_{\tP}-a_{\tG}}
\qquad\text{avec $\Q\subset P\subset R$ et $t\in W^P$} \ptf\label{eq13.2f}
\]
\newnot{eta(Q,R;t)@$\tvedQRt$}{tvedqrt}%
L'ensemble des $\tP$ satisfaisant ces conditions
peut être vide. S'il est non vide, il existe $\tP_{1}\subset \tP_{2}$ tel que ce soit l'ensemble des $\tP$
tels que $\tP_{1}\subset \tP\subset \tP_{2}$.
Remarquons en passant que $P _{2}=\Rm$.
On en déduit que $\tvedQRt$ est non nul si et seulement s'il existe un unique $\tP$ vérifiant
les conditions de~\eqref{eq13.2f}.
Et dans ce cas on a:
\[
\tvedQRt=(-1)^{a_{\tRm}-a_{\tG}}\ptf
\]

On a défini à la section~\ref{Apure} une expression $\mathbf{A}^T_{s,t}(B)$ pour
\[
t\in \weyl^{\pQ}(\ga_S,\Qdo)\Qquad{et}s\in \weyl^Q\bigl(\theto (\ga_S),t(\ga_S)\bigr)\ptf
\]
On peut aussi bien la définir, au moins formellement, pour $t\in \weyl^\G(\ga_S,Q)$
et $s$ comme précédemment.

\begin{proposition}\label{W2.9.1bis}
Soient $t\in \mathbf{W}^G(\ga_{S},Q)$ et $s\in \mathbf{W}^Q\bigl(\theto(\ga_S),t(\ga_{S})\bigr)$.
On suppose $\tvedQRt\ne 0$.
\begin{enumerate}[(i)]
\item L'expression $\mathbf{A}^T_{s,t}(B)$ est convergente dans l'ordre indiqué.
\item Supposons $t\notin \mathbf{W}^{\pQ}(\ga_{S},Q_{0})$. Alors on a une majoration
\[
\lvert \mathbf{A}_{s,t}^T(B)\rvert \ll \mathbf{d}_{\PO}(T)^{-r}
\]
pour tout réel $r$.
\end{enumerate}
\end{proposition}

\addtocounter{equation}{2}
\begin{proof} 
On reprend la preuve de la proposition~\ref{W2.5}.
Le lemme~\ref{W2.6}(ii) s'applique, avec les mêmes définitions. On obtient
\[
\mathbf{A}_{s,t}^T(B)=\lvert \det(\iota)\rvert ^{-1}\int_{\ga_{Q}^G}\tilde{\sigma}_{Q}^R(H)
\int_{\ima\mathfrak{b}_{S}^*}\int_{\ima\mathfrak{c}_{_{t}S}^*}\varphi_{_{t}M}^{T+H,Q}(\lambda;\Lambda)\dd \Lambda\dd \lambda\dd H,
\]
puis, par inversion de Fourier,
\[
\mathbf{A}_{s,t}^T(B)=\lvert \det(\iota)\rvert ^{-1}\int_{\ga_{Q}^G}\tilde{\sigma}_{Q}^R(H)
\int_{\ima\mathfrak{b}_{S}^*}\int_{\mathfrak{b}_{_{t}S}}\hvf_{_{t}M}^{T+H,Q}(\lambda;U)\dd U\dd \lambda\dd H\ptf
\]
Plus précisément, la convergence absolue du membre de droite ci-dessus entraîne
la convergence dans l'ordre indiqué de $\mathbf{A}_{s,t}^T(B)$ et la validité de l'égalité précédente.
On doit donc prouver que ce membre de droite est absolument convergent et le majorer sous l'hypothèse de~(ii).
Maintenant,
on reprend la preuve du lemme~\ref{W2.7} consacrée à l'expression~\eqref{eq13.4d}.
Le début de cette preuve vaut aussi pour $t\in \mathbf{W}^G$: ce n'est qu'à partir de la relation \eqref{eq13.10d} du lemme~\ref{W2.7}
qu'était utilisée l'hypothèse $t\in \mathbf{W}^{\pQ}$. On obtient que l'expression déduite du membre de droite
de l'égalité ci-dessus en rempla\c{c}ant la fonction à intégrer par sa valeur absolue est essentiellement majorée par
\[
\sum_{P'\in \mathcal{F}^Q(_{t}M)}I^T(P'),
\]
avec
\[
I^T(P')\ll \int_{\ga_{Q}^G}\tilde{\sigma}_{Q}^R(H)\int_{\mathfrak{b}_{\sharp}}
\Gamma_{_{t}M}^{P'}\bigl(U^{P'},\YY(T)\bigr)\phi_{e}(U_{P',e}-Y_{P',e}(T)-H_{e})\dd U\dd H,
\]
\cf lemme~\ref{W2.7}\eqref{eq13.7d}. Fixons $P'$. On va prouver que
\begin{subequations}\label{eq13.3f}
\begin{equation}
\lVert H\rVert \ll 1+\lVert U_{P',e}-Y_{P',e}(T)-H_{e}\rVert \label{eq13.3fi}
\end{equation}
pour tout $T$, tout $H$ tel que $\tilde{\sigma}_{Q}^R(H)=1$ et tout $U\in \mathfrak{b}_{\sharp}$
tel que $\Gamma_{_{t}M}^{P'}\bigl(U^{P'},\YY(T)\bigr)=1$ et que, d'autre part,
si $t\notin \mathbf{W}^{\pQ}(\ga_{S},Q_{0})$, on a
\begin{equation}
\lVert T\rVert +\lVert H\rVert \ll 1+\lVert U_{P',e}-Y_{P',e}(T)-H_{e}\rVert \label{eq13.3fii}
\end{equation}
\end{subequations}
pour tout $T$, $H$ et $U$ comme ci-dessus.
En admettant cela et puisque $\phi_{e}$ est de Schwartz, on a une majoration
\[
I^T(P')\ll C_{r}\int_{\ga_{Q}^G}(1+\lVert H\rVert)^{-r}\int_{\mathfrak{b}_{\sharp}}
\Gamma_{_{t}M}^{P'}\bigl(U^{P'},\YY(T)\bigr)\dd U\dd H
\]
pour tout réel $r$, avec $C_{r}=1$ sans hypothèse sur $t$ et, par contre, $C_{r}=\lVert T\rVert ^{-r}$
sous l'hypothèse de~(ii). Puisque $U\mapsto U^{P'}$ est injective, l'intégrale en $U$ est essentiellement majorée par
$\lVert T\rVert ^D$ pour un $D$ convenable. L'intégrale en $H$ est convergente,
ce qui démontre le~(i) de l'énoncé. Sous l'hypothèse de~(ii), on obtient
\[
I^T(P')\ll \lVert T\rVert ^{-r}
\]
pour tout réel $r$, ce qui démontre le~(ii).
Démontrons les assertions~\eqref{eq13.3f}. Posons
\[
w=s\theta_{0}(t)^{-1}
\]
et introduisons l'application linéaire
\[
q_{w}\colon V\mapsto \bigl((1-w\theta_{0})(V)\bigr)_{Q}
\]
de $\ga_{0}^G$ dans lui-même.
Les hypothèses sur $s$ et $t$ impliquent que $w\theta_{0}$ conservent $\ga_{0}^{_{t}S}$.
Puisque cet espace est contenu dans $\ga_{0}^Q$, il est aussi contenu dans le noyau de $q_{w}$.
La preuve de l'équation~\eqref{eq13.10d} du lemme~\ref{W2.7} montre que $\mathfrak{b}_{_{t}S}$ est le noyau de $1-w\theta_{0}$
dans $t(\ga_S)^G$, donc est contenu dans celui de $q_{w}$. On a évidemment
\[
\lVert q_{w}(-U_{P',e}+Y_{P',e}(T)+H_{e})\rVert \ll \lVert U_{P',e}-Y_{P',e}(T)-H_{e}\rVert
\]
et il suffit de prouver la majoration
\begin{subequations}\label{eqref13.4f}
\begin{equation}
\lVert H\rVert \ll 1+\lVert q_{w}(-U_{P',e}+Y_{P',e}(T)+H_{e})\rVert
\label{eq13.4fi}
\end{equation}
sans hypothèse sur $t$, respectivement
\begin{equation}
\lVert T\rVert +\lVert H\rVert \ll 1+\lVert q_{w}(-U_{P',e}+Y_{P',e}(T)+H_{e})\rVert
\label{eq13.4fii}
\end{equation}
\end{subequations}
si $t\notin \mathbf{W}^{\pQ}(\ga_{S},Q_{0})$.
Parce que $\mathfrak{d}\subset \mathfrak{b}_{_{t}S}$ est contenu dans le noyau de $q_{w}$, on a
\[
q_{w}(-U_{P',e}+Y_{P',e}(T)+H_{e})=q_{w}(-U_{P'}+Y_{P'}(T)+H)\ptf
\]
Puisque $U\in \mathfrak{b}_{\sharp}$
appartient aussi au noyau de $q_{w}$, on a
\[
q_{w}(-U_{P'})=q_{w}(U^{P'})\ptf
\]
L'hypothèse
$\Gamma_{_{t}M}^{P'}\bigl(U^{P'},\YY(T)\bigr)=1$ signifie que $U^{P'}$ est dans l'enveloppe convexe des
$Y_{S'}(T)^{P'}$ pour $S'\in\mathcal{P}^{P'}(_{t}M)$. On peut donc écrire
\[
U^{P'}=\sum_{S'\in \mathcal{P}^{P'}(_{t}M)}x_{S'}Y_{S'}(T)^{P'}\;,
\]
avec des réels $x_{S'}\geq0$ tels que $\sum_{S'\in \mathcal{P}^{P'}(_{t}M)}x_{S'}=1$. Alors
\[
q_{w}(-U_{P',e}+Y_{P',e}(T)+H_{e})=q_{w}\bigl(H+Y_{\mathbf{x}}(T)\bigr)\;,
\]
où
\[
Y_{\mathbf{x}}(T)=\sum_{S'\in \mathcal{P}^{P'}(_{t}M)}x_{S'}Y_{S'}(T)\ptf
\]
D'après l'hypothèse $\tvedQRt\ne 0$, il existe un unique espace parabolique, que l'on note
$\tP$, qui vérifie la condition~\eqref{eq13.1f}. Parce que $w\in W^P$, on a
\[
q_{w}(V)_P =q_{w}(V_P)=(1-\theta_{0})(V_P)
\]
pour tout $V$. Parce que les $S'$ sont contenus dans $Q$, a fortiori dans $P$, on a $Y_{\mathbf{x}}(T)_P =T_P $, d'où
$q_{w}(Y_{\mathbf{x}}(T)_P)=0$ puisque $T=\theta_{0}(T)$. On en déduit
\[
\Bigl(q_{w}\bigl(H+Y_{\mathbf{x}}(T)\bigr)\Bigr)_P =(1-\theta_{0})(H_P).
\]
Comme dans le corollaire~\ref{recroissb}
on écrit $H$ comme une somme de vecteurs deux à deux orthogonaux
\[
H=X_0+X_1+X_2
\]
avec $X_0=H^P$ et $H_P=X_1+X_2$ où $X_1$ est dans l'image de l'application $(1-\theto)$ et
$X_2$ est la projection de $H_P$ sur le sous-espace des $\theto$\hyph invariants.
Comme $P=\Rm$ il résulte du lemme~\ref{croiss}(ii) que l'on a
\[
\lVert X_2\rVert \ll \lVert X_0\rVert +\lVert X_1\rVert
\]
et donc aussi
\[
\lVert X_1+X_2\rVert\ll \lVert X_0\rVert +\lVert X_1\rVert
\]
et comme
\[
\lVert X_1\rVert\ll \lVert(1-\theto)X_1\rVert=\lVert(1-\theto)H_P\rVert
\]
on a
\[
\lVert H_P \rVert \ll \lVert (1-\theta_{0})(H_P)\rVert +\lVert H^P\rVert \ptf
\]
D'où
\[
\lVert H_P \rVert \ll \lVert \Bigr(q_{w}\bigr(H+Y_{\mathbf{x}}(T)\bigr)\Bigr)_P \rVert +\lVert H^P\rVert \ptf
\]
Puisque d'autre part on a aussi $q_{w}(V)^P=q_{w}(V^P)$ pour tout $V$, il suffit donc de prouver une majoration
\begin{subequations}\label{eqref13.5f}
\begin{equation}
\lVert H^P\rVert \ll 1+\lVert q_{w}(H^P+Y_{\mathbf{x}}(T)^P)\rVert \label{eq13.5fi}
\end{equation}
respectivement
\begin{equation}
\lVert T\rVert +\lVert H^P\rVert \ll 1+\lVert q_{w}(H^P+Y_{\mathbf{x}}(T)^P)\rVert \ptf\label{eq13.5fii}
\end{equation}
\end{subequations}
Pour $S'\in \mathcal{P}^{P'}(_{t}M)$, fixons $u_{S'}\in \mathbf{W}^Q$ tel que $u_{S'}(S')$ soit standard.
On a
\[
Y_{S'}(T)^P=(u_{S'}^{-1}T+T_{0}-u_{S'}^{-1}T_{0})^P_{_{t}S}\ptf
\]
Puisque $\ga^{_{t}S}_{0}$ est contenu dans le noyau de $q_{w}$, on peut supprimer l'indice $_{t}S$:
\[
q_{w}(Y_{S'}(T)^P)=q_{w}(u_{S'}^{-1}T^P+T_{0}^P-u_{S'}^{-1}T_{0}^P).
\]
L'élément
\[
\sum_{S'\in \mathcal{P}^{P'}(_{t}M)}x_{S'}q_{w}(T_{0}^P-u_{S'}^{-1}T_{0}^P)
\]
reste borné, on peut le négliger. Posons $v_{S'}=\theta_{0}(u_{S'}^{-1})$.
On a
\[
(1-w\theta_{0})u_{S'}^{-1}T^P=u_{S'}^{-1}T^P-T^P+X^T(S')
\]
où
\[
X^T(S')=T^P-wv_{S'}\theta_{0} T^P=T^P-wv_{S'}T^P\ptf
\]
Puisque $u_{S'}\in W^Q$, on a $(u_{S'}^{-1}T^P-T^P)_{Q}=0$ et on obtient $q_{w}(u_{S'}^{-1}T^P)=X^T(S')_{Q}$.
De ces calculs résulte la majoration
\[
\Biggl\lVert q_{w}(H^P)+\sum_{S'\in \mathcal{P}^{P'}(_{t}M)}x_{S'}X^T(S')_{Q}\Biggr\rVert \ll 1+\lVert q_{w}(H^P+Y_{\mathbf{x}}(T)^P)\rVert \ptf
\]
Il nous suffit donc de prouver la majoration
\begin{subequations}\label{eqref13.6f}
\begin{equation}
\lVert H^P\rVert \ll \Biggl\lVert q_{w}(H^P)+\sum_{S'\in \mathcal{P}^{P'}(_{t}M)}x_{S'}X^T(S')_{Q}\Biggr\rVert\;,
\label{eq13.6fi}
\end{equation}
sans hypothèse sur $t$, respectivement
\begin{equation}
\lVert T\rVert +\lVert H^P\rVert \ll \lVert q_{w}(H^P)+\sum_{S'\in \mathcal{P}^{P'}(_{t}M)}x_{S'}X^T(S')_{Q}\rVert\;,
\label{eq13.6fii}
\end{equation}
\end{subequations}
si $t\not\in \mathbf{W}^{\pQ}(\ga_{S},Q_{0})$.
L'élément $H^P$ appartient à $\ga_{Q}^P$ et vérifie la condition
\[
\tau_{Q}^P(H^P)=1\ptf
\]
Donc $q_{w}(H^P)$
appartient au cône engendré par les $q_{w}(\varpi^P)$ pour $\varpi\in \hDelta_{Q}-\hDelta_P $.
La définition
\[
X^T(S')=T^P-wv_{S'}T^P
\]
et les propriétés habituelles montrent que $X^T(S')_{Q}$
appartient au cône engendré par les $\check{\alpha}_{Q}$ pour $\alpha\in \Delta_{0}^P-\Delta_{0}^Q$.
Montrons que:
\begin{enumerate}\setcounter{enumi}{6}
\item\label{eq13.7f}\itshape
le cône engendré par les $q_{w}(\varpi^P)$ pour $\varpi\in \hDelta_{Q}-\hDelta_P$
et les $\check{\alpha}_{Q}$ pour $\alpha\in \Delta_{0}^P-\Delta_{0}^Q$ est un \emph{\og \emph{vrai}\fg} cône\textup,
c'est-à-dire ne contient pas d'espace vectoriel non nul.
\end{enumerate}
Soient
\[
Y=\sum_{\varpi\in \hDelta_{Q}-\hDelta_P }y_{\varpi}\varpi^P\Qquad{et}
Z=\sum_{\alpha\in \Delta_{0}^P-\Delta_{0}^Q}z_{\alpha}\check{\alpha}_{Q}
\]
avec des coefficients positifs ou nuls.
Il faut voir que l'égalité $q_{w}(Y)+Z=0$ entraîne $Y=0$ et $Z=0$. L'argument est le même que dans la preuve du lemme~\ref{recroissa}.
Puisque le produit scalaire $(Y,Z)$ est positif ou nul, l'égalité $q_{w}(Y)+Z=0$ entraîne,
par produit scalaire avec $Y$, l'inégalité $\bigl(Y,q_{w}(Y)\bigr)\leq0$. On a
\[
q_{w}(Y)=Y-(w\theta_{0}Y)_{Q}
\]
d'où $(Y,Y)\leq (Y,(w\theta_{0}Y)_{Q})$.
Par Cauchy-Schwartz, cela entraîne les égalités
\[
(w\theta_{0}Y)_{Q}=w\theta_{0}Y=Y\ptf
\]
Introduisons le parabolique standard $P_{1}$
tel que
\[
\hDelta_{P_{1}}=\hDelta_P \cup\{\varpi\in \hDelta_{Q}-\hDelta_P\mid y_{\varpi}\ne 0\}\ptf
\]
L'élément $Y$ appartient à la chambre positive associée au parabolique $P_{1}\cap M_P $ du Levi standard $M_P $ de $P$.
L'égalité $w\theta_{0}Y=Y$ entraîne $w\theta_{0}(P_{1})=P_{1}$. Puisque $P_{1}$ et $\theta_{0}(P_{1})$
sont standard, cela entraîne $P_{1}=\theta_{0}(P_{1})$ puis $w\in \mathbf{W}^{P_{1}}$. On a aussi $Q\subset P_{1}$,
donc $s\in \mathbf{W}^{P_{1}}$. Les relations $w\in \mathbf{W}^{P_{1}}$ et $\theta_{0}(P_{1})=P_{1}$ entraînent alors $t\in\weyl^{P_{1}}$.
D'après l'unicité de $\tP$, on a alors $P_{1}=P$. D'après la définition de $P_{1}$, cela entraîne $Y=0$.
Mais alors l'égalité $q_{w}(Y)+Z=0$ entraîne aussi $Z=0$. Cela prouve~\eqref{eq13.7f}.
On a donc une majoration
\[
\lVert H^P\rVert +\sum_{S'\in \mathcal{P}^{P'}(_{t}M)}x_{S'}\lVert X^T(S')_{Q}\rVert
\ll \Biggl\lVert q_{w}(H^P)+\sum_{S'\in \mathcal{P}^{P'}(_{t}M)}x_{S'}X^T(S')_{Q}\Biggr\rVert \ptf
\]
\emph{A fortiori}, on a la majoration~\eqref{eq13.6fi}.
Pour obtenir~\eqref{eq13.6fii}, il suffit maintenant de prouver que, si $t\notin \mathbf{W}^{\pQ}(\ga_{S},Q_{0})$,
on a une majoration
\[
\lVert T\rVert \ll \lVert X^T(S')\rVert
\]
pour tout $S'$. On reprend l'argument utilisé dans la preuve du lemme~\ref{W2.8.2}.
L'élément $X^T(S')$ est combinaison linéaire de tous les $\check{\alpha}$ pour
$\alpha\in \Rac(v_{S'}^{-1}w^{-1})$, avec des coefficients supérieurs ou égaux à $\mathbf{d}_{\PO}(T)$.
La majoration ci-dessus s'ensuit pourvu qu'il y ait au moins un $\alpha\in\Rac(v_{S'}^{-1}w^{-1})$ tel que $\check{\alpha}_{Q}\ne 0$.
S'il n'en est pas ainsi, on a $v_{S'}^{-1}w^{-1}\in\weyl^Q$, c'est-à-dire $\theta(u_{S'})\theta(t)s^{-1}\in\weyl^Q$, d'où $t'=u_{S'}t\in\weyl^{\pQ}$.
Soit $t''$ l'élément de longueur minimale dans la classe $\weyl^{\Qdo} t'$, il résulte du lemme~\ref{inclusi}
que
\[
t''\in\mathbf{W}^{\pQ}(\ga_{S},Q)\subset\mathbf{W}^G(\ga_{S},Q)\ptf
\]
Mais comme par hypothèse on a aussi $t\in\mathbf{W}^G(\ga_{S},Q)$ il en résulte que $t=t''$ et
donc $t\in\weyl^{\pQ}(\ga_{S},Q_{0})$ contrairement à l'hypothèse. Cela achève la démonstration.
\end{proof}

On pose
\[
\mathbf{A}^T(B)=\sum_{t\in \mathbf{W}^G(\ga_{S},Q)}\tvedQRt
\sum_{s\in \mathbf{W}^Q\left(\theto(\ga_S),t(\ga_S)\right)}\mathbf{A}_{s,t}^T(B)\ptf
\]

\begin{corollaire}\label{W2.9.2bis}
\begin{enumerate}[(i)]
\item Supposons $\tvedQR\ne 0$. Alors on a la majoration
\[
\lvert \tvedQR A^T(B)-\mathbf{A}^T(B)\rvert \ll \mathbf{d}_{\PO}(T)^{-r}
\]
pour tout réel $r$.

\item Supposons $\tvedQR=0$. Alors on a la majoration
\[
\lvert \mathbf{A}^T(B)\rvert \ll \mathbf{d}_{\PO}(T)^{-r}
\]
pour tout réel $r$.
\end{enumerate}
\end{corollaire}

\begin{proof}
La proposition nous dit qu'à des termes négligeables près, on peut imposer la condition
$t\in \mathbf{W}^{\pQ}(\ga_{S},Q_{0})$ dans la définition de $\mathbf{A}^T(B)$.
Mais, pour $t$ dans cet ensemble, on a par définition l'égalité $\tvedQRt=\tvedQR$.
Le~(ii) devient clair tandis que le~(i) résulte du corollaire~\ref{Wbiscor}.
\end{proof}


Pour des sous-groupes paraboliques $Q,S$ vérifiant $\PO\subset\Q$ et $\PO\subset S\subset Q'$,
pour une représentation $\sigma\in \Pi_{\mathrm{disc}}(M_{S})$, pour des éléments
$t\in \mathbf{W}^G(\ga_{S},Q)$ et $s\in \mathbf{W}^Q\bigl(\theto(\ga_S),t(\ga_S)\bigr)$,
pour des éléments $\nu\in \ima(\ga_{S}^G)^*$ et $\tmu\in \ima(\theto(\ga_S)^G)^*$, on considère l'opérateur
{\multlinegap0pt\begin{multline*}
\gmomega^{T,Q}_{s,t}(S,\sigma,f,\omega; \tmu,\nu)=\sum_{\Sp\in \mathcal{P}^Q(\tdM)}\ee^{\langle s\tmu-t\nu,Y_{\Sp}(T)\rangle}
\epsilon^Q_{\Sp}(s\tmu-t\nu)\\
{}\times\Mint(t,\nu)^{-1}\Mint_{\Sp\rest \tdS}(t\nu)^{-1}\Mint_{\Sp\rest \tdS}(s\tmu)\Mint(s,\tmu)
\treg_{S,\sigma,\nu}(f,\omega)\ptf
\end{multline*}}%
D'après la proposition~\ref{GMspec}, cet opérateur est une fonction lisse en $\tmu$ et $\nu$.
On peut donc imposer $\tmu=\theta_{0}\nu$. Introduisons
{\multlinegap0pt\begin{multline*}
\JJJ^T_{\chi}(B,\fff,\omega)\\
\begin{aligned}
&=\sum_{\{Q,R\mid \PO\subset Q\subset R\}}\sum_{\{S\mid \PO\subset S\subset\pQ\}}
\smash[t]{\frac{1}{n^{\pQ}(S)}}\\
&\quad{}\times\sum_{\sigma\in \Pi_{\mathrm{disc}}(M_S)}
\sum_{\WPsi\in \base^{\pQ}_\chi(\sigma)}\sum_{t\in \weyl^\G(\ga_S,Q)}\sum_{s\in \weyl^Q\left(\theto(\ga_S),t(\ga_S)\right)}
\tved(\Q,\R;t)\\
&\qquad{}\times\int_{\ga_Q^G}\tsQR(H) \biggl(\int_{\ima(\ga_S^\G)^*}\mspace{-10mu}
\ee^{\langle (s\theto-t)\mu,H\rangle}\langle
\gmomega^{T,Q}_{s,t}(S,\sigma,f,\omega
;\theto\mu,\mu)\WPsi,\WPsi\rangle
B_{\sigma}(\mu)\dd \mu\biggr)\dd H\ptf
\end{aligned}
\end{multline*}}%
On\hspace*{-1.4pt} dispose\hspace*{-1.4pt} de\hspace*{-1.4pt} l'expression\hspace*{-1.4pt} $J^T_{\chi}(B,f,\omega)$\hspace*{-1.4pt} du\hspace*{-1.4pt} théorème~\ref{thWA}\hspace*{-1.4pt} et\hspace*{-1.4pt}
l'expression\hspace*{-1.4pt} $\JJJ^T_{\chi}(B,f,\omega)$ en est une approximation:

\begin{proposition}\label{2.10bis}
L'expression $\JJJ^T_{\chi}(B,f,\omega)$ est convergente.
On a une majoration
\[
\lvert J^T_{\chi}(B,f,\omega)-\JJJ^T_{\chi}(B,f,\omega)\rvert \ll \mathbf{d}_{\PO}(T)^{-r}
\]
pour tout réel $r$.
\end{proposition}

\begin{proof}
Cela résulte du corollaire~\ref{W2.9.2bis} et du fait que chacune des quadruples intégrales
intervenant dans $J^T_{\chi}(B,\fff,\omega)$ est de la
forme $A^T(\varphi B_{\sigma})$ étudiées ci-dessus, pour
des fonctions $\varphi$ qui sont des coefficients des opérateurs $\treg_{S,\sigma,\mu}(f,\omega)$.
\end{proof}
\chapter{Formules explicites}

\section{Combinatoire finale}\label{sec14.1}

Soient $S$, $\oS$ et $\Q$ trois sous-groupes
paraboliques standard et on suppose que $\oS=\theto S\subset\Q$.
Considérons
\[
t\in\weyl^G(\ga_S,Q)\Qquad{et}s\in\weyl^Q\bigl(\ga_{\oS},t(\ga_S)\bigr)
\]
alors on pose $\us=t\moins s$ et
$\tus=\us\theto$ appartient à $\weyl^{\tG}(\ga_S,\ga_S)$.
On a la réciproque:

\begin{lemme}\label{tst}
Tout $\tus\in\weyl^{\tG}(\ga_S,\ga_S)$ s'écrit d'une façon et d'une seule
sous de la forme
\[
\tus=t\moins s\theto
\]
avec $s$ et $t$ comme ci-dessus.
\end{lemme}

\begin{proof}
On écrit $\tus=\us\theto$ et on
observe que
$\us\moins$ appartient à $\weyl^\G(\ga_S)$.
D'après le lemme~\ref{echange} il existe un unique élément
$t$ de longueur minimale dans la classe $\weyl^\Q\us\moins$. On pose $s=t\us$ et on a donc
$s\in\weyl^\Q$. D'après le lemme~\ref{weylga}
le couple $(s,t)$ a les propriétés requises.
\end{proof}

Soit $\Sp$ le sous-groupe
parabolique standard tel que $\ga_\Sp=t(\ga_S)$.
On introduit
\[
\uS=t\moins\Sp\ptf
\]
Considérons des paramètres $\mu$ et $\nu$ dans $\ima\ga_S^*$ et
posons
\[
\Lambda=\tus\mu-\nu\ptf
\]
On rappelle que l'on a introduit
\[
Y_\us=\TK-\us\moins\TK=\HO(w_\us\moins)\ptf
\]
Introduisons de plus
\[
Y_{\tus}=\theta_0\moins Y_{\us}
=\theta_0\moins\TK-\tus\moins\TK
\Qquad{et}
a_S(\mu,\tus)=\ee^{\langle\mu+\demisom_S,Y_{\tus}\rangle}\ptf
\]

\begin{lemme}\label{st} 
Avec les notations de la section~\ref{gmspec}\textup,
\begin{equation}
\Mint_{\Sp\rest S}(t,\nu)\moins\Mint_{\Sp\rest \oS}(s,\ttmu)\label{eq14.1}
\end{equation}
est égal à
\begin{equation}
a_S(\mu,\tus) \ee^{\langle\Lambda,Y_t\rangle}
\Mint_{\uS\rest S}(\nu)\moins
\Mint_{\uS\rest S}(\nu+\Lambda)
\Mint_{S\rest \tus S}(\nu+\Lambda)\mathbf{\us}\ptf\label{eq14.2}
\end{equation}
\end{lemme}

\begin{proof}
Le lemme~\ref{sta} montre que l'expression~\eqref{eq14.1} est égale à
\begin{equation}
\ee^{\langle\Lambda,Y_t\rangle}\Mint_{\uS\rest S}(1,\nu)\moins
\Mint_{\uS\rest S}(1,\nu+\Lambda)\Mint_{S\rest \oS}(\us,\ttmu)\ptf\label{eq14.3}
\end{equation}
Maintenant
\[
\Mint_{S\rest \oS}(\us,\ttmu)=
\Mint_{S\rest \us\oS}(1,\nu+\Lambda)\Mint_{\us\oS\rest \oS}(\us,\ttmu)
\]
et, d'après le lemme~\ref{Mintsimple}
\[
\Mint_{\us\oS\rest \oS}(\us,\tmu)=\ee^{\langle\tmu+\demisom_{\oS},Y_{\us}\rangle}\mathbf{\us}
\]
on a donc
\[
\Mint_{S\rest \oS}(\us,\ttmu)=
\ee^{\langle\ttmu+\demisom_{\oS},Y_{\us}\rangle}
\Mint_{S\rest \us\oS}(\nu+\Lambda)\mathbf{\us}\ptf
\]
On a donc montré que~\eqref{eq14.1}
est égal à
\begin{equation}
\ee^{\langle\Lambda, Y_t\rangle +\langle\ttmu+\demisom_{\oS},Y_{\us}\rangle}
\Mint_{\uS\rest S}(\nu)\moins
\Mint_{\uS\rest S}(\nu+\Lambda)
\Mint_{S\rest \us\oS}(\nu+\Lambda)\mathbf{\us}\label{eq14.4}
\end{equation}
et on observe que
\[
\ee^{\langle\ttmu+\demisom_{\oS},Y_{\us}\rangle}=\ee^{\langle\mu+\demisom_S,Y_{\tus}\rangle}=
a_S(\mu,\tus)\ptf
\]
Enfin, on remarque que $\us\oS=\tus S$.
\end{proof}

On va donner une nouvelle expression pour l'opérateur $\gmomega^{T,Q}_{s,t}$
introduit à la section~\ref{provi}.
On suppose $S$ (et donc aussi $\oS$) standard. On note $\M$ le sous-groupe de Levi de $S$
et on pose
\[
\uQ=t\moins\Q=\tus\pQ\ptf
\]
On a $\uQ\in\FF(\M)$.
Suivant la proposition~\ref{GMspec} on pose pour $\uS\in\Parab(\M)$:
\[
\MM(S,\T,\nu;\Lambda,\uS)=\ee^{\langle\Lambda,Y_{\uS}(T)\rangle}\,
\Mint_{\uS\rest S}(\nu)\moins\Mint_{\uS\rest S}(\nu+\Lambda)
\]
et
\[
\MM_{\M}^{\uQ}(S,\T,\nu;\Lambda)=
\sum_{\uS\in\Parab^{\uQ}(\M)}\epsilon_{\uS}^{\uQ}(\Lambda)\MM(S,\T,\nu;\Lambda,\uS)\ptf
\]

\begin{proposition}\label{page47}
Si $\Lambda=\tus\mu-\nu$
\[
\gmomega^{T,Q}_{s,t}(S,\sigma,f,\omega; \ttmu,\nu)\]
est égal à
\[
\frac{a_S(\mu,\tus)}{a_S(\nu,\tus)}\MM_{\M}^{\uQ}(S,\T,\nu;\Lambda)
\Mint_{S\rest \tus S}(\nu+\Lambda)\reg_{S,\sigma,\nu}(\tus,\fff,\omega)
\ptf
\]
\end{proposition}

\begin{proof}
On rappelle que, par définition,
\[
\gmomega^{T,Q}_{s,t}(S,\sigma,f,\omega; \tmu,\nu)
\]
est égal à
{\multlinegap0pt\begin{multline*}
\smash[b]{\sum_{\Sp\in {\mathcal P}^Q(\tdM)}}\ee^{\langle s\tmu-t\nu,Y_{\Sp}(T)\rangle}\epsilon^Q_{\Sp}(s\tmu-t\nu)\\
{}\times\Mint(t,\nu)^{-1}\Mint_{\Sp\rest \tdS}(t\nu)^{-1}\Mint_{\Sp\rest \tdS}(s\tmu)\Mint(s,\tmu)
\treg_{S,\sigma,\nu}(f,\omega)
\ptf
\end{multline*}}%
On commence par observer que
\[
\Mint_{\Sp\rest \tdS}(s\tmu)\Mint(s,\tmu)=
\Mint_{\Sp\rest \tdS}(1,s\tmu)\Mint_{\tdS\rest \oS}(s,\tmu)=\Mint_{\Sp\rest \oS}(s,\tmu)
\]
et
\[
\Mint(t,\nu)^{-1}\Mint_{\Sp\rest \tdS}(t\nu)^{-1}=
\Mint_{\tdS\rest S}(t,\nu)^{-1}\Mint_{\Sp\rest \tdS}(1,t\nu)^{-1}=
\Mint_{\Sp\rest S}(t,\nu)^{-1}\ptf
\]
On en déduit que
{\multlinegap0pt
\begin{multline*}
\gmomega^{T,Q}_{s,t}(S,\sigma,f,\omega; \tmu,\nu)\\
= \sum_{\Sp\in {\mathcal P}^Q(\tdM)}\ee^{\langle s\tmu-t\nu,Y_{\Sp}(T)\rangle}
\epsilon^Q_{\Sp}(s\tmu-t\nu)
\Mint_{\Sp\rest S}(t,\nu)^{-1}\Mint_{\Sp\rest \oS}(s,\tmu)
\treg_{S,\sigma,\nu}(f,\omega)
\ptf
\end{multline*}}%
On va utiliser ceci pour $\tmu=\ttmu$.
On rappelle que
\[
\Lambda=\tus\mu-\nu=t\moins(s\ttmu-t\nu)\ptf
\]
On a donc
\[
\ee^{\langle s\ttmu-t\nu,Y_{\Sp}(T)\rangle}
\epsilon^Q_{\Sp}(s\ttmu-t\nu)=\ee^{\langle \Lambda,t\moins Y_{\Sp}(T)\rangle}
\epsilon^{\uQ}_{\uS}(\Lambda)
\]
D'après le lemme~\ref{asu} on a
\[
\reg_{S,\sigma,\nu}(\tus,\fff,\omega)=
\ee^{\langle\nu+\demisom_S,\theta_0\moins\HO(w_\us\moins)\rangle}
\mathbf{\us}\treg_{S,\sigma,\nu}(\fff,\omega)
=a_S(\nu,\tus)\mathbf{\us}\treg_{S,\sigma,\nu}(\fff,\omega)
\ptf
\]
Il résulte alors du lemme~\ref{st} que
\[
\Mint_{\Sp\rest S}(t,\nu)\moins\Mint_{\Sp\rest \oS}(s,\ttmu)
\treg_{S,\sigma,\nu}(\fff,\omega)
\]
est égal à
\[
\frac{a_S(\mu,\tus)}{a_S(\nu,\tus)}\ee^{\langle\Lambda, Y_t\rangle}
\Mint_{\uS\rest S}(\nu)\moins
\Mint_{\uS\rest S}(\nu+\Lambda)
\Mint_{S\rest \tus S}(\nu+\Lambda)\reg_{S,\sigma,\nu}(\tus,\fff,\omega)\ptf
\]
Enfin, on observe que
\[
t\moins Y_\Sp(T)+Y_t=Y_{\uS}(T)\ptf
\]
On obtient alors que
\[
\gmomega^{T,Q}_{s,t}(S,\sigma,f,\omega; \ttmu,\nu)
\]
est égal à
{\multlinegap0pt\begin{multline*}
\smash[b]{\frac{a_S(\mu,\tus)}{a_S(\nu,\tus)}\sum_{\uS\in {\mathcal P}^{\uQ}(\M)}}
\ee^{\langle\Lambda,Y_{\uS}(T)\rangle}
\epsilon^{\uQ}_{\uS}(\Lambda)\\
{}\times\Mint_{\uS\rest S}(\nu)\moins
\Mint_{\uS\rest S}(\nu+\Lambda)
\Mint_{S\rest \tus S}(\nu+\Lambda)\reg_{S,\sigma,\nu}(\tus,\fff,\omega)
\ptf
\end{multline*}}%
En introduisant les opérateurs $\MM$,
ceci se récrit
\[
\frac{a_S(\mu,\tus)}{a_S(\nu,\tus)}
\sum_{\uS\in {\mathcal P}^{\uQ}(\M)}
\epsilon^{\uQ}_{\uS}(\Lambda)\MM(S,\T,\nu;\Lambda,\uS)
\Mint_{S\rest \tus S}(\nu+\Lambda)\reg_{S,\sigma,\nu}(\tus,\fff,\omega)
\]
soit encore
\[
\frac{a_S(\mu,\tus)}{a_S(\nu,\tus)}\MM_{\M}^{\uQ}(S,\T,\nu;\Lambda)
\Mint_{S\rest \tus S}(\nu+\Lambda)\reg_{S,\sigma,\nu}(\tus,\fff,\omega)
\ptf\qedhere
\]
\end{proof}

On introduit, pour $\Lambda=(\tus-1)\mu$
\[
 \AAA_\M^{\uQ}(\T,\sigma,\tus,\mu)=
\tr\bigl(\MM_{\M}^{\uQ}(S,\T,\mu;\Lambda)
\Mint_{S\rest \tus S}(\mu+\Lambda)\reg_{S,\sigma,\mu}(\tus,\fff,\omega)\bigr)\ptf
\]

\begin{lemme}\label{bluval}
\begin{enumerate}[(i)]
\item On a\textup:
\[
\AAA_\M^{\uQ}(\T,\sigma,\tus,\mu)=
\sum_{\WPsi\in \base^{\pQ}_\chi(\sigma)}
\langle\gmomega^{T,Q}_{s,t}(S,\sigma,f,\omega;\ttmu,\mu)\WPsi,\WPsi\rangle
\ptf\]
\item L'expression pour $\AAA_\M^{\uQ}(\T,\sigma,\tus,\mu)$
est invariante si on remplace $\tus$\textup, $\M$\textup, $\uQ$\textup, $S$ et $\sigma$
par des conjugués sous un élément du groupe de Weyl\textup, et simultanément
$\mu$ par le translaté par le même élément.
\end{enumerate}
\end{lemme}

\begin{proof}
D'après la proposition~\ref{page47} et compte tenu du lemme~\ref{tst}, on sait que pour $\Lambda=\tus\mu-\nu$
\[
\gmomega^{T,Q}_{s,t}(S,\sigma,f,\omega; \ttmu,\nu)=
\MM_{\M}^{\uQ}(S,\T,\nu;\Lambda)
\Mint_{S\rest \tus S}(\nu+\Lambda)\reg_{S,\sigma,\mu}(\tus,\fff,\omega)\ptf
\]
D'après la proposition~\ref{GMspec} la fonction $\MM_{\M}^{\uQ}$
est lisse pour les valeurs imaginaires pures des variables $\nu$ et $\Lambda$.
On peut donc imposer $\mu=\nu$.
Ceci prouve~(i).
L'assertion~(ii) est une conséquence directe des équations fonctionnelles
satisfaites par les opérateurs d'entrelacement.
\end{proof}

La fonction $B_\sigma$ et le sous-groupe parabolique $S$ étant fixés, on note
\[
\JJJ_{\uQ}^{\uR}(\M,\T,\sigma,\tus)
\]
l'intégrale double itérée
\[
\int_{\ga_{\uQ}^\G}\tsuQuR(H_1)
\biggl(\int_{\ima(\ga_\M^\G)^*} \ee^{\langle(\tus-1)\mu,H_1\rangle}
\AAA_\M^{\uQ}(\T,\sigma,\tus,\mu)
B_{\sigma}(\mu)\dd\mu\biggr)\dd H_1\ptf
\]
Ici $\M$ est le sous-groupe de Levi de $S$ (qui est standard), on rappelle que $S\subset\pQ\subset\R$ et on a posé
\[
\uQ=\tus\pQ=t\moins\Q\Quad{et} \uR=t\moins\R\ptf
\]
On observera que plus généralement, $\M$ étant donné, cette expression est bien définie pour
tout $S\in\Parab(\M)$ et tout
\[
\tus\in\weyl^{\tG}(\ga_\M,\ga_\M)\ptf
\]
Rappelons que pour énoncer la proposition~\ref{2.10bis} nous avons introduit l'expression~$\JJJ^T_{\chi}$.

\begin{proposition}\label{symetrisation}
L'expression pour $\JJJ^T_{\chi}$ se récrit
\[
\JJJ^T_{\chi}(B,\fff,\omega)=\sum_{\M\in\mathcali L^\G/\weyl^\G}
\frac{1}{w^\G(\M)}\sum_{\sigma\in \Pi_{\mathrm{disc}}(M)_\chi}
\sum_{\tus\in\weyl^{\tG}(\ga_\M,\ga_\M)}\JJJ^T_{\M}(B,\sigma,\fff,\omega,\tus)
\]
avec
\[
\JJJ^T_{\M}(B,\sigma,\fff,\omega,\tus)=
\sum_{\{\uQ,\uR\mid \M\subset\uQ\subset\uR\}}
\tved(\uQ,\uR;\us)
\JJJ_{\uQ}^{\uR}(\M,\T,\sigma,\tus)
\ptf\]
\end{proposition}

\begin{proof}
Par définition
{\multlinegap0pt\begin{multline*}
\begin{aligned}
&\JJJ^T_{\chi}(B,\fff,\omega)=\sum_{\{Q,R\mid \PO\subset Q\subset R\}}\sum_{\{S\mid \PO\subset S\subset\pQ\}}
\frac{1}{n^{\pQ}(S)}\\
&\quad{}\times\sum_{\sigma\in \Pi_{\mathrm{disc}}(M_S)}
\sum_{\WPsi\in \base^{\pQ}_\chi(\sigma)}\sum_{t\in \weyl^\G(\ga_S,Q)}
\sum_{s\in \weyl^Q\left(\theto(\ga_S),t(\ga_S)\right)}\tved(\Q,\R;t)\\
&\qquad{}\times\int_{\ga_Q^G}\tsQR(H)
\biggl(\int_{\ima(\ga_S^\G)^*}\mspace{-10mu}\ee^{\langle(s\theto-t)\mu,H\rangle}\langle
\gmomega^{T,Q}_{s,t}(S,\sigma,f,\omega;\theto\mu,\mu)\WPsi,\WPsi\rangle
B_{\sigma}(\mu)\dd\mu\biggr)\dd H \ptf
\end{aligned}
\end{multline*}}%
D'après le lemme~\ref{bluval}(i) et en remarquant que
\[
\tved(\Q,\R;t)=\tved(\uQ,\uR;\us)
\]
on a
\begin{multline*}
\JJJ^T_{\chi}(B,\fff,\omega)=\sum_{\{S\mid \PO\subset S\}}
\sum_{\sigma\in \Pi_{\mathrm{disc}}(M)_\chi}
\sum_{\{\pQ,\R\mid S\subset \pQ\subset R\}}\frac{1}{n^{\pQ}(S)}\\
{}\times\sum_{\tus\in\weyl^{\tG}(\ga_\M,\ga_\M)}\tved(\uQ,\uR;\us)\JJJ_{\uQ}^{\uR}(\M,\T,\sigma,\tus)\ptf
\end{multline*}
Par symétrisation, en invoquant le lemme~\ref{bluval}(ii),
on obtient la formule souhaitée.
\end{proof}

\begin{lemme}\label{prolong}
Il existe un ensemble de \GM-familles dépendant
de $\tus$\textup, $T$ et d'un paramètre $\mu\in\ima\ga_\M^*$\textup:
\[
(\Lambda,\uS)\mapsto\ebf(\tus,\T,\mu;\Lambda,\uS)
\]
à support compact en $\LL$ et $\mu$ et telles que si
$\Lambda=(\tus-1)\mu$ alors
\[
\AAA_\M^{\uQ}(\T,\sigma,\tus,\mu)B_{\sigma}(\mu)=\ebf_\M^{\uQ}(\tus,\T,\mu;\Lambda)\ptf
\]
\end{lemme}

\begin{proof}
Soit $h$ une fonction lisse à support compact telle que $h(0)=1$.
La \GM-famille
\[
\ebf(\tus,\T,\mu;\Lambda,\uS)
\]
égale, par définition, à
\[
\tr\bigl(\MM(S,\T,\mu;\Lambda,\uS)\phantom{\reg_\mu}
\Mint_{S\rest \tus S}(\mu+\Lambda)\reg_{S,\sigma,\mu}(\tus,\fff,\omega)\bigr)
B_{\sigma}(\mu)h(\Lambda-(\tus-1)\mu)
\]
où $\mu$ et $\Lambda$ sont des variables indépendantes,
est à support compact en $\LL$ et $\mu$\,\footnote{Pour ne pas surcharger les notations nous n'avons pas fait apparaître,
dans
\[
\ebf(\tus,\T,\mu;\Lambda,\uS)
\]
la fonction $B_\sigma$, ni la fonction $h$, ni la donnée $\chi$,
qui resteront fixes dans toute cette section.}.
Avec les notations usuelles on pose
\[
\ebf_\M^{\uQ}(\tus,\T,\mu;\Lambda)=
\sum_{\uS\in\Parab^{\uQ}(\M)}\epsilon_{\uS}^{\uQ}(\Lambda)\ebf(\tus,\T,\mu;\Lambda,\uS)\ptf
\]
Pour conclure on rappelle que, pour $\Lambda=(\tus-1)\mu$, on a
\[
 \AAA_\M^{\uQ}(\T,\sigma,\tus,\mu)=
\tr\bigl(
\MM_{\M}^{\uQ}(S,\T,\mu;\Lambda)
\Mint_{S\rest \tus S}(\mu+\Lambda)\reg_{S,\sigma,\mu}(\tus,\fff,\omega)\bigr)\ptf\qedhere
\]
\end{proof}

Notons $\tL$ le sous-ensemble de Levi minimal contenant l'ensemble $\M\tus$.
Le sous-espace $\gats$ est l'espace des points fixes sous
$\tus$ dans $\gasn$.
On a en particulier
\[
\ga_{\tnL}^{\tG}\subset\ga_{\nL}^\G\subset\ga_\M^\G\ptf
\]
On décompose l'espace $\gasn$ comme somme de $\gats$ et de son orthogonal $\gb_{\tus}$:
\[
\gasn=\gats\oplus\gb_{\tus}
\]
Cette décomposition est stable sous l'action de $(\tus-1)$ et l'espace $\gats$ est l'image de cette application.
Notons $\tL$ le sous-ensemble de Levi minimal contenant l'ensemble $\M\tus$.

\begin{proposition}\label{simpli} 
Il existe une fonction
\[
\UU\mapsto \vf(\tus,\T,\mu;\UU)\]
dans l'espace de Schwartz sur $\gHM$ et à support compact en $\mu$\textup, telle que
\[
\JJJ^T_{\M}(B,\sigma,\fff,\omega,\tus)
\]
soit égal à
\[
\int_{H\in\ga_\M^\G}
\biggl(\int_{\UU\in\gHM}\int_{\ima(\ga_\M^\G)^*}
\ee^{\langle(\tus-1)\mu,H+U_\G\rangle}\,\Gamma_{\tL}^{\tG}(H,\UU)
\vf(\tus,\T,\mu;\UU)\dd\mu\dd\UU\biggr)\dd H\ptf
\]
\end{proposition}

\begin{proof}
On rappelle que, par définition,
\[
\JJJ_{\uQ}^{\uR}(\M,\T,\sigma,\tus)
\]
est l'intégrale itérée
\[
\int_{\ga_{\uQ}^\G}\tsuQuR(H_1)
\biggl(\int_{\ima(\ga_\M^\G)^*} \ee^{\langle(\tus-1)\mu,H_1\rangle}
\AAA_\M^{\uQ}(\T,\sigma,\tus,\mu)
B_{\sigma}(\mu)\dd\mu\biggr)\dd H_1\ptf
\]
On sait que si $\Lambda=(\tus-1)\mu$, on a, d'après le lemme~\ref{prolong}
\[
\AAA_\M^{\uQ}(\T,\sigma,\tus,\mu)B_{\sigma}(\mu)=\ebf_\M^{\uQ}(\tus,\T,\mu;\Lambda)\ptf
\]
D'après le corollaire~\ref{transfour},
il existe une fonction $\vf$
dans l'espace de Schwartz sur $\gHM$
fournissant la \GM-famille par transformation de
Fourier:
\[
\ebf(\tus,\T,\mu;\Lambda,\uS)=\int_{\gHM}\ee^{\Lambda(\UU)}
\vf(\tus,\T,\mu;\UU)\dd\UU\ptf
\]
On sait d'après le lemme~\ref{gmorth} que
\[
\ebf_\M^{\uQ}(\tus,\T,\mu;\Lambda)=\int_{H^1\in\ga_\M^{\uQ}}\int_{\UU\in\gHM}
\ee^{\Lambda(H^1+U_{\uQ})}\,\Gamma_\M^{\uQ}(H^1,\UU)
\vf(\tus,\T,\mu;\UU)\dd H\dd\UU
\]
cette intégrale double étant absolument convergente. En effet, l'intégrale en $H$
de la valeur absolue de $\Gamma_\M^{\uQ}(H,\UU)$
est majorée pas un polynôme en $\UU$, et par ailleurs $\vf$ est à décroissance rapide en $\UU$.
Le résultat est à support compact comme fonction de $\mu$.
Donc
\[
\JJJ_{\uQ}^{\uR}(\M,\T,\sigma,\tus)
\]
est donné par l'intégrale itérée
\[
\int_{H_1\in\ga_{\uQ}^\G}\biggl(\int_{\UU\in\gHM}\int_{H^1\in\ga_\M^{\uQ}}
\int_{\ima(\ga_\M^\G)^*}
F_{\uQ}^{\uR}(\tus;H_1+H^1,\mu,\UU)\dd\mu\dd\UU\dd H^1\biggr)\dd H_1
\]
avec
\[
F_{\uQ}^{\uR}(\tus;H,\mu,\UU)=\ee^{\langle(\tus-1)\mu,H+U_{\uQ}\rangle}\,
\tsuQuR(H)\,\Gamma_\M^{\uQ}(H,\UU)
\vf(\tus,\T,\mu;\UU)
\]
et l'intégrale triple étant absolument convergente.
Ceci peut se récrire comme une intégrale itérée
\[
\JJJ_{\uQ}^{\uR}(\M,\T,\sigma,\tus)=
\int_{H\in\ga_\M^\G}\int_{\UU\in\gHM}
\biggl(\int_{\ima(\ga_\M^\G)^*}
F_{\uQ}^{\uR}(\tus;H,\mu,\UU)\dd\mu\biggr)\dd\UU\dd H\ptf
\]
Pour le voir il faut démontrer la convergence d'une intégrale de la forme
\[
\int_{H\in\ga_\M^\G}\int_{\UU\in\gHM}
\tsuQuR(H)\,\Gamma_\M^{\uQ}(H,\UU)\psi\bigl((\tus^*-1)(H+U_{\uQ}),\UU\bigr)
\dd\UU\dd H
\]
où $\tus^*$ est l'adjoint de $\tus$ et
$\psi(H,\UU)$ est la fonction dans l'espace de Schwartz définie par
\[
\psi(H,\UU)=\int_{\ima(\ga_\M^\G)^*}\ee^{\langle\mu,H\rangle}\,\vf(\tus,\T,\mu;\UU)\dd\mu\ptf
\]
Compte tenu de ce que $\Gamma_\M^{\uQ}(H^1,\UU)$ est à support compact
dans $\ga_\M^{\uQ}$ avec volume polynomial en $\UU$
il suffit d'établir la convergence d'une intégrale du type
\[
\int_{H\in\ga_{\uQ}^\G}\tsuQuR(H)\xi\bigl((\tus^*-1)H\bigr)\dd H
\]
où $\xi$ est une fonction dans l'espace de Schwartz et ceci résulte du lemme~\ref{rebigron}.
On peut alors faire le changement de variable $H_1\mapsto H_1-\U_{\uQ}^\G$ et donc
\[
\JJJ_{\uQ}^{\uR}(\M,\T,\sigma,\tus)=
\int_{H\in\ga_\M^\G}\biggl(\int_{\UU\in\gHM}\int_{\ima(\ga_\M^\G)^*}G_{\uQ}^{\uR}(\tus;H,\mu,\UU)\dd\mu\dd\UU\biggr)\dd H
\]
avec
\[
G_{\uQ}^{\uR}(\tus;H,\mu,\UU)=\ee^{\langle(\tus-1)\mu,H+U_\G\rangle}\,
\tsuQuR(H-U_{\uQ})\,\Gamma_\M^{\uQ}(H,\UU)
\vf(\tus,\T,\mu;\UU)\ptf
\]
Mais
\[
\tved(\uQ,\uR;\us)
\JJJ_{\uQ}^{\uR}(\M,\T,\sigma,\tus)=
\sum_{\{\tP\mid \uQ\subset P\subset\uR,\tus\in\weyl^{\tP}\}}(-1)^{a_{\tP}-a_{\tG}}\,
\JJJ_{\uQ}^{\uR}(\M,\T,\sigma,\tus)\ptf\]
Donc
\begin{equation}
\JJJ^T_{\M}(B,\sigma,\fff,\omega,\tus)=
\sum_{\{\uQ,\uR\mid \M\subset\uQ\subset\uR\}}\tved(\uQ,\uR;\us)
\JJJ_{\uQ}^{\uR}(\M,\T,\sigma,\tus)\label{eq14.1a}
\end{equation}
est égal à l'intégrale itérée en $\mu$, $\UU$ et $H^1$ puis en $H_1$ et de
\begin{equation}
\sum_{\{\tP\mid \M\subset P,\tus\in\weyl^{\tP}\}}
\sum_{\{\uQ,\uR\mid \M\subset\uQ\subset P\subset\uR\}}(-1)^{a_{\tP}-a_{\tG}}\,
G_{\uQ}^{\uR}(\tus;H,\mu,\UU)\ptf\label{eq14.2a}
\end{equation}
Mais, en effectuant d'abord la somme sur $\uR$ on voit que, d'après le lemme~\ref{repart},
\begin{equation}
\sum_{\{\tP\mid \M\subset P,\tus\in\weyl^{\tP}\}}
\sum_{\{\uQ,\uR\mid \M\subset\uQ\subset P\subset\uR\}}
(-1)^{a_{\tP}-a_{\tG}}\,\tsuQuR(H-U_{\uQ})\,\Gamma_\M^{\uQ}(H,\UU)\label{eq14.3a}
\end{equation}
est égal à
\[
\sum_{\{\tP\mid \M\subset P,\tus\in\weyl^{\tP}\}}
\sum_{\{\uQ\mid \M\subset\uQ\subset P\}}
(-1)^{a_{\tP}-a_{\tG}}\,\htau_{\tP}(H-U_{\tP})\tau_{\uQ}^P(H-U_{\uQ})\,\Gamma_\M^{\uQ}(H,\UU)
\]
qui d'après le lemme~\ref{decomp}\eqref{eq1.3c}, par sommation sur $\uQ$, est égal à
\begin{equation}
\sum_{\{\tP\mid \M\subset P,\tus\in\weyl^{\tP}\}}
(-1)^{a_{\tP}-a_{\tG}}\,\,\htau_{\tP}(H-U_{\tP})\label{eq14.5a}
\end{equation}
qui, à son tour, d'après la proposition~\ref{envconvt} est égal à
\begin{equation}
\Gamma_{\tL}^{\tG}(H,\UU)\ptf\label{eq14.6a}
\end{equation}
L'égalité de~\eqref{eq14.3a} et~\eqref{eq14.6a} montre que~\eqref{eq14.1a} est égal à l'intégrale itérée de
\begin{equation}
\ee^{\langle (\tus-1)\mu,H+U_\G\rangle}\,\Gamma_{\tL}^{\tG}(H,\UU)\vf(\tus,\T,\mu;\UU)\ptf\qedhere\label{eq14.7a}
\end{equation}
\end{proof}

Posons
\[
\MM_{\tnL}^{\tG}(S,\T,\nu)=\MM_{\tnL}^{\tG}(S,\T,\nu;0)\ptf
\]

\begin{proposition}\label{final}
L'expression
\begin{equation}
\JJJ^T_{\M}(B,\sigma,\fff,\omega,\tus)\label{eq14.1b}
\end{equation}
est égale au produit de
\[
\frac{1}{\lvert\det(\tus-1\rest  \ga_{\M}^\G /\gats)\rvert}
\]
et de
\[
\int_{\ima(\gats)^*}
\tr\bigr(\MM_{\tnL}^{\tG}(S,\T,\nu)\Mint_{S\rest \tus S}(0)
\reg_{S,\sigma,\nu}(\tus,\fff,\omega)B_\sigma(\nu)\bigr)\dd\nu\ptf
\]
\end{proposition}

\begin{proof} 
On sait d'après la proposition~\ref{simpli} que l'expression
\eqref{eq14.1b} est égale à la triple intégrale itérée
\[
\int_{H\in\ga_\M^\G}\biggl(\int_{\UU\in\gHM}
\int_{\ima(\ga_\M^\G)^*}\ee^{\langle(\tus-1)\mu,H+U_\G\rangle}\,\Gamma_{\tL}^{\tG}(H,\UU)
\vf(\tus,\T,\mu;\UU)\dd\mu\dd\UU\biggr)\dd H\ptf
\]
On rappelle que l'on peut écrire
\[
\mu=\nu+\eta
\qquad\text{avec $\nu\in \ima(\gats)^*$ et $\eta\in (i\gb_{\tus})^*$}
\]
et on pose comme d'habitude
\[
\Lambda=(\tus-1)\mu=(\tus-1)\eta\ptf
\]
Le changement de variable $(\nu,\Lambda)\mapsto\mu$ dans une intégrale sur $\mu$ s'écrit:
\[
\int_{\mu\in\ima\gasns}\phi(\mu)\dd\mu=\frac{1}{\lvert\det(\tus-1\rest  \ga_{\M} /\gats)\rvert}
\int_{\nu\in i(\gats)^*}\int_{\Lambda\in \ima(\gb_{\tus})^*}
\phi\bigl(\mu(\nu,\Lambda)\bigr)\dd\nu\dd\Lambda\ptf
\]
Maintenant, on peut aussi décomposer $H\in\gasn$ en
\[
H=X+Y\in\gb_{\tus}\oplus\gats\ptf
\]
On en déduit que la triple intégrale itérée peut encore s'écrire comme le produit de
\[
\frac{1}{\lvert\det(\tus-1\rest  \ga_{\M}^\G /\gats)\rvert}
\]
et de l'intégrale en $X\in\gb_{\tus}$ de
\begin{multline*}
\smash[b]{\int_{\UU\in\gHM}\int_{Y\in\gats}
\int_{\nu\in \ima(\gats)^*}\int_{\Lambda\in i(\gb_{\tus})^*}}
\ee^{\langle \LL,X+Y+U_\G\rangle}\Gamma_{\tL}^{\tG}(Y,\UU)\\ 
{}\times
\vf(\tus,\T,\mu(\nu,\LL);\UU)\dd\LL \dd\nu\dd\UU \dd Y
\end{multline*}
qui est encore égale, d'après les lemmes~\ref{gmorth} et~\ref{tgmfour},
à l'intégrale itérée
\[
\int_{X\in\gb_{\tus}}\biggl(\int_{\nu\in \ima(\gats)^*}\int_{\Lambda\in \ima(\gb_{\tus})^*}
\ee^{\langle\LL,X\rangle}\,\ebf_{\tL}^{\tG}(\tus,\T,\nu;\LL)\dd\LL\dd\nu\biggr)\dd X
\]
qui, par inversion de Fourier, se récrit
\[
\int_{\ima(\gats)^*} \ebf_{\tL}^{\tG}(\tus,\T,\nu;0)\dd\nu
\]
soit encore, par définition de la \GM-famille $\ebf$:
\[
\int_{\ima(\gats)^*}
\tr\bigl(\MM_{\tnL}^{\tG}(S,\T,\nu)\Mint_{S\rest \tus S}(\nu)
\reg_{S,\sigma,\nu}(\tus,\fff,\omega)B_\sigma(\nu)\bigr)\dd\nu\ptf
\]
Enfin on observe de plus que
\[
\Mint_{S\rest \tus S}(\nu)=\Mint_{S\rest \tus S}(0)
\]
lorsque $\nu\in{i(\gats)^*}$.
\end{proof}

En résumé on a obtenu la

\begin{proposition}\label{maintha}
\[
\JJJ^\T_\chi(B,\fff,\omega)
=\sum_{\M\in\mathcali L^\G/\weyl^\G}
\frac{1}{w^\G(\M)} \sum_{\sigma\in \Pi_{\mathrm{disc}}(M)_\chi}\sum_{\tus\in\weyl^{\tG}(\M)}
\JJJ^\T_{\M,\sigma,\tus}(B,\fff,\omega)
\]
avec
{\multlinegap0pt
\begin{multline*}
\JJJ^\T_{\M,\sigma,\tus}(B,\fff,\omega)\\ =
\frac{1}{\lvert\det(\tus-1\rest  \ga_{\M}^\G /\gats)\rvert} 
\int_{\ima(\gats)^*}\mspace{-10mu}
\tr\bigl(\MM_{\tnL}^{\tG}(S,\T,\nu)\Mint_{S\rest \tus S}(0)
\reg_{S,\sigma,\nu}(\tus,\fff,\omega)B_\sigma(\nu)\bigr)\dd\nu\ptf
\end{multline*}}%
où on note $\tL$ le sous-ensemble de Levi minimal contenant l'ensemble $\M\tus$.
\end{proposition}

\begin{lemme}\label{polygm}
La fonction
\[
\T\mapsto\MM_{\tnL}^{\tG}(S,\T,\nu)
\]
est un polynôme à valeurs opérateurs.
\end{lemme}

\begin{proof}
La \GM-famille
\[
\MM(S,\T,\nu;\Lambda,\uS)=\ee^{\langle\Lambda,Y_{\uS}(T)\rangle}\,
\Mint_{\uS\rest S}(\nu)\moins\Mint_{\uS\rest S}(\nu+\Lambda)
\]
est un produit de deux \GM-familles et d'après la proposition~\ref{GammaHXt}
on a
\[
\MM_{\tnL}^{\tG}(S,\T,\nu;0)=\sum_{\tP\in\FF(\tL)}c_{\tL}^{\tP}(S,\T;0)\,d_{\tP}^{\tG}(S,\nu;0)
\]
avec
\[
c(S,\T;\LL,\uS)=\ee^{\langle\Lambda,Y_{\uS}(T)\rangle}
\]
mais le lemme~\ref{gmorth} montre que
\[
c_{\tL}^{\tP}(S,\T;0)=\int_{\ga_{\tL}^{\tP}} \Gamma_{\tL}^{\tP}\bigl(H,\mathcali Y(\T)\bigr)\dd H
\]
qui est un polynôme en $\T$ d'après le lemme~\ref{LaplaceGamma}.
\end{proof}

\begin{proposition}\label{mainthb}
Supposons $B(0)=1$\textup, alors il existe $c>0$ tel que\textup,
si $\dPO(T)\geq c$\textup, on a
{\multlinegap0pt\begin{multline*}
J^T_{\chi}(\fff,\omega)=\lim_{\epsilon\to 0}\sum_{\M\in\mathcali L^\G/\weyl^\G}
\frac{1}{w^\G(\M)} \sum_{\sigma\in \Pi_{\mathrm{disc}}(M)_\chi}\sum_{\tus\in\weyl^{\tG}(\M)}
\frac{1}{\lvert\det(\tus-1\rest  \ga_{\M}^\G /\gats)\rvert}\\
{}\times\int_{\ima(\gats)^*}
\tr\bigl(\MM_{\tnL}^{\tG}(S,\T,\nu)\Mint_{S\rest \tus S}(0)
\reg_{S,\sigma,\nu}(\tus,\fff,\omega)B^\epsilon_\sigma(\nu)\bigr)\dd\nu\ptf
\end{multline*}}%
\end{proposition}

\begin{proof} 
Rappelons que d'après le théorème~\ref{thWA},
il existe un polynôme $p^\T_\chi(B,\fff,\omega)$ en $T$ tel que
\[
\lim_{\dPO(T)\to \infty} \bigl(J_{\chi}^T(B,\fff,\omega)-p^\T_\chi(B,\fff,\omega)\bigr)=0\ptf
\]
C'est dire que $p^\T_\chi(B,\fff,\omega)$ et $J_{\chi}^T(B,\fff,\omega)$
sont asymptotes, quand $\dPO(T)$ tend vers l'infini. De plus, la proposition~\ref{2.10bis} montre que
$J_{\chi}^T(B,\fff,\omega)$
et $\JJJ_{\chi}^T(B,\fff,\omega)$ sont également asymptotes.
Mais $\JJJ_{\chi}^T(B,\fff,\omega)$ est aussi un polynôme d'après la proposition~\ref{maintha} et le lemme~\ref{polygm}.
Maintenant, deux polynômes asymptotes sont nécessairement égaux:
\[
p^\T_\chi(B,\fff,\omega)=\JJJ^\T_\chi(B,\fff,\omega)\ptf
\]
Enfin, d'après le théorème~\ref{thWA}, si nous
supposons $B(0)=1$, alors il existe $c>0$ tel que
\[
J^\T_\chi(\fff,\omega)=\lim_{\epsilon\to 0}\,\JJJ^\T_\chi(B^{\epsilon},\fff,\omega)
\]
si $\dPO(T)\geq c$.
\end{proof}

\section[Élimination de la fonction $B$]{\mathversion{bold}Élimination de la fonction $B$}\label{sec14.2}

\begin{theoreme}\label{mainthc}
{\multlinegap0pt\begin{multline*}
J^\T_\chi(\fff,\omega)=\sum_{\M\in\mathcali L^\G/\weyl^\G}
\frac{1}{w^\G(\M)} \sum_{\sigma\in \Pi_{\mathrm{disc}}(M)_\chi}\sum_{\tus\in\weyl^{\tG}(\M)}
\frac{1}{\lvert\det(\tus-1\rest  \ga_{\M}^\G /\gats)\rvert}\\
{}\times\int_{\ima(\gats)^*}
\tr\bigl(\MM_{\tnL}^{\tG}(S,\T,\nu)\Mint_{S\rest \tus S}(0)
\reg_{S,\sigma,\nu}(\tus,\fff,\omega)\bigr)\dd\nu\ptf
\end{multline*}}%
\end{theoreme}

\begin{proof}
La normalisation des opérateurs d'entrelacement
a été établie pour la première fois par Langlands dans \cite{MS}*{Lecture~15},
puis reprise par Arthur dans \cite{Ainter}. Ceci étant établi
on peut reprendre essentiellement mot à mot la preuve d'Arthur
dans les sections~6 à~9 de \cite{AeisII}
(qui elle était conditionnelle à l'existence d'une telle normalisation)
pour montrer que l'expression
{\multlinegap0pt\begin{multline*}
\sum_{\M\in\mathcali L^\G/\weyl^\G}
\frac{1}{w^\G(\M)} \sum_{\sigma\in \Pi_{\mathrm{disc}}(M)_\chi}\sum_{\tus\in\weyl^{\tG}(\M)}
\frac{1}{\lvert\det(\tus-1\rest  \ga_{\M}^\G /\gats)\rvert}\\
{}\times\int_{\ima(\gats)^*}
\tr\bigl(\MM_{\tnL}^{\tG}(S,\T,\nu)\Mint_{S\rest \tus S}(0)
\reg_{S,\sigma,\nu}(\tus,\fff,\omega)\bigr)\dd\nu\ptf
\end{multline*}}%
est absolument convergente.
La seule étape non évidente est
l'extension au cas tordu des résultats de la section~7 de \cite{AeisII}.
Cela fait l'objet du corollaire~\ref{gmrad}.
Le théorème de convergence dominée montre alors que, si $B(0)=1$
cette expression est la limite pour $\epsilon\to0$ de
{\multlinegap0pt\begin{multline*}
\sum_{\M\in\mathcali L^\G/\weyl^\G}
\frac{1}{w^\G(\M)} \sum_{\sigma\in \Pi_{\mathrm{disc}}(M)_\chi}\sum_{\tus\in\weyl^{\tG}(\M)}
\frac{1}{\lvert\det(\tus-1\rest  \ga_{\M}^\G /\gats)\rvert}\\
{}\times\int_{\ima(\gats)^*}
\tr\bigl(\MM_{\tnL}^{\tG}(S,\T,\nu)\Mint_{S\rest \tus S}(0)
\reg_{S,\sigma,\nu}(\tus,\fff,\omega)B^\epsilon_\sigma(\nu)\bigr)\dd\nu\ptf
\end{multline*}}%
Si nous supposons $\dPO(T)\geq c$, l'égalité cherchée résulte alors de la proposition \ref{mainthb}.
On observe enfin que les deux membres sont des polynômes. L'égalité est donc toujours vraie.
\end{proof}

\section{Développement spectral fin}\label{spectral}

Le développement spectral grossier de la formule des traces est, par définition,
la valeur en $T=\TK$ de la série des $J^\T_\chi$:
\[
J^{\tG}(\fff,\omega)=\sum_\chi J^{\TK}_\chi(\fff,\omega)\ptf
\]
En le combinant avec développement spectral des termes $J^{\TK}_\chi$
on obtient le développement spectral fin.
On peut le formuler au moyen de la \GM-famille
\[
\MM(S,\nu;\Lambda,\uS)=
\Mint_{\uS\rest S}(\nu)\moins\Mint_{\uS\rest S}(\nu+\Lambda)
\]
qui donne naissance à l'opérateur
$\MM_{\tnL}^{\tG}(S,\nu;\LL)$ et on note $\MM_{\tnL}^{\tG}(S,\nu)$ sa valeur en $\LL=0$.

\begin{theoreme}\label{mainthd}
\[
J^{\tG}(\fff,\omega)=\sum_{\M\in\mathcali L^\G/\weyl^\G}
\frac{1}{w^\G(\M)} J^{\tG}_\M(\fff,\omega)
\]
avec
\[
J^{\tG}_\M(\fff,\omega)= \sum_{\sigma\in \Pi_{\mathrm{disc}}(M)}\sum_{\tus\in\weyl^{\tG}(\M)}
\frac{1}{\lvert\det(\tus-1\rest  \ga_{\M}^\G /\gats)\rvert}J^{\tG}_{\M,\sigma}(\fff,\omega,\tus)
\]
et
\[
J^{\tG}_{\M,\sigma}(\fff,\omega,\tus)=
\int_{\ima(\gats)^*}
\tr\bigl(\MM_{\tnL}^{\tG}(S,\nu)\Mint_{S\rest \tus S}(0)\reg_{S,\sigma,\nu}(\tus,\fff,\omega)\bigr)\,d\nu\ptf
\]
\end{theoreme}

\begin{proof}
Il résulte de la proposition~\ref{mainthc} que
\[
J^{\tG}(\fff,\omega)=\sum_\chi\Biggl(\sum_{\M\in\mathcali L^\G/\weyl^\G}
\frac{1}{w^\G(\M)} J^{\tG}_{\M,\chi}(\fff,\omega)\Biggr)
\]
avec
{\multlinegap0pt\begin{multline*}
J^{\tG}_{\M,\chi}(\fff,\omega)= \sum_{\sigma\in \Pi_{\mathrm{disc}}(M)_\chi}\sum_{\tus\in\weyl^{\tG}(\M)}
\frac{1}{\lvert\det(\tus-1\rest  \ga_{\M}^\G /\gats)\rvert}\\
{}\times\int_{\ima(\gats)^*}
\tr\bigl(\MM_{\tnL}^{\tG}(S,\TK,\nu)\Mint_{S\rest \tus S}(0)
\reg_{S,\sigma,\nu}(\tus,\fff,\omega)\bigr)\dd\nu\ptf
\end{multline*}}%
Maintenant d'après l'équation~\eqref{eq5.4b} de la proposition~\ref{GMspec} on a
\[
\MM_{\tnL}^{\tG}(S,\TK,\nu)=\MM_{\tnL}^{\tG}(S,\nu)\ptf
\]
De plus, grâce aux travaux récents de
Finis, Lapid et Müller \citelist{\cite{FL}\cite{FLM}},
on sait maintenant que le développement spectral
est absolument convergent. Leurs travaux ne concernent que le cas classique
(non tordu) mais ils s'étendent sans modification au cas général.
On peut donc omettre les sommations partielles suivant les~$\chi$.
\end{proof}

Notons $\treg_{\mathrm{disc}}$ la restriction de $\treg$ au spectre discret pour $\G$.
On pose
\[
J_{\G,\mathrm{disc}}^{\tG}(\fff,\omega)=\tr\treg_{\mathrm{disc}}(\fff,\omega)=
\sum_{\pi\in\Pi_{\mathrm{disc}}(\tG,\omega)}\mtG(\pi,\tpi) \tr\tpi(\fff,\omega)
\]
où $\Pi_{\mathrm{disc}}(\tG,\omega)$ est l'ensemble des classes d'équivalence de
représentations irréductibles $\pi$ qui sont les
restrictions à $\Gadef $ de représentations $\tpi$ de $(\Gadef,\omega)$ qui
interviennent dans le spectre discret
\[
L^2_{\mathrm{disc}}(\XG)\ptf
\]
On sait que pour $\delta\in\tGadef$ on a
\[
\pi\circ\theta\simeq\pi\otimes\omega\qquad\text{avec $\theta=\tAd(\delta)$}\ptf
\]
Enfin $\mtG(\pi,\tpi)$ est la multiplicité tordue de
$\tpi$ dans le spectre discret. Rappelons que c'est un
torseur à valeurs dans $\CM^\times$;
cette notion de multiplicité tordue a été discutée dans la section~\ref{multi}.

Nous avons omis la sommation partielle --- utilisée chez Arthur --- suivant les
modules des caractères infinitésimaux à l'infini désormais inutile
puisque, comme observé plus haut, nous savons que le développement spectral est absolument convergent.
On remarquera que pour le spectre discret il suffit d'invoquer~\cite{Mu}.

La partie discrète du développement spectral de la formule des traces est une distribution
\[
J_{\mathrm{disc}}^{\tG}
\]
qui
est une somme de termes
parmi lesquels on a $J_{\G,\mathrm{disc}}^{\tG}$ la trace dans le spectre discret.
Cependant, d'autres termes discrets, \cad ne faisant pas apparaître d'intégrale dans leur expression,
quoique provenant du spectre continu, contribuent à l'expression spectrale de la formule des traces;
nous allons les décrire.
Soit $M\in\Levi^\G$ un sous-groupe de Levi de $G$.
On notera $\weyl^{\tG}(M)_{\mathrm{reg}}$ le sous-ensemble des $\tus\in \weyl^{\tG}(M)$ réguliers, \cad tels que
\[
\det(\tus-1\rest \ga_M/\ga_\G)\ne0\ptf
\]
On pose
\[
J_{M,\mathrm{disc}}^{\tG}(\fff,\omega)=\sum_{\tus\in \weyl^{\tG}(M)_{\mathrm{reg}}}
\frac{1}{\lvert\det(\tus-1\rest \ga_M/\ga_\G)\rvert}
\tr\bigl(\Mint_{S\rest \tus(S)}(0)\reg_{S,\mathrm{disc},0}(\tus,\fff,\omega)\bigr)
\]
où $S$ est un sous-groupe parabolique de Levi $\M$.
L'expression est indépendante du choix de $S$.

\begin{proposition}\label{spectrediscret}
La partie discrète de la formule des traces peut s'écrire\text:
\[
J_{\mathrm{disc}}^{\tG}(\fff,\omega)=\sum_{M\in\Levi^\G} \frac{\lvert\weyl^M\rvert}{\lvert\weyl^G\rvert}\,
J_{M,\mathrm{disc}}^{\tG}(\fff,\omega)
\]
ou si on préfère
\[
J_{\mathrm{disc}}^{\tG}(\fff,\omega)=\sum_{M\in\Levi^\G/\weyl^G}\frac{1}{\lvert\weyl^G(M)\rvert}\,
J_{M,\mathrm{disc}}^{\tG}(\fff,\omega)
\]
où la somme porte sur un ensemble de représentants
des orbites de $\weyl^G$ dans $\Levi^\G$.
\end{proposition}

\begin{proof}
Les termes discrets sont ceux qui dans le théorème~\ref{mainthd}
ne font pas apparaître d'intégrale, \cad les termes où $\gats$ est réduit à $0$.
Ce sont donc ceux pour lesquels $\tus$ est régulier.
\end{proof}

Plus généralement, posons
\[
J^{\tG}_{\M,\tus}(\fff,\omega)=
\int_{\ima(\gats)^*}
\tr\bigl(\MM_{\tnL}^{\tG}(S,\nu)\Mint_{S\rest \tus S}(0)
\reg_{S,\mathrm{disc},\nu}(\tus,\fff,\omega)\bigr)\dd\nu\ptf
\]
Soit $\tL$ un sous-espace de Levi semi-standard (\ie $\MO\subset L$). Définissons
\[ 
J_{\tL}^{\tG}(\fff,\omega)=\sum_{\M\in\Levi ^L}
\frac{\lvert\weyl^\M\rvert}{\lvert\weyl^L\rvert}\sum_{\tus\in\weyl^{\tL}(\M)_{\mathrm{reg}}}\frac{1}{\lvert\det(\tus-1\rest  \ga_\M^L)\rvert}
J^{\tG}_{\M,\tus}(\fff,\omega)\ptf
\]
On remarquera que
\[
J_{\tG}^{\tG}(\fff,\omega)=J_{\mathrm{disc}}^{\tG}(\fff,\omega)\ptf
\]
On note $\td\weyl^L$ le quotient par $\MO$ du normalisateur de $\tMO$ dans $L$.
Soit $\Levi ^{\tG}$ l'ensemble des $\tL$ contenant $\tMO$.
On note enfin $\theta_L$ l'automorphisme induit sur $\ga_L$
par un quelconque élément $\tus$ de $\tL(F)$.
Avec ces notations, on a le théorème suivant:

\begin{theoreme}\label{finale}
\[
J^{\tG}(\fff,\omega)=\sum_{\tL\in\Levi ^{\tG}}
\frac{\lvert\td\weyl^L\rvert}{\lvert\td\weyl^\G\rvert}
\frac{1}{\lvert\det(\theta_L-1\rest \ga_L^\G/\ga_{\tL}^{\tG})\lvert}J_{\tL}^{\tG}(\fff,\omega)\ptf
\]
\end{theoreme}

\begin{proof} 
Il suffit d'observer que l'on peut
écrire $J^{\tG}_\M(\fff,\omega)$ sous la forme
\[
J^{\tG}_\M(\fff,\omega)=\sum_{\tus\in\weyl^{\tG}(\M)}
\frac{1}{\lvert\det(\tus-1\rest  \ga_{\M}^\G /\gats)\rvert}J^{\tG}_{\M,\tus}(\fff,\omega)
\]
et que si $\tL$ est défini au moyen de $\tus$ on a $\tus\in\weyl^{\tL}(\M)_{\mathrm{reg}}$,
et enfin que
\[
\det(\tus-1\rest  \ga_\M^\G /\ga_L^\G)=\det(\tus-1\rest  \ga_\M^L)\ptf
\]
\end{proof}

En vue de la stabilisation de la formule des traces tordue il est utile
de reformuler ce théorème en renormalisant les distributions comme suit.
On pose
\[
\JJ_{\tL}^{\tG}=\jtL\moins J_{\tL}^{\tG}
\Qquad{et}
\JJ^{\tG}=j(\tG)\moins J^{\tG}\]
avec
\[
\jtL=\lvert\det(\theta_L-1\rest \ga_L/\ga_{\tL})\rvert\ptf
\]
Le théorème~\ref{finale} se récrit alors

\begin{corollaire}
\[
\JJ^{\tG}(\fff,\omega)=\sum_{\tL\in\Levi ^{\tG}}
\frac{\lvert\td\weyl^L\rvert}{\lvert\td\weyl^\G\rvert}
\JJ_{\tL}^{\tG}(\fff,\omega)\ptf
\]
\end{corollaire}

\backmatter

\renewcommand{\indexname}{Index des notations}
\begin{theindex}

  \item $\ll $, 55

  \indexspace

  \item $\mathfrak  {A}_G$, 3
  \item $a^G(\delta )$, 117
  \item $\mathfrak  {a}_P$, 3
  \item $a_P^Q$, 4
  \item $a(s,\chm )$, 17
  \item $A_{s,t}^T$, 198
  \item $\mathbf  {A}_{s,t}^T$, 198

  \indexspace

  \item $\Delta (M,s)$, 12
  \item $\delta _0$, 44
  \item $\Delta _P^Q$, 5
  \item $\Delta _{\tP }^{\tQ }$, 45
  \item $\hDelta _P^Q$, 5
  \item $\mathbf  {d}_{P_0}$, 8

  \indexspace

  \item $\tve(Q,R)$, 133
  \item $\epsilon _P^Q$, 28
  \item $\mathaccentV {hat}05E\epsilon _P^Q$, 28
  \item $\mathaccentV {widetilde}672{\eta }(Q,R)$, 144
  \item $\mathaccentV {widetilde}672{\eta }(Q,R;t)$, 211

  \indexspace

  \item $\mathaccentV {widetilde}672{f}$, 117
  \item $f^1$, 105
  \item $F_{P_0}^{Q}$, 72
  \item $\mathcal  {F}^Q(M)$, 12

  \indexspace

  \item $G$, 3
  \item $\tG $, 43
  \item $G(\mathbb  {A})^1$, 3
  \item $\Gamma _M^Q(H,\mathcal  {X})$, 24
  \item $\Gamma _M^Q(H,\mathcal  {X},\chm )$, 17
  \item $\Gamma _P^Q(H,X)$, 23
  \item $\gamma _P^R(\Lambda ,X)$, 28

  \indexspace

  \item $\mathbf  {H}_0$, 62
  \item $\mathbf  {H}_{G}$, 3

  \indexspace

  \item $j(\mathaccentV {widetilde}672{G})$, 44

  \indexspace

  \item $K_{\tG }(x,y)$, 105
  \item $K_{\tP }$, 120

  \indexspace

  \item $\Lambda ^{T,Q}$, 81
  \item $\mathcal  {L}^Q(M)$, 12

  \indexspace

  \item $\Mint(s,\lambda )$, 95
  \item $M_0$, 44
  \item $\mathbf  {M}_{Q\rest{}P}(s,\lambda)$, 94

  \indexspace

  \item $\mathcal  {O}_\delta (f,\omega )$, 117
  \item $\omega ^{T,Q}_{s,t}$, 198

  \indexspace

  \item $\tP $, 45
  \item $P_0$, 44
  \item $\phi _{M,s}^\chm $, 17
  \item $\phi _P^Q$, 21
  \item $\phi _P^{Q,R}$, 21
  \item $\mathcal  {P}^Q(M)$, 12

  \indexspace

  \item $Q^{+}$, 52
  \item $Q^s$, 13
  \item $Q_s$, 13

  \indexspace

  \item $\mathcal  {R}$, 4
  \item $\mathcal  {R}(s)$, 8
  \item $\mathcal  {R}(s,t)$, 8
  \item $R^{-}$, 52
  \item $\treg (f,\omega )$, 103
  \item $\boldsymbol  {\rho }_{P,\sigma ,\mu }(\delta ,f,\omega )$, 104
  \item $\mathaccentV {widetilde}672{\boldsymbol  {\rho }}_{P,\sigma ,\mu }(f,\omega )$, 
		104

  \indexspace

  \item $\mathbf  {s}$, 94
  \item $\sigma _{Q}^{R}$, 53
  \item $\tsQR $, 53
  \item ${}_tS$, 197
  \item $\boldsymbol  {\mathfrak  {S}}_{t,\Omega }$, 68

  \indexspace

  \item $\htau _P^Q$, 19
  \item $\tau _P^Q$, 19
  \item $\theta _0$, 44

  \indexspace

  \item $\mathbf  {W}$, 8
  \item $\mathbf  {W}(\mathfrak  {a}_M)$, 12
  \item $\mathbf  {W}(\mathfrak  {a}_P,\mathfrak  {a}_Q)$, 10
  \item $\mathbf  {W}(\mathfrak  {a}_P,R)$, 11

  \indexspace

  \item $X_F$, 3
  \item $\mathbf  {X}_{G}$, 68
  \item $\mathbf  {X}_{P,G}$, 68

  \indexspace

  \item $\mathbf  {Y}_P$, 68
  \item $Y_S$, 97
  \item $Y_s$, 97
  \item $Y_S(T)$, 67
  \item $Y_s(T)$, 67
  \item $Y_s(x,T)$, 66

\end{theindex}

\printindex

\appendix\bigskip
\newcount\Cor\Cor=0
\def\newcor{\global\advance\Cor by 1
\par\medskip\noindent (\romannumeral\Cor) - }

\def\obs#1{{\color{blue}#1}}
\def\footobs#1{\footnote{\obs{#1}}}
\def\add#1{{\color{magenta}#1}}
\def\supp#1{{\color{orange}#1}}

\catcode`\@=11
\catcode`\Ž=13
\catcode`\=13
\catcode`\ˆ=13
\catcode`\=13
\catcode`\=13
\catcode`\‰=13
\catcode`\™=13
\catcode`\"=13
\catcode`\•=13
\catcode`\ž=13
\catcode`\=13
\catcode`\'=13
\catcode`\Ë=13
\def Ž{\'e}
\def {\`e}
\def ˆ{\`a}
\def {\`u}
\def {\^e}
\def ‰{\^a}
\def ™{\^o}
\def "{\^\i}
\def •{\"\i}
\def ž{\^u}
\def {\c c}
\def '{\"e}
\def Ë{\`A}
\catcode`\:=13
\def :{~\string:}
\catcode`\;=13
\def ;{~\string;}
\catcode`\!=13
\def !{~\string!}
\def\cad{c'est-\`a-dire\ }

\def\xrightarrow#1{\,\,\smash{\mathop{-\!\!\!-\!\!\!\longrightarrow}^#1}\,\,\,}
\def\bydef{\buildrel {\rm d\acute{e}f}\over{=}}

\def\ES#1{\EuScript{#1}}
\def\ck#1{\check{#1}}
\def\wt#1{\widetilde{#1}}
\def\bs#1{\boldsymbol{#1}}
\def\wh{\widehat}
\def\mun{^{-1}}

\def\Lbra{[\![}
\def\Rbra{]\!]}

\def\ni{\noindent}
\def\pni{\par\noindent}
\def\ptf{\,.}
\def\vg{\,,}
\def\vgq{\,,\quad}
\def\pv{\,;}
\def\pvq{\,;\quad}

\def\A{A}
\def\G{G}
\def\L{L}
\def\M{M}
\def\Y{Y}

\def\AM{\mathbb A}
\def\CM{\mathbb C}
\def\FM{\mathbb F}
\def\GM{\mathbb G}
\def\NM{\mathbb N}
\def\RM{\mathbb R}
\def\QM{\mathbb Q}
\def\UM{\mathbb U}
\def\ZM{\mathbb Z}
\def\bsbbc{\mathbbm{c}}
\def\mbmb{\mathbbm{b}}

\def\bfHom{{\rm Hom}}
\def\vol{{\rm vol}}

\def\AA{\mathfrak A}
\def\BB{\mathfrak B}
\def\EE{\mathfrak E}
\def\FF{\mathfrak F}
\def\HH{\mathfrak H}
\def\LL{\mathfrak L}
\def\NN{\mathfrak N}
\def\OO{\mathfrak O}
\def\UU{\mathfrak U}
\def\VV{\mathfrak V}
\def\TT{\mathfrak T}
\def\XX{\mathfrak X}
\def\YY{\mathfrak Y}
\def\ZZ{\mathfrak Z}

\def\Siegel{\bs{\mathfrak{S}}}

\def\ag{\mathfrak{a}}
\def\agp{\ag_{P}}
\def\cc{\mathfrak c}
\def\uu{\mathfrak u}
\def\oo{\mathfrak o}
\def\pp{\mathfrak p}

\def\bfA{{\bf A}}
\def\bfC{{\bf C}}
\def\bfD{{\bf D}}
\def\bfH{{\bf H}}
\def\bfM{{\bf M}}
\def\bfW{{\bf W}}

\def\Rho{\bs\rho}
\def\tRho{\wt{\Rho}}
\def\bsmu{\bs{\mu}}
\def\bsPi{\bs{\Pi}}
\def\bsigma{\bs{\sigma}}
\def\tbsrho{\wt{\bs{\rho}}}

\def\bsA{\bs{A}}
\def\bsJ{\bs{J}}
\def\bsK{\bs{K}}
\def\bsS{\bs{S}}
\def\bsX{\bs{X}}
\def\ovbsX{\overline{\bsX}}
\def\bsY{\bs{Y}}
\def\bsc{\bs{c}}
\def\bsd{\bs{d}}
\def\bscMF{\bsc_{M,F}}
\def\bscM{\bsc_{M}}
\def\bsO{\boldsymbol{\Omega}}
\def\bso{\boldsymbol{\omega}}
\def\brabsO{[\bsO]}
\def\brabso{[\bso]}
\def\bchi{\boldsymbol{\chi}}
\def\bESD{\boldsymbol{\ESD}}
\def\bESP{\boldsymbol\ESP}
\def\bPI{\boldsymbol{\Pi}}

\def\bsmu{\boldsymbol{\mu}}

\def\fun{f}
\def\tdf{f}
\def\tG{{\widetilde G}}
\def\tH{{\widetilde H}}
\def\tL{{\widetilde L}}
\def\tM{{\widetilde M}}
\def\tP{{\widetilde P}}
\def\tQ{{\widetilde Q}}
\def\tR{{\widetilde R}}
\def\tS{{\widetilde S}}
\def\tZ{{\widetilde Z}}
\def\tu{{\widetilde u}}

\def\Dcal{\mathcal{D}}

\def\brT_#1{[T]_{#1}}
\def\brTo_#1{[T_0]_{#1}}
\def\brTX_#1{[T-X]_{#1}}

\def\iipi{2i\pi}
\def\iipiZ{\iipi\ZM}
\def\hag{\wh{\ag}}
\def\hagp{\wh{\ag}_P}
\def\hagpf{\wh{\ESA}_P}

\def\logq{\log_q}

\def\bsl{\backslash}
\def\dd{\,{\mathrm d}}

\def\cM{|\wh\bsbbc_M|}
\def\Stab{{\rm Stab}_M(\sigma)}
\def\stab{\vert \Stab\vert}
\def\unstab {\frac{\cM}{\stab}}
\def\frGM{\frac{1}{w^G(M)}}
\def\frQpM{\frac{1}{w^{Q'}(M)}}

\def\gammaMF{\gamma_{M,F}}
\def\overag{\bar{\ag}}

\def\Jres{\mathfrak J}
\def\Jvec{J}
\def\bJres{\bs{\Jres}}
\def\bJvec{\bs{Jvec}}

\def\ttM{{_tM}}
\def\ttS{{_tS}}
\def\mf{m}

\def\ZTX{H_Z^{T-X}} 
\def\ZX{H_Z^X} 
\def\ZTSs{H^T_{Z,S''}}

\def\nP{P} 
\def\Pp{P'}
\def\Qo{{Q_0}}

\def\Kappa#1{\kappa^{#1}}
\def\kopa{\chi}
\def\cMsig{\wh c_M(\sigma)}
\def\cMSsig{\wh{c}_{M_S}(\sigma)}
\def\bsESEsigma{\bs{\ESE}(\sigma)}

\def\disc{{\rm disc}}
\def\cusp{{\rm cusp}}
\def\spec{{\rm spec}}
\def\trace{{\rm trace}}
\def\spur{\mathfrak{Sp}}
\def\spurs{\spur_\sigma}
\def\Autom{\bs{\mathcal{A}}}
\def\Automd{\Autom}
\def\Base{\bs{\mathcal{B}}}

\def\fetu{\ESL(\sigma,\tu)}
\def\Fetu{\bs{\mathfrak{l}}(\sigma,\tu)}
\def\Lie{{\rm{Lie}\, }}
\def\st{{\rm{st}}}
\def\prim{{\rm{prim}}}

\def\bmMtunu{\bs{\mu}_{M,\tu}(\nu)}
\def\bmMtunuo{\bs{\mu}_{M,\tu}(\nu,0)}
\def\azx{a_{\tZ,\mu},}
\def\dmu{\dot\mu}

\def\oD{\Delta}
\def\hD{\widehat\Delta}
\def\vp{\varpi}
\def\V{V}
\def\com#1{\quad\hbox{#1}\quad}
\def\coma#1{\,\,\hbox{#1}\,\,}
\def\comm#1{\qquad\hbox{#1}\qquad}
\def\lb{\langle}
\def\rb{\rangle}

\chapter{Erratum}\label{erratum}

\def\citeLW#1{#1}

\newcor 
Il a ŽtŽ observŽ par P.-H. Chaudouard que la preuve de  1.8.1  pages~22-23
n'est valable que si les  sous-espaces $\ag_Q^R$ et $\ag_S^R$ sont orthogonaux.
Un argument diffŽrent est donnŽ en \ref{cpqrx} ci-dessous.
On commence par Žtablir un rŽsultat auxiliaire. 
Soit $V$ un espace vectoriel rŽel euclidien de dimension finie
On note $(\,,\,)$ le produit scalaire, $\|\,\|$ la norme et $V^*$ l'espace dual (que l'on peut identifier ˆ $V$). 
Soit $\oD$ une base de $V^*$ (que l'on ne suppose pas nŽcessairement obtuse) 
et soit $\oD^*$ la base de $V^*$ duale de $\Delta$. 
On note  $\alpha\mapsto\xi_\alpha$ la bijection naturelle entre $\oD$ et 
$\oD^*$. 
\begin{proposition} \label{ocd}\label{cor}
Soit $X\in V$. On note $C(\oD,X)$ l'ensemble des
$H\in V$ tels que pour tout $\alpha \in \oD$ on ait 
$$\alpha(H)\ge0 \com{et}
\xi_\alpha(H-X) \le0\ptf$$
Il existe une constante $c>0$ telle que
$$||H||\le c||X||\ptf$$
En particulier, l'ensemble $C(\oD,X)$ est compact et si $X=0$ alors $C(\oD,X)$ est rŽduit ˆ $\{0\}$.
\end{proposition}

\begin{proof}Pour $H\in C(\oD,X)$ on a
$$ ( H,H )=\sum_{\alpha\in\oD} \alpha(H) \xi_\alpha(H)\le \sum_{\alpha\in\oD} \alpha(H) \xi_\alpha(X)$$
ce qui implique l'existence d'une constante $c>0$ telle que
$$||H||^2\le c||H||\,||X||
\comm{et donc}||H||\le c||X||\ptf$$
\end{proof}

On considre maintenant trois sous-groupes paraboliques $P\subset Q\subset R$ 
dans $G$. Soit $V=\ag_P^R$. On dispose de la base $\oD_P^R$ de $V^*=\ag_P^{R,*}$.
On note encore $\alpha\mapsto\xi_\alpha$ la bijection naturelle entre $\oD_P^R$ et 
sa base duale\footnote{On prendra garde que cette base duale n'est pas en gŽnŽral 
la base $\hat\Delta_P^R$: leurs ŽlŽments sont colinŽaires, les coefficients de proportionnalitŽ Žtant positifs
mais non nŽcessairement Žgaux ˆ 1.  Plus gŽnŽralement, $\hat{\oD}_Q^R$ est une base du dual de $\ag_Q^R$ 
dont les ŽlŽments sont colinŽaires aux $\xi_\alpha$ pour $\alpha\in \Delta_P^S$ et 
$\hat{\Delta}_S^R=\hat{\Delta}_P^R \smallsetminus \hat{\Delta}_Q^R$ est une base du dual de $\ag_S^R$ dont 
les ŽlŽments sont colinŽaires aux $\xi_\alpha$ pour $\alpha \in \oD_P^Q$, les coefficients de proportionnalitŽ Žtant tous positifs. 
L'ensemble $C(P,Q,R,X)$ dŽfini ici est donc le mme (plus prŽcisŽment, c'est sa projection sur $V$) 
que celui de  1.8].}. D'aprs  1.2.2], il existe un sous-groupe parabolique $S$ tel 
que $\oD_P^R$ soit l'union disjointe  
de $\oD_P^Q$ et $\oD_P^S$. Soit $X\in\V$. On note $C(P,Q,R,X)$ l'ensemble des $H\in V$ tels que
$$\alpha(H)>0 \com{pour tout} \alpha\in\oD_P^Q\;,\quad
\alpha(H)\le0 \com{pour tout} \alpha\in\oD_P^S$$
et
$$\xi_\alpha(H-X)>0 \com{pour tout} \alpha\in\oD_P^S\;, \quad
\xi_\alpha(H-X)\le0 \com{pour tout} \alpha\in\oD_P^Q\ptf$$
Le lemme  1.8.1  affirme que

\begin{lemme} \label{cpqrx}
\begin{enumerate}
\item[(i)] L'ensemble $C(P,Q,R,X)$ est relativement compact
et  conte\-nu dans une boule de rayon $c||X||$ pour une constante $c>0$. 
\item[(ii)] $C(P,Q,R,0)$ est vide sauf si $P=R$.
\item[(iii)] Si $X$ est dans l'adhŽrence la chambre de Weyl positive, 
i.e. si $\alpha(X)\ge0$ pour tout $\alpha\in\oD_P^R$, alors $C(P,Q,R,X)$ est vide si $Q\subsetneq R$.
\end{enumerate}
\end{lemme}

\begin{proof}  
On observe que l'adhŽrence de $C(P,Q,R,X)$ est un fermŽ dans
$$C(\oD,X)\subset\V=\ag_P^R$$ o $\oD$ est la base de $V$ formŽe des $\alpha\in\oD_P^Q$
et des $-\alpha$ pour $\alpha\in\oD_P^S$.
L'assertion (i) rŽsulte de \ref{ocd}. Lorsque $X=0$ il rŽsulte de \ref{cor} que $C(P,Q,R,0)\subset\{0\}$ et
il y a des inŽgalitŽs strictes qui excluent $H=0$ si $P\ne R$. Pour Žtablir (iii)
 on observe que
$$\alpha(H)= \alpha(X)+
\sum_{\beta\in\oD_P^Q}(\alpha,\beta)\xi_\beta(H-X) 
+\sum_{\beta\in\oD_P^S}(\alpha,\beta)\xi_\beta(H-X)\ptf
$$
ConsidŽrons $H\in C(P,Q,R,X)$,
$\alpha\in\oD_P^S$ et $X$  dans l'adhŽrence la chambre de Weyl positive.
Comme $\oD_P^R$ est une base obtuse, on a
$$\alpha(H) \le0\;,\quad \alpha(X)\ge0\com{et}
\sum_{\beta\in\oD_P^Q}(\alpha,\beta)\xi_\beta(H-X)\ge0$$
ce qui implique 
$$ \sum_{\beta\in\oD_P^S}(\alpha,\beta)\langle \xi_\beta,H-X\rangle\le0\ptf$$
Pour toute famille de scalaires positifs $\{c_\alpha\,\vert\, \alpha \in \Delta_P^S\}$ on a donc
$$\sum_{\alpha\in\oD_P^S}(\eta,\alpha) \xi_\alpha(H-X)\le0\com{avec}
\eta=\sum_{\alpha\in\oD_P^S}c_\alpha\alpha\ptf$$
Les restrictions $\overline\xi_\alpha$ des $\xi_\alpha$ au sous-espace
$\ag_P^S\subset \ag_P^R$ forment la base duale  de $\oD_P^S$ qui est aig\"ue puisque 
$\oD_P^S$ est obtuse  1.2.5]. On  choisit
$$\eta=\sum_{\alpha\in\oD_P^S}\overline\xi_\alpha \com{\cad}
c_\alpha=\sum_{\beta\in\oD_P^S}(\overline\xi_\alpha,\overline\xi_\beta)>0\ptf$$
Dans ce cas
 $(\eta,\alpha)=1$ pour tout $\alpha\in\oD_P^S$, ce qui fournit l'inŽgalitŽ
 $$\sum_{\alpha\in\oD_P^S}\xi_\alpha(H-X)\le0\ptf$$
Comme $\xi_\alpha(H-X)>0$ pour tout $\alpha\in\oD_P^S$, 
c'est une contradiction si $\oD_P^S$ est non vide.  
\end{proof}

\newcor Dans l'ŽnoncŽ de  4.2.2  page 86 appara"t un signe $(-1)^{a_Q-a_G}$ qui est erronŽ:
 il faut supprimer ce signe. L'erreur est engendrŽe
 par une faute de frappe dans la preuve: dans la formule (1) en haut de la page 87 
 on doit lire $(-1)^{a_R-a_Q}$ au lieu de $(-1)^{a_R-a_G}$.

\newcor 
Il a ŽtŽ observŽ par Delorme que
dans la proposition  4.3.3  il faut ajouter une condition sur $T$. Non seulement
$T$ doit tre \og rŽgulier\fg mais il doit tre \og loin des murs\fg.
De fait, il est supposŽ implicitement que ${\bf d}_{P_0}(T)> 0$, mais il faut de plus fixer 
une constante $c$ (positive assez grande) et se limiter ˆ des $T$ vŽrifiant $$\| T \| \le c\, {\bf d}_{P_0} (T )$$ ou, ce qui 
revient au mme, imposer (pour une autre constante)
$$\langle\alpha,T \rangle \le c \langle \beta,T \rangle$$ 
pour tous les $\alpha$ et $\beta \in \Delta$. %
L'argument qui suit corrigeant la preuve de 4.3.3  est empruntŽ ˆ un message de 
Waldspurger ˆ Delorme.

\begin{proof}On reprend la preuve et les notations de la proposition  4.3.2]. 
On doit majorer $$|x|^B |A_{Q,R}(x)|$$ pour $Q \subsetneq R$. L'expression $A_{Q,R}(x)$ est donnŽe 
par une somme sur $\xi\in Q(F)\backslash G(F)$, mais pour $x$ fixŽ il y a au plus un $\xi$ modulo $Q(F)$ tel que
$$F^Q_{P_0}(\xi x, T)\sigma_Q^R(\bfH_0(\xi x -T)=1\ptf$$
Quitte ˆ remplacer $x$ par $\xi x$, on peut supposer
$$F_{P_0} ^Q (x, T ) = 1\vg \leqno{(1)}$$
$$\sigma_Q^R(\bfH_0(x) - T) = 1\ptf\leqno{(2)}$$
D'aprs  3.6.2], on peut aussi supposer $x = n_1ak$ avec $n_1 \in N_Q(\AM)$, 
$a \in \AA^G_0$ et $k \in \Omega$ un compact 
(qui ne dŽpend que de du choix du sous-groupe parabolique $Q$). 
Pour les majorations, on peut aussi bien supposer $k = 1$. On impose aussi que $n_1$ est dans un compact, 
ce que l'on peut assurer en le multipliant ˆ gauche par un ŽlŽment de $N_Q(F)$. 
D'aprs  4.3.1], pour tout entier $r \ge 1$, il existe un opŽrateur diffŽrentiel $X_{Q,R}$ de degrŽ $r$, qui ne
dŽpend que de $Q$, $R$, $r$, tel que
$$|A_{Q,R}(x)| \le \sup_{n'_1\in N_Q(\AM)} |\bs{\rho}({\rm Ad}(a\mun)X_{Q,R})\varphi(n'_1a)|\ptf$$
L'opŽrateur $X_{Q,R}$ est de la forme $X_{Q,R}= \sum_m Z_m$ avec
$${\rm Ad}(e^{H_1})Z_m = e^{\lambda_m (H_1)} Z_m \quad \hbox{pour} \quad H_1\in \ag_Q^R$$
o $\lambda_m$ est une combinaison linŽaire ˆ coefficients entiers $\geq r$ des racines dans $\Delta_Q^R$.
On en dŽduit bien une majoration
$$|A_{Q,R}(x)|\ll \sup_{n'_1,m}e^{-r \sum_{\alpha\in \Delta_Q^R} \langle \alpha,\bfH_0(a)\rangle }
|\bs{\rho}(Z_m)\varphi(n'_1a)|$$
o les $Z_m$ ne dŽpendent que de $Q$, $R$, $r$. 
D'o 
$$\qquad |x|^B |A_{Q,R}(x)| \ll e^{C(a)}\sup_{n'_1,m}|\bs{\rho}(Z_m)\varphi(n'_1a)||x|^{-A}\leqno(3)$$ 
\pni
o on a posŽ
$$C(a)=c(B+A) \| H_0(a) \| -r \sum_{\alpha\in \Delta_Q^R} \langle \alpha,H_0(a)\rangle\vg$$
$c$ Žtant une constante telle que $|a| \le e^{c\| H_0(a)\| }$. La proposition 4.3.3 rŽsulte de (3)
ˆ condition de prouver que 
$$C(a) \le - A\cdot{\bf d}_{P_0}(T)\ptf$$
Ecrivons
$$H_0(a)=H_Q +H_Q^R +H_R$$
avec $H_Q \in \AA_Q$, $H_Q^R \in \AA^R_Q$ et $H_R \in \AA_R$. La condition (1) impose $\| H_Q\| \ll \| T\| $. La
condition (2) impose (d'aprs  2.11.6]) que $\| H_R - T_R\| \ll \| H_Q^R - T_Q^R\| $, 
avec des dŽfinitions Žvidentes de $T_R$ et $T_Q^R$. D'o $\| H_R\| \ll \| H_Q \|+\| T \|$. Elle impose
aussi que $\langle\alpha,H_Q^R -T_Q \rangle>0$ pour tout $\alpha\in \Delta_Q^R $. On en dŽduit
$$C(a) \le c_1(A+B)\| H_Q^R \| +c_2(A+B)\| T \| -\frac{r}{2}
\sum_{\alpha\in \Delta_Q^R} \langle \alpha, H_Q+T_Q \rangle$$
o $c_1$, $c_2$ sont des constantes positives qui ne dŽpendent de rien. 
En choisissant $r \ge c_3(A+B)$ pour une constante $c_3$ convenable, on a
$$c_1(A+B)\| H_Q \| - \sum_{\alpha\in \Delta_Q^R} <\langle\alpha, H_Q \rangle\le 0\ptf$$
Pour $\alpha \in \Delta^R_Q$, soit $\beta \in \Delta$ se projetant sur $\alpha$. On sait que $\alpha$ est la somme de 
$\beta$ et d'une combinaison linŽaire des $\gamma \in \Delta_Q$ 
ˆ coefficients positifs ou nuls. Donc $$\langle \alpha,T_Q \rangle\ge \langle \beta,T \rangle\ge {\bf d}_{P_0}(T)\ptf$$
 Alors $$C(a) \le c_2(A+B)\| T\| -|\Delta^R_Q|\frac{r}{2}{\bf d}_{P_0}(T)$$
Si on impose, comme dit plus haut, une condition $\|T\| \le c\,{\bf d}_{P_0}(T)$, alors pour 
$r \ge c4(2A+B)$, on a bien
$$C(a) \le - A\cdot{\bf d}_{P_0}(T)\ptf$$
\end{proof} 

\newcor La formule pour le produit scalaire $\langle \Phi,\Psi \rangle_P$,
en bas de la page 93, donnŽe par une intŽgrale sur 
${\bf X}_P=\AA_PP(F)N_P(\AM)\backslash G(\AM)$ n'a pas de sens si $P\neq G$. En effet le groupe 
$\AA_PP(F)N_P(\AM)$ n'est pas unimodulaire et il n'y a pas de mesure $\G(\AM)$-invariante ˆ droite sur 
${\bf X}_{P}$. 
Il convient de poser
$$\langle \Phi,\Psi \rangle_P
= \int_{\bf K}\int_{{\bf X}_M} \Phi(mk) \overline{\Psi(mk)} \dd m \dd k\ptf$$
En revanche sur ${\bf X}_{P,G}= \AA_G P(F)N_P(\AM)\backslash G(\AM)$ on dispose de la mesure quotient
et on a la formule d'intŽgration
$$\int_{{\bf X}_{P,G}}e^{\langle 2\rho_P,\bfH_P(x) \rangle} \phi(x)\dd x 
= \int_{\bf K}\int_{{\bf X}_{M,G}} \phi(mk) \dd m \dd k$$
 qui est utilisŽe dans la preuve du thŽorme  5.4.2]. 

\newcor Page 104.
Dans la dŽfinition de l'opŽrateur $\bs\rho_{P,\sigma,\mu}(\delta,y,\omega)$
il manque le caractre $\omega$ dans la dŽfinition de l'espace d'arrivŽe.

\newcor Page 106, ligne 16 il faut lire
$$K_\tG(f,\omega;x,y)= \sum_i \int_{{\bf X}_G}K_\tG(g_i,\omega;x,z) 
K_G^*(\omega h_i;z,y)\dd z \ptf$$

\newcor Page 111.
Dans l'ŽnoncŽ de la proposition  7.2.2 
il faut prendre pour $K_1$ et $K_2$ des noyaux sur $\bsX_{\theta(P),G}\times\bsX_{P,G}$ et 
sur ${\bf X}_{P,G}\times{\bf X}_{P,G}$ et remplacer les intŽgrales sur ${\bf X}_P$ par des intŽgrales sur ${\bf X}_{P,G}$. 

\newcor Page 112. Il faut remplacer $\mathcal{B}^P(\sigma)$ par $\mathcal{B}^S(\sigma)$ dans 
la dŽfinition de $H_\sigma$: ici $M$ est un sous-groupe de Levi de $S$ avec $S\subset P$. 
Cette erreur due ˆ un copier-coller intempestif s'est propagŽe dans toute la suite.

\newcor Page 126, ligne 2. Il faut remplacer $\mathfrak{n}(\AM)$ par $\mathfrak{n}_P(\AM)$ o $\mathfrak{n}_P$ est l'algbre de Lie 
de $N=N_P$ (rappel: on a notŽ $\mathfrak{n}$ l'algbre de Lie de $N_{Q^+}$).

\newcor Page 128, ligne -15. La dŽfinition de $j_{\tP\!,\oo}(x)$ n'a pas de sens. 
Il faut la remplacer par 
$$j_{\tP\!,\oo}(x)= 
\sum_{\delta \in \oo \cap \tM_P(F)} \int_{N(\delta,\AM)}\omega(x)f^1(x^{-1}\delta n x) \dd n\ptf$$
De mme, dans l'ŽnoncŽ du lemme 9.2.2 page 129 il faut remplacer la premire expression par
$$\int_{N(F)\backslash N(\AM)} \sum_{\delta \in \oo\cap \tM_P(F)} 
\int_{N(\delta,\AM)} \phi(n^{-1} \delta n_1 n)\dd n_1 \dd n \ptf$$

\newcor Page 128, ligne -13. Lire:
{\it la partie quasi semi-simple $\delta_s$ de $\delta$.}

\newcor Page 129. Dans la preuve du lemme 9.2.2, il convient de dire (mme si c'est implicite) 
que l'application du lemme  9.2.1 induit une application surjective sur les groupes des points adŽliques.

\newcor Page 130, ligne -10. Lire $s(\delta)$ au lieu de $\delta_s$.

\newcor Page 130, ligne -11. Lire $I_{s(\delta)}(F)$ au lieu de $I_{\delta_s}$.

\newcor Page 130, ligne -6. Pour une classe de conjugaison 
quasi semi-simple $\cc$ 
on doit prŽciser que $k_\cc^T$ est 
dŽfinie en remplaant, dans la dŽfinition de $k_{\oo}^T$  9.2], la somme sur $P\supset P_0$
 par la somme sur les $P\supset P_0$ tels que $\tM_P$ 
contienne un conjuguŽ de $\tM_\delta$ et l'ensemble $\oo$ par l'orbite $\cc$.

\newcor Page 150. Dans la formule pour $\hat{\tau}_\tP(H-X)$ 
il faut sommer sur les $\tQ$ tels que $\tP\subset\tQ$, i.e. prendre $\tR=\tG$.

\newcor Page 168, ligne 2. Lire: {\it et sera ˆ l'\oe uvre...}

\newcor Page 182. Dans l'ŽnoncŽ du point $(ii)$ du lemme 13.1.1, il faut lire: {\it Pour tout $\lambda \in \ag_0^*$, etc}.

\newcor Pages 183--186: l'indice $Q_0$ devrait tre en exposant dans les termes ˆ droite des Žquations (2) 
page 183, (7) page 185 et (8) page 186. 

\newcor Page 194, ligne 11, il faut lire: {\it holomorphe en $\lambda$ et $\nu$}.

\newcor Page 197, lignes 2 et 6: les $+Y$ devraient tre des $-Y$.

\newcor Page 203. Dans l'inŽgalitŽ ligne~12 et dans celle ligne~-3, il faut remplacer $Y_{P'\!,e}(T)$ par $-Y_{P'\!,e}(T)$.

\newcor Page 203, ligne -9. L'inŽgalitŽ est stricte: il faut lire $\| U+V \| > \rho \| T \|$.

\newcor Page 205, majoration (1): il faut lire $\vert A^T_{s,t} - E^T_1\vert$.
 
\newcor Page 207: $\wt{Y}_{S'}(X)$ est la projection de $u^{-1}X$ sur $t(\ag_S)^Q$.

\newcor Page 219, ligne -6: l'expression pour $\bfA_M^{Q_1}(T,\sigma,\tu,\mu)$ faisant intervenir une trace
n'a de sens que si l'opŽrateur est un endomorphisme de l'espace ${\mathcal A}({\bf X}_S,\sigma)$. 
La formulation correcte est celle du point ($i$) du lemme  14.1.4].

\end{document}